%% file: Paper.tex
\documentclass[10pt,reqno]{NumPDEsArticle}

\usepackage[basicdelimiters]{NumPDEsMacros}

\usepackage[english]{babel}
\usepackage{csquotes}
\usepackage{enumitem}
\usepackage{amssymb}
\usepackage{algorithm}
\usepackage{graphicx}
\usepackage{subcaption}
\usepackage{pgfplots}
\pgfplotsset{compat=1.9}
\usepackage[export]{adjustbox}%
\usepackage{pgfplotstable}
\usepackage{orcidlink}
\usepackage{todonotes}

\usepackage{tikz}
\usetikzlibrary{%
    calc,%
    hobby,%
    matrix,%
    shapes,%
    arrows,%
    fadings,%
    fpu,%
    cd,%
    spy,%
    3d,%
    positioning,%
    shapes.multipart,%
    shadows,%
    shapes.symbols%
}
\usepackage{pgfplots}
\usepackage{pgfplotstable}
\usepackage{tikz-3dplot}
\pgfplotsset{%
    compat=newest,%
    every axis/.style={scale only axis},%
    grid style={densely dotted, semithick},%
}
\definecolor{TUblue}{rgb}{0,0.4,0.6}
\definecolor{TUgray}{rgb}{0.3922,0.3882,0.3882}
\definecolor{TUgreen}{rgb}{0,0.4941,0.4431}
\definecolor{TUmagenta}{rgb}{0.7294,0.2745,0.5098}
\definecolor{TUyellow}{rgb}{0.8824,0.5373,0.1333}
\newcommand\drawslopetriangle[4][ST]{
    \pgfplotsextra
    {
        \pgfkeys{/pgf/fpu=true}
        \pgfmathsetmacro\leftcoord{#3}
        \pgfmathsetmacro\rightcoord{10*#3}
        \pgfmathsetmacro\bottomcoord{#4}
        \pgfmathsetmacro\topcoord{10^(#2)*#4}
        \pgfkeys{/pgf/fpu=false}
        \coordinate (#1-BL) at (axis cs:\leftcoord,\bottomcoord);
        \coordinate (#1-BR) at (axis cs:\rightcoord,\bottomcoord);
        \coordinate (#1-TL) at (axis cs:\leftcoord,\topcoord);
        \shadedraw[%
            bottom color = black!20,%
            middle color = black!5,%
            top color    = white,%
            draw         = black,%
            font         = \scriptsize%
        ]
        (#1-TL) -- (#1-BL) node[midway, left] {\(#2\)} -- (#1-BR) node[midway, below] {\(1\)} -- (#1-TL);
    }
}
\newcommand\drawswappedslopetriangle[4][SST]{
    \pgfplotsextra
    {
        \pgfkeys{/pgf/fpu=true}
        \pgfmathsetmacro\leftcoord{#3/10}
        \pgfmathsetmacro\rightcoord{#3}
        \pgfmathsetmacro\topcoord{#4}
        \pgfmathsetmacro\bottomcoord{10^(-#2)*#4}
        \pgfkeys{/pgf/fpu=false}
        \coordinate (#1-TR) at (axis cs:\rightcoord,\topcoord);
        \coordinate (#1-BR) at (axis cs:\rightcoord,\bottomcoord);
        \coordinate (#1-TL) at (axis cs:\leftcoord,\topcoord);
        \shadedraw[%
            bottom color = black!20,%
            middle color = black!5,%
            top color    = white,%
            draw         = black,%
            font         = \scriptsize%
        ]
        (#1-BR) -- (#1-TR) node[midway, right] {\(#2\)} -- (#1-TL) node[midway, above] {\(1\)} -- (#1-BR);
    }
}
\newcommand\drawslopetriangleup[4][STU]{
    \pgfplotsextra
    {
        \pgfkeys{/pgf/fpu=true}
        \pgfmathsetmacro\leftcoord{#3}
        \pgfmathsetmacro\rightcoord{10*#3}
        \pgfmathsetmacro\bottomcoord{#4}
        \pgfmathsetmacro\topcoord{10^(#2)*#4}
        \pgfkeys{/pgf/fpu=false}
        \coordinate (#1-BL) at (axis cs:\leftcoord,\bottomcoord);
        \coordinate (#1-BR) at (axis cs:\rightcoord,\bottomcoord);
        \coordinate (#1-TR) at (axis cs:\rightcoord,\topcoord);
        \shadedraw[%
            bottom color = black!20,%
            middle color = black!5,%
            top color    = white,%
            draw         = black%
        ]
        (#1-BL) -- (#1-BR) node[midway, below] {\(1\)} -- (#1-TR) node[midway, right] {\(#2\)} -- (#1-BL);
    }
}

\pgfplotsset{
    colormap={parula}{
        rgb=(0.2081,0.1663,0.5292)
        rgb=(0.2116,0.1898,0.5777)
        rgb=(0.2123,0.2138,0.627)
        rgb=(0.2081,0.2386,0.6771)
        rgb=(0.1959,0.2645,0.7279)
        rgb=(0.1707,0.2919,0.7792)
        rgb=(0.1253,0.3242,0.8303)
        rgb=(0.0591,0.3598,0.8683)
        rgb=(0.0117,0.3875,0.882)
        rgb=(0.006,0.4086,0.8828)
        rgb=(0.0165,0.4266,0.8786)
        rgb=(0.0329,0.443,0.872)
        rgb=(0.0498,0.4586,0.8641)
        rgb=(0.0629,0.4737,0.8554)
        rgb=(0.0723,0.4887,0.8467)
        rgb=(0.0779,0.504,0.8384)
        rgb=(0.0793,0.52,0.8312)
        rgb=(0.0749,0.5375,0.8263)
        rgb=(0.0641,0.557,0.824)
        rgb=(0.0488,0.5772,0.8228)
        rgb=(0.0343,0.5966,0.8199)
        rgb=(0.0265,0.6137,0.8135)
        rgb=(0.0239,0.6287,0.8038)
        rgb=(0.0231,0.6418,0.7913)
        rgb=(0.0228,0.6535,0.7768)
        rgb=(0.0267,0.6642,0.7607)
        rgb=(0.0384,0.6743,0.7436)
        rgb=(0.059,0.6838,0.7254)
        rgb=(0.0843,0.6928,0.7062)
        rgb=(0.1133,0.7015,0.6859)
        rgb=(0.1453,0.7098,0.6646)
        rgb=(0.1801,0.7177,0.6424)
        rgb=(0.2178,0.725,0.6193)
        rgb=(0.2586,0.7317,0.5954)
        rgb=(0.3022,0.7376,0.5712)
        rgb=(0.3482,0.7424,0.5473)
        rgb=(0.3953,0.7459,0.5244)
        rgb=(0.442,0.7481,0.5033)
        rgb=(0.4871,0.7491,0.484)
        rgb=(0.53,0.7491,0.4661)
        rgb=(0.5709,0.7485,0.4494)
        rgb=(0.6099,0.7473,0.4337)
        rgb=(0.6473,0.7456,0.4188)
        rgb=(0.6834,0.7435,0.4044)
        rgb=(0.7184,0.7411,0.3905)
        rgb=(0.7525,0.7384,0.3768)
        rgb=(0.7858,0.7356,0.3633)
        rgb=(0.8185,0.7327,0.3498)
        rgb=(0.8507,0.7299,0.336)
        rgb=(0.8824,0.7274,0.3217)
        rgb=(0.9139,0.7258,0.3063)
        rgb=(0.945,0.7261,0.2886)
        rgb=(0.9739,0.7314,0.2666)
        rgb=(0.9938,0.7455,0.2403)
        rgb=(0.999,0.7653,0.2164)
        rgb=(0.9955,0.7861,0.1967)
        rgb=(0.988,0.8066,0.1794)
        rgb=(0.9789,0.8271,0.1633)
        rgb=(0.9697,0.8481,0.1475)
        rgb=(0.9626,0.8705,0.1309)
        rgb=(0.9589,0.8949,0.1132)
        rgb=(0.9598,0.9218,0.0948)
        rgb=(0.9661,0.9514,0.0755)
        rgb=(0.9763,0.9831,0.0538)
    }
}

\definecolor{col1}{HTML}{332288}
\definecolor{col2}{HTML}{88CCEE}
\definecolor{col3}{HTML}{44AA99}
\definecolor{col4}{HTML}{117733}
\definecolor{col5}{HTML}{999933}
\definecolor{col6}{HTML}{DDCC77}
\definecolor{col7}{HTML}{CC6677}
\definecolor{col8}{HTML}{882255}
\definecolor{col9}{HTML}{AA4499}
\pgfplotsset{%
	degdefault/.style = {%
		mark = *,%
		mark size = 2pt,%
		every mark/.append style = {solid},%
		gray,%
		every mark/.append style = {fill = gray!60!white}%
	},%
	marker1/.style = {%
		degdefault,%
		mark size = 1.6pt,%
		col1,%
		every mark/.append style = {fill = col1!60!white}%
	},%
	marker2/.style = {%
		degdefault,%
		mark = pentagon*,%
		mark size = 1.8pt,%
		every mark/.append style = {rotate = 180},
		col3,%
		every mark/.append style = {fill = col3!60!white}%
	},%
	marker3/.style = {%
		degdefault,%
		mark = triangle*,%
		mark size = 1.8pt,%
		col4,%
		every mark/.append style = {fill = col4!60!white}%
	},%
	marker4/.style = {%
		degdefault,%
		mark = halfdiamond*,%
		every mark/.append style = {rotate = 90},%
		mark size = 2pt,%
		col5,%
		every mark/.append style = {fill = col5!60!white}%
	},%
	marker5/.style = {%
		degdefault,%
		mark = halfsquare*,%
		every mark/.append style = {rotate = 135},%
		mark size = 2pt,%
		col7,%
		every mark/.append style = {fill = col7!60!white}%
	},%
	marker6/.style = {%
		degdefault,%
		mark = halfsquare*,%
		every mark/.append style = {rotate = 315},%
		mark size = 2pt,%
		col8,%
		every mark/.append style = {fill = col8!60!white}%
	},%
	marker7/.style = {%
		degdefault,%
		mark = halfdiamond*,%
		every mark/.append style = {rotate = 180},%
		mark size = 1.8pt,%
		col9,%
		every mark/.append style = {fill = col8!60!white}%
	},%
	uniform/.style = {%
		dashed,%
		every mark/.append style = {%
			black!50!white,%
			fill = black!20!white
		}%
	},%
	adaptive/.style = {%
		solid%
	},%
}

\pgfplotstableset{%
	create on use/estimator/.style={%
		create col/expr={\thisrow{eta} + \thisrow{mu}}
	},
	create on use/cumulativeNdof/.style={%
		create col/expr={%
				\pgfmathaccuma + \thisrow{ndof}
		}
	}
}

\newcommand{\exact}{\star}

\newcommand{\vvvert}{|\mkern-1.5mu|\mkern-1.5mu|}

\renewcommand\D{{\textup{D}}}

\DeclareMathOperator{\sign}{sign}

\makeatletter
\newcommand{\labeltext}[2]{%
    \@bsphack
    \csname phantomsection\endcsname 
    \def\@currentlabel{#1}{\label{#2}}%
    \@esphack
}
\makeatother

\newcommand\subfigref[1]{(\textsc{\subref{#1}})}

\title[Adaptive nonlinear LSFEM]{Global convergence of adaptive least-squares finite element methods for nonlinear PDEs}
\address{TU Wien, Institute of Analysis and Scientific Computing, Wiedner Hauptstr. 8--10/E101/4, 1040 Vienna, Austria}
\author{Philipp Bringmann~\orcidlink{0000-0002-4546-5165}}
\email{philipp.bringmann@asc.tuwien.ac.at (corresponding author)}
\author{Dirk Praetorius~\orcidlink{0000-0002-1977-9830}}
\email{dirk.praetorius@asc.tuwien.ac.at}
\thanks{%
	This research was funded in whole or in part by the Austrian Science Fund (FWF)
	[\href{https://www.fwf.ac.at/en/research-radar/10.55776/I6802}{10.55776/I6802},
	\href{https://www.fwf.ac.at/en/research-radar/10.55776/P33216}{10.55776/P33216}, and
	\href{https://www.fwf.ac.at/en/research-radar/10.55776/PAT3699424}{10.55776/PAT3699424}].
}

\begin{document}

\maketitle
\thispagestyle{fancy}

\begin{abstract}
The Zarantonello fixed-point iteration is an established linearization scheme for quasilinear PDEs with strongly monotone and Lipschitz continuous nonlinearity in Hilbert spaces. This paper presents a weighted least-squares minimization for the computation of the update of this scheme. The resulting formulation allows for a conforming least-squares finite element discretization of the primal and dual variable of the PDE with arbitrary polynomial degree. The least-squares functional provides a built-in a posteriori discretization error estimator in each linearization step motivating an adaptive Uzawa-type algorithm with an outer linearization loop and an inner adaptive mesh-refinement loop. For quasilinear PDEs in divergence form satisfying a 2-growth condition, we prove global R-linear convergence of the computed linearization iterates for arbitrary initial guesses. Particular focus is on the role of the weights in the least-squares functional of the linearized problem and their influence on the robustness of the Zarantonello damping parameter. Numerical experiments illustrate the performance of the proposed algorithm.
\end{abstract}

%
%
\section{Introduction}
\label{sec:introduction}

\subsection{Motivation}
Least-squares methods have enjoyed ongoing attention in the numerical solution of partial differential equations (PDEs) for several decades.
Besides the formulation-inherent well-posedness, this is primarily due to their built-in a~posteriori error estimation which directly enables the application in adaptive mesh-refinement algorithms;
see \cite{Bringmann2024} for a recent literature review on adaptive least-squares finite element methods (LSFEMs)
and \cite{Bringmann2023} for the convergence analysis with rates of adaptive LSFEMs for linear problems.
Moreover, their intrinsic symmetrization and stabilization motivated the application to space-time formulations of parabolic and hyperbolic PDEs;
see, e.g.,
\cite{FuhrerKarkulik2021,GantnerStevenson2021,GantnerStevenson2024_rates,FuhrerGonzalezKarkulik2025,HoonhoutLoscherSteinbachUrzuaTorres2026,KotheSteinbach2026,KotheLoscherSteinbach2023}.
Further advantages include the versatile weak enforcement of boundary conditions \cite{MonsuurSmeetsStevenson2025} as well as
the flexible choice of the discretization encouraging the use of least-squares cost functionals in the context of physics-informed neural networks
\cite{RaissiPerdikarisKarniadakis2019,CaiChenLiuLiu2020,MeissnerHuynhKuberryBochev2025}.
The equal-order approximation of primal and dual variable
is particularly attractive for applications in computational mechanics; 
see, e.g., \cite{MullerStarkeSchwarzSchroder2014}.

\subsection{Literature}
\label{sec:intro:literature}
For nonlinear PDEs, however, least-squares formulations are less prevalent in the literature.
The main reason is the possible lack of convexity of the least-squares functional which consists of the sum
of the nonlinear residuals of the (first-order system of) PDEs in squared Lebesgue norms.
Nevertheless, least-squares approaches have been successfully applied to a wide range of applications, e.g.,
various formulations of the Navier--Stokes equations 
\cite{BochevGunzburger1993,BochevCaiManteuffelMcCormick1998},
the nonlinear Stokes equation \cite{MonneslandLeeGunzburgerYoon2016},
the geometrically nonlinear elasticity problem \cite{ManteuffelMcCormickSchmidtWestphal2006},
the hyperelasticity problem \cite{MullerStarkeSchwarzSchroder2014},
sea-ice models \cite{BertrandSchneider2024},
and the Monge--Amp\`ere equation \cite{Westphal2019,BrennerSungTanZhang2024}.
These references follow different approaches.
The methods in \cite{BochevGunzburger1993,BochevCaiManteuffelMcCormick1998}
employ least-squares minimization of the nonlinear residuals to be solved with Newton's method,
but this approach is tailored to the Navier--Stokes equations;
see also the discussion in~\cite[Section~8.4]{BochevGunzburger2009}.
In most of the cases, the formulations
are based on a Gauss--Newton method which first linearizes the residuals with the Newton method and
then applies a least-squares minimization to compute the update direction.
This can be interpreted as an inexact Newton method and the recent work \cite{BertrandBrodbeckRickenSchneider2025}
applied the local convergence theory for such schemes to Gauss--Newton least-squares methods for nonlinear PDEs.
The relation between linearization and minimization in LSFEMs for nonlinear problems
is discussed in~\cite{PayetteReddy2011}.
The least-squares functional may also be used as an error estimator and refinement indicator for other discretizations
of nonlinear PDEs \cite{LiZhang2025}.

The discontinuous Petrov--Galerkin method (DPG) is a minimal residual method for primal, dual, or ultraweak variational formulations.
The flexibility of their broken test spaces leads to improved stability properties.
A DPG method for a quasi-linear model problem has been presented in \cite{CarstensenBringmannHellwigWriggers2018}.
It aims to minimize the residual of the nonlinear variational formulation.
The authors employ the close relation to least-squares methods to prove the existence of discrete minimizers while the uniqueness remained open.
Nevertheless, a sufficient a~posteriori criterion for the uniqueness is given in \cite[Theorem~4.4]{CarstensenBringmannHellwigWriggers2018}.
The reader is referred to \cite[Section~3.1]{CarstensenBringmannHellwigWriggers2018} and \cite[Section~4]{BringmannCarstensenTran2022} for a discussion of the problems of nonlinear residual minimization methods.
An alternative approach from \cite{CantinHeuer2018} establishes a DPG method based on a first-order formulation
where the nonlinear residual is posed as a side constraint to the minimization of the remaining linear residual.

In this paper, we will assume a 2-growth condition and the divergence form of the PDE.
In particular, we restrict to a Hilbert space setting excluding, e.g., the \(p\)-Laplace problem for \(p \neq 2\).
Beyond this setting, several publications investigate minimal residual methods in Banach spaces; see, e.g.,
\cite{Guermond2004,MugaVanDerZee2020,HoustonRoggendorfVanDerZee2022,LiDemkowicz2022}.
Relaxed Ka\v{c}anov schemes from \cite{DieningFornasierTomasiWank2020,BalciDieningStorn2023} 
enable the efficient solution of \(p\)-Laplace problems as employed in \cite{Storn2024}
to solve linear problems in \(W^{-1,p'}\) spaces for large \(p\).
For minimal residual methods for problems in nondivergence form, we refer, e.g., to \cite{Gallistl2017,QiuZhang2020,Fuhrer2021}
and the already mentioned references \cite{Westphal2019,BrennerSungTanZhang2024} for the Monge-Amp\`ere equation.
A novel minimal residual method in \(L^p(\Omega)\) norms has been recently introduced in \cite{GallistlTran2025}
for a class of fully nonlinear PDEs.

\subsection{Approach in this paper}
\label{sec:intro:approach}
We consider the model problem of strongly monotone and Lipschitz continuous quasilinear PDEs in divergence form
on a bounded polyhedral Lipschitz domain \(\Omega \subset \R^d\) with \(d \in \N\).
For general right-hand sides \(f_1 - \div f_2 \in H^{-1}(\Omega)\) 
with \(f_1 \in L^2(\Omega)\) and \(f_2 \in L^2(\Omega; \R^d)\),
it seeks \(u^\star \in H^1_0(\Omega)\) satisfying
\begin{equation}
	\label{eq:intro:nonlinear_primal_pde}
	- \div (\sigma(\nabla u^\star))
	=
	f_1
	- 
	\div(f_2)
	\quad \text{in } \Omega.
\end{equation}
The reader is referred to Section~\ref{sec:model_problem} for the detailed assumptions on the nonlinear mapping \(\sigma\colon \R^d \to \R^d\).
Our goal is to analyze an adaptive LSFEM for nonlinear PDEs with guaranteed global convergence,
i.e., for arbitrary initial guess on an arbitrary (and potentially coarse) initial mesh. 
The derivation of the discrete linearized problem can be summarized in the following three key steps:
\begin{enumerate}[label=(S\arabic*)]
	\item
		\label{step:fos}
		Rewriting the PDE as a nonlinear first-order system.
	\item
		\label{step:zarantonello}
		Application of the Zarantonello iteration to obtain a linearized first-order system of PDEs.
	\item
		\label{step:lsfem}
		Discrete solution of the linearized system with an LSFEM.
\end{enumerate}

As the first step~\ref{step:fos}, the nonlinear PDE~\eqref{eq:intro:nonlinear_primal_pde} is reformulated into a
first-order system with the solution \((p^\star, u^\star) \in H(\div, \Omega) \times H^1_0(\Omega)\) to
\begin{equation}
	\label{eq:intro:nonlinear_pde}
	- \div p^\star
	=
	f_1
	\quad\text{and}\quad
	p^\star - \sigma(\nabla u^\star)
	=
	- f_2
	\quad \text{in } \Omega.
\end{equation}

In order to linearize this system of PDEs in step~\ref{step:zarantonello},
we employ the Zarantonello fixed-point iteration \cite{Zarantonello1960} 
instead of the Newton method used in \cite{BertrandBrodbeckRickenSchneider2025}.
The Zarantonello iteration is also employed
for the iterative solution of nonlinear finite element discretizations
in the context of adaptive iterative linearized FEMs
\cite{CongreveWihler2017,HeidWihler2020_convergence,HeidWihler2020_linearization,HeidPraetoriusWihler2021,HaberlPraetoriusSchimankoVohralik2021,EggerEngertsbergerDomenigRoppertKaltenbacher2025}
as well as for the iterative symmetrization of nonsymmetric problems
\cite{BrunnerInnerbergerMiraciPraetoriusStreitbergerHeid2024}.
Given some previous iterate \((p^{k-1}, u^{k-1}) \in H(\div, \Omega) \times H^1_0(\Omega)\),
a damping parameter \(\delta > 0\), and positive weights \(\omega_1, \omega_2 > 0\), 
the Zarantonello linearization \cite{Zarantonello1960} of the nonlinear problem~\eqref{eq:intro:nonlinear_pde} seeks
the solution \((p^k_\star, u^k_\star) \in H(\div, \Omega) \times H^1_0(\Omega)\) to
\begin{equation}
	\label{eq:intro:zarantonello_fos}
	\begin{aligned}
		- \omega_1 \div p^k_\star
		&=
		- \omega_1 \div p^{k-1} + \delta \omega_1 \, \big[ f_1 + \div p^{k-1} \big],
		\\
		p^k_\star - \omega_2^2 \, \nabla u^k_\star
		&=
		p^{k-1} - \omega_2^2 \, \nabla u^{k-1} - \delta \, \big[ f_2 + p^{k-1} - \sigma(\nabla u^{k-1}) \big].
	\end{aligned}
\end{equation}

For the practical solution of the \(k\)-th iterates in step~\ref{step:lsfem},
this paper investigates the approximate solution of the linear first-order system~\eqref{eq:intro:zarantonello_fos} of PDEs
using a weighted least-squares approach with the functional 
\(Z_k(f_1, f_2; \,\cdot\,,\,\cdot\,)\colon H(\div, \Omega) \times H^1_0(\Omega) \to \R\) defined by
\begin{align*}
	Z_k(f_1, f_2; p, u)
	&\coloneqq
	\omega_1^2 C_{\textup{F}}^2 \,
	\Vert \div (p - p^{k-1}) + \delta \, [f_1 + \div p^{k-1}] \Vert_{L^2(\Omega)}^2
	\\
	&\phantom{{}\coloneqq{}}
	+
	\Vert p - p^{k-1} - \omega_2^2 \, \nabla (u - u^{k-1}) 
	+ \delta \, [f_2 + p^{k-1} - \sigma(\nabla u^{k-1})] \Vert_{L^2(\Omega)}^2.
\end{align*}
Therein, the Friedrichs constant \(C_{\textup{F}} > 0\) ensures the robustness with respect to the size of the domain \(\Omega\).
The solution of the PDE system~\eqref{eq:intro:zarantonello_fos} is formulated as the minimization problem 
of finding \((p^k_\star, u^k_\star) \in H(\div, \Omega) \times H^1_0(\Omega)\) such that
\begin{equation}
	\label{eq:intro:minimization}
	Z_k(f_1, f_2; p^k_\star, u^k_\star)
	=
	\min_{(p, u) \in H(\div, \Omega) \times H^1_0(\Omega)}
	Z_k(f_1, f_2; p, u).
\end{equation}
The damping parameter \(\delta > 0\) and the weights \(\omega_1, \omega_2 > 0\) 
can be chosen in terms of the monotonicity and the Lipschitz constant of the primal PDE~\eqref{eq:intro:nonlinear_primal_pde}
to guarantee a well-posed and convergent Zarantonello iteration.
See Sections~\ref{sec:Zarantonello_least_squares}--\ref{sec:alternative_weightings} 
for an analysis and discussion of different sufficient choices.
The minimization over finite element spaces will eventually lead to computable discrete iterates.
This completes step~\ref{step:lsfem}.

We highlight that the residual in~\eqref{eq:intro:nonlinear_pde} leads to the direct inclusion of general right-hand sides
\(f_1 - \div f_2 \in H^{-1}(\Omega)\) into the least-squares minimization~\eqref{eq:intro:minimization}.
An alternative treatment of rough right-hand sides in the context of minimal residual methods
consists of the application of regularizing operators as introduced in \cite[Section~3]{FuhrerHeuerKarkulik2022} for the lowest-order case and
in~\cite[Section~4.3]{DieningStornTscherpel2023} for arbitrary polynomial degree.

The least-squares functional \(Z_k(f_1, f_2;\,\cdot\,,\,\cdot\,)\) provides a built-in a~posteriori estimator
for the discretization error of the linearized problem.
This motivates its application in an adaptive mesh-refinement algorithm for determining the approximate Zarantonello update
resulting in an adaptive Uzawa-type algorithm with an outer linearization loop and
an inner adaptive mesh-refinement loop; cf.~\cite{BanschMorinNochetto2002}.
The combination of linearization and adaptive mesh refinement follows the analysis from~\cite{FuhrerPraetorius2018}
guaranteeing global convergence of the overall algorithm 
for arbitrary initial guesses \((p^0, u^0) \in H(\div, \Omega) \times H^1_0(\Omega)\).
To the best of our knowledge, this is the first global convergence result for an adaptive LSFEM
for nonlinear PDEs.
While fixed-point iterations are typically slower than the Newton method, 
the presented algorithm is particularly attractive due to the guaranteed convergence and
may be used to compute a good initial guess for a subsequent Gauss--Newton iteration.

\subsection{Outline}

The outline of the paper reads as follows.
In Section~\ref{sec:preliminaries}, we introduce the theoretical foundation of the Zarantonello iteration.
Section~\ref{sec:weighted_least_squares} presents the weighted least-squares minimization for the linearized problem
with general right-hand sides in \(H^{-1}(\Omega)\) in step~\ref{step:lsfem}.
The strongly monotone and globally Lipschitz continuous first-order system from step~\ref{step:fos} 
is introduced in Section~\ref{sec:model_problem} followed by a discussion of possible related least-squares formulations.
Section~\ref{sec:Zarantonello_least_squares} establishes the well-posedness 
of the Zarantonello-linearized least-squares formulation in step~\ref{step:zarantonello}.
Alternative weightings are discussed in Section~\ref{sec:alternative_weightings}
whereas the corresponding proofs are deferred to the Appendices~\ref{app:balanced_weighting}--\ref{app:split_weighting}.
Suitable a~posteriori error estimates allow formulating an adaptive Uzawa-type algorithm
with the adaptive LSFEM in Section~\ref{sec:adaptive_Zarantonello_LS}.
The main result of this paper is the global convergence of the adaptive algorithm in Theorem~\ref{thm:ls_plain_convergence}.
Numerical experiments in Section~\ref{sec:applications} illustrate the performance of the proposed algorithm.

%
%
\section{Preliminaries}
\label{sec:preliminaries}

%
%
\subsection{Zarantonello iteration}
\label{sec:Zarantonello}
\noindent
Consider a Hilbert space \(X\) with scalar product \(\mathcal{A}(\,\cdot\,;\,\cdot\,)\colon X \times X \to \R\)
and induced norm \(\vvvert \cdot \vvvert_{\mathcal{A}}\).
Let \(\mathcal{B}(\,\cdot\,;\,\cdot\,)\colon X \times X \to \R\) denote a nonlinear mapping which is linear in the second component.
Given a right-hand side \(\mathcal{F} \in X^*\), the corresponding nonlinear problem seeks \(x^\star \in X\) with
\begin{equation}
	\label{eq:nonlinear_weak}
	\mathcal{B}(x^\star; y)
	=
	\mathcal{F}(y)
	\quad\text{for all }
	y \in X.
\end{equation}
The well-posedness of this formulation follows from the strong monotonicity and Lipschitz continuity of \(\mathcal{B}\)
with respect to the norm \(\vvvert \cdot \vvvert_{\mathcal{A}}\), 
i.e., there exist constants \(\alpha, L > 0\) such that, for all \(x, y, z \in X\),
\begin{equation}
	\label{eq:Zarantonello:assumptions}
    \alpha\,
	\vvvert x - y \vvvert_{\mathcal{A}}^2
    \leq
	\mathcal{B}(x; x - y) - \mathcal{B}(y; x - y)
    \text{ and }
	\mathcal{B}(x; z) - \mathcal{B}(y; z)
    \leq
    L\,
	\vvvert x - y \vvvert_{\mathcal{A}}\,
	\vvvert z \vvvert_{\mathcal{A}}.
\end{equation}
In this case, the Browder--Minty theorem provides existence and uniqueness of the solution \(x^\star \in X\) 
to the nonlinear problem~\eqref{eq:nonlinear_weak}; see~\cite[Section~25.4]{Zeidler1990}.
The proof in~\cite{Zarantonello1960} employs a fixed-point iteration \(\Psi\colon X \to X\)
defined, for a damping parameter \(\delta > 0\) and a given iterate \(x^{k-1} \in X\), by
\begin{equation}
	\label{eq:Zarantonello:update}
	\mathcal{A}(\Psi(x^{k-1}); y)
	=
	\mathcal{A}(x^{k-1}; y)
	+
	\delta[ \mathcal{F}(y) + \mathcal{B}(x^{k-1}; y) ].
\end{equation}
The iteration \(x^k \coloneqq \Psi(x^{k-1})\) is well-defined
by the Riesz representation theorem for the scalar product \(\mathcal{A}\).
We will use it to linearize the nonlinear first-order system in step~\ref{step:zarantonello} of Subsection~\ref{sec:intro:approach}.
Throughout this paper, the number \(k \in \N_0\) denotes the Zarantonello iteration index.
For any small \(0 < \delta < 2\alpha/L^2\),
it is well-known from \cite[Theorem~25.B]{Zeidler1990} that \(\Psi\) is a contraction in the norm \(\vvvert \cdot \vvvert_{\mathcal{A}}\)
with factor \(0 < \rho_{\textup{Z}} \coloneqq [1 - \delta (2\alpha + L^2 \delta)]^{1/2} < 1\) such that
\begin{equation}
	\label{eq:Zarantonello:contraction}
	\vvvert  x^\star - x^k \vvvert_{\mathcal{A}}
	\leq
	\rho_{\textup{Z}}\,
	\vvvert x^\star - x^{k-1} \vvvert_{\mathcal{A}}.
\end{equation}

\subsection{Sobolev spaces}
\noindent
The nonlinear PDE~\eqref{eq:intro:nonlinear_pde} is formulated on 
a bounded Lipschitz domain \(\Omega \subset \R^d\) with polyhedral boundary \(\partial\Omega\) 
in arbitrary spatial dimension \(d \in \N\).
This paper employs standard notation for Sobolev and Lebesgue spaces
\(H^1_0(\Omega)\), \(H(\div, \Omega)\), \(L^2(\Omega)\), and \(L^2(\Omega; \R^d)\).
The \(L^2\) scalar products and norms on scalar- and vector-valued functions are denoted by the same index
in \((\cdot, \cdot)_{L^2(\Omega)}\) and \(\Vert \cdot \Vert_{L^2(\Omega)}\).
The domain-dependent Friedrichs constant is uniquely determined as the smallest possible constant \(C_{\textup{F}} > 0\)
satisfying the Friedrichs inequality
\begin{equation}
	\label{eq:friedrichs}
	\Vert v \Vert_{L^2(\Omega)} 
	\leq
	C_{\textup{F}} \,
	\Vert \nabla v \Vert_{L^2(\Omega)}
    \quad \text{for all } v \in H^1_0(\Omega).
\end{equation}
The upper bound \(C_{\textup{F}} \leq \operatorname{width}(\Omega) / \pi\) is sharp with
the width of the domain \(\Omega\) defined as the smallest possible distance 
of two parallel hyperplanes (lines in 2D, planes in 3D) enclosing \(\Omega\) in
\[
    \operatorname{width}(\Omega)
    \coloneqq
    \inf
    \left\{
        \ell > 0
        \;:\;
        \begin{gathered}
			\exists H_1, H_2 \subseteq \R^{d} \text{ hyperplanes with }
            \Omega \subseteq \conv(H_1 \cup H_2) \text{ and }
            \\
            \operatorname{dist}(H_1, H_2) \coloneqq \inf\{\vert x_1 - x_2 \vert : x_1 \in H_1, x_2 \in H_2\} = \ell,\,
        \end{gathered}
    \right\}.
\]

\subsection{Triangulations and refinement}
\noindent
The discretization will be based on conforming triangulations \(\TT\) of the bounded polyhedral domain \(\Omega\) into compact simplices.
The local mesh refinement employs a newest-vertex bisection (NVB) algorithm such as~\cite{Maubach1995,Traxler1997}.
Classical refinement algorithms require admissibility conditions on an initial conforming triangulation \(\TT_0\) 
to ensure conformity of the refined triangulation and to control the number of newly created simplices in the closure step of the refinement;
see, e.g., \cite[Section~6]{Traxler1997} and \cite[Section~4]{Stevenson2008}.
For \(d = 2\), these restrictive assumptions do not apply~\cite{kpp2013} 
and, for $d \ge 2$, the recent work~\cite{DieningGehringStorn2025} designed a novel initialization routine for arbitrary \(\TT_0\).
Concerning the case $d=1$, we refer to~\cite{affkp2013}.
For each triangulation \(\TT_H\) and marked elements \(\MM_H \subseteq \TT_H\),
let \(\TT_h \coloneqq \texttt{refine}(\TT_H, \MM_H)\) be the coarsest refinement of \(\TT_H\) 
such that at least all elements \(T \in \MM_H\) have been refined, i.e.,
\(\MM_H  \subseteq \TT_H \setminus \TT_h\).
We write \(\TT_h \in \T(\TT_H)\) if \(\TT_h\) can be obtained from \(\TT_H\) by finitely many steps of NVB,
and abbreviate \(\T \coloneqq \T(\TT_0)\).

\subsection{Finite element spaces}
\label{sec:fespaces}
\noindent
Let \(P^{m}(K)\) denote the space of polynomials on the subset \(K \subset \overline\Omega\)
of degree at most \(m \in \N_0\).
Throughout the paper, we employ conforming finite element spaces of Raviart--Thomas and Lagrange type
defined, for any \(\TT \in \T\), by
\begin{equation}
	\label{eq:finite_element_spaces}
	\begin{split}
		RT^{m}(\TT)
		&\coloneqq
		\big\{
			q_h \in H(\div, \Omega)
			\;:\;
			\forall T \in \TT, \: q_h \vert_T \in P^{m}(T; \R^d) + P^{m}(T) \cdot \operatorname{id}
		\big\},
		\\
		S^{m+1}_0(\TT)
		&\coloneqq
		\big\{
			v_h \in  H^1_0(\Omega)
			\;:\;
			\forall T \in \TT, \: v_h \vert_T \in P^{m+1}(T)
		\big\}.
	\end{split}
\end{equation}
The reader is referred to \cite[Chapter~2]{BoffiBrezziFortin2013} for a comprehensive introduction of these spaces.

\begin{remark}[other discretizations]
	For the ease of the presentation, we restrict ourselves to simplicial triangulations \(\TT\).
	However, all proofs in this paper can be generalized to \emph{any} conforming discretization
	of the Sobolev spaces \(H(\div, \Omega)\) and \(H^1_0(\Omega)\).
	In fact, the crucial result is the plain convergence result of the linear LSFEM
	in Theorem~\ref{thm:ls_plain_convergence} below.
	We refer to \cite[Section~2]{FuhrerPraetorius2020} for a detailed presentation of the sufficient conditions
	for this result.
\end{remark}

%
%
\section{Weighted least-squares minimization for linear problems}
\label{sec:weighted_least_squares}
\noindent
This section recalls the analysis of a weighted least-squares discretization on adaptively generated meshes
as we intend to use it in step~\ref{step:lsfem} of Subsection~\ref{sec:intro:approach}.
The reader is referred to~\cite{BochevGunzburger2009} for a comprehensive introduction of 
the least-squares finite element method for linear problems.
The Zarantonello iteration for the nonlinear PDE
results in a linear diffusion problem to be solved by a weighted least-squares method.
Given right-hand sides \(g_1 \in L^2(\Omega)\) and \(g_2 \in L^2(\Omega; \R^d)\) 
and positive weights \(\omega_1, \omega_2 > 0\),
the linearized PDE seeks the solution \((p^\star, u^\star) \in H(\div, \Omega) \times H^1_0(\Omega)\) to
\begin{equation}
	\label{eq:linear:pde}
	- \omega_1 \, \div p^\star 
	=
	g_1
	\quad\text{and}\quad
	p^\star - \omega_2^2 \, \nabla u^\star 
	=
	- g_2
	\quad \text{in } \Omega.
\end{equation}
For suitable right-hand sides \(g_1\) and \(g_2\), this first-order system takes the form of the Zarantonello
linearization~\eqref{eq:intro:zarantonello_fos} from step~\ref{step:zarantonello} in Subsection~\ref{sec:intro:approach}.
The weighting in the second residual is one particular choice 
in the linearized problem in Section~\ref{sec:Zarantonello_least_squares} below.
Further alternative weightings are discussed in the subsequent Section~\ref{sec:alternative_weightings}.
With the Friedrichs constant \(C_{\textup{F}} > 0\) from~\eqref{eq:friedrichs},
define the least-squares functional \(LS(g_1, g_2; \,\cdot\,,\,\cdot\,)\colon H(\div, \Omega) \times H^1_0(\Omega) \to \R\)
for the solution of the linear problem~\eqref{eq:linear:pde} as
\begin{equation}
	\label{eq:linear:least_squares_functional}
	LS(g_1, g_2; p, u)
	\coloneqq
	C_{\textup{F}}^2 \, 
	\Vert g_1 + \omega_1 \, \div p \Vert_{L^2(\Omega)}^2
	+
	\Vert g_2 + p - \omega_2^2 \, \nabla u \Vert_{L^2(\Omega)}^2.
\end{equation}
The first variation of this quadratic functional leads to the bilinear form
\(\mathcal{A}(\,\cdot\,,\,\cdot\,; \,\cdot\,,\,\cdot\,)\colon [H(\div, \Omega) \times H^1_0(\Omega)] \times
[H(\div, \Omega) \times H^1_0(\Omega)] \to \R\) with
\begin{equation}
	\label{eq:linear:ls_bilinear_form}
	\mathcal{A}(p, u; q, v)
	\coloneqq
	\omega_1^2 C_{\textup{F}}^2 \,
	(\div p,\: \div q)_{L^2(\Omega)}
	+
	(p - \omega_2^2 \, \nabla u,\: q - \omega_2^2 \, \nabla v)_{L^2(\Omega)}.
\end{equation}
The fundamental equivalence in Theorem~\ref{thm:fundamental_equivalence} below
ensures that this defines a scalar product on \(H(\div, \Omega) \times H^1_0(\Omega)\)
inducing the norm \(\vvvert \cdot \vvvert_{\mathcal{A}}\) defined by
\begin{equation}
	\label{eq:linear:ls_norm}
	\vvvert (p, u) \vvvert_{\mathcal{A}}^2
	\coloneqq
	C_{\textup{F}}^2 \,
	\Vert \omega_1 \div p \Vert_{L^2(\Omega)}^2
	+
	\Vert p - \omega_2^2 \, \nabla u \Vert_{L^2(\Omega)}^2.
\end{equation}
In Section~\ref{sec:Zarantonello_least_squares} below, 
we will use \(\mathcal{A}\) as a scalar product of a Zarantonello linearization 
on the Hilbert space \(X \coloneqq H(\div, \Omega) \times H^1_0(\Omega)\)
in the spirit of Subsection~\ref{sec:Zarantonello}.

The unique exact minimizer \((p^\star, u^\star) \in H(\div, \Omega) \times H^1_0(\Omega)\) 
of the functional~\eqref{eq:linear:least_squares_functional} with
\[
	LS(g_1, g_2; p^\star, u^\star)
	=
	\min_{(p,u) \in H(\div, \Omega) \times H^1_0(\Omega)}
	LS(g_1, g_2; p, u)
\]
is characterized by the Euler--Lagrange equation, for all \((q, v) \in H(\div, \Omega) \times H^1_0(\Omega)\),
\begin{equation}
	\label{eq:linear:ls_variational}
	\mathcal{A}(p^\star, u^\star; q, v)
	=
	- C_{\textup{F}}^2 \,
	(g_1, \omega_1 \, \div q)_{L^2(\Omega)}
	-
	(g_2, q - \omega_2^2 \, \nabla v)_{L^2(\Omega)}.
\end{equation}
The well-posedness of this formulation follows from the equivalence of the norm \(\vvvert \cdot \vvvert_{\mathcal{A}}\) with 
the weighted norm on the product space \(H(\div, \Omega) \times H^1_0(\Omega)\) given by
\begin{equation}
	\label{eq:weighted_norm}
	\vvvert (p, u) \vvvert^2
	\coloneqq
	C_{\textup{F}}^2 \,
	\Vert \omega_1 \, \div p \Vert_{L^2(\Omega)}^2
	+
	\Vert p \Vert_{L^2(\Omega)}^2
	+
	\Vert \omega_2^2 \, \nabla u \Vert_{L^2(\Omega)}^2.
\end{equation}
This norm is equivalent to the unweighted norm on \(H(\div, \Omega) \times H^1_0(\Omega)\) defined by
\begin{equation}
	\label{eq:unweighted_norm}
	\vvvert (p, u) \vvvert_{\textup{uw}}^2
	\coloneqq
	C_{\textup{F}}^2 \,
	\Vert \div p \Vert_{L^2(\Omega)}^2
	+
	\Vert p \Vert_{L^2(\Omega)}^2
	+
	\Vert \nabla u \Vert_{L^2(\Omega)}^2.
\end{equation}
In fact, for all \((p, u) \in H(\div, \Omega) \times H^1_0(\Omega)\), it holds that
\[
	\min\big\{
		1,\:
		\omega_1^2,\:
		\omega_2^4
	\big\} \,
	\vvvert (p, u) \vvvert_{\textup{uw}}^2
	\leq
	\vvvert (p, u) \vvvert^2
	\leq
	\max\big\{
		1,\:
		\omega_1^2,\:
		\omega_2^4
	\big\} \,
	\vvvert (p, u) \vvvert_{\textup{uw}}^2.
\]
Here, the consistent weighting with the Friedrichs constant \(C_{\textup{F}}\) 
ensures that the fundamental equivalence constants are independent of the size of the domain \(\Omega\) 
(and even the spatial dimension \(d \in \N\)).
The authors assume the following result to be well-known; see, e.g., 
the original contributions \cite[Lemma~4.3]{Jespersen1977} for the Poisson model problem
and \cite[Theorem~2.1]{PehlivanovCareyLazarov1994} for second-order elliptic diffusion problems
and the more recent publication \cite[Proof of Lemma~6.2]{ErnGudiSmearsVohralik2022}
for a similar weighting as in~\eqref{eq:linear:ls_bilinear_form}.
However, the proof is given here in detail for the sake of explicit constants in terms of \(\omega_1\) and \(\omega_2\).
\begin{theorem}[fundamental equivalence]
    \label{thm:fundamental_equivalence}
	For any \(q \in H(\div, \Omega)\) and \(v \in H^1_0(\Omega)\),
    \begin{equation}
		\label{eq:fundamental_equivalence:weighted}
		\min\bigg\{
			\frac12,\:
			\bigg(1 + \frac{4}{\omega_1^2}\bigg)^{-1}
		\bigg\} \,
		\vvvert (q, v) \vvvert^2
		\leq
		\vvvert (q, v) \vvvert_{\mathcal{A}}^2
		\leq
		2 \,
		\vvvert (q, v) \vvvert^2.
    \end{equation}
\end{theorem}
\begin{proof}
    \emph{Step~1.}\,
	The proof of the \emph{ellipticity} of the least-squares functional
	(i.e., the lower bound in~\eqref{eq:fundamental_equivalence:weighted}) departs from the
    binomial formula followed by an integration by parts
    \begin{align*}
		\Vert q \Vert_{L^2(\Omega)}^2
		+
		\Vert \omega_2^2 \, \nabla v \Vert_{L^2(\Omega)}^2
		&=
		\Vert q - \omega_2^2 \, \nabla v \Vert_{L^2(\Omega)}^2
		+
		2\omega_2^2 \, (q, \nabla v)_{L^2(\Omega)}
		\\
		&=
		\Vert q - \omega_2^2 \, \nabla v \Vert_{L^2(\Omega)}^2
		-
		2 \omega_2^2 \, (\div q, v)_{L^2(\Omega)}.
    \end{align*}
	The Cauchy--Schwarz, the Friedrichs, and a weighted Young inequality imply
    \begin{align*}
		- 2 \omega_2^2 \,
		(\div q, v)_{L^2(\Omega)}
		&\leq
		2 \omega_2^2 \,
		\Vert \div q \Vert_{L^2(\Omega)} \, \Vert v \Vert_{L^2(\Omega)}
		\leq
		2 C_{\textup{F}} \,
		\Vert \div q \Vert_{L^2(\Omega)} \, \Vert \omega_2^2 \, \nabla v \Vert_{L^2(\Omega)}
		\\
		&\leq
		\frac{2C_{\textup{F}}^2}{\omega_1^2} \,
		\Vert \omega_1 \, \div q \Vert_{L^2(\Omega)}^2
		+
		\frac12 \,
		\Vert \omega_2^2 \, \nabla v \Vert_{L^2(\Omega)}^2.
    \end{align*}
	The combination of the two previous displayed formulas
	and the absorption of \(\frac12 \, \Vert \omega_2^2 \, \nabla v \Vert_{L^2(\Omega)}^2\) into the left-hand side read
    \[
		2 \,
		\Vert q \Vert_{L^2(\Omega)}^2
        +
		\Vert \omega_2^2 \, \nabla v \Vert_{L^2(\Omega)}^2
        \leq
		\frac{4C_{\textup{F}}^2}{\omega_1^2} \,
		\Vert \omega_1 \, \div q \Vert_{L^2(\Omega)}^2
        +
		2 \, \Vert q - \omega_2^2 \, \nabla v \Vert_{L^2(\Omega)}^2.
    \]
	The addition of \(C_{\textup{F}}^2 \, \Vert \omega_1 \, \div q \Vert_{L^2(\Omega)}^2\) results in
    \begin{align*}
		\vvvert (q, v) \vvvert^2
		&\leq
		C_{\textup{F}}^2 \,
		\Vert \omega_1 \, \div q \Vert_{L^2(\Omega)}^2
		+
		2 \,
		\Vert q \Vert_{L^2(\Omega)}^2
        +
		\Vert \omega_2^2 \, \nabla v \Vert_{L^2(\Omega)}^2
		\\
		&\leq
		\bigg(1 + \frac{4}{\omega_1^2}\bigg) \, C_{\textup{F}}^2 \, 
		\Vert \omega_1 \, \div q \Vert_{L^2(\Omega)}^2
        +
		2\, \Vert q - \omega_2^2 \, \nabla v \Vert_{L^2(\Omega)}^2.
    \end{align*}
	This concludes the proof of the lower bound with
    \[
		\vvvert (q, v) \vvvert^2
        \leq
		\max\bigg\{ 1 + \frac{4}{\omega_1^2},\: 2 \bigg\} \, 
		\vvvert (q, v) \vvvert_{\mathcal{A}}^2.
    \]

	\emph{Step~2.}\,
    The proof of the \emph{boundedness} of the least-squares functional
	(i.e., the upper bound in the estimate~\eqref{eq:fundamental_equivalence:weighted})
	employs the triangle inequality and the Young inequality to establish
	\begin{align*}
		\Vert q - \omega_2^2 \, \nabla v \Vert_{L^2(\Omega)}^2
		&\leq
		2 \,
		\big[
			\Vert q \Vert_{L^2(\Omega)}^2
			+
			\Vert \omega_2^2 \, \nabla v \Vert_{L^2(\Omega)}^2
		\big].
	\end{align*}
	The addition of \(C_{\textup{F}}^2 \, \Vert \omega_1 \, \div q \Vert_{L^2(\Omega)}^2\) concludes the proof with
	\(\vvvert (q, v) \vvvert_{\mathcal{A}}^2 \leq 2 \, \vvvert (q, v) \vvvert^2\).
\end{proof}

The fundamental equivalence ensures well-posedness of the continuous least-squares problem~\eqref{eq:linear:ls_variational}
as well as of the discrete LSFEM of piecewise polynomial degree \(m \in \N_0\).
The latter seeks the discrete minimizers \((p_h, u_h) \in RT^{m}(\TT) \times S^{m+1}_0(\TT)\) 
of the functional~\eqref{eq:linear:least_squares_functional} satisfying
\[
	LS(g_1, g_2; p_h, u_h)
	=
	\min_{(q_h, v_h) \in RT^m(\TT) \times S^{m+1}_0(\TT)}
	LS(g_1, g_2; q_h, v_h).
\]
The unique discrete minimizer is characterized by, for all \((q_h, v_h) \in RT^{m}(\TT) \times S^{m+1}_0(\TT)\),
\begin{equation}
	\label{eq:linear:lsfem_variational}
	\mathcal{A}(p_h, u_h; q_h, v_h)
	=
	- C_{\textup{F}}^2 \,
	(g_1, \omega_1 \, \div q_h)_{L^2(\Omega)}
	-
	(g_2, q_h - \omega_2^2 \ \nabla v_h)_{L^2(\Omega)}.
\end{equation}
Another immediate consequence of the fundamental equivalence~\eqref{eq:fundamental_equivalence:weighted}
is the built-in a~posteriori error estimate for \emph{every conforming} approximation
\(q \in H(\div, \Omega)\) and \(v \in H^1_0(\Omega)\) to the exact solution \((p^\star, u^\star)\) of the least-squares problem
\begin{equation}
	\label{eq:linear:aposteriori}
	LS(g_1, g_2; q, v)
	\eqsim
	\vvvert (p^\star - q, u^\star - v) \vvvert^2.
\end{equation}
This estimate is even asymptotically exact in the sense that the ratio of least-squares functional
and error of the exact discrete solution tend to one under uniform refinement \cite[Theorem~3.1]{CarstensenStorn2018}.
The a~posteriori estimate~\eqref{eq:linear:aposteriori} motivates 
the definition of an a~posteriori error estimator by the local contributions to the least-squares functional
\begin{equation}
	\label{eq:linear:eta}
	\eta(T; q, v)^2
	\coloneqq
	C_{\textup{F}}^2 \,
	\Vert g_1 + \omega_1 \, \div q \Vert_{L^2(T)}^2
	+
	\Vert g_2 + q - \omega_2^2 \, \nabla v \Vert_{L^2(T)}^2
\end{equation}
with the full contribution abbreviated as
\[
	LS(g_1, g_2; q, v)
	=
	\vvvert (p^\star - q, u^\star - v) \vvvert_{\mathcal{A}}^2
	=
	\eta(q, v)^2
	\coloneqq
	\sum_{T \in \TT} \eta(T; q, v)^2.
\]
The local contributions \(\eta(T; q, v)\) are used to steer the adaptive mesh refinement in Algorithm~\ref{alg:ALSFEM}.
\begin{algorithm}[t]
	\caption{Adaptive least-squares FEM (ALSFEM) for linear problem~\eqref{eq:linear:lsfem_variational}}
	\label{alg:ALSFEM}
	\flushleft
	\noindent
	{\bfseries Input:}
	Initial mesh \(\TT_0\),
	marking parameter \(0 < \theta \le 1\), tolerance \(\tau \geq 0\).
	\medskip

	\noindent
	\textbf{for} \(\ell = 0, 1, 2, \dots\) \textbf{do}
	\begin{enumerate}[label=(\alph*),leftmargin=2em]
		\item
			\textbf{Solve.}
			Compute the discrete solutions
			\((p_\ell, u_\ell) \in RT^{m}(\TT_\ell) \times S^{m+1}_0(\TT_\ell)\)
			to~\eqref{eq:linear:lsfem_variational}.
		\item
			\textbf{Estimate.}
			Compute the refinement indicators \(\eta_\ell(T; p_\ell, u_\ell)\) from~\eqref{eq:linear:eta}
			for all \(T \in \TT_\ell\).
		\item
			\textbf{If} \(\eta_\ell(p_\ell, u_\ell) \leq \tau\),
			\textbf{then break} the \(\ell\) loop and terminate.
		\item
			\textbf{Mark.}
			Determine a set \(\MM_\ell \subseteq \TT_\ell\) of minimal cardinality satisfying
			\[
				\theta \, \eta_\ell(p_\ell, u_\ell)^{2}
				\leq
				\sum_{T \in \MM_\ell} \eta_\ell(T; p_\ell, u_\ell)^{2}.
			\]
		\item
			\textbf{Refine.}
			Generate the refined mesh \(\TT_{\ell+1} \coloneqq \refine(\TT_\ell, \MM_\ell)\) by NVB.
	\end{enumerate}
	\textbf{end for}
	\medskip

	\noindent
	{\bfseries Output:}
	Sequence of successively refined triangulations \(\TT_\ell\) with
	corresponding discrete solutions 
	\((p_\ell, u_\ell) \in RT^{m}(\TT_\ell) \times S^{m+1}_0(\TT_\ell)\).
\end{algorithm}
The plain convergence analysis from~\cite{Siebert2011} applies to Algorithm~\ref{alg:ALSFEM}
and provides the following result
independently proven in \cite[Theorem~2]{FuhrerPraetorius2020} and~\cite[Theorem~3.3]{GantnerStevenson2021}.
\begin{theorem}[plain convergence]
	\label{thm:ls_plain_convergence}
	For \(\tau = 0\), the output \((p_\ell, u_\ell)_{\ell \in \N_0}\) of Algorithm~\ref{alg:ALSFEM} satisfies
	\[
		\vvvert (p^\star - p_\ell, u^\star - u_\ell) \vvvert^2
		+
		LS(g_1, g_2; p_\ell, u_\ell)
		\to
		0
		\quad\text{as }
		\ell \to \infty.
	\]
	For positive tolerance \(\tau > 0\), the algorithm thus terminates after finitely many steps. \hfill \qed
\end{theorem}
We refer to~\cite[Theorem~4.1]{CarstensenParkBringmann2017} for the stronger result of Q-linear convergence
under the additional assumptions of sufficient data approximation and sufficiently large bulk parameter \(\theta\).
Note that the convergence analysis with rates of minimal residual methods 
\cite{CarstensenPark2015,CarstensenHellwig2018,Carstensen2020,CarstensenMa2021,Bringmann2023,Bringmann2024} has so far only been established 
for alternative residual-based error estimators instead of the built-in estimator used in Algorithm~\ref{alg:ALSFEM}.

%
%
\section{Strongly monotone model problem}
\label{sec:model_problem}
\noindent
This section discusses the detailed assumptions guaranteeing the well-posedness of the nonlinear problem~\eqref{eq:intro:nonlinear_primal_pde}.
We also introduce the first-order reformulation as required for step~\ref{step:fos} in Subsection~\ref{sec:intro:approach}.
The remaining part of this section discusses two straight-forward, but less favorable applications of the least-squares method
to solve the nonlinear problem.

\subsection{Nonlinear first-order system}
Let the right-hand sides \(f_1 \in L^2(\Omega)\) and \(f_2 \in L^2(\Omega; \R^d)\) be given.
Throughout the paper, consider a nonlinear flux mapping \(\sigma\colon \R^d \to \R^d\)
in the nonlinear elliptic PDE with homogeneous Dirichlet boundary conditions.
It seeks \(u^\star \in H^1_0(\Omega)\) such that
\begin{equation}
    \label{eq:nonlinear}
    - \div(\sigma(\nabla u^\star))
    =
	f_1 - \div f_2
    \quad\text{in } \Omega.
\end{equation}
The corresponding weak formulation takes the form of problem~\eqref{eq:nonlinear_weak}
with \(X \coloneqq H^1_0(\Omega)\) as well as
the nonlinear mapping \(\widehat{\mathcal{B}}(\,\cdot\,;\,\cdot\;)\colon H^1_0(\Omega) \times H^1_0(\Omega) \to \R\)
and the right-hand side \(\widehat{\mathcal{F}} \in H^{-1}(\Omega)\) defined by
\[
	\widehat{\mathcal{B}}(u; v)
    \coloneqq
	(\sigma(\nabla u), \nabla v )_{L^2(\Omega)}
	\quad\text{and}\quad
	\widehat{\mathcal{F}}(v)
	\coloneqq
	(f_1, v)_{L^2(\Omega)} + (f_2, \nabla v)_{L^2(\Omega)}.
\]
The primal formulation of the Zarantonello iteration employs
the scalar product \(\widehat{\mathcal{A}}(\,\cdot\,;\,\cdot\,)\colon H^1_0(\Omega) \times H^1_0(\Omega) \to \R\)
with \(\widehat{\mathcal{A}}(u; v) \coloneqq (\nabla u, \nabla v)_{L^2(\Omega)}\) for all \(u, v \in H^1_0(\Omega)\).
Given any \(u^{k-1} \in H^1_0(\Omega)\) and a damping parameter \(\delta > 0\),
it seeks the exact solution \(u^k_\exact \in H^1_0(\Omega)\) satisfying,
for all \(v \in H^1_0(\Omega)\),
\begin{equation}
	\label{eq:Zarantonello:primal}
	\widehat{\mathcal{A}}(u^k_\star; v)
	=
	\widehat{\mathcal{A}}(u^{k-1}; v)
	+
	\delta \,
	[\widehat{\mathcal{F}}(v) - \widehat{\mathcal{B}}(u^{k-1}; v)].
\end{equation}
In order to verify the assumptions~\eqref{eq:Zarantonello:assumptions} from Subsection~\ref{sec:Zarantonello},
suppose that the nonlinear mapping \(\sigma \in C^1(\R^d; \R^d)\) is differentiable
and that the derivative \(\D\sigma\colon \R^d \to \R^{d \times d}\) is uniformly elliptic and bounded,
i.e., there exist constants \(0 < \Lambda_1 < \Lambda_2 < \infty\) such that:
\begin{description}
	\item[\bfseries (N1) ellipticity]\labeltext{N1}{assum:ellipticity} 
		For all \(\xi, a \in \R^d\), it holds
		\(\Lambda_1 \, \vert a \vert^2 \leq (\D\sigma(\xi) \, a) \cdot a\).
	\item[\bfseries (N2) boundedness]\labeltext{N2}{assum:boudedness} 
		For all \(\xi, a, b \in \R^d\), it holds
		\(\vert (\D\sigma(\xi) \, a) \cdot b \vert \leq \Lambda_2 \, \vert a \vert \, \vert b \vert\).
\end{description}
These two assumptions formalize a 2-growth condition on the nonlinearity \(\sigma\) and, thereby, 
restrict our analysis to the Hilbert space setting.
We refer to the literature presented in the introductory Subsection~\ref{sec:intro:literature} for methods beyond it.

The fundamental theorem of calculus ensures that, for any \(u, v \in H^1_0(\Omega)\),
the componentwise integrated matrix \(M \coloneqq \int_0^1 \D\sigma(\nabla (u + s(v - u))) \d{s} \in L^\infty(\Omega; \R^{d \times d})\) satisfies
\begin{equation}
	\label{eq:stress_difference}
	\sigma(\nabla u) - \sigma(\nabla v)
	=
	\int_0^1 \frac{\textup{d}}{\textup{d} s} \sigma(\nabla (u + s(v - u))) \d{s}
	=
	M \, \nabla(v - u)
\end{equation}
almost everywhere in \(\Omega\).
Under the assumptions~\eqref{assum:ellipticity}--\eqref{assum:boudedness}, the relation
\[
	\widehat{\mathcal{B}}(u; z) - \widehat{\mathcal{B}}(v; z)
	=
	(\sigma(\nabla u) - \sigma(\nabla v), \nabla z)_{L^2(\Omega)}
	=
	(M \, \nabla(v - u), \nabla z)_{L^2(\Omega)}
\]
for all \(u, v, z \in H^1_0(\Omega)\) reveals that the mapping \(\widehat{\mathcal{B}}\)
satisfies~\eqref{eq:Zarantonello:assumptions} with \(\alpha = \Lambda_1\) and \(L = \Lambda_2\).
In particular, the nonlinear PDE~\eqref{eq:nonlinear} is well-posed
and the iteration~\eqref{eq:Zarantonello:primal} is contractive
for any \(0 < \delta < \delta^\star \coloneqq 2\Lambda_1 / \Lambda_2^2\).

The introduction of the additional flux-like variable \(p^\star \coloneqq \sigma(\nabla u^\star) - f_2 \in H(\div, \Omega)\)
for the exact solution \(u^\star \in H^1_0(\Omega)\) from~\eqref{eq:nonlinear} leads to
the nonlinear first-order system of step~\ref{step:fos} in Subsection~\ref{sec:intro:approach}
\begin{equation}
    \label{eq:nonlinear_fos}
	- \div p^\star
	=
	f_1 
	\quad\text{and}\quad
	- p^\star + \sigma (\nabla u^\star)
	=
	f_2
	\quad
	\text{in } \Omega.
\end{equation}

\subsection{Nonlinear least-squares minimization}
\label{sec:ls_nonlinear}
As an alternative to the steps~\ref{step:zarantonello}--\ref{step:lsfem} in Subsection~\ref{sec:intro:approach},
one can directly apply the least-squares minimization to the nonlinear first-order system~\eqref{eq:nonlinear_fos} of PDEs.
To this end, define the (weighted) nonlinear residual mapping
\(\mathcal{R}(f_1, f_2; \,\cdot\,,\,\cdot\,)\colon H(\div, \Omega) \times H^1_0(\Omega) \to L^2(\Omega) \times L^2(\Omega; \R^d)\)
by, for all \((p, u) \in H(\div, \Omega) \times H^1_0(\Omega)\),
\[
	\mathcal{R}(f_1, f_2; p, u)
	\coloneqq
	\big( C_{\textup{F}}\, (f_1 + \div p),\: f_2 + p - \sigma (\nabla u) \big).
\]
The corresponding least-squares functional 
\(N(f_1, f_2;\,\cdot\,,\,\cdot\,)\colon H(\div, \Omega) \times H^1_0(\Omega) \to \R\) reads
\begin{equation}
	\label{eq:nonlinear_ls_functional}
	N(f_1, f_2; p, u)
	\coloneqq
	\Vert \mathcal{R}(f_1, f_2; p, u) \Vert_{L^2(\Omega)}^2
	=
	C_{\textup{F}}^2 \, \Vert f_1 + \div p \Vert_{L^2(\Omega)}^2
	+
	\Vert f_2 + p - \sigma (\nabla u) \Vert_{L^2(\Omega)}^2.
\end{equation}
The nonlinear least-squares formulation seeks 
minimizers \((p^\star, u^\star) \in H(\div, \Omega) \times H^1_0(\Omega)\) satisfying
\begin{equation}
	\label{eq:nonlinear_minimization}
	N(f_1, f_2; p^\star, u^\star)
	=
	\min_{(p, u) \in H(\div, \Omega) \times H^1_0(\Omega)}
	N(f_1, f_2; p, u).
\end{equation}
The unique solution~\(u^\star \in H^1_0(\Omega)\) to~\eqref{eq:nonlinear_weak} and 
\(p^\star \coloneqq \sigma(\nabla u^\star) - f_2\) obviously minimize
the non-negative functional~\eqref{eq:nonlinear_ls_functional}.
The following result shows that the minimization of the least-squares functional~\eqref{eq:nonlinear_ls_functional}
is justified and provides a solution to the nonlinear PDE~\eqref{eq:nonlinear}.
\begin{lemma}[nonlinear fundamental equivalence]
	\label{lem:nonlinear_fundamental_equivalence}
	Suppose the derivative \(\D\sigma\) satisfies~\eqref{assum:ellipticity}--\eqref{assum:boudedness}
	and that it is pointwise symmetric, i.e., \(\D\sigma(\xi) = \D\sigma(\xi)^\top\) for all \(\xi \in \R^d\).
	Then, for all \((p, u), (q, v) \in H(\div, \Omega) \times H^1_0(\Omega)\), there holds the equivalence
	\begin{equation}
		\label{eq:nonlinear_fundamental_equivalence}
		\Vert \mathcal{R}(f_1, f_2; p, u) 
		- \mathcal{R}(f_1, f_2; q, v) \Vert_{L^2(\Omega)}^2
		\eqsim
		C_{\textup{F}}^2 \,
		\Vert \div (p - q) \Vert_{L^2(\Omega)}^2
		+
		\Vert p - q \Vert_{L^2(\Omega)}^2
		+
		\Vert \nabla (u - v) \Vert_{L^2(\Omega)}^2.
	\end{equation}
	The hidden equivalence constants depend only on \(\Lambda_1\) and \(\Lambda_2\).
	In particular, they are independent of the size of the domain \(\Omega\).
\end{lemma}
\begin{proof}
	The proof follows verbatim the proof of \cite[Lemma~4.2]{CarstensenBringmannHellwigWriggers2018}
	with a straight-forward modification for robust constants independent of \(C_{\textup{F}}\).
	Since it is only presented for the convex energy minimization problem therein, 
	the proof makes use of the pointwise symmetry of \(\D\sigma\).
\end{proof}
This equivalence implies that the exact minimizers \((p^\star, u^\star) \in H(\div, \Omega) \times H^1_0(\Omega)\)
in~\eqref{eq:nonlinear_minimization} are indeed unique.
The direct method of calculus of variations proves the existence of a discrete minimizer
in any finite dimensional subspace.
The uniqueness of such discrete minimizers, however, is not guaranteed in general
because the nonlinear least-squares functional~\eqref{eq:nonlinear_ls_functional} might not be strictly convex.
As for the linear case in~\eqref{eq:linear:aposteriori},
the nonlinear least-squares functional~\eqref{eq:nonlinear_ls_functional} provides a built-in a~posteriori error estimator.
\begin{proposition}[a~posteriori error estimates]
	\label{prop:nonlinear:aposteriori}
	Under the assumptions of Lemma~\ref{lem:nonlinear_fundamental_equivalence},
	for any approximation \((q, v) \in H(\div, \Omega) \times H^1_0(\Omega)\)
	to the exact solution \((p^\star, u^\star) \in H(\div, \Omega) \times H^1_0(\Omega)\)
	with vanishing residual \(\mathcal{R}(f_1, f_2; p^\star, u^\star) = 0\) in \(L^2(\Omega) \times L^2(\Omega; \R^d)\),
	the nonlinear fundamental equivalence~\eqref{eq:nonlinear_fundamental_equivalence} implies
	\begin{align}
		\label{eq:nonlinear:aposteriori}
		\vvvert (p^\star - q, u^\star - v) \vvvert_{\textup{uw}}^2
		\eqsim
		\Vert \mathcal{R}(f_1, f_2; q, v) \Vert_{L^2(\Omega)}^2
		=
		C_{\textup{F}}^2 \,
		\Vert f_1 + \div q \Vert_{L^2(\Omega)}^2
		+
		\Vert f_2 + q - \sigma(\nabla v) \Vert_{L^2(\Omega)}^2.
		\qed
	\end{align}
\end{proposition}%

The Euler--Lagrange equations for the minimization problem~\eqref{eq:nonlinear_minimization} 
of the least-squares functional~\eqref{eq:nonlinear_ls_functional} 
seek \((p^\star, u^\star) \in H(\div, \Omega) \times H^1_0(\Omega)\) with,
for all \((q, v) \in H(\div, \Omega) \times H^1_0(\Omega)\),
\begin{equation}
	\label{eq:nonlinear:least_squares_variational}
	C_{\textup{F}}^2 \, 
	(f_1 + \div p^\star, \div q)_{L^2(\Omega)}
	+
	(f_2 + p^\star - \sigma(\nabla u^\star), q - \D\sigma(\nabla u^\star) \nabla v)_{L^2(\Omega)}
	=
	0
\end{equation}
However, this formulation does not fit into the framework of the Zarantonello iteration
as both Lipschitz continuity and strong monotonicity of the associated mapping \(\widehat{\mathcal{B}}\) are unclear.
This is due to the fact that, in applications, the second derivative \(\D^2\sigma\) is not necessarily bounded.
The reader is referred to the discussion in \cite[Section~4]{BringmannCarstensenTran2022} for further details.

\subsection{Least-squares minimization for primal Zarantonello iteration}
\label{sec:ls_primal_zarantonello}
A direct approach in line with the step~\ref{step:zarantonello} from Subsection~\ref{sec:intro:approach} 
applies the least-squares discretization to the primal Zarantonello iteration~\eqref{eq:Zarantonello:primal}.
It determines the weak solution \(\widehat{u}^k_\star \in H^1_0(\Omega)\) to the linear PDE
\begin{equation}
	\label{eq:Zarantonello:primal_strong}
	- \Delta \widehat{u}^k_\star
	=
	- \Delta u^{k-1}
	+
	\delta [f_1 - \div f_2 + \div \sigma(\nabla u^{k-1})]
	\quad\text{in } \Omega.
\end{equation}
As usual for the least-squares formulation for the Poisson problem~\cite{PehlivanovCareyLazarov1994},
a first-order reformulation allows avoiding \(C^1\) conforming discretizations.
For the additional variable 
\(\widehat{p}^k_\star \coloneqq \nabla \widehat{u}^k_\star - \nabla u^{k-1} - \delta [f_2 - \sigma(\nabla u^{k-1})]\),
a first-order system equivalent to \eqref{eq:Zarantonello:primal_strong} reads
\[
	- \div \widehat{p}^k_\star
	=
	\delta f_1
	\quad\text{and}\quad
	- \widehat{p}^k_\star + \nabla \widehat{u}^k_\star
	=
	\nabla u^{k-1} + \delta [f_2 - \sigma(\nabla u^{k-1})]
	\quad\text{in } \Omega.
\]
The corresponding least-squares formulation seeks minimizers 
\((\widehat{p}^k_\star, \widehat{u}^k_\star) \in H(\div, \Omega) \times H^1_0(\Omega)\)
of
\begin{align*}
	(p, u)
	\mapsto
	C_{\textup{F}}^2 \,
	\Vert \delta f_1 + \div p \Vert_{L^2(\Omega)}^2
	+
	\Vert p - \nabla u + \nabla u^{k-1} + \delta [f_2 - \sigma(\nabla u^{k-1})] \Vert_{L^2(\Omega)}^2.
\end{align*}
They are characterized by the Euler--Lagrange equations,
for all \(q \in H(\div, \Omega)\) and \(v \in H^1_0(\Omega)\),
\begin{multline*}
	C_{\textup{F}}^2 \,
	(\div \widehat{p}^k_\star, \div q)_{L^2(\Omega)}
	+
	(\widehat{p}^k_\star - \nabla \widehat{u}^k_\star, q - \nabla v)_{L^2(\Omega)}
	\\
	=
	- C_{\textup{F}}^2 \,
	(\delta f_1, \div q)_{L^2(\Omega)}
	-
	(\nabla u^{k-1} + \delta [f_2 - \sigma(\nabla u^{k-1})], q - \nabla v)_{L^2(\Omega)}.
\end{multline*}
The standard conforming finite element spaces from Subsection~\ref{sec:fespaces}
allow for the discrete solution of the Zarantonello iteration~\eqref{eq:Zarantonello:primal_strong}.
While this is a perfectly justified discretization of the linearized problem 
to realize step~\ref{step:lsfem} in Subsection~\ref{sec:intro:approach},
it might be less appealing to explicitly approximate the physically meaningless variable \(\widehat{p}^k_\star\).

%
%
\section{Zarantonello least-squares formulation}
\label{sec:Zarantonello_least_squares}
\noindent
The two approaches from the Subsections~\ref{sec:ls_nonlinear}--\ref{sec:ls_primal_zarantonello} appear disadvantageous.
Instead, the minimal residual methods in the literature usually employ a linearize-first approach
such as a Gauss--Newton scheme; cf.\ \cite{BertrandBrodbeckRickenSchneider2025}.
The remaining part of the paper is devoted to the development of a least-squares formulation of the nonlinear problem
employing the fixed-point iteration of Zarantonello from Section~\ref{sec:Zarantonello} 
as the step~\ref{step:zarantonello} from Subsection~\ref{sec:intro:approach}.
The formulation will be the basis for the finite element discretization in the final step~\ref{step:lsfem}.
This section derives the least-squares formulation for the linearized problem and 
establishes its well-posedness by providing sufficient choices of the damping parameter \(\delta\) 
and of the weights \(\omega_1\) and \(\omega_2\) in~\eqref{eq:intro:zarantonello_fos}.

For any given iterate \((p^{k-1}, u^{k-1}) \in H(\div, \Omega) \times H^1_0(\Omega)\) and damping parameter \(\delta > 0\),
the formal application of the Zarantonello iteration \eqref{eq:Zarantonello:primal}
to the nonlinear system of PDEs~\eqref{eq:nonlinear_fos}
seeks \((p^k_\star, u^k_\star) \in H(\div, \Omega) \times H^1_0(\Omega)\) such that
\begin{align*}
	- \omega_1 \div p^k_\star
	&=
	- \omega_1 \div p^{k-1}
	+
	\delta \omega_1 \, \big[ f_1 + \div p^{k-1} \big],
	\\
	p^k_\star
	-
	\omega_2^2 \, \nabla u^k_\star
	&=
	p^{k-1}
	-
	\omega_2^2 \, \nabla u^{k-1}
	-
	\delta \, \big[ f_2 + p^{k-1} - \sigma(\nabla u^{k-1}) \big].
\end{align*}
This is a linear system of the form~\eqref{eq:linear:pde} and can be treated as in Section~\ref{sec:weighted_least_squares}.
For this reason, we apply the least-squares approach and define the corresponding
Zarantonello least-squares functional \(Z_k(f_1, f_2; \,\cdot\,;\,\cdot\,)\colon H(\div, \Omega) \times H^1_0(\Omega) \to \R\)
with some scalar weights \(\omega_1, \omega_2 > 0\) by
\begin{equation}
	\label{eq:linearized_LS:functional}
	\begin{aligned}
		Z_k(f_1, f_2; p, u)
		&\coloneqq
		\omega_1^2 C_{\textup{F}}^2 \,
		\Vert \div (p - p^{k-1}) + \delta\, [f_1 + \div p^{k-1}] \Vert_{L^2(\Omega)}^2
		\\
		&\phantom{{}\coloneqq{}}
		+
		\Vert
			p - p^{k-1} - \omega_2^2 \, \nabla (u - u^{k-1}) 
			+
			\delta\, [f_2 + p^{k-1} - \sigma(\nabla u^{k-1})]
		\Vert_{L^2(\Omega)}^2.
	\end{aligned}
\end{equation}
The least-squares formulation seeks the exact minimizer \((p^k_\star, u^k_\star) \in H(\div, \Omega) \times H^1_0(\Omega)\) satisfying
\begin{equation}
	\label{eq:linearized_LS:minimization}
	Z_k(f_1, f_2; p^k_\star, u^k_\star)
	=
	\min_{(p,u) \in H(\div, \Omega) \times H^1_0(\Omega)}
	Z_k(f_1, f_2; p, u).
\end{equation}
The first variation of the quadratic functional~\eqref{eq:linearized_LS:functional} leads to
the scalar product \(\mathcal{A}(\,\cdot\,,\,\cdot\,;\,\cdot\,,\,\cdot\,)\) from~\eqref{eq:linear:ls_bilinear_form}
with induced norm \(\vvvert \cdot \vvvert_{\mathcal{A}}\) in~\eqref{eq:linear:ls_norm}
as well as the nonlinear mapping \(\mathcal{B}(\,\cdot\,,\,\cdot\,;\,\cdot\,,\,\cdot\,)\colon 
[H(\div, \Omega) \times H^1_0(\Omega)] \times [H(\div, \Omega) \times H^1_0(\Omega)] \to \R\)
and the right-hand side \(\mathcal{F} \in [H(\div, \Omega) \times H^1_0(\Omega)]^*\) with
\begin{subequations}
\label{eq:nonlinear:operator_and_rhs}
\begin{align}
	\label{eq:nonlinear:least_squares_operator}
	\mathcal{B}(p, u; q, v)
	&\coloneqq
	\omega_1^2 C_{\textup{F}}^2 \,
	(\div p,\: \div q)_{L^2(\Omega)}
	+
	(p - \sigma(\nabla u),\: q - \omega_2^2 \, \nabla v)_{L^2(\Omega)},
	\\
	\mathcal{F}(q, v) 
	&\coloneqq
	- \omega_1^2 C_{\textup{F}}^2 \, (f_1, \div q)_{L^2(\Omega)}
	- (f_2, q - \omega_2^2 \, \nabla v)_{L^2(\Omega)}.
\end{align}
\end{subequations}
The Euler--Lagrange equation for the minimization~\eqref{eq:linearized_LS:minimization} 
of the Zarantonello least-squares functional~\eqref{eq:linearized_LS:functional}
seeks \((p^k_\star, u^k_\star) \in H(\div, \Omega) \times H^1_0(\Omega)\) satisfying,
for all \((q, v) \in H(\div, \Omega) \times H^1_0(\Omega)\),
\begin{equation}
	\label{eq:linearized_LS:iteration}
	\mathcal{A}(p^k_\star, u^k_\star; q, v)
	=
	\mathcal{A}(p^{k-1}, u^{k-1}; q, v)
	+
	\delta [ \mathcal{F}(q, v) - \mathcal{B}(p^{k-1}, u^{k-1}; q, v) ].
\end{equation}
In explicit terms, this reads
\begin{align*}
	&\hspace{-0.5em}
	\omega_1^2 C_{\textup{F}}^2 \,
	(\div p_\star^k,\: \div q)_{L^2(\Omega)}
	+
	(p_\star^k - \omega_2^2 \, \nabla u_\star^k,\: q - \omega_2^2 \, \nabla v)_{L^2(\Omega)}
	\\
	&=
	\omega_1^2 C_{\textup{F}}^2 \, (\div p^{k-1},\: \div q)_{L^2(\Omega)}
	+
	(p^{k-1} - \omega_2^2 \, \nabla u^{k-1},\: q - \omega_2^2 \, \nabla v)_{L^2(\Omega)}
	\\
	&\phantom{{}\leq{}}
	-
	\delta \,
	\big[
		\omega_1^2 C_{\textup{F}} \,
		(f_1 + \div p^{k-1},\: \div q)_{L^2(\Omega)} 
		+
		(f_2 + p^{k-1} - \sigma(\nabla u^{k-1}),\: q - \omega_2^2 \, \nabla v)_{L^2(\Omega)}
	\big].
\end{align*}
We highlight that the least-squares method for the inexact solution of the Zarantonello-linearized system of PDEs
takes the form~\eqref{eq:Zarantonello:update} of a Zarantonello iteration itself
with the space \(X \coloneqq H(\div, \Omega) \times H^1_0(\Omega)\).
If it is a contraction (which is confirmed in Corollary~\ref{cor:well_posedness:ls_norm} below),
the iteration converges to the unique solution \((p^\star, u^\star) \in H(\div, \Omega) \times H^1_0(\Omega)\) of the nonlinear equation,
for all \((q, v) \in H(\div, \Omega) \times H^1_0(\Omega)\),
\begin{equation}
	\label{eq:nonlinear:least_squares_weighted}
	\mathcal{B}(p^\star, u^\star; q, v)
	=
	\mathcal{F}(q, v).
\end{equation}
The linear weighting with \(\omega_2\) in the definition~\eqref{eq:nonlinear:least_squares_operator} of \(\mathcal{B}\)
can be interpreted as a remedy for the possibly nonmonotone contribution \(\D\sigma(\nabla u)\)
in formulation~\eqref{eq:nonlinear:least_squares_variational} 
to enable the proofs of strong monotonicity and Lipschitz continuity.
A similar idea is also used in \cite{Riveros2023}.
It is important to note that the formulation~\eqref{eq:nonlinear:least_squares_weighted} is equivalent 
to the nonlinear least-squares problem~\eqref{eq:nonlinear:least_squares_variational}.
\begin{lemma}
	Every solution~\((p^\star, u^\star) \in H(\div, \Omega) \times H^1_0(\Omega)\) to~\eqref{eq:nonlinear:least_squares_weighted}
	also solves the nonlinear least-squares problem~\eqref{eq:nonlinear:least_squares_variational}.
	In particular, it minimizes the nonlinear least-squares functional in~\eqref{eq:nonlinear_minimization}.
\end{lemma}
\begin{proof}
	\emph{Step~1.}
	Given \(\varphi \in L^2(\Omega; \R^d)\), there exist \(q \in H(\div, \Omega)\) with \(\div q = 0\) and \(\widetilde v \in H^1_0(\Omega)\)
	such that \(\varphi = \nabla \widetilde v + q\) \cite[Theorem~I.2.7]{GiraulRaviart1986}.
	The variational formulation~\eqref{eq:nonlinear:least_squares_weighted} 
	tested with \(q\) and \(v = - \omega_2^{-2} \, \widetilde v\) yields
	\[
		0
		=
		(f_2 + p^\star - \sigma(\nabla u^\star),\: q - \omega_2^2 \, \nabla v)_{L^2(\Omega)}
		=
		(f_2 + p^\star - \sigma(\nabla u^\star),\: \varphi)_{L^2(\Omega)}.
	\]
	Since this holds for all \(\varphi \in L^2(\Omega; \R^d)\),
	it follows that \(p^\star = \sigma(\nabla u^\star) - f_2\) in \(L^2(\Omega; \R^d)\).

	\emph{Step~2.}
	Given \(\psi \in L^2(\Omega)\),
	the surjectivity of the weak divergence operator \(\div\colon H(\div, \Omega) \to L^2(\Omega)\)
	implies the existence of \(q \in H(\div, \Omega)\) such that \(\div q = \psi\) \cite[Section~7.1.2]{BoffiBrezziFortin2013}.
	With the equality \(p^\star = \sigma(\nabla u^\star) - f_2\),
	the variational formulation~\eqref{eq:nonlinear:least_squares_weighted} tested with \(q\) and \(v = 0\) proves
	\[
		0
		=
		(f_{1} + \div p^\star,\: \div q)_{L^2(\Omega)}
		=
		(f_{1} + \div p^\star,\: \psi)_{L^2(\Omega)}.
	\]
	Since this holds for all \(\psi \in L^2(\Omega)\),
	it follows that \(- \div p^\star = f_{1}\) in \(L^2(\Omega)\).
	In particular, this verifies~\eqref{eq:nonlinear:least_squares_variational} and concludes the proof.
\end{proof}

For suitable choices of the weights \(\omega_1, \omega_2 > 0\) 
in~\eqref{eq:linearized_LS:functional}, the following result asserts that 
the Zarantonello iteration~\eqref{eq:linearized_LS:iteration} is indeed a contraction.
\begin{theorem}[well-posedness]
	\label{thm:well_posedness}
	Assume that the mapping \(\sigma\colon \R^d \to \R^d\) is differentiable
	and its (not necessarily symmetric) derivative \(\D\sigma\) satisfies~\eqref{assum:ellipticity}--\eqref{assum:boudedness}.
	We use the constants \(\Lambda_1, \Lambda_2 > 0\) from~\eqref{assum:ellipticity} and~\eqref{assum:boudedness}
	to choose the weights \(\omega_1, \omega_2 > 0\) as
	\begin{equation}
		\label{eq:choice_omega_weighting}
		\omega_1^2
		\coloneqq
		\frac{2 \omega_2^2}{\Lambda_1}
		=
		\frac{2 \Lambda_2^2}{\Lambda_1^2}
		\quad\text{and}\quad
		\omega_2^2
		\coloneqq
		\frac{\Lambda_2^2}{\Lambda_1}.
	\end{equation}
	This choice ensures strong monotonicity and Lipschitz continuity of the mapping \(\mathcal{B}\) 
	with respect to the weighted norm \(\vvvert \cdot \vvvert\) from~\eqref{eq:weighted_norm}, i.e.,
	for all \((p, u), (q, v), (r, z) \in H(\div, \Omega) \times H^1_0(\Omega)\),
	\begin{align*}
		\frac{\Lambda_1^2}{4\Lambda_2^2} \,
		\vvvert (p - q, u - v) \vvvert^2
		&\leq
		\mathcal{B}(p, u;p - q, u - v) - \mathcal{B}(q, v;p - q, u - v),
		\\
		\mathcal{B}(p, u; r, z) - \mathcal{B}(q, v; r, z)
		&\leq
		2 \,
		\vvvert (p - q, u - v) \vvvert \,
		\vvvert (r, z) \vvvert.
	\end{align*}
\end{theorem}

\begin{proof}
	\emph{Step~1 (strong monotonicity).}
	For any \(u, v \in H^1_0(\Omega)\), recall the matrix
	\[
		M = \int_0^1 \D\sigma (\nabla (u + s(v - u))) \textup{d}s \in L^\infty(\Omega; \R^{d \times d})
	\]
	from~\eqref{eq:stress_difference}
	satisfying \(\sigma(\nabla u) - \sigma(\nabla v) = M \, \nabla(v - u)\) almost everywhere in \(\Omega\).
	With an integration by parts, this allows rewriting
	\begin{equation}
		\label{eq:monotonicity:split}
		\begin{aligned}
			&\hspace{-2em}
			(p - q - [\sigma(\nabla u) - \sigma(\nabla v)],\: p - q - \omega_2^2 \, \nabla(u - v))_{L^2(\Omega)}
			\\
			&=
			(p - q - M\nabla (u - v),\: p - q - \omega_2^2 \, \nabla(u - v))_{L^2(\Omega)}
			\\
			&=
			\Vert p - q \Vert_{L^2(\Omega)}^2
			+
			\omega_2^2 \,
			(M \nabla(u - v),\: \nabla(u - v) )_{L^2(\Omega)}
			\\
			&\phantom{{}\leq{}}
			+
			\omega_2^2 \,
			(\div (p - q),\: u - v)_{L^2(\Omega)}
			-
			(p - q,\: M \nabla(u - v))_{L^2(\Omega)}.
		\end{aligned}
	\end{equation}
	The ellipticity~\eqref{assum:ellipticity} of \(\D\sigma\) guarantees
	\begin{equation}
		\label{eq:monotonicity:ellipticity}
		\Lambda_1 \, \Vert \nabla(u - v) \Vert_{L^2(\Omega)}^2
		\leq
		(M \nabla(u - v),\: \nabla(u - v) )_{L^2(\Omega)}.
	\end{equation}
	A Cauchy--Schwarz inequality, the Friedrichs inequality, and a weighted Young inequality yield
	\begin{equation}
		\label{eq:monotonicity:Friedrichs}
		\begin{aligned}
			- (\div (p - q),\: u - v)_{L^2(\Omega)}
			&\leq
			C_{\textup{F}} \,
			\Vert \div (p - q) \Vert_{L^2(\Omega)}
			\Vert \nabla (u - v) \Vert_{L^2(\Omega)}.
			\\
			&\leq
			\frac{C_{\textup{F}}^2}{\Lambda_1} \,
			\Vert \div (p - q) \Vert_{L^2(\Omega)}^2
			+
			\frac{\Lambda_1}{4} \, \Vert \nabla (u - v) \Vert_{L^2(\Omega)}^2.
		\end{aligned}
	\end{equation}
	The boundedness~\eqref{assum:boudedness} of \(\D\sigma\), the Cauchy--Schwarz inequality,
	and an unweighted Young inequality show
	\begin{equation}
		\label{eq:monotonicity:boudedness}
		\begin{aligned}
			(p - q,\: M \nabla (u - v))_{L^2(\Omega)}
			&\leq
			\Lambda_2 \,
			\Vert p - q \Vert_{L^2(\Omega)} \,
			\Vert \nabla (u - v) \Vert_{L^2(\Omega)}
			\\
			&\leq
			\frac12 \,
			\Vert p - q \Vert_{L^2(\Omega)}^2
			+
			\frac{\Lambda_2^2}{2} \, \Vert \nabla (u - v) \Vert_{L^2(\Omega)}^2.
		\end{aligned}
	\end{equation}
	The combination of the three previous displayed formulas~\eqref{eq:monotonicity:ellipticity}--\eqref{eq:monotonicity:boudedness} 
	with the initial split~\eqref{eq:monotonicity:split}
	and adding \(C_{\textup{F}}^2 \Vert \omega_1 \div (p - q) \Vert_{L^2(\Omega)}^2\) to both sides results in
	\begin{equation}
		\label{eq:monotonicity:lower_bound}
		\begin{aligned}
			&\hspace{-2em}
			\Big(\omega_1^2 - \frac{\omega_2^2}{\Lambda_1} \Big) C_{\textup{F}}^2\,
			\Vert \div (p - q) \Vert_{L^2(\Omega)}^2
			+
			\frac{1}{2} \, \Vert p - q \Vert_{L^2(\Omega)}^2
			+
			\frac{3\omega_2^2 \Lambda_1 - 2 \Lambda_2^2}4
			\, 
			\Vert \nabla(u - v) \Vert_{L^2(\Omega)}^2
			\\
			&\leq
			C_{\textup{F}}^2 \,
			\Vert \omega_1 \div (p - q) \Vert_{L^2(\Omega)}^2
			+
			(p - q - [\sigma(\nabla u) - \sigma(\nabla v)],\: p - q - \omega_2^2 \, \nabla(u - v))_{L^2(\Omega)}
			\\
			&\eqreff*{eq:nonlinear:least_squares_operator}=\;
			\mathcal{B}(p, u; p - q, u - v) - \mathcal{B}(q, v; p - q, u - v).
		\end{aligned}
	\end{equation}
	The definition of the weighted norm \(\vvvert \cdot \vvvert\) from~\eqref{eq:weighted_norm} and
	the weights from~\eqref{eq:choice_omega_weighting} thus lead to
	\begin{align*}
		&\hspace{-2em}
		\min\bigg\{
			\frac12,\:
			\frac{\Lambda_1}{4\omega_2^2}
		\bigg\} \,
		\vvvert (p - q, u - v) \vvvert^2
		\\
		&=
		\min\bigg\{
			\frac12,\:
			\frac{\Lambda_1}{4\omega_2^2}
		\bigg\} \,
		\bigg[
			\omega_1^2 C_{\textup{F}}^2 \,
			\Vert \div (p - q) \Vert_{L^2(\Omega)}^2
			+
			\Vert p - q \Vert_{L^2(\Omega)}^2
			+
			\Vert \omega_2^2 \nabla(u - v) \Vert_{L^2(\Omega)}^2
		\bigg]
		\\
		&\leq
		\frac{\omega_1^2 C_{\textup{F}}^2}{2} \,
		\Vert \div (p - q) \Vert_{L^2(\Omega)}^2
		+
		\frac12 \, 
		\Vert p - q \Vert_{L^2(\Omega)}^2
		+
		\frac{\omega_2^2 \, \Lambda_1}{4} \, 
		\Vert \nabla(u - v) \Vert_{L^2(\Omega)}^2
		\\
		&\leq
		\mathcal{B}(p, u; p - q, u - v) - \mathcal{B}(q, v; p - q, u - v).
	\end{align*}
	This and the relation \(\Lambda_1 \leq \Lambda_2\) conclude the proof of strong monotonicity with constant
	\[
		\min\bigg\{
			\frac12,\:
			\frac{\Lambda_1}{4\omega_2^2}
		\bigg\}
		\eqreff{eq:choice_omega_weighting}=
		\min\bigg\{
			\frac12,\:
			\frac{\Lambda_1^2}{4\Lambda_2^2}
		\bigg\}
		=
		\frac{\Lambda_1^2}{4\Lambda_2^2}.
	\]

	\emph{Step~2 (Lipschitz continuity).}
	For any \((p, u), (q, v), (r, z) \in H(\div, \Omega) \times H^1_0(\Omega)\),
	an analogous computation as for~\eqref{eq:monotonicity:split} in Step~1 with
	the matrix \(M \in L^\infty(\Omega; \R^{d \times d})\) shows
	\begin{align*}
		\mathcal{B}(p, u; r, z) - \mathcal{B}(q, v; r, z)
		\;&\eqreff*{eq:nonlinear:least_squares_operator}=\;
		\omega_1^2 C_{\textup{F}}^2 \, (\div (p - q),\: \div r)_{L^2(\Omega)}
		+
		(p - q - [\sigma(\nabla u) - \sigma(\nabla v)],\: r - \omega_2^2 \, \nabla z)_{L^2(\Omega)}
		\\
		&=
		\omega_1^2 C_{\textup{F}}^2 \,
		(\div (p - q),\: \div r)_{L^2(\Omega)}
		+
		(p - q,\: r)_{L^2(\Omega)}
		+
		\omega_2^2 \, (M \nabla(u - v),\: \nabla z)_{L^2(\Omega)}
		\\
		&\phantom{{}\leq{}}
		+
		\omega_2^2 \,
		(\div (p - q),\: z)_{L^2(\Omega)}
		-
		(M \nabla(u - v),\: r)_{L^2(\Omega)}.
	\end{align*}
	The boundedness of \(\D\sigma\) from~\eqref{assum:boudedness} and
	the Cauchy--Schwarz inequality in \(L^2(\Omega)\) verify
	\begin{align*}
		&\hspace{-0.5em}
		\mathcal{B}(p, u; r, z) - \mathcal{B}(q, v; r, z)
		\\
		&\leq
		\omega_1^2 C_{\textup{F}}^2 \,
		\Vert \div (p - q) \Vert_{L^2(\Omega)} \, \Vert \div r \Vert_{L^2(\Omega)}
		+
		\Vert p - q \Vert_{L^2(\Omega)} \, \Vert r \Vert_{L^2(\Omega)}
		+
		\omega_2^2 \Lambda_2 \,
		\Vert \nabla(u - v) \Vert_{L^2(\Omega)} \, \Vert \nabla z \Vert_{L^2(\Omega)}
		\\
		&\phantom{{}\leq{}}
		+
		\omega_2^2 C_{\textup{F}} \,
		\Vert \div (p - q) \Vert_{L^2(\Omega)} \, \Vert \nabla z \Vert_{L^2(\Omega)}
		+
		\Lambda_2 \,
		\Vert \nabla(u - v) \Vert_{L^2(\Omega)} \, \Vert r \Vert_{L^2(\Omega)}
		\\
		&=
		C_{\textup{F}}^2 \,
		\Vert \omega_1 \div (p - q) \Vert_{L^2(\Omega)} \, \Vert \omega_1 \div r \Vert_{L^2(\Omega)}
		+
		\Vert p - q \Vert_{L^2(\Omega)} \, \Vert r \Vert_{L^2(\Omega)}
		\\
		&\phantom{{}\leq{}}
		+
		\frac{\Lambda_2}{\omega_2^2} \,
		\Vert \omega_2^2 \, \nabla(u - v) \Vert_{L^2(\Omega)} \, \Vert \omega_2^2 \, \nabla z \Vert_{L^2(\Omega)}
		+
		\frac{1}{\omega_1} C_{\textup{F}} \,
		\Vert \omega_1 \div (p - q) \Vert_{L^2(\Omega)} \, \Vert \omega_2^2 \nabla z \Vert_{L^2(\Omega)}
		\\
		&\phantom{{}\leq{}}
		+
		\frac{\Lambda_2}{\omega_2^2} \,
		\Vert \omega_2^2 \, \nabla(u - v) \Vert_{L^2(\Omega)} \, \Vert r \Vert_{L^2(\Omega)}.
	\end{align*}
	A Cauchy--Schwarz inequality in \(\R^5\) results in
	\begin{align*}
		&\hspace{-1em}
		\mathcal{B}(p, u; r, z) - \mathcal{B}(q, v; r, z)
		\\
		&\leq
		\max\bigg\{
			1,\:
			\frac{\Lambda_2}{\omega_2^2},\:
			\frac{1}{\omega_1}
		\bigg\}
		\Big[
			2 C_{\textup{F}}^2 \,
			\Vert \omega_1 \div (p - q) \Vert_{L^2(\Omega)}^2
			+
			\Vert p - q \Vert_{L^2(\Omega)}^2
			+
			2 \Vert \omega_2^2 \, \nabla(u - v) \Vert_{L^2(\Omega)}^2
		\Big]^{1/2}
		\\
		&\phantom{{}\leq{}}
		\times
		\Big[
			C_{\textup{F}}^2 \,
			\Vert \omega_1 \div r \Vert_{L^2(\Omega)}^2
			+
			2 \,
			\Vert r \Vert_{L^2(\Omega)}^2
			+
			2 \,
			\Vert \omega_2^2 \, \nabla z \Vert_{L^2(\Omega)}^2
		\Big]^{1/2}
		\\
		&\leq
		2
		\max\bigg\{
			1,\:
			\frac{\Lambda_2}{\omega_2^2},\:
			\frac{1}{\omega_1}
		\bigg\}\,
		\vvvert (p - q, u - v) \vvvert \,
		\vvvert (r, z) \vvvert.
	\end{align*}
	This and the relation \(\Lambda_1 \leq \Lambda_2\) conclude
	the proof of the Lipschitz continuity with constant
	\[
		2
		\max\bigg\{
			1,\:
			\frac{\Lambda_2}{\omega_2^2},\:
			\frac{1}{\omega_1}
		\bigg\}\,
		\eqreff{eq:choice_omega_weighting}=
		2
		\max\bigg\{
			1,\:
			\frac{\Lambda_1}{\Lambda_2},\:
			\frac{\Lambda_1}{\sqrt2\Lambda_2}
		\bigg\}
		=
		2.
		\qedhere
	\]
\end{proof}

We emphasize that the direct proof of strong monotonicity and Lipschitz continuity
with respect to the least-squares norm appears impossible,
because the constitutive residual must be split in order to bound the matrix \(M\). 
However, the fundamental equivalence from Theorem~\ref{thm:fundamental_equivalence} 
relates the weighted least-squares norm \(\vvvert \cdot \vvvert_{\mathcal{A}}\) 
with the weighted product norm \(\vvvert \cdot \vvvert\).
Inserting the weights from~\eqref{eq:choice_omega_weighting}
into the equivalence~\eqref{eq:fundamental_equivalence:weighted}
and the estimate \(\Lambda_1 \leq \Lambda_2\) verify
\[
	\min\bigg\{
		\frac12,\:
		\bigg(1 + \frac{2 \Lambda_1^2}{\Lambda_2^2}\bigg)^{-1}
	\bigg\} \,
	\vvvert (q, v) \vvvert^2
	\leq
	\vvvert (q, v) \vvvert_{\mathcal{A}}^2
	\leq
	2 \, \vvvert (q, v) \vvvert^2.
\]
The combination of this and Theorem~\ref{thm:well_posedness} result in the following corollary.
\begin{corollary}
	\label{cor:well_posedness:ls_norm}
	Under the assumptions of Theorem~\ref{thm:well_posedness} and the choice of the weights~\eqref{eq:choice_omega_weighting},
	the nonlinear operator \(\mathcal{B}\) from~\eqref{eq:nonlinear:least_squares_operator}
	is strongly monotone and Lipschitz continuous 
	with respect to the weighted least-squares norm \(\vvvert \cdot \vvvert_{\mathcal{A}}\) from~\eqref{eq:linear:ls_norm},
	i.e.,
	for all \((p, u), (q, v), (r, z) \in H(\div, \Omega) \times H^1_0(\Omega)\),
	\begin{align*}
		\alpha_{\textup{LS}} \,
		\vvvert (p - q, u - v) \vvvert_{\mathcal{A}}^2
		&\leq
		\mathcal{B}(p, u; p - q, u - v) - \mathcal{B}(q, v; p - q, u - v),
		\\
		\mathcal{B}(p, u; r, z) - \mathcal{B}(q, v; r, z)
		&\leq
		L_{\textup{LS}} \,
		\vvvert (p - q, u - v) \vvvert_{\mathcal{A}} \,
		\vvvert (r, z) \vvvert_{\mathcal{A}},
	\end{align*}
	with the constants
	\[
		\alpha_{\textup{LS}}
		\coloneqq
		\frac{\Lambda_1^2}{8\Lambda_2^2} \,
		\quad\text{and}\quad
		L_{\textup{LS}}
		\coloneqq
		2
		\max\bigg\{
			2,\:
			1 + \frac{2 \Lambda_1^2}{\Lambda_2^2}
		\bigg\}.
	\]
	This justifies the application of the Zarantonello iteration from Section~\ref{sec:Zarantonello} in the sense that,
	for all damping parameters \(0 < \delta < \delta_{\textup{LS}}^\star \coloneqq 2\alpha_{\textup{LS}} / L_{\textup{LS}}^2\),
	the iteration~\eqref{eq:linearized_LS:iteration} is a contraction in the norm \(\vvvert \cdot \vvvert_{\mathcal{A}}\).
	In particular, the iterates converge to the solution
	of the first-order optimality condition~\eqref{eq:nonlinear:least_squares_variational}.
	\hfill\qed
\end{corollary}

\begin{remark}
	\label{rem:damping_parameter}
	The estimate \(\Lambda_1 \leq \Lambda_2\) allows bounding the Lipschitz constant 
	in Corollary~\ref{cor:well_posedness:ls_norm} by \(L_{\textup{LS}} \leq 6\).
	Hence, the sufficient damping parameter can be bounded from below by
	\[
		\delta_{\textup{LS}}^\star
		\coloneqq
		\frac{2 \alpha_{\textup{LS}}}{L_{\textup{LS}}^2}
		\geq
		\frac{\Lambda_1^2}{144\Lambda_2^2}.
	\]
	The quality of this parameter depends on the ratio of the problem-dependent constants \(\Lambda_1, \Lambda_2 > 0\)
	from~\eqref{assum:ellipticity} and~\eqref{assum:boudedness}.
	In the case of a small constant \(0 < \Lambda_1 \ll 1\),
	the damping parameter \(\delta_{\textup{LS}}^\star\) scales moderately worse
	than the damping parameter \(\delta^\star = 2\Lambda_1/\Lambda_2^2\)
	for the primal formulation~\eqref{eq:Zarantonello:primal}.
	However, if \(1 \ll \Lambda_1 \leq \Lambda_2\) are large but the ratio \(0 \ll \Lambda_1 / \Lambda_2 \leq 1\) is close to one, 
	then the damping parameter might be even larger (i.e., better) than \(\delta^\star = 2 \Lambda_1 / \Lambda_2^2\)
	from the primal iteration in Section~\ref{sec:model_problem}.
	The following section presents results for other weightings
	in the Zarantonello least-squares functional~\eqref{eq:linearized_LS:functional}.
\end{remark}

%
%
\section{Alternative weightings}
\label{sec:alternative_weightings}
\noindent
Section~\ref{sec:Zarantonello_least_squares} investigated the constitutive residual \(p - \omega_2^2 \, \nabla u\) 
with emphasis on the gradient term by the prefactor \(\omega_2^2\).
It guarantees to find a lower bound for the left-hand side in~\eqref{eq:monotonicity:lower_bound}
in the monotonicity proof for Theorem~\ref{thm:well_posedness}.
However, there is some flexibility in how to realize this lower bound by placing the inverse weight \(\omega_2^{-2}\) 
in front of the flux variable or even split it.

To illustrate that the emphasized gradient weighting results in the most robust constants and damping parameter
(cf.~Remark~\ref{rem:damping_parameter}), this
section presents the following alternative weightings and
the appropriate choices of the (possibly different) weights \(\widetilde\omega_1, \widetilde\omega_2 > 0\)
as well as the resulting monotonicity and Lipschitz constants:
\begin{itemize}
	\item
		Balanced weighting with residual \(\widetilde\omega_2^{-1} \, p - \widetilde\omega_2 \, \nabla u\)
		in Subsection~\ref{sec:balanced_weighting}.
	\item
		Downscaled flux variable with residual \(\widetilde\omega_2^{-1} \, p - \nabla u\)
		in Subsection~\ref{sec:downscaled_flux_variable}.
	\item
		Split weighting with residual \(\Lambda_1 \, p - \Lambda_2^2 \, \nabla u\)
		in Subsection~\ref{sec:split_weighting}.
\end{itemize}
A modification of the weighting affects the Zarantonello least-squares functional \(Z_k\) from~\eqref{eq:linearized_LS:functional} and,
thus, also the mappings \(\mathcal{B}\) and \(\mathcal{F}\) from~\eqref{eq:nonlinear:operator_and_rhs}
and the bilinear form \(\mathcal{A}\) from~\eqref{eq:linear:ls_bilinear_form}.
The mappings and norms with alternative weightings are indicated by a tilde
in order to distinguish them from the rest of the paper.
With a little abuse of notation, they denote different mappings and norms in the following subsections
depending on the choice of the weighting.

The investigation in this section focusses on the monotonicity and Lipschitz constants 
\(\widetilde\alpha_{\textup{LS}}, \widetilde{L}_{\textup{LS}} > 0\) 
of the mapping \(\widetilde{\mathcal{B}}(\,\cdot\,,\,\cdot\,;\,\cdot\,;\,\cdot\,)\colon
[H(\div, \Omega) \times H^1_0(\Omega)] \times [H(\div, \Omega) \times H^1_0(\Omega)] \to \R\)
with respect to the weighted least-squares norm \(\vvvert \cdot \vvvert_{\widetilde{\mathcal{A}}}\) in,
for all \((p, u), (q, v), (r, z) \in H(\div, \Omega) \times H^1_0(\Omega)\),
\begin{equation}
	\label{eq:alternative_weightings:monotonicity_lipschitz}
	\begin{aligned}
		\widetilde\alpha_{\textup{LS}} \,
		\vvvert (p - q, u - v) \vvvert_{\widetilde{\mathcal{A}}}^2
		&\leq
		\widetilde{\mathcal{B}}(p, u; p - q, u - v) - \widetilde{\mathcal{B}}(q, v; p - q, u - v),
		\\
		\widetilde{\mathcal{B}}(p, u; r, z) - \widetilde{\mathcal{B}}(q, v; r, z)
		&\leq
		\widetilde L_{\textup{LS}} \,
		\vvvert (p - q, u - v) \vvvert_{\widetilde{\mathcal{A}}} \,
		\vvvert (r, z) \vvvert_{\widetilde{\mathcal{A}}}.
	\end{aligned}
\end{equation}%
For the sake of a concise presentation, we only state the results 
and refer to the Appendices~\ref{app:balanced_weighting}--\ref{app:split_weighting}
for the detailed proofs.
In comparison with the emphasized-gradient weighting from Section~\ref{sec:Zarantonello_least_squares}, 
all presented alternative weightings in this section exhibit inferior constants in the sense that they
scale worse with respect to the problem-dependent constants \(\Lambda_1, \Lambda_2 > 0\)
from~\eqref{assum:ellipticity} and~\eqref{assum:boudedness} and 
thus tend to lead to a smaller damping parameter \(\widetilde{\delta}_{\textup{LS}}\).
From a theoretical perspective, the former appears to be the most favorable choice among the four considered weightings;
see Section~\ref{sec:applications} for a numerical comparison.

\subsection{Balanced weighting}
\label{sec:balanced_weighting}
\noindent
In this subsection, we discuss the weighted least-squares functional,
for \((p, u), (q, v) \in H(\div, \Omega) \times H^1_0(\Omega)\),
\begin{align*}
	\widetilde Z_k(p, u; q, v)
	&\coloneqq
	\widetilde\omega_1^2 C_{\textup{F}}^2 \,
	\Vert \div (p - p^{k-1}) + \delta \, [f_1 + \div p^{k-1}] \Vert_{L^2(\Omega)}^2
	\\
	&\phantom{{}\coloneqq{}}
	+
	\Vert \widetilde\omega_2^{-1} \, (p - p^{k-1}) - \widetilde\omega_2 \, \nabla (u - u^{k-1}) 
	+ \delta \, [f_2 + p^{k-1} - \sigma(\nabla u^{k-1})] \Vert_{L^2(\Omega)}^2.
\end{align*}
The first variation of the functional \(\widetilde{Z}_k\) leads 
to the following nonlinear mapping \(\widetilde{\mathcal{B}}\) in the formulation~\eqref{eq:nonlinear:least_squares_weighted}
as well as the norms \(\vvvert \cdot \vvvert_{\widetilde{\mathcal{A}}}\) and \(\vvvert \cdot \vvvert_{\widetilde\omega}\),
for all \((p, u), (q, v) \in H(\div, \Omega) \times H^1_0(\Omega)\),
\begin{subequations}
\label{eq:weighting2}
\begin{align}
	\label{eq:weighting2:nonlinear_operator}
	\widetilde{\mathcal{B}}(p, u; q, v)
	&\coloneqq
	\widetilde\omega_1^2 C_{\textup{F}}^2 \,
	(\div q, \div v)_{L^2(\Omega)}
	+
	(p - \sigma(\nabla u), \widetilde\omega_2^{-1} \, q - \widetilde\omega_2 \, \nabla v)_{L^2(\Omega)},
	\\
	\label{eq:weighting2:ls_norm}
	\vvvert (q, v) \vvvert_{\widetilde{\mathcal{A}}}^2
	&\coloneqq
	C_{\textup{F}}^2 \,
	\Vert \widetilde\omega_1 \div q \Vert_{L^2(\Omega)}^2
	+
	\Vert \widetilde\omega_2^{-1} \, q - \widetilde\omega_2 \, \nabla v \Vert_{L^2(\Omega)}^2,
	\\
	\vvvert (q, v) \vvvert_{\widetilde\omega}^2
	&\coloneqq
	C_{\textup{F}}^2 \,
	\Vert \widetilde\omega_1 \div q \Vert_{L^2(\Omega)}^2
	+
	\Vert \widetilde\omega_2^{-1} \, q \Vert_{L^2(\Omega)}^2
	+
	\Vert \widetilde\omega_2 \, \nabla v \Vert_{L^2(\Omega)}^2.
\end{align}
\end{subequations}
The choice of the weights
\begin{equation}
	\label{eq:weighting2:choice_omega}
	\widetilde\omega_1^2
	\coloneqq
	\frac{2 \widetilde\omega_2}{\Lambda_1}
	=
	\frac{2 \Lambda_2}{\Lambda_1^{3/2}}
	>
	0
	\quad\text{and}\quad
	\widetilde\omega_2^2
	\coloneqq
	\frac{\Lambda_2^2}{\Lambda_1}
	>
	0,
\end{equation}
ensures the strong monotonicity and Lipschitz continuity
in~\eqref{eq:alternative_weightings:monotonicity_lipschitz} with the constants
\begin{equation}
	\label{eq:weighting2:constants}
	\widetilde\alpha_{\textup{LS}}
	\coloneqq
	\frac12
	\min\bigg\{
		\frac12,\:
		\frac{\Lambda_2}{\Lambda_1^{1/2}},\:
		\frac{\Lambda_1^{3/2}}{4\Lambda_2}
	\bigg\},
	\quad
	\widetilde L_{\textup{LS}}
	\coloneqq
	2
	\max\bigg\{
		1,\:
		\frac{\Lambda_2}{\Lambda_1^{1/2}},\:
		\frac{\Lambda_1^{3/2}}{2\Lambda_2}
	\bigg\}\,
	\max\bigg\{
		2,\:
		1 + \frac{2 \Lambda_1^{5/2}}{\Lambda_2^3}
	\bigg\}.
\end{equation}
The estimate \(\Lambda_1 \leq \Lambda_2\) leads to the following bounds of the constants
\begin{align*}
	\widetilde\alpha_{\textup{LS}}
	&\geq
	\frac12
	\min\bigg\{
		\frac12,\:
		\Lambda_2^{1/2},\:
		\frac{\Lambda_1^{3/2}}{4\Lambda_2}
	\bigg\},
	\quad
	\widetilde L_{\textup{LS}}
	\leq
	2
	\max\bigg\{
		1,\:
		\frac{\Lambda_2}{\Lambda_1^{1/2}},\:
		\frac{\Lambda_1^{1/2}}{2}
	\bigg\}\,
	\max\bigg\{
		2,\:
		1 + \frac{2}{\Lambda_2^{1/2}}
	\bigg\}.
\end{align*}
In order to compare with the discussion in Remark~\ref{rem:damping_parameter},
the damping parameter scales differently depending on the constants \(\Lambda_1, \Lambda_2 > 0\).
In the first case \(0 < \Lambda_1 < \Lambda_2 \ll 1\), the bounds simplify further to
\[
	\widetilde\alpha_{\textup{LS}}
	\geq
	\frac{\Lambda_1^{3/2}}{8\Lambda_2},
	\quad
	\widetilde L_{\textup{LS}}
	\leq
	6
	\max\bigg\{
		\frac1{\Lambda_2^{1/2}},\:
		\frac{\Lambda_2^{1/2}}{\Lambda_1^{1/2}}
	\bigg\}
	\quad
	\implies
	\quad
	\widetilde\delta_{\textup{LS}}^\star
	\geq 
	\frac{\Lambda_1^{3/2}}{144}
	\min\bigg\{
		1,\:
		\frac{\Lambda_1}{\Lambda_2^{2}}
	\bigg\}.
\]
In the second case \(0 < \Lambda_1 \ll 1 \ll \Lambda_2\), we obtain
\[
	\widetilde\alpha_{\textup{LS}}
	\geq
	\frac{\Lambda_1^{3/2}}{8\Lambda_2},
	\quad
	\widetilde L_{\textup{LS}}
	\leq
	6
	\frac{\Lambda_2}{\Lambda_1^{1/2}}
	\quad
	\implies
	\quad
	\widetilde\delta_{\textup{LS}}^\star
	\geq 
	\frac{\Lambda_1^{5/2}}{144\Lambda_2^3}.
\]
In the remaining case \(1 \ll \Lambda_1 < \Lambda_2\),
\begin{gather*}
	\widetilde\alpha_{\textup{LS}}
	\geq
	\frac14
	\min\bigg\{
		1,\:
		\frac{\Lambda_1^{3/2}}{2\Lambda_2}
	\bigg\},
	\;
	\widetilde L_{\textup{LS}}
	\leq
	\frac{3}{\Lambda_1^{1/2}}
	\max\big\{
		2\Lambda_2,\:
		\Lambda_1
	\big\}
	\;\implies\;
	\widetilde\delta_{\textup{LS}}^\star
	\geq 
	\frac{1}{36\Lambda_1^3}
	\min\bigg\{
		1,\:
		\frac{\Lambda_1^{3/2}}{2\Lambda_2},\:
		\frac{\Lambda_1^{2}}{4\Lambda_2^2},\:
		\frac{\Lambda_1^{7/2}}{8\Lambda_2^3}
	\bigg\}.
\end{gather*}
In the all three cases, the lower bound for the damping parameter \(\widetilde\delta_{\textup{LS}}^\star\)
scales worse than the bound \(\delta_{\textup{LS}}^\star\) from Remark~\ref{rem:damping_parameter}
for the emphasized-gradient weighting in Section~\ref{sec:Zarantonello_least_squares}.

\subsection{Downscaled flux variable}
\label{sec:downscaled_flux_variable}
\noindent
The second alternatively weighted least-squares functional reads,
for \((p, u), (q, v) \in H(\div, \Omega) \times H^1_0(\Omega)\),
\begin{align*}
	\widetilde Z_k(p, u; q, v)
	&\coloneqq
	\widetilde\omega_1^2 C_{\textup{F}}^2 \,
	\Vert \div (p - p^{k-1}) + \delta \, [f_1 + \div p^{k-1}] \Vert_{L^2(\Omega)}^2
	\\
	&\phantom{{}\coloneqq{}}
	+
	\Vert \widetilde\omega_2^{-2} \, (p - p^{k-1}) - \nabla (u - u^{k-1}) 
	+ \delta \, [f_2 + p^{k-1} - \sigma(\nabla u^{k-1})] \Vert_{L^2(\Omega)}^2.
\end{align*}
This functional induces the nonlinear mapping \(\widetilde{\mathcal{B}}\) in the formulation~\eqref{eq:nonlinear:least_squares_weighted}
and the norms \(\vvvert \cdot \vvvert_{\widetilde{\mathcal{A}}}\) and \(\vvvert \cdot \vvvert_{\widetilde\omega}\) 
as follows:
\begin{subequations}
\label{eq:weighting3}
\begin{align}
	\label{eq:weighting3:nonlinear_operator}
	\widetilde{\mathcal{B}}(p, u; q, v)
	&\coloneqq
	\widetilde\omega_1^2 C_{\textup{F}}^2 \,
	(\div q, \div v)_{L^2(\Omega)}
	+
	(p - \sigma(\nabla u), \widetilde\omega_2^{-2} \, q - \nabla v)_{L^2(\Omega)},
	\\
	\label{eq:weighting3:ls_norm}
	\vvvert (q, v) \vvvert_{\widetilde{\mathcal{A}}}^2
	&\coloneqq
	C_{\textup{F}}^2 \,
	\Vert \widetilde\omega_1 \div q \Vert_{L^2(\Omega)}^2
	+
	\Vert \widetilde\omega_2^{-2} \, q - \nabla v \Vert_{L^2(\Omega)}^2,
	\\
	\vvvert (q, v) \vvvert_{\widetilde\omega}^2
	&\coloneqq
	C_{\textup{F}}^2 \,
	\Vert \widetilde\omega_1 \div q \Vert_{L^2(\Omega)}^2
	+
	\Vert \widetilde\omega_2^{-2} \, q \Vert_{L^2(\Omega)}^2
	+
	\Vert \nabla v \Vert_{L^2(\Omega)}^2.
\end{align}
\end{subequations}
For the weights given by
\begin{equation}
	\label{eq:weighting3:choice_omega}
	\widetilde\omega_1^2
	\coloneqq
	\frac{2}{\Lambda_1}
	>
	0
	\quad\text{and}\quad
	\widetilde\omega_2^2
	\coloneqq
	\frac{\Lambda_2^2}{\Lambda_1}
	>
	0,
\end{equation}
the strong monotonicity and Lipschitz continuity in~\eqref{eq:alternative_weightings:monotonicity_lipschitz} hold
with the constants
\begin{equation}
	\label{eq:weighting3:constants}
	\widetilde\alpha_{\textup{LS}}
	\coloneqq
	\frac12
	\min\bigg\{
		\frac12,\:
		\frac{\Lambda_2^2}{2\Lambda_1},\:
		\frac{\Lambda_1}{4}
	\bigg\},
	\quad
	\widetilde L_{\textup{LS}}
	\coloneqq
	2
	\max\bigg\{
		1,\:
		\frac{\Lambda_2^2}{\Lambda_1},\:
		\Lambda_2,\:
		\frac{\Lambda_1^{1/2}}{\sqrt2}
	\bigg\}
	\max\bigg\{
		2,\:
		1 + \frac{2 \Lambda_1^3}{\Lambda_2^4}
	\bigg\} \, 
\end{equation}
With \(\Lambda_1 \leq \Lambda_2\), these constants can be bounded by
\[
	\widetilde\alpha_{\textup{LS}}
	\geq
	\frac12
	\min\bigg\{
		\frac12,\:
		\frac{\Lambda_1}{4}
	\bigg\},
	\quad
	\widetilde L_{\textup{LS}}
	\leq
	2
	\max\bigg\{
		1,\:
		\frac{\Lambda_2^2}{\Lambda_1},\:
		\frac{\Lambda_1^{1/2}}{\sqrt2}
	\bigg\}
	\max\bigg\{
		2,\:
		1 + \frac{2}{\Lambda_2}
	\bigg\}.
\]
To compare with Remark~\ref{rem:damping_parameter},
we investigate the scaling of the damping parameter 
with respect to the constants \(\Lambda_1, \Lambda_2 > 0\).
In the first case \(0 < \Lambda_1 < \Lambda_2 \ll 1\)
(with \(\Lambda_1 \leq 1/2\)), the bounds simplify further to
\[
	\widetilde\alpha_{\textup{LS}}
	\geq
	\frac{\Lambda_1}{8},
	\quad
	\widetilde L_{\textup{LS}}
	\leq
	\frac{6}{\Lambda_2}
	\max\bigg\{
		1,\:
		\frac{\Lambda_2^2}{\Lambda_1}
	\bigg\}
	\quad
	\implies
	\quad
	\widetilde\delta_{\textup{LS}}^\star
	\geq 
	\frac{\Lambda_1\Lambda_2^2}{144}
	\min\bigg\{
		1,\:
		\frac{\Lambda_1^2}{\Lambda_2^4}
	\bigg\}.
\]
In the second case \(0 < \Lambda_1 \ll 1 \ll \Lambda_2\) (with \(\Lambda_2 \geq 2\)), we obtain
\[
	\widetilde\alpha_{\textup{LS}}
	\geq
	\frac{\Lambda_1}{8},
	\quad
	\widetilde L_{\textup{LS}}
	\leq
	\frac{4 \Lambda_2^2}{\Lambda_1}
	\quad
	\implies
	\quad
	\widetilde\delta_{\textup{LS}}^\star
	\geq 
	\frac{\Lambda_1^3}{64 \Lambda_2^4}.
\]
In the remaining case \(1 \ll \Lambda_1 < \Lambda_2\)
(with \(\Lambda_1 \geq 2\)),
it holds that
\[
	\widetilde\alpha_{\textup{LS}}
	\geq
	\frac14,
	\quad
	\widetilde L_{\textup{LS}}
	\leq
	\frac{4\Lambda_2^2}{\Lambda_1}
	\quad
	\implies
	\quad
	\widetilde\delta_{\textup{LS}}^\star
	\geq 
	\frac{\Lambda_1^2}{32\Lambda_2^4}.
\]
Again in all three cases, the lower bound for the damping parameter \(\widetilde\delta_{\textup{LS}}^\star\)
scales worse than the bound \(\delta_{\textup{LS}}^\star\) from Remark~\ref{rem:damping_parameter}
for the emphasized-gradient weighting in Section~\ref{sec:Zarantonello_least_squares}.

\subsection{Split weighting}
\label{sec:split_weighting}
\noindent
In this subsection, we present the weighted least-squares functional,
for \((p, u), (q, v) \in H(\div, \Omega) \times H^1_0(\Omega)\),
\begin{align*}
	\widetilde Z_k(p, u; q, v)
	&\coloneqq
	\widetilde\omega_1^2 C_{\textup{F}}^2 \,
	\Vert \div (p - p^{k-1}) + \delta \, [f_1 + \div p^{k-1}] \Vert_{L^2(\Omega)}^2
	\\
	&\phantom{{}\coloneqq{}}
	+
	\Vert \Lambda_1 \, (p - p^{k-1}) - \Lambda_2^2 \, \nabla (u - u^{k-1}) 
	+ \delta \, [f_2 + p^{k-1} - \sigma(\nabla u^{k-1})] \Vert_{L^2(\Omega)}^2.
\end{align*}
The resulting nonlinear mapping \(\widetilde{\mathcal{B}}\) in the formulation~\eqref{eq:nonlinear:least_squares_weighted}
and norms \(\vvvert \cdot \vvvert_{\widetilde{\mathcal{A}}}\) and \(\vvvert \cdot \vvvert_{\widetilde\omega}\) read
\begin{subequations}
\label{eq:weighting4}
\begin{align}
	\label{eq:weighting4:nonlinear_operator}
	\widetilde{\mathcal{B}}(p, u; q, v)
	&\coloneqq
	\widetilde\omega_1^2 C_{\textup{F}}^2 \,
	(\div q, \div v)_{L^2(\Omega)}
	+
	(p - \sigma(\nabla u), \Lambda_1 \, q - \Lambda_2^2 \, \nabla v)_{L^2(\Omega)},
	\\
	\label{eq:weighting4:ls_norm}
	\vvvert (q, v) \vvvert_{\widetilde{\mathcal{A}}}^2
	&\coloneqq
	C_{\textup{F}}^2 \,
	\Vert \widetilde\omega_1 \div q \Vert_{L^2(\Omega)}^2
	+
	\Vert \Lambda_1 \, q - \Lambda_2^2 \, \nabla v \Vert_{L^2(\Omega)}^2,
	\\
	\vvvert (q, v) \vvvert_{\widetilde\omega}^2
	&\coloneqq
	C_{\textup{F}}^2 \,
	\Vert \widetilde\omega_1 \div q \Vert_{L^2(\Omega)}^2
	+
	\Vert \Lambda_1 \, q \Vert_{L^2(\Omega)}^2
	+
	\Vert \Lambda_2^2 \, \nabla v \Vert_{L^2(\Omega)}^2.
\end{align}
\end{subequations}
The weight
\begin{equation}
	\label{eq:weighting4:choice_omega}
	\widetilde\omega_1^2
	\coloneqq
	\frac{2\Lambda_2^2}{\Lambda_1}
	>
	0
\end{equation}
provides the strong monotonicity and Lipschitz continuity in~\eqref{eq:alternative_weightings:monotonicity_lipschitz}
with constants
\begin{equation}
	\label{eq:weighting4:constants}
	\widetilde\alpha_{\textup{LS}}
	\coloneqq
	\frac12
	\min\bigg\{
		\frac12,\:
		\frac1{2\Lambda_1},\:
		\frac{\Lambda_1}{4\Lambda_2^2}
	\bigg\},
	\quad
	\widetilde L_{\textup{LS}}
	\coloneqq
	2
	\max\bigg\{
		1,\:
		\frac1{\Lambda_1},\:
		\frac1{\Lambda_2},\:
		\frac{\Lambda_1}{2\Lambda_2^2}
	\bigg\}\,
	\max\bigg\{
		2,\:
		1 + \frac{2\Lambda_1^3}{\Lambda_2^2}
	\bigg\}.
\end{equation}
A simplification of these constants with \(\Lambda_1 \leq \Lambda_2\) yields the bounds
\[
	\widetilde\alpha_{\textup{LS}}
	\geq
	\frac12
	\min\bigg\{
		\frac12,\:
		\frac{\Lambda_1}{4\Lambda_2^2}
	\bigg\},
	\quad
	\widetilde L_{\textup{LS}}
	\leq
	2
	\max\bigg\{
		1,\:
		\frac{1}{\Lambda_1}
	\bigg\}\,
	\max\big\{
		2,\:
		1 + 2 \Lambda_1
	\big\}.
\]
For comparison with Remark~\ref{rem:damping_parameter}, 
we again consider the scaling of the damping parameter with respect to \(\Lambda_1, \Lambda_2 > 0\).
In the first case \(0 < \Lambda_1 < \Lambda_2 \ll 1\) (with \(\Lambda_1 \leq 1/2\)), the bounds simplify further to
\[
	\widetilde\alpha_{\textup{LS}}
	\geq
	\frac14
	\min\bigg\{
		1,\:
		\frac{\Lambda_1}{2\Lambda_2^2}
	\bigg\},
	\quad
	\widetilde L_{\textup{LS}}
	\leq
	\frac{4}{\Lambda_1}
	\quad
	\implies
	\quad
	\widetilde\delta_{\textup{LS}}^\star
	\geq 
	\frac{\Lambda_1^2}{32}
	\min\bigg\{
		1,\:
		\frac{\Lambda_1}{2\Lambda_2^2}
	\bigg\}.
\]
In the second case \(0 < \Lambda_1 \ll 1 \ll \Lambda_2\) (with \(\Lambda_1 \leq 1/2\)), it holds that
\[
	\widetilde\alpha_{\textup{LS}}
	\geq
	\frac{\Lambda_1}{8\Lambda_2^2},
	\quad
	\widetilde L_{\textup{LS}}
	\leq
	\frac{4}{\Lambda_1}
	\quad
	\implies
	\quad
	\widetilde\delta_{\textup{LS}}^\star
	\geq 
	\frac{\Lambda_1^3}{64 \Lambda_2^2}.
\]
In the remaining case \(1 \ll \Lambda_1 < \Lambda_2\), we get
\[
	\widetilde\alpha_{\textup{LS}}
	\geq
	\frac{\Lambda_1}{8\Lambda_2^2},
	\quad
	\widetilde L_{\textup{LS}}
	\leq
	6\Lambda_1
	\quad
	\implies
	\quad
	\widetilde\delta_{\textup{LS}}^\star
	\geq 
	\frac{1}{144 \Lambda_1 \Lambda_2^2}.
\]
In all three cases,
the lower bound for the damping parameter \(\widetilde\delta_{\textup{LS}}^\star\)
scales worse than the bound \(\delta_{\textup{LS}}^\star\) from Remark~\ref{rem:damping_parameter}
for the emphasized-gradient weighting in Section~\ref{sec:Zarantonello_least_squares}.

%
%
\section{Adaptive Zarantonello least-squares FEM}
\label{sec:adaptive_Zarantonello_LS}
\noindent
This section is devoted to  step~\ref{step:lsfem} of the procedure described in Subsection~\ref{sec:intro:approach}.
It consists of the discretization of the Zarantonello iteration~\eqref{eq:linearized_LS:iteration} using
the conforming finite element subspaces from \eqref{eq:finite_element_spaces} 
in Subsection~\ref{sec:fespaces} with fixed polynomial degree \(m \in \N_0\).
The underlying meshes are generated using an adaptive refinement strategy 
driven by the built-in least-squares discretization error estimator.

\subsection{Finite element discretization}
Let the optimal damping parameter \(\delta_{\textup{LS}}^\star \coloneqq 2 \alpha_{\textup{LS}}/ L_{\textup{LS}}^2\)
be chosen with the constants from Corollary~\ref{cor:well_posedness:ls_norm}.
For any Zarantonello linearization index \(k \in \N\), we
assume that the previous iterate \((p^{k-1}_H, u^{k-1}_H) \in RT^m(\TT_H) \times S^{m+1}_0(\TT_H)\)
consists of discrete functions on a coarse mesh \(\TT_H \in \T\).
Let \(\TT_h \in \T(\TT_H)\) be a conforming refinement of \(\TT_H\) ensuring the nestedness of the discrete spaces 
\(RT^{m}(\TT_H) \times S^{m+1}_0(\TT_H) \subseteq RT^{m}(\TT_h) \times S^{m+1}_0(\TT_h)\).
Given \(0 < \delta < \delta_{\textup{LS}}^\star\),
the Zarantonello LSFEM seeks the next iterate 
\((p^k_h, u^k_h) \in RT^{m}(\TT_{h}) \times S^{m+1}_0(\TT_{h})\)
such that, for all \((q_h, v_h) \in RT^{m}(\TT_{h}) \times S^{m+1}_0(\TT_{h})\),
\begin{equation}
	\label{eq:linearized_LS:discretization}
	\mathcal{A}(p^k_h, u^k_h; q_h, v_h)
	=
	\mathcal{A}(p^{k-1}_{H}, u^{k-1}_{H}; q_h, v_h)
	+
	\delta \, [ \mathcal{F}(q_h, v_h) - \mathcal{B}(p^{k-1}_{H}, u^{k-1}_{H}; q_h, v_h) ].
\end{equation}
As in the continuous case, this variational formulation characterizes 
the discrete minimizers 
\begin{equation}
	\label{eq:linearized_LS:discrete_minimization}
	Z_k(f_1, f_2; p^k_h, u^k_h)
	=
	\min_{(q_h, v_h) \in RT^m(\TT_h) \times S^{m+1}_0(\TT_h)}
	Z_k(f_1, f_2; q_h, v_h)
\end{equation}
of the least-squares functional \(Z_k\) from~\eqref{eq:linearized_LS:functional}.

\subsection{A~posteriori error estimation}
\label{sec:aposteriori}
Two sources of error need to be controlled in order to provide an upper bound of the overall approximation error: 
the error of the Zarantonello linearization and the discretization error.
The latter one will be measured by the natural built-in least-squares error estimator 
in terms of the least-squares functional for the linearized problem.

For the weights~\eqref{eq:choice_omega_weighting} and the right-hand sides 
\begin{align*}
	g_1^{k-1}
	&\coloneqq
	- \omega_1 \, \div p^{k-1}_{H} 
	+ \delta \omega_1\, [f_1 + \div p^{k-1}_{H}],
	\\
	g_2^{k-1}
	&\coloneqq
	- p^{k-1}_{H} + \omega_2^2 \, \nabla u^{k-1}_{H}
	+ \delta \, [f_2 + p^{k-1}_{H} - \sigma(\nabla u^{k-1}_{H})],
\end{align*}
the functional \(Z_k\) coincides with the weighted least-squares functional \(LS\) defined in~\eqref{eq:linear:least_squares_functional} 
for the linear problem, i.e., for all \((p, u) \in H(\div, \Omega) \times H^1_0(\Omega)\),
\[
	Z_k(f_1, f_2; p, u)
	=
	LS(g_1^{k-1}, g_2^{k-1}; p, u).
\]
The fundamental equivalence~\eqref{eq:fundamental_equivalence:weighted} in Theorem~\ref{thm:fundamental_equivalence}
guarantees that \(\vvvert \cdot \vvvert_{\mathcal{A}}\) is a norm on \(H(\div, \Omega) \times H^1_0(\Omega)\).
For any approximation \((q, v) \in H(\div, \Omega) \times H^1_0(\Omega)\) to the exact minimizers \((p_\star^k, u_\star^k)\)
from~\eqref{eq:linearized_LS:minimization},
the exact built-in error estimate for the discretization error in this least-squares norm reads
\begin{equation}
	\label{eq:zarantonello_ls_built_in_error}
	\vvvert (p^k_\star - q, u^k_\star - v) \vvvert_{\mathcal{A}}^2
	=
	LS(g_1^{k-1}, g_2^{k-1}; q, v)
	=
	Z_k(f_1, f_2; q, v).
\end{equation}
This motivates the definition of the local contributions on any simplex \(T \in \TT_{h}\) by
\begin{equation}
	\label{eq:nonlinear:eta}
	\begin{aligned}
		\eta_k(T; q, v)^2
		&\coloneqq
		\omega_1^2 C_{\textup{F}}^2 \,
		\Vert \div (q - p^{k-1}_{H}) + \delta\, [f_1 + \div p^{k-1}_{H}] \Vert_{L^2(T)}^2
		\\
		&\phantom{{}\coloneqq{}}
		+
		\Vert
		q - p^{k-1}_{H} - \omega_2^2 \, \nabla (v - u^{k-1}_{H}) 
			+
			\delta\, [f_2 + p^{k-1}_{H} - \sigma(\nabla u^{k-1}_{H})]
		\Vert_{L^2(T)}^2.
	\end{aligned}
\end{equation}

The error of contractive linearization schemes is typically measured by the norm of the difference between two consecutive iterates;
cf.~\cite[Lemma~1]{ghps2021}.
This leads to the definition of the local contributions
\begin{equation}
	\label{eq:nonlinear:mu}
	\mu_k(T; q, v)^2
	\coloneqq
	\omega_1^2 C_{\textup{F}}^2 \,
	\Vert \div (q - p^{k-1}_{H}) \Vert_{L^2(T)}^2
	+
	\Vert q - p^{k-1}_{H} - \omega_2^2 \, \nabla (v - u^{k-1}_{H}) \Vert_{L^2(T)}^2.
\end{equation}
In the following, we employ the abbreviations
\begin{align*}
	\eta_k(q, v)^2
	&\coloneqq
	\sum_{T \in \TT_{h}} \eta_k(T; q, v)^2
	=
	Z_k(f_1, f_2; q, v),
	\;
	\mu_k(q, v)^2
	\coloneqq
	\sum_{T \in \TT_{h}} \mu_k(T; q, v)^2
	=
	\vvvert (q - p^{k-1}_{H}, v - u^{k-1}_{H}) \vvvert_{\mathcal{A}}^2.
\end{align*}
The sum of the discretization error estimator \(\eta_k(q, v)\) and the linearization error estimator \(\mu_k(q, v)\)
provides a reliable and efficient error estimator for the total error of the Zarantonello LSFEM.
\begin{proposition}[a posteriori error estimates]
	\label{prop:a_posteriori}
	Let \((p^\star, u^\star) \in H(\div, \Omega) \times H^1_0(\Omega)\) denote the exact solution
	to the nonlinear first-order system~\eqref{eq:nonlinear_fos}.
	For any \((q, v) \in H(\div, \Omega) \times H^1_0(\Omega)\),
	there holds the reliability estimate
	\[
		\vvvert (p^\star - q, u^\star - v) \vvvert_{\mathcal{A}}
		\lesssim
		\eta_k(q, v) + \mu_k(q, v).
	\]
	Moreover, the exact discrete minimizers \((p^k_h, u^k_h) \in RT^{m}(\TT_h) \times S^{m+1}_0(\TT_h)\) of 
	the Zarantonello least-squares functional \(Z_k(f_1, f_2; \,\cdot\,,\,\cdot\,)\) from~\eqref{eq:linearized_LS:discrete_minimization}
	satisfy the efficiency estimate
	\[
		\eta_k(p^k_h, u^k_h) + \mu_k(p^k_h, u^k_h)
		\lesssim
		\vvvert (p^\star - p^{k-1}_{H}, u^\star - u^{k-1}_{H}) \vvvert_{\mathcal{A}}.
	\]
	The hidden constants depend only on the contraction factor \(\rho_{\textup{Z}}\) 
	of the Zarantonello iteration~\eqref{eq:linearized_LS:iteration}; see Corollary~\ref{cor:well_posedness:ls_norm}.
	The fundamental equivalence~\eqref{eq:fundamental_equivalence:weighted} from Theorem~\ref{thm:fundamental_equivalence}
	extends these results to the weighted norm \(\vvvert \cdot \vvvert\) from~\eqref{eq:weighted_norm}.
\end{proposition}

\begin{proof}
	\emph{Step~1 (reliability).}
	For the previous iterate \((p^{k-1}_H, u^{k-1}_H) \in RT^m(\TT_H) \times S^{m+1}_0(\TT_H)\),
	the exact Zarantonello iteration generates \((p^k_\star, u^k_\star) \in H(\div, \Omega) \times H^1_0(\Omega)\)
	by solving the minimization problem~\eqref{eq:linearized_LS:minimization}.
	Theorem~\ref{thm:well_posedness} guarantees existence of the contraction factor \(0 < \rho_{\textup{Z}} < 1\) 
	in the estimate~\eqref{eq:Zarantonello:contraction} which, in the situation at hand, reads as
	\begin{equation}
		\label{eq:a_posteriori:zarantonello_ls_contraction}
		\vvvert (p^\star - p^k_\star, u^\star - u^k_\star) \vvvert_{\mathcal{A}}
		\leq 
		\rho_{\textup{Z}} \,
		\vvvert (p^\star - p^{k-1}_{H}, u^\star - u^{k-1}_{H}) \vvvert_{\mathcal{A}}.
	\end{equation}
	This, two triangle inequalities, and the equality~\eqref{eq:zarantonello_ls_built_in_error} yield,
	for any \((q, v) \in H(\div, \Omega) \times H^1_0(\Omega)\),
	\begin{align*}
		\vvvert (p^\star - q, u^\star - v) \vvvert_{\mathcal{A}}
		&\leq
		\vvvert (p^\star - p^k_\star, u^\star - u^k_\star) \vvvert_{\mathcal{A}}
		+
		\vvvert (p^k_\star - q, u^k_\star - v) \vvvert_{\mathcal{A}}
		\\
		&\eqreff*{eq:a_posteriori:zarantonello_ls_contraction}\leq
		\rho_{\textup{Z}} \,
		\vvvert (p^\star - p^{k-1}_{H}, u^\star - u^{k-1}_{H}) \vvvert_{\mathcal{A}}
		+
		\vvvert (p^k_\star - q, u^k_\star - v) \vvvert_{\mathcal{A}}
		\\
		&\leq
		\rho_{\textup{Z}} \,
		\vvvert (p^\star - q, u^\star - v) \vvvert_{\mathcal{A}}
		+
		\rho_{\textup{Z}} \,
		\vvvert (q - p^{k-1}_{H}, v - u^{k-1}_{H}) \vvvert_{\mathcal{A}}
		+
		\vvvert (p^k_\star - q, u^k_\star - v) \vvvert_{\mathcal{A}}
		\\
		&\eqreff*{eq:zarantonello_ls_built_in_error}{=}
		\rho_{\textup{Z}} \,
		\vvvert (p^\star - q, u^\star - v) \vvvert_{\mathcal{A}}
		+
		\rho_{\textup{Z}} \,
		\mu_k(q, v)
		+
		\eta_k(q, v).
	\end{align*}
	The absorption of the first summand on the right-hand side concludes the proof of the reliability estimate
	\[
		\vvvert (p^\star - q, u^\star - v) \vvvert_{\mathcal{A}}
		\leq
		\frac{1}{1 - \rho_{\textup{Z}}} \,
		\big[
		\eta_k(q, v) 
		+
		\mu_k(q, v) \big].
	\]

	\emph{Step~2 (efficiency).}
	The exact solution \((p^k_h, u^k_h) \in RT^{m}(\TT_{h}) \times S^{m+1}_0(\TT_{h})\)
	to the discrete problem~\eqref{eq:linearized_LS:discretization} 
	solves the minimization problem~\eqref{eq:linearized_LS:discrete_minimization} of
	the Zarantonello least-squares functional \(Z_k(f_1, f_2;\,\cdot\,,\,\cdot\,)\)
	over the spaces \(RT^{m}(\TT_{h}) \times S^{m+1}_0(\TT_{h})\).
	Since the nestedness of the discrete spaces ensures that
	\((p^{k-1}_{H}, u^{k-1}_{H}) \in RT^m(\TT_H) \times S^{m+1}_0(\TT_H) \subseteq RT^{m}(\TT_h) \times S^{m+1}_0(\TT_h)\),
	this implies with~\eqref{eq:zarantonello_ls_built_in_error} that
	\begin{equation}
		\label{eq:a_posteriori:nestedness}
		\eta_k(p^k_h, u^k_h)
		=
		Z_k(f_1, f_2; p^k_h, u^k_h)^{1/2}
		\leq
		Z_k(f_1, f_2; p^{k-1}_{H}, u^{k-1}_{H})^{1/2}
		\eqreff{eq:zarantonello_ls_built_in_error}=
		\vvvert (p^k_\star - p^{k-1}_{H}, u^k_\star - u^{k-1}_{H}) \vvvert_{\mathcal{A}}.
	\end{equation}
	This, a triangle inequality, and the estimate~\eqref{eq:a_posteriori:zarantonello_ls_contraction} yield
	\begin{align*}
		\eta_k(p^k_h, u^k_h)
		&\eqreff*{eq:a_posteriori:nestedness}{\leq}
		\vvvert (p^k_\star - p^{k-1}_{H}, u^k_\star - u^{k-1}_{H}) \vvvert_{\mathcal{A}}
		\\
		&\leq
		\vvvert (p^\star - p^k_\star, u^\star - u^k_\star) \vvvert_{\mathcal{A}}
		+
		\vvvert (p^\star - p^{k-1}_{H}, u^\star - u^{k-1}_{H}) \vvvert_{\mathcal{A}}
		\\
		&\eqreff*{eq:a_posteriori:zarantonello_ls_contraction}\leq
		(1 + \rho_{\textup{Z}})
		\vvvert (p^\star - p^{k-1}_{H}, u^\star - u^{k-1}_{H}) \vvvert_{\mathcal{A}}.
	\end{align*}
	The same arguments establish
	\begin{align*}
		\mu_k(p^k_h, u^k_h)
		&=
		\vvvert (p^k_h - p^{k-1}_{H}, u^k_h - u^{k-1}_{H}) \vvvert_{\mathcal{A}}
		\\
		&\leq
		\vvvert (p^\star - p^{k-1}_{H}, u^\star - u^{k-1}_{H}) \vvvert_{\mathcal{A}}
		+
		\vvvert (p^\star - p^k_\star, u^\star - u^k_\star) \vvvert_{\mathcal{A}}
		+
		\vvvert (p^k_\star - p^k_h, u^k_\star - u^k_h) \vvvert_{\mathcal{A}}
		\\
		&\eqreff*{eq:a_posteriori:zarantonello_ls_contraction}\leq
		2(1 + \rho_{\textup{Z}}) \,
		\vvvert (p^\star - p^{k-1}_{H}, u^\star - u^{k-1}_{H}) \vvvert_{\mathcal{A}}.
	\end{align*}
	The sum of the two previous displayed formulas concludes the proof of the efficiency estimate
	\[
		\eta_k(p^k_h, u^k_h) + \mu_k(p^k_h, u^k_h)
		\leq
		3(1 + \rho_{\textup{Z}}) \,
		\vvvert (p^\star - p^{k-1}_{H}, u^\star - u^{k-1}_{H}) \vvvert_{\mathcal{A}}.
		\qedhere
	\]
\end{proof}

In combination with Theorem~\ref{thm:fundamental_equivalence} and Proposition~\ref{prop:nonlinear:aposteriori},
we immediately deduce the following equivalence.
\begin{corollary}
	\label{cor:aposteriori}
	Using the notation of Proposition~\ref{prop:a_posteriori}, it holds that
	\[
		\vvvert (p^\star - p^k_h, u^\star - u^k_h) \vvvert_{\textup{uw}}
		\eqsim
		N(f_1, f_2; p^k_h, u^k_h)^{1/2}
		\lesssim
		\eta_k(p^k_h, u^k_h) + \mu_k(p^k_h, u^k_h).
		\qquad\qed
	\]
\end{corollary}
This justifies that both, the nonlinear least-squares functional \(N(f_1, f_2; \,\cdot\,;\,\cdot\,)\) and 
the sum \(\eta_k(\,\cdot\,,\,\cdot\,) + \mu_k(\,\cdot\,,\,\cdot\,)\), are upper bounds for the error of the approximation
in the unweighted norm \(\vvvert \,\cdot\, \vvvert_{\textup{uw}}\) on \(H(\div, \Omega) \times H^1_0(\Omega)\)
from~\eqref{eq:unweighted_norm}.
However, we consider the estimator \(\eta_k(\,\cdot\,,\,\cdot\,) + \mu_k(\,\cdot\,,\,\cdot\,)\) 
from Proposition~\ref{prop:a_posteriori} as advantageous because
it follows without assuming the symmetry \(\D\sigma = \D\sigma^\top\).
Moreover, it contains the built-in discretization error estimator \(\eta_k\)
which can be computed without an additional (possibly expensive) quadrature of the nonlinear least-squares functional.
In the concluding Section~\ref{sec:applications} below,
we will compare both error measures as part of the numerical investigation.
To summarize, the various quantities are considered as measures for the following error contributions arising in the corresponding steps of the
procedure described in Subsection~\ref{sec:intro:approach}:
\begin{itemize}
	\item
		\(\eta_k(\,\cdot\,,\,\cdot\,)\): Discretization error in step~\ref{step:lsfem}.
	\item
		\(\mu_k(\,\cdot\,,\,\cdot\,)\): Linearization error in step~\ref{step:zarantonello}.
	\item
		\(N(f_1, f_2; \,\cdot\,;\,\cdot\,)\): Overall approximation error.
\end{itemize}

\subsection{Adaptive algorithm}
The key idea for the adaptive mesh-refinement algorithm is to apply the established adaptive LSFEM driven by the built-in error estimator
\cite{CarstensenParkBringmann2017,FuhrerPraetorius2020,GantnerStevenson2021,Bringmann2024}
for the solution of the linearized system of PDEs in step~\ref{step:lsfem} of Subsection~\ref{sec:intro:approach}.
The inner adaptive algorithm is stopped as soon as the desired accuracy for the discrete solution is reached.
In order to determine this point and to steer the adaptive mesh refinement,
we employ the estimator \(\eta_k(\,\cdot\,,\,\cdot\,)\) for the discretization error
of the discrete problem~\eqref{eq:linearized_LS:discretization}.
The simplices with large estimator contributions are selected by the D\"orfler marking criterion \cite[Section~5]{Dorfler1996}.
This results in the combined Algorithm~\ref{alg:adaptive} for linearization and adaptive LSFEM.
Note that the linearization error estimator \(\mu_k(\,\cdot\,,\,\cdot\,)\) is not employed therein
as the convergence of the outer linearization loop is fully guaranteed by the Zarantonello iteration;
see Theorem~\ref{thm:global_convergence} below.
\begin{algorithm}[t]
	\caption{Adaptive Zarantonello least-squares FEM}
	\label{alg:adaptive}
	\flushleft
	\noindent
	\textbf{Input:}
	Initial mesh \(\TT_0^1 \coloneqq \TT_0\),
	initial iterates \(p^0_0 \coloneqq p^0_{\underline\ell} \in RT^{m}(\TT_0^1)\)
	and \(u^0_0 \coloneqq u^0_{\underline\ell} \in S^{m+1}_0(\TT_0^1)\),
	marking parameter \(0 < \theta \leq 1\),
	stopping parameter \(0 < \gamma < 1\).
	\medskip

	\noindent
	\textbf{for all} $k = 1, 2, 3, \dots$ \textbf{do}
	\hfill \textit{\% linearization loop}
	\begin{enumerate}[label=(\roman*),leftmargin=2em]
		\item 
			\textbf{for all} $\ell = 0, 1, 2, \dots$ \textbf{do}
			\hfill \textit{\% refinement loop}
			\begin{enumerate}[label=(\alph*),leftmargin=2em]
				\item 
					\textbf{Solve.} 
					Compute the exact solution 
					\((p^k_\ell, u^k_\ell) \in RT^{m}(\TT^k_\ell) \times S^{m+1}_0(\TT^k_\ell)\)
					to the discrete linear problem~\eqref{eq:linearized_LS:discretization} with respect to
					\(p^{k-1}_{H} = p^{k-1}_{\underline\ell}\) and \(u^{k-1}_{H} = u^{k-1}_{\underline\ell}\).
				\item 
					\textbf{Estimate.}
					Compute \(\eta_k(T; p^k_\ell, u^k_\ell)^2\) from~\eqref{eq:nonlinear:eta} for all \(T \in \TT^k_\ell\).
				\item
					\textbf{If} \(\eta_k(p^k_\ell, u^k_\ell) \leq \gamma^k\),
					\textbf{then break} the \(\ell\) loop.
					\hfill \textit{\% stopping criterion}
				\item 
					\textbf{Mark.}
					Determine a minimal subset \(\MM^k_\ell \subseteq \TT^k_\ell\) such that 
					\[
						\theta \, \eta_k(p^k_\ell, u^k_\ell)^2
						\leq
						\sum_{T \in \MM^k_\ell} \eta_k(T; p^k_\ell, u^k_\ell)^2
					\]
				\item 
					\textbf{Refine.}
					Generate refined mesh $\TT^k_{\ell+1} \coloneqq \refine(\TT^k_\ell, \MM^k_\ell)$ by NVB.
			\end{enumerate}
		\item[] \textbf{end for}
		\item
			Define \(\underline\ell[k] \coloneqq \ell\),
			\(\TT^{k+1}_0 \coloneqq \TT^k_{\underline\ell} \coloneqq \TT^k_\ell\),
			\(p^k_{\underline\ell} \coloneqq p^k_\ell\), and \(u^k_{\underline\ell} \coloneqq u^k_\ell\)
			(nested iteration).
		\end{enumerate}
		\textbf{end for}
		\medskip

		\noindent
		\textbf{Output:} Sequentially ordered meshes \(\TT^k_\ell\) with corresponding discrete functions
		\((p^k_\ell, u^k_\ell) \in RT^{m}(\TT^k_\ell) \times S^{m+1}_0(\TT^k_\ell)\).
\end{algorithm}

For a clear presentation, Algorithm~\ref{alg:adaptive} is formulated with a double index.
The upper index \(k\) refers to the outer loop of the Zarantonello iteration and 
the lower index \(\ell\) to the inner loop performing the adaptive mesh refinement.
The latter is restarted for every step of the linearization loop.
The final mesh index \(\underline\ell[k]\) depends on the linearization index \(k \in \N\),
but this dependency is omitted in the notation whenever it is clear from the context, e.g., 
for \(\TT_{\underline\ell}^k\) replacing \(\TT_{\underline\ell[k]}^k\) and
analogously for \(p_{\underline\ell}^k\) and \(u_{\underline\ell}^k\).
Nevertheless, Algorithm~\ref{alg:adaptive} can be realized with a single index in practice.
This leads to the sequentially ordered index set
\[
	\mathcal{Q}
	\coloneqq
	\big\{
		(k, \ell) \in \N_0^2
		\;:\;
		\TT^k_\ell \text{ is output by Algorithm~\ref{alg:adaptive} }
	\big\}.
\]

The convergence of the adaptive LSFEM from Theorem~\ref{thm:ls_plain_convergence} ensures
that the inner \(\ell\)-loop always terminates.
The accepted solutions converge R-linearly as stated in the following main result.
\begin{theorem}[global R-linear convergence]
	\label{thm:global_convergence}
	The sequence \((p^k_{\underline\ell}, u^k_{\underline\ell})_{k \in \N_0}\) of final iterates of the inner mesh-refinement loop
	of Algorithm~\ref{alg:adaptive} converges R-linearly to
	the exact solution \((p^\star, u^\star) \in H(\div, \Omega) \times H^1_0(\Omega)\),
	i.e., there exists constants \(C_{\textup{lin}} > 0\) and \(0 < \rho < 1\) such that
	\[
		\vvvert (p^\star - p^k_{\underline\ell}, u^\star - u^k_{\underline\ell}) \vvvert_{\mathcal{A}}
		\leq
		C_{\textup{lin}}\,
		\rho^k.
	\]
\end{theorem}

\begin{proof}
	\emph{Step~1.}\,
	With the contraction factor \(\rho_{\textup{Z}} < 1\) from~\eqref{eq:Zarantonello:contraction}
	and the parameter \(\gamma < 1\) from Algorithm~\ref{alg:adaptive},
	let \(0 < \rho_\star \coloneqq \max\{\rho_{\textup{Z}}, \gamma \} < 1\) and choose \(\rho_\star < \rho < 1\).
	Let \(k_\star \in \N\) denote the smallest integer such that \(k_\star \leq (\rho / \rho_\star)^{k_\star}\).
	Hence, for any \(k \in \N\), it holds either that \(k < k_\star\) or that \(k (\rho_\star / \rho)^k \leq 1\).
	This ensures, for \(C_\star \coloneqq k_\star (\rho_\star / \rho)^{k_\star} > 0\) and all \(k \in \N\),
	\begin{equation}
		\label{eq:contraction_factor}
		k \, \rho_\star^k
		=
		k\, \frac{\rho_\star^k}{\rho^k} \,
		\rho^k
		\leq
		\begin{cases}
			C_\star \, \rho^k & \text{if } k < k_\star, \\
			\rho^k & \text{if } k \geq k_\star.
		\end{cases}
	\end{equation}
	
	\emph{Step~2.}\,
	The contraction~\eqref{eq:Zarantonello:contraction} of the exact Zarantonello iteration reads
	\[
		\vvvert (p^\star - p^k_\star, u^\star - u^k_\star) \vvvert_{\mathcal{A}}
		\leq
		\rho_{\textup{Z}} \,
		\vvvert (p^\star - p^{k-1}_{\underline\ell}, u^\star - u^{k-1}_{\underline\ell}) \vvvert_{\mathcal{A}}.
	\]
	The error equality~\eqref{eq:zarantonello_ls_built_in_error} of the least-squares functional
	and the stopping criterion~(c) of the adaptive mesh-refinement loop
	in Algorithm~\ref{alg:adaptive} imply
	\[
		\vvvert (p^k_\star - p^k_{\underline\ell}, u^k_\star - u^k_{\underline\ell}) \vvvert_{\mathcal{A}}
		=
		\eta_k(p^k_{\underline\ell}, u^k_{\underline\ell})
		\leq
		\gamma^k.
	\]
	The combination of the two previous inequalities with a triangle inequality proves
	\begin{align*}
		\vvvert (p^\star - p^k_{\underline\ell}, u^\star - u^k_{\underline\ell}) \vvvert_{\mathcal{A}}
		&\leq
		\vvvert (p^\star - p^k_\star, u^\star - u^k_\star) \vvvert_{\mathcal{A}}
		+
		\vvvert (p^k_\star - p^k_{\underline\ell}, u^k_\star - u^k_{\underline\ell}) \vvvert_{\mathcal{A}}
		\\
		&\leq
		\rho_{\textup{Z}} \,
		\vvvert (p^\star - p^{k-1}_{\underline\ell}, u^\star - u^{k-1}_{\underline\ell}) \vvvert_{\mathcal{A}}
		+
		\gamma^k.
	\end{align*}
	By induction on \(k \in \N\), this results in
	\begin{align*}
		\vvvert (p^\star - p^k_{\underline\ell}, u^\star - u^k_{\underline\ell}) \vvvert_{\mathcal{A}}
		&\leq
		\rho_{\textup{Z}}^k \,
		\vvvert (p^\star - p^0_0, u^\star - u^0_0) \vvvert_{\mathcal{A}}
		+
		\sum_{j=0}^{k-1} \rho_{\textup{Z}}^j \, \gamma^{k-j}
		\leq
		\rho^k \,
		\vvvert (p^\star - p^0_0, u^\star - u^0_0) \vvvert_{\mathcal{A}}
		+
		k \, \rho_\star^k.
	\end{align*}
	This and the estimate~\eqref{eq:contraction_factor} with the constant \(C_\star = k_\star (\rho_\star/\rho)^{k_\star}\)
	conclude the proof with the generic constant
	\[
		C_{\textup{lin}}
		\coloneqq
		\vvvert (p^\star - p^0_0, u^\star - u^0_0) \vvvert_{\mathcal{A}}
		+
		C_\star.
		\qedhere
	\]
\end{proof}

%
%
\section{Applications}
\label{sec:applications}
\noindent
This section is devoted to the application of the adaptive Zarantonello LSFEM of Algorithm~\ref{alg:adaptive}
to some quasilinear PDEs.
We discuss the applicability of the theory and investigate the performance of the adaptive method in numerical computations.
The implementation is based on the octAFEM software package that was also used in~\cite{Bringmann2024} for adaptive linear LSFEMs.
The complete code for reproducing the numerical experiments 
is published as a code capsule on the Code Ocean platform \cite{BringmannPraetorius2025_code}.

\subsection{Convex energy minimization}
\label{sec:convex_energy_minimization}
\noindent
Problems of the form of the model problem~\eqref{eq:nonlinear} typically arise
in the context of the minimization of convex energy functionals.
This subsection presents a general framework from~\cite[Chapter~25]{Zeidler1990} for such type of problems.
Some practical applications are described in~\cite[Section~2.2]{gmz2012}.
Given a function \(\phi \in C^2(0, \infty)\) with
\begin{equation}
	\label{eq:alc:assumptions}
	\Lambda_1 \leq \phi(t) \leq \Lambda_2
	\quad \text{and} \quad
	\Lambda_1 \leq \phi(t) + t \phi'(t) \leq \Lambda_2
	\quad \text{for all } t > 0,
\end{equation}
define the convex potential function \(\Phi(t) \coloneqq \int_0^t s\,\phi(s) \,\textrm{d}s\) for \(t \geq 0\).
The minimizer \(u^\star \in H^1_0(\Omega)\) of the energy functional
\[
	\mathcal{E}(u)
	\coloneqq
	\int_\Omega
	\Phi(\vert \nabla u \vert) \d{x}
	-
	\int_\Omega (f_1 u + f_2 \cdot \nabla u) \d{x}
\]
is characterized by the Euler--Lagrange equation~\eqref{eq:nonlinear} 
including the nonlinear mapping \(\sigma\colon \R^d \to \R^d\) with \(\sigma(\xi) = \phi(\vert \xi \vert) \xi\).
The function \(\sign\colon \R^d \to \R^d\) with \(\sign(\xi) \coloneqq \xi / \vert \xi \vert\)
for \(\xi \in \R^d \setminus \{0\}\) allows calculating the derivative
\[
	\D \sigma(\xi) 
	=
	\phi(\vert \xi \vert) \, I_{d \times d} 
	+
	\phi'(\vert \xi \vert) \, \vert \xi \vert \sign(\xi) \otimes \sign(\xi).
\]
In the experiment below, the prefactor \(\phi'(\vert \xi \vert) \, \vert \xi \vert\) makes \(\D\sigma\) a continuous function in the full \(\R^d\).
Hence, the assumptions~\eqref{eq:alc:assumptions} on the function \(\phi\) guarantee
that \(\sigma \in C^1(\R^d; \R^d)\) and that \(\D\sigma \in C^0(\R^d; \R^{d \times d}_{\textup{sym}})\) 
satisfies the conditions~\eqref{assum:ellipticity}--\eqref{assum:boudedness}.
The reader is referred to \cite[Section~3]{CarstensenBringmannHellwigWriggers2018} and \cite[Section~4]{BringmannCarstensenTran2022}
for a more detailed discussion of the application of minimal residual methods to this problem class.

As a benchmark example for the Zarantonello LSFEM from Section~\ref{sec:Zarantonello_least_squares},
we consider the coefficient function \(\phi \in C^\infty(0, \infty)\) with
\(\phi(t) \coloneqq 2 + (1 + t)^{-1}\) from~\cite[Section~VII]{CarstensenStephan1995}.
This function and its derivative \(\phi'(t) = - (1 + t)^{-2}\) satisfy,
for all \(t > 0\),
\[
	2 \leq \phi(t) \leq 3
	\quad\text{and}\quad
	2
	\leq
	\phi(t) + t \phi'(t)
	=
	2
	+
	\frac{1}{(1 + t)^2}
	\leq
	3
\]
verifying the assumptions~\eqref{eq:alc:assumptions} for the constants \(\Lambda_1 = 2\) and \(\Lambda_2 = 3\).
The choice~\eqref{eq:choice_omega_weighting} of the weights \(\omega_1, \omega_2 > 0\)
for the least-squares functional~\eqref{eq:linearized_LS:functional} reads
\[
	\omega_1^2
	=
	\frac{2\Lambda_2^2}{\Lambda_1^2}
	=
	\frac{9}{2}
	=
	\frac{\Lambda_2^2}{\Lambda_1}
	=
	\omega_2^2.
\]
Let \(\Omega \coloneqq (-1, 1)^2 \setminus [0, 1)^2\)
be the L-shaped domain with approximated Friedrichs constant
\[
	\lambda_1^{-1/2} \leq C_{\textup{F}} \coloneqq 0.32208292665417854
\]
computed from guaranteed lower bounds for the first Dirichlet eigenvalue \(\lambda_1 \geq 9.639723838973880\)
of the Laplace operator; see, e.g., \cite[Section~6.3]{CarstensenGedicke2014} and \cite[Section~7.2]{CarstensenErnPuttkammer2021}.
Moreover, let \(f_1 \equiv 1 \in L^2(\Omega)\) and
\(f_2 \equiv 0 \in L^2(\Omega; \R^2)\)
be the right-hand side of the Euler--Lagrange equation~\eqref{eq:nonlinear}.
The following experiments consider lowest-order discretizations with polynomial degree \(m = 0\).

For one run of Algorithm~\ref{alg:adaptive} with parameters \(\delta = 1\), \(\gamma = 0.9\), and \(\theta = 0.3\),
the convergence history plot in Figure~\ref{fig:nonlinear:errors} compares the three error quantities
from~\eqref{eq:nonlinear:eta}, \eqref{eq:nonlinear:mu}, and~\eqref{eq:nonlinear_ls_functional}
(cf.\ the a~posteriori estimates from Propositions~\ref{prop:nonlinear:aposteriori} and~\ref{prop:a_posteriori}
and Corollary~\ref{cor:aposteriori}) 
with the abbreviations, for all \((\ell, k) \in \mathcal{Q}\),
\begin{equation}
	\label{eq:applications:errors}
	\eta^k_\ell \coloneqq \eta_k(p^k_\ell, u^k_\ell),
	\quad
	\mu^k_\ell \coloneqq \mu_k(p^k_\ell, u^k_\ell),
	\quad \text{and} \quad
	N^k_\ell \coloneqq N(f_1, f_2; p^k_\ell, u^k_\ell)^{1/2}.
\end{equation}
Due to the cumulative nature of Algorithm~\ref{alg:adaptive},
these quantities are plotted against the cumulative number of degrees of freedom defined by
\begin{equation}
	\label{eq:cumulative_ndof}
	\sum_{(k, \ell) \in \mathcal{Q}}
	\big[
		\dim(RT^0(\TT^k_\ell))
		+
		\dim(S^1_0(\TT^k_\ell))
	\big].
\end{equation}
This value roughly represents the overall computational effort
to determine the approximation \((p^k_\ell, u^k_\ell) \in RT^0(\TT^k_\ell) \times S^1_0(\TT^k_\ell)\)
because all intermediate meshes need to be computed as well to obtain it.
Note that~\eqref{eq:cumulative_ndof} does only reflect the overall computational cost under the idealistic assumption
that the arising linear systems are solved in linear complexity, which, based on an inexact iterative solver,
is not the focus of the present work.

The experiment confirms that the linearized least-squares estimator \(\eta^k_\ell\)
and the nonlinear least-squares estimator \(N^k_\ell\) match very well and essentially measure the discretization error
whereas the linearization error \(\mu^k_\ell\) is smaller and gets significantly reduced in every Zarantonello update 
(nearly vertical steps of the graphs).
Figure~\ref{fig:nonlinear:mesh} shows the adaptively generated mesh \(\TT^k_{1}\) of the Zarantonello iterate \(k = 29\) 
from this computation.
It exhibits an increased refinement towards the re-entrant corner of the L-shaped domain.
The discrete solution \((p^{29}_{1}, u^{29}_{1}) \in RT^0(\TT^{29}_{1}) \times S^1_{0}(\TT^{29}_{1})\)
from the same experiment are depicted in Figures~\ref{fig:nonlinear:potential} and~\ref{fig:nonlinear:flux}.
For improved visualization, the discrete flux variable \(p^{29}_{1}\)
is only evaluated at 208 equidistributed points in the domain \(\Omega\).

The convergence plot in Figure~\ref{fig:nonlinear:convergence_history:gamma}
displays the full error estimator \(\eta^k_\ell + \mu^k_\ell\)
for various choices of the reduction parameter \(0 < \gamma < 1\).
The figure shows that smaller values of \(\gamma\) emphasize the focus on the mesh refinement
to the detriment of the linearization error.
For \(\gamma = 0.1\), the algorithm performs almost only mesh refinement steps (nearly horizontal steps of the graphs).
The results are comparable for both parameter selections \(\delta \in \{0.5, 1\}\).
Consequently, comparably large values like \(\gamma = 0.9\) are advisable.

Figure~\ref{fig:nonlinear:convergence_history:theta} investigates the influence of the bulk parameter \(0 < \theta \leq 1\).
Numerical experiments, e.g., in \cite[Section~4.3]{BringmannCarstensenTran2022}, \cite[Section~6.2]{ghps2021}, 
and \cite[Section~7]{HaberlPraetoriusSchimankoVohralik2021}
on adaptive methods for related examples suggest an expected optimal convergence rate of \(0.5\)
with respect to both, the number of degrees of freedom and the cumulative number of degrees of freedom~\eqref{eq:cumulative_ndof}.
As usual for adaptive mesh-refinement algorithms, the expected optimal rate is achieved 
for small and moderate choices of \(0.3 \leq \theta \leq 0.7\).
However, smaller values of \(\theta\) lead to a more adaptive behavior in that
there are more mesh refinement steps with fewer linearization steps in between.
The uniform refinement with \(\theta = 1\) and adaptive refinement with large bulk parameter \(\theta = 0.9\)
exhibit suboptimal convergence rates.
Both employed damping parameters \(\delta \in \{0.5, 1\}\) result in a similar behavior.

\begin{figure}
	\centering
	\input{figures/Fig04_Nonlinear_convergence_damping.tex}
	\caption{%
		Convergence history plot for Algorithm~\ref{alg:adaptive} 
		applied to the convex energy minimization problem from Subsection~\ref{sec:convex_energy_minimization}
		and various choices of the damping parameter \(0 < \delta \leq 1\).
		The remaining parameters read \(\theta = 0.3\) and \(\gamma = 0.9\).
	}
	\label{fig:nonlinear:convergence_history:damping}
\end{figure}
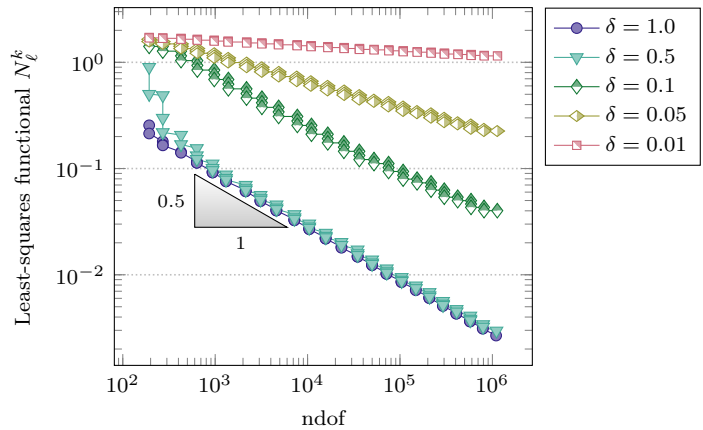

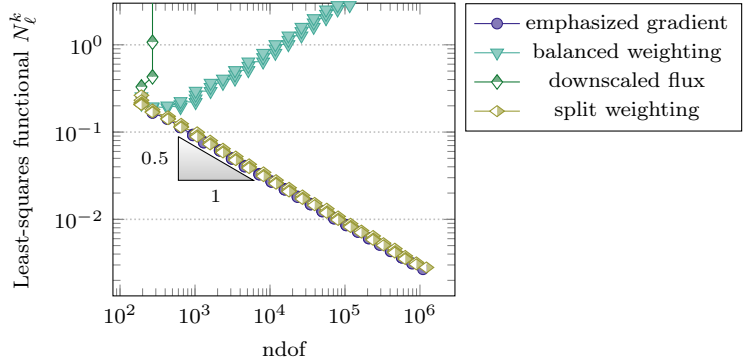
\begin{figure}
	\centering
	\subcaptionbox{\label{fig:nonlinear:convergence_history:weighting:residual} Error measure}[.48\textwidth]{
		\input{figures/Fig05a_Nonlinear_convergence_weighting.tex}
	}
	\hfil
	\subcaptionbox{\label{fig:nonlinear:convergence_history:weighting:condition}
		Condition number}[.48\textwidth]{
		\input{figures/Fig05b_Nonlinear_condition_weighting.tex}
	}
	\caption{%
		Convergence history and condition number plot 
		for Algorithm~\ref{alg:adaptive} applied to the convex energy minimization problem from Subsection~\ref{sec:convex_energy_minimization}
		with the weightings from Sections~\ref{sec:Zarantonello_least_squares} and~\ref{sec:alternative_weightings}.
		The chosen parameters read \(\delta = 1\), \(\theta = 0.3\), and \(\gamma = 0.9\).
	}
	\label{fig:nonlinear:convergence_history:weighting}
\end{figure}

A study of the damping parameter \(0 < \delta \leq 1\) is presented in
Figure~\ref{fig:nonlinear:convergence_history:damping}.
Since the damping parameter influences the value of the error estimators \(\eta_k\) and \(\mu_k\),
the plot only shows the values 
of the nonlinear least-squares functional \(N^k_\ell = N(f_1, f_2; p^k_\ell, u^k_\ell)^{1/2}\)
as a unified measure of the overall error; cf.~Subsection~\ref{sec:aposteriori}.
It is well-known that the damping parameter is crucial for the performance of the Zarantonello iteration
\cite[Section~6.4]{BringmannMiraciPraetorius2024}.
While large parameters like \(\delta \in \{0.5, 1\}\) lead to small estimator values and an optimal convergence rate,
the choices \(\delta \in \{0.01, 0.05, 0.1\}\) result in significantly larger estimator values and 
even a reduced convergence rate.
With the monotonicity and Lipschitz constants from Corollary~\ref{cor:well_posedness:ls_norm},
a theoretically justified sufficient value 
\[
	\delta_{\textup{LS}}^\star
	=
	\frac{2\alpha_{\textup{LS}}}{L_{\textup{LS}}^2} 
	\geq
	\frac{\Lambda_1^2}{144\Lambda_2^2}
	=
	\frac{1}{324}
	\approx
	3.086 \times 10^{-3}
	\geq \delta
\]
would lead to a very slow convergence in practice.
The undamped iteration with \(\delta = 1\) is most efficient in this example.
For this choice of \(\delta\),
Figure~\ref{fig:nonlinear:convergence_history:weighting} compares the weighting
of the least-squares functional \(Z_k\) as in~\eqref{eq:linearized_LS:functional} 
with the three alternative weightings presented in Section~\ref{sec:alternative_weightings}.
The choice of the weights \(\omega_1, \omega_2 > 0\) follows the corresponding theoretical values
from~\eqref{eq:choice_omega_weighting},
\eqref{eq:weighting2:choice_omega},
\eqref{eq:weighting3:choice_omega}, and
\eqref{eq:weighting4:choice_omega}.
The performance of the methods differs significantly for the four weightings.
While the adaptive algorithm with emphasized-gradient weighting
and with split weighting from Subsection~\ref{sec:split_weighting} converges with optimal rate,
the scheme with the balanced weighting from Subsection~\ref{sec:balanced_weighting}
or downscaled flux from Subsection~\ref{sec:downscaled_flux_variable} does not converge at all.
This empirically supports the better robustness of the weighting with emphasized gradient; 
see Remark~\ref{rem:damping_parameter}.
The stiffness matrix \(A_\ell \in \R^{M \times M}\) for the discrete problem~\eqref{eq:linearized_LS:discretization}
is independent of the damping parameter \(\delta > 0\) and the Zarantonello index \(k \in \N_0\) as both only belong to the right-hand side.
It depends, however, on the triangulation \(\TT_\ell\) and 
the weightings from Sections~\ref{sec:Zarantonello_least_squares} and~\ref{sec:alternative_weightings}.
Figure~\ref{fig:nonlinear:convergence_history:weighting:condition} reveals that all weighting strategies
result in comparable condition numbers with respect to the spectral radius
\begin{equation}
	\label{eq:condition}
	\operatorname{cond}_2(A_\ell)
	\coloneqq
	\Vert A \Vert_2 \Vert A^{-1} \Vert_2 
	=
	\frac{\lambda_{\textup{max}}(A_\ell)}{\lambda_{\textup{min}}(A_\ell)}.
\end{equation}

\begin{figure}
	\centering
	\subcaptionbox{\label{fig:nonlinear:errors} Convergence history of estimators from~\eqref{eq:applications:errors}}[.48\textwidth]{
		\input{figures/Fig01_Nonlinear_errors.tex}
	}
	\hfil
	\subcaptionbox{\label{fig:nonlinear:mesh} Adaptive mesh \(\TT_0^{15}\) (2\,424 triangles)}[.48\textwidth]{
		\input{figures/Fig01_Nonlinear_mesh.tex}
	}
	\bigskip

	\subcaptionbox{\label{fig:nonlinear:potential} Discrete function \(u_1^{40} \in S^1(\TT^{40}_1)\)}[.45\textwidth]{
		\input{figures/Fig01_Nonlinear_potential.tex}
	}
	\hfil
	\subcaptionbox{\label{fig:nonlinear:flux} Discrete function \(p_1^{40} \in RT^0(\TT^{40}_1)\)}[.45\textwidth]{
		\input{figures/Fig01_Nonlinear_flux.tex}
	}
	\caption{%
		Convergence and mesh plot as well as discrete solutions on the final mesh \(\TT^{40}_1\) (with \(\#\TT^{40}_1 = 548\,798\))
		for Algorithm~\ref{alg:adaptive} applied to the convex energy minimization problem from Subsection~\ref{sec:convex_energy_minimization}.
		The chosen parameters read \(\delta = 1\), \(\gamma = 0.9\), and \(\theta = 0.3\).
		(Figure~\subfigref{fig:nonlinear:flux} was created using the MATLAB function \texttt{quiver2.m} \cite{quiver2}.)
	}
	\label{fig:nonlinear:mesh_solution}
\end{figure}

\begin{figure}
	\centering
	\subcaptionbox{\label{fig:nonlinear:convergence_history:gamma:estimator} 
		Built-in estimator \(\eta^k_\ell + \mu^k_\ell\) for \(\delta = 1\)}[.48\textwidth]{
		\input{figures/Fig02a_Nonlinear_convergence_gamma.tex}
	}
	\hfil
	\subcaptionbox{\label{fig:nonlinear:convergence_history:gamma:residual} 
		Least-squares functional \(N^k_\ell\) for \(\delta = 1\)}[.48\textwidth]{
		\input{figures/Fig02b_Nonlinear_convergence_gamma.tex}
	}
	\bigskip

	\subcaptionbox{\label{fig:nonlinear:convergence_history:gamma:estimator:delta05}
		Built-in estimator \(\eta^k_\ell + \mu^k_\ell\) for \(\delta = 0.5\)}[.48\textwidth]{
		\input{figures/Fig02c_Nonlinear_convergence_gamma.tex}
	}
	\hfil
	\subcaptionbox{\label{fig:nonlinear:convergence_history:gamma:residual:delta05}
		Least-squares functional \(N^k_\ell\) for \(\delta = 0.5\)}[.48\textwidth]{
		\input{figures/Fig02d_Nonlinear_convergence_gamma.tex}
	}
	\caption{%
		Convergence history plot of the error estimators~\eqref{eq:applications:errors}
		for Algorithm~\ref{alg:adaptive} applied to the convex energy minimization problem from Subsection~\ref{sec:convex_energy_minimization}
		with various reduction parameters \(0 < \gamma < 1\).
		The remaining parameters read \(\delta \in \{0.5, 1\}\) and \(\theta = 0.3\).
	}
	\label{fig:nonlinear:convergence_history:gamma}
\end{figure}
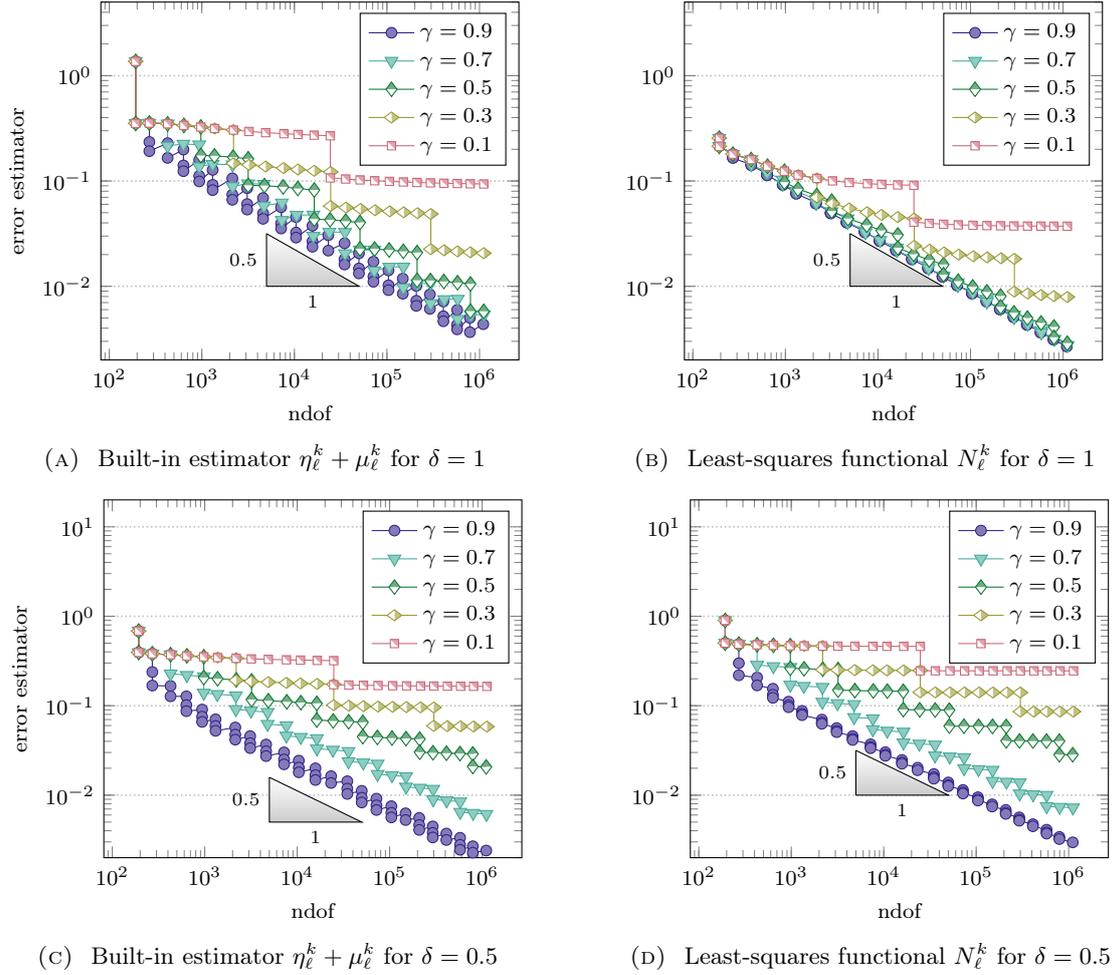

\begin{figure}
	\centering
	\subcaptionbox{\label{fig:nonlinear:convergence_history:theta:estimator} 
		Built-in estimator \(\eta^k_\ell + \mu^k_\ell\) for \(\delta = 1\)}[.48\textwidth]{
		\input{figures/Fig03a_Nonlinear_convergence_theta.tex}
	}
	\hfil
	\subcaptionbox{\label{fig:nonlinear:convergence_history:theta:residual}
		Least-squares functional \(N^k_\ell\) for \(\delta = 1\)}[.48\textwidth]{
		\input{figures/Fig03b_Nonlinear_convergence_theta.tex}
	}
	\bigskip

	\subcaptionbox{\label{fig:nonlinear:convergence_history:theta:estimator:delta05}
		Built-in estimator \(\eta^k_\ell + \mu^k_\ell\) for \(\delta = 0.5\)}[.48\textwidth]{
		\input{figures/Fig03c_Nonlinear_convergence_theta.tex}
	}
	\hfil
	\subcaptionbox{\label{fig:nonlinear:convergence_history:theta:residual:delta05}
		Least-squares functional \(N^k_\ell\) for \(\delta = 0.5\)}[.48\textwidth]{
		\input{figures/Fig03d_Nonlinear_convergence_theta.tex}
	}
	\caption{%
		Convergence history plot of the estimators~\eqref{eq:applications:errors}
		for Algorithm~\ref{alg:adaptive} applied to the convex energy minimization problem from Subsection~\ref{sec:convex_energy_minimization}
		with various choices of the bulk parameter \(0 < \theta \leq 1\).
		The remaining parameters read \(\delta = 1\) and \(\gamma = 0.9\).
	}
	\label{fig:nonlinear:convergence_history:theta}
\end{figure}
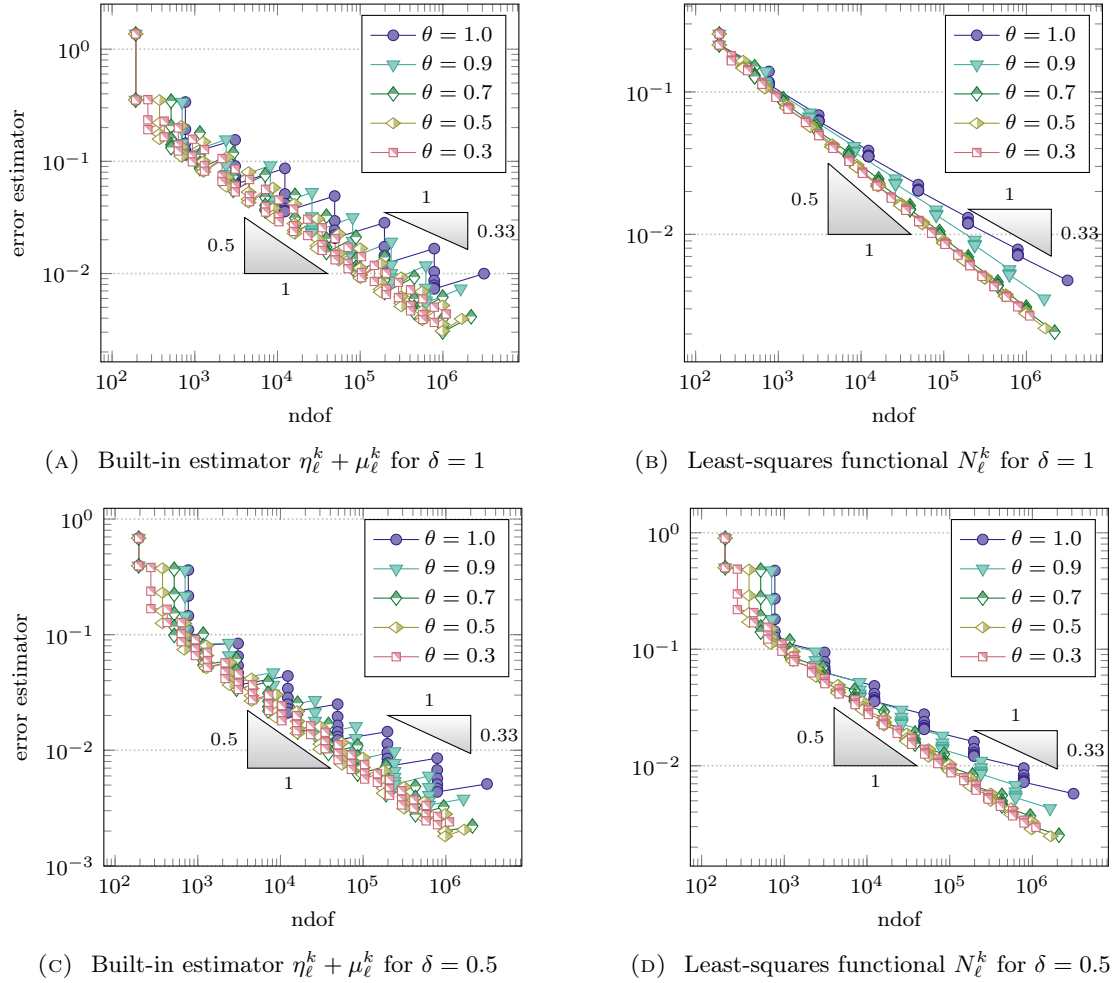

\subsection{Porous media flow}
\label{sec:porous_media_flow}
\noindent
This subsection is devoted to a model for the flow of a fluid 
through a porous medium \(\Omega \subset \R^2\) without gravity
which is not fully covered by our theory.
We choose material parameters \(k_1 = 0.2\) and \(k_2 = 20\)
in the nonlinear mapping \(\sigma\colon \R^d \to \R^d\) with
\begin{equation}
	\label{eq:porousmedium:sigma}
	\sigma(\xi)
	\coloneqq
	\frac{2\xi}{k_1 + \sqrt{k_1^2 + k_2 \vert \xi \vert}}.
\end{equation}
Since \(\vert \D\sigma(\xi) \vert \to 0\) as \(\vert \xi \vert \to \infty\),
its derivative does not satisfy the monotonicity assumption~\eqref{assum:ellipticity}.

The variable \(u\colon \R^2 \to \R\) describes the pressure of the fluid.
The relation between the pressure gradient \(\nabla u\) and the fluid velocity \(p\colon \R^2 \to \R^2\)
is given by Forchheimer's law
\[
	- p
	=
	\sigma(\nabla u)
	=
	\frac{2\nabla u}{k_1 + \sqrt{k_1^2 + k_2 \vert \nabla u \vert}}
	\quad \text{in } \Omega.
\]
Given an external mass flow rate \(f \in L^2(\Omega)\), this law is complemented by the mass conservation equation
\(\div p = f\) in \(\Omega\); see, e.g., \cite{DouglasPaesLemeGiorgi1993,Park1995} and the references therein.
Note that the different sign convention for \(p\) does not affect the analysis in the previous sections.
The coefficient function \(\phi\colon [0, \infty) \to \R\) and its derivative read
\[
	\phi(t)
	=
	\frac{2}{k_1 + \sqrt{k_1^2 + k_2 t}}
	\quad \text{and} \quad
	\phi'(t)
	=
	-\frac{k_2}{(k_1 + \sqrt{k_1^2 + k_2 t})^2 \sqrt{k_1^2 + k_2 t}}.
\]
Given any upper bound \(T > 0\), it holds, for \(0 \leq t \leq T\),
\begin{gather*}
	\frac{2}{k_1 + \sqrt{k_1^2 + k_2 T}}
	< 
	\phi(t)
	\leq
	\frac{1}{k_1},
	\quad
	\varphi(T)
	<
	\varphi(t)
	\coloneqq
	\phi(t) + t \phi'(t)
	=
	\frac{2 k_1}{(k_1 + \sqrt{k_1^2 + k_2 t}) \sqrt{k_1^2 + k_2 t}}
	\leq
	\frac{1}{k_1}.
\end{gather*}
Hence, the assumptions~\eqref{assum:ellipticity}--\eqref{assum:boudedness} are satisfied
with \(\Lambda_1 = \varphi(T)\) and \(\Lambda_2 = 1 / k_1\) on the bounded set \(\{\xi \in \R^2 : \vert \xi \vert \leq T\}\).
Under the assumption that the gradient \(\nabla u^\star\) of the exact solution is uniformly bounded
\(\vert u^\star \vert \leq T\) almost everywhere in \(\Omega\),
the analysis of the Zarantonello LSFEM from Section~\ref{sec:Zarantonello_least_squares} applies.
This assumption has been made, e.g., in \cite[Section~1]{Park1995} and 
has been rigorously proven in \cite[Theorem~1.4]{CianchiMazya2011} for convex domains in three spatial dimensions.
However, we expect that the solution in the setting at hand with a nonconvex domain is less regular and thus does not satisfy such a uniform bound.
Nevertheless, Algorithm~\ref{alg:adaptive} performs well in the experiments and,
heuristically, we choose \(T = 10^{-2}\) leading to the constants
\(\Lambda_1 = \varphi(T) \approx 1.1835\) and \(\Lambda_2 = k_1^{-1} = 5\) 
in the conditions~\eqref{assum:ellipticity}--\eqref{assum:boudedness}.
This motivates the selection of the weights and the damping parameter as
\[
	\omega_1^{2} = \frac{2\Lambda_2^2}{\Lambda_1^2} \approx 35.6969
	\quad \text{and} \quad
	\omega_2^{2} = \frac{\Lambda_2^2}{\Lambda_1} \approx 21.1237.
\]
Figure~\ref{fig:porousmedium:damping:grad_norm} supports this choice empirically by showing 
a (generously chosen) upper bound \(T = 10^{-2}\) of the discrete gradient norm \(\Vert \nabla u^k_\ell \Vert_{L^\infty(\Omega)}\).

In the remaining part of this section, we consider a benchmark problem
on the L-shaped domain \(\Omega = (-1, 1)^2 \setminus [0, 1)^2\)
with Friedrichs constant \(C_{\textup{F}} \leq 0.3221\) 
from Subsection~\ref{sec:convex_energy_minimization}.
The given right-hand side \(f \in L^2(\Omega)\) 
with local support \(\operatorname{supp}(f) = [-0.6, -0.4] \times [0.4, 0.6]\)
is illustrated in Figure~\ref{fig:porousmedium:rhs} and reads
\[
	f(x)
	\coloneqq
	\begin{cases}
		1 & \text{if } -0.6 < x_1 < -0.4 \text{ and } 0.4 < x_2 < 0.6, \\
		0 & \text{otherwise.}
	\end{cases}
\]
The following experiments consider lowest-order discretizations with \(m = 0\).

The adaptively generated mesh \(\TT^k_\ell\) with \(k = 17\) and \(\ell = 1\) in Figure~\ref{fig:porousmedium:mesh}
displays a strong refinement towards the re-entrant corner of the L-shaped domain
as well as at the support of the right-hand side \(f\).
Figures~\ref{fig:porousmedium:potential}--\ref{fig:porousmedium:flux} show the discrete solution 
\((p^{46}_{1}, u^{46}_{1}) \in RT^0(\TT^{46}_{1}) \times S^1_0(\TT^{46}_{1})\)
where discrete flux variable \(p^{46}_{1}\) is evaluated at equidistributed 208 points in the domain \(\Omega\).
The plots illustrate that the fluid is transported away from the region with high mass flow rate \(f\)
where the pressure \(u^\star\) is high.
The fact that both variables are physically relevant quantities
make the least-squares approach particularly attractive for this problem
as it provides equal approximation quality for both variables.

The flux mapping \(\sigma\) from~\eqref{eq:porousmedium:sigma} resembles the nonlinearity 
of the \(\textup{p}\)-Laplace problem for small exponent \(\textup{p} = 3/2\);
see \cite[Section~6.1]{BelenkiDieningKreuzer2012} for a discussion of the regularity of solutions in that context.
Related numerical experiments in \cite[Section~6.2]{BelenkiDieningKreuzer2012} and~\cite[Section~6]{DieningFornasierTomasiWank2020}
indicate an expected optimal convergence rate of \(0.5\) with respect to (cumulative) number of degrees of freedom.

\begin{figure}
	\centering
	\subcaptionbox{\label{fig:porousmedium:weighting:error} Error measure}[.48\textwidth]{
		\input{figures/Fig08a_PorousMedium_convergence_weighting.tex}
	}
	\hfil
	\subcaptionbox{\label{fig:porousmedium:weighting:condition} Condition number}[.48\textwidth]{
		\input{figures/Fig08b_PorousMedium_condition_weighting.tex}
	}
	\caption{%
		Convergence history and condition number plots for Algorithm~\ref{alg:adaptive}
		applied to the porous medium flow problem from Subsection~\ref{sec:porous_media_flow}
		with the weightings from Sections~\ref{sec:Zarantonello_least_squares} and~\ref{sec:alternative_weightings}.
		The chosen parameters read \(\delta = 1\), \(\gamma = 0.9\), and \(\theta = 0.3\).
	}
	\label{fig:porousmedium:weighting}
\end{figure}
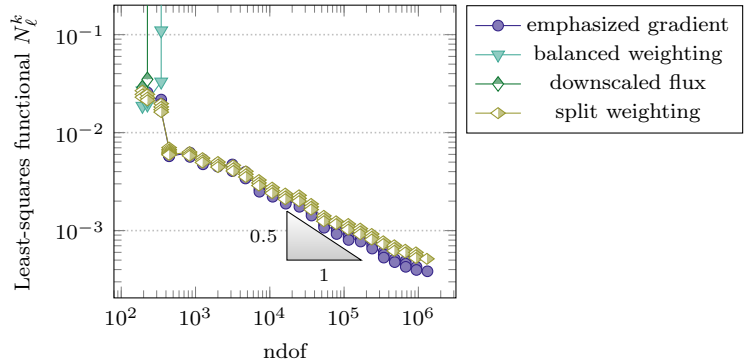

The investigation of the damping parameter \(0 < \delta \leq 1\) in Figure~\ref{fig:porousmedium:damping:error}
confirms the observations from Subsection~\ref{sec:convex_energy_minimization}
and shows best performance with the expected optimal rate \(0.5\) for the undamped iteration with \(\delta = 1\).
Figure~\ref{fig:porousmedium:weighting:error} shows that, in the present example again,
only the emphasized-gradient and the split weighting converge.
They also turn out to be more favorable concerning the condition number~\eqref{eq:condition} of the discrete linearized system
with respect to the spectral radius as displayed in Figure~\ref{fig:porousmedium:weighting:condition}.

\begin{figure}
	\centering
	\subcaptionbox{\label{fig:porousmedium:rhs} Flow mass rate \(f \in L^2(\Omega)\)}[.45\textwidth]{
		\input{figures/Fig06_PorousMedium_rhs.tex}
	}
	\hfill
	\subcaptionbox{\label{fig:porousmedium:mesh} Adaptive mesh \(\TT_{1}^{17}\) (1\,568 triangles)}[.45\textwidth]{
		\input{figures/Fig06_PorousMedium_mesh.tex}
	}
	\bigskip

	\subcaptionbox{\label{fig:porousmedium:potential} Discrete pressure \(u^{46}_{1}\) with \(\#\TT^{46}_1 = 662\,621\)}[.45\textwidth]{
		\input{figures/Fig06_PorousMedium_potential.tex}
	}
	\hfil
	\subcaptionbox{\label{fig:porousmedium:flux} Discrete flux \(p^{46}_{1}\) with \(\#\TT^{46}_1 = 662\,621\)}[.45\textwidth]{
		\input{figures/Fig06_PorousMedium_flux.tex}
	}
	\caption{%
		Right-hand side, mesh and solution plots for the porous medium flow problem from Subsection~\ref{sec:porous_media_flow}.
		The chosen parameters read \(\delta = 1\), \(\gamma = 0.9\), and \(\theta = 0.3\).
		(Figure~\subfigref{fig:porousmedium:flux} was created using the MATLAB function \texttt{quiver2.m} \cite{quiver2}.)
	}
	\label{fig:porousmedium:mesh_solution}
\end{figure}
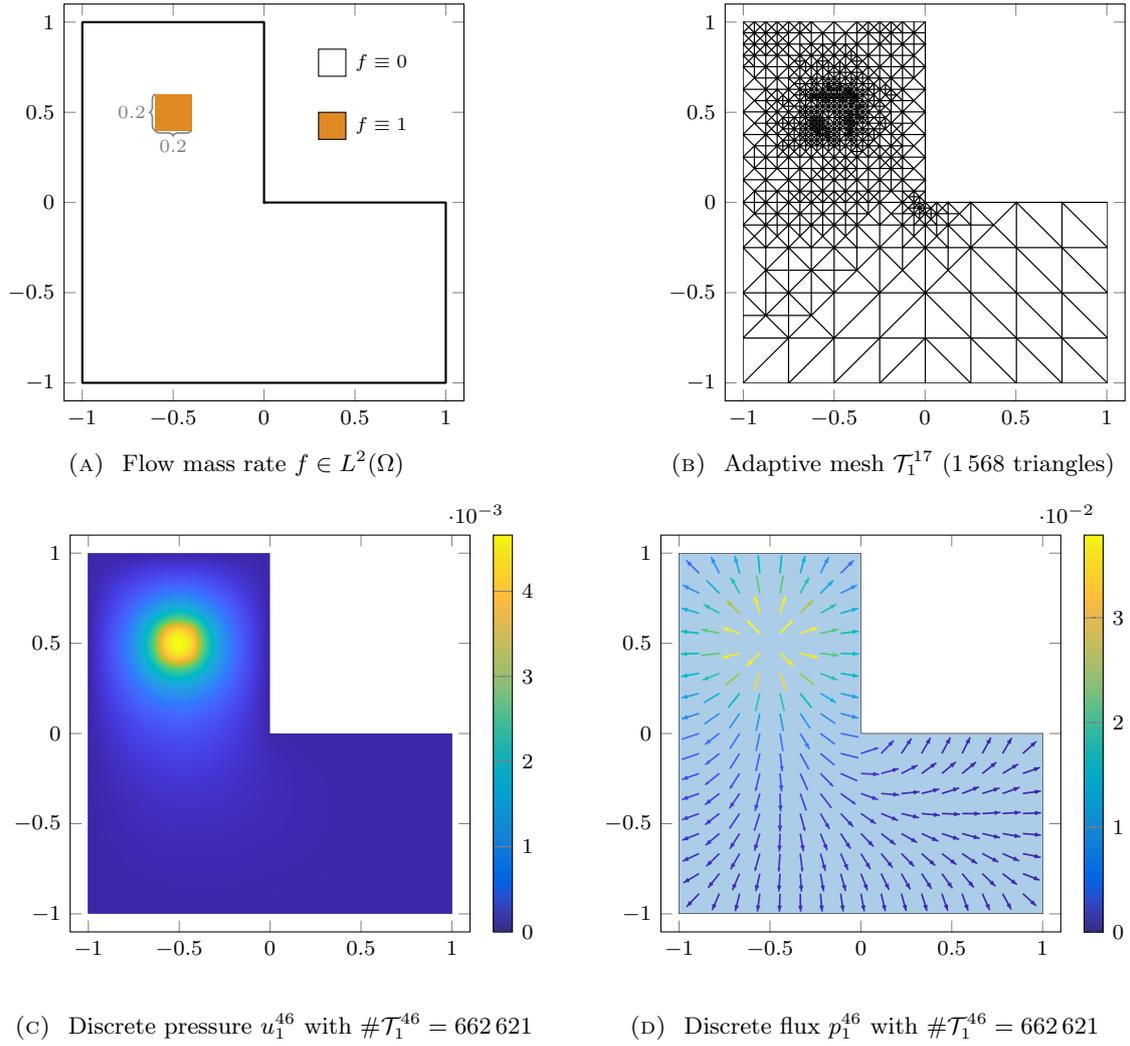

\begin{figure}
	\subcaptionbox{\label{fig:porousmedium:damping:error} Error measure}[.48\textwidth]{
		\input{figures/Fig07a_PorousMedium_convergence_damping.tex}
	}
	\hfil
	\subcaptionbox{\label{fig:porousmedium:damping:grad_norm} Gradient norm}[.48\textwidth]{
		\input{figures/Fig07b_PorousMedium_gradient_norm_damping.tex}
	}
	\caption{%
		Convergence history and norm plots for Algorithm~\ref{alg:adaptive}
		applied to the porous medium flow problem from Subsection~\ref{sec:porous_media_flow}
		for various choices of the damping parameter \(0 < \delta \leq 1\).
		The remaining parameters read \(\gamma = 0.9\) and \(\theta = 0.3\).
	}
	\label{fig:porousmedium:damping}
\end{figure}
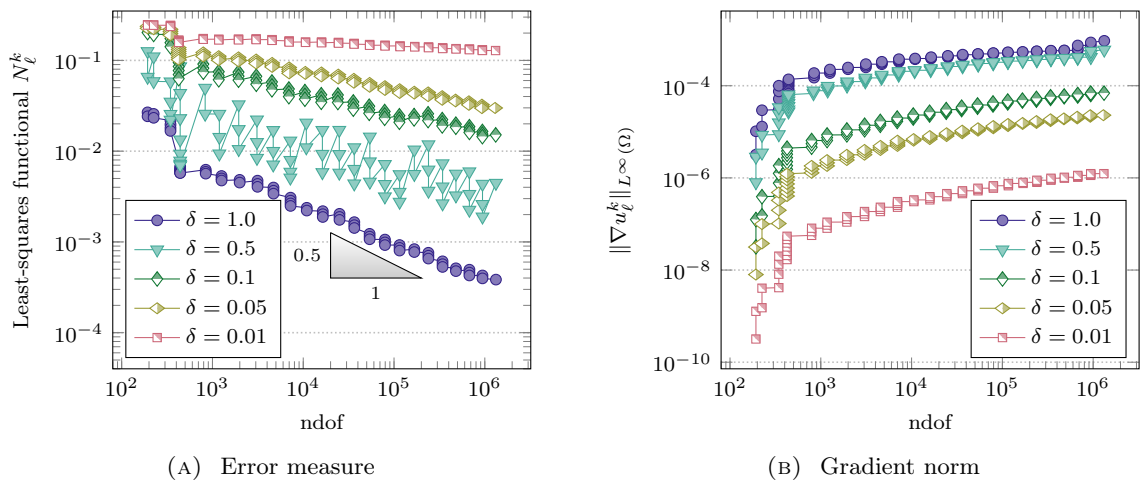

{
    \sloppy
	\printbibliography  
}

\appendix

\section{Weighting 2: Balanced weighting}
\label{app:balanced_weighting}
\noindent
This appendix is devoted to the proofs of the results from Subsection~\ref{sec:balanced_weighting}.
They concern the nonlinear mapping \(\widetilde{\mathcal{B}}\) and the norms 
\(\vvvert \cdot \vvvert_{\widetilde{\mathcal{A}}}\) and \(\vvvert \cdot \vvvert_{\widetilde\omega}\)
as introduced in~\eqref{eq:weighting2} with weights \(\omega_1, \omega_2 > 0\)
chosen according to~\eqref{eq:weighting2:choice_omega}.
The strong monotonicity and Lipschitz continuity~\eqref{eq:alternative_weightings:monotonicity_lipschitz}
with the constants from~\eqref{eq:weighting2:constants}
are a direct consequence of the following estimates:
For all \((p, u), (q, v), (r, z) \in H(\div, \Omega) \times H^1_0(\Omega)\),
it holds that
\begin{gather}
	\label{eq:weighting2:fundamental_equivalence}
	\min\bigg\{
		\frac12,\:
		\bigg(1 + \frac{2 \Lambda_1^{5/2}}{\Lambda_2^3}\bigg)^{-1}
	\bigg\} \,
	\vvvert (q, v) \vvvert_{\widetilde\omega}^2
	\leq
	\vvvert (q, v) \vvvert_{\widetilde{\mathcal{A}}}^2
	\leq
	2\,
	\vvvert (q, v) \vvvert_{\widetilde\omega}^2,
	\\
	\label{eq:weighting2:monotonicity_weighted_norm}
	\min\bigg\{
		\frac12,\:
		\frac{\Lambda_2}{\Lambda_1^{1/2}},\:
		\frac{\Lambda_1^{3/2}}{4 \Lambda_2}
	\bigg\}\,
	\vvvert (p - q, u - v) \vvvert_{\widetilde\omega}^2
	\leq
	\widetilde{\mathcal{B}}(p, u; p - q, u - v) - \widetilde{\mathcal{B}}(q, v; p - q, u - v),
	\\
	\label{eq:weighting2:Lipschitz_weighted_norm}
	\widetilde{\mathcal{B}}(p, u; r, z) - \widetilde{\mathcal{B}}(q, v; r, z)
	\leq
	2
	\max\bigg\{
		1,\:
		\frac{\Lambda_2}{\Lambda_1^{1/2}},\:
		\frac{\Lambda_1^{3/2}}{2\Lambda_2}
	\bigg\}\,
	\vvvert (p - q, u - v) \vvvert_{\widetilde\omega} \,
	\vvvert (r, z) \vvvert_{\widetilde\omega}.
\end{gather}
The remaining part of this section contains the proofs of these three estimates.

\begin{proof}[Proof of equivalence~\eqref{eq:weighting2:fundamental_equivalence}]
	The proof follows the same steps as the proof of Theorem~\ref{thm:fundamental_equivalence}.
	It is given here in full detail for the ease of reading.

	\emph{Step~1 (lower bound).}\,
	The binomial formula and an integration by parts show
    \begin{align*}
		\Vert \widetilde\omega_2^{-1} \, q \Vert_{L^2(\Omega)}^2
		+
		\Vert \widetilde\omega_2 \, \nabla v \Vert_{L^2(\Omega)}^2
		&=
		\Vert \widetilde\omega_2^{-1} \, q - \omega_2 \, \nabla v \Vert_{L^2(\Omega)}^2
		+
		2\, (q, \nabla v)_{L^2(\Omega)}
		\\
		&=
		\Vert \widetilde\omega_2^{-1} \, q - \widetilde\omega_2 \, \nabla v \Vert_{L^2(\Omega)}^2
		-
		2 \, (\div q, v)_{L^2(\Omega)}.
    \end{align*}
	The Cauchy--Schwarz, Friedrichs and weighted Young inequality prove
    \begin{align*}
		- 2 \,
		(\div q, v)_{L^2(\Omega)}
		&\leq
		2 \,
		\Vert \div q \Vert_{L^2(\Omega)} \, \Vert v \Vert_{L^2(\Omega)}
		\leq
		2 C_{\textup{F}} \,
		\Vert \div q \Vert_{L^2(\Omega)} \, \Vert \nabla v \Vert_{L^2(\Omega)}
		\\
		&\leq
		\frac{2C_{\textup{F}}}{\widetilde\omega_1 \widetilde\omega_2} \,
		\Vert \widetilde\omega_1 \, \div q \Vert_{L^2(\Omega)} \, \Vert \widetilde\omega_2 \, \nabla v \Vert_{L^2(\Omega)}
		\\
		&\leq
		\frac{2C_{\textup{F}}^2}{\widetilde\omega_1^2 \widetilde\omega_2^2} \,
		\Vert \widetilde\omega_1 \, \div q \Vert_{L^2(\Omega)}^2
		+
		\frac12 \,
		\Vert \widetilde\omega_2 \, \nabla v \Vert_{L^2(\Omega)}^2.
    \end{align*}
	The combination of the two previous formulas
	and the absorption of \(\frac12 \, \Vert \widetilde\omega_2 \, \nabla v \Vert_{L^2(\Omega)}^2\) into the left-hand side yield
    \[
		2 \,
		\Vert \widetilde\omega_2^{-1} \, q \Vert_{L^2(\Omega)}^2
        +
		\Vert \widetilde\omega_2 \, \nabla v \Vert_{L^2(\Omega)}^2
        \leq
		\frac{4C_{\textup{F}}^2}{\widetilde\omega_1^2 \widetilde\omega_2^2} \,
		\Vert \widetilde\omega_1 \, \div q \Vert_{L^2(\Omega)}^2
        +
		2 \, \Vert \widetilde\omega_2^{-1} \, q - \widetilde\omega_2 \, \nabla v \Vert_{L^2(\Omega)}^2.
    \]
	The addition of \(C_{\textup{F}}^2 \, \Vert \widetilde\omega_1 \, \div q \Vert_{L^2(\Omega)}^2\)
	concludes the proof of the lower bound with
    \begin{align*}
		\vvvert (q, v) \vvvert_{\widetilde\omega}^2
		&\leq
		C_{\textup{F}}^2 \,
		\Vert \widetilde\omega_1 \, \div q \Vert_{L^2(\Omega)}^2
		+
		\Vert \widetilde\omega_2^{-1} \, q \Vert_{L^2(\Omega)}^2
        +
		\Vert \widetilde\omega_2 \, \nabla v \Vert_{L^2(\Omega)}^2
		\\
		&\leq
		\bigg(1 + \frac{4}{\widetilde\omega_1^2\widetilde\omega_2^2}\bigg) C_{\textup{F}}^2 \, 
		\Vert \widetilde\omega_1 \, \div q \Vert_{L^2(\Omega)}^2
        +
		2\, \Vert \widetilde\omega_2^{-1} \, q - \widetilde\omega_2 \, \nabla v \Vert_{L^2(\Omega)}^2
		\\
		&\eqreff*{eq:weighting2:ls_norm}\leq
		\max\bigg\{
			2,\:
			1 + \frac{4}{\widetilde\omega_1^2\widetilde\omega_2^2}
		\bigg\} \, 
		\vvvert (q, v) \vvvert_{\widetilde{\mathcal{A}}}^2
		\eqreff{eq:weighting2:choice_omega}=
		\max\bigg\{
			2,\:
			1 + \frac{2 \Lambda_1^{5/2}}{\Lambda_2^3}
		\bigg\} \,
		\vvvert (q, v) \vvvert_{\widetilde{\mathcal{A}}}^2.
    \end{align*}

	\emph{Step~2 (upper bound).}\,
	As in the proof of Theorem~\ref{thm:fundamental_equivalence},
	the triangle and the Young inequality verify the upper bound
	\(
		\vvvert (q, v) \vvvert_{\widetilde{\mathcal{A}}}^2
        \leq
		2 \,
		\vvvert (q, v) \vvvert_{\widetilde\omega}^2.
	\)
\end{proof}

\begin{proof}[Proof of monotonicity~\eqref{eq:weighting2:monotonicity_weighted_norm}]
	The proof follows the same steps as for Theorem~\ref{thm:well_posedness}.
	With the changed weighting, the equality~\eqref{eq:monotonicity:split} reads
	\begin{align*}
		&\hspace{-2em}(p - q - [\sigma(\nabla u) - \sigma(\nabla v)],\: \widetilde\omega_2^{-1} \, (p - q) - \widetilde\omega_2 \, \nabla(u -
		v))_{L^2(\Omega)}
		\\
		&=
		(p - q - M\nabla (u - v),\: \widetilde\omega_2^{-1} \, (p - q) - \widetilde\omega_2 \, \nabla(u - v))_{L^2(\Omega)}
		\\
		&=
		\widetilde\omega_2^{-1} \, 
		\Vert p - q \Vert_{L^2(\Omega)}^2
		+
		\widetilde\omega_2 \,
		(M \nabla(u - v),\: \nabla(u - v) )_{L^2(\Omega)}
		\\
		&\phantom{{}\leq{}}
		+
		\widetilde\omega_2 \,
		(\div (p - q),\: u - v)_{L^2(\Omega)}
		-
		\widetilde\omega_2^{-1} \, 
		(p - q,\: M \nabla(u - v))_{L^2(\Omega)}.
	\end{align*}
	The combination with~\eqref{eq:monotonicity:ellipticity}--\eqref{eq:monotonicity:boudedness}
	and adding \(C_{\textup{F}}^2 \Vert \widetilde\omega_1 \div (p - q) \Vert_{L^2(\Omega)}^2\) to both sides show
	\begin{align*}
		&\hspace{-0.5em}
		\Big(\widetilde\omega_1^2 - \frac{\widetilde\omega_2}{\Lambda_1} \Big) C_{\textup{F}}^2\,
		\Vert \div (p - q) \Vert_{L^2(\widetilde\Omega)}^2
		+
		\frac{\widetilde\omega_2^{-1}}{2} \, \Vert p - q \Vert_{L^2(\Omega)}^2
		+
		\frac{3\widetilde\omega_2 \Lambda_1 - 2 \widetilde\omega_2^{-1} \Lambda_2^2}4
		\, 
		\Vert \nabla(u - v) \Vert_{L^2(\Omega)}^2
		\\
		&\leq
		C_{\textup{F}}^2 \,
		\Vert \widetilde\omega_1 \div (p - q) \Vert_{L^2(\Omega)}^2
		+
		(p - q - [\sigma(\nabla u) - \sigma(\nabla v)],\: \widetilde\omega_2^{-1} \, (p - q) - \widetilde\omega_2 \, \nabla(u - v))_{L^2(\Omega)}
		\\
		&=
		\widetilde{\mathcal{B}}(p, u; p - q, u - v) - \widetilde{\mathcal{B}}(q, v; p - q, u - v).
	\end{align*}
	The choice of the weights in~\eqref{eq:weighting2:choice_omega} leads to 
	\begin{multline*}
		\frac{C_{\textup{F}}^2}{2} \,
		\Vert \widetilde\omega_1 \div (p - q) \Vert_{L^2(\Omega)}^2
		+
		\frac{\widetilde\omega_2}2 \, 
		\Vert \widetilde\omega_2^{-1} \, (p - q) \Vert_{L^2(\Omega)}^2
		+
		\frac{\Lambda_1}{4\widetilde\omega_2} \, 
		\Vert \widetilde\omega_2 \, \nabla(u - v) \Vert_{L^2(\Omega)}^2
		\\
		\leq
		\widetilde{\mathcal{B}}(p, u; p - q, u - v) - \widetilde{\mathcal{B}}(q, v; p - q, u - v)
	\end{multline*}
	and concludes the proof of
	\begin{align*}
		\min\bigg\{
			\frac12,\:
			\frac{\Lambda_2}{\Lambda_1^{1/2}},\:
			\frac{\Lambda_1^{3/2}}{4 \Lambda_2}
		\bigg\}\,
		\vvvert (p - q, u - v) \vvvert_{\widetilde\omega}^2
		&=
		\min\bigg\{
			\frac12,\:
			\frac{\widetilde\omega_2}{2},\:
			\frac{\Lambda_1}{4\widetilde\omega_2}
		\bigg\} \,
		\vvvert (p - q, u - v) \vvvert_{\widetilde\omega}^2
		\\
		&\leq
		\widetilde{\mathcal{B}}(p, u; p - q, u - v) - \widetilde{\mathcal{B}}(q, v; p - q, u - v).
		\qedhere
	\end{align*}
\end{proof}

\begin{proof}[Proof of Lipschitz continuity~\eqref{eq:weighting2:Lipschitz_weighted_norm}]
	An analogous computation as in the proof of~\eqref{eq:weighting2:monotonicity_weighted_norm} establishes
	\begin{align*}
		&\hspace{-2em}
		\widetilde{\mathcal{B}}(p, u; r, z) - \widetilde{\mathcal{B}}(q, v; r, z)
		\\
		&=
		\widetilde\omega_1^2 C_{\textup{F}}^2 \, (\div (p - q),\: \div r)_{L^2(\Omega)}
		+
		(p - q - [\sigma(\nabla u) - \sigma(\nabla v)],\: \widetilde\omega_2^{-1} \, r - \widetilde\omega_2 \, \nabla z)_{L^2(\Omega)}
		\\
		&=
		\widetilde\omega_1^2 C_{\textup{F}}^2 \,
		(\div (p - q),\: \div r)_{L^2(\Omega)}
		+
		\widetilde\omega_2^{-1} \, (p - q,\: r)_{L^2(\Omega)}
		+
		\widetilde\omega_2 \, (M \nabla(u - v),\: \nabla z)_{L^2(\Omega)}
		\\
		&\phantom{{}\leq{}}
		+
		\widetilde\omega_2 \,
		(\div (p - q),\: z)_{L^2(\Omega)}
		-
		\widetilde\omega_2^{-1} \, (r,\: M \nabla(u - v))_{L^2(\Omega)}.
	\end{align*}
	The boundedness of \(\D\sigma\) from~\eqref{assum:boudedness} and a
	Cauchy--Schwarz inequality in \(L^2(\Omega)\) imply
	\begin{align*}
		&\hspace{-1em}
		\widetilde{\mathcal{B}}(p, u; r, z) - \widetilde{\mathcal{B}}(q, v; r, z)
		\\
		&\leq
		\widetilde\omega_1^2 C_{\textup{F}}^2 \,
		\Vert \div (p - q) \Vert_{L^2(\Omega)} \, \Vert \div r \Vert_{L^2(\Omega)}
		+
		\widetilde\omega_2^{-1} \, \Vert p - q \Vert_{L^2(\Omega)} \, \Vert r \Vert_{L^2(\Omega)}
		\\
		&\phantom{{}\leq{}}
		+
		\widetilde\omega_2 \Lambda_2 \,
		\Vert \nabla(u - v) \Vert_{L^2(\Omega)} \, \Vert \nabla z \Vert_{L^2(\Omega)}
		+
		\widetilde\omega_2 C_{\textup{F}} \,
		\Vert \div (p - q) \Vert_{L^2(\Omega)} \, \Vert \nabla z \Vert_{L^2(\Omega)}
		\\
		&\phantom{{}\leq{}}
		+
		\widetilde\omega_2^{-1} \Lambda_2 \,
		\Vert r \Vert_{L^2(\Omega)} \, \Vert \nabla(u - v) \Vert_{L^2(\Omega)}
		\\
		&=
		C_{\textup{F}}^2 \,
		\Vert \widetilde\omega_1 \div (p - q) \Vert_{L^2(\Omega)} \, \Vert \widetilde\omega_1 \div r \Vert_{L^2(\Omega)}
		+
		\widetilde\omega_2 \,
		\Vert \widetilde\omega_2^{-1} \, (p - q) \Vert_{L^2(\Omega)} \, \Vert \widetilde\omega_2^{-1} \, r \Vert_{L^2(\Omega)}
		\\
		&\phantom{{}\leq{}}
		+
		\frac{\Lambda_2}{\widetilde\omega_2} \,
		\Vert \widetilde\omega_2 \, \nabla(u - v) \Vert_{L^2(\Omega)} \, \Vert \widetilde\omega_2 \, \nabla z \Vert_{L^2(\Omega)}
		+
		\frac{1}{\widetilde\omega_1} C_{\textup{F}} \,
		\Vert \widetilde\omega_1 \div (p - q) \Vert_{L^2(\Omega)} \, \Vert \widetilde\omega_2 \, \nabla z \Vert_{L^2(\Omega)}
		\\
		&\phantom{{}\leq{}}
		+
		\frac{\Lambda_2}{\widetilde\omega_2} \,
		\Vert \widetilde\omega_2^{-1} \, r \Vert_{L^2(\Omega)} \, \Vert \widetilde\omega_2 \, \nabla(u - v) \Vert_{L^2(\Omega)}.
	\end{align*}
	A Cauchy--Schwarz inequality in \(\R^5\) results in
	\begin{align*}
		&\hspace{-2em}
		\widetilde{\mathcal{B}}(p, u; r, z) - \widetilde{\mathcal{B}}(q, v; r, z)
		\\
		&\leq
		\max\bigg\{
			1,\:
			\widetilde\omega_2,\:
			\frac{\Lambda_2}{\widetilde\omega_{2}},\:
			\frac{1}{\widetilde\omega_1}
		\bigg\}
		\\
		&\phantom{{}\leq{}}
		\times
		\Big[
			2 C_{\textup{F}}^2 \,
			\Vert \widetilde\omega_1 \div (p - q) \Vert_{L^2(\Omega)}^2
			+
			\Vert \widetilde\omega_2^{-1} (p - q) \Vert_{L^2(\Omega)}^2
			+
			2 \Vert \widetilde\omega_2 \, \nabla(u - v) \Vert_{L^2(\Omega)}^2
		\Big]^{1/2}
		\\
		&\phantom{{}\leq{}}
		\times
		\Big[
			C_{\textup{F}}^2 \,
			\Vert \widetilde\omega_1 \div r \Vert_{L^2(\Omega)}^2
			+
			2 \,
			\Vert \widetilde\omega_2^{-1} \, r \Vert_{L^2(\Omega)}^2
			+
			2 \,
			\Vert \widetilde\omega_2 \, \nabla z \Vert_{L^2(\Omega)}^2
		\Big]^{1/2}.
	\end{align*}
	This and the estimate \(\Lambda_1 \leq \Lambda_2\) conclude the proof of
	\begin{align*}
		\widetilde{\mathcal{B}}(p, u; r, z) - \widetilde{\mathcal{B}}(q, v; r, z)
		&\leq
		2
		\max\bigg\{
			1,\:
			\widetilde\omega_2,\:
			\frac{\Lambda_2}{\widetilde\omega_{2}},\:
			\frac{1}{\widetilde\omega_1}
		\bigg\}
		\vvvert (p - q, u - v) \vvvert_{\widetilde\omega} \,
		\vvvert (r, z) \vvvert_{\widetilde\omega}
		\\
		&\eqreff*{eq:weighting2:choice_omega}=
		2
		\max\bigg\{
			1,\:
			\frac{\Lambda_2}{\Lambda_1^{1/2}},\:
			\frac{\Lambda_1}{\Lambda_2},\:
			\frac{\Lambda_1^{3/2}}{2\Lambda_2}
		\bigg\}\,
		\vvvert (p - q, u - v) \vvvert_{\widetilde\omega} \,
		\vvvert (r, z) \vvvert_{\widetilde\omega}
		\\
		&=
		2
		\max\bigg\{
			1,\:
			\frac{\Lambda_2}{\Lambda_1^{1/2}},\:
			\frac{\Lambda_1^{3/2}}{2\Lambda_2}
		\bigg\}\,
		\vvvert (p - q, u - v) \vvvert_{\widetilde\omega} \,
		\vvvert (r, z) \vvvert_{\widetilde\omega}.
		\qedhere
	\end{align*}
\end{proof}

%
%
\section{Weighting 3: Downscaled flux}
\label{app:downscale_flux}
\noindent
This section is devoted to the proofs for Subsection~\ref{sec:balanced_weighting}
guaranteeing the strong monotonicity and Lipschitz continuity~\eqref{eq:alternative_weightings:monotonicity_lipschitz}
with the constants~\eqref{eq:weighting3:constants}.
The assertion for the nonlinear mapping \(\widetilde{\mathcal{B}}\) and the norms
\(\vvvert \cdot \vvvert_{\widetilde{\mathcal{A}}}\) and \(\vvvert \cdot \vvvert_{\widetilde\omega}\)
from~\eqref{eq:weighting3} with the choice~\eqref{eq:weighting3:choice_omega} of the weights \(\widetilde\omega_1, \widetilde\omega_2 > 0\)
immediately follows from the estimates:
For all \((p, u), (q, v), (r, z) \in H(\div, \Omega) \times H^1_0(\Omega)\),
it holds that
\begin{gather}
	\label{eq:weighting3:fundamental_equivalence}
	\min\bigg\{
		\frac12,\:
		\bigg(1 + \frac{2\Lambda_1^3}{\Lambda_2^4}\bigg)^{-1}
	\bigg\} \,
	\vvvert (q, v) \vvvert_{\widetilde\omega}^2
	\leq
	\vvvert (q, v) \vvvert_{\widetilde{\mathcal{A}}}^2
	\leq
	2 \,
	\vvvert (q, v) \vvvert_{\widetilde\omega}^2,
	\\
	\label{eq:weighting3:monotonicity_weighted_norm}
	\min\bigg\{
		\frac12,\:
		\frac{\Lambda_2^2}{2\Lambda_1},\:
		\frac{\Lambda_1}{4}
	\bigg\}\,
	\vvvert (p - q, u - v) \vvvert_{\widetilde\omega}^2
	\leq
	\widetilde{\mathcal{B}}(p, u; p - q, u - v) - \widetilde{\mathcal{B}}(q, v; p - q, u - v),
	\\
	\label{eq:weighting3:Lipschitz_weighted_norm}
	\widetilde{\mathcal{B}}(p, u; r, z) - \widetilde{\mathcal{B}}(q, v; r, z)
	\leq
	4
	\max\bigg\{
		1,\:
		\frac{\Lambda_2^2}{\Lambda_1},\:
		\Lambda_2,\:
		\frac{\Lambda_1^{1/2}}{\sqrt2}
	\bigg\}\,
	\vvvert (p - q, u - v) \vvvert_{\widetilde\omega} \,
	\vvvert (r, z) \vvvert_{\widetilde\omega}.
\end{gather}

\begin{proof}[Proof of fundamental equivalence~\eqref{eq:weighting3:fundamental_equivalence}]
	The proof follows the argumentation in the proof of Theorem~\ref{thm:fundamental_equivalence}.
	It is given here in full detail for the ease of reading.

	\emph{Step~1 (lower bound).}\,
	With the binomial formula followed by an integration by parts, it follows that
	\begin{equation}
		\label{eq:weighting3:proof_lower_bound}
		\begin{aligned}
			\Vert \widetilde\omega_2^{-2} \, q \Vert_{L^2(\Omega)}^2
			+
			\Vert \nabla v \Vert_{L^2(\Omega)}^2
			&=
			\Vert \widetilde\omega_2^{-2} \, q - \nabla v \Vert_{L^2(\Omega)}^2
			+
			2 \widetilde\omega_2^{-2} \, (q, \nabla v)_{L^2(\Omega)}
			\\
			&=
			\Vert \widetilde\omega_2^{-2} \, q - \nabla v \Vert_{L^2(\Omega)}^2
			-
			2 \widetilde\omega_2^{-2} \, (\div q, v)_{L^2(\Omega)}.
		\end{aligned}
	\end{equation}
	The Cauchy--Schwarz, Friedrichs and weighted Young inequality show
    \begin{align*}
		- 2 \widetilde\omega_2^{-2} \,
		(\div q, v)_{L^2(\Omega)}
		&\leq
		2 \widetilde\omega_2^{-2} \,
		\Vert \div q \Vert_{L^2(\Omega)} \, \Vert v \Vert_{L^2(\Omega)}
		\leq
		\frac{2 C_{\textup{F}}}{\widetilde\omega_2^2} \,
		\Vert \div q \Vert_{L^2(\Omega)} \, \Vert \nabla v \Vert_{L^2(\Omega)}
		\\
		&\leq
		\frac{2C_{\textup{F}}}{\widetilde\omega_1\widetilde\omega_2^2} \,
		\Vert \widetilde\omega_1 \, \div q \Vert_{L^2(\Omega)} \, \Vert \nabla v \Vert_{L^2(\Omega)}
		\\
		&\leq
		\frac{2C_{\textup{F}}^2}{\widetilde\omega_1^2 \widetilde\omega_2^4} \,
		\Vert \widetilde\omega_1 \, \div q \Vert_{L^2(\Omega)}^2
		+
		\frac12 \,
		\Vert \nabla v \Vert_{L^2(\Omega)}^2.
    \end{align*}
	The combination with~\eqref{eq:weighting3:proof_lower_bound} and
	the absorption of \(\frac12 \, \Vert \nabla v \Vert_{L^2(\Omega)}^2\) yield
    \[
		2 \,
		\Vert \widetilde\omega_2^{-2} \, q \Vert_{L^2(\Omega)}^2
        +
		\Vert \nabla v \Vert_{L^2(\Omega)}^2
        \leq
		\frac{4C_{\textup{F}}^2}{\widetilde\omega_1^2\widetilde\omega_2^4} \,
		\Vert \widetilde\omega_1 \, \div q \Vert_{L^2(\Omega)}^2
        +
		2 \, \Vert \widetilde\omega_2^{-2} \, q - \nabla v \Vert_{L^2(\Omega)}^2.
    \]
	Adding \(C_{\textup{F}}^2 \, \Vert \widetilde\omega_1 \, \div q \Vert_{L^2(\Omega)}^2\) 
	concludes the proof of the lower bound with
    \begin{align*}
		\vvvert (q, v) \vvvert_{\widetilde\omega}^2
		&\leq
		C_{\textup{F}}^2 \,
		\Vert \widetilde\omega_1 \, \div q \Vert_{L^2(\Omega)}^2
		+
		2 \,
		\Vert \widetilde\omega_2^{-2} \, q \Vert_{L^2(\Omega)}^2
        +
		\Vert \nabla v \Vert_{L^2(\Omega)}^2
		\\
		&\leq
		\bigg(1 + \frac4{\widetilde\omega_1^2\widetilde\omega_2^4}\bigg) C_{\textup{F}}^2 \, 
		\Vert \widetilde\omega_1 \, \div q \Vert_{L^2(\Omega)}^2
        +
		2\, \Vert \widetilde\omega_2^{-2} \, q - \nabla v \Vert_{L^2(\Omega)}^2
		\\
		&\eqreff*{eq:weighting3:ls_norm}\leq
		\max\bigg\{
			2,\:
			1 + \frac4{\widetilde\omega_1^2\widetilde\omega_2^4}
		\bigg\} \, 
		\vvvert (q, v) \vvvert_{\widetilde{\mathcal{A}}}^2
		\eqreff{eq:weighting3:choice_omega}=
		\max\bigg\{
			2,\:
			1 + \frac{2 \Lambda_1^3}{\Lambda_2^4}
		\bigg\} \,
		\vvvert (q, v) \vvvert_{\widetilde{\mathcal{A}}}^2.
    \end{align*}

	\emph{Step~2 (upper bound).}\,
	The upper bound follows immediately from the triangle and the Young inequality
	\(
		\vvvert (q, v) \vvvert_{\widetilde{\mathcal{A}}}^2
        \leq
		2\,
		\vvvert (q, v) \vvvert_{\widetilde\omega}^2.
	\)
\end{proof}

\begin{proof}[Proof of monotonicity~\eqref{eq:weighting3:monotonicity_weighted_norm}]
	The proof proceeds analogously to the one of Theorem~\ref{thm:well_posedness}.
	With the modified weighting, the equality~\eqref{eq:monotonicity:split} reads
	\begin{align*}
		&\hspace{-2em}(p - q - [\sigma(\nabla u) - \sigma(\nabla v)],\: \widetilde\omega_2^{-2} \, (p - q) - \nabla(u - v))_{L^2(\Omega)}
		\\
		&=
		(p - q - M\nabla (u - v),\: \widetilde\omega_2^{-2} \, (p - q) - \nabla(u - v))_{L^2(\Omega)}
		\\
		&=
		\widetilde\omega_2^{-2} \, 
		\Vert p - q \Vert_{L^2(\Omega)}^2
		+
		(M \nabla(u - v),\: \nabla(u - v) )_{L^2(\Omega)}
		\\
		&\phantom{{}\leq{}}
		+
		(\div (p - q),\: u - v)_{L^2(\Omega)}
		-
		\widetilde\omega_2^{-2} \, 
		(p - q,\: M \nabla(u - v))_{L^2(\Omega)}.
	\end{align*}
	Applying the estimates~\eqref{eq:monotonicity:ellipticity}--\eqref{eq:monotonicity:boudedness}
	and adding \(C_{\textup{F}}^2 \Vert \widetilde\omega_1 \div (p - q) \Vert_{L^2(\Omega)}^2\) to both sides result in
	\begin{align*}
		&\hspace{-1em}
		\Big(\widetilde\omega_1^2 - \frac1{\Lambda_1} \Big) C_{\textup{F}}^2\,
		\Vert \div (p - q) \Vert_{L^2(\Omega)}^2
		+
		\frac{\widetilde\omega_2^{-2}}{2} \, \Vert p - q \Vert_{L^2(\Omega)}^2
		+
		\frac{3\Lambda_1 - 2 \widetilde\omega_2^{-2} \Lambda_2^2}4
		\, 
		\Vert \nabla(u - v) \Vert_{L^2(\Omega)}^2
		\\
		&\leq
		C_{\textup{F}}^2 \,
		\Vert \widetilde\omega_1 \div (p - q) \Vert_{L^2(\Omega)}^2
		+
		(p - q - [\sigma(\nabla u) - \sigma(\nabla v)],\: \widetilde\omega_2^{-2} \, (p - q) - \nabla(u - v))_{L^2(\Omega)}
		\\
		&=
		\widetilde{\mathcal{B}}(p, u; p - q, u - v) - \widetilde{\mathcal{B}}(q, v; p - q, u - v).
	\end{align*}
	The weights from~\eqref{eq:weighting3:choice_omega} ensure
	\begin{multline*}
		\frac{C_{\textup{F}}^2}{2} \,
		\Vert \widetilde\omega_1 \div (p - q) \Vert_{L^2(\Omega)}^2
		+
		\frac{\widetilde\omega_2^2}2 \, 
		\Vert \widetilde\omega_2^{-2} \, (p - q) \Vert_{L^2(\Omega)}^2
		+
		\frac{\Lambda_1}{4} \, 
		\Vert \nabla(u - v) \Vert_{L^2(\Omega)}^2
		\\
		\leq
		\widetilde{\mathcal{B}}(p, u; p - q, u - v) - \widetilde{\mathcal{B}}(q, v; p - q, u - v)
	\end{multline*}
	and conclude the proof of
	\begin{align*}
		\min\bigg\{
			\frac12,\:
			\frac{\Lambda_2^2}{2\Lambda_1},\:
			\frac{\Lambda_1}{4}
		\bigg\}\,
		\vvvert (p - q, u - v) \vvvert_{\widetilde\omega}^2
		&=
		\min\bigg\{
			\frac12,\:
			\frac{\widetilde\omega_2^2}{2},\:
			\frac{\Lambda_1}{4}
		\bigg\} \,
		\vvvert (p - q, u - v) \vvvert_{\widetilde\omega}^2
		\\
		&\leq
		\widetilde{\mathcal{B}}(p, u; p - q, u - v) - \widetilde{\mathcal{B}}(q, v; p - q, u - v).
		\qedhere
	\end{align*}
\end{proof}

\begin{proof}[Proof of Lipschitz continuity~\eqref{eq:weighting3:Lipschitz_weighted_norm}]
	Analogously to proof of the monotonicity \eqref{eq:weighting3:monotonicity_weighted_norm}, it follows that
	\begin{align*}
		&\hspace{-2em}
		\widetilde{\mathcal{B}}(p, u; r, z) - \widetilde{\mathcal{B}}(q, v; r, z)
		\\
		&=
		\widetilde\omega_1^2 C_{\textup{F}}^2 \, (\div (p - q),\: \div r)_{L^2(\Omega)}
		+
		(p - q - [\sigma(\nabla u) - \sigma(\nabla v)],\: \widetilde\omega_2^{-2} \, r - \nabla z)_{L^2(\Omega)}
		\\
		&=
		\widetilde\omega_1^2 C_{\textup{F}}^2 \,
		(\div (p - q),\: \div r)_{L^2(\Omega)}
		+
		\widetilde\omega_2^{-2} \, (p - q,\: r)_{L^2(\Omega)}
		+
		(M \nabla(u - v),\: \nabla z)_{L^2(\Omega)}
		\\
		&\phantom{{}\leq{}}
		+
		(\div (p - q),\: z)_{L^2(\Omega)}
		-
		\widetilde\omega_2^{-2} \, 
		(r,\: M \nabla(u - v))_{L^2(\Omega)}.
	\end{align*}
	The boundedness of \(\D\sigma\) from~\eqref{assum:boudedness} and a
	Cauchy--Schwarz inequalities in \(L^2(\Omega)\) prove
	\begin{align*}
		&\hspace{-1em}
		\widetilde{\mathcal{B}}(p, u; r, z) - \widetilde{\mathcal{B}}(q, v; r, z)
		\\
		&\leq
		\widetilde\omega_1^2 C_{\textup{F}}^2 \,
		\Vert \div (p - q) \Vert_{L^2(\Omega)} \, \Vert \div r \Vert_{L^2(\Omega)}
		+
		\widetilde\omega_2^{-2} \, \Vert p - q \Vert_{L^2(\Omega)} \, \Vert r \Vert_{L^2(\Omega)}
		\\
		&\phantom{{}\leq{}}
		+
		\Lambda_2 \,
		\Vert \nabla(u - v) \Vert_{L^2(\Omega)} \, \Vert \nabla z \Vert_{L^2(\Omega)}
		+
		C_{\textup{F}} \,
		\Vert \div (p - q) \Vert_{L^2(\Omega)} \, \Vert \nabla z \Vert_{L^2(\Omega)}
		\\
		&\phantom{{}\leq{}}
		+
		\widetilde\omega_2^{-2} \Lambda_2 \,
		\Vert r \Vert_{L^2(\Omega)} \, \Vert \nabla(u - v) \Vert_{L^2(\Omega)}
		\\
		&=
		C_{\textup{F}}^2 \,
		\Vert \widetilde\omega_1 \div (p - q) \Vert_{L^2(\Omega)} \, \Vert \widetilde\omega_1 \div r \Vert_{L^2(\Omega)}
		+
		\widetilde\omega_2^2 \,
		\Vert \widetilde\omega_2^{-2} \, (p - q) \Vert_{L^2(\Omega)} \,
		\Vert \widetilde\omega_2^{-2} \, r \Vert_{L^2(\Omega)}
		\\
		&\phantom{{}\leq{}}
		+
		\Lambda_2 \,
		\Vert \nabla(u - v) \Vert_{L^2(\Omega)} \, \Vert \nabla z \Vert_{L^2(\Omega)}
		+
		\frac{1}{\widetilde\omega_1}
		C_{\textup{F}} \,
		\Vert \widetilde\omega_1 \div (p - q) \Vert_{L^2(\Omega)} \, \Vert \nabla z \Vert_{L^2(\Omega)}
		\\
		&\phantom{{}\leq{}}
		+
		\Lambda_2 \,
		\Vert \widetilde\omega_2^{-2} \, r \Vert_{L^2(\Omega)} \, \Vert \nabla(u - v) \Vert_{L^2(\Omega)}.
	\end{align*}
	A Cauchy--Schwarz inequality in \(\R^5\) results in
	\begin{align*}
		&\hspace{-2em}
		\widetilde{\mathcal{B}}(p, u; r, z) - \widetilde{\mathcal{B}}(q, v; r, z)
		\\
		&\leq
		\max\bigg\{
			1,\:
			\widetilde\omega_2^2,\:
			\Lambda_2,\:
			\frac{1}{\widetilde\omega_1}
		\bigg\}
		\\
		&\phantom{{}\leq{}}
		\times
		\Big[
			2 C_{\textup{F}}^2 \,
			\Vert \widetilde\omega_1 \div (p - q) \Vert_{L^2(\Omega)}^2
			+
			\Vert \widetilde\omega_2^{-2} (p - q) \Vert_{L^2(\Omega)}^2
			+
			2 \Vert \nabla(u - v) \Vert_{L^2(\Omega)}^2
		\Big]^{1/2}
		\\
		&\phantom{{}\leq{}}
		\times
		\Big[
			C_{\textup{F}}^2 \,
			\Vert \widetilde\omega_1 \div r \Vert_{L^2(\Omega)}^2
			+
			2 \,
			\Vert \widetilde\omega_2^{-2} \, r \Vert_{L^2(\Omega)}^2
			+
			2 \,
			\Vert \nabla z \Vert_{L^2(\Omega)}^2
		\Big]^{1/2}.
	\end{align*}
	This concludes the proof of
	\begin{align*}
		\widetilde{\mathcal{B}}(p, u; r, z) - \widetilde{\mathcal{B}}(q, v; r, z)
		&\leq
		2
		\max\bigg\{
			1,\:
			\widetilde\omega_2^2,\:
			\Lambda_2,\:
			\frac{1}{\widetilde\omega_1}
		\bigg\}\,
		\vvvert (p - q, u - v) \vvvert_{\widetilde\omega} \,
		\vvvert (r, z) \vvvert_{\widetilde\omega}
		\\
		&\eqreff*{eq:weighting3:choice_omega}=
		2
		\max\bigg\{
			1,\:
			\frac{\Lambda_2^2}{\Lambda_1},\:
			\Lambda_2,\:
			\frac{\Lambda_1^{1/2}}{\sqrt{2}}
		\bigg\} \,
		\vvvert (p - q, u - v) \vvvert_{\widetilde\omega} \,
		\vvvert (r, z) \vvvert_{\widetilde\omega}.
		\qedhere
	\end{align*}
\end{proof}

\section{Weighting 4: Split weighting}
\label{app:split_weighting}
\noindent
This appendix proves the strong monotonicity and Lipschitz continuity~\eqref{eq:alternative_weightings:monotonicity_lipschitz}
presented in Subsection~\ref{sec:split_weighting} for the nonlinear mapping \(\widetilde{\mathcal{B}}\) 
and the norms \(\vvvert \cdot \vvvert_{\widetilde{\mathcal{A}}}\) and \(\vvvert \cdot \vvvert_{\widetilde\omega}\) from~\eqref{eq:weighting4}.
The choice~\eqref{eq:weighting4:choice_omega} ensures~\eqref{eq:alternative_weightings:monotonicity_lipschitz}
with the constants~\eqref{eq:weighting4:constants}.
This immediately follows from the estimates:
For all \((p, u), (q, v), (r, z) \in H(\div, \Omega) \times H^1_0(\Omega)\), it holds that
\begin{gather}
	\label{eq:weighting4:fundamental_equivalence}
	\min\bigg\{
		\frac12,\:
		\bigg(1 + \frac{2 \Lambda_1^3}{\Lambda_2^2}\bigg)^{-1}
	\bigg\} \,
	\vvvert (q, v) \vvvert_{\widetilde\omega}^2
	\leq
	\vvvert (q, v) \vvvert_{\widetilde{\mathcal{A}}}^2
	\leq
	2\,
	\vvvert (q, v) \vvvert_{\widetilde\omega}^2,
	\\
	\label{eq:weighting4:monotonicity_weighted_norm}
	\min\bigg\{
		\frac12,\:
		\frac1{2\Lambda_1},\:
		\frac{\Lambda_1}{4\Lambda_2^2}
	\bigg\}\,
	\vvvert (p - q, u - v) \vvvert_{\widetilde\omega}^2
	\leq
	\widetilde{\mathcal{B}}(p, u; p - q, u - v) - \widetilde{\mathcal{B}}(q, v; p - q, u - v),
	\\
	\label{eq:weighting4:Lipschitz_weighted_norm}
	\widetilde{\mathcal{B}}(p, u; r, z) - \widetilde{\mathcal{B}}(q, v; r, z)
	\leq
	4
	\max\bigg\{
		1,\:
		\frac1{\Lambda_1},\:
		\frac1{\Lambda_2},\:
		\frac{\Lambda_1}{2\Lambda_2^2}
	\bigg\} \,
	\vvvert (p - q, u - v) \vvvert_{\widetilde\omega} \,
	\vvvert (r, z) \vvvert_{\widetilde\omega}.
\end{gather}

\begin{proof}[Proof of fundamental equivalence~\eqref{eq:weighting4:fundamental_equivalence}]
	The proof employs the arguments the proof of Theorem~\ref{thm:fundamental_equivalence}.
	They are given here in full detail for the ease of reading.

	\emph{Step~1 (lower bound).}\,
	To begin with, the binomial formula and an integration by parts provide
	\begin{equation}
		\label{eq:weighting4:proof_lower_bound}
		\begin{aligned}
			\Vert \Lambda_1 \, q \Vert_{L^2(\Omega)}^2
			+
			\Vert \Lambda_2^2 \, \nabla v \Vert_{L^2(\Omega)}^2
			&=
			\Vert \Lambda_1 \, q - \Lambda_2^2 \, \nabla v \Vert_{L^2(\Omega)}^2
			+
			2\Lambda_1\Lambda_2^2 \,
			(q, \nabla v)_{L^2(\Omega)}
			\\
			&=
			\Vert \Lambda_1 \, q - \Lambda_2^2 \, \nabla v \Vert_{L^2(\Omega)}^2
			-
			2\Lambda_1\Lambda_2^2 \,
			(\div q, v)_{L^2(\Omega)}.
		\end{aligned}
	\end{equation}
	The Cauchy--Schwarz, Friedrichs and weighted Young inequality establish
    \begin{align*}
		- 2 \Lambda_1 \Lambda_2^2 \,
		(\div q, v)_{L^2(\Omega)}
		&\leq
		2 \Lambda_1 \Lambda_2^2 \,
		\Vert \div q \Vert_{L^2(\Omega)} \, \Vert v \Vert_{L^2(\Omega)}
		\leq
		2 C_{\textup{F}} \,
		\Vert \div q \Vert_{L^2(\Omega)} \, \Vert \Lambda_2^2 \, \nabla v \Vert_{L^2(\Omega)}
		\\
		&\leq
		\frac{2\Lambda_1^2 C_{\textup{F}}^2}{\widetilde\omega_1^2} \,
		\Vert \widetilde\omega_1 \, \div q \Vert_{L^2(\Omega)}^2
		+
		\frac12 \,
		\Vert \Lambda_2^2 \, \nabla v \Vert_{L^2(\Omega)}^2.
    \end{align*}
	This, the binomial formula from~\eqref{eq:weighting4:proof_lower_bound},
	and the absorption of \(\frac12 \, \Vert \Lambda_2^2 \, \nabla v \Vert_{L^2(\Omega)}^2\) lead to
    \[
		2 \,
		\Vert \Lambda_1 \, q \Vert_{L^2(\Omega)}^2
        +
		\Vert \Lambda_2^2 \, \nabla v \Vert_{L^2(\Omega)}^2
        \leq
		\frac{4\Lambda_1^2 C_{\textup{F}}^2}{\widetilde\omega_1^2} \,
		\Vert \widetilde\omega_1 \, \div q \Vert_{L^2(\Omega)}^2
        +
		2 \, \Vert \Lambda_1 \, q - \Lambda_2^2 \, \nabla v \Vert_{L^2(\Omega)}^2.
    \]
	The addition of \(C_{\textup{F}}^2 \, \Vert \widetilde\omega_1 \, \div q \Vert_{L^2(\Omega)}^2\) 
	concludes the proof of the lower bound via
    \begin{align*}
		\vvvert (q, v) \vvvert^2
		&\leq
		C_{\textup{F}}^2 \,
		\Vert \widetilde\omega_1 \, \div q \Vert_{L^2(\Omega)}^2
		+
		\Vert \Lambda_1 \, q \Vert_{L^2(\Omega)}^2
        +
		\Vert \Lambda_2^2 \, \nabla v \Vert_{L^2(\Omega)}^2
		\\
		&\leq
		\bigg(1 + \frac{4\Lambda_1^2}{\widetilde\omega_1^2}\bigg) \, C_{\textup{F}}^2 \, 
		\Vert \widetilde\omega_1 \, \div q \Vert_{L^2(\Omega)}^2
        +
		2\, \Vert \Lambda_1 \, q - \Lambda_2^2 \, \nabla v \Vert_{L^2(\Omega)}^2
		\\
		&\eqreff*{eq:weighting4:ls_norm}\leq
		\max\bigg\{ 
			2,\:
			1 + \frac{4\Lambda_1^2}{\widetilde\omega_1^2}
		\bigg\} \, 
		\vvvert (q, v) \vvvert_{\mathcal{A}}^2
		\eqreff{eq:weighting4:choice_omega}=
		\max\bigg\{
			2,\:
			1 + \frac{2 \Lambda_1^3}{\Lambda_2^2}
		\bigg\} \,
		\vvvert (q, v) \vvvert_{\mathcal{A}}^2.
    \end{align*}

	\emph{Step~2 (upper bound).}\,
	Analogously to the proof of Theorem~\ref{thm:fundamental_equivalence}, the triangle and the Young inequality
	prove
	\(
		\vvvert (q, v) \vvvert_{\mathcal{A}}^2
        \leq
		2\,
		\vvvert (q, v) \vvvert^2.
	\)
\end{proof}

\begin{proof}[Proof of monotonicity~\eqref{eq:weighting4:monotonicity_weighted_norm}]
	The proof proceeds analogously to Theorem~\ref{thm:well_posedness}.
	With the changed weights, the equality~\eqref{eq:monotonicity:split} reads
	\begin{align*}
		&\hspace{-2em}(p - q - [\sigma(\nabla u) - \sigma(\nabla v)],\: \Lambda_1 \, (p - q) - \Lambda_2^2 \, \nabla(u - v))_{L^2(\Omega)}
		\\
		&=
		(p - q - M\nabla (u - v),\: \Lambda_1 \, (p - q) - \Lambda_2^2 \, \nabla(u - v))_{L^2(\Omega)}
		\\
		&=
		\Lambda_1 \, \Vert p - q \Vert_{L^2(\Omega)}^2
		+
		\Lambda_2^2 \,
		(M \nabla(u - v),\: \nabla(u - v) )_{L^2(\Omega)}
		\\
		&\phantom{{}\leq{}}
		+
		\Lambda_2^2 \,
		(\div (p - q),\: u - v)_{L^2(\Omega)}
		-
		\Lambda_1 \, 
		(p - q,\: M \nabla(u - v))_{L^2(\Omega)}.
	\end{align*}
	This, the estimates~\eqref{eq:monotonicity:ellipticity}--\eqref{eq:monotonicity:boudedness},
	and adding \(C_{\textup{F}}^2 \Vert \widetilde\omega_1 \div (p  - q) \Vert_{L^2(\Omega)}^2\) results in
	\begin{align*}
		&\hspace{-0.5em}
		\bigg(\widetilde\omega_1^2 - \frac{\Lambda_2^2}{\Lambda_1} \bigg) C_{\textup{F}}^2\,
		\Vert \div (p - q) \Vert_{L^2(\Omega)}^2
		+
		\frac{\Lambda_1}{2} \, \Vert p - q \Vert_{L^2(\Omega)}^2
		+
		\frac{\Lambda_1 \Lambda_2^2}4
		\, 
		\Vert \nabla(u - v) \Vert_{L^2(\Omega)}^2
		\\
		&\leq
		C_{\textup{F}}^2 \,
		\Vert \widetilde\omega_1 \div (p - q) \Vert_{L^2(\Omega)}^2
		+
		(p - q - [\sigma(\nabla u) - \sigma(\nabla v)],\: \Lambda_1 (p - q) - \Lambda_2^2 \, \nabla(u - v))_{L^2(\Omega)}
		\\
		&=
		\widetilde{\mathcal{B}}(p, u; p - q, u - v) - \widetilde{\mathcal{B}}(q, v; p - q, u - v).
	\end{align*}
	The choice of the weights in~\eqref{eq:weighting4:choice_omega} proves
	\begin{multline*}
		\frac{C_{\textup{F}}^2}{2} \,
		\Vert \widetilde\omega_1 \div (p - q) \Vert_{L^2(\Omega)}^2
		+
		\frac1{2\Lambda_1} \, 
		\Vert \Lambda_1 \, (p - q) \Vert_{L^2(\Omega)}^2
		+
		\frac{\Lambda_1}{4\Lambda_2^2} \, 
		\Vert \Lambda_2^2 \, \nabla(u - v) \Vert_{L^2(\Omega)}^2
		\\
		\leq
		\widetilde{\mathcal{B}}(p, u; p - q, u - v) - \widetilde{\mathcal{B}}(q, v; p - q, u - v)
	\end{multline*}
	and concludes the proof of
	\[
		\min\bigg\{
			\frac12,\:
			\frac1{2\Lambda_1},\:
			\frac{\Lambda_1}{4\Lambda_2^2}
		\bigg\} \,
		\vvvert (p - q, u - v) \vvvert^2
		\leq
		\widetilde{\mathcal{B}}(p, u; p - q, u - v) - \widetilde{\mathcal{B}}(q, v; p - q, u - v).
		\qedhere
	\]
\end{proof}

\begin{proof}[Proof of Lipschitz continuity~\eqref{eq:weighting4:Lipschitz_weighted_norm}]
	As in the proof of the monotonicity~\eqref{eq:weighting4:monotonicity_weighted_norm},
	it holds that
	\begin{align*}
		&\hspace{-2em}
		\widetilde{\mathcal{B}}(p, u; r, z) - \widetilde{\mathcal{B}}(q, v; r, z)
		\\
		&=
		\widetilde\omega_1^2 C_{\textup{F}}^2 \, (\div (p - q),\: \div r)_{L^2(\Omega)}
		+
		(p - q - [\sigma(\nabla u) - \sigma(\nabla v)],\: \Lambda_1 \, r - \Lambda_2^2 \, \nabla z)_{L^2(\Omega)}
		\\
		&=
		\widetilde\omega_1^2 C_{\textup{F}}^2 \,
		(\div (p - q),\: \div r)_{L^2(\Omega)}
		+
		\Lambda_1 \,
		(p - q,\: r)_{L^2(\Omega)}
		+
		\Lambda_2^2 \,
		(M \nabla(u - v),\: \nabla z)_{L^2(\Omega)}
		\\
		&\phantom{{}\leq{}}
		+
		\Lambda_2^2 \,
		(\div (p - q),\: z)_{L^2(\Omega)}
		-
		\Lambda_1 \,
		(r,\: M \nabla(u - v))_{L^2(\Omega)}.
	\end{align*}
	The boundedness of \(\D\sigma\) from~\eqref{assum:boudedness} and a
	Cauchy--Schwarz inequalities in \(L^2(\Omega)\) verify
	\begin{align*}
		&\hspace{-1em}
		\widetilde{\mathcal{B}}(p, u; r, z) - \widetilde{\mathcal{B}}(q, v; r, z)
		\\
		&\leq
		\widetilde\omega_1^2 C_{\textup{F}}^2 \,
		\Vert \div (p - q) \Vert_{L^2(\Omega)} \, \Vert \div r \Vert_{L^2(\Omega)}
		+
		\Lambda_1 \,
		\Vert p - q \Vert_{L^2(\Omega)} \, \Vert r \Vert_{L^2(\Omega)}
		\\
		&\phantom{{}\leq{}}
		+
		\Lambda_2^3 \,
		\Vert \nabla(u - v) \Vert_{L^2(\Omega)} \, \Vert \nabla z \Vert_{L^2(\Omega)}
		+
		\Lambda_2^2 C_{\textup{F}} \,
		\Vert \div (p - q) \Vert_{L^2(\Omega)} \, \Vert \nabla z \Vert_{L^2(\Omega)}
		\\
		&\phantom{{}\leq{}}
		+
		\Lambda_1 \Lambda_2 \,
		\Vert r \Vert_{L^2(\Omega)} \, \Vert \nabla(u - v) \Vert_{L^2(\Omega)}
		\\
		&=
		C_{\textup{F}}^2 \,
		\Vert \widetilde\omega_1 \, \div (p - q) \Vert_{L^2(\Omega)} \, \Vert \widetilde\omega_1 \, \div r \Vert_{L^2(\Omega)}
		+
		\frac1{\Lambda_1} \,
		\Vert \Lambda_1 \, (p - q) \Vert_{L^2(\Omega)} \, \Vert \Lambda_1 \, r \Vert_{L^2(\Omega)}
		\\
		&\phantom{{}\leq{}}
		+
		\frac1{\Lambda_2} \,
		\Vert \Lambda_2^2 \, \nabla(u - v) \Vert_{L^2(\Omega)} \, \Vert \Lambda_2^2 \, \nabla z \Vert_{L^2(\Omega)}
		+
		\frac{C_{\textup{F}}}{\widetilde\omega_1} \,
		\Vert \widetilde\omega_1 \, \div (p - q) \Vert_{L^2(\Omega)} \, \Vert \Lambda_2^2 \, \nabla z \Vert_{L^2(\Omega)}
		\\
		&\phantom{{}\leq{}}
		+
		\frac1{\Lambda_2} \,
		\Vert \Lambda_1 \, r \Vert_{L^2(\Omega)} \, \Vert \Lambda_2^2 \, \nabla(u - v) \Vert_{L^2(\Omega)}.
	\end{align*}
	A Cauchy--Schwarz inequality in \(\R^5\) results in
	\begin{align*}
		&\hspace{-1em}
		\widetilde{\mathcal{B}}(p, u; r, z) - \widetilde{\mathcal{B}}(q, v; r, z)
		\\
		&\leq
		\max\bigg\{
			1,\:
			\frac1{\Lambda_1},\:
			\frac1{\Lambda_2},\:
			\frac{1}{\widetilde\omega_1}
		\bigg\}
		\\
		&\phantom{{}\leq{}}
		\times
		\Big[
			2 C_{\textup{F}}^2 \,
			\Vert \widetilde\omega_1 \div (p - q) \Vert_{L^2(\Omega)}^2
			+
			\Vert p - q \Vert_{L^2(\Omega)}^2
			+
			2 \Vert \widetilde\omega_2^2 \, \nabla(u - v) \Vert_{L^2(\Omega)}^2
		\Big]^{1/2}
		\\
		&\phantom{{}\leq{}}
		\times
		\Big[
			C_{\textup{F}}^2 \,
			\Vert \widetilde\omega_1 \div r \Vert_{L^2(\Omega)}^2
			+
			2 \,
			\Vert r \Vert_{L^2(\Omega)}^2
			+
			2 \,
			\Vert \widetilde\omega_2^2 \, \nabla z \Vert_{L^2(\Omega)}^2
		\Big]^{1/2}.
	\end{align*}
	This concludes the proof of
	\begin{align*}
		\widetilde{\mathcal{B}}(p, u; r, z) - \widetilde{\mathcal{B}}(q, v; r, z)
		&\leq
		2
		\max\bigg\{
			1,\:
			\frac1{\Lambda_1},\:
			\frac1{\Lambda_2},\:
			\frac{1}{\widetilde\omega_1}
		\bigg\}\,
		\vvvert (p - q, u - v) \vvvert \,
		\vvvert (r, z) \vvvert
		\\
		&\eqreff*{eq:weighting4:choice_omega}=
		2
		\max\bigg\{
			1,\:
			\frac1{\Lambda_1},\:
			\frac1{\Lambda_2},\:
			\frac{\Lambda_1^{1/2}}{\sqrt{2}\Lambda_2}
		\bigg\}\,
		\vvvert (p - q, u - v) \vvvert \,
		\vvvert (r, z) \vvvert.
		\qedhere
	\end{align*}
\end{proof}

\end{document}

%% file: figures/Fig01_Nonlinear_errors.tex
\begin{tikzpicture}[>=stealth]
    \begin{loglogaxis}[%
            width            = 5.6cm,%
            xlabel           = {cumulative ndof},%
            ylabel           = {error estimator},%
            ymajorgrids      = true,%
            font             = \footnotesize,%
            grid style       = {%
                densely dotted,%
                semithick%
            },%
            legend style     = {%
                legend pos = north east,%
                font = \footnotesize%
            },%
        ]

        \addlegendimage{marker1}
        \addlegendentry{\(\eta^k_\ell\)}
        \addlegendimage{marker2}
        \addlegendentry{\(\mu^k_\ell\)}
        \addlegendimage{marker3}
        \addlegendentry{\(N^k_\ell\)}


        \pgfplotstableread[col sep=comma]{
       k,     ell,    ndof,   nElem,        eta,         mu,        res,    case,   maxGradU
       0,       0,     193,      96,2.01791e-01,1.16608e+00,2.55303e-01,       Z,8.62011e-03
       1,       0,     193,      96,2.01405e-01,1.51358e-01,2.13723e-01,       R,1.60043e-02
       1,       1,     271,     135,1.63298e-01,1.91851e-01,1.79901e-01,       Z,2.95074e-02
       2,       0,     271,     135,1.63037e-01,7.16045e-02,1.67435e-01,       Z,4.14307e-02
       3,       0,     271,     135,1.62901e-01,2.90517e-02,1.65264e-01,       R,4.69924e-02
       3,       1,     427,     213,1.39329e-01,8.92647e-02,1.42920e-01,       Z,7.06245e-02
       4,       0,     427,     213,1.39178e-01,2.65099e-02,1.40792e-01,       R,8.57288e-02
       4,       1,     633,     316,1.12356e-01,8.63095e-02,1.16536e-01,       Z,1.20423e-01
       5,       0,     633,     316,1.12328e-01,2.82196e-02,1.13647e-01,       Z,1.47185e-01
       6,       0,     633,     316,1.12315e-01,1.15980e-02,1.13149e-01,       R,1.60334e-01
       6,       1,     937,     468,9.12346e-02,6.65240e-02,9.38424e-02,       Z,1.62853e-01
       7,       0,     937,     468,9.12420e-02,2.03783e-02,9.19383e-02,       Z,1.70365e-01
       8,       0,     937,     468,9.12472e-02,7.78354e-03,9.16595e-02,       R,1.73925e-01
       8,       1,    1309,     654,7.54564e-02,5.18938e-02,7.74173e-02,       Z,2.33313e-01
       9,       0,    1309,     654,7.54574e-02,1.63962e-02,7.59310e-02,       Z,2.67113e-01
      10,       0,    1309,     654,7.54574e-02,6.39634e-03,7.57034e-02,       R,2.84215e-01
      10,       1,    2137,    1068,6.08165e-02,4.51230e-02,6.24039e-02,       Z,3.76092e-01
      11,       0,    2137,    1068,6.08435e-02,1.33715e-02,6.12043e-02,       Z,4.36950e-01
      12,       0,    2137,    1068,6.08539e-02,5.40803e-03,6.10093e-02,       R,4.68905e-01
      12,       1,    3089,    1544,4.92090e-02,3.62065e-02,5.05179e-02,       Z,6.00485e-01
      13,       0,    3089,    1544,4.92085e-02,1.11692e-02,4.94582e-02,       Z,6.94138e-01
      14,       0,    3089,    1544,4.92082e-02,4.33988e-03,4.92970e-02,       R,7.44808e-01
      14,       1,    4589,    2294,4.01562e-02,2.87710e-02,4.11967e-02,       Z,9.64786e-01
      15,       0,    4589,    2294,4.01500e-02,9.09263e-03,4.03372e-02,       Z,1.12408e+00
      16,       0,    4589,    2294,4.01479e-02,3.57172e-03,4.02034e-02,       R,1.21258e+00
      16,       1,    7187,    3593,3.25030e-02,2.38362e-02,3.33021e-02,       Z,1.56129e+00
      17,       0,    7187,    3593,3.25053e-02,7.16600e-03,3.26503e-02,       Z,1.82386e+00
      18,       0,    7187,    3593,3.25067e-02,2.87589e-03,3.25452e-02,       R,1.97272e+00
      18,       1,   10419,    5209,2.68242e-02,1.85855e-02,2.74027e-02,       Z,2.52575e+00
      19,       0,   10419,    5209,2.68268e-02,5.54342e-03,2.69252e-02,       Z,2.95377e+00
      20,       0,   10419,    5209,2.68280e-02,2.17301e-03,2.68520e-02,       R,3.20017e+00
      20,       1,   15673,    7836,2.18321e-02,1.57425e-02,2.23725e-02,       Z,3.24298e+00
      21,       0,   15673,    7836,2.18306e-02,4.87122e-03,2.19192e-02,       Z,3.40512e+00
      22,       0,   15673,    7836,2.18303e-02,1.91498e-03,2.18486e-02,       R,3.49461e+00
      22,       1,   23249,   11624,1.78925e-02,1.26526e-02,1.83225e-02,       Z,4.42420e+00
      23,       0,   23249,   11624,1.78942e-02,3.92574e-03,1.79647e-02,       R,5.04005e+00
      23,       1,   34811,   17405,1.47452e-02,1.08718e-02,1.51273e-02,       Z,6.43375e+00
      24,       0,   34811,   17405,1.47492e-02,3.35415e-03,1.48135e-02,       Z,7.70067e+00
      25,       0,   34811,   17405,1.47509e-02,1.34287e-03,1.47634e-02,       R,8.44703e+00
      25,       1,   49807,   24903,1.22872e-02,8.27143e-03,1.25374e-02,       Z,1.06775e+01
      26,       0,   49807,   24903,1.22882e-02,2.48284e-03,1.23277e-02,       Z,1.25320e+01
      27,       0,   49807,   24903,1.22886e-02,9.69309e-04,1.22958e-02,       R,1.36208e+01
      27,       1,   71459,   35729,1.01768e-02,6.95582e-03,1.04136e-02,       Z,1.69306e+01
      28,       0,   71459,   35729,1.01764e-02,2.20749e-03,1.02143e-02,       Z,1.96559e+01
      29,       0,   71459,   35729,1.01763e-02,8.74291e-04,1.01829e-02,       R,2.12546e+01
      29,       1,  103867,   51933,8.51327e-03,5.64327e-03,8.68291e-03,       Z,2.69515e+01
      30,       0,  103867,   51933,8.51369e-03,1.70441e-03,8.54125e-03,       Z,3.14958e+01
      31,       0,  103867,   51933,8.51393e-03,6.78939e-04,8.51880e-03,       R,3.41693e+01
      31,       1,  148269,   74134,7.12601e-03,4.70829e-03,7.25922e-03,       Z,3.45545e+01
      32,       0,  148269,   74134,7.12693e-03,1.37864e-03,7.14859e-03,       R,3.63128e+01
      32,       1,  207617,  103808,6.02299e-03,4.05184e-03,6.15480e-03,       Z,4.56736e+01
      33,       0,  207617,  103808,6.02330e-03,1.26496e-03,6.04426e-03,       Z,5.25239e+01
      34,       0,  207617,  103808,6.02347e-03,4.99684e-04,6.02701e-03,       R,5.65371e+01
      34,       1,  290993,  145496,5.05900e-03,3.30732e-03,5.16659e-03,       Z,6.98422e+01
      35,       0,  290993,  145496,5.05879e-03,1.04959e-03,5.07584e-03,       R,8.05348e+01
      35,       1,  407169,  203584,4.29572e-03,2.87049e-03,4.38946e-03,       Z,1.02460e+02
      36,       0,  407169,  203584,4.29585e-03,9.01526e-04,4.31101e-03,       Z,1.23282e+02
      37,       0,  407169,  203584,4.29592e-03,3.59901e-04,4.29850e-03,       R,1.35684e+02
      37,       1,  572869,  286434,3.62967e-03,2.32595e-03,3.69334e-03,       Z,1.37730e+02
      38,       0,  572869,  286434,3.63005e-03,6.80692e-04,3.64041e-03,       Z,1.46105e+02
      39,       0,  572869,  286434,3.63022e-03,2.71917e-04,3.63200e-03,       R,1.50865e+02
      39,       1,  787563,  393781,3.09451e-03,1.91740e-03,3.14701e-03,       Z,1.87968e+02
      40,       0,  787563,  393781,3.09452e-03,5.72286e-04,3.10273e-03,       R,2.13507e+02
      40,       1, 1097597,  548798,2.61659e-03,1.74844e-03,2.67739e-03,       Z,2.65291e+02
        }\tableData

        %
        %

        \addplot+ [marker1, adaptive, forget plot]
        table [col sep=comma, x=cumulativeNdof, y=eta] {\tableData};

        \addplot+ [marker2, adaptive, forget plot]
        table [col sep=comma, x=cumulativeNdof, y=mu] {\tableData};

        \addplot+ [marker3, adaptive, forget plot]
        table [col sep=comma, x=cumulativeNdof, y=res] {\tableData};

        \drawslopetriangle[ST1]{0.5}{5e4}{7e-4} 
        \drawswappedslopetriangle[ST2]{0.5}{6e6}{1.5e-2} 
    \end{loglogaxis}
\end{tikzpicture}

%% file: figures/Fig01_Nonlinear_mesh.tex
\begin{tikzpicture}
    \begin{axis}[%
        axis equal image,%
        width=6.1cm,%
        xmin=-1.1, xmax=1.1,%
        ymin=-1.1, ymax=1.1,%
        font=\footnotesize%
    ]
        \addplot graphics [xmin=-1, xmax=1, ymin=-1, ymax=1]
            {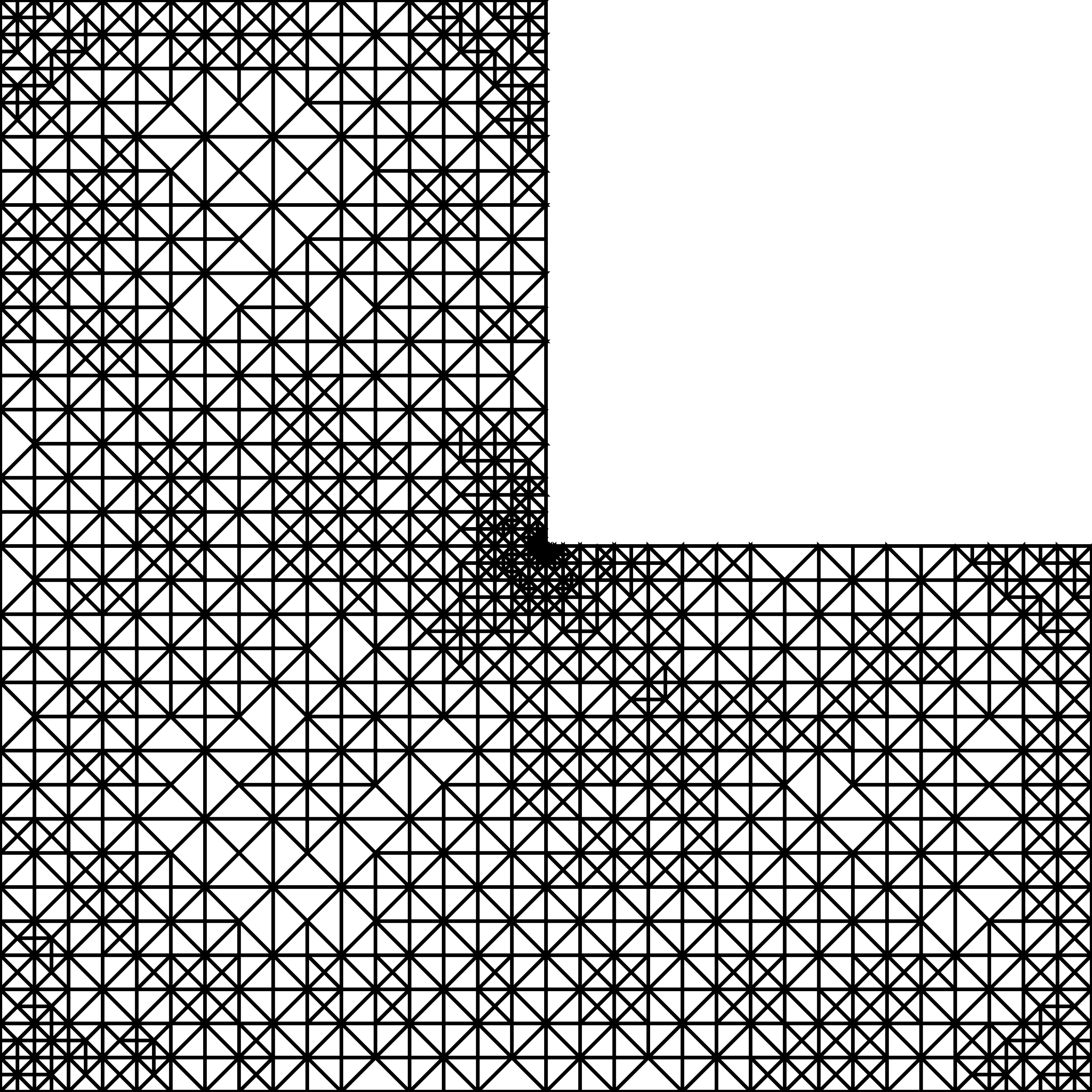};
    \end{axis}
\end{tikzpicture}

%% file: figures/Fig01_Nonlinear_potential.tex
\begin{tikzpicture}
    \begin{axis}[%
        axis equal image,%
        width=6.1cm,%
        xmin=-1.1, xmax=1.1,%
        ymin=-1.1, ymax=1.1,%
        font=\footnotesize,%
        point meta min=0.0,%
        point meta max=5.77470135e-04,%
        colorbar,%
        colorbar style={%
            font=\footnotesize,%
            width=2.5mm,%
            title style={yshift=-2mm},%
        },%
    ]
        \addplot graphics [xmin=-1, xmax=1, ymin=-1, ymax=1] {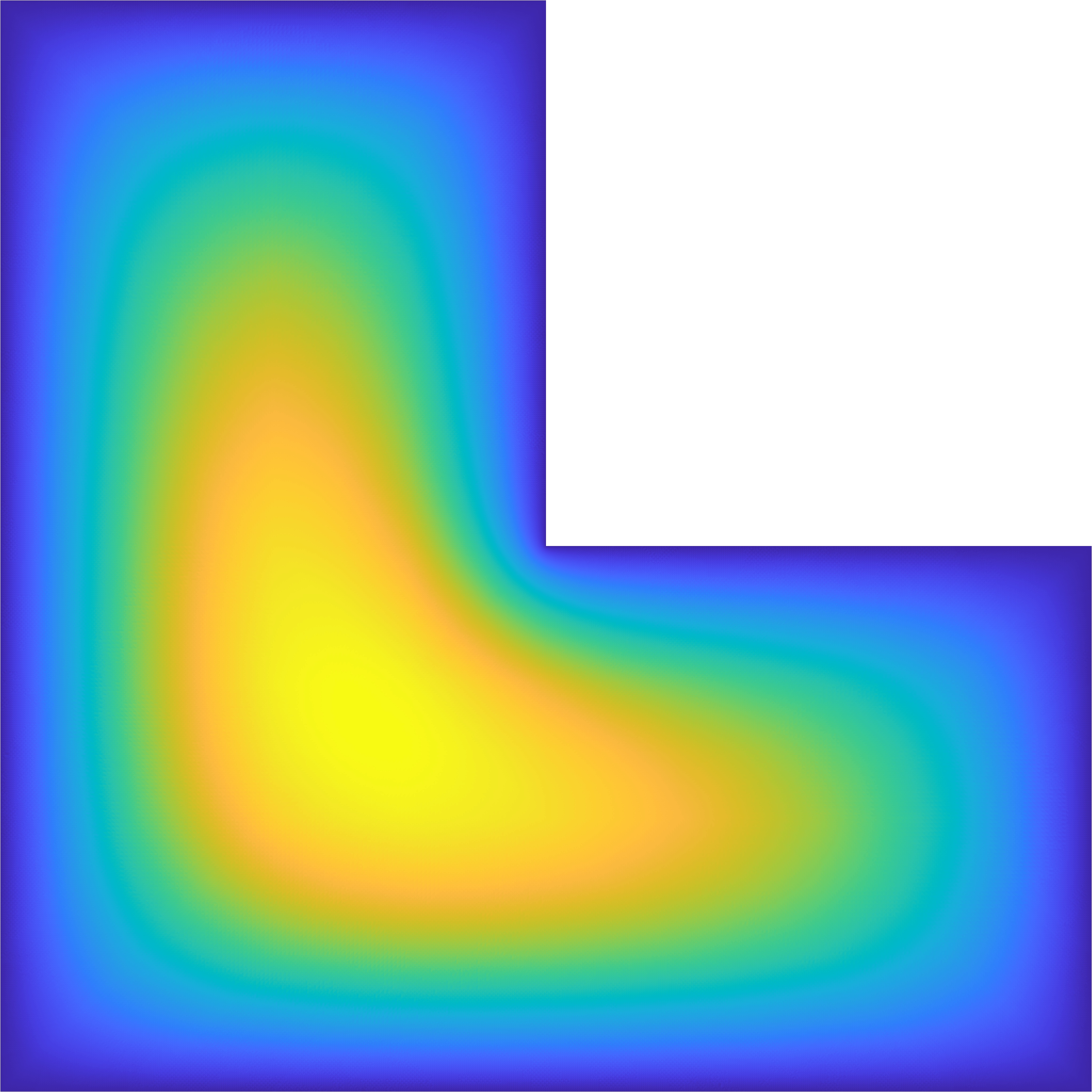};
    \end{axis}
\end{tikzpicture}

%% file: figures/Fig01_Nonlinear_flux.tex
\begin{tikzpicture}
    \begin{axis}[%
        axis equal image,%
        width=6.1cm,%
        xmin=-1.1, xmax=1.1,%
        ymin=-1.1, ymax=1.1,%
        font=\footnotesize,%
        point meta min=0.0,%
        point meta max=8.57313746e-03,%
        colorbar,%
        colorbar style={%
            font=\footnotesize,%
            width=2.5mm,%
            title style={yshift=-2mm},%
        },%
    ]
        \addplot graphics [xmin=-1, xmax=1, ymin=-1, ymax=1]
{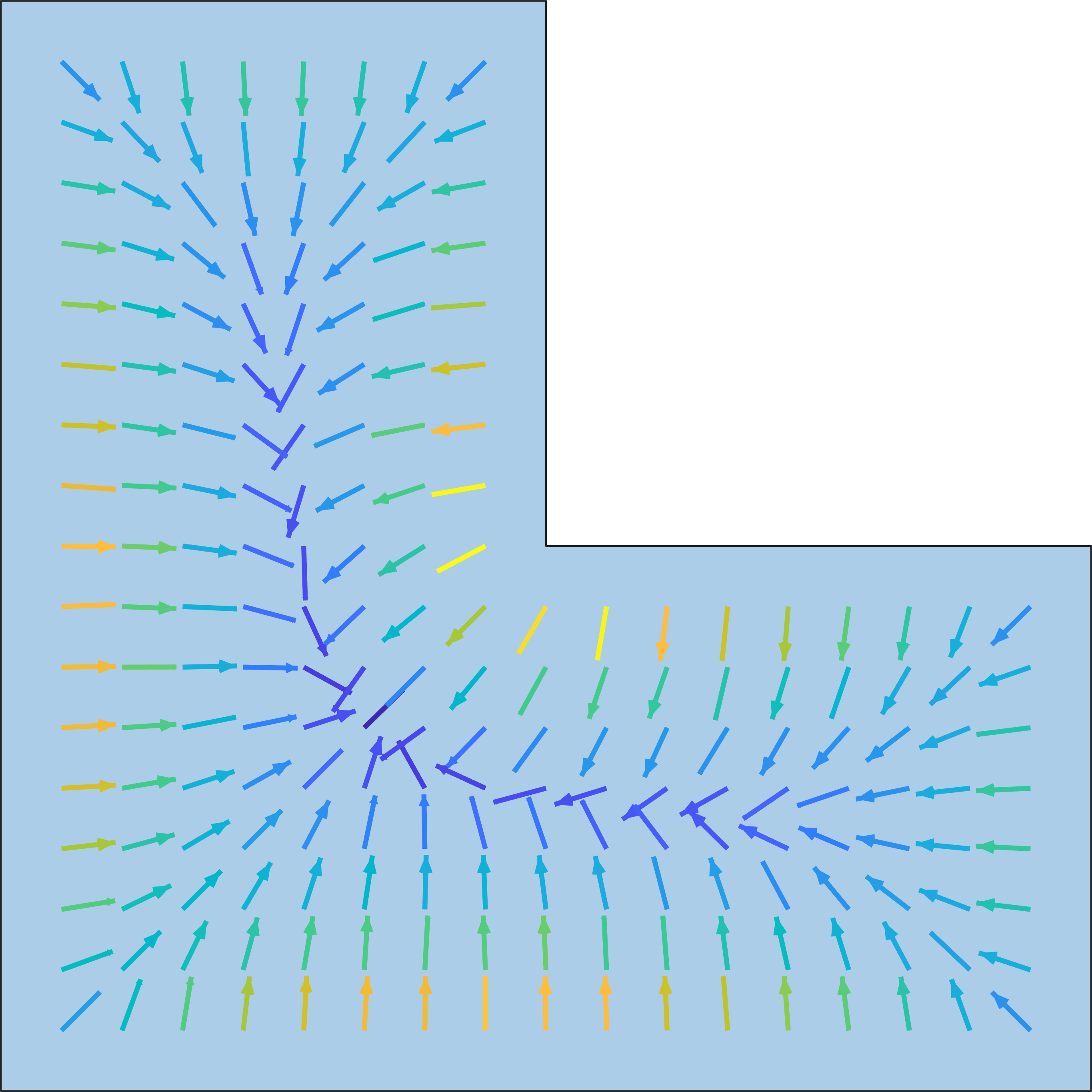};
    \end{axis}
\end{tikzpicture}

%% file: figures/Fig02a_Nonlinear_convergence_gamma.tex
\begin{tikzpicture}[>=stealth]
    \begin{loglogaxis}[%
            width            = 5.5cm,%
            xlabel           = {cumulative ndof},%
            ylabel           = {error estimator},%
            ymajorgrids      = true,%
            ymin             = 2e-3,%
            ymax             = 5,%
            font             = \footnotesize,%
            grid style       = {%
                densely dotted,%
                semithick%
            },%
            legend style     = {%
                legend pos = north east,%
                font = \footnotesize%
            },%
        ]

        \addlegendimage{marker1}
        \addlegendentry{\(\gamma = 0.9\)}
        \addlegendimage{marker2}
        \addlegendentry{\(\gamma = 0.7\)}
        \addlegendimage{marker3}
        \addlegendentry{\(\gamma = 0.5\)}
        \addlegendimage{marker4}
        \addlegendentry{\(\gamma = 0.3\)}
        \addlegendimage{marker5}
        \addlegendentry{\(\gamma = 0.1\)}


        \pgfplotstableread[col sep=comma]{
       k,     ell,    ndof,   nElem,        eta,         mu,        res,    case,   maxGradU
       0,       0,     193,      96,2.01791e-01,1.16608e+00,2.55303e-01,       Z,8.62011e-03
       1,       0,     193,      96,2.01405e-01,1.51358e-01,2.13723e-01,       R,1.60043e-02
       1,       1,     271,     135,1.63298e-01,1.91851e-01,1.79901e-01,       Z,2.95074e-02
       2,       0,     271,     135,1.63037e-01,7.16045e-02,1.67435e-01,       Z,4.14307e-02
       3,       0,     271,     135,1.62901e-01,2.90517e-02,1.65264e-01,       R,4.69924e-02
       3,       1,     427,     213,1.39329e-01,8.92647e-02,1.42920e-01,       Z,7.06245e-02
       4,       0,     427,     213,1.39178e-01,2.65099e-02,1.40792e-01,       R,8.57288e-02
       4,       1,     633,     316,1.12356e-01,8.63095e-02,1.16536e-01,       Z,1.20423e-01
       5,       0,     633,     316,1.12328e-01,2.82196e-02,1.13647e-01,       Z,1.47185e-01
       6,       0,     633,     316,1.12315e-01,1.15980e-02,1.13149e-01,       R,1.60334e-01
       6,       1,     937,     468,9.12346e-02,6.65240e-02,9.38424e-02,       Z,1.62853e-01
       7,       0,     937,     468,9.12420e-02,2.03783e-02,9.19383e-02,       Z,1.70365e-01
       8,       0,     937,     468,9.12472e-02,7.78354e-03,9.16595e-02,       R,1.73925e-01
       8,       1,    1309,     654,7.54564e-02,5.18938e-02,7.74173e-02,       Z,2.33313e-01
       9,       0,    1309,     654,7.54574e-02,1.63962e-02,7.59310e-02,       Z,2.67113e-01
      10,       0,    1309,     654,7.54574e-02,6.39634e-03,7.57034e-02,       R,2.84215e-01
      10,       1,    2137,    1068,6.08165e-02,4.51230e-02,6.24039e-02,       Z,3.76092e-01
      11,       0,    2137,    1068,6.08435e-02,1.33715e-02,6.12043e-02,       Z,4.36950e-01
      12,       0,    2137,    1068,6.08539e-02,5.40803e-03,6.10093e-02,       R,4.68905e-01
      12,       1,    3089,    1544,4.92090e-02,3.62065e-02,5.05179e-02,       Z,6.00485e-01
      13,       0,    3089,    1544,4.92085e-02,1.11692e-02,4.94582e-02,       Z,6.94138e-01
      14,       0,    3089,    1544,4.92082e-02,4.33988e-03,4.92970e-02,       R,7.44808e-01
      14,       1,    4589,    2294,4.01562e-02,2.87710e-02,4.11967e-02,       Z,9.64786e-01
      15,       0,    4589,    2294,4.01500e-02,9.09263e-03,4.03372e-02,       Z,1.12408e+00
      16,       0,    4589,    2294,4.01479e-02,3.57172e-03,4.02034e-02,       R,1.21258e+00
      16,       1,    7187,    3593,3.25030e-02,2.38362e-02,3.33021e-02,       Z,1.56129e+00
      17,       0,    7187,    3593,3.25053e-02,7.16600e-03,3.26503e-02,       Z,1.82386e+00
      18,       0,    7187,    3593,3.25067e-02,2.87589e-03,3.25452e-02,       R,1.97272e+00
      18,       1,   10419,    5209,2.68242e-02,1.85855e-02,2.74027e-02,       Z,2.52575e+00
      19,       0,   10419,    5209,2.68268e-02,5.54342e-03,2.69252e-02,       Z,2.95377e+00
      20,       0,   10419,    5209,2.68280e-02,2.17301e-03,2.68520e-02,       R,3.20017e+00
      20,       1,   15673,    7836,2.18321e-02,1.57425e-02,2.23725e-02,       Z,3.24298e+00
      21,       0,   15673,    7836,2.18306e-02,4.87122e-03,2.19192e-02,       Z,3.40512e+00
      22,       0,   15673,    7836,2.18303e-02,1.91498e-03,2.18486e-02,       R,3.49461e+00
      22,       1,   23249,   11624,1.78925e-02,1.26526e-02,1.83225e-02,       Z,4.42420e+00
      23,       0,   23249,   11624,1.78942e-02,3.92574e-03,1.79647e-02,       R,5.04005e+00
      23,       1,   34811,   17405,1.47452e-02,1.08718e-02,1.51273e-02,       Z,6.43375e+00
      24,       0,   34811,   17405,1.47492e-02,3.35415e-03,1.48135e-02,       Z,7.70067e+00
      25,       0,   34811,   17405,1.47509e-02,1.34287e-03,1.47634e-02,       R,8.44703e+00
      25,       1,   49807,   24903,1.22872e-02,8.27143e-03,1.25374e-02,       Z,1.06775e+01
      26,       0,   49807,   24903,1.22882e-02,2.48284e-03,1.23277e-02,       Z,1.25320e+01
      27,       0,   49807,   24903,1.22886e-02,9.69309e-04,1.22958e-02,       R,1.36208e+01
      27,       1,   71459,   35729,1.01768e-02,6.95582e-03,1.04136e-02,       Z,1.69306e+01
      28,       0,   71459,   35729,1.01764e-02,2.20749e-03,1.02143e-02,       Z,1.96559e+01
      29,       0,   71459,   35729,1.01763e-02,8.74291e-04,1.01829e-02,       R,2.12546e+01
      29,       1,  103867,   51933,8.51327e-03,5.64327e-03,8.68291e-03,       Z,2.69515e+01
      30,       0,  103867,   51933,8.51369e-03,1.70441e-03,8.54125e-03,       Z,3.14958e+01
      31,       0,  103867,   51933,8.51393e-03,6.78939e-04,8.51880e-03,       R,3.41693e+01
      31,       1,  148269,   74134,7.12601e-03,4.70829e-03,7.25922e-03,       Z,3.45545e+01
      32,       0,  148269,   74134,7.12693e-03,1.37864e-03,7.14859e-03,       R,3.63128e+01
      32,       1,  207617,  103808,6.02299e-03,4.05184e-03,6.15480e-03,       Z,4.56736e+01
      33,       0,  207617,  103808,6.02330e-03,1.26496e-03,6.04426e-03,       Z,5.25239e+01
      34,       0,  207617,  103808,6.02347e-03,4.99684e-04,6.02701e-03,       R,5.65371e+01
      34,       1,  290993,  145496,5.05900e-03,3.30732e-03,5.16659e-03,       Z,6.98422e+01
      35,       0,  290993,  145496,5.05879e-03,1.04959e-03,5.07584e-03,       R,8.05348e+01
      35,       1,  407169,  203584,4.29572e-03,2.87049e-03,4.38946e-03,       Z,1.02460e+02
      36,       0,  407169,  203584,4.29585e-03,9.01526e-04,4.31101e-03,       Z,1.23282e+02
      37,       0,  407169,  203584,4.29592e-03,3.59901e-04,4.29850e-03,       R,1.35684e+02
      37,       1,  572869,  286434,3.62967e-03,2.32595e-03,3.69334e-03,       Z,1.37730e+02
      38,       0,  572869,  286434,3.63005e-03,6.80692e-04,3.64041e-03,       Z,1.46105e+02
      39,       0,  572869,  286434,3.63022e-03,2.71917e-04,3.63200e-03,       R,1.50865e+02
      39,       1,  787563,  393781,3.09451e-03,1.91740e-03,3.14701e-03,       Z,1.87968e+02
      40,       0,  787563,  393781,3.09452e-03,5.72286e-04,3.10273e-03,       R,2.13507e+02
      40,       1, 1097597,  548798,2.61659e-03,1.74844e-03,2.67739e-03,       Z,2.65291e+02
        }\tableDataNine

        \pgfplotstableread[col sep=comma]{
       k,     ell,    ndof,   nElem,        eta,         mu,        res,    case,   maxGradU
       0,       0,     193,      96,2.01791e-01,1.16608e+00,2.55303e-01,       Z,8.62011e-03
       1,       0,     193,      96,2.01405e-01,1.51358e-01,2.13723e-01,       R,1.60043e-02
       1,       1,     271,     135,1.63298e-01,1.91851e-01,1.79901e-01,       R,2.95074e-02
       1,       2,     427,     213,1.39639e-01,2.09700e-01,1.59908e-01,       Z,4.90105e-02
       2,       0,     427,     213,1.39332e-01,7.61752e-02,1.43841e-01,       R,7.50225e-02
       2,       1,     633,     316,1.12460e-01,1.12110e-01,1.20289e-01,       R,1.08018e-01
       2,       2,     937,     468,9.13725e-02,1.29873e-01,1.02550e-01,       Z,1.10158e-01
       3,       0,     937,     468,9.13069e-02,4.59695e-02,9.35695e-02,       R,1.44640e-01
       3,       1,    1309,     654,7.55073e-02,6.89114e-02,7.97201e-02,       R,1.99973e-01
       3,       2,    2153,    1076,6.04047e-02,8.24706e-02,6.69198e-02,       Z,2.80392e-01
       4,       0,    2153,    1076,6.04488e-02,2.84836e-02,6.16732e-02,       R,3.86472e-01
       4,       1,    3153,    1576,4.86624e-02,4.57967e-02,5.12906e-02,       R,5.19901e-01
       4,       2,    4689,    2344,3.99149e-02,5.35927e-02,4.39412e-02,       Z,7.12169e-01
       5,       0,    4689,    2344,3.99048e-02,1.83312e-02,4.06100e-02,       R,9.98491e-01
       5,       1,    7325,    3662,3.22427e-02,2.98133e-02,3.38350e-02,       Z,1.30277e+00
       6,       0,    7325,    3662,3.22403e-02,1.02134e-02,3.25208e-02,       R,1.67540e+00
       6,       1,   10569,    5284,2.64970e-02,2.10158e-02,2.73840e-02,       R,2.20376e+00
       6,       2,   16031,    8015,2.16501e-02,2.59812e-02,2.32130e-02,       Z,2.95671e+00
       7,       0,   16031,    8015,2.16474e-02,8.36772e-03,2.19066e-02,       R,4.10615e+00
       7,       1,   23841,   11920,1.77103e-02,1.49992e-02,1.84257e-02,       R,4.16675e+00
       7,       2,   35227,   17613,1.45974e-02,1.80429e-02,1.57432e-02,       Z,5.46391e+00
       8,       0,   35227,   17613,1.46008e-02,5.88381e-03,1.47925e-02,       R,7.13634e+00
       8,       1,   51059,   25529,1.21291e-02,1.00343e-02,1.25928e-02,       R,9.25386e+00
       8,       2,   72969,   36484,1.00920e-02,1.20812e-02,1.08526e-02,       Z,1.19822e+01
       9,       0,   72969,   36484,1.00913e-02,3.99199e-03,1.02164e-02,       R,1.66263e+01
       9,       1,  105873,   52936,8.42995e-03,6.83416e-03,8.74110e-03,       R,2.17044e+01
       9,       2,  150399,   75199,7.06158e-03,8.24038e-03,7.55037e-03,       Z,2.20959e+01
      10,       0,  150399,   75199,7.06261e-03,2.66930e-03,7.14387e-03,       R,2.89801e+01
      10,       1,  212649,  106324,5.95229e-03,4.64497e-03,6.15504e-03,       R,3.73907e+01
      10,       2,  295965,  147982,5.02110e-03,5.63863e-03,5.35470e-03,       Z,4.82329e+01
      11,       0,  295965,  147982,5.02070e-03,1.86138e-03,5.07547e-03,       R,6.70297e+01
      11,       1,  415273,  207636,4.25512e-03,3.25056e-03,4.39339e-03,       R,8.71670e+01
      11,       2,  582005,  291002,3.59627e-03,3.96724e-03,3.81665e-03,       Z,8.87665e+01
      12,       0,  582005,  291002,3.59679e-03,1.27659e-03,3.63336e-03,       R,1.16515e+02
      12,       1,  805989,  402994,3.06038e-03,2.28050e-03,3.15392e-03,       R,1.49615e+02
      12,       2, 1118667,  559333,2.59318e-03,2.80036e-03,2.75084e-03,       Z,1.93440e+02
        }\tableDataSeven

        \pgfplotstableread[col sep=comma]{
       k,     ell,    ndof,   nElem,        eta,         mu,        res,    case,   maxGradU
       0,       0,     193,      96,2.01791e-01,1.16608e+00,2.55303e-01,       Z,8.62011e-03
       1,       0,     193,      96,2.01405e-01,1.51358e-01,2.13723e-01,       R,1.60043e-02
       1,       1,     271,     135,1.63298e-01,1.91851e-01,1.79901e-01,       R,2.95074e-02
       1,       2,     427,     213,1.39639e-01,2.09700e-01,1.59908e-01,       R,4.90105e-02
       1,       3,     641,     320,1.11407e-01,2.25968e-01,1.38294e-01,       R,7.67314e-02
       1,       4,     977,     488,9.07399e-02,2.35031e-01,1.23960e-01,       Z,1.18293e-01
       2,       0,     977,     488,9.06315e-02,8.41899e-02,9.72697e-02,       R,2.01611e-01
       2,       1,    1369,     684,7.45102e-02,9.87433e-02,8.38655e-02,       R,2.06406e-01
       2,       2,    2179,    1089,6.02900e-02,1.08014e-01,7.27452e-02,       R,2.88715e-01
       2,       3,    3169,    1584,4.86591e-02,1.13729e-01,6.44483e-02,       Z,4.09290e-01
       3,       0,    3169,    1584,4.86385e-02,4.22191e-02,5.17022e-02,       R,6.03994e-01
       3,       1,    4745,    2372,3.93809e-02,5.09637e-02,4.40339e-02,       R,8.06870e-01
       3,       2,    7559,    3779,3.18376e-02,5.59868e-02,3.80748e-02,       R,1.11757e+00
       3,       3,   10879,    5439,2.61074e-02,5.88775e-02,3.38894e-02,       R,1.53329e+00
       3,       4,   16375,    8187,2.13214e-02,6.07746e-02,3.07304e-02,       Z,2.20133e+00
       4,       0,   16375,    8187,2.13090e-02,2.21379e-02,2.31777e-02,       R,3.64712e+00
       4,       1,   24255,   12127,1.76206e-02,2.51728e-02,2.01731e-02,       R,4.86656e+00
       4,       2,   35827,   17913,1.45003e-02,2.70905e-02,1.77650e-02,       R,4.95881e+00
       4,       3,   51117,   25558,1.21334e-02,2.82301e-02,1.60738e-02,       Z,6.77818e+00
       5,       0,   51117,   25558,1.21354e-02,1.05393e-02,1.28955e-02,       R,1.02212e+01
       5,       1,   73513,   36756,1.00609e-02,1.25349e-02,1.11727e-02,       R,1.31404e+01
       5,       2,  105679,   52839,8.43582e-03,1.36815e-02,9.87292e-03,       R,1.77532e+01
       5,       3,  150327,   75163,7.06787e-03,1.44358e-02,8.83834e-03,       R,1.80421e+01
       5,       4,  211713,  105856,5.96778e-03,1.49242e-02,8.06792e-03,       Z,2.48980e+01
       6,       0,  211713,  105856,5.96807e-03,5.42884e-03,6.37407e-03,       R,3.95315e+01
       6,       1,  294863,  147431,5.03214e-03,6.30616e-03,5.60011e-03,       R,5.08938e+01
       6,       2,  411915,  205957,4.27453e-03,6.84240e-03,4.99389e-03,       R,6.87772e+01
       6,       3,  578355,  289177,3.60926e-03,7.21550e-03,4.48789e-03,       R,7.01036e+01
       6,       4,  797931,  398965,3.07301e-03,7.45967e-03,4.10851e-03,       Z,9.67506e+01
       7,       0,  797931,  398965,3.07333e-03,2.72661e-03,3.27263e-03,       R,1.56262e+02
       7,       1, 1110073,  555036,2.60296e-03,3.17874e-03,2.88244e-03,       Z,2.01218e+02
        }\tableDataFive

        \pgfplotstableread[col sep=comma]{
       k,     ell,    ndof,   nElem,        eta,         mu,        res,    case,   maxGradU
       0,       0,     193,      96,2.01791e-01,1.16608e+00,2.55303e-01,       Z,8.62011e-03
       1,       0,     193,      96,2.01405e-01,1.51358e-01,2.13723e-01,       R,1.60043e-02
       1,       1,     271,     135,1.63298e-01,1.91851e-01,1.79901e-01,       R,2.95074e-02
       1,       2,     427,     213,1.39639e-01,2.09700e-01,1.59908e-01,       R,4.90105e-02
       1,       3,     641,     320,1.11407e-01,2.25968e-01,1.38294e-01,       R,7.67314e-02
       1,       4,     977,     488,9.07399e-02,2.35031e-01,1.23960e-01,       R,1.18293e-01
       1,       5,    1377,     688,7.45122e-02,2.40668e-01,1.13849e-01,       R,1.76056e-01
       1,       6,    2199,    1099,5.98309e-02,2.44731e-01,1.05794e-01,       Z,2.42559e-01
       2,       0,    2199,    1099,5.97773e-02,8.72164e-02,6.97431e-02,       R,4.41010e-01
       2,       1,    3237,    1618,4.80609e-02,9.41816e-02,6.10484e-02,       R,6.03934e-01
       2,       2,    4893,    2446,3.89686e-02,9.82927e-02,5.49449e-02,       R,8.47235e-01
       2,       3,    7649,    3824,3.16881e-02,1.00876e-01,5.05376e-02,       R,1.21745e+00
       2,       4,   10949,    5474,2.61538e-02,1.02450e-01,4.75975e-02,       R,1.21955e+00
       2,       5,   16395,    8197,2.13444e-02,1.03559e-01,4.54069e-02,       R,1.77374e+00
       2,       6,   24667,   12333,1.75445e-02,1.04270e-01,4.39225e-02,       Z,2.62359e+00
       3,       0,   24667,   12333,1.75210e-02,4.02755e-02,2.42263e-02,       R,4.71555e+00
       3,       1,   36207,   18103,1.44053e-02,4.14920e-02,2.22840e-02,       R,6.33533e+00
       3,       2,   51851,   25925,1.20542e-02,4.22350e-02,2.09830e-02,       R,8.70095e+00
       3,       3,   74547,   37273,9.99424e-03,4.27693e-02,1.99866e-02,       R,8.69752e+00
       3,       4,  107247,   53623,8.39046e-03,4.31126e-02,1.93076e-02,       R,1.21834e+01
       3,       5,  151831,   75915,7.03540e-03,4.33544e-02,1.88096e-02,       R,1.74713e+01
       3,       6,  214041,  107020,5.94238e-03,4.35176e-02,1.84656e-02,       R,2.54542e+01
       3,       7,  297875,  148937,5.01422e-03,4.36343e-02,1.82169e-02,       Z,2.54922e+01
       4,       0,  297875,  148937,5.01116e-03,1.75141e-02,8.90149e-03,       R,4.70471e+01
       4,       1,  415469,  207734,4.26177e-03,1.77114e-02,8.53948e-03,       R,6.25420e+01
       4,       2,  579321,  289660,3.60606e-03,1.78565e-02,8.25937e-03,       R,8.56216e+01
       4,       3,  797191,  398595,3.07297e-03,1.79559e-02,8.06114e-03,       R,1.19452e+02
       4,       4, 1112677,  556338,2.60201e-03,1.80301e-02,7.91071e-03,       Z,1.19417e+02
        }\tableDataThree

        \pgfplotstableread[col sep=comma]{
       k,     ell,    ndof,   nElem,        eta,         mu,        res,    case,   maxGradU
       0,       0,     193,      96,2.01791e-01,1.16608e+00,2.55303e-01,       Z,8.62011e-03
       1,       0,     193,      96,2.01405e-01,1.51358e-01,2.13723e-01,       R,1.60043e-02
       1,       1,     271,     135,1.63298e-01,1.91851e-01,1.79901e-01,       R,2.95074e-02
       1,       2,     427,     213,1.39639e-01,2.09700e-01,1.59908e-01,       R,4.90105e-02
       1,       3,     641,     320,1.11407e-01,2.25968e-01,1.38294e-01,       R,7.67314e-02
       1,       4,     977,     488,9.07399e-02,2.35031e-01,1.23960e-01,       R,1.18293e-01
       1,       5,    1377,     688,7.45122e-02,2.40668e-01,1.13849e-01,       R,1.76056e-01
       1,       6,    2199,    1099,5.98309e-02,2.44731e-01,1.05794e-01,       R,2.42559e-01
       1,       7,    3237,    1618,4.80872e-02,2.47307e-01,1.00353e-01,       R,3.73499e-01
       1,       8,    4901,    2450,3.89676e-02,2.48907e-01,9.68075e-02,       R,5.81609e-01
       1,       9,    7713,    3856,3.16253e-02,2.49946e-01,9.43938e-02,       R,9.11302e-01
       1,      10,   10967,    5483,2.61304e-02,2.50580e-01,9.28934e-02,       R,9.12945e-01
       1,      11,   16361,    8180,2.14089e-02,2.51027e-01,9.18272e-02,       R,1.42526e+00
       1,      12,   24439,   12219,1.76262e-02,2.51321e-01,9.11152e-02,       Z,2.24771e+00
       2,       0,   24439,   12219,1.75952e-02,8.93997e-02,4.07109e-02,       R,4.55854e+00
       2,       1,   35823,   17911,1.45070e-02,8.99524e-02,3.95964e-02,       R,6.18251e+00
       2,       2,   51215,   25607,1.21460e-02,9.03015e-02,3.88761e-02,       R,8.57713e+00
       2,       3,   73801,   36900,1.00617e-02,9.05575e-02,3.83416e-02,       R,8.55828e+00
       2,       4,  106409,   53204,8.42932e-03,9.07240e-02,3.79862e-02,       R,1.21156e+01
       2,       5,  150595,   75297,7.07308e-03,9.08398e-02,3.77356e-02,       R,1.75240e+01
       2,       6,  211607,  105803,5.98027e-03,9.09183e-02,3.75656e-02,       R,2.56967e+01
       2,       7,  295083,  147541,5.04428e-03,9.09750e-02,3.74431e-02,       R,2.57073e+01
       2,       8,  413065,  206532,4.27721e-03,9.10143e-02,3.73572e-02,       R,3.83690e+01
       2,       9,  577189,  288594,3.61787e-03,9.10429e-02,3.72942e-02,       R,5.75245e+01
       2,      10,  795753,  397876,3.08210e-03,9.10626e-02,3.72509e-02,       R,8.75013e+01
       2,      11, 1111411,  555705,2.60834e-03,9.10774e-02,3.72187e-02,       Z,8.75714e+01
        }\tableDataOne

        %
        %

        \addplot+ [marker1, adaptive, forget plot]
        table [col sep=comma, x=cumulativeNdof, y=estimator] {\tableDataNine};

        \addplot+ [marker2, adaptive, forget plot]
        table [col sep=comma, x=cumulativeNdof, y=estimator] {\tableDataSeven};

        \addplot+ [marker3, adaptive, forget plot]
        table [col sep=comma, x=cumulativeNdof, y=estimator] {\tableDataFive};

        \addplot+ [marker4, adaptive, forget plot]
        table [col sep=comma, x=cumulativeNdof, y=estimator] {\tableDataThree};

        \addplot+ [marker5, adaptive, forget plot]
        table [col sep=comma, x=cumulativeNdof, y=estimator] {\tableDataOne};

        \drawslopetriangle[ST1]{0.5}{5e4}{1e-2} 
    \end{loglogaxis}
\end{tikzpicture}

%% file: figures/Fig02b_Nonlinear_convergence_gamma.tex
\begin{tikzpicture}[>=stealth]
    \begin{loglogaxis}[%
            width            = 5.5cm,%
            xlabel           = {cumulative ndof},%
            ymajorgrids      = true,%
            ymax             = 5,%
            ymin             = 2e-3,%
            font             = \footnotesize,%
            grid style       = {%
                densely dotted,%
                semithick%
            },%
            legend style     = {%
                legend pos = north east,%
                font = \footnotesize%
            },%
        ]

        \addlegendimage{marker1}
        \addlegendentry{\(\gamma = 0.9\)}
        \addlegendimage{marker2}
        \addlegendentry{\(\gamma = 0.7\)}
        \addlegendimage{marker3}
        \addlegendentry{\(\gamma = 0.5\)}
        \addlegendimage{marker4}
        \addlegendentry{\(\gamma = 0.3\)}
        \addlegendimage{marker5}
        \addlegendentry{\(\gamma = 0.1\)}


        \pgfplotstableread[col sep=comma]{
       k,     ell,    ndof,   nElem,        eta,         mu,        res,    case,   maxGradU
       0,       0,     193,      96,2.01791e-01,1.16608e+00,2.55303e-01,       Z,8.62011e-03
       1,       0,     193,      96,2.01405e-01,1.51358e-01,2.13723e-01,       R,1.60043e-02
       1,       1,     271,     135,1.63298e-01,1.91851e-01,1.79901e-01,       Z,2.95074e-02
       2,       0,     271,     135,1.63037e-01,7.16045e-02,1.67435e-01,       Z,4.14307e-02
       3,       0,     271,     135,1.62901e-01,2.90517e-02,1.65264e-01,       R,4.69924e-02
       3,       1,     427,     213,1.39329e-01,8.92647e-02,1.42920e-01,       Z,7.06245e-02
       4,       0,     427,     213,1.39178e-01,2.65099e-02,1.40792e-01,       R,8.57288e-02
       4,       1,     633,     316,1.12356e-01,8.63095e-02,1.16536e-01,       Z,1.20423e-01
       5,       0,     633,     316,1.12328e-01,2.82196e-02,1.13647e-01,       Z,1.47185e-01
       6,       0,     633,     316,1.12315e-01,1.15980e-02,1.13149e-01,       R,1.60334e-01
       6,       1,     937,     468,9.12346e-02,6.65240e-02,9.38424e-02,       Z,1.62853e-01
       7,       0,     937,     468,9.12420e-02,2.03783e-02,9.19383e-02,       Z,1.70365e-01
       8,       0,     937,     468,9.12472e-02,7.78354e-03,9.16595e-02,       R,1.73925e-01
       8,       1,    1309,     654,7.54564e-02,5.18938e-02,7.74173e-02,       Z,2.33313e-01
       9,       0,    1309,     654,7.54574e-02,1.63962e-02,7.59310e-02,       Z,2.67113e-01
      10,       0,    1309,     654,7.54574e-02,6.39634e-03,7.57034e-02,       R,2.84215e-01
      10,       1,    2137,    1068,6.08165e-02,4.51230e-02,6.24039e-02,       Z,3.76092e-01
      11,       0,    2137,    1068,6.08435e-02,1.33715e-02,6.12043e-02,       Z,4.36950e-01
      12,       0,    2137,    1068,6.08539e-02,5.40803e-03,6.10093e-02,       R,4.68905e-01
      12,       1,    3089,    1544,4.92090e-02,3.62065e-02,5.05179e-02,       Z,6.00485e-01
      13,       0,    3089,    1544,4.92085e-02,1.11692e-02,4.94582e-02,       Z,6.94138e-01
      14,       0,    3089,    1544,4.92082e-02,4.33988e-03,4.92970e-02,       R,7.44808e-01
      14,       1,    4589,    2294,4.01562e-02,2.87710e-02,4.11967e-02,       Z,9.64786e-01
      15,       0,    4589,    2294,4.01500e-02,9.09263e-03,4.03372e-02,       Z,1.12408e+00
      16,       0,    4589,    2294,4.01479e-02,3.57172e-03,4.02034e-02,       R,1.21258e+00
      16,       1,    7187,    3593,3.25030e-02,2.38362e-02,3.33021e-02,       Z,1.56129e+00
      17,       0,    7187,    3593,3.25053e-02,7.16600e-03,3.26503e-02,       Z,1.82386e+00
      18,       0,    7187,    3593,3.25067e-02,2.87589e-03,3.25452e-02,       R,1.97272e+00
      18,       1,   10419,    5209,2.68242e-02,1.85855e-02,2.74027e-02,       Z,2.52575e+00
      19,       0,   10419,    5209,2.68268e-02,5.54342e-03,2.69252e-02,       Z,2.95377e+00
      20,       0,   10419,    5209,2.68280e-02,2.17301e-03,2.68520e-02,       R,3.20017e+00
      20,       1,   15673,    7836,2.18321e-02,1.57425e-02,2.23725e-02,       Z,3.24298e+00
      21,       0,   15673,    7836,2.18306e-02,4.87122e-03,2.19192e-02,       Z,3.40512e+00
      22,       0,   15673,    7836,2.18303e-02,1.91498e-03,2.18486e-02,       R,3.49461e+00
      22,       1,   23249,   11624,1.78925e-02,1.26526e-02,1.83225e-02,       Z,4.42420e+00
      23,       0,   23249,   11624,1.78942e-02,3.92574e-03,1.79647e-02,       R,5.04005e+00
      23,       1,   34811,   17405,1.47452e-02,1.08718e-02,1.51273e-02,       Z,6.43375e+00
      24,       0,   34811,   17405,1.47492e-02,3.35415e-03,1.48135e-02,       Z,7.70067e+00
      25,       0,   34811,   17405,1.47509e-02,1.34287e-03,1.47634e-02,       R,8.44703e+00
      25,       1,   49807,   24903,1.22872e-02,8.27143e-03,1.25374e-02,       Z,1.06775e+01
      26,       0,   49807,   24903,1.22882e-02,2.48284e-03,1.23277e-02,       Z,1.25320e+01
      27,       0,   49807,   24903,1.22886e-02,9.69309e-04,1.22958e-02,       R,1.36208e+01
      27,       1,   71459,   35729,1.01768e-02,6.95582e-03,1.04136e-02,       Z,1.69306e+01
      28,       0,   71459,   35729,1.01764e-02,2.20749e-03,1.02143e-02,       Z,1.96559e+01
      29,       0,   71459,   35729,1.01763e-02,8.74291e-04,1.01829e-02,       R,2.12546e+01
      29,       1,  103867,   51933,8.51327e-03,5.64327e-03,8.68291e-03,       Z,2.69515e+01
      30,       0,  103867,   51933,8.51369e-03,1.70441e-03,8.54125e-03,       Z,3.14958e+01
      31,       0,  103867,   51933,8.51393e-03,6.78939e-04,8.51880e-03,       R,3.41693e+01
      31,       1,  148269,   74134,7.12601e-03,4.70829e-03,7.25922e-03,       Z,3.45545e+01
      32,       0,  148269,   74134,7.12693e-03,1.37864e-03,7.14859e-03,       R,3.63128e+01
      32,       1,  207617,  103808,6.02299e-03,4.05184e-03,6.15480e-03,       Z,4.56736e+01
      33,       0,  207617,  103808,6.02330e-03,1.26496e-03,6.04426e-03,       Z,5.25239e+01
      34,       0,  207617,  103808,6.02347e-03,4.99684e-04,6.02701e-03,       R,5.65371e+01
      34,       1,  290993,  145496,5.05900e-03,3.30732e-03,5.16659e-03,       Z,6.98422e+01
      35,       0,  290993,  145496,5.05879e-03,1.04959e-03,5.07584e-03,       R,8.05348e+01
      35,       1,  407169,  203584,4.29572e-03,2.87049e-03,4.38946e-03,       Z,1.02460e+02
      36,       0,  407169,  203584,4.29585e-03,9.01526e-04,4.31101e-03,       Z,1.23282e+02
      37,       0,  407169,  203584,4.29592e-03,3.59901e-04,4.29850e-03,       R,1.35684e+02
      37,       1,  572869,  286434,3.62967e-03,2.32595e-03,3.69334e-03,       Z,1.37730e+02
      38,       0,  572869,  286434,3.63005e-03,6.80692e-04,3.64041e-03,       Z,1.46105e+02
      39,       0,  572869,  286434,3.63022e-03,2.71917e-04,3.63200e-03,       R,1.50865e+02
      39,       1,  787563,  393781,3.09451e-03,1.91740e-03,3.14701e-03,       Z,1.87968e+02
      40,       0,  787563,  393781,3.09452e-03,5.72286e-04,3.10273e-03,       R,2.13507e+02
      40,       1, 1097597,  548798,2.61659e-03,1.74844e-03,2.67739e-03,       Z,2.65291e+02
        }\tableDataNine

        \pgfplotstableread[col sep=comma]{
       k,     ell,    ndof,   nElem,        eta,         mu,        res,    case,   maxGradU
       0,       0,     193,      96,2.01791e-01,1.16608e+00,2.55303e-01,       Z,8.62011e-03
       1,       0,     193,      96,2.01405e-01,1.51358e-01,2.13723e-01,       R,1.60043e-02
       1,       1,     271,     135,1.63298e-01,1.91851e-01,1.79901e-01,       R,2.95074e-02
       1,       2,     427,     213,1.39639e-01,2.09700e-01,1.59908e-01,       Z,4.90105e-02
       2,       0,     427,     213,1.39332e-01,7.61752e-02,1.43841e-01,       R,7.50225e-02
       2,       1,     633,     316,1.12460e-01,1.12110e-01,1.20289e-01,       R,1.08018e-01
       2,       2,     937,     468,9.13725e-02,1.29873e-01,1.02550e-01,       Z,1.10158e-01
       3,       0,     937,     468,9.13069e-02,4.59695e-02,9.35695e-02,       R,1.44640e-01
       3,       1,    1309,     654,7.55073e-02,6.89114e-02,7.97201e-02,       R,1.99973e-01
       3,       2,    2153,    1076,6.04047e-02,8.24706e-02,6.69198e-02,       Z,2.80392e-01
       4,       0,    2153,    1076,6.04488e-02,2.84836e-02,6.16732e-02,       R,3.86472e-01
       4,       1,    3153,    1576,4.86624e-02,4.57967e-02,5.12906e-02,       R,5.19901e-01
       4,       2,    4689,    2344,3.99149e-02,5.35927e-02,4.39412e-02,       Z,7.12169e-01
       5,       0,    4689,    2344,3.99048e-02,1.83312e-02,4.06100e-02,       R,9.98491e-01
       5,       1,    7325,    3662,3.22427e-02,2.98133e-02,3.38350e-02,       Z,1.30277e+00
       6,       0,    7325,    3662,3.22403e-02,1.02134e-02,3.25208e-02,       R,1.67540e+00
       6,       1,   10569,    5284,2.64970e-02,2.10158e-02,2.73840e-02,       R,2.20376e+00
       6,       2,   16031,    8015,2.16501e-02,2.59812e-02,2.32130e-02,       Z,2.95671e+00
       7,       0,   16031,    8015,2.16474e-02,8.36772e-03,2.19066e-02,       R,4.10615e+00
       7,       1,   23841,   11920,1.77103e-02,1.49992e-02,1.84257e-02,       R,4.16675e+00
       7,       2,   35227,   17613,1.45974e-02,1.80429e-02,1.57432e-02,       Z,5.46391e+00
       8,       0,   35227,   17613,1.46008e-02,5.88381e-03,1.47925e-02,       R,7.13634e+00
       8,       1,   51059,   25529,1.21291e-02,1.00343e-02,1.25928e-02,       R,9.25386e+00
       8,       2,   72969,   36484,1.00920e-02,1.20812e-02,1.08526e-02,       Z,1.19822e+01
       9,       0,   72969,   36484,1.00913e-02,3.99199e-03,1.02164e-02,       R,1.66263e+01
       9,       1,  105873,   52936,8.42995e-03,6.83416e-03,8.74110e-03,       R,2.17044e+01
       9,       2,  150399,   75199,7.06158e-03,8.24038e-03,7.55037e-03,       Z,2.20959e+01
      10,       0,  150399,   75199,7.06261e-03,2.66930e-03,7.14387e-03,       R,2.89801e+01
      10,       1,  212649,  106324,5.95229e-03,4.64497e-03,6.15504e-03,       R,3.73907e+01
      10,       2,  295965,  147982,5.02110e-03,5.63863e-03,5.35470e-03,       Z,4.82329e+01
      11,       0,  295965,  147982,5.02070e-03,1.86138e-03,5.07547e-03,       R,6.70297e+01
      11,       1,  415273,  207636,4.25512e-03,3.25056e-03,4.39339e-03,       R,8.71670e+01
      11,       2,  582005,  291002,3.59627e-03,3.96724e-03,3.81665e-03,       Z,8.87665e+01
      12,       0,  582005,  291002,3.59679e-03,1.27659e-03,3.63336e-03,       R,1.16515e+02
      12,       1,  805989,  402994,3.06038e-03,2.28050e-03,3.15392e-03,       R,1.49615e+02
      12,       2, 1118667,  559333,2.59318e-03,2.80036e-03,2.75084e-03,       Z,1.93440e+02
        }\tableDataSeven

        \pgfplotstableread[col sep=comma]{
       k,     ell,    ndof,   nElem,        eta,         mu,        res,    case,   maxGradU
       0,       0,     193,      96,2.01791e-01,1.16608e+00,2.55303e-01,       Z,8.62011e-03
       1,       0,     193,      96,2.01405e-01,1.51358e-01,2.13723e-01,       R,1.60043e-02
       1,       1,     271,     135,1.63298e-01,1.91851e-01,1.79901e-01,       R,2.95074e-02
       1,       2,     427,     213,1.39639e-01,2.09700e-01,1.59908e-01,       R,4.90105e-02
       1,       3,     641,     320,1.11407e-01,2.25968e-01,1.38294e-01,       R,7.67314e-02
       1,       4,     977,     488,9.07399e-02,2.35031e-01,1.23960e-01,       Z,1.18293e-01
       2,       0,     977,     488,9.06315e-02,8.41899e-02,9.72697e-02,       R,2.01611e-01
       2,       1,    1369,     684,7.45102e-02,9.87433e-02,8.38655e-02,       R,2.06406e-01
       2,       2,    2179,    1089,6.02900e-02,1.08014e-01,7.27452e-02,       R,2.88715e-01
       2,       3,    3169,    1584,4.86591e-02,1.13729e-01,6.44483e-02,       Z,4.09290e-01
       3,       0,    3169,    1584,4.86385e-02,4.22191e-02,5.17022e-02,       R,6.03994e-01
       3,       1,    4745,    2372,3.93809e-02,5.09637e-02,4.40339e-02,       R,8.06870e-01
       3,       2,    7559,    3779,3.18376e-02,5.59868e-02,3.80748e-02,       R,1.11757e+00
       3,       3,   10879,    5439,2.61074e-02,5.88775e-02,3.38894e-02,       R,1.53329e+00
       3,       4,   16375,    8187,2.13214e-02,6.07746e-02,3.07304e-02,       Z,2.20133e+00
       4,       0,   16375,    8187,2.13090e-02,2.21379e-02,2.31777e-02,       R,3.64712e+00
       4,       1,   24255,   12127,1.76206e-02,2.51728e-02,2.01731e-02,       R,4.86656e+00
       4,       2,   35827,   17913,1.45003e-02,2.70905e-02,1.77650e-02,       R,4.95881e+00
       4,       3,   51117,   25558,1.21334e-02,2.82301e-02,1.60738e-02,       Z,6.77818e+00
       5,       0,   51117,   25558,1.21354e-02,1.05393e-02,1.28955e-02,       R,1.02212e+01
       5,       1,   73513,   36756,1.00609e-02,1.25349e-02,1.11727e-02,       R,1.31404e+01
       5,       2,  105679,   52839,8.43582e-03,1.36815e-02,9.87292e-03,       R,1.77532e+01
       5,       3,  150327,   75163,7.06787e-03,1.44358e-02,8.83834e-03,       R,1.80421e+01
       5,       4,  211713,  105856,5.96778e-03,1.49242e-02,8.06792e-03,       Z,2.48980e+01
       6,       0,  211713,  105856,5.96807e-03,5.42884e-03,6.37407e-03,       R,3.95315e+01
       6,       1,  294863,  147431,5.03214e-03,6.30616e-03,5.60011e-03,       R,5.08938e+01
       6,       2,  411915,  205957,4.27453e-03,6.84240e-03,4.99389e-03,       R,6.87772e+01
       6,       3,  578355,  289177,3.60926e-03,7.21550e-03,4.48789e-03,       R,7.01036e+01
       6,       4,  797931,  398965,3.07301e-03,7.45967e-03,4.10851e-03,       Z,9.67506e+01
       7,       0,  797931,  398965,3.07333e-03,2.72661e-03,3.27263e-03,       R,1.56262e+02
       7,       1, 1110073,  555036,2.60296e-03,3.17874e-03,2.88244e-03,       Z,2.01218e+02
        }\tableDataFive

        \pgfplotstableread[col sep=comma]{
       k,     ell,    ndof,   nElem,        eta,         mu,        res,    case,   maxGradU
       0,       0,     193,      96,2.01791e-01,1.16608e+00,2.55303e-01,       Z,8.62011e-03
       1,       0,     193,      96,2.01405e-01,1.51358e-01,2.13723e-01,       R,1.60043e-02
       1,       1,     271,     135,1.63298e-01,1.91851e-01,1.79901e-01,       R,2.95074e-02
       1,       2,     427,     213,1.39639e-01,2.09700e-01,1.59908e-01,       R,4.90105e-02
       1,       3,     641,     320,1.11407e-01,2.25968e-01,1.38294e-01,       R,7.67314e-02
       1,       4,     977,     488,9.07399e-02,2.35031e-01,1.23960e-01,       R,1.18293e-01
       1,       5,    1377,     688,7.45122e-02,2.40668e-01,1.13849e-01,       R,1.76056e-01
       1,       6,    2199,    1099,5.98309e-02,2.44731e-01,1.05794e-01,       Z,2.42559e-01
       2,       0,    2199,    1099,5.97773e-02,8.72164e-02,6.97431e-02,       R,4.41010e-01
       2,       1,    3237,    1618,4.80609e-02,9.41816e-02,6.10484e-02,       R,6.03934e-01
       2,       2,    4893,    2446,3.89686e-02,9.82927e-02,5.49449e-02,       R,8.47235e-01
       2,       3,    7649,    3824,3.16881e-02,1.00876e-01,5.05376e-02,       R,1.21745e+00
       2,       4,   10949,    5474,2.61538e-02,1.02450e-01,4.75975e-02,       R,1.21955e+00
       2,       5,   16395,    8197,2.13444e-02,1.03559e-01,4.54069e-02,       R,1.77374e+00
       2,       6,   24667,   12333,1.75445e-02,1.04270e-01,4.39225e-02,       Z,2.62359e+00
       3,       0,   24667,   12333,1.75210e-02,4.02755e-02,2.42263e-02,       R,4.71555e+00
       3,       1,   36207,   18103,1.44053e-02,4.14920e-02,2.22840e-02,       R,6.33533e+00
       3,       2,   51851,   25925,1.20542e-02,4.22350e-02,2.09830e-02,       R,8.70095e+00
       3,       3,   74547,   37273,9.99424e-03,4.27693e-02,1.99866e-02,       R,8.69752e+00
       3,       4,  107247,   53623,8.39046e-03,4.31126e-02,1.93076e-02,       R,1.21834e+01
       3,       5,  151831,   75915,7.03540e-03,4.33544e-02,1.88096e-02,       R,1.74713e+01
       3,       6,  214041,  107020,5.94238e-03,4.35176e-02,1.84656e-02,       R,2.54542e+01
       3,       7,  297875,  148937,5.01422e-03,4.36343e-02,1.82169e-02,       Z,2.54922e+01
       4,       0,  297875,  148937,5.01116e-03,1.75141e-02,8.90149e-03,       R,4.70471e+01
       4,       1,  415469,  207734,4.26177e-03,1.77114e-02,8.53948e-03,       R,6.25420e+01
       4,       2,  579321,  289660,3.60606e-03,1.78565e-02,8.25937e-03,       R,8.56216e+01
       4,       3,  797191,  398595,3.07297e-03,1.79559e-02,8.06114e-03,       R,1.19452e+02
       4,       4, 1112677,  556338,2.60201e-03,1.80301e-02,7.91071e-03,       Z,1.19417e+02
        }\tableDataThree

        \pgfplotstableread[col sep=comma]{
       k,     ell,    ndof,   nElem,        eta,         mu,        res,    case,   maxGradU
       0,       0,     193,      96,2.01791e-01,1.16608e+00,2.55303e-01,       Z,8.62011e-03
       1,       0,     193,      96,2.01405e-01,1.51358e-01,2.13723e-01,       R,1.60043e-02
       1,       1,     271,     135,1.63298e-01,1.91851e-01,1.79901e-01,       R,2.95074e-02
       1,       2,     427,     213,1.39639e-01,2.09700e-01,1.59908e-01,       R,4.90105e-02
       1,       3,     641,     320,1.11407e-01,2.25968e-01,1.38294e-01,       R,7.67314e-02
       1,       4,     977,     488,9.07399e-02,2.35031e-01,1.23960e-01,       R,1.18293e-01
       1,       5,    1377,     688,7.45122e-02,2.40668e-01,1.13849e-01,       R,1.76056e-01
       1,       6,    2199,    1099,5.98309e-02,2.44731e-01,1.05794e-01,       R,2.42559e-01
       1,       7,    3237,    1618,4.80872e-02,2.47307e-01,1.00353e-01,       R,3.73499e-01
       1,       8,    4901,    2450,3.89676e-02,2.48907e-01,9.68075e-02,       R,5.81609e-01
       1,       9,    7713,    3856,3.16253e-02,2.49946e-01,9.43938e-02,       R,9.11302e-01
       1,      10,   10967,    5483,2.61304e-02,2.50580e-01,9.28934e-02,       R,9.12945e-01
       1,      11,   16361,    8180,2.14089e-02,2.51027e-01,9.18272e-02,       R,1.42526e+00
       1,      12,   24439,   12219,1.76262e-02,2.51321e-01,9.11152e-02,       Z,2.24771e+00
       2,       0,   24439,   12219,1.75952e-02,8.93997e-02,4.07109e-02,       R,4.55854e+00
       2,       1,   35823,   17911,1.45070e-02,8.99524e-02,3.95964e-02,       R,6.18251e+00
       2,       2,   51215,   25607,1.21460e-02,9.03015e-02,3.88761e-02,       R,8.57713e+00
       2,       3,   73801,   36900,1.00617e-02,9.05575e-02,3.83416e-02,       R,8.55828e+00
       2,       4,  106409,   53204,8.42932e-03,9.07240e-02,3.79862e-02,       R,1.21156e+01
       2,       5,  150595,   75297,7.07308e-03,9.08398e-02,3.77356e-02,       R,1.75240e+01
       2,       6,  211607,  105803,5.98027e-03,9.09183e-02,3.75656e-02,       R,2.56967e+01
       2,       7,  295083,  147541,5.04428e-03,9.09750e-02,3.74431e-02,       R,2.57073e+01
       2,       8,  413065,  206532,4.27721e-03,9.10143e-02,3.73572e-02,       R,3.83690e+01
       2,       9,  577189,  288594,3.61787e-03,9.10429e-02,3.72942e-02,       R,5.75245e+01
       2,      10,  795753,  397876,3.08210e-03,9.10626e-02,3.72509e-02,       R,8.75013e+01
       2,      11, 1111411,  555705,2.60834e-03,9.10774e-02,3.72187e-02,       Z,8.75714e+01
        }\tableDataOne

        %
        %

        \addplot+ [marker1, adaptive, forget plot]
        table [col sep=comma, x=cumulativeNdof, y=res] {\tableDataNine};

        \addplot+ [marker2, adaptive, forget plot]
        table [col sep=comma, x=cumulativeNdof, y=res] {\tableDataSeven};

        \addplot+ [marker3, adaptive, forget plot]
        table [col sep=comma, x=cumulativeNdof, y=res] {\tableDataFive};

        \addplot+ [marker4, adaptive, forget plot]
        table [col sep=comma, x=cumulativeNdof, y=res] {\tableDataThree};

        \addplot+ [marker5, adaptive, forget plot]
        table [col sep=comma, x=cumulativeNdof, y=res] {\tableDataOne};

        \drawslopetriangle[ST1]{0.5}{2e4}{1e-2} 
    \end{loglogaxis}
\end{tikzpicture}

%% file: figures/Fig02c_Nonlinear_convergence_gamma.tex
\begin{tikzpicture}[>=stealth]
    \begin{loglogaxis}[%
            width            = 5.5cm,%
            xlabel           = {cumulative ndof},%
            ylabel           = {error estimator},%
            ymajorgrids      = true,%
            ymin             = 2e-3,%
            ymax             = 20,%
            font             = \footnotesize,%
            grid style       = {%
                densely dotted,%
                semithick%
            },%
            legend style     = {%
                legend pos = north east,%
                font = \footnotesize%
            },%
        ]

        \addlegendimage{marker1}
        \addlegendentry{\(\gamma = 0.9\)}
        \addlegendimage{marker2}
        \addlegendentry{\(\gamma = 0.7\)}
        \addlegendimage{marker3}
        \addlegendentry{\(\gamma = 0.5\)}
        \addlegendimage{marker4}
        \addlegendentry{\(\gamma = 0.3\)}
        \addlegendimage{marker5}
        \addlegendentry{\(\gamma = 0.1\)}


        \pgfplotstableread[col sep=comma]{
       k,     ell,    ndof,   nElem,        eta,         mu,        res,    case,   maxGradU
       0,       0,     193,      96,1.00896e-01,5.83038e-01,8.99738e-01,       Z,2.15503e-03
       1,       0,     193,      96,1.00842e-01,2.93805e-01,5.01852e-01,       R,6.04181e-03
       1,       1,     271,     135,8.17297e-02,2.99684e-01,4.88116e-01,       Z,1.00747e-02
       2,       0,     271,     135,8.16617e-02,1.56242e-01,2.98824e-01,       Z,1.89306e-02
       3,       0,     271,     135,8.15955e-02,8.64786e-02,2.18991e-01,       R,2.67471e-02
       3,       1,     427,     213,6.97909e-02,9.62580e-02,2.06544e-01,       Z,3.54307e-02
       4,       0,     427,     213,6.97253e-02,5.77212e-02,1.68723e-01,       R,5.01653e-02
       4,       1,     633,     316,5.62695e-02,7.09020e-02,1.54874e-01,       Z,6.32253e-02
       5,       0,     633,     316,5.62384e-02,4.53859e-02,1.32507e-01,       Z,8.79335e-02
       6,       0,     633,     316,5.62149e-02,3.06981e-02,1.22503e-01,       R,1.08203e-01
       6,       1,     937,     468,4.56701e-02,4.49080e-02,1.10924e-01,       Z,1.09260e-01
       7,       0,     937,     468,4.56572e-02,2.97143e-02,1.00952e-01,       Z,1.26100e-01
       8,       0,     937,     468,4.56485e-02,2.03060e-02,9.62142e-02,       R,1.39063e-01
       8,       1,    1309,     654,3.77501e-02,3.27270e-02,8.74973e-02,       Z,1.65001e-01
       9,       0,    1309,     654,3.77440e-02,2.15617e-02,8.13071e-02,       Z,1.96791e-01
      10,       0,    1309,     654,3.77398e-02,1.46174e-02,7.83987e-02,       R,2.21918e-01
      10,       1,    2141,    1070,3.03214e-02,2.68063e-02,6.99016e-02,       Z,2.61717e-01
      11,       0,    2141,    1070,3.03293e-02,1.71424e-02,6.49729e-02,       Z,3.15697e-01
      12,       0,    2141,    1070,3.03347e-02,1.13973e-02,6.27711e-02,       R,3.59023e-01
      12,       1,    3101,    1550,2.45310e-02,2.11736e-02,5.62089e-02,       Z,4.20140e-01
      13,       0,    3101,    1550,2.45329e-02,1.36017e-02,5.23790e-02,       Z,5.07914e-01
      14,       0,    3101,    1550,2.45344e-02,9.04632e-03,5.06600e-02,       R,5.79139e-01
      14,       1,    4611,    2305,2.00855e-02,1.67434e-02,4.57119e-02,       Z,6.66814e-01
      15,       0,    4611,    2305,2.00832e-02,1.08548e-02,4.27472e-02,       Z,8.01035e-01
      16,       0,    4611,    2305,2.00817e-02,7.25705e-03,4.13983e-02,       R,9.10722e-01
      16,       1,    7209,    3604,1.62690e-02,1.38296e-02,3.70372e-02,       Z,1.06319e+00
      17,       0,    7209,    3604,1.62696e-02,8.81547e-03,3.45919e-02,       Z,1.28711e+00
      18,       0,    7209,    3604,1.62701e-02,5.83506e-03,3.35101e-02,       R,1.47191e+00
      18,       1,   10405,    5202,1.34099e-02,1.09060e-02,3.02270e-02,       Z,1.70595e+00
      19,       0,   10405,    5202,1.34113e-02,6.95138e-03,2.83665e-02,       Z,2.06611e+00
      20,       0,   10405,    5202,1.34124e-02,4.59347e-03,2.75480e-02,       R,2.36538e+00
      20,       1,   15749,    7874,1.09239e-02,9.03663e-03,2.47425e-02,       Z,2.38567e+00
      21,       0,   15749,    7874,1.09231e-02,5.79918e-03,2.31691e-02,       Z,2.64613e+00
      22,       0,   15749,    7874,1.09225e-02,3.85009e-03,2.24666e-02,       R,2.85506e+00
      22,       1,   23243,   11621,8.96472e-03,7.33204e-03,2.02599e-02,       Z,3.26478e+00
      23,       0,   23243,   11621,8.96461e-03,4.71221e-03,1.89964e-02,       R,3.78608e+00
      23,       1,   34815,   17407,7.36648e-03,6.95012e-03,1.72182e-02,       Z,4.37836e+00
      24,       0,   34815,   17407,7.36743e-03,4.45108e-03,1.58784e-02,       Z,5.35046e+00
      25,       0,   34815,   17407,7.36813e-03,2.95363e-03,1.52754e-02,       R,6.16870e+00
      25,       1,   49897,   24948,6.13029e-03,5.04309e-03,1.38672e-02,       Z,7.11315e+00
      26,       0,   49897,   24948,6.13083e-03,3.23682e-03,1.29947e-02,       Z,8.65007e+00
      27,       0,   49897,   24948,6.13123e-03,2.14839e-03,1.26058e-02,       R,9.94367e+00
      27,       1,   71927,   35963,5.07506e-03,4.05603e-03,1.14271e-02,       Z,1.13337e+01
      28,       0,   71927,   35963,5.07490e-03,2.62395e-03,1.07382e-02,       Z,1.36648e+01
      29,       0,   71927,   35963,5.07481e-03,1.75212e-03,1.04271e-02,       R,1.56234e+01
      29,       1,  104369,   52184,4.24766e-03,3.28344e-03,9.47892e-03,       Z,1.79842e+01
      30,       0,  104369,   52184,4.24777e-03,2.10161e-03,8.94193e-03,       Z,2.17464e+01
      31,       0,  104369,   52184,4.24785e-03,1.39413e-03,8.70459e-03,       R,2.49139e+01
      31,       1,  148793,   74396,3.55369e-03,2.71278e-03,7.89503e-03,       Z,2.51079e+01
      32,       0,  148793,   74396,3.55397e-03,1.71787e-03,7.46003e-03,       R,2.78847e+01
      32,       1,  208883,  104441,2.99862e-03,2.56711e-03,6.83621e-03,       Z,3.18377e+01
      33,       0,  208883,  104441,2.99873e-03,1.64025e-03,6.37899e-03,       Z,3.76191e+01
      34,       0,  208883,  104441,2.99882e-03,1.08604e-03,6.17594e-03,       R,4.24448e+01
      34,       1,  292283,  146141,2.52507e-03,1.94844e-03,5.64490e-03,       Z,4.81223e+01
      35,       0,  292283,  146141,2.52500e-03,1.26107e-03,5.32361e-03,       Z,5.71812e+01
      36,       0,  292283,  146141,2.52495e-03,8.42382e-04,5.17882e-03,       R,6.47636e+01
      36,       1,  410411,  205205,2.14108e-03,1.58138e-03,4.73709e-03,       Z,7.43316e+01
      37,       0,  410411,  205205,2.14110e-03,1.01269e-03,4.48825e-03,       R,8.92206e+01
      37,       1,  577419,  288709,1.80705e-03,1.53115e-03,4.10657e-03,       Z,8.99622e+01
      38,       0,  577419,  288709,1.80718e-03,9.74680e-04,3.83738e-03,       Z,1.03147e+02
      39,       0,  577419,  288709,1.80727e-03,6.44270e-04,3.71868e-03,       R,1.14021e+02
      39,       1,  793523,  396761,1.54019e-03,1.14417e-03,3.41008e-03,       Z,1.29910e+02
      40,       0,  793523,  396761,1.54024e-03,7.31283e-04,3.22924e-03,       R,1.52856e+02
      40,       1, 1106953,  553476,1.30207e-03,1.10079e-03,2.96734e-03,       Z,1.74297e+02
        }\tableDataNine

        \pgfplotstableread[col sep=comma]{
       k,     ell,    ndof,   nElem,        eta,         mu,        res,    case,   maxGradU
       0,       0,     193,      96,1.00896e-01,5.83038e-01,8.99738e-01,       Z,2.15503e-03
       1,       0,     193,      96,1.00842e-01,2.93805e-01,5.01852e-01,       R,6.04181e-03
       1,       1,     271,     135,8.17297e-02,2.99684e-01,4.88116e-01,       R,1.00747e-02
       1,       2,     427,     213,6.99069e-02,3.02660e-01,4.80820e-01,       Z,1.56773e-02
       2,       0,     427,     213,6.98299e-02,1.58073e-01,2.85109e-01,       R,3.23139e-02
       2,       1,     641,     320,5.57020e-02,1.63586e-01,2.75120e-01,       R,4.30450e-02
       2,       2,     977,     488,4.53764e-02,1.66746e-01,2.69120e-01,       Z,5.79765e-02
       3,       0,     977,     488,4.53452e-02,9.35437e-02,1.70930e-01,       R,1.03840e-01
       3,       1,    1369,     684,3.72854e-02,9.70382e-02,1.65432e-01,       R,1.05584e-01
       3,       2,    2175,    1087,3.02469e-02,9.94572e-02,1.61321e-01,       Z,1.34751e-01
       4,       0,    2175,    1087,3.02307e-02,6.07891e-02,1.10045e-01,       R,2.09471e-01
       4,       1,    3167,    1583,2.42966e-02,6.33946e-02,1.06195e-01,       R,2.58962e-01
       4,       2,    4785,    2392,1.96794e-02,6.49763e-02,1.03720e-01,       Z,3.27499e-01
       5,       0,    4785,    2392,1.96687e-02,4.25774e-02,7.41567e-02,       R,5.21308e-01
       5,       1,    7635,    3817,1.58791e-02,4.41310e-02,7.18007e-02,       Z,6.42873e-01
       6,       0,    7635,    3817,1.58740e-02,3.01536e-02,5.42181e-02,       R,9.44379e-01
       6,       1,   10929,    5464,1.30444e-02,3.14812e-02,5.23073e-02,       R,1.12779e+00
       6,       2,   16349,    8174,1.06930e-02,3.23556e-02,5.09874e-02,       Z,1.40312e+00
       7,       0,   16349,    8174,1.06898e-02,2.24178e-02,3.86262e-02,       R,2.16632e+00
       7,       1,   24281,   12140,8.80195e-03,2.32240e-02,3.74432e-02,       R,2.63082e+00
       7,       2,   35881,   17940,7.22839e-03,2.37609e-02,3.65944e-02,       Z,2.66497e+00
       8,       0,   35881,   17940,7.22785e-03,1.65527e-02,2.76365e-02,       R,3.88803e+00
       8,       1,   51773,   25886,6.03988e-03,1.70222e-02,2.69259e-02,       R,4.66707e+00
       8,       2,   74033,   37016,5.02121e-03,1.73500e-02,2.64190e-02,       Z,5.66099e+00
       9,       0,   74033,   37016,5.02029e-03,1.21277e-02,1.99072e-02,       R,8.70114e+00
       9,       1,  105997,   52998,4.21124e-03,1.24318e-02,1.94411e-02,       R,1.05330e+01
       9,       2,  150561,   75280,3.53042e-03,1.26420e-02,1.91014e-02,       Z,1.06623e+01
      10,       0,  150561,   75280,3.53037e-03,8.84217e-03,1.43573e-02,       R,1.55651e+01
      10,       1,  212419,  106209,2.98265e-03,9.04164e-03,1.40496e-02,       R,1.86211e+01
      10,       2,  295901,  147950,2.51498e-03,9.18272e-03,1.38288e-02,       Z,2.25018e+01
      11,       0,  295901,  147950,2.51462e-03,6.42955e-03,1.03807e-02,       R,3.46095e+01
      11,       1,  413067,  206533,2.13394e-03,6.56573e-03,1.01691e-02,       R,4.18007e+01
      11,       2,  578621,  289310,1.80331e-03,6.66412e-03,1.00085e-02,       Z,4.23118e+01
      12,       0,  578621,  289310,1.80334e-03,4.66396e-03,7.50456e-03,       R,6.18576e+01
      12,       1,  799557,  399778,1.53616e-03,4.75865e-03,7.35694e-03,       R,7.39057e+01
      12,       2, 1113237,  556618,1.30095e-03,4.82826e-03,7.24763e-03,       Z,8.92067e+01
        }\tableDataSeven

        \pgfplotstableread[col sep=comma]{
       k,     ell,    ndof,   nElem,        eta,         mu,        res,    case,   maxGradU
       0,       0,     193,      96,1.00896e-01,5.83038e-01,8.99738e-01,       Z,2.15503e-03
       1,       0,     193,      96,1.00842e-01,2.93805e-01,5.01852e-01,       R,6.04181e-03
       1,       1,     271,     135,8.17297e-02,2.99684e-01,4.88116e-01,       R,1.00747e-02
       1,       2,     427,     213,6.99069e-02,3.02660e-01,4.80820e-01,       R,1.56773e-02
       1,       3,     641,     320,5.57489e-02,3.05585e-01,4.73651e-01,       R,2.34373e-02
       1,       4,     977,     488,4.54090e-02,3.07292e-01,4.69418e-01,       Z,3.48510e-02
       2,       0,     977,     488,4.53788e-02,1.61072e-01,2.62712e-01,       R,8.07874e-02
       2,       1,    1377,     688,3.72684e-02,1.63139e-01,2.58708e-01,       R,1.03494e-01
       2,       2,    2193,    1096,2.99938e-02,1.64632e-01,2.55696e-01,       R,1.27786e-01
       2,       3,    3211,    1605,2.41396e-02,1.65592e-01,2.53770e-01,       Z,1.72957e-01
       3,       0,    3211,    1605,2.41245e-02,9.30548e-02,1.49481e-01,       R,3.36926e-01
       3,       1,    4873,    2436,1.95290e-02,9.41266e-02,1.47635e-01,       R,4.26751e-01
       3,       2,    7715,    3857,1.58063e-02,9.48228e-02,1.46383e-01,       R,5.56299e-01
       3,       3,   10945,    5472,1.30380e-02,9.52429e-02,1.45624e-01,       R,5.56938e-01
       3,       4,   16281,    8140,1.07043e-02,9.55333e-02,1.45101e-01,       Z,7.41440e-01
       4,       0,   16281,    8140,1.06964e-02,5.85311e-02,9.04627e-02,       R,1.40914e+00
       4,       1,   24481,   12240,8.79100e-03,5.88475e-02,8.99279e-02,       R,1.76715e+00
       4,       2,   35987,   17993,7.24273e-03,5.90580e-02,8.95642e-02,       R,2.27521e+00
       4,       3,   51237,   25618,6.06929e-03,5.91901e-02,8.93368e-02,       Z,3.00384e+00
       5,       0,   51237,   25618,6.06727e-03,3.90671e-02,5.92547e-02,       R,5.68832e+00
       5,       1,   73687,   36843,5.02982e-03,3.92142e-02,5.90141e-02,       R,5.68141e+00
       5,       2,  106079,   53039,4.21968e-03,3.93096e-02,5.88535e-02,       R,7.07788e+00
       5,       3,  150283,   75141,3.54106e-03,3.93765e-02,5.87394e-02,       R,9.05153e+00
       5,       4,  211377,  105688,2.99199e-03,3.94221e-02,5.86625e-02,       Z,1.18727e+01
       6,       0,  211377,  105688,2.99122e-03,2.72608e-02,4.04444e-02,       R,2.25062e+01
       6,       1,  294385,  147192,2.52243e-03,2.73082e-02,4.03668e-02,       R,2.24792e+01
       6,       2,  411375,  205687,2.14130e-03,2.73407e-02,4.03121e-02,       R,2.79402e+01
       6,       3,  576061,  288030,1.81075e-03,2.73646e-02,4.02713e-02,       R,3.56156e+01
       6,       4,  793317,  396658,1.54200e-03,2.73810e-02,4.02435e-02,       Z,4.66137e+01
       7,       0,  793317,  396658,1.54172e-03,1.93760e-02,2.83445e-02,       R,8.86124e+01
       7,       1, 1105589,  552794,1.30524e-03,1.93934e-02,2.83159e-02,       Z,8.85758e+01
        }\tableDataFive

        \pgfplotstableread[col sep=comma]{
       k,     ell,    ndof,   nElem,        eta,         mu,        res,    case,   maxGradU
       0,       0,     193,      96,1.00896e-01,5.83038e-01,8.99738e-01,       Z,2.15503e-03
       1,       0,     193,      96,1.00842e-01,2.93805e-01,5.01852e-01,       R,6.04181e-03
       1,       1,     271,     135,8.17297e-02,2.99684e-01,4.88116e-01,       R,1.00747e-02
       1,       2,     427,     213,6.99069e-02,3.02660e-01,4.80820e-01,       R,1.56773e-02
       1,       3,     641,     320,5.57489e-02,3.05585e-01,4.73651e-01,       R,2.34373e-02
       1,       4,     977,     488,4.54090e-02,3.07292e-01,4.69418e-01,       R,3.48510e-02
       1,       5,    1381,     690,3.73448e-02,3.08376e-01,4.66697e-01,       R,5.23952e-02
       1,       6,    2189,    1094,2.99783e-02,3.09179e-01,4.64642e-01,       Z,5.40313e-02
       2,       0,    2189,    1094,2.99573e-02,1.62286e-01,2.52940e-01,       R,1.31048e-01
       2,       1,    3195,    1597,2.41318e-02,1.63254e-01,2.51004e-01,       R,1.71958e-01
       2,       2,    4879,    2439,1.96252e-02,1.63857e-01,2.49790e-01,       R,2.30039e-01
       2,       3,    7651,    3825,1.58929e-02,1.64261e-01,2.48955e-01,       R,3.17242e-01
       2,       4,   10867,    5433,1.31219e-02,1.64505e-01,2.48451e-01,       R,4.33678e-01
       2,       5,   16469,    8234,1.06653e-02,1.64683e-01,2.48087e-01,       R,6.22685e-01
       2,       6,   24835,   12417,8.76087e-03,1.64795e-01,2.47854e-01,       Z,9.06923e-01
       3,       0,   24835,   12417,8.75270e-03,9.26565e-02,1.41164e-01,       R,2.16615e+00
       3,       1,   36229,   18114,7.22051e-03,9.27885e-02,1.40921e-01,       R,2.20083e+00
       3,       2,   52175,   26087,6.01421e-03,9.28745e-02,1.40764e-01,       R,2.82780e+00
       3,       3,   74999,   37499,5.00137e-03,9.29345e-02,1.40656e-01,       R,3.74416e+00
       3,       4,  107787,   53893,4.19510e-03,9.29744e-02,1.40582e-01,       R,4.94246e+00
       3,       5,  152169,   76084,3.51567e-03,9.30026e-02,1.40530e-01,       R,6.88397e+00
       3,       6,  216029,  108014,2.96418e-03,9.30218e-02,1.40495e-01,       R,7.01351e+00
       3,       7,  300569,  150284,2.50281e-03,9.30354e-02,1.40471e-01,       Z,9.93719e+00
       4,       0,  300569,  150284,2.50061e-03,5.70285e-02,8.60139e-02,       R,2.29403e+01
       4,       1,  420709,  210354,2.12225e-03,5.70438e-02,8.59874e-02,       R,2.86502e+01
       4,       2,  583813,  291906,1.79808e-03,5.70549e-02,8.59679e-02,       R,2.91304e+01
       4,       3,  807913,  403956,1.52906e-03,5.70628e-02,8.59543e-02,       R,3.82146e+01
       4,       4, 1127527,  563763,1.29441e-03,5.70686e-02,8.59445e-02,       Z,4.98876e+01
        }\tableDataThree

        \pgfplotstableread[col sep=comma]{
       k,     ell,    ndof,   nElem,        eta,         mu,        res,    case,   maxGradU
       0,       0,     193,      96,1.00896e-01,5.83038e-01,8.99738e-01,       Z,2.15503e-03
       1,       0,     193,      96,1.00842e-01,2.93805e-01,5.01852e-01,       R,6.04181e-03
       1,       1,     271,     135,8.17297e-02,2.99684e-01,4.88116e-01,       R,1.00747e-02
       1,       2,     427,     213,6.99069e-02,3.02660e-01,4.80820e-01,       R,1.56773e-02
       1,       3,     641,     320,5.57489e-02,3.05585e-01,4.73651e-01,       R,2.34373e-02
       1,       4,     977,     488,4.54090e-02,3.07292e-01,4.69418e-01,       R,3.48510e-02
       1,       5,    1381,     690,3.73448e-02,3.08376e-01,4.66697e-01,       R,5.23952e-02
       1,       6,    2189,    1094,2.99783e-02,3.09179e-01,4.64642e-01,       R,5.40313e-02
       1,       7,    3195,    1597,2.41491e-02,3.09689e-01,4.63347e-01,       R,8.17904e-02
       1,       8,    4879,    2439,1.96379e-02,3.10007e-01,4.62537e-01,       R,1.23970e-01
       1,       9,    7655,    3827,1.59002e-02,3.10222e-01,4.61985e-01,       R,1.91111e-01
       1,      10,   10859,    5429,1.31339e-02,3.10351e-01,4.61653e-01,       R,2.85329e-01
       1,      11,   16477,    8238,1.06776e-02,3.10445e-01,4.61412e-01,       R,4.45174e-01
       1,      12,   24851,   12425,8.76662e-03,3.10505e-01,4.61258e-01,       Z,6.93721e-01
       2,       0,   24851,   12425,8.75886e-03,1.63148e-01,2.45857e-01,       R,1.93822e+00
       2,       1,   36199,   18099,7.22864e-03,1.63223e-01,2.45700e-01,       R,1.97182e+00
       2,       2,   52107,   26053,6.01776e-03,1.63272e-01,2.45598e-01,       R,2.57913e+00
       2,       3,   74999,   37499,5.00252e-03,1.63307e-01,2.45528e-01,       R,3.47031e+00
       2,       4,  107965,   53982,4.19669e-03,1.63329e-01,2.45481e-01,       R,4.66001e+00
       2,       5,  152609,   76304,3.51382e-03,1.63346e-01,2.45447e-01,       R,6.57992e+00
       2,       6,  216267,  108133,2.96338e-03,1.63356e-01,2.45424e-01,       R,9.37509e+00
       2,       7,  300871,  150435,2.50260e-03,1.63364e-01,2.45409e-01,       R,9.59850e+00
       2,       8,  421901,  210950,2.12337e-03,1.63370e-01,2.45398e-01,       R,1.36190e+01
       2,       9,  585557,  292778,1.79846e-03,1.63373e-01,2.45389e-01,       R,2.03088e+01
       2,      10,  809309,  404654,1.53008e-03,1.63376e-01,2.45384e-01,       R,2.07997e+01
       2,      11, 1128989,  564494,1.29601e-03,1.63378e-01,2.45380e-01,       Z,3.02703e+01
        }\tableDataOne

        %
        %

        \addplot+ [marker1, adaptive, forget plot]
        table [col sep=comma, x=cumulativeNdof, y=estimator] {\tableDataNine};

        \addplot+ [marker2, adaptive, forget plot]
        table [col sep=comma, x=cumulativeNdof, y=estimator] {\tableDataSeven};

        \addplot+ [marker3, adaptive, forget plot]
        table [col sep=comma, x=cumulativeNdof, y=estimator] {\tableDataFive};

        \addplot+ [marker4, adaptive, forget plot]
        table [col sep=comma, x=cumulativeNdof, y=estimator] {\tableDataThree};

        \addplot+ [marker5, adaptive, forget plot]
        table [col sep=comma, x=cumulativeNdof, y=estimator] {\tableDataOne};

        \drawslopetriangle[ST1]{0.5}{7e4}{5e-3} 
    \end{loglogaxis}
\end{tikzpicture}

%% file: figures/Fig02d_Nonlinear_convergence_gamma.tex
\begin{tikzpicture}[>=stealth]
    \begin{loglogaxis}[%
            width            = 5.5cm,%
            xlabel           = {cumulative ndof},%
            ymajorgrids      = true,%
            ymax             = 20,%
            ymin             = 2e-3,%
            font             = \footnotesize,%
            grid style       = {%
                densely dotted,%
                semithick%
            },%
            legend style     = {%
                legend pos = north east,%
                font = \footnotesize%
            },%
        ]

        \addlegendimage{marker1}
        \addlegendentry{\(\gamma = 0.9\)}
        \addlegendimage{marker2}
        \addlegendentry{\(\gamma = 0.7\)}
        \addlegendimage{marker3}
        \addlegendentry{\(\gamma = 0.5\)}
        \addlegendimage{marker4}
        \addlegendentry{\(\gamma = 0.3\)}
        \addlegendimage{marker5}
        \addlegendentry{\(\gamma = 0.1\)}


        \pgfplotstableread[col sep=comma]{
       k,     ell,    ndof,   nElem,        eta,         mu,        res,    case,   maxGradU
       0,       0,     193,      96,1.00896e-01,5.83038e-01,8.99738e-01,       Z,2.15503e-03
       1,       0,     193,      96,1.00842e-01,2.93805e-01,5.01852e-01,       R,6.04181e-03
       1,       1,     271,     135,8.17297e-02,2.99684e-01,4.88116e-01,       Z,1.00747e-02
       2,       0,     271,     135,8.16617e-02,1.56242e-01,2.98824e-01,       Z,1.89306e-02
       3,       0,     271,     135,8.15955e-02,8.64786e-02,2.18991e-01,       R,2.67471e-02
       3,       1,     427,     213,6.97909e-02,9.62580e-02,2.06544e-01,       Z,3.54307e-02
       4,       0,     427,     213,6.97253e-02,5.77212e-02,1.68723e-01,       R,5.01653e-02
       4,       1,     633,     316,5.62695e-02,7.09020e-02,1.54874e-01,       Z,6.32253e-02
       5,       0,     633,     316,5.62384e-02,4.53859e-02,1.32507e-01,       Z,8.79335e-02
       6,       0,     633,     316,5.62149e-02,3.06981e-02,1.22503e-01,       R,1.08203e-01
       6,       1,     937,     468,4.56701e-02,4.49080e-02,1.10924e-01,       Z,1.09260e-01
       7,       0,     937,     468,4.56572e-02,2.97143e-02,1.00952e-01,       Z,1.26100e-01
       8,       0,     937,     468,4.56485e-02,2.03060e-02,9.62142e-02,       R,1.39063e-01
       8,       1,    1309,     654,3.77501e-02,3.27270e-02,8.74973e-02,       Z,1.65001e-01
       9,       0,    1309,     654,3.77440e-02,2.15617e-02,8.13071e-02,       Z,1.96791e-01
      10,       0,    1309,     654,3.77398e-02,1.46174e-02,7.83987e-02,       R,2.21918e-01
      10,       1,    2141,    1070,3.03214e-02,2.68063e-02,6.99016e-02,       Z,2.61717e-01
      11,       0,    2141,    1070,3.03293e-02,1.71424e-02,6.49729e-02,       Z,3.15697e-01
      12,       0,    2141,    1070,3.03347e-02,1.13973e-02,6.27711e-02,       R,3.59023e-01
      12,       1,    3101,    1550,2.45310e-02,2.11736e-02,5.62089e-02,       Z,4.20140e-01
      13,       0,    3101,    1550,2.45329e-02,1.36017e-02,5.23790e-02,       Z,5.07914e-01
      14,       0,    3101,    1550,2.45344e-02,9.04632e-03,5.06600e-02,       R,5.79139e-01
      14,       1,    4611,    2305,2.00855e-02,1.67434e-02,4.57119e-02,       Z,6.66814e-01
      15,       0,    4611,    2305,2.00832e-02,1.08548e-02,4.27472e-02,       Z,8.01035e-01
      16,       0,    4611,    2305,2.00817e-02,7.25705e-03,4.13983e-02,       R,9.10722e-01
      16,       1,    7209,    3604,1.62690e-02,1.38296e-02,3.70372e-02,       Z,1.06319e+00
      17,       0,    7209,    3604,1.62696e-02,8.81547e-03,3.45919e-02,       Z,1.28711e+00
      18,       0,    7209,    3604,1.62701e-02,5.83506e-03,3.35101e-02,       R,1.47191e+00
      18,       1,   10405,    5202,1.34099e-02,1.09060e-02,3.02270e-02,       Z,1.70595e+00
      19,       0,   10405,    5202,1.34113e-02,6.95138e-03,2.83665e-02,       Z,2.06611e+00
      20,       0,   10405,    5202,1.34124e-02,4.59347e-03,2.75480e-02,       R,2.36538e+00
      20,       1,   15749,    7874,1.09239e-02,9.03663e-03,2.47425e-02,       Z,2.38567e+00
      21,       0,   15749,    7874,1.09231e-02,5.79918e-03,2.31691e-02,       Z,2.64613e+00
      22,       0,   15749,    7874,1.09225e-02,3.85009e-03,2.24666e-02,       R,2.85506e+00
      22,       1,   23243,   11621,8.96472e-03,7.33204e-03,2.02599e-02,       Z,3.26478e+00
      23,       0,   23243,   11621,8.96461e-03,4.71221e-03,1.89964e-02,       R,3.78608e+00
      23,       1,   34815,   17407,7.36648e-03,6.95012e-03,1.72182e-02,       Z,4.37836e+00
      24,       0,   34815,   17407,7.36743e-03,4.45108e-03,1.58784e-02,       Z,5.35046e+00
      25,       0,   34815,   17407,7.36813e-03,2.95363e-03,1.52754e-02,       R,6.16870e+00
      25,       1,   49897,   24948,6.13029e-03,5.04309e-03,1.38672e-02,       Z,7.11315e+00
      26,       0,   49897,   24948,6.13083e-03,3.23682e-03,1.29947e-02,       Z,8.65007e+00
      27,       0,   49897,   24948,6.13123e-03,2.14839e-03,1.26058e-02,       R,9.94367e+00
      27,       1,   71927,   35963,5.07506e-03,4.05603e-03,1.14271e-02,       Z,1.13337e+01
      28,       0,   71927,   35963,5.07490e-03,2.62395e-03,1.07382e-02,       Z,1.36648e+01
      29,       0,   71927,   35963,5.07481e-03,1.75212e-03,1.04271e-02,       R,1.56234e+01
      29,       1,  104369,   52184,4.24766e-03,3.28344e-03,9.47892e-03,       Z,1.79842e+01
      30,       0,  104369,   52184,4.24777e-03,2.10161e-03,8.94193e-03,       Z,2.17464e+01
      31,       0,  104369,   52184,4.24785e-03,1.39413e-03,8.70459e-03,       R,2.49139e+01
      31,       1,  148793,   74396,3.55369e-03,2.71278e-03,7.89503e-03,       Z,2.51079e+01
      32,       0,  148793,   74396,3.55397e-03,1.71787e-03,7.46003e-03,       R,2.78847e+01
      32,       1,  208883,  104441,2.99862e-03,2.56711e-03,6.83621e-03,       Z,3.18377e+01
      33,       0,  208883,  104441,2.99873e-03,1.64025e-03,6.37899e-03,       Z,3.76191e+01
      34,       0,  208883,  104441,2.99882e-03,1.08604e-03,6.17594e-03,       R,4.24448e+01
      34,       1,  292283,  146141,2.52507e-03,1.94844e-03,5.64490e-03,       Z,4.81223e+01
      35,       0,  292283,  146141,2.52500e-03,1.26107e-03,5.32361e-03,       Z,5.71812e+01
      36,       0,  292283,  146141,2.52495e-03,8.42382e-04,5.17882e-03,       R,6.47636e+01
      36,       1,  410411,  205205,2.14108e-03,1.58138e-03,4.73709e-03,       Z,7.43316e+01
      37,       0,  410411,  205205,2.14110e-03,1.01269e-03,4.48825e-03,       R,8.92206e+01
      37,       1,  577419,  288709,1.80705e-03,1.53115e-03,4.10657e-03,       Z,8.99622e+01
      38,       0,  577419,  288709,1.80718e-03,9.74680e-04,3.83738e-03,       Z,1.03147e+02
      39,       0,  577419,  288709,1.80727e-03,6.44270e-04,3.71868e-03,       R,1.14021e+02
      39,       1,  793523,  396761,1.54019e-03,1.14417e-03,3.41008e-03,       Z,1.29910e+02
      40,       0,  793523,  396761,1.54024e-03,7.31283e-04,3.22924e-03,       R,1.52856e+02
      40,       1, 1106953,  553476,1.30207e-03,1.10079e-03,2.96734e-03,       Z,1.74297e+02
        }\tableDataNine

        \pgfplotstableread[col sep=comma]{
       k,     ell,    ndof,   nElem,        eta,         mu,        res,    case,   maxGradU
       0,       0,     193,      96,1.00896e-01,5.83038e-01,8.99738e-01,       Z,2.15503e-03
       1,       0,     193,      96,1.00842e-01,2.93805e-01,5.01852e-01,       R,6.04181e-03
       1,       1,     271,     135,8.17297e-02,2.99684e-01,4.88116e-01,       R,1.00747e-02
       1,       2,     427,     213,6.99069e-02,3.02660e-01,4.80820e-01,       Z,1.56773e-02
       2,       0,     427,     213,6.98299e-02,1.58073e-01,2.85109e-01,       R,3.23139e-02
       2,       1,     641,     320,5.57020e-02,1.63586e-01,2.75120e-01,       R,4.30450e-02
       2,       2,     977,     488,4.53764e-02,1.66746e-01,2.69120e-01,       Z,5.79765e-02
       3,       0,     977,     488,4.53452e-02,9.35437e-02,1.70930e-01,       R,1.03840e-01
       3,       1,    1369,     684,3.72854e-02,9.70382e-02,1.65432e-01,       R,1.05584e-01
       3,       2,    2175,    1087,3.02469e-02,9.94572e-02,1.61321e-01,       Z,1.34751e-01
       4,       0,    2175,    1087,3.02307e-02,6.07891e-02,1.10045e-01,       R,2.09471e-01
       4,       1,    3167,    1583,2.42966e-02,6.33946e-02,1.06195e-01,       R,2.58962e-01
       4,       2,    4785,    2392,1.96794e-02,6.49763e-02,1.03720e-01,       Z,3.27499e-01
       5,       0,    4785,    2392,1.96687e-02,4.25774e-02,7.41567e-02,       R,5.21308e-01
       5,       1,    7635,    3817,1.58791e-02,4.41310e-02,7.18007e-02,       Z,6.42873e-01
       6,       0,    7635,    3817,1.58740e-02,3.01536e-02,5.42181e-02,       R,9.44379e-01
       6,       1,   10929,    5464,1.30444e-02,3.14812e-02,5.23073e-02,       R,1.12779e+00
       6,       2,   16349,    8174,1.06930e-02,3.23556e-02,5.09874e-02,       Z,1.40312e+00
       7,       0,   16349,    8174,1.06898e-02,2.24178e-02,3.86262e-02,       R,2.16632e+00
       7,       1,   24281,   12140,8.80195e-03,2.32240e-02,3.74432e-02,       R,2.63082e+00
       7,       2,   35881,   17940,7.22839e-03,2.37609e-02,3.65944e-02,       Z,2.66497e+00
       8,       0,   35881,   17940,7.22785e-03,1.65527e-02,2.76365e-02,       R,3.88803e+00
       8,       1,   51773,   25886,6.03988e-03,1.70222e-02,2.69259e-02,       R,4.66707e+00
       8,       2,   74033,   37016,5.02121e-03,1.73500e-02,2.64190e-02,       Z,5.66099e+00
       9,       0,   74033,   37016,5.02029e-03,1.21277e-02,1.99072e-02,       R,8.70114e+00
       9,       1,  105997,   52998,4.21124e-03,1.24318e-02,1.94411e-02,       R,1.05330e+01
       9,       2,  150561,   75280,3.53042e-03,1.26420e-02,1.91014e-02,       Z,1.06623e+01
      10,       0,  150561,   75280,3.53037e-03,8.84217e-03,1.43573e-02,       R,1.55651e+01
      10,       1,  212419,  106209,2.98265e-03,9.04164e-03,1.40496e-02,       R,1.86211e+01
      10,       2,  295901,  147950,2.51498e-03,9.18272e-03,1.38288e-02,       Z,2.25018e+01
      11,       0,  295901,  147950,2.51462e-03,6.42955e-03,1.03807e-02,       R,3.46095e+01
      11,       1,  413067,  206533,2.13394e-03,6.56573e-03,1.01691e-02,       R,4.18007e+01
      11,       2,  578621,  289310,1.80331e-03,6.66412e-03,1.00085e-02,       Z,4.23118e+01
      12,       0,  578621,  289310,1.80334e-03,4.66396e-03,7.50456e-03,       R,6.18576e+01
      12,       1,  799557,  399778,1.53616e-03,4.75865e-03,7.35694e-03,       R,7.39057e+01
      12,       2, 1113237,  556618,1.30095e-03,4.82826e-03,7.24763e-03,       Z,8.92067e+01
        }\tableDataSeven

        \pgfplotstableread[col sep=comma]{
       k,     ell,    ndof,   nElem,        eta,         mu,        res,    case,   maxGradU
       0,       0,     193,      96,1.00896e-01,5.83038e-01,8.99738e-01,       Z,2.15503e-03
       1,       0,     193,      96,1.00842e-01,2.93805e-01,5.01852e-01,       R,6.04181e-03
       1,       1,     271,     135,8.17297e-02,2.99684e-01,4.88116e-01,       R,1.00747e-02
       1,       2,     427,     213,6.99069e-02,3.02660e-01,4.80820e-01,       R,1.56773e-02
       1,       3,     641,     320,5.57489e-02,3.05585e-01,4.73651e-01,       R,2.34373e-02
       1,       4,     977,     488,4.54090e-02,3.07292e-01,4.69418e-01,       Z,3.48510e-02
       2,       0,     977,     488,4.53788e-02,1.61072e-01,2.62712e-01,       R,8.07874e-02
       2,       1,    1377,     688,3.72684e-02,1.63139e-01,2.58708e-01,       R,1.03494e-01
       2,       2,    2193,    1096,2.99938e-02,1.64632e-01,2.55696e-01,       R,1.27786e-01
       2,       3,    3211,    1605,2.41396e-02,1.65592e-01,2.53770e-01,       Z,1.72957e-01
       3,       0,    3211,    1605,2.41245e-02,9.30548e-02,1.49481e-01,       R,3.36926e-01
       3,       1,    4873,    2436,1.95290e-02,9.41266e-02,1.47635e-01,       R,4.26751e-01
       3,       2,    7715,    3857,1.58063e-02,9.48228e-02,1.46383e-01,       R,5.56299e-01
       3,       3,   10945,    5472,1.30380e-02,9.52429e-02,1.45624e-01,       R,5.56938e-01
       3,       4,   16281,    8140,1.07043e-02,9.55333e-02,1.45101e-01,       Z,7.41440e-01
       4,       0,   16281,    8140,1.06964e-02,5.85311e-02,9.04627e-02,       R,1.40914e+00
       4,       1,   24481,   12240,8.79100e-03,5.88475e-02,8.99279e-02,       R,1.76715e+00
       4,       2,   35987,   17993,7.24273e-03,5.90580e-02,8.95642e-02,       R,2.27521e+00
       4,       3,   51237,   25618,6.06929e-03,5.91901e-02,8.93368e-02,       Z,3.00384e+00
       5,       0,   51237,   25618,6.06727e-03,3.90671e-02,5.92547e-02,       R,5.68832e+00
       5,       1,   73687,   36843,5.02982e-03,3.92142e-02,5.90141e-02,       R,5.68141e+00
       5,       2,  106079,   53039,4.21968e-03,3.93096e-02,5.88535e-02,       R,7.07788e+00
       5,       3,  150283,   75141,3.54106e-03,3.93765e-02,5.87394e-02,       R,9.05153e+00
       5,       4,  211377,  105688,2.99199e-03,3.94221e-02,5.86625e-02,       Z,1.18727e+01
       6,       0,  211377,  105688,2.99122e-03,2.72608e-02,4.04444e-02,       R,2.25062e+01
       6,       1,  294385,  147192,2.52243e-03,2.73082e-02,4.03668e-02,       R,2.24792e+01
       6,       2,  411375,  205687,2.14130e-03,2.73407e-02,4.03121e-02,       R,2.79402e+01
       6,       3,  576061,  288030,1.81075e-03,2.73646e-02,4.02713e-02,       R,3.56156e+01
       6,       4,  793317,  396658,1.54200e-03,2.73810e-02,4.02435e-02,       Z,4.66137e+01
       7,       0,  793317,  396658,1.54172e-03,1.93760e-02,2.83445e-02,       R,8.86124e+01
       7,       1, 1105589,  552794,1.30524e-03,1.93934e-02,2.83159e-02,       Z,8.85758e+01
        }\tableDataFive

        \pgfplotstableread[col sep=comma]{
       k,     ell,    ndof,   nElem,        eta,         mu,        res,    case,   maxGradU
       0,       0,     193,      96,1.00896e-01,5.83038e-01,8.99738e-01,       Z,2.15503e-03
       1,       0,     193,      96,1.00842e-01,2.93805e-01,5.01852e-01,       R,6.04181e-03
       1,       1,     271,     135,8.17297e-02,2.99684e-01,4.88116e-01,       R,1.00747e-02
       1,       2,     427,     213,6.99069e-02,3.02660e-01,4.80820e-01,       R,1.56773e-02
       1,       3,     641,     320,5.57489e-02,3.05585e-01,4.73651e-01,       R,2.34373e-02
       1,       4,     977,     488,4.54090e-02,3.07292e-01,4.69418e-01,       R,3.48510e-02
       1,       5,    1381,     690,3.73448e-02,3.08376e-01,4.66697e-01,       R,5.23952e-02
       1,       6,    2189,    1094,2.99783e-02,3.09179e-01,4.64642e-01,       Z,5.40313e-02
       2,       0,    2189,    1094,2.99573e-02,1.62286e-01,2.52940e-01,       R,1.31048e-01
       2,       1,    3195,    1597,2.41318e-02,1.63254e-01,2.51004e-01,       R,1.71958e-01
       2,       2,    4879,    2439,1.96252e-02,1.63857e-01,2.49790e-01,       R,2.30039e-01
       2,       3,    7651,    3825,1.58929e-02,1.64261e-01,2.48955e-01,       R,3.17242e-01
       2,       4,   10867,    5433,1.31219e-02,1.64505e-01,2.48451e-01,       R,4.33678e-01
       2,       5,   16469,    8234,1.06653e-02,1.64683e-01,2.48087e-01,       R,6.22685e-01
       2,       6,   24835,   12417,8.76087e-03,1.64795e-01,2.47854e-01,       Z,9.06923e-01
       3,       0,   24835,   12417,8.75270e-03,9.26565e-02,1.41164e-01,       R,2.16615e+00
       3,       1,   36229,   18114,7.22051e-03,9.27885e-02,1.40921e-01,       R,2.20083e+00
       3,       2,   52175,   26087,6.01421e-03,9.28745e-02,1.40764e-01,       R,2.82780e+00
       3,       3,   74999,   37499,5.00137e-03,9.29345e-02,1.40656e-01,       R,3.74416e+00
       3,       4,  107787,   53893,4.19510e-03,9.29744e-02,1.40582e-01,       R,4.94246e+00
       3,       5,  152169,   76084,3.51567e-03,9.30026e-02,1.40530e-01,       R,6.88397e+00
       3,       6,  216029,  108014,2.96418e-03,9.30218e-02,1.40495e-01,       R,7.01351e+00
       3,       7,  300569,  150284,2.50281e-03,9.30354e-02,1.40471e-01,       Z,9.93719e+00
       4,       0,  300569,  150284,2.50061e-03,5.70285e-02,8.60139e-02,       R,2.29403e+01
       4,       1,  420709,  210354,2.12225e-03,5.70438e-02,8.59874e-02,       R,2.86502e+01
       4,       2,  583813,  291906,1.79808e-03,5.70549e-02,8.59679e-02,       R,2.91304e+01
       4,       3,  807913,  403956,1.52906e-03,5.70628e-02,8.59543e-02,       R,3.82146e+01
       4,       4, 1127527,  563763,1.29441e-03,5.70686e-02,8.59445e-02,       Z,4.98876e+01
        }\tableDataThree

        \pgfplotstableread[col sep=comma]{
       k,     ell,    ndof,   nElem,        eta,         mu,        res,    case,   maxGradU
       0,       0,     193,      96,1.00896e-01,5.83038e-01,8.99738e-01,       Z,2.15503e-03
       1,       0,     193,      96,1.00842e-01,2.93805e-01,5.01852e-01,       R,6.04181e-03
       1,       1,     271,     135,8.17297e-02,2.99684e-01,4.88116e-01,       R,1.00747e-02
       1,       2,     427,     213,6.99069e-02,3.02660e-01,4.80820e-01,       R,1.56773e-02
       1,       3,     641,     320,5.57489e-02,3.05585e-01,4.73651e-01,       R,2.34373e-02
       1,       4,     977,     488,4.54090e-02,3.07292e-01,4.69418e-01,       R,3.48510e-02
       1,       5,    1381,     690,3.73448e-02,3.08376e-01,4.66697e-01,       R,5.23952e-02
       1,       6,    2189,    1094,2.99783e-02,3.09179e-01,4.64642e-01,       R,5.40313e-02
       1,       7,    3195,    1597,2.41491e-02,3.09689e-01,4.63347e-01,       R,8.17904e-02
       1,       8,    4879,    2439,1.96379e-02,3.10007e-01,4.62537e-01,       R,1.23970e-01
       1,       9,    7655,    3827,1.59002e-02,3.10222e-01,4.61985e-01,       R,1.91111e-01
       1,      10,   10859,    5429,1.31339e-02,3.10351e-01,4.61653e-01,       R,2.85329e-01
       1,      11,   16477,    8238,1.06776e-02,3.10445e-01,4.61412e-01,       R,4.45174e-01
       1,      12,   24851,   12425,8.76662e-03,3.10505e-01,4.61258e-01,       Z,6.93721e-01
       2,       0,   24851,   12425,8.75886e-03,1.63148e-01,2.45857e-01,       R,1.93822e+00
       2,       1,   36199,   18099,7.22864e-03,1.63223e-01,2.45700e-01,       R,1.97182e+00
       2,       2,   52107,   26053,6.01776e-03,1.63272e-01,2.45598e-01,       R,2.57913e+00
       2,       3,   74999,   37499,5.00252e-03,1.63307e-01,2.45528e-01,       R,3.47031e+00
       2,       4,  107965,   53982,4.19669e-03,1.63329e-01,2.45481e-01,       R,4.66001e+00
       2,       5,  152609,   76304,3.51382e-03,1.63346e-01,2.45447e-01,       R,6.57992e+00
       2,       6,  216267,  108133,2.96338e-03,1.63356e-01,2.45424e-01,       R,9.37509e+00
       2,       7,  300871,  150435,2.50260e-03,1.63364e-01,2.45409e-01,       R,9.59850e+00
       2,       8,  421901,  210950,2.12337e-03,1.63370e-01,2.45398e-01,       R,1.36190e+01
       2,       9,  585557,  292778,1.79846e-03,1.63373e-01,2.45389e-01,       R,2.03088e+01
       2,      10,  809309,  404654,1.53008e-03,1.63376e-01,2.45384e-01,       R,2.07997e+01
       2,      11, 1128989,  564494,1.29601e-03,1.63378e-01,2.45380e-01,       Z,3.02703e+01
        }\tableDataOne

        %
        %

        \addplot+ [marker1, adaptive, forget plot]
        table [col sep=comma, x=cumulativeNdof, y=res] {\tableDataNine};

        \addplot+ [marker2, adaptive, forget plot]
        table [col sep=comma, x=cumulativeNdof, y=res] {\tableDataSeven};

        \addplot+ [marker3, adaptive, forget plot]
        table [col sep=comma, x=cumulativeNdof, y=res] {\tableDataFive};

        \addplot+ [marker4, adaptive, forget plot]
        table [col sep=comma, x=cumulativeNdof, y=res] {\tableDataThree};

        \addplot+ [marker5, adaptive, forget plot]
        table [col sep=comma, x=cumulativeNdof, y=res] {\tableDataOne};

        \drawslopetriangle[ST1]{0.5}{4e4}{1e-2} 
    \end{loglogaxis}
\end{tikzpicture}

%% file: figures/Fig03a_Nonlinear_convergence_theta.tex
\begin{tikzpicture}[>=stealth]
    \begin{loglogaxis}[%
            width            = 5.5cm,%
            xlabel           = {cumulative ndof},%
            ylabel           = {error estimator},%
            ymajorgrids      = true,%
            font             = \footnotesize,%
            grid style       = {%
                densely dotted,%
                semithick%
            },%
            legend style     = {%
                legend pos = north east,%
                font = \footnotesize%
            },%
        ]

        \addlegendimage{marker1}
        \addlegendentry{\(\theta = 1.0\)}
        \addlegendimage{marker2}
        \addlegendentry{\(\theta = 0.9\)}
        \addlegendimage{marker3}
        \addlegendentry{\(\theta = 0.7\)}
        \addlegendimage{marker4}
        \addlegendentry{\(\theta = 0.5\)}
        \addlegendimage{marker5}
        \addlegendentry{\(\theta = 0.3\)}

        \pgfplotstableread[col sep=comma]{
       k,     ell,    ndof,   nElem,        eta,         mu,        res,    case,   maxGradU
       0,       0,     193,      96,2.01791e-01,1.16608e+00,2.55303e-01,       Z,8.62011e-03
       1,       0,     193,      96,2.01405e-01,1.51358e-01,2.13723e-01,       R,1.60043e-02
       1,       1,     769,     384,1.12784e-01,2.25284e-01,1.39191e-01,       Z,2.99741e-02
       2,       0,     769,     384,1.12303e-01,8.10536e-02,1.17627e-01,       Z,4.22288e-02
       3,       0,     769,     384,1.12054e-01,3.30142e-02,1.13609e-01,       Z,4.79857e-02
       4,       0,     769,     384,1.11941e-01,1.36903e-02,1.12875e-01,       Z,5.05454e-02
       5,       0,     769,     384,1.11891e-01,5.73025e-03,1.12733e-01,       Z,5.16611e-02
       6,       0,     769,     384,1.11870e-01,2.41213e-03,1.12702e-01,       R,5.21436e-02
       6,       1,    3073,    1536,6.26217e-02,9.27318e-02,6.87761e-02,       Z,7.66793e-02
       7,       0,    3073,    1536,6.24245e-02,2.84779e-02,6.35550e-02,       Z,8.93124e-02
       8,       0,    3073,    1536,6.23220e-02,1.15357e-02,6.26293e-02,       Z,9.52050e-02
       9,       0,    3073,    1536,6.22729e-02,4.74016e-03,6.24579e-02,       Z,9.78745e-02
      10,       0,    3073,    1536,6.22503e-02,1.97414e-03,6.24219e-02,       Z,9.90706e-02
      11,       0,    3073,    1536,6.22402e-02,8.32388e-04,6.24127e-02,       Z,9.96045e-02
      12,       0,    3073,    1536,6.22356e-02,3.54877e-04,6.24098e-02,       R,9.98426e-02
      12,       1,   12289,    6144,3.54291e-02,5.11681e-02,3.87598e-02,       Z,1.37360e-01
      13,       0,   12289,    6144,3.53046e-02,1.59085e-02,3.59017e-02,       Z,1.56870e-01
      14,       0,   12289,    6144,3.52393e-02,6.65822e-03,3.53607e-02,       Z,1.66335e-01
      15,       0,   12289,    6144,3.52072e-02,2.83805e-03,3.52528e-02,       Z,1.70812e-01
      16,       0,   12289,    6144,3.51919e-02,1.23043e-03,3.52282e-02,       Z,1.72909e-01
      17,       0,   12289,    6144,3.51847e-02,5.41589e-04,3.52215e-02,       R,1.73887e-01
      17,       1,   49153,   24576,2.04452e-02,2.86400e-02,2.23215e-02,       Z,2.31860e-01
      18,       0,   49153,   24576,2.03650e-02,9.11792e-03,2.07177e-02,       Z,2.63285e-01
      19,       0,   49153,   24576,2.03228e-02,3.97903e-03,2.03882e-02,       Z,2.79122e-01
      20,       0,   49153,   24576,2.03015e-02,1.77528e-03,2.03163e-02,       Z,2.86919e-01
      21,       0,   49153,   24576,2.02911e-02,8.07364e-04,2.02986e-02,       Z,2.90720e-01
      22,       0,   49153,   24576,2.02860e-02,3.73056e-04,2.02935e-02,       R,2.92566e-01
      22,       1,  196609,   98304,1.20156e-02,1.63489e-02,1.31172e-02,       Z,3.83252e-01
      23,       0,  196609,   98304,1.19665e-02,5.36752e-03,1.21918e-02,       Z,4.34103e-01
      24,       0,  196609,   98304,1.19406e-02,2.45144e-03,1.19844e-02,       Z,4.60631e-01
      25,       0,  196609,   98304,1.19273e-02,1.14588e-03,1.19353e-02,       Z,4.74158e-01
      26,       0,  196609,   98304,1.19206e-02,5.45600e-04,1.19223e-02,       Z,4.80986e-01
      27,       0,  196609,   98304,1.19172e-02,2.63524e-04,1.19184e-02,       R,4.84418e-01
      27,       1,  786433,  393216,7.17220e-03,9.52099e-03,7.84206e-03,       Z,6.27218e-01
      28,       0,  786433,  393216,7.14367e-03,3.23389e-03,7.29427e-03,       Z,7.09811e-01
      29,       0,  786433,  393216,7.12862e-03,1.54320e-03,7.16094e-03,       Z,7.54197e-01
      30,       0,  786433,  393216,7.12077e-03,7.52340e-04,7.12700e-03,       Z,7.77516e-01
      31,       0,  786433,  393216,7.11672e-03,3.72639e-04,7.11762e-03,       Z,7.89639e-01
      32,       0,  786433,  393216,7.11464e-03,1.86736e-04,7.11473e-03,       R,7.95913e-01
      32,       1, 3145729, 1572864,4.33557e-03,5.64409e-03,4.75296e-03,       Z,1.02160e+00
        }\tableDataUniform

        \pgfplotstableread[col sep=comma]{
       k,     ell,    ndof,   nElem,        eta,         mu,        res,    case,   maxGradU
       0,       0,     193,      96,2.01791e-01,1.16608e+00,2.55303e-01,       Z,8.62011e-03
       1,       0,     193,      96,2.01405e-01,1.51358e-01,2.13723e-01,       R,1.60043e-02
       1,       1,     703,     351,1.12733e-01,2.25309e-01,1.39288e-01,       Z,2.98463e-02
       2,       0,     703,     351,1.12324e-01,8.12089e-02,1.17677e-01,       Z,4.20513e-02
       3,       0,     703,     351,1.12101e-01,3.30720e-02,1.13659e-01,       Z,4.77874e-02
       4,       0,     703,     351,1.12000e-01,1.37114e-02,1.12928e-01,       Z,5.03366e-02
       5,       0,     703,     351,1.11956e-01,5.73726e-03,1.12788e-01,       Z,5.14466e-02
       6,       0,     703,     351,1.11936e-01,2.41407e-03,1.12758e-01,       R,5.19260e-02
       6,       1,    2391,    1195,6.61138e-02,9.03576e-02,7.15130e-02,       Z,7.65045e-02
       7,       0,    2391,    1195,6.59581e-02,2.71769e-02,6.69728e-02,       Z,8.91836e-02
       8,       0,    2391,    1195,6.58669e-02,1.10473e-02,6.61665e-02,       Z,9.51254e-02
       9,       0,    2391,    1195,6.58220e-02,4.55063e-03,6.60159e-02,       Z,9.78282e-02
      10,       0,    2391,    1195,6.58012e-02,1.89834e-03,6.59839e-02,       Z,9.90429e-02
      11,       0,    2391,    1195,6.57917e-02,8.01219e-04,6.59756e-02,       R,9.95861e-02
      11,       1,    8199,    4099,3.81656e-02,5.35965e-02,4.15922e-02,       Z,1.36968e-01
      12,       0,    8199,    4099,3.80611e-02,1.66748e-02,3.86640e-02,       Z,1.56670e-01
      13,       0,    8199,    4099,3.80013e-02,6.90029e-03,3.81241e-02,       Z,1.66295e-01
      14,       0,    8199,    4099,3.79709e-02,2.91237e-03,3.80186e-02,       Z,1.70878e-01
      15,       0,    8199,    4099,3.79562e-02,1.25237e-03,3.79949e-02,       Z,1.73036e-01
      16,       0,    8199,    4099,3.79493e-02,5.47694e-04,3.79884e-02,       R,1.74048e-01
      16,       1,   26125,   13062,2.27153e-02,3.04049e-02,2.46314e-02,       Z,2.31822e-01
      17,       0,   26125,   13062,2.26457e-02,9.66573e-03,2.29947e-02,       Z,2.63323e-01
      18,       0,   26125,   13062,2.26064e-02,4.15435e-03,2.26701e-02,       Z,2.79312e-01
      19,       0,   26125,   13062,2.25861e-02,1.83238e-03,2.26009e-02,       Z,2.87237e-01
      20,       0,   26125,   13062,2.25759e-02,8.26653e-04,2.25840e-02,       Z,2.91123e-01
      21,       0,   26125,   13062,2.25710e-02,3.79942e-04,2.25791e-02,       R,2.93020e-01
      21,       1,   81129,   40564,1.37386e-02,1.79121e-02,1.48914e-02,       Z,3.83372e-01
      22,       0,   81129,   40564,1.36958e-02,5.84043e-03,1.39212e-02,       Z,4.34327e-01
      23,       0,   81129,   40564,1.36718e-02,2.61078e-03,1.37142e-02,       Z,4.61091e-01
      24,       0,   81129,   40564,1.36592e-02,1.20111e-03,1.36667e-02,       Z,4.74831e-01
      25,       0,   81129,   40564,1.36527e-02,5.65658e-04,1.36544e-02,       Z,4.81809e-01
      26,       0,   81129,   40564,1.36494e-02,2.71254e-04,1.36507e-02,       R,4.85336e-01
      26,       1,  237083,  118541,8.43651e-03,1.07333e-02,9.13725e-03,       Z,6.27592e-01
      27,       0,  237083,  118541,8.41244e-03,3.56515e-03,8.56004e-03,       Z,7.10314e-01
      28,       0,  237083,  118541,8.39885e-03,1.65040e-03,8.42905e-03,       Z,7.55055e-01
      29,       0,  237083,  118541,8.39153e-03,7.88194e-04,8.39712e-03,       Z,7.78714e-01
      30,       0,  237083,  118541,8.38769e-03,3.85422e-04,8.38847e-03,       Z,7.91089e-01
      31,       0,  237083,  118541,8.38570e-03,1.91715e-04,8.38581e-03,       R,7.97529e-01
      31,       1,  622943,  311471,5.23192e-03,6.55621e-03,5.67583e-03,       Z,1.02237e+00
      32,       0,  622943,  311471,5.21971e-03,2.22893e-03,5.31944e-03,       Z,1.15654e+00
      33,       0,  622943,  311471,5.21271e-03,1.05952e-03,5.23495e-03,       Z,1.23085e+00
      34,       0,  622943,  311471,5.20885e-03,5.20729e-04,5.21351e-03,       Z,1.27109e+00
      35,       0,  622943,  311471,5.20677e-03,2.62109e-04,5.20754e-03,       R,1.29264e+00
      35,       1, 1641709,  820854,3.25188e-03,4.07486e-03,3.53803e-03,       Z,1.64747e+00
        }\tableDataNine

        \pgfplotstableread[col sep=comma]{
       k,     ell,    ndof,   nElem,        eta,         mu,        res,    case,   maxGradU
       0,       0,     193,      96,2.01791e-01,1.16608e+00,2.55303e-01,       Z,8.62011e-03
       1,       0,     193,      96,2.01405e-01,1.51358e-01,2.13723e-01,       R,1.60043e-02
       1,       1,     517,     258,1.27356e-01,2.17379e-01,1.50284e-01,       Z,2.98864e-02
       2,       0,     517,     258,1.27082e-01,7.85000e-02,1.31914e-01,       Z,4.21105e-02
       3,       0,     517,     258,1.26931e-01,3.19091e-02,1.28602e-01,       Z,4.78501e-02
       4,       0,     517,     258,1.26862e-01,1.31929e-02,1.28010e-01,       Z,5.03978e-02
       5,       0,     517,     258,1.26832e-01,5.50208e-03,1.27898e-01,       R,5.15057e-02
       5,       1,    1153,     576,8.44757e-02,9.47649e-02,8.95817e-02,       Z,7.56992e-02
       6,       0,    1153,     576,8.43011e-02,2.93916e-02,8.53637e-02,       Z,8.85102e-02
       7,       0,    1153,     576,8.42127e-02,1.18800e-02,8.46309e-02,       Z,9.45159e-02
       8,       0,    1153,     576,8.41713e-02,4.86653e-03,8.44958e-02,       Z,9.72486e-02
       9,       0,    1153,     576,8.41525e-02,2.01885e-03,8.44674e-02,       R,9.84768e-02
       9,       1,    2919,    1459,5.38574e-02,6.46924e-02,5.74642e-02,       Z,1.35598e-01
      10,       0,    2919,    1459,5.37818e-02,2.00117e-02,5.44459e-02,       Z,1.55551e-01
      11,       0,    2919,    1459,5.37398e-02,8.20002e-03,5.39122e-02,       Z,1.65305e-01
      12,       0,    2919,    1459,5.37187e-02,3.41994e-03,5.38129e-02,       Z,1.69952e-01
      13,       0,    2919,    1459,5.37084e-02,1.45094e-03,5.37920e-02,       R,1.72141e-01
      13,       1,    7097,    3548,3.53435e-02,4.04665e-02,3.74625e-02,       Z,2.29919e-01
      14,       0,    7097,    3548,3.53033e-02,1.24673e-02,3.56747e-02,       Z,2.62086e-01
      15,       0,    7097,    3548,3.52805e-02,5.12783e-03,3.53586e-02,       Z,2.78413e-01
      16,       0,    7097,    3548,3.52685e-02,2.16481e-03,3.52987e-02,       Z,2.86503e-01
      17,       0,    7097,    3548,3.52625e-02,9.38419e-04,3.52857e-02,       R,2.90469e-01
      17,       1,   16307,    8153,2.30790e-02,2.66774e-02,2.44940e-02,       Z,3.80420e-01
      18,       0,   16307,    8153,2.30492e-02,8.26915e-03,2.32996e-02,       Z,4.32367e-01
      19,       0,   16307,    8153,2.30334e-02,3.46768e-03,2.30808e-02,       Z,4.59662e-01
      20,       0,   16307,    8153,2.30253e-02,1.50041e-03,2.30377e-02,       Z,4.73675e-01
      21,       0,   16307,    8153,2.30212e-02,6.69445e-04,2.30280e-02,       R,4.80791e-01
      21,       1,   39093,   19546,1.49609e-02,1.75098e-02,1.58944e-02,       Z,6.22794e-01
      22,       0,   39093,   19546,1.49524e-02,5.38533e-03,1.51205e-02,       Z,7.07533e-01
      23,       0,   39093,   19546,1.49470e-02,2.27216e-03,1.49788e-02,       Z,7.53377e-01
      24,       0,   39093,   19546,1.49438e-02,9.96672e-04,1.49507e-02,       Z,7.77617e-01
      25,       0,   39093,   19546,1.49420e-02,4.54381e-04,1.49444e-02,       R,7.90293e-01
      25,       1,   88995,   44497,9.87849e-03,1.12198e-02,1.04873e-02,       Z,1.01427e+00
      26,       0,   88995,   44497,9.87129e-03,3.53985e-03,9.98154e-03,       Z,1.15182e+00
      27,       0,   88995,   44497,9.86735e-03,1.50206e-03,9.88825e-03,       Z,1.22804e+00
      28,       0,   88995,   44497,9.86523e-03,6.66675e-04,9.86947e-03,       Z,1.26933e+00
      29,       0,   88995,   44497,9.86410e-03,3.09466e-04,9.86515e-03,       R,1.29144e+00
      29,       1,  197781,   98890,6.56978e-03,7.36439e-03,6.95204e-03,       Z,1.64627e+00
      30,       0,  197781,   98890,6.56691e-03,2.28134e-03,6.63825e-03,       Z,1.86919e+00
      31,       0,  197781,   98890,6.56528e-03,9.80535e-04,6.57939e-03,       Z,1.99499e+00
      32,       0,  197781,   98890,6.56437e-03,4.42053e-04,6.56733e-03,       Z,2.06439e+00
      33,       0,  197781,   98890,6.56386e-03,2.08874e-04,6.56454e-03,       R,2.10223e+00
      33,       1,  450789,  225394,4.36746e-03,4.90440e-03,4.63479e-03,       Z,2.66504e+00
      34,       0,  450789,  225394,4.36653e-03,1.55379e-03,4.41611e-03,       Z,3.02466e+00
      35,       0,  450789,  225394,4.36596e-03,6.63383e-04,4.37582e-03,       Z,3.23027e+00
      36,       0,  450789,  225394,4.36561e-03,2.98074e-04,4.36773e-03,       Z,3.34518e+00
      37,       0,  450789,  225394,4.36540e-03,1.41075e-04,4.36589e-03,       R,3.40866e+00
      37,       1,  995387,  497693,2.91583e-03,3.25186e-03,3.08660e-03,       Z,4.30238e+00
      38,       0,  995387,  497693,2.91512e-03,1.01444e-03,2.94708e-03,       Z,4.88014e+00
      39,       0,  995387,  497693,2.91475e-03,4.35240e-04,2.92119e-03,       Z,5.21334e+00
      40,       0,  995387,  497693,2.91455e-03,1.96692e-04,2.91596e-03,       Z,5.40120e+00
      41,       0,  995387,  497693,2.91444e-03,9.36892e-05,2.91477e-03,       R,5.50588e+00
      41,       1, 2195593, 1097796,1.94762e-03,2.17015e-03,2.06239e-03,       Z,6.92622e+00
        }\tableDataSeven

        \pgfplotstableread[col sep=comma]{
       k,     ell,    ndof,   nElem,        eta,         mu,        res,    case,   maxGradU
       0,       0,     193,      96,2.01791e-01,1.16608e+00,2.55303e-01,       Z,8.62011e-03
       1,       0,     193,      96,2.01405e-01,1.51358e-01,2.13723e-01,       R,1.60043e-02
       1,       1,     375,     187,1.45838e-01,2.05437e-01,1.65286e-01,       Z,3.01275e-02
       2,       0,     375,     187,1.45491e-01,7.55836e-02,1.49989e-01,       Z,4.24912e-02
       3,       0,     375,     187,1.45309e-01,3.07517e-02,1.47267e-01,       Z,4.82856e-02
       4,       0,     375,     187,1.45226e-01,1.27174e-02,1.46776e-01,       R,5.08524e-02
       4,       1,     687,     343,1.05500e-01,1.00609e-01,1.10111e-01,       Z,7.44749e-02
       5,       0,     687,     343,1.05429e-01,2.97934e-02,1.06656e-01,       Z,8.78115e-02
       6,       0,     687,     343,1.05385e-01,1.21745e-02,1.06055e-01,       Z,9.40806e-02
       7,       0,     687,     343,1.05362e-01,5.04232e-03,1.05945e-01,       R,9.69384e-02
       7,       1,    1273,     636,7.71608e-02,7.19220e-02,8.06007e-02,       Z,1.33102e-01
       8,       0,    1273,     636,7.71021e-02,2.27474e-02,7.78147e-02,       Z,1.53452e-01
       9,       0,    1273,     636,7.70727e-02,8.98311e-03,7.73643e-02,       Z,1.63403e-01
      10,       0,    1273,     636,7.70581e-02,3.61714e-03,7.72869e-02,       R,1.68140e-01
      10,       1,    2431,    1215,5.66725e-02,5.23380e-02,5.88568e-02,       Z,2.25717e-01
      11,       0,    2431,    1215,5.66569e-02,1.56181e-02,5.71018e-02,       Z,2.59233e-01
      12,       0,    2431,    1215,5.66470e-02,6.44591e-03,5.67957e-02,       Z,2.76250e-01
      13,       0,    2431,    1215,5.66413e-02,2.70930e-03,5.67398e-02,       R,2.84678e-01
      13,       1,    4477,    2238,4.10719e-02,3.90983e-02,4.28571e-02,       Z,3.73885e-01
      14,       0,    4477,    2238,4.10544e-02,1.21856e-02,4.13600e-02,       Z,4.28056e-01
      15,       0,    4477,    2238,4.10457e-02,4.80906e-03,4.11212e-02,       Z,4.56533e-01
      16,       0,    4477,    2238,4.10413e-02,1.94910e-03,4.10809e-02,       R,4.71147e-01
      16,       1,    8873,    4436,2.95398e-02,2.85585e-02,3.07760e-02,       Z,6.11323e-01
      17,       0,    8873,    4436,2.95405e-02,8.58805e-03,2.97591e-02,       Z,6.99538e-01
      18,       0,    8873,    4436,2.95400e-02,3.49628e-03,2.95884e-02,       Z,7.47273e-01
      19,       0,    8873,    4436,2.95394e-02,1.45941e-03,2.95586e-02,       R,7.72495e-01
      19,       1,   16057,    8028,2.17448e-02,2.00467e-02,2.25988e-02,       Z,9.95871e-01
      20,       0,   16057,    8028,2.17383e-02,6.15893e-03,2.18797e-02,       Z,1.14093e+00
      21,       0,   16057,    8028,2.17355e-02,2.46427e-03,2.17629e-02,       Z,1.22139e+00
      22,       0,   16057,    8028,2.17342e-02,1.01479e-03,2.17429e-02,       R,1.26495e+00
      22,       1,   30435,   15217,1.59042e-02,1.48480e-02,1.65359e-02,       Z,1.61854e+00
      23,       0,   30435,   15217,1.59081e-02,4.50618e-03,1.60142e-02,       Z,1.85455e+00
      24,       0,   30435,   15217,1.59097e-02,1.80946e-03,1.59298e-02,       Z,1.98794e+00
      25,       0,   30435,   15217,1.59102e-02,7.47428e-04,1.59156e-02,       R,2.06153e+00
      25,       1,   55039,   27519,1.17607e-02,1.07415e-02,1.22058e-02,       Z,2.61861e+00
      26,       0,   55039,   27519,1.17598e-02,3.26674e-03,1.18330e-02,       Z,2.99823e+00
      27,       0,   55039,   27519,1.17594e-02,1.31011e-03,1.17725e-02,       R,3.21573e+00
      27,       1,  100409,   50204,8.70079e-03,8.01853e-03,9.05578e-03,       Z,4.08631e+00
      28,       0,  100409,   50204,8.70167e-03,2.50649e-03,8.75957e-03,       Z,4.75456e+00
      29,       0,  100409,   50204,8.70209e-03,9.99064e-04,8.71220e-03,       Z,5.14279e+00
      30,       0,  100409,   50204,8.70227e-03,4.09005e-04,8.70439e-03,       R,5.36239e+00
      30,       1,  177329,   88664,6.53937e-03,5.75616e-03,6.76266e-03,       Z,6.77039e+00
      31,       0,  177329,   88664,6.53962e-03,1.72195e-03,6.57653e-03,       Z,7.76521e+00
      32,       0,  177329,   88664,6.53976e-03,6.93132e-04,6.54620e-03,       Z,8.34385e+00
      33,       0,  177329,   88664,6.53983e-03,2.85590e-04,6.54112e-03,       R,8.67229e+00
      33,       1,  314545,  157272,4.89062e-03,4.35117e-03,5.07479e-03,       Z,1.09005e+01
      34,       0,  314545,  157272,4.89042e-03,1.35527e-03,4.91961e-03,       Z,1.24659e+01
      35,       0,  314545,  157272,4.89037e-03,5.35089e-04,4.89519e-03,       Z,1.33791e+01
      36,       0,  314545,  157272,4.89036e-03,2.16046e-04,4.89124e-03,       R,1.38991e+01
      36,       1,  562481,  281240,3.67352e-03,3.23535e-03,3.79885e-03,       Z,1.74569e+01
      37,       0,  562481,  281240,3.67401e-03,9.65760e-04,3.69466e-03,       Z,1.99567e+01
      38,       0,  562481,  281240,3.67423e-03,3.87712e-04,3.67779e-03,       Z,2.14172e+01
      39,       0,  562481,  281240,3.67433e-03,1.58900e-04,3.67501e-03,       R,2.22503e+01
      39,       1,  990365,  495182,2.77914e-03,2.40880e-03,2.87515e-03,       Z,2.79037e+01
      40,       0,  990365,  495182,2.77902e-03,7.37203e-04,2.79434e-03,       Z,3.18813e+01
      41,       0,  990365,  495182,2.77899e-03,2.92358e-04,2.78151e-03,       R,3.42065e+01
      41,       1, 1692197,  846098,2.11711e-03,1.82375e-03,2.19089e-03,       Z,4.30213e+01
        }\tableDataFive

        \pgfplotstableread[col sep=comma]{
       k,     ell,    ndof,   nElem,        eta,         mu,        res,    case,   maxGradU
       0,       0,     193,      96,2.01791e-01,1.16608e+00,2.55303e-01,       Z,8.62011e-03
       1,       0,     193,      96,2.01405e-01,1.51358e-01,2.13723e-01,       R,1.60043e-02
       1,       1,     271,     135,1.63298e-01,1.91851e-01,1.79901e-01,       Z,2.95074e-02
       2,       0,     271,     135,1.63037e-01,7.16045e-02,1.67435e-01,       Z,4.14307e-02
       3,       0,     271,     135,1.62901e-01,2.90517e-02,1.65264e-01,       R,4.69924e-02
       3,       1,     427,     213,1.39329e-01,8.92647e-02,1.42920e-01,       Z,7.06245e-02
       4,       0,     427,     213,1.39178e-01,2.65099e-02,1.40792e-01,       R,8.57288e-02
       4,       1,     633,     316,1.12356e-01,8.63095e-02,1.16536e-01,       Z,1.20423e-01
       5,       0,     633,     316,1.12328e-01,2.82196e-02,1.13647e-01,       Z,1.47185e-01
       6,       0,     633,     316,1.12315e-01,1.15980e-02,1.13149e-01,       R,1.60334e-01
       6,       1,     937,     468,9.12346e-02,6.65240e-02,9.38424e-02,       Z,1.62853e-01
       7,       0,     937,     468,9.12420e-02,2.03783e-02,9.19383e-02,       Z,1.70365e-01
       8,       0,     937,     468,9.12472e-02,7.78354e-03,9.16595e-02,       R,1.73925e-01
       8,       1,    1309,     654,7.54564e-02,5.18938e-02,7.74173e-02,       Z,2.33313e-01
       9,       0,    1309,     654,7.54574e-02,1.63962e-02,7.59310e-02,       Z,2.67113e-01
      10,       0,    1309,     654,7.54574e-02,6.39634e-03,7.57034e-02,       R,2.84215e-01
      10,       1,    2137,    1068,6.08165e-02,4.51230e-02,6.24039e-02,       Z,3.76092e-01
      11,       0,    2137,    1068,6.08435e-02,1.33715e-02,6.12043e-02,       Z,4.36950e-01
      12,       0,    2137,    1068,6.08539e-02,5.40803e-03,6.10093e-02,       R,4.68905e-01
      12,       1,    3089,    1544,4.92090e-02,3.62065e-02,5.05179e-02,       Z,6.00485e-01
      13,       0,    3089,    1544,4.92085e-02,1.11692e-02,4.94582e-02,       Z,6.94138e-01
      14,       0,    3089,    1544,4.92082e-02,4.33988e-03,4.92970e-02,       R,7.44808e-01
      14,       1,    4589,    2294,4.01562e-02,2.87710e-02,4.11967e-02,       Z,9.64786e-01
      15,       0,    4589,    2294,4.01500e-02,9.09263e-03,4.03372e-02,       Z,1.12408e+00
      16,       0,    4589,    2294,4.01479e-02,3.57172e-03,4.02034e-02,       R,1.21258e+00
      16,       1,    7187,    3593,3.25030e-02,2.38362e-02,3.33021e-02,       Z,1.56129e+00
      17,       0,    7187,    3593,3.25053e-02,7.16600e-03,3.26503e-02,       Z,1.82386e+00
      18,       0,    7187,    3593,3.25067e-02,2.87589e-03,3.25452e-02,       R,1.97272e+00
      18,       1,   10419,    5209,2.68242e-02,1.85855e-02,2.74027e-02,       Z,2.52575e+00
      19,       0,   10419,    5209,2.68268e-02,5.54342e-03,2.69252e-02,       Z,2.95377e+00
      20,       0,   10419,    5209,2.68280e-02,2.17301e-03,2.68520e-02,       R,3.20017e+00
      20,       1,   15673,    7836,2.18321e-02,1.57425e-02,2.23725e-02,       Z,3.24298e+00
      21,       0,   15673,    7836,2.18306e-02,4.87122e-03,2.19192e-02,       Z,3.40512e+00
      22,       0,   15673,    7836,2.18303e-02,1.91498e-03,2.18486e-02,       R,3.49461e+00
      22,       1,   23249,   11624,1.78925e-02,1.26526e-02,1.83225e-02,       Z,4.42420e+00
      23,       0,   23249,   11624,1.78942e-02,3.92574e-03,1.79647e-02,       R,5.04005e+00
      23,       1,   34811,   17405,1.47452e-02,1.08718e-02,1.51273e-02,       Z,6.43375e+00
      24,       0,   34811,   17405,1.47492e-02,3.35415e-03,1.48135e-02,       Z,7.70067e+00
      25,       0,   34811,   17405,1.47509e-02,1.34287e-03,1.47634e-02,       R,8.44703e+00
      25,       1,   49807,   24903,1.22872e-02,8.27143e-03,1.25374e-02,       Z,1.06775e+01
      26,       0,   49807,   24903,1.22882e-02,2.48284e-03,1.23277e-02,       Z,1.25320e+01
      27,       0,   49807,   24903,1.22886e-02,9.69309e-04,1.22958e-02,       R,1.36208e+01
      27,       1,   71459,   35729,1.01768e-02,6.95582e-03,1.04136e-02,       Z,1.69306e+01
      28,       0,   71459,   35729,1.01764e-02,2.20749e-03,1.02143e-02,       Z,1.96559e+01
      29,       0,   71459,   35729,1.01763e-02,8.74291e-04,1.01829e-02,       R,2.12546e+01
      29,       1,  103867,   51933,8.51327e-03,5.64327e-03,8.68291e-03,       Z,2.69515e+01
      30,       0,  103867,   51933,8.51369e-03,1.70441e-03,8.54125e-03,       Z,3.14958e+01
      31,       0,  103867,   51933,8.51393e-03,6.78939e-04,8.51880e-03,       R,3.41693e+01
      31,       1,  148269,   74134,7.12601e-03,4.70829e-03,7.25922e-03,       Z,3.45545e+01
      32,       0,  148269,   74134,7.12693e-03,1.37864e-03,7.14859e-03,       R,3.63128e+01
      32,       1,  207617,  103808,6.02299e-03,4.05184e-03,6.15480e-03,       Z,4.56736e+01
      33,       0,  207617,  103808,6.02330e-03,1.26496e-03,6.04426e-03,       Z,5.25239e+01
      34,       0,  207617,  103808,6.02347e-03,4.99684e-04,6.02701e-03,       R,5.65371e+01
      34,       1,  290993,  145496,5.05900e-03,3.30732e-03,5.16659e-03,       Z,6.98422e+01
      35,       0,  290993,  145496,5.05879e-03,1.04959e-03,5.07584e-03,       R,8.05348e+01
      35,       1,  407169,  203584,4.29572e-03,2.87049e-03,4.38946e-03,       Z,1.02460e+02
      36,       0,  407169,  203584,4.29585e-03,9.01526e-04,4.31101e-03,       Z,1.23282e+02
      37,       0,  407169,  203584,4.29592e-03,3.59901e-04,4.29850e-03,       R,1.35684e+02
      37,       1,  572869,  286434,3.62967e-03,2.32595e-03,3.69334e-03,       Z,1.37730e+02
      38,       0,  572869,  286434,3.63005e-03,6.80692e-04,3.64041e-03,       Z,1.46105e+02
      39,       0,  572869,  286434,3.63022e-03,2.71917e-04,3.63200e-03,       R,1.50865e+02
      39,       1,  787563,  393781,3.09451e-03,1.91740e-03,3.14701e-03,       Z,1.87968e+02
      40,       0,  787563,  393781,3.09452e-03,5.72286e-04,3.10273e-03,       R,2.13507e+02
      40,       1, 1097597,  548798,2.61659e-03,1.74844e-03,2.67739e-03,       Z,2.65291e+02
        }\tableDataThree

        %
        %

        \addplot+ [marker1, adaptive, forget plot]
        table [col sep=comma, x=cumulativeNdof, y=estimator] {\tableDataUniform};

        \addplot+ [marker2, adaptive, forget plot]
        table [col sep=comma, x=cumulativeNdof, y=estimator] {\tableDataNine};

        \addplot+ [marker3, adaptive, forget plot]
        table [col sep=comma, x=cumulativeNdof, y=estimator] {\tableDataSeven};

        \addplot+ [marker4, adaptive, forget plot]
        table [col sep=comma, x=cumulativeNdof, y=estimator] {\tableDataFive};

        \addplot+ [marker5, adaptive, forget plot]
        table [col sep=comma, x=cumulativeNdof, y=estimator] {\tableDataThree};

        \drawslopetriangle[ST1]{0.5}{4.5e4}{1e-2} 
        \drawswappedslopetriangle[ST2]{0.38}{5e6}{3.5e-2} 
    \end{loglogaxis}
\end{tikzpicture}

%% file: figures/Fig03b_Nonlinear_convergence_theta.tex
\begin{tikzpicture}[>=stealth]
    \begin{loglogaxis}[%
            width            = 5.5cm,%
            xlabel           = {cumulative ndof},%
            xmax             = 4e7,%
            ymajorgrids      = true,%
            font             = \footnotesize,%
            grid style       = {%
                densely dotted,%
                semithick%
            },%
            legend style     = {%
                legend pos = north east,%
                font = \footnotesize%
            },%
        ]

        \addlegendimage{marker1}
        \addlegendentry{\(\theta = 1.0\)}
        \addlegendimage{marker2}
        \addlegendentry{\(\theta = 0.9\)}
        \addlegendimage{marker3}
        \addlegendentry{\(\theta = 0.7\)}
        \addlegendimage{marker4}
        \addlegendentry{\(\theta = 0.5\)}
        \addlegendimage{marker5}
        \addlegendentry{\(\theta = 0.3\)}

        \pgfplotstableread[col sep=comma]{
       k,     ell,    ndof,   nElem,        eta,         mu,        res,    case,   maxGradU
       0,       0,     193,      96,2.01791e-01,1.16608e+00,2.55303e-01,       Z,8.62011e-03
       1,       0,     193,      96,2.01405e-01,1.51358e-01,2.13723e-01,       R,1.60043e-02
       1,       1,     769,     384,1.12784e-01,2.25284e-01,1.39191e-01,       Z,2.99741e-02
       2,       0,     769,     384,1.12303e-01,8.10536e-02,1.17627e-01,       Z,4.22288e-02
       3,       0,     769,     384,1.12054e-01,3.30142e-02,1.13609e-01,       Z,4.79857e-02
       4,       0,     769,     384,1.11941e-01,1.36903e-02,1.12875e-01,       Z,5.05454e-02
       5,       0,     769,     384,1.11891e-01,5.73025e-03,1.12733e-01,       Z,5.16611e-02
       6,       0,     769,     384,1.11870e-01,2.41213e-03,1.12702e-01,       R,5.21436e-02
       6,       1,    3073,    1536,6.26217e-02,9.27318e-02,6.87761e-02,       Z,7.66793e-02
       7,       0,    3073,    1536,6.24245e-02,2.84779e-02,6.35550e-02,       Z,8.93124e-02
       8,       0,    3073,    1536,6.23220e-02,1.15357e-02,6.26293e-02,       Z,9.52050e-02
       9,       0,    3073,    1536,6.22729e-02,4.74016e-03,6.24579e-02,       Z,9.78745e-02
      10,       0,    3073,    1536,6.22503e-02,1.97414e-03,6.24219e-02,       Z,9.90706e-02
      11,       0,    3073,    1536,6.22402e-02,8.32388e-04,6.24127e-02,       Z,9.96045e-02
      12,       0,    3073,    1536,6.22356e-02,3.54877e-04,6.24098e-02,       R,9.98426e-02
      12,       1,   12289,    6144,3.54291e-02,5.11681e-02,3.87598e-02,       Z,1.37360e-01
      13,       0,   12289,    6144,3.53046e-02,1.59085e-02,3.59017e-02,       Z,1.56870e-01
      14,       0,   12289,    6144,3.52393e-02,6.65822e-03,3.53607e-02,       Z,1.66335e-01
      15,       0,   12289,    6144,3.52072e-02,2.83805e-03,3.52528e-02,       Z,1.70812e-01
      16,       0,   12289,    6144,3.51919e-02,1.23043e-03,3.52282e-02,       Z,1.72909e-01
      17,       0,   12289,    6144,3.51847e-02,5.41589e-04,3.52215e-02,       R,1.73887e-01
      17,       1,   49153,   24576,2.04452e-02,2.86400e-02,2.23215e-02,       Z,2.31860e-01
      18,       0,   49153,   24576,2.03650e-02,9.11792e-03,2.07177e-02,       Z,2.63285e-01
      19,       0,   49153,   24576,2.03228e-02,3.97903e-03,2.03882e-02,       Z,2.79122e-01
      20,       0,   49153,   24576,2.03015e-02,1.77528e-03,2.03163e-02,       Z,2.86919e-01
      21,       0,   49153,   24576,2.02911e-02,8.07364e-04,2.02986e-02,       Z,2.90720e-01
      22,       0,   49153,   24576,2.02860e-02,3.73056e-04,2.02935e-02,       R,2.92566e-01
      22,       1,  196609,   98304,1.20156e-02,1.63489e-02,1.31172e-02,       Z,3.83252e-01
      23,       0,  196609,   98304,1.19665e-02,5.36752e-03,1.21918e-02,       Z,4.34103e-01
      24,       0,  196609,   98304,1.19406e-02,2.45144e-03,1.19844e-02,       Z,4.60631e-01
      25,       0,  196609,   98304,1.19273e-02,1.14588e-03,1.19353e-02,       Z,4.74158e-01
      26,       0,  196609,   98304,1.19206e-02,5.45600e-04,1.19223e-02,       Z,4.80986e-01
      27,       0,  196609,   98304,1.19172e-02,2.63524e-04,1.19184e-02,       R,4.84418e-01
      27,       1,  786433,  393216,7.17220e-03,9.52099e-03,7.84206e-03,       Z,6.27218e-01
      28,       0,  786433,  393216,7.14367e-03,3.23389e-03,7.29427e-03,       Z,7.09811e-01
      29,       0,  786433,  393216,7.12862e-03,1.54320e-03,7.16094e-03,       Z,7.54197e-01
      30,       0,  786433,  393216,7.12077e-03,7.52340e-04,7.12700e-03,       Z,7.77516e-01
      31,       0,  786433,  393216,7.11672e-03,3.72639e-04,7.11762e-03,       Z,7.89639e-01
      32,       0,  786433,  393216,7.11464e-03,1.86736e-04,7.11473e-03,       R,7.95913e-01
      32,       1, 3145729, 1572864,4.33557e-03,5.64409e-03,4.75296e-03,       Z,1.02160e+00
        }\tableDataUniform

        \pgfplotstableread[col sep=comma]{
       k,     ell,    ndof,   nElem,        eta,         mu,        res,    case,   maxGradU
       0,       0,     193,      96,2.01791e-01,1.16608e+00,2.55303e-01,       Z,8.62011e-03
       1,       0,     193,      96,2.01405e-01,1.51358e-01,2.13723e-01,       R,1.60043e-02
       1,       1,     703,     351,1.12733e-01,2.25309e-01,1.39288e-01,       Z,2.98463e-02
       2,       0,     703,     351,1.12324e-01,8.12089e-02,1.17677e-01,       Z,4.20513e-02
       3,       0,     703,     351,1.12101e-01,3.30720e-02,1.13659e-01,       Z,4.77874e-02
       4,       0,     703,     351,1.12000e-01,1.37114e-02,1.12928e-01,       Z,5.03366e-02
       5,       0,     703,     351,1.11956e-01,5.73726e-03,1.12788e-01,       Z,5.14466e-02
       6,       0,     703,     351,1.11936e-01,2.41407e-03,1.12758e-01,       R,5.19260e-02
       6,       1,    2391,    1195,6.61138e-02,9.03576e-02,7.15130e-02,       Z,7.65045e-02
       7,       0,    2391,    1195,6.59581e-02,2.71769e-02,6.69728e-02,       Z,8.91836e-02
       8,       0,    2391,    1195,6.58669e-02,1.10473e-02,6.61665e-02,       Z,9.51254e-02
       9,       0,    2391,    1195,6.58220e-02,4.55063e-03,6.60159e-02,       Z,9.78282e-02
      10,       0,    2391,    1195,6.58012e-02,1.89834e-03,6.59839e-02,       Z,9.90429e-02
      11,       0,    2391,    1195,6.57917e-02,8.01219e-04,6.59756e-02,       R,9.95861e-02
      11,       1,    8199,    4099,3.81656e-02,5.35965e-02,4.15922e-02,       Z,1.36968e-01
      12,       0,    8199,    4099,3.80611e-02,1.66748e-02,3.86640e-02,       Z,1.56670e-01
      13,       0,    8199,    4099,3.80013e-02,6.90029e-03,3.81241e-02,       Z,1.66295e-01
      14,       0,    8199,    4099,3.79709e-02,2.91237e-03,3.80186e-02,       Z,1.70878e-01
      15,       0,    8199,    4099,3.79562e-02,1.25237e-03,3.79949e-02,       Z,1.73036e-01
      16,       0,    8199,    4099,3.79493e-02,5.47694e-04,3.79884e-02,       R,1.74048e-01
      16,       1,   26125,   13062,2.27153e-02,3.04049e-02,2.46314e-02,       Z,2.31822e-01
      17,       0,   26125,   13062,2.26457e-02,9.66573e-03,2.29947e-02,       Z,2.63323e-01
      18,       0,   26125,   13062,2.26064e-02,4.15435e-03,2.26701e-02,       Z,2.79312e-01
      19,       0,   26125,   13062,2.25861e-02,1.83238e-03,2.26009e-02,       Z,2.87237e-01
      20,       0,   26125,   13062,2.25759e-02,8.26653e-04,2.25840e-02,       Z,2.91123e-01
      21,       0,   26125,   13062,2.25710e-02,3.79942e-04,2.25791e-02,       R,2.93020e-01
      21,       1,   81129,   40564,1.37386e-02,1.79121e-02,1.48914e-02,       Z,3.83372e-01
      22,       0,   81129,   40564,1.36958e-02,5.84043e-03,1.39212e-02,       Z,4.34327e-01
      23,       0,   81129,   40564,1.36718e-02,2.61078e-03,1.37142e-02,       Z,4.61091e-01
      24,       0,   81129,   40564,1.36592e-02,1.20111e-03,1.36667e-02,       Z,4.74831e-01
      25,       0,   81129,   40564,1.36527e-02,5.65658e-04,1.36544e-02,       Z,4.81809e-01
      26,       0,   81129,   40564,1.36494e-02,2.71254e-04,1.36507e-02,       R,4.85336e-01
      26,       1,  237083,  118541,8.43651e-03,1.07333e-02,9.13725e-03,       Z,6.27592e-01
      27,       0,  237083,  118541,8.41244e-03,3.56515e-03,8.56004e-03,       Z,7.10314e-01
      28,       0,  237083,  118541,8.39885e-03,1.65040e-03,8.42905e-03,       Z,7.55055e-01
      29,       0,  237083,  118541,8.39153e-03,7.88194e-04,8.39712e-03,       Z,7.78714e-01
      30,       0,  237083,  118541,8.38769e-03,3.85422e-04,8.38847e-03,       Z,7.91089e-01
      31,       0,  237083,  118541,8.38570e-03,1.91715e-04,8.38581e-03,       R,7.97529e-01
      31,       1,  622943,  311471,5.23192e-03,6.55621e-03,5.67583e-03,       Z,1.02237e+00
      32,       0,  622943,  311471,5.21971e-03,2.22893e-03,5.31944e-03,       Z,1.15654e+00
      33,       0,  622943,  311471,5.21271e-03,1.05952e-03,5.23495e-03,       Z,1.23085e+00
      34,       0,  622943,  311471,5.20885e-03,5.20729e-04,5.21351e-03,       Z,1.27109e+00
      35,       0,  622943,  311471,5.20677e-03,2.62109e-04,5.20754e-03,       R,1.29264e+00
      35,       1, 1641709,  820854,3.25188e-03,4.07486e-03,3.53803e-03,       Z,1.64747e+00
        }\tableDataNine

        \pgfplotstableread[col sep=comma]{
       k,     ell,    ndof,   nElem,        eta,         mu,        res,    case,   maxGradU
       0,       0,     193,      96,2.01791e-01,1.16608e+00,2.55303e-01,       Z,8.62011e-03
       1,       0,     193,      96,2.01405e-01,1.51358e-01,2.13723e-01,       R,1.60043e-02
       1,       1,     517,     258,1.27356e-01,2.17379e-01,1.50284e-01,       Z,2.98864e-02
       2,       0,     517,     258,1.27082e-01,7.85000e-02,1.31914e-01,       Z,4.21105e-02
       3,       0,     517,     258,1.26931e-01,3.19091e-02,1.28602e-01,       Z,4.78501e-02
       4,       0,     517,     258,1.26862e-01,1.31929e-02,1.28010e-01,       Z,5.03978e-02
       5,       0,     517,     258,1.26832e-01,5.50208e-03,1.27898e-01,       R,5.15057e-02
       5,       1,    1153,     576,8.44757e-02,9.47649e-02,8.95817e-02,       Z,7.56992e-02
       6,       0,    1153,     576,8.43011e-02,2.93916e-02,8.53637e-02,       Z,8.85102e-02
       7,       0,    1153,     576,8.42127e-02,1.18800e-02,8.46309e-02,       Z,9.45159e-02
       8,       0,    1153,     576,8.41713e-02,4.86653e-03,8.44958e-02,       Z,9.72486e-02
       9,       0,    1153,     576,8.41525e-02,2.01885e-03,8.44674e-02,       R,9.84768e-02
       9,       1,    2919,    1459,5.38574e-02,6.46924e-02,5.74642e-02,       Z,1.35598e-01
      10,       0,    2919,    1459,5.37818e-02,2.00117e-02,5.44459e-02,       Z,1.55551e-01
      11,       0,    2919,    1459,5.37398e-02,8.20002e-03,5.39122e-02,       Z,1.65305e-01
      12,       0,    2919,    1459,5.37187e-02,3.41994e-03,5.38129e-02,       Z,1.69952e-01
      13,       0,    2919,    1459,5.37084e-02,1.45094e-03,5.37920e-02,       R,1.72141e-01
      13,       1,    7097,    3548,3.53435e-02,4.04665e-02,3.74625e-02,       Z,2.29919e-01
      14,       0,    7097,    3548,3.53033e-02,1.24673e-02,3.56747e-02,       Z,2.62086e-01
      15,       0,    7097,    3548,3.52805e-02,5.12783e-03,3.53586e-02,       Z,2.78413e-01
      16,       0,    7097,    3548,3.52685e-02,2.16481e-03,3.52987e-02,       Z,2.86503e-01
      17,       0,    7097,    3548,3.52625e-02,9.38419e-04,3.52857e-02,       R,2.90469e-01
      17,       1,   16307,    8153,2.30790e-02,2.66774e-02,2.44940e-02,       Z,3.80420e-01
      18,       0,   16307,    8153,2.30492e-02,8.26915e-03,2.32996e-02,       Z,4.32367e-01
      19,       0,   16307,    8153,2.30334e-02,3.46768e-03,2.30808e-02,       Z,4.59662e-01
      20,       0,   16307,    8153,2.30253e-02,1.50041e-03,2.30377e-02,       Z,4.73675e-01
      21,       0,   16307,    8153,2.30212e-02,6.69445e-04,2.30280e-02,       R,4.80791e-01
      21,       1,   39093,   19546,1.49609e-02,1.75098e-02,1.58944e-02,       Z,6.22794e-01
      22,       0,   39093,   19546,1.49524e-02,5.38533e-03,1.51205e-02,       Z,7.07533e-01
      23,       0,   39093,   19546,1.49470e-02,2.27216e-03,1.49788e-02,       Z,7.53377e-01
      24,       0,   39093,   19546,1.49438e-02,9.96672e-04,1.49507e-02,       Z,7.77617e-01
      25,       0,   39093,   19546,1.49420e-02,4.54381e-04,1.49444e-02,       R,7.90293e-01
      25,       1,   88995,   44497,9.87849e-03,1.12198e-02,1.04873e-02,       Z,1.01427e+00
      26,       0,   88995,   44497,9.87129e-03,3.53985e-03,9.98154e-03,       Z,1.15182e+00
      27,       0,   88995,   44497,9.86735e-03,1.50206e-03,9.88825e-03,       Z,1.22804e+00
      28,       0,   88995,   44497,9.86523e-03,6.66675e-04,9.86947e-03,       Z,1.26933e+00
      29,       0,   88995,   44497,9.86410e-03,3.09466e-04,9.86515e-03,       R,1.29144e+00
      29,       1,  197781,   98890,6.56978e-03,7.36439e-03,6.95204e-03,       Z,1.64627e+00
      30,       0,  197781,   98890,6.56691e-03,2.28134e-03,6.63825e-03,       Z,1.86919e+00
      31,       0,  197781,   98890,6.56528e-03,9.80535e-04,6.57939e-03,       Z,1.99499e+00
      32,       0,  197781,   98890,6.56437e-03,4.42053e-04,6.56733e-03,       Z,2.06439e+00
      33,       0,  197781,   98890,6.56386e-03,2.08874e-04,6.56454e-03,       R,2.10223e+00
      33,       1,  450789,  225394,4.36746e-03,4.90440e-03,4.63479e-03,       Z,2.66504e+00
      34,       0,  450789,  225394,4.36653e-03,1.55379e-03,4.41611e-03,       Z,3.02466e+00
      35,       0,  450789,  225394,4.36596e-03,6.63383e-04,4.37582e-03,       Z,3.23027e+00
      36,       0,  450789,  225394,4.36561e-03,2.98074e-04,4.36773e-03,       Z,3.34518e+00
      37,       0,  450789,  225394,4.36540e-03,1.41075e-04,4.36589e-03,       R,3.40866e+00
      37,       1,  995387,  497693,2.91583e-03,3.25186e-03,3.08660e-03,       Z,4.30238e+00
      38,       0,  995387,  497693,2.91512e-03,1.01444e-03,2.94708e-03,       Z,4.88014e+00
      39,       0,  995387,  497693,2.91475e-03,4.35240e-04,2.92119e-03,       Z,5.21334e+00
      40,       0,  995387,  497693,2.91455e-03,1.96692e-04,2.91596e-03,       Z,5.40120e+00
      41,       0,  995387,  497693,2.91444e-03,9.36892e-05,2.91477e-03,       R,5.50588e+00
      41,       1, 2195593, 1097796,1.94762e-03,2.17015e-03,2.06239e-03,       Z,6.92622e+00
        }\tableDataSeven

        \pgfplotstableread[col sep=comma]{
       k,     ell,    ndof,   nElem,        eta,         mu,        res,    case,   maxGradU
       0,       0,     193,      96,2.01791e-01,1.16608e+00,2.55303e-01,       Z,8.62011e-03
       1,       0,     193,      96,2.01405e-01,1.51358e-01,2.13723e-01,       R,1.60043e-02
       1,       1,     375,     187,1.45838e-01,2.05437e-01,1.65286e-01,       Z,3.01275e-02
       2,       0,     375,     187,1.45491e-01,7.55836e-02,1.49989e-01,       Z,4.24912e-02
       3,       0,     375,     187,1.45309e-01,3.07517e-02,1.47267e-01,       Z,4.82856e-02
       4,       0,     375,     187,1.45226e-01,1.27174e-02,1.46776e-01,       R,5.08524e-02
       4,       1,     687,     343,1.05500e-01,1.00609e-01,1.10111e-01,       Z,7.44749e-02
       5,       0,     687,     343,1.05429e-01,2.97934e-02,1.06656e-01,       Z,8.78115e-02
       6,       0,     687,     343,1.05385e-01,1.21745e-02,1.06055e-01,       Z,9.40806e-02
       7,       0,     687,     343,1.05362e-01,5.04232e-03,1.05945e-01,       R,9.69384e-02
       7,       1,    1273,     636,7.71608e-02,7.19220e-02,8.06007e-02,       Z,1.33102e-01
       8,       0,    1273,     636,7.71021e-02,2.27474e-02,7.78147e-02,       Z,1.53452e-01
       9,       0,    1273,     636,7.70727e-02,8.98311e-03,7.73643e-02,       Z,1.63403e-01
      10,       0,    1273,     636,7.70581e-02,3.61714e-03,7.72869e-02,       R,1.68140e-01
      10,       1,    2431,    1215,5.66725e-02,5.23380e-02,5.88568e-02,       Z,2.25717e-01
      11,       0,    2431,    1215,5.66569e-02,1.56181e-02,5.71018e-02,       Z,2.59233e-01
      12,       0,    2431,    1215,5.66470e-02,6.44591e-03,5.67957e-02,       Z,2.76250e-01
      13,       0,    2431,    1215,5.66413e-02,2.70930e-03,5.67398e-02,       R,2.84678e-01
      13,       1,    4477,    2238,4.10719e-02,3.90983e-02,4.28571e-02,       Z,3.73885e-01
      14,       0,    4477,    2238,4.10544e-02,1.21856e-02,4.13600e-02,       Z,4.28056e-01
      15,       0,    4477,    2238,4.10457e-02,4.80906e-03,4.11212e-02,       Z,4.56533e-01
      16,       0,    4477,    2238,4.10413e-02,1.94910e-03,4.10809e-02,       R,4.71147e-01
      16,       1,    8873,    4436,2.95398e-02,2.85585e-02,3.07760e-02,       Z,6.11323e-01
      17,       0,    8873,    4436,2.95405e-02,8.58805e-03,2.97591e-02,       Z,6.99538e-01
      18,       0,    8873,    4436,2.95400e-02,3.49628e-03,2.95884e-02,       Z,7.47273e-01
      19,       0,    8873,    4436,2.95394e-02,1.45941e-03,2.95586e-02,       R,7.72495e-01
      19,       1,   16057,    8028,2.17448e-02,2.00467e-02,2.25988e-02,       Z,9.95871e-01
      20,       0,   16057,    8028,2.17383e-02,6.15893e-03,2.18797e-02,       Z,1.14093e+00
      21,       0,   16057,    8028,2.17355e-02,2.46427e-03,2.17629e-02,       Z,1.22139e+00
      22,       0,   16057,    8028,2.17342e-02,1.01479e-03,2.17429e-02,       R,1.26495e+00
      22,       1,   30435,   15217,1.59042e-02,1.48480e-02,1.65359e-02,       Z,1.61854e+00
      23,       0,   30435,   15217,1.59081e-02,4.50618e-03,1.60142e-02,       Z,1.85455e+00
      24,       0,   30435,   15217,1.59097e-02,1.80946e-03,1.59298e-02,       Z,1.98794e+00
      25,       0,   30435,   15217,1.59102e-02,7.47428e-04,1.59156e-02,       R,2.06153e+00
      25,       1,   55039,   27519,1.17607e-02,1.07415e-02,1.22058e-02,       Z,2.61861e+00
      26,       0,   55039,   27519,1.17598e-02,3.26674e-03,1.18330e-02,       Z,2.99823e+00
      27,       0,   55039,   27519,1.17594e-02,1.31011e-03,1.17725e-02,       R,3.21573e+00
      27,       1,  100409,   50204,8.70079e-03,8.01853e-03,9.05578e-03,       Z,4.08631e+00
      28,       0,  100409,   50204,8.70167e-03,2.50649e-03,8.75957e-03,       Z,4.75456e+00
      29,       0,  100409,   50204,8.70209e-03,9.99064e-04,8.71220e-03,       Z,5.14279e+00
      30,       0,  100409,   50204,8.70227e-03,4.09005e-04,8.70439e-03,       R,5.36239e+00
      30,       1,  177329,   88664,6.53937e-03,5.75616e-03,6.76266e-03,       Z,6.77039e+00
      31,       0,  177329,   88664,6.53962e-03,1.72195e-03,6.57653e-03,       Z,7.76521e+00
      32,       0,  177329,   88664,6.53976e-03,6.93132e-04,6.54620e-03,       Z,8.34385e+00
      33,       0,  177329,   88664,6.53983e-03,2.85590e-04,6.54112e-03,       R,8.67229e+00
      33,       1,  314545,  157272,4.89062e-03,4.35117e-03,5.07479e-03,       Z,1.09005e+01
      34,       0,  314545,  157272,4.89042e-03,1.35527e-03,4.91961e-03,       Z,1.24659e+01
      35,       0,  314545,  157272,4.89037e-03,5.35089e-04,4.89519e-03,       Z,1.33791e+01
      36,       0,  314545,  157272,4.89036e-03,2.16046e-04,4.89124e-03,       R,1.38991e+01
      36,       1,  562481,  281240,3.67352e-03,3.23535e-03,3.79885e-03,       Z,1.74569e+01
      37,       0,  562481,  281240,3.67401e-03,9.65760e-04,3.69466e-03,       Z,1.99567e+01
      38,       0,  562481,  281240,3.67423e-03,3.87712e-04,3.67779e-03,       Z,2.14172e+01
      39,       0,  562481,  281240,3.67433e-03,1.58900e-04,3.67501e-03,       R,2.22503e+01
      39,       1,  990365,  495182,2.77914e-03,2.40880e-03,2.87515e-03,       Z,2.79037e+01
      40,       0,  990365,  495182,2.77902e-03,7.37203e-04,2.79434e-03,       Z,3.18813e+01
      41,       0,  990365,  495182,2.77899e-03,2.92358e-04,2.78151e-03,       R,3.42065e+01
      41,       1, 1692197,  846098,2.11711e-03,1.82375e-03,2.19089e-03,       Z,4.30213e+01
        }\tableDataFive

        \pgfplotstableread[col sep=comma]{
       k,     ell,    ndof,   nElem,        eta,         mu,        res,    case,   maxGradU
       0,       0,     193,      96,2.01791e-01,1.16608e+00,2.55303e-01,       Z,8.62011e-03
       1,       0,     193,      96,2.01405e-01,1.51358e-01,2.13723e-01,       R,1.60043e-02
       1,       1,     271,     135,1.63298e-01,1.91851e-01,1.79901e-01,       Z,2.95074e-02
       2,       0,     271,     135,1.63037e-01,7.16045e-02,1.67435e-01,       Z,4.14307e-02
       3,       0,     271,     135,1.62901e-01,2.90517e-02,1.65264e-01,       R,4.69924e-02
       3,       1,     427,     213,1.39329e-01,8.92647e-02,1.42920e-01,       Z,7.06245e-02
       4,       0,     427,     213,1.39178e-01,2.65099e-02,1.40792e-01,       R,8.57288e-02
       4,       1,     633,     316,1.12356e-01,8.63095e-02,1.16536e-01,       Z,1.20423e-01
       5,       0,     633,     316,1.12328e-01,2.82196e-02,1.13647e-01,       Z,1.47185e-01
       6,       0,     633,     316,1.12315e-01,1.15980e-02,1.13149e-01,       R,1.60334e-01
       6,       1,     937,     468,9.12346e-02,6.65240e-02,9.38424e-02,       Z,1.62853e-01
       7,       0,     937,     468,9.12420e-02,2.03783e-02,9.19383e-02,       Z,1.70365e-01
       8,       0,     937,     468,9.12472e-02,7.78354e-03,9.16595e-02,       R,1.73925e-01
       8,       1,    1309,     654,7.54564e-02,5.18938e-02,7.74173e-02,       Z,2.33313e-01
       9,       0,    1309,     654,7.54574e-02,1.63962e-02,7.59310e-02,       Z,2.67113e-01
      10,       0,    1309,     654,7.54574e-02,6.39634e-03,7.57034e-02,       R,2.84215e-01
      10,       1,    2137,    1068,6.08165e-02,4.51230e-02,6.24039e-02,       Z,3.76092e-01
      11,       0,    2137,    1068,6.08435e-02,1.33715e-02,6.12043e-02,       Z,4.36950e-01
      12,       0,    2137,    1068,6.08539e-02,5.40803e-03,6.10093e-02,       R,4.68905e-01
      12,       1,    3089,    1544,4.92090e-02,3.62065e-02,5.05179e-02,       Z,6.00485e-01
      13,       0,    3089,    1544,4.92085e-02,1.11692e-02,4.94582e-02,       Z,6.94138e-01
      14,       0,    3089,    1544,4.92082e-02,4.33988e-03,4.92970e-02,       R,7.44808e-01
      14,       1,    4589,    2294,4.01562e-02,2.87710e-02,4.11967e-02,       Z,9.64786e-01
      15,       0,    4589,    2294,4.01500e-02,9.09263e-03,4.03372e-02,       Z,1.12408e+00
      16,       0,    4589,    2294,4.01479e-02,3.57172e-03,4.02034e-02,       R,1.21258e+00
      16,       1,    7187,    3593,3.25030e-02,2.38362e-02,3.33021e-02,       Z,1.56129e+00
      17,       0,    7187,    3593,3.25053e-02,7.16600e-03,3.26503e-02,       Z,1.82386e+00
      18,       0,    7187,    3593,3.25067e-02,2.87589e-03,3.25452e-02,       R,1.97272e+00
      18,       1,   10419,    5209,2.68242e-02,1.85855e-02,2.74027e-02,       Z,2.52575e+00
      19,       0,   10419,    5209,2.68268e-02,5.54342e-03,2.69252e-02,       Z,2.95377e+00
      20,       0,   10419,    5209,2.68280e-02,2.17301e-03,2.68520e-02,       R,3.20017e+00
      20,       1,   15673,    7836,2.18321e-02,1.57425e-02,2.23725e-02,       Z,3.24298e+00
      21,       0,   15673,    7836,2.18306e-02,4.87122e-03,2.19192e-02,       Z,3.40512e+00
      22,       0,   15673,    7836,2.18303e-02,1.91498e-03,2.18486e-02,       R,3.49461e+00
      22,       1,   23249,   11624,1.78925e-02,1.26526e-02,1.83225e-02,       Z,4.42420e+00
      23,       0,   23249,   11624,1.78942e-02,3.92574e-03,1.79647e-02,       R,5.04005e+00
      23,       1,   34811,   17405,1.47452e-02,1.08718e-02,1.51273e-02,       Z,6.43375e+00
      24,       0,   34811,   17405,1.47492e-02,3.35415e-03,1.48135e-02,       Z,7.70067e+00
      25,       0,   34811,   17405,1.47509e-02,1.34287e-03,1.47634e-02,       R,8.44703e+00
      25,       1,   49807,   24903,1.22872e-02,8.27143e-03,1.25374e-02,       Z,1.06775e+01
      26,       0,   49807,   24903,1.22882e-02,2.48284e-03,1.23277e-02,       Z,1.25320e+01
      27,       0,   49807,   24903,1.22886e-02,9.69309e-04,1.22958e-02,       R,1.36208e+01
      27,       1,   71459,   35729,1.01768e-02,6.95582e-03,1.04136e-02,       Z,1.69306e+01
      28,       0,   71459,   35729,1.01764e-02,2.20749e-03,1.02143e-02,       Z,1.96559e+01
      29,       0,   71459,   35729,1.01763e-02,8.74291e-04,1.01829e-02,       R,2.12546e+01
      29,       1,  103867,   51933,8.51327e-03,5.64327e-03,8.68291e-03,       Z,2.69515e+01
      30,       0,  103867,   51933,8.51369e-03,1.70441e-03,8.54125e-03,       Z,3.14958e+01
      31,       0,  103867,   51933,8.51393e-03,6.78939e-04,8.51880e-03,       R,3.41693e+01
      31,       1,  148269,   74134,7.12601e-03,4.70829e-03,7.25922e-03,       Z,3.45545e+01
      32,       0,  148269,   74134,7.12693e-03,1.37864e-03,7.14859e-03,       R,3.63128e+01
      32,       1,  207617,  103808,6.02299e-03,4.05184e-03,6.15480e-03,       Z,4.56736e+01
      33,       0,  207617,  103808,6.02330e-03,1.26496e-03,6.04426e-03,       Z,5.25239e+01
      34,       0,  207617,  103808,6.02347e-03,4.99684e-04,6.02701e-03,       R,5.65371e+01
      34,       1,  290993,  145496,5.05900e-03,3.30732e-03,5.16659e-03,       Z,6.98422e+01
      35,       0,  290993,  145496,5.05879e-03,1.04959e-03,5.07584e-03,       R,8.05348e+01
      35,       1,  407169,  203584,4.29572e-03,2.87049e-03,4.38946e-03,       Z,1.02460e+02
      36,       0,  407169,  203584,4.29585e-03,9.01526e-04,4.31101e-03,       Z,1.23282e+02
      37,       0,  407169,  203584,4.29592e-03,3.59901e-04,4.29850e-03,       R,1.35684e+02
      37,       1,  572869,  286434,3.62967e-03,2.32595e-03,3.69334e-03,       Z,1.37730e+02
      38,       0,  572869,  286434,3.63005e-03,6.80692e-04,3.64041e-03,       Z,1.46105e+02
      39,       0,  572869,  286434,3.63022e-03,2.71917e-04,3.63200e-03,       R,1.50865e+02
      39,       1,  787563,  393781,3.09451e-03,1.91740e-03,3.14701e-03,       Z,1.87968e+02
      40,       0,  787563,  393781,3.09452e-03,5.72286e-04,3.10273e-03,       R,2.13507e+02
      40,       1, 1097597,  548798,2.61659e-03,1.74844e-03,2.67739e-03,       Z,2.65291e+02
        }\tableDataThree

        %
        %

        \addplot+ [marker1, adaptive, forget plot]
        table [col sep=comma, x=cumulativeNdof, y=res] {\tableDataUniform};

        \addplot+ [marker2, adaptive, forget plot]
        table [col sep=comma, x=cumulativeNdof, y=res] {\tableDataNine};

        \addplot+ [marker3, adaptive, forget plot]
        table [col sep=comma, x=cumulativeNdof, y=res] {\tableDataSeven};

        \addplot+ [marker4, adaptive, forget plot]
        table [col sep=comma, x=cumulativeNdof, y=res] {\tableDataFive};

        \addplot+ [marker5, adaptive, forget plot]
        table [col sep=comma, x=cumulativeNdof, y=res] {\tableDataThree};

        \drawslopetriangle[ST1]{0.5}{3e4}{1e-2} 
        \drawswappedslopetriangle[ST2]{0.38}{8e6}{1.8e-2} 
    \end{loglogaxis}
\end{tikzpicture}

%% file: figures/Fig03c_Nonlinear_convergence_theta.tex
\begin{tikzpicture}[>=stealth]
    \begin{loglogaxis}[%
            width            = 5.5cm,%
            xlabel           = {cumulative ndof},%
            ylabel           = {error estimator},%
            ymajorgrids      = true,%
            font             = \footnotesize,%
            grid style       = {%
                densely dotted,%
                semithick%
            },%
            legend style     = {%
                legend pos = north east,%
                font = \footnotesize%
            },%
        ]

        \addlegendimage{marker1}
        \addlegendentry{\(\theta = 1.0\)}
        \addlegendimage{marker2}
        \addlegendentry{\(\theta = 0.9\)}
        \addlegendimage{marker3}
        \addlegendentry{\(\theta = 0.7\)}
        \addlegendimage{marker4}
        \addlegendentry{\(\theta = 0.5\)}
        \addlegendimage{marker5}
        \addlegendentry{\(\theta = 0.3\)}

        \pgfplotstableread[col sep=comma]{
       k,     ell,    ndof,   nElem,        eta,         mu,        res,    case,   maxGradU
       0,       0,     193,      96,1.00896e-01,5.83038e-01,8.99738e-01,       Z,2.15503e-03
       1,       0,     193,      96,1.00842e-01,2.93805e-01,5.01852e-01,       R,6.04181e-03
       1,       1,     769,     384,5.64967e-02,3.05448e-01,4.74067e-01,       Z,1.02133e-02
       2,       0,     769,     384,5.63841e-02,1.59918e-01,2.71882e-01,       Z,1.92483e-02
       3,       0,     769,     384,5.62719e-02,8.92128e-02,1.81029e-01,       Z,2.72373e-02
       4,       0,     769,     384,5.61810e-02,5.41496e-02,1.42649e-01,       Z,3.36899e-02
       5,       0,     769,     384,5.61119e-02,3.55249e-02,1.26584e-01,       Z,3.86706e-02
       6,       0,     769,     384,5.60608e-02,2.45154e-02,1.19483e-01,       R,4.24166e-02
       6,       1,    3073,    1536,3.14079e-02,5.25105e-02,9.39606e-02,       Z,5.30379e-02
       7,       0,    3073,    1536,3.13361e-02,3.37825e-02,7.81187e-02,       Z,6.47204e-02
       8,       0,    3073,    1536,3.12796e-02,2.26926e-02,7.03315e-02,       Z,7.38215e-02
       9,       0,    3073,    1536,3.12365e-02,1.56564e-02,6.64459e-02,       Z,8.07343e-02
      10,       0,    3073,    1536,3.12042e-02,1.09626e-02,6.44843e-02,       Z,8.59021e-02
      11,       0,    3073,    1536,3.11804e-02,7.73507e-03,6.34857e-02,       Z,8.97253e-02
      12,       0,    3073,    1536,3.11628e-02,5.47919e-03,6.29738e-02,       R,9.25339e-02
      12,       1,   12289,    6144,1.77500e-02,2.61931e-02,4.84744e-02,       Z,1.09930e-01
      13,       0,   12289,    6144,1.77111e-02,1.64596e-02,4.15607e-02,       Z,1.25945e-01
      14,       0,   12289,    6144,1.76806e-02,1.08460e-02,3.83363e-02,       Z,1.38312e-01
      15,       0,   12289,    6144,1.76572e-02,7.38545e-03,3.67898e-02,       Z,1.47702e-01
      16,       0,   12289,    6144,1.76395e-02,5.13163e-03,3.60274e-02,       Z,1.54752e-01
      17,       0,   12289,    6144,1.76263e-02,3.60781e-03,3.56430e-02,       R,1.60006e-01
      17,       1,   49153,   24576,1.02494e-02,1.47869e-02,2.77041e-02,       Z,1.86912e-01
      18,       0,   49153,   24576,1.02232e-02,9.32789e-03,2.38844e-02,       Z,2.12573e-01
      19,       0,   49153,   24576,1.02028e-02,6.18994e-03,2.20925e-02,       Z,2.32576e-01
      20,       0,   49153,   24576,1.01872e-02,4.25413e-03,2.12207e-02,       Z,2.47931e-01
      21,       0,   49153,   24576,1.01753e-02,2.98734e-03,2.07822e-02,       Z,2.59598e-01
      22,       0,   49153,   24576,1.01665e-02,2.12432e-03,2.05560e-02,       R,2.68398e-01
      22,       1,  196609,   98304,6.02592e-03,8.45920e-03,1.61169e-02,       Z,3.10525e-01
      23,       0,  196609,   98304,6.00933e-03,5.36453e-03,1.39865e-02,       Z,3.51489e-01
      24,       0,  196609,   98304,5.99653e-03,3.59427e-03,1.29793e-02,       Z,3.83736e-01
      25,       0,  196609,   98304,5.98676e-03,2.50093e-03,1.24802e-02,       Z,4.08760e-01
      26,       0,  196609,   98304,5.97934e-03,1.78043e-03,1.22229e-02,       Z,4.27987e-01
      27,       0,  196609,   98304,5.97375e-03,1.28429e-03,1.20865e-02,       R,4.42660e-01
      27,       1,  786433,  393216,3.59759e-03,4.93887e-03,9.55579e-03,       Z,5.08958e-01
      28,       0,  786433,  393216,3.58765e-03,3.15440e-03,8.34046e-03,       Z,5.74636e-01
      29,       0,  786433,  393216,3.58004e-03,2.13778e-03,7.75950e-03,       Z,6.26830e-01
      30,       0,  786433,  393216,3.57424e-03,1.50809e-03,7.46534e-03,       Z,6.67738e-01
      31,       0,  786433,  393216,3.56985e-03,1.08935e-03,7.30969e-03,       Z,6.99493e-01
      32,       0,  786433,  393216,3.56651e-03,7.97329e-04,7.22490e-03,       R,7.23975e-01
      32,       1, 3145729, 1572864,2.17466e-03,2.93711e-03,5.75569e-03,       Z,8.28612e-01
        }\tableDataUniform

        \pgfplotstableread[col sep=comma]{
       k,     ell,    ndof,   nElem,        eta,         mu,        res,    case,   maxGradU
       0,       0,     193,      96,1.00896e-01,5.83038e-01,8.99738e-01,       Z,2.15503e-03
       1,       0,     193,      96,1.00842e-01,2.93805e-01,5.01852e-01,       R,6.04181e-03
       1,       1,     703,     351,5.64663e-02,3.05453e-01,4.74074e-01,       Z,1.01724e-02
       2,       0,     703,     351,5.63683e-02,1.59935e-01,2.71898e-01,       Z,1.91661e-02
       3,       0,     703,     351,5.62692e-02,8.92360e-02,1.81053e-01,       Z,2.71207e-02
       4,       0,     703,     351,5.61885e-02,5.41734e-02,1.42680e-01,       Z,3.35467e-02
       5,       0,     703,     351,5.61270e-02,3.55455e-02,1.26621e-01,       Z,3.85073e-02
       6,       0,     703,     351,5.60814e-02,2.45315e-02,1.19524e-01,       R,4.22385e-02
       6,       1,    2383,    1191,3.31752e-02,5.14424e-02,9.52261e-02,       Z,5.28840e-02
       7,       0,    2383,    1191,3.31156e-02,3.29410e-02,8.03548e-02,       Z,6.45744e-02
       8,       0,    2383,    1191,3.30664e-02,2.20730e-02,7.31928e-02,       Z,7.36873e-02
       9,       0,    2383,    1191,3.30281e-02,1.52166e-02,6.96646e-02,       Z,8.06156e-02
      10,       0,    2383,    1191,3.29990e-02,1.06559e-02,6.78947e-02,       Z,8.58001e-02
      11,       0,    2383,    1191,3.29773e-02,7.52269e-03,6.69957e-02,       R,8.96395e-02
      11,       1,    8179,    4089,1.91331e-02,2.78930e-02,5.21184e-02,       Z,1.06761e-01
      12,       0,    8179,    4089,1.90972e-02,1.76306e-02,4.47827e-02,       Z,1.23466e-01
      13,       0,    8179,    4089,1.90681e-02,1.16561e-02,4.13348e-02,       Z,1.36419e-01
      14,       0,    8179,    4089,1.90454e-02,7.94494e-03,3.96754e-02,       Z,1.46287e-01
      15,       0,    8179,    4089,1.90280e-02,5.51721e-03,3.88576e-02,       Z,1.53718e-01
      16,       0,    8179,    4089,1.90149e-02,3.87323e-03,3.84466e-02,       R,1.59268e-01
      16,       1,   26179,   13089,1.13670e-02,1.57277e-02,3.01960e-02,       Z,1.86043e-01
      17,       0,   26179,   13089,1.13436e-02,9.94288e-03,2.62283e-02,       Z,2.11884e-01
      18,       0,   26179,   13089,1.13248e-02,6.59483e-03,2.43746e-02,       Z,2.32067e-01
      19,       0,   26179,   13089,1.13102e-02,4.52110e-03,2.34789e-02,       Z,2.47594e-01
      20,       0,   26179,   13089,1.12991e-02,3.16343e-03,2.30323e-02,       Z,2.59415e-01
      21,       0,   26179,   13089,1.12906e-02,2.24052e-03,2.28041e-02,       R,2.68350e-01
      21,       1,   81731,   40865,6.87298e-03,9.23360e-03,1.80467e-02,       Z,3.10330e-01
      22,       0,   81731,   40865,6.85840e-03,5.85879e-03,1.57782e-02,       Z,3.51341e-01
      23,       0,   81731,   40865,6.84685e-03,3.91467e-03,1.47142e-02,       Z,3.83683e-01
      24,       0,   81731,   40865,6.83789e-03,2.71040e-03,1.41932e-02,       Z,4.08828e-01
      25,       0,   81731,   40865,6.83101e-03,1.91811e-03,1.39284e-02,       Z,4.28187e-01
      26,       0,   81731,   40865,6.82577e-03,1.37512e-03,1.37901e-02,       R,4.42989e-01
      26,       1,  239123,  119561,4.22392e-03,5.53539e-03,1.09900e-02,       Z,5.09066e-01
      27,       0,  239123,  119561,4.21549e-03,3.52365e-03,9.66221e-03,       Z,5.74754e-01
      28,       0,  239123,  119561,4.20886e-03,2.37068e-03,9.03767e-03,       Z,6.27040e-01
      29,       0,  239123,  119561,4.20371e-03,1.65683e-03,8.72823e-03,       Z,6.68095e-01
      30,       0,  239123,  119561,4.19975e-03,1.18517e-03,8.56827e-03,       Z,7.00024e-01
      31,       0,  239123,  119561,4.19672e-03,8.59448e-04,8.48305e-03,       R,7.24687e-01
      31,       1,  625537,  312768,2.62074e-03,3.38863e-03,6.79931e-03,       Z,8.28984e-01
      32,       0,  625537,  312768,2.61624e-03,2.17064e-03,5.99822e-03,       Z,9.34509e-01
      33,       0,  625537,  312768,2.61271e-03,1.47229e-03,5.61801e-03,       Z,1.01919e+00
      34,       0,  625537,  312768,2.60997e-03,1.03829e-03,5.42711e-03,       Z,1.08624e+00
      35,       0,  625537,  312768,2.60785e-03,7.49694e-04,5.32693e-03,       R,1.13884e+00
      35,       1, 1635367,  817683,1.63056e-03,2.16892e-03,4.30183e-03,       Z,1.30116e+00
        }\tableDataNine

        \pgfplotstableread[col sep=comma]{
       k,     ell,    ndof,   nElem,        eta,         mu,        res,    case,   maxGradU
       0,       0,     193,      96,1.00896e-01,5.83038e-01,8.99738e-01,       Z,2.15503e-03
       1,       0,     193,      96,1.00842e-01,2.93805e-01,5.01852e-01,       R,6.04181e-03
       1,       1,     519,     259,6.38055e-02,3.04005e-01,4.77590e-01,       Z,1.02252e-02
       2,       0,     519,     259,6.37283e-02,1.58959e-01,2.78836e-01,       Z,1.92879e-02
       3,       0,     519,     259,6.36515e-02,8.84681e-02,1.91235e-01,       Z,2.73038e-02
       4,       0,     519,     259,6.35897e-02,5.35346e-02,1.55177e-01,       Z,3.37777e-02
       5,       0,     519,     259,6.35431e-02,3.50246e-02,1.40405e-01,       R,3.87732e-02
       5,       1,    1165,     582,4.20711e-02,5.91140e-02,1.17933e-01,       Z,4.87928e-02
       6,       0,    1165,     582,4.20080e-02,3.83894e-02,1.00551e-01,       Z,6.11941e-02
       7,       0,    1165,     582,4.19582e-02,2.60325e-02,9.22421e-02,       Z,7.09445e-02
       8,       0,    1165,     582,4.19203e-02,1.81007e-02,8.81552e-02,       Z,7.83964e-02
       9,       0,    1165,     582,4.18920e-02,1.27465e-02,8.61034e-02,       R,8.39922e-02
       9,       1,    2941,    1470,2.68518e-02,3.45889e-02,6.97040e-02,       Z,1.00633e-01
      10,       0,    2941,    1470,2.68224e-02,2.20155e-02,6.12483e-02,       Z,1.18523e-01
      11,       0,    2941,    1470,2.67990e-02,1.46315e-02,5.73077e-02,       Z,1.32468e-01
      12,       0,    2941,    1470,2.67810e-02,1.00060e-02,5.54206e-02,       Z,1.43128e-01
      13,       0,    2941,    1470,2.67674e-02,6.96035e-03,5.44948e-02,       R,1.51171e-01
      13,       1,    7087,    3543,1.76907e-02,2.12599e-02,4.45065e-02,       Z,1.77426e-01
      14,       0,    7087,    3543,1.76738e-02,1.34714e-02,3.96055e-02,       Z,2.05035e-01
      15,       0,    7087,    3543,1.76604e-02,8.91304e-03,3.73632e-02,       Z,2.26686e-01
      16,       0,    7087,    3543,1.76502e-02,6.07257e-03,3.63046e-02,       Z,2.43381e-01
      17,       0,    7087,    3543,1.76424e-02,4.21301e-03,3.57903e-02,       R,2.56108e-01
      17,       1,   16103,    8051,1.15985e-02,1.39455e-02,2.91248e-02,       Z,2.97202e-01
      18,       0,   16103,    8051,1.15880e-02,8.80567e-03,2.59301e-02,       Z,3.40819e-01
      19,       0,   16103,    8051,1.15798e-02,5.81820e-03,2.44781e-02,       Z,3.75337e-01
      20,       0,   16103,    8051,1.15736e-02,3.96673e-03,2.37942e-02,       Z,4.02230e-01
      21,       0,   16103,    8051,1.15689e-02,2.75828e-03,2.34616e-02,       R,4.22961e-01
      21,       1,   38099,   19049,7.58359e-03,9.16167e-03,1.90655e-02,       Z,4.87840e-01
      22,       0,   38099,   19049,7.57915e-03,5.77802e-03,1.69687e-02,       Z,5.57789e-01
      23,       0,   38099,   19049,7.57564e-03,3.81662e-03,1.60186e-02,       Z,6.13667e-01
      24,       0,   38099,   19049,7.57290e-03,2.60385e-03,1.55721e-02,       Z,6.57636e-01
      25,       0,   38099,   19049,7.57078e-03,1.81339e-03,1.53550e-02,       R,6.91876e-01
      25,       1,   85689,   42844,5.02520e-03,5.94580e-03,1.25684e-02,       Z,7.94082e-01
      26,       0,   85689,   42844,5.02242e-03,3.77600e-03,1.12238e-02,       Z,9.06422e-01
      27,       0,   85689,   42844,5.02029e-03,2.50698e-03,1.06089e-02,       Z,9.96914e-01
      28,       0,   85689,   42844,5.01865e-03,1.71646e-03,1.03178e-02,       Z,1.06873e+00
      29,       0,   85689,   42844,5.01740e-03,1.19864e-03,1.01757e-02,       R,1.12515e+00
      29,       1,  188823,   94411,3.34892e-03,3.92375e-03,8.32964e-03,       Z,1.28669e+00
      30,       0,  188823,   94411,3.34771e-03,2.47730e-03,7.45448e-03,       Z,1.46728e+00
      31,       0,  188823,   94411,3.34678e-03,1.64025e-03,7.05919e-03,       Z,1.61375e+00
      32,       0,  188823,   94411,3.34606e-03,1.12303e-03,6.87326e-03,       Z,1.73081e+00
      33,       0,  188823,   94411,3.34551e-03,7.85732e-04,6.78258e-03,       R,1.82340e+00
      33,       1,  429041,  214520,2.23149e-03,2.61347e-03,5.56736e-03,       Z,2.07910e+00
      34,       0,  429041,  214520,2.23092e-03,1.66480e-03,4.98149e-03,       Z,2.36912e+00
      35,       0,  429041,  214520,2.23049e-03,1.10838e-03,4.71283e-03,       Z,2.60553e+00
      36,       0,  429041,  214520,2.23016e-03,7.60869e-04,4.58539e-03,       Z,2.79546e+00
      37,       0,  429041,  214520,2.22991e-03,5.32769e-04,4.52308e-03,       R,2.94647e+00
      37,       1,  943037,  471518,1.49213e-03,1.74066e-03,3.70745e-03,       Z,3.35174e+00
      38,       0,  943037,  471518,1.49184e-03,1.10030e-03,3.32063e-03,       Z,3.81653e+00
      39,       0,  943037,  471518,1.49162e-03,7.29146e-04,3.14583e-03,       Z,4.19678e+00
      40,       0,  943037,  471518,1.49146e-03,4.99572e-04,3.06364e-03,       R,4.50337e+00
      40,       1, 2096743, 1048371,9.99384e-04,1.21460e-03,2.53295e-03,       Z,5.13068e+00
        }\tableDataSeven

        \pgfplotstableread[col sep=comma]{
       k,     ell,    ndof,   nElem,        eta,         mu,        res,    case,   maxGradU
       0,       0,     193,      96,1.00896e-01,5.83038e-01,8.99738e-01,       Z,2.15503e-03
       1,       0,     193,      96,1.00842e-01,2.93805e-01,5.01852e-01,       R,6.04181e-03
       1,       1,     375,     187,7.30049e-02,3.01928e-01,4.82719e-01,       Z,1.02544e-02
       2,       0,     375,     187,7.29210e-02,1.57686e-01,2.88665e-01,       Z,1.93538e-02
       3,       0,     375,     187,7.28373e-02,8.75717e-02,2.05067e-01,       Z,2.74022e-02
       4,       0,     375,     187,7.27696e-02,5.28474e-02,1.71519e-01,       R,3.39013e-02
       4,       1,     687,     343,5.28379e-02,7.27765e-02,1.50321e-01,       Z,4.32609e-02
       5,       0,     687,     343,5.27983e-02,4.60496e-02,1.26514e-01,       Z,5.67190e-02
       6,       0,     687,     343,5.27662e-02,3.08746e-02,1.15892e-01,       Z,6.74484e-02
       7,       0,     687,     343,5.27415e-02,2.14333e-02,1.10901e-01,       R,7.57145e-02
       7,       1,    1273,     636,3.86311e-02,4.18173e-02,9.55972e-02,       Z,9.12754e-02
       8,       0,    1273,     636,3.86036e-02,2.73493e-02,8.61269e-02,       Z,1.10755e-01
       9,       0,    1273,     636,3.85826e-02,1.85027e-02,8.16252e-02,       Z,1.26075e-01
      10,       0,    1273,     636,3.85669e-02,1.27766e-02,7.94403e-02,       R,1.37847e-01
      10,       1,    2427,    1213,2.83838e-02,2.90689e-02,6.78536e-02,       Z,1.63262e-01
      11,       0,    2427,    1213,2.83750e-02,1.83994e-02,6.18905e-02,       Z,1.93581e-01
      12,       0,    2427,    1213,2.83680e-02,1.21767e-02,5.92240e-02,       Z,2.17526e-01
      13,       0,    2427,    1213,2.83625e-02,8.30375e-03,5.79809e-02,       R,2.36064e-01
      13,       1,    4475,    2237,2.05880e-02,2.12018e-02,4.94050e-02,       Z,2.75763e-01
      14,       0,    4475,    2237,2.05814e-02,1.35953e-02,4.49953e-02,       Z,3.23475e-01
      15,       0,    4475,    2237,2.05763e-02,9.03279e-03,4.29828e-02,       Z,3.61466e-01
      16,       0,    4475,    2237,2.05725e-02,6.14325e-03,4.20417e-02,       R,3.91168e-01
      16,       1,    8869,    4434,1.47608e-02,1.55912e-02,3.55043e-02,       Z,4.53656e-01
      17,       0,    8869,    4434,1.47595e-02,9.83978e-03,3.22543e-02,       Z,5.29742e-01
      18,       0,    8869,    4434,1.47585e-02,6.47761e-03,3.08090e-02,       Z,5.90869e-01
      19,       0,    8869,    4434,1.47576e-02,4.38897e-03,3.01432e-02,       R,6.39124e-01
      19,       1,   16071,    8035,1.08862e-02,1.08877e-02,2.58301e-02,       Z,7.38544e-01
      20,       0,   16071,    8035,1.08841e-02,6.94290e-03,2.36358e-02,       Z,8.61415e-01
      21,       0,   16071,    8035,1.08826e-02,4.59816e-03,2.26480e-02,       Z,9.61000e-01
      22,       0,   16071,    8035,1.08816e-02,3.12329e-03,2.21898e-02,       R,1.04033e+00
      22,       1,   30395,   15197,7.93599e-03,8.07361e-03,1.88890e-02,       Z,1.19725e+00
      23,       0,   30395,   15197,7.93651e-03,5.11601e-03,1.72509e-02,       Z,1.39518e+00
      24,       0,   30395,   15197,7.93688e-03,3.37330e-03,1.65214e-02,       Z,1.55677e+00
      25,       0,   30395,   15197,7.93712e-03,2.28494e-03,1.61859e-02,       R,1.68644e+00
      25,       1,   55561,   27780,5.86180e-03,5.81878e-03,1.38608e-02,       Z,1.93289e+00
      26,       0,   55561,   27780,5.86173e-03,3.69590e-03,1.27001e-02,       Z,2.24962e+00
      27,       0,   55561,   27780,5.86169e-03,2.44060e-03,1.21823e-02,       Z,2.50963e+00
      28,       0,   55561,   27780,5.86166e-03,1.65446e-03,1.19438e-02,       R,2.71944e+00
      28,       1,   99747,   49873,4.35531e-03,4.25765e-03,1.02691e-02,       Z,3.11041e+00
      29,       0,   99747,   49873,4.35545e-03,2.71858e-03,9.42652e-03,       Z,3.61836e+00
      30,       0,   99747,   49873,4.35556e-03,1.80048e-03,9.04787e-03,       R,4.03706e+00
      30,       1,  177753,   88876,3.26724e-03,3.39673e-03,7.83911e-03,       Z,4.63412e+00
      31,       0,  177753,   88876,3.26744e-03,2.16464e-03,7.13623e-03,       Z,5.50647e+00
      32,       0,  177753,   88876,3.26758e-03,1.43304e-03,6.81919e-03,       Z,6.23357e+00
      33,       0,  177753,   88876,3.26769e-03,9.73115e-04,6.67216e-03,       R,6.82767e+00
      33,       1,  313711,  156855,2.45004e-03,2.37108e-03,5.76884e-03,       Z,7.79692e+00
      34,       0,  313711,  156855,2.44997e-03,1.52224e-03,5.30117e-03,       Z,9.13022e+00
      35,       0,  313711,  156855,2.44993e-03,1.01155e-03,5.08914e-03,       Z,1.02373e+01
      36,       0,  313711,  156855,2.44990e-03,6.87505e-04,4.99061e-03,       R,1.11400e+01
      36,       1,  556573,  278286,1.84423e-03,1.75314e-03,4.30128e-03,       Z,1.26941e+01
      37,       0,  556573,  278286,1.84433e-03,1.10616e-03,3.96507e-03,       Z,1.47852e+01
      38,       0,  556573,  278286,1.84441e-03,7.27001e-04,3.81756e-03,       R,1.65197e+01
      38,       1,  980655,  490327,1.39764e-03,1.40605e-03,3.32444e-03,       Z,1.88977e+01
      39,       0,  980655,  490327,1.39765e-03,8.99735e-04,3.03948e-03,       Z,2.24431e+01
      40,       0,  980655,  490327,1.39766e-03,5.96780e-04,2.91049e-03,       Z,2.54096e+01
      41,       0,  980655,  490327,1.39766e-03,4.05292e-04,2.85061e-03,       R,2.78427e+01
      41,       1, 1665503,  832751,1.06640e-03,9.90213e-04,2.48002e-03,       Z,3.17260e+01
        }\tableDataFive

        \pgfplotstableread[col sep=comma]{
       k,     ell,    ndof,   nElem,        eta,         mu,        res,    case,   maxGradU
       0,       0,     193,      96,1.00896e-01,5.83038e-01,8.99738e-01,       Z,2.15503e-03
       1,       0,     193,      96,1.00842e-01,2.93805e-01,5.01852e-01,       R,6.04181e-03
       1,       1,     271,     135,8.17297e-02,2.99684e-01,4.88116e-01,       Z,1.00747e-02
       2,       0,     271,     135,8.16617e-02,1.56242e-01,2.98824e-01,       Z,1.89306e-02
       3,       0,     271,     135,8.15955e-02,8.64786e-02,2.18991e-01,       R,2.67471e-02
       3,       1,     427,     213,6.97909e-02,9.62580e-02,2.06544e-01,       Z,3.54307e-02
       4,       0,     427,     213,6.97253e-02,5.77212e-02,1.68723e-01,       R,5.01653e-02
       4,       1,     633,     316,5.62695e-02,7.09020e-02,1.54874e-01,       Z,6.32253e-02
       5,       0,     633,     316,5.62384e-02,4.53859e-02,1.32507e-01,       Z,8.79335e-02
       6,       0,     633,     316,5.62149e-02,3.06981e-02,1.22503e-01,       R,1.08203e-01
       6,       1,     937,     468,4.56701e-02,4.49080e-02,1.10924e-01,       Z,1.09260e-01
       7,       0,     937,     468,4.56572e-02,2.97143e-02,1.00952e-01,       Z,1.26100e-01
       8,       0,     937,     468,4.56485e-02,2.03060e-02,9.62142e-02,       R,1.39063e-01
       8,       1,    1309,     654,3.77501e-02,3.27270e-02,8.74973e-02,       Z,1.65001e-01
       9,       0,    1309,     654,3.77440e-02,2.15617e-02,8.13071e-02,       Z,1.96791e-01
      10,       0,    1309,     654,3.77398e-02,1.46174e-02,7.83987e-02,       R,2.21918e-01
      10,       1,    2141,    1070,3.03214e-02,2.68063e-02,6.99016e-02,       Z,2.61717e-01
      11,       0,    2141,    1070,3.03293e-02,1.71424e-02,6.49729e-02,       Z,3.15697e-01
      12,       0,    2141,    1070,3.03347e-02,1.13973e-02,6.27711e-02,       R,3.59023e-01
      12,       1,    3101,    1550,2.45310e-02,2.11736e-02,5.62089e-02,       Z,4.20140e-01
      13,       0,    3101,    1550,2.45329e-02,1.36017e-02,5.23790e-02,       Z,5.07914e-01
      14,       0,    3101,    1550,2.45344e-02,9.04632e-03,5.06600e-02,       R,5.79139e-01
      14,       1,    4611,    2305,2.00855e-02,1.67434e-02,4.57119e-02,       Z,6.66814e-01
      15,       0,    4611,    2305,2.00832e-02,1.08548e-02,4.27472e-02,       Z,8.01035e-01
      16,       0,    4611,    2305,2.00817e-02,7.25705e-03,4.13983e-02,       R,9.10722e-01
      16,       1,    7209,    3604,1.62690e-02,1.38296e-02,3.70372e-02,       Z,1.06319e+00
      17,       0,    7209,    3604,1.62696e-02,8.81547e-03,3.45919e-02,       Z,1.28711e+00
      18,       0,    7209,    3604,1.62701e-02,5.83506e-03,3.35101e-02,       R,1.47191e+00
      18,       1,   10405,    5202,1.34099e-02,1.09060e-02,3.02270e-02,       Z,1.70595e+00
      19,       0,   10405,    5202,1.34113e-02,6.95138e-03,2.83665e-02,       Z,2.06611e+00
      20,       0,   10405,    5202,1.34124e-02,4.59347e-03,2.75480e-02,       R,2.36538e+00
      20,       1,   15749,    7874,1.09239e-02,9.03663e-03,2.47425e-02,       Z,2.38567e+00
      21,       0,   15749,    7874,1.09231e-02,5.79918e-03,2.31691e-02,       Z,2.64613e+00
      22,       0,   15749,    7874,1.09225e-02,3.85009e-03,2.24666e-02,       R,2.85506e+00
      22,       1,   23243,   11621,8.96472e-03,7.33204e-03,2.02599e-02,       Z,3.26478e+00
      23,       0,   23243,   11621,8.96461e-03,4.71221e-03,1.89964e-02,       R,3.78608e+00
      23,       1,   34815,   17407,7.36648e-03,6.95012e-03,1.72182e-02,       Z,4.37836e+00
      24,       0,   34815,   17407,7.36743e-03,4.45108e-03,1.58784e-02,       Z,5.35046e+00
      25,       0,   34815,   17407,7.36813e-03,2.95363e-03,1.52754e-02,       R,6.16870e+00
      25,       1,   49897,   24948,6.13029e-03,5.04309e-03,1.38672e-02,       Z,7.11315e+00
      26,       0,   49897,   24948,6.13083e-03,3.23682e-03,1.29947e-02,       Z,8.65007e+00
      27,       0,   49897,   24948,6.13123e-03,2.14839e-03,1.26058e-02,       R,9.94367e+00
      27,       1,   71927,   35963,5.07506e-03,4.05603e-03,1.14271e-02,       Z,1.13337e+01
      28,       0,   71927,   35963,5.07490e-03,2.62395e-03,1.07382e-02,       Z,1.36648e+01
      29,       0,   71927,   35963,5.07481e-03,1.75212e-03,1.04271e-02,       R,1.56234e+01
      29,       1,  104369,   52184,4.24766e-03,3.28344e-03,9.47892e-03,       Z,1.79842e+01
      30,       0,  104369,   52184,4.24777e-03,2.10161e-03,8.94193e-03,       Z,2.17464e+01
      31,       0,  104369,   52184,4.24785e-03,1.39413e-03,8.70459e-03,       R,2.49139e+01
      31,       1,  148793,   74396,3.55369e-03,2.71278e-03,7.89503e-03,       Z,2.51079e+01
      32,       0,  148793,   74396,3.55397e-03,1.71787e-03,7.46003e-03,       R,2.78847e+01
      32,       1,  208883,  104441,2.99862e-03,2.56711e-03,6.83621e-03,       Z,3.18377e+01
      33,       0,  208883,  104441,2.99873e-03,1.64025e-03,6.37899e-03,       Z,3.76191e+01
      34,       0,  208883,  104441,2.99882e-03,1.08604e-03,6.17594e-03,       R,4.24448e+01
      34,       1,  292283,  146141,2.52507e-03,1.94844e-03,5.64490e-03,       Z,4.81223e+01
      35,       0,  292283,  146141,2.52500e-03,1.26107e-03,5.32361e-03,       Z,5.71812e+01
      36,       0,  292283,  146141,2.52495e-03,8.42382e-04,5.17882e-03,       R,6.47636e+01
      36,       1,  410411,  205205,2.14108e-03,1.58138e-03,4.73709e-03,       Z,7.43316e+01
      37,       0,  410411,  205205,2.14110e-03,1.01269e-03,4.48825e-03,       R,8.92206e+01
      37,       1,  577419,  288709,1.80705e-03,1.53115e-03,4.10657e-03,       Z,8.99622e+01
      38,       0,  577419,  288709,1.80718e-03,9.74680e-04,3.83738e-03,       Z,1.03147e+02
      39,       0,  577419,  288709,1.80727e-03,6.44270e-04,3.71868e-03,       R,1.14021e+02
      39,       1,  793523,  396761,1.54019e-03,1.14417e-03,3.41008e-03,       Z,1.29910e+02
      40,       0,  793523,  396761,1.54024e-03,7.31283e-04,3.22924e-03,       R,1.52856e+02
      40,       1, 1106953,  553476,1.30207e-03,1.10079e-03,2.96734e-03,       Z,1.74297e+02
        }\tableDataThree

        %
        %

        \addplot+ [marker1, adaptive, forget plot]
        table [col sep=comma, x=cumulativeNdof, y=estimator] {\tableDataUniform};

        \addplot+ [marker2, adaptive, forget plot]
        table [col sep=comma, x=cumulativeNdof, y=estimator] {\tableDataNine};

        \addplot+ [marker3, adaptive, forget plot]
        table [col sep=comma, x=cumulativeNdof, y=estimator] {\tableDataSeven};

        \addplot+ [marker4, adaptive, forget plot]
        table [col sep=comma, x=cumulativeNdof, y=estimator] {\tableDataFive};

        \addplot+ [marker5, adaptive, forget plot]
        table [col sep=comma, x=cumulativeNdof, y=estimator] {\tableDataThree};

        \drawslopetriangle[ST1]{0.5}{4e4}{7e-3} 
        \drawswappedslopetriangle[ST2]{0.38}{4e6}{2e-2} 
    \end{loglogaxis}
\end{tikzpicture}

%% file: figures/Fig03d_Nonlinear_convergence_theta.tex
\begin{tikzpicture}[>=stealth]
    \begin{loglogaxis}[%
            width            = 5.5cm,%
            xlabel           = {cumulative ndof},%
            ymajorgrids      = true,%
            font             = \footnotesize,%
            grid style       = {%
                densely dotted,%
                semithick%
            },%
            legend style     = {%
                legend pos = north east,%
                font = \footnotesize%
            },%
        ]

        \addlegendimage{marker1}
        \addlegendentry{\(\theta = 1.0\)}
        \addlegendimage{marker2}
        \addlegendentry{\(\theta = 0.9\)}
        \addlegendimage{marker3}
        \addlegendentry{\(\theta = 0.7\)}
        \addlegendimage{marker4}
        \addlegendentry{\(\theta = 0.5\)}
        \addlegendimage{marker5}
        \addlegendentry{\(\theta = 0.3\)}

        \pgfplotstableread[col sep=comma]{
       k,     ell,    ndof,   nElem,        eta,         mu,        res,    case,   maxGradU
       0,       0,     193,      96,1.00896e-01,5.83038e-01,8.99738e-01,       Z,2.15503e-03
       1,       0,     193,      96,1.00842e-01,2.93805e-01,5.01852e-01,       R,6.04181e-03
       1,       1,     769,     384,5.64967e-02,3.05448e-01,4.74067e-01,       Z,1.02133e-02
       2,       0,     769,     384,5.63841e-02,1.59918e-01,2.71882e-01,       Z,1.92483e-02
       3,       0,     769,     384,5.62719e-02,8.92128e-02,1.81029e-01,       Z,2.72373e-02
       4,       0,     769,     384,5.61810e-02,5.41496e-02,1.42649e-01,       Z,3.36899e-02
       5,       0,     769,     384,5.61119e-02,3.55249e-02,1.26584e-01,       Z,3.86706e-02
       6,       0,     769,     384,5.60608e-02,2.45154e-02,1.19483e-01,       R,4.24166e-02
       6,       1,    3073,    1536,3.14079e-02,5.25105e-02,9.39606e-02,       Z,5.30379e-02
       7,       0,    3073,    1536,3.13361e-02,3.37825e-02,7.81187e-02,       Z,6.47204e-02
       8,       0,    3073,    1536,3.12796e-02,2.26926e-02,7.03315e-02,       Z,7.38215e-02
       9,       0,    3073,    1536,3.12365e-02,1.56564e-02,6.64459e-02,       Z,8.07343e-02
      10,       0,    3073,    1536,3.12042e-02,1.09626e-02,6.44843e-02,       Z,8.59021e-02
      11,       0,    3073,    1536,3.11804e-02,7.73507e-03,6.34857e-02,       Z,8.97253e-02
      12,       0,    3073,    1536,3.11628e-02,5.47919e-03,6.29738e-02,       R,9.25339e-02
      12,       1,   12289,    6144,1.77500e-02,2.61931e-02,4.84744e-02,       Z,1.09930e-01
      13,       0,   12289,    6144,1.77111e-02,1.64596e-02,4.15607e-02,       Z,1.25945e-01
      14,       0,   12289,    6144,1.76806e-02,1.08460e-02,3.83363e-02,       Z,1.38312e-01
      15,       0,   12289,    6144,1.76572e-02,7.38545e-03,3.67898e-02,       Z,1.47702e-01
      16,       0,   12289,    6144,1.76395e-02,5.13163e-03,3.60274e-02,       Z,1.54752e-01
      17,       0,   12289,    6144,1.76263e-02,3.60781e-03,3.56430e-02,       R,1.60006e-01
      17,       1,   49153,   24576,1.02494e-02,1.47869e-02,2.77041e-02,       Z,1.86912e-01
      18,       0,   49153,   24576,1.02232e-02,9.32789e-03,2.38844e-02,       Z,2.12573e-01
      19,       0,   49153,   24576,1.02028e-02,6.18994e-03,2.20925e-02,       Z,2.32576e-01
      20,       0,   49153,   24576,1.01872e-02,4.25413e-03,2.12207e-02,       Z,2.47931e-01
      21,       0,   49153,   24576,1.01753e-02,2.98734e-03,2.07822e-02,       Z,2.59598e-01
      22,       0,   49153,   24576,1.01665e-02,2.12432e-03,2.05560e-02,       R,2.68398e-01
      22,       1,  196609,   98304,6.02592e-03,8.45920e-03,1.61169e-02,       Z,3.10525e-01
      23,       0,  196609,   98304,6.00933e-03,5.36453e-03,1.39865e-02,       Z,3.51489e-01
      24,       0,  196609,   98304,5.99653e-03,3.59427e-03,1.29793e-02,       Z,3.83736e-01
      25,       0,  196609,   98304,5.98676e-03,2.50093e-03,1.24802e-02,       Z,4.08760e-01
      26,       0,  196609,   98304,5.97934e-03,1.78043e-03,1.22229e-02,       Z,4.27987e-01
      27,       0,  196609,   98304,5.97375e-03,1.28429e-03,1.20865e-02,       R,4.42660e-01
      27,       1,  786433,  393216,3.59759e-03,4.93887e-03,9.55579e-03,       Z,5.08958e-01
      28,       0,  786433,  393216,3.58765e-03,3.15440e-03,8.34046e-03,       Z,5.74636e-01
      29,       0,  786433,  393216,3.58004e-03,2.13778e-03,7.75950e-03,       Z,6.26830e-01
      30,       0,  786433,  393216,3.57424e-03,1.50809e-03,7.46534e-03,       Z,6.67738e-01
      31,       0,  786433,  393216,3.56985e-03,1.08935e-03,7.30969e-03,       Z,6.99493e-01
      32,       0,  786433,  393216,3.56651e-03,7.97329e-04,7.22490e-03,       R,7.23975e-01
      32,       1, 3145729, 1572864,2.17466e-03,2.93711e-03,5.75569e-03,       Z,8.28612e-01
        }\tableDataUniform

        \pgfplotstableread[col sep=comma]{
       k,     ell,    ndof,   nElem,        eta,         mu,        res,    case,   maxGradU
       0,       0,     193,      96,1.00896e-01,5.83038e-01,8.99738e-01,       Z,2.15503e-03
       1,       0,     193,      96,1.00842e-01,2.93805e-01,5.01852e-01,       R,6.04181e-03
       1,       1,     703,     351,5.64663e-02,3.05453e-01,4.74074e-01,       Z,1.01724e-02
       2,       0,     703,     351,5.63683e-02,1.59935e-01,2.71898e-01,       Z,1.91661e-02
       3,       0,     703,     351,5.62692e-02,8.92360e-02,1.81053e-01,       Z,2.71207e-02
       4,       0,     703,     351,5.61885e-02,5.41734e-02,1.42680e-01,       Z,3.35467e-02
       5,       0,     703,     351,5.61270e-02,3.55455e-02,1.26621e-01,       Z,3.85073e-02
       6,       0,     703,     351,5.60814e-02,2.45315e-02,1.19524e-01,       R,4.22385e-02
       6,       1,    2383,    1191,3.31752e-02,5.14424e-02,9.52261e-02,       Z,5.28840e-02
       7,       0,    2383,    1191,3.31156e-02,3.29410e-02,8.03548e-02,       Z,6.45744e-02
       8,       0,    2383,    1191,3.30664e-02,2.20730e-02,7.31928e-02,       Z,7.36873e-02
       9,       0,    2383,    1191,3.30281e-02,1.52166e-02,6.96646e-02,       Z,8.06156e-02
      10,       0,    2383,    1191,3.29990e-02,1.06559e-02,6.78947e-02,       Z,8.58001e-02
      11,       0,    2383,    1191,3.29773e-02,7.52269e-03,6.69957e-02,       R,8.96395e-02
      11,       1,    8179,    4089,1.91331e-02,2.78930e-02,5.21184e-02,       Z,1.06761e-01
      12,       0,    8179,    4089,1.90972e-02,1.76306e-02,4.47827e-02,       Z,1.23466e-01
      13,       0,    8179,    4089,1.90681e-02,1.16561e-02,4.13348e-02,       Z,1.36419e-01
      14,       0,    8179,    4089,1.90454e-02,7.94494e-03,3.96754e-02,       Z,1.46287e-01
      15,       0,    8179,    4089,1.90280e-02,5.51721e-03,3.88576e-02,       Z,1.53718e-01
      16,       0,    8179,    4089,1.90149e-02,3.87323e-03,3.84466e-02,       R,1.59268e-01
      16,       1,   26179,   13089,1.13670e-02,1.57277e-02,3.01960e-02,       Z,1.86043e-01
      17,       0,   26179,   13089,1.13436e-02,9.94288e-03,2.62283e-02,       Z,2.11884e-01
      18,       0,   26179,   13089,1.13248e-02,6.59483e-03,2.43746e-02,       Z,2.32067e-01
      19,       0,   26179,   13089,1.13102e-02,4.52110e-03,2.34789e-02,       Z,2.47594e-01
      20,       0,   26179,   13089,1.12991e-02,3.16343e-03,2.30323e-02,       Z,2.59415e-01
      21,       0,   26179,   13089,1.12906e-02,2.24052e-03,2.28041e-02,       R,2.68350e-01
      21,       1,   81731,   40865,6.87298e-03,9.23360e-03,1.80467e-02,       Z,3.10330e-01
      22,       0,   81731,   40865,6.85840e-03,5.85879e-03,1.57782e-02,       Z,3.51341e-01
      23,       0,   81731,   40865,6.84685e-03,3.91467e-03,1.47142e-02,       Z,3.83683e-01
      24,       0,   81731,   40865,6.83789e-03,2.71040e-03,1.41932e-02,       Z,4.08828e-01
      25,       0,   81731,   40865,6.83101e-03,1.91811e-03,1.39284e-02,       Z,4.28187e-01
      26,       0,   81731,   40865,6.82577e-03,1.37512e-03,1.37901e-02,       R,4.42989e-01
      26,       1,  239123,  119561,4.22392e-03,5.53539e-03,1.09900e-02,       Z,5.09066e-01
      27,       0,  239123,  119561,4.21549e-03,3.52365e-03,9.66221e-03,       Z,5.74754e-01
      28,       0,  239123,  119561,4.20886e-03,2.37068e-03,9.03767e-03,       Z,6.27040e-01
      29,       0,  239123,  119561,4.20371e-03,1.65683e-03,8.72823e-03,       Z,6.68095e-01
      30,       0,  239123,  119561,4.19975e-03,1.18517e-03,8.56827e-03,       Z,7.00024e-01
      31,       0,  239123,  119561,4.19672e-03,8.59448e-04,8.48305e-03,       R,7.24687e-01
      31,       1,  625537,  312768,2.62074e-03,3.38863e-03,6.79931e-03,       Z,8.28984e-01
      32,       0,  625537,  312768,2.61624e-03,2.17064e-03,5.99822e-03,       Z,9.34509e-01
      33,       0,  625537,  312768,2.61271e-03,1.47229e-03,5.61801e-03,       Z,1.01919e+00
      34,       0,  625537,  312768,2.60997e-03,1.03829e-03,5.42711e-03,       Z,1.08624e+00
      35,       0,  625537,  312768,2.60785e-03,7.49694e-04,5.32693e-03,       R,1.13884e+00
      35,       1, 1635367,  817683,1.63056e-03,2.16892e-03,4.30183e-03,       Z,1.30116e+00
        }\tableDataNine

        \pgfplotstableread[col sep=comma]{
       k,     ell,    ndof,   nElem,        eta,         mu,        res,    case,   maxGradU
       0,       0,     193,      96,1.00896e-01,5.83038e-01,8.99738e-01,       Z,2.15503e-03
       1,       0,     193,      96,1.00842e-01,2.93805e-01,5.01852e-01,       R,6.04181e-03
       1,       1,     519,     259,6.38055e-02,3.04005e-01,4.77590e-01,       Z,1.02252e-02
       2,       0,     519,     259,6.37283e-02,1.58959e-01,2.78836e-01,       Z,1.92879e-02
       3,       0,     519,     259,6.36515e-02,8.84681e-02,1.91235e-01,       Z,2.73038e-02
       4,       0,     519,     259,6.35897e-02,5.35346e-02,1.55177e-01,       Z,3.37777e-02
       5,       0,     519,     259,6.35431e-02,3.50246e-02,1.40405e-01,       R,3.87732e-02
       5,       1,    1165,     582,4.20711e-02,5.91140e-02,1.17933e-01,       Z,4.87928e-02
       6,       0,    1165,     582,4.20080e-02,3.83894e-02,1.00551e-01,       Z,6.11941e-02
       7,       0,    1165,     582,4.19582e-02,2.60325e-02,9.22421e-02,       Z,7.09445e-02
       8,       0,    1165,     582,4.19203e-02,1.81007e-02,8.81552e-02,       Z,7.83964e-02
       9,       0,    1165,     582,4.18920e-02,1.27465e-02,8.61034e-02,       R,8.39922e-02
       9,       1,    2941,    1470,2.68518e-02,3.45889e-02,6.97040e-02,       Z,1.00633e-01
      10,       0,    2941,    1470,2.68224e-02,2.20155e-02,6.12483e-02,       Z,1.18523e-01
      11,       0,    2941,    1470,2.67990e-02,1.46315e-02,5.73077e-02,       Z,1.32468e-01
      12,       0,    2941,    1470,2.67810e-02,1.00060e-02,5.54206e-02,       Z,1.43128e-01
      13,       0,    2941,    1470,2.67674e-02,6.96035e-03,5.44948e-02,       R,1.51171e-01
      13,       1,    7087,    3543,1.76907e-02,2.12599e-02,4.45065e-02,       Z,1.77426e-01
      14,       0,    7087,    3543,1.76738e-02,1.34714e-02,3.96055e-02,       Z,2.05035e-01
      15,       0,    7087,    3543,1.76604e-02,8.91304e-03,3.73632e-02,       Z,2.26686e-01
      16,       0,    7087,    3543,1.76502e-02,6.07257e-03,3.63046e-02,       Z,2.43381e-01
      17,       0,    7087,    3543,1.76424e-02,4.21301e-03,3.57903e-02,       R,2.56108e-01
      17,       1,   16103,    8051,1.15985e-02,1.39455e-02,2.91248e-02,       Z,2.97202e-01
      18,       0,   16103,    8051,1.15880e-02,8.80567e-03,2.59301e-02,       Z,3.40819e-01
      19,       0,   16103,    8051,1.15798e-02,5.81820e-03,2.44781e-02,       Z,3.75337e-01
      20,       0,   16103,    8051,1.15736e-02,3.96673e-03,2.37942e-02,       Z,4.02230e-01
      21,       0,   16103,    8051,1.15689e-02,2.75828e-03,2.34616e-02,       R,4.22961e-01
      21,       1,   38099,   19049,7.58359e-03,9.16167e-03,1.90655e-02,       Z,4.87840e-01
      22,       0,   38099,   19049,7.57915e-03,5.77802e-03,1.69687e-02,       Z,5.57789e-01
      23,       0,   38099,   19049,7.57564e-03,3.81662e-03,1.60186e-02,       Z,6.13667e-01
      24,       0,   38099,   19049,7.57290e-03,2.60385e-03,1.55721e-02,       Z,6.57636e-01
      25,       0,   38099,   19049,7.57078e-03,1.81339e-03,1.53550e-02,       R,6.91876e-01
      25,       1,   85689,   42844,5.02520e-03,5.94580e-03,1.25684e-02,       Z,7.94082e-01
      26,       0,   85689,   42844,5.02242e-03,3.77600e-03,1.12238e-02,       Z,9.06422e-01
      27,       0,   85689,   42844,5.02029e-03,2.50698e-03,1.06089e-02,       Z,9.96914e-01
      28,       0,   85689,   42844,5.01865e-03,1.71646e-03,1.03178e-02,       Z,1.06873e+00
      29,       0,   85689,   42844,5.01740e-03,1.19864e-03,1.01757e-02,       R,1.12515e+00
      29,       1,  188823,   94411,3.34892e-03,3.92375e-03,8.32964e-03,       Z,1.28669e+00
      30,       0,  188823,   94411,3.34771e-03,2.47730e-03,7.45448e-03,       Z,1.46728e+00
      31,       0,  188823,   94411,3.34678e-03,1.64025e-03,7.05919e-03,       Z,1.61375e+00
      32,       0,  188823,   94411,3.34606e-03,1.12303e-03,6.87326e-03,       Z,1.73081e+00
      33,       0,  188823,   94411,3.34551e-03,7.85732e-04,6.78258e-03,       R,1.82340e+00
      33,       1,  429041,  214520,2.23149e-03,2.61347e-03,5.56736e-03,       Z,2.07910e+00
      34,       0,  429041,  214520,2.23092e-03,1.66480e-03,4.98149e-03,       Z,2.36912e+00
      35,       0,  429041,  214520,2.23049e-03,1.10838e-03,4.71283e-03,       Z,2.60553e+00
      36,       0,  429041,  214520,2.23016e-03,7.60869e-04,4.58539e-03,       Z,2.79546e+00
      37,       0,  429041,  214520,2.22991e-03,5.32769e-04,4.52308e-03,       R,2.94647e+00
      37,       1,  943037,  471518,1.49213e-03,1.74066e-03,3.70745e-03,       Z,3.35174e+00
      38,       0,  943037,  471518,1.49184e-03,1.10030e-03,3.32063e-03,       Z,3.81653e+00
      39,       0,  943037,  471518,1.49162e-03,7.29146e-04,3.14583e-03,       Z,4.19678e+00
      40,       0,  943037,  471518,1.49146e-03,4.99572e-04,3.06364e-03,       R,4.50337e+00
      40,       1, 2096743, 1048371,9.99384e-04,1.21460e-03,2.53295e-03,       Z,5.13068e+00
        }\tableDataSeven

        \pgfplotstableread[col sep=comma]{
       k,     ell,    ndof,   nElem,        eta,         mu,        res,    case,   maxGradU
       0,       0,     193,      96,1.00896e-01,5.83038e-01,8.99738e-01,       Z,2.15503e-03
       1,       0,     193,      96,1.00842e-01,2.93805e-01,5.01852e-01,       R,6.04181e-03
       1,       1,     375,     187,7.30049e-02,3.01928e-01,4.82719e-01,       Z,1.02544e-02
       2,       0,     375,     187,7.29210e-02,1.57686e-01,2.88665e-01,       Z,1.93538e-02
       3,       0,     375,     187,7.28373e-02,8.75717e-02,2.05067e-01,       Z,2.74022e-02
       4,       0,     375,     187,7.27696e-02,5.28474e-02,1.71519e-01,       R,3.39013e-02
       4,       1,     687,     343,5.28379e-02,7.27765e-02,1.50321e-01,       Z,4.32609e-02
       5,       0,     687,     343,5.27983e-02,4.60496e-02,1.26514e-01,       Z,5.67190e-02
       6,       0,     687,     343,5.27662e-02,3.08746e-02,1.15892e-01,       Z,6.74484e-02
       7,       0,     687,     343,5.27415e-02,2.14333e-02,1.10901e-01,       R,7.57145e-02
       7,       1,    1273,     636,3.86311e-02,4.18173e-02,9.55972e-02,       Z,9.12754e-02
       8,       0,    1273,     636,3.86036e-02,2.73493e-02,8.61269e-02,       Z,1.10755e-01
       9,       0,    1273,     636,3.85826e-02,1.85027e-02,8.16252e-02,       Z,1.26075e-01
      10,       0,    1273,     636,3.85669e-02,1.27766e-02,7.94403e-02,       R,1.37847e-01
      10,       1,    2427,    1213,2.83838e-02,2.90689e-02,6.78536e-02,       Z,1.63262e-01
      11,       0,    2427,    1213,2.83750e-02,1.83994e-02,6.18905e-02,       Z,1.93581e-01
      12,       0,    2427,    1213,2.83680e-02,1.21767e-02,5.92240e-02,       Z,2.17526e-01
      13,       0,    2427,    1213,2.83625e-02,8.30375e-03,5.79809e-02,       R,2.36064e-01
      13,       1,    4475,    2237,2.05880e-02,2.12018e-02,4.94050e-02,       Z,2.75763e-01
      14,       0,    4475,    2237,2.05814e-02,1.35953e-02,4.49953e-02,       Z,3.23475e-01
      15,       0,    4475,    2237,2.05763e-02,9.03279e-03,4.29828e-02,       Z,3.61466e-01
      16,       0,    4475,    2237,2.05725e-02,6.14325e-03,4.20417e-02,       R,3.91168e-01
      16,       1,    8869,    4434,1.47608e-02,1.55912e-02,3.55043e-02,       Z,4.53656e-01
      17,       0,    8869,    4434,1.47595e-02,9.83978e-03,3.22543e-02,       Z,5.29742e-01
      18,       0,    8869,    4434,1.47585e-02,6.47761e-03,3.08090e-02,       Z,5.90869e-01
      19,       0,    8869,    4434,1.47576e-02,4.38897e-03,3.01432e-02,       R,6.39124e-01
      19,       1,   16071,    8035,1.08862e-02,1.08877e-02,2.58301e-02,       Z,7.38544e-01
      20,       0,   16071,    8035,1.08841e-02,6.94290e-03,2.36358e-02,       Z,8.61415e-01
      21,       0,   16071,    8035,1.08826e-02,4.59816e-03,2.26480e-02,       Z,9.61000e-01
      22,       0,   16071,    8035,1.08816e-02,3.12329e-03,2.21898e-02,       R,1.04033e+00
      22,       1,   30395,   15197,7.93599e-03,8.07361e-03,1.88890e-02,       Z,1.19725e+00
      23,       0,   30395,   15197,7.93651e-03,5.11601e-03,1.72509e-02,       Z,1.39518e+00
      24,       0,   30395,   15197,7.93688e-03,3.37330e-03,1.65214e-02,       Z,1.55677e+00
      25,       0,   30395,   15197,7.93712e-03,2.28494e-03,1.61859e-02,       R,1.68644e+00
      25,       1,   55561,   27780,5.86180e-03,5.81878e-03,1.38608e-02,       Z,1.93289e+00
      26,       0,   55561,   27780,5.86173e-03,3.69590e-03,1.27001e-02,       Z,2.24962e+00
      27,       0,   55561,   27780,5.86169e-03,2.44060e-03,1.21823e-02,       Z,2.50963e+00
      28,       0,   55561,   27780,5.86166e-03,1.65446e-03,1.19438e-02,       R,2.71944e+00
      28,       1,   99747,   49873,4.35531e-03,4.25765e-03,1.02691e-02,       Z,3.11041e+00
      29,       0,   99747,   49873,4.35545e-03,2.71858e-03,9.42652e-03,       Z,3.61836e+00
      30,       0,   99747,   49873,4.35556e-03,1.80048e-03,9.04787e-03,       R,4.03706e+00
      30,       1,  177753,   88876,3.26724e-03,3.39673e-03,7.83911e-03,       Z,4.63412e+00
      31,       0,  177753,   88876,3.26744e-03,2.16464e-03,7.13623e-03,       Z,5.50647e+00
      32,       0,  177753,   88876,3.26758e-03,1.43304e-03,6.81919e-03,       Z,6.23357e+00
      33,       0,  177753,   88876,3.26769e-03,9.73115e-04,6.67216e-03,       R,6.82767e+00
      33,       1,  313711,  156855,2.45004e-03,2.37108e-03,5.76884e-03,       Z,7.79692e+00
      34,       0,  313711,  156855,2.44997e-03,1.52224e-03,5.30117e-03,       Z,9.13022e+00
      35,       0,  313711,  156855,2.44993e-03,1.01155e-03,5.08914e-03,       Z,1.02373e+01
      36,       0,  313711,  156855,2.44990e-03,6.87505e-04,4.99061e-03,       R,1.11400e+01
      36,       1,  556573,  278286,1.84423e-03,1.75314e-03,4.30128e-03,       Z,1.26941e+01
      37,       0,  556573,  278286,1.84433e-03,1.10616e-03,3.96507e-03,       Z,1.47852e+01
      38,       0,  556573,  278286,1.84441e-03,7.27001e-04,3.81756e-03,       R,1.65197e+01
      38,       1,  980655,  490327,1.39764e-03,1.40605e-03,3.32444e-03,       Z,1.88977e+01
      39,       0,  980655,  490327,1.39765e-03,8.99735e-04,3.03948e-03,       Z,2.24431e+01
      40,       0,  980655,  490327,1.39766e-03,5.96780e-04,2.91049e-03,       Z,2.54096e+01
      41,       0,  980655,  490327,1.39766e-03,4.05292e-04,2.85061e-03,       R,2.78427e+01
      41,       1, 1665503,  832751,1.06640e-03,9.90213e-04,2.48002e-03,       Z,3.17260e+01
        }\tableDataFive

        \pgfplotstableread[col sep=comma]{
       k,     ell,    ndof,   nElem,        eta,         mu,        res,    case,   maxGradU
       0,       0,     193,      96,1.00896e-01,5.83038e-01,8.99738e-01,       Z,2.15503e-03
       1,       0,     193,      96,1.00842e-01,2.93805e-01,5.01852e-01,       R,6.04181e-03
       1,       1,     271,     135,8.17297e-02,2.99684e-01,4.88116e-01,       Z,1.00747e-02
       2,       0,     271,     135,8.16617e-02,1.56242e-01,2.98824e-01,       Z,1.89306e-02
       3,       0,     271,     135,8.15955e-02,8.64786e-02,2.18991e-01,       R,2.67471e-02
       3,       1,     427,     213,6.97909e-02,9.62580e-02,2.06544e-01,       Z,3.54307e-02
       4,       0,     427,     213,6.97253e-02,5.77212e-02,1.68723e-01,       R,5.01653e-02
       4,       1,     633,     316,5.62695e-02,7.09020e-02,1.54874e-01,       Z,6.32253e-02
       5,       0,     633,     316,5.62384e-02,4.53859e-02,1.32507e-01,       Z,8.79335e-02
       6,       0,     633,     316,5.62149e-02,3.06981e-02,1.22503e-01,       R,1.08203e-01
       6,       1,     937,     468,4.56701e-02,4.49080e-02,1.10924e-01,       Z,1.09260e-01
       7,       0,     937,     468,4.56572e-02,2.97143e-02,1.00952e-01,       Z,1.26100e-01
       8,       0,     937,     468,4.56485e-02,2.03060e-02,9.62142e-02,       R,1.39063e-01
       8,       1,    1309,     654,3.77501e-02,3.27270e-02,8.74973e-02,       Z,1.65001e-01
       9,       0,    1309,     654,3.77440e-02,2.15617e-02,8.13071e-02,       Z,1.96791e-01
      10,       0,    1309,     654,3.77398e-02,1.46174e-02,7.83987e-02,       R,2.21918e-01
      10,       1,    2141,    1070,3.03214e-02,2.68063e-02,6.99016e-02,       Z,2.61717e-01
      11,       0,    2141,    1070,3.03293e-02,1.71424e-02,6.49729e-02,       Z,3.15697e-01
      12,       0,    2141,    1070,3.03347e-02,1.13973e-02,6.27711e-02,       R,3.59023e-01
      12,       1,    3101,    1550,2.45310e-02,2.11736e-02,5.62089e-02,       Z,4.20140e-01
      13,       0,    3101,    1550,2.45329e-02,1.36017e-02,5.23790e-02,       Z,5.07914e-01
      14,       0,    3101,    1550,2.45344e-02,9.04632e-03,5.06600e-02,       R,5.79139e-01
      14,       1,    4611,    2305,2.00855e-02,1.67434e-02,4.57119e-02,       Z,6.66814e-01
      15,       0,    4611,    2305,2.00832e-02,1.08548e-02,4.27472e-02,       Z,8.01035e-01
      16,       0,    4611,    2305,2.00817e-02,7.25705e-03,4.13983e-02,       R,9.10722e-01
      16,       1,    7209,    3604,1.62690e-02,1.38296e-02,3.70372e-02,       Z,1.06319e+00
      17,       0,    7209,    3604,1.62696e-02,8.81547e-03,3.45919e-02,       Z,1.28711e+00
      18,       0,    7209,    3604,1.62701e-02,5.83506e-03,3.35101e-02,       R,1.47191e+00
      18,       1,   10405,    5202,1.34099e-02,1.09060e-02,3.02270e-02,       Z,1.70595e+00
      19,       0,   10405,    5202,1.34113e-02,6.95138e-03,2.83665e-02,       Z,2.06611e+00
      20,       0,   10405,    5202,1.34124e-02,4.59347e-03,2.75480e-02,       R,2.36538e+00
      20,       1,   15749,    7874,1.09239e-02,9.03663e-03,2.47425e-02,       Z,2.38567e+00
      21,       0,   15749,    7874,1.09231e-02,5.79918e-03,2.31691e-02,       Z,2.64613e+00
      22,       0,   15749,    7874,1.09225e-02,3.85009e-03,2.24666e-02,       R,2.85506e+00
      22,       1,   23243,   11621,8.96472e-03,7.33204e-03,2.02599e-02,       Z,3.26478e+00
      23,       0,   23243,   11621,8.96461e-03,4.71221e-03,1.89964e-02,       R,3.78608e+00
      23,       1,   34815,   17407,7.36648e-03,6.95012e-03,1.72182e-02,       Z,4.37836e+00
      24,       0,   34815,   17407,7.36743e-03,4.45108e-03,1.58784e-02,       Z,5.35046e+00
      25,       0,   34815,   17407,7.36813e-03,2.95363e-03,1.52754e-02,       R,6.16870e+00
      25,       1,   49897,   24948,6.13029e-03,5.04309e-03,1.38672e-02,       Z,7.11315e+00
      26,       0,   49897,   24948,6.13083e-03,3.23682e-03,1.29947e-02,       Z,8.65007e+00
      27,       0,   49897,   24948,6.13123e-03,2.14839e-03,1.26058e-02,       R,9.94367e+00
      27,       1,   71927,   35963,5.07506e-03,4.05603e-03,1.14271e-02,       Z,1.13337e+01
      28,       0,   71927,   35963,5.07490e-03,2.62395e-03,1.07382e-02,       Z,1.36648e+01
      29,       0,   71927,   35963,5.07481e-03,1.75212e-03,1.04271e-02,       R,1.56234e+01
      29,       1,  104369,   52184,4.24766e-03,3.28344e-03,9.47892e-03,       Z,1.79842e+01
      30,       0,  104369,   52184,4.24777e-03,2.10161e-03,8.94193e-03,       Z,2.17464e+01
      31,       0,  104369,   52184,4.24785e-03,1.39413e-03,8.70459e-03,       R,2.49139e+01
      31,       1,  148793,   74396,3.55369e-03,2.71278e-03,7.89503e-03,       Z,2.51079e+01
      32,       0,  148793,   74396,3.55397e-03,1.71787e-03,7.46003e-03,       R,2.78847e+01
      32,       1,  208883,  104441,2.99862e-03,2.56711e-03,6.83621e-03,       Z,3.18377e+01
      33,       0,  208883,  104441,2.99873e-03,1.64025e-03,6.37899e-03,       Z,3.76191e+01
      34,       0,  208883,  104441,2.99882e-03,1.08604e-03,6.17594e-03,       R,4.24448e+01
      34,       1,  292283,  146141,2.52507e-03,1.94844e-03,5.64490e-03,       Z,4.81223e+01
      35,       0,  292283,  146141,2.52500e-03,1.26107e-03,5.32361e-03,       Z,5.71812e+01
      36,       0,  292283,  146141,2.52495e-03,8.42382e-04,5.17882e-03,       R,6.47636e+01
      36,       1,  410411,  205205,2.14108e-03,1.58138e-03,4.73709e-03,       Z,7.43316e+01
      37,       0,  410411,  205205,2.14110e-03,1.01269e-03,4.48825e-03,       R,8.92206e+01
      37,       1,  577419,  288709,1.80705e-03,1.53115e-03,4.10657e-03,       Z,8.99622e+01
      38,       0,  577419,  288709,1.80718e-03,9.74680e-04,3.83738e-03,       Z,1.03147e+02
      39,       0,  577419,  288709,1.80727e-03,6.44270e-04,3.71868e-03,       R,1.14021e+02
      39,       1,  793523,  396761,1.54019e-03,1.14417e-03,3.41008e-03,       Z,1.29910e+02
      40,       0,  793523,  396761,1.54024e-03,7.31283e-04,3.22924e-03,       R,1.52856e+02
      40,       1, 1106953,  553476,1.30207e-03,1.10079e-03,2.96734e-03,       Z,1.74297e+02
        }\tableDataThree

        %
        %

        \addplot+ [marker1, adaptive, forget plot]
        table [col sep=comma, x=cumulativeNdof, y=res] {\tableDataUniform};

        \addplot+ [marker2, adaptive, forget plot]
        table [col sep=comma, x=cumulativeNdof, y=res] {\tableDataNine};

        \addplot+ [marker3, adaptive, forget plot]
        table [col sep=comma, x=cumulativeNdof, y=res] {\tableDataSeven};

        \addplot+ [marker4, adaptive, forget plot]
        table [col sep=comma, x=cumulativeNdof, y=res] {\tableDataFive};

        \addplot+ [marker5, adaptive, forget plot]
        table [col sep=comma, x=cumulativeNdof, y=res] {\tableDataThree};

        \drawslopetriangle[ST1]{0.5}{4e4}{1e-2} 
        \drawswappedslopetriangle[ST2]{0.38}{6e6}{2e-2} 
    \end{loglogaxis}
\end{tikzpicture}

%% file: figures/Fig04_Nonlinear_convergence_damping.tex
\begin{tikzpicture}[>=stealth]
    \begin{loglogaxis}[%
            width            = 5.5cm,%
            xlabel           = {cumulative ndof},%
            ylabel           = {Least-squares functional \(N^k_\ell\)},%
            ymajorgrids      = true,%
            font             = \footnotesize,%
            grid style       = {%
                densely dotted,%
                semithick%
            },%
            legend style     = {%
                legend pos = south west,%
                font = \footnotesize%
            },%
        ]

        \addlegendimage{marker1}
        \addlegendentry{\(\delta = 1.0\phantom{0}\)}
        \addlegendimage{marker2}
        \addlegendentry{\(\delta = 0.5\phantom{0}\)}
        \addlegendimage{marker3}
        \addlegendentry{\(\delta = 0.1\phantom{0}\)}
        \addlegendimage{marker4}
        \addlegendentry{\(\delta = 0.05\)}
        \addlegendimage{marker5}
        \addlegendentry{\(\delta = 0.01\)}

        \pgfplotstableread[col sep=comma]{
       k,     ell,    ndof,   nElem,        eta,         mu,        res,    case,   maxGradU
       0,       0,     193,      96,2.01791e-01,1.16608e+00,2.55303e-01,       Z,8.62011e-03
       1,       0,     193,      96,2.01405e-01,1.51358e-01,2.13723e-01,       R,1.60043e-02
       1,       1,     271,     135,1.63298e-01,1.91851e-01,1.79901e-01,       Z,2.95074e-02
       2,       0,     271,     135,1.63037e-01,7.16045e-02,1.67435e-01,       Z,4.14307e-02
       3,       0,     271,     135,1.62901e-01,2.90517e-02,1.65264e-01,       R,4.69924e-02
       3,       1,     427,     213,1.39329e-01,8.92647e-02,1.42920e-01,       Z,7.06245e-02
       4,       0,     427,     213,1.39178e-01,2.65099e-02,1.40792e-01,       R,8.57288e-02
       4,       1,     633,     316,1.12356e-01,8.63095e-02,1.16536e-01,       Z,1.20423e-01
       5,       0,     633,     316,1.12328e-01,2.82196e-02,1.13647e-01,       Z,1.47185e-01
       6,       0,     633,     316,1.12315e-01,1.15980e-02,1.13149e-01,       R,1.60334e-01
       6,       1,     937,     468,9.12346e-02,6.65240e-02,9.38424e-02,       Z,1.62853e-01
       7,       0,     937,     468,9.12420e-02,2.03783e-02,9.19383e-02,       Z,1.70365e-01
       8,       0,     937,     468,9.12472e-02,7.78354e-03,9.16595e-02,       R,1.73925e-01
       8,       1,    1309,     654,7.54564e-02,5.18938e-02,7.74173e-02,       Z,2.33313e-01
       9,       0,    1309,     654,7.54574e-02,1.63962e-02,7.59310e-02,       Z,2.67113e-01
      10,       0,    1309,     654,7.54574e-02,6.39634e-03,7.57034e-02,       R,2.84215e-01
      10,       1,    2137,    1068,6.08165e-02,4.51230e-02,6.24039e-02,       Z,3.76092e-01
      11,       0,    2137,    1068,6.08435e-02,1.33715e-02,6.12043e-02,       Z,4.36950e-01
      12,       0,    2137,    1068,6.08539e-02,5.40803e-03,6.10093e-02,       R,4.68905e-01
      12,       1,    3089,    1544,4.92090e-02,3.62065e-02,5.05179e-02,       Z,6.00485e-01
      13,       0,    3089,    1544,4.92085e-02,1.11692e-02,4.94582e-02,       Z,6.94138e-01
      14,       0,    3089,    1544,4.92082e-02,4.33988e-03,4.92970e-02,       R,7.44808e-01
      14,       1,    4589,    2294,4.01562e-02,2.87710e-02,4.11967e-02,       Z,9.64786e-01
      15,       0,    4589,    2294,4.01500e-02,9.09263e-03,4.03372e-02,       Z,1.12408e+00
      16,       0,    4589,    2294,4.01479e-02,3.57172e-03,4.02034e-02,       R,1.21258e+00
      16,       1,    7187,    3593,3.25030e-02,2.38362e-02,3.33021e-02,       Z,1.56129e+00
      17,       0,    7187,    3593,3.25053e-02,7.16600e-03,3.26503e-02,       Z,1.82386e+00
      18,       0,    7187,    3593,3.25067e-02,2.87589e-03,3.25452e-02,       R,1.97272e+00
      18,       1,   10419,    5209,2.68242e-02,1.85855e-02,2.74027e-02,       Z,2.52575e+00
      19,       0,   10419,    5209,2.68268e-02,5.54342e-03,2.69252e-02,       Z,2.95377e+00
      20,       0,   10419,    5209,2.68280e-02,2.17301e-03,2.68520e-02,       R,3.20017e+00
      20,       1,   15673,    7836,2.18321e-02,1.57425e-02,2.23725e-02,       Z,3.24298e+00
      21,       0,   15673,    7836,2.18306e-02,4.87122e-03,2.19192e-02,       Z,3.40512e+00
      22,       0,   15673,    7836,2.18303e-02,1.91498e-03,2.18486e-02,       R,3.49461e+00
      22,       1,   23249,   11624,1.78925e-02,1.26526e-02,1.83225e-02,       Z,4.42420e+00
      23,       0,   23249,   11624,1.78942e-02,3.92574e-03,1.79647e-02,       R,5.04005e+00
      23,       1,   34811,   17405,1.47452e-02,1.08718e-02,1.51273e-02,       Z,6.43375e+00
      24,       0,   34811,   17405,1.47492e-02,3.35415e-03,1.48135e-02,       Z,7.70067e+00
      25,       0,   34811,   17405,1.47509e-02,1.34287e-03,1.47634e-02,       R,8.44703e+00
      25,       1,   49807,   24903,1.22872e-02,8.27143e-03,1.25374e-02,       Z,1.06775e+01
      26,       0,   49807,   24903,1.22882e-02,2.48284e-03,1.23277e-02,       Z,1.25320e+01
      27,       0,   49807,   24903,1.22886e-02,9.69309e-04,1.22958e-02,       R,1.36208e+01
      27,       1,   71459,   35729,1.01768e-02,6.95582e-03,1.04136e-02,       Z,1.69306e+01
      28,       0,   71459,   35729,1.01764e-02,2.20749e-03,1.02143e-02,       Z,1.96559e+01
      29,       0,   71459,   35729,1.01763e-02,8.74291e-04,1.01829e-02,       R,2.12546e+01
      29,       1,  103867,   51933,8.51327e-03,5.64327e-03,8.68291e-03,       Z,2.69515e+01
      30,       0,  103867,   51933,8.51369e-03,1.70441e-03,8.54125e-03,       Z,3.14958e+01
      31,       0,  103867,   51933,8.51393e-03,6.78939e-04,8.51880e-03,       R,3.41693e+01
      31,       1,  148269,   74134,7.12601e-03,4.70829e-03,7.25922e-03,       Z,3.45545e+01
      32,       0,  148269,   74134,7.12693e-03,1.37864e-03,7.14859e-03,       R,3.63128e+01
      32,       1,  207617,  103808,6.02299e-03,4.05184e-03,6.15480e-03,       Z,4.56736e+01
      33,       0,  207617,  103808,6.02330e-03,1.26496e-03,6.04426e-03,       Z,5.25239e+01
      34,       0,  207617,  103808,6.02347e-03,4.99684e-04,6.02701e-03,       R,5.65371e+01
      34,       1,  290993,  145496,5.05900e-03,3.30732e-03,5.16659e-03,       Z,6.98422e+01
      35,       0,  290993,  145496,5.05879e-03,1.04959e-03,5.07584e-03,       R,8.05348e+01
      35,       1,  407169,  203584,4.29572e-03,2.87049e-03,4.38946e-03,       Z,1.02460e+02
      36,       0,  407169,  203584,4.29585e-03,9.01526e-04,4.31101e-03,       Z,1.23282e+02
      37,       0,  407169,  203584,4.29592e-03,3.59901e-04,4.29850e-03,       R,1.35684e+02
      37,       1,  572869,  286434,3.62967e-03,2.32595e-03,3.69334e-03,       Z,1.37730e+02
      38,       0,  572869,  286434,3.63005e-03,6.80692e-04,3.64041e-03,       Z,1.46105e+02
      39,       0,  572869,  286434,3.63022e-03,2.71917e-04,3.63200e-03,       R,1.50865e+02
      39,       1,  787563,  393781,3.09451e-03,1.91740e-03,3.14701e-03,       Z,1.87968e+02
      40,       0,  787563,  393781,3.09452e-03,5.72286e-04,3.10273e-03,       R,2.13507e+02
      40,       1, 1097597,  548798,2.61659e-03,1.74844e-03,2.67739e-03,       Z,2.65291e+02
        }\tableDataOne

        \pgfplotstableread[col sep=comma]{
       k,     ell,    ndof,   nElem,        eta,         mu,        res,    case,   maxGradU
       0,       0,     193,      96,1.00896e-01,5.83038e-01,8.99738e-01,       Z,2.15503e-03
       1,       0,     193,      96,1.00842e-01,2.93805e-01,5.01852e-01,       R,6.04181e-03
       1,       1,     271,     135,8.17297e-02,2.99684e-01,4.88116e-01,       Z,1.00747e-02
       2,       0,     271,     135,8.16617e-02,1.56242e-01,2.98824e-01,       Z,1.89306e-02
       3,       0,     271,     135,8.15955e-02,8.64786e-02,2.18991e-01,       R,2.67471e-02
       3,       1,     427,     213,6.97909e-02,9.62580e-02,2.06544e-01,       Z,3.54307e-02
       4,       0,     427,     213,6.97253e-02,5.77212e-02,1.68723e-01,       R,5.01653e-02
       4,       1,     633,     316,5.62695e-02,7.09020e-02,1.54874e-01,       Z,6.32253e-02
       5,       0,     633,     316,5.62384e-02,4.53859e-02,1.32507e-01,       Z,8.79335e-02
       6,       0,     633,     316,5.62149e-02,3.06981e-02,1.22503e-01,       R,1.08203e-01
       6,       1,     937,     468,4.56701e-02,4.49080e-02,1.10924e-01,       Z,1.09260e-01
       7,       0,     937,     468,4.56572e-02,2.97143e-02,1.00952e-01,       Z,1.26100e-01
       8,       0,     937,     468,4.56485e-02,2.03060e-02,9.62142e-02,       R,1.39063e-01
       8,       1,    1309,     654,3.77501e-02,3.27270e-02,8.74973e-02,       Z,1.65001e-01
       9,       0,    1309,     654,3.77440e-02,2.15617e-02,8.13071e-02,       Z,1.96791e-01
      10,       0,    1309,     654,3.77398e-02,1.46174e-02,7.83987e-02,       R,2.21918e-01
      10,       1,    2141,    1070,3.03214e-02,2.68063e-02,6.99016e-02,       Z,2.61717e-01
      11,       0,    2141,    1070,3.03293e-02,1.71424e-02,6.49729e-02,       Z,3.15697e-01
      12,       0,    2141,    1070,3.03347e-02,1.13973e-02,6.27711e-02,       R,3.59023e-01
      12,       1,    3101,    1550,2.45310e-02,2.11736e-02,5.62089e-02,       Z,4.20140e-01
      13,       0,    3101,    1550,2.45329e-02,1.36017e-02,5.23790e-02,       Z,5.07914e-01
      14,       0,    3101,    1550,2.45344e-02,9.04632e-03,5.06600e-02,       R,5.79139e-01
      14,       1,    4611,    2305,2.00855e-02,1.67434e-02,4.57119e-02,       Z,6.66814e-01
      15,       0,    4611,    2305,2.00832e-02,1.08548e-02,4.27472e-02,       Z,8.01035e-01
      16,       0,    4611,    2305,2.00817e-02,7.25705e-03,4.13983e-02,       R,9.10722e-01
      16,       1,    7209,    3604,1.62690e-02,1.38296e-02,3.70372e-02,       Z,1.06319e+00
      17,       0,    7209,    3604,1.62696e-02,8.81547e-03,3.45919e-02,       Z,1.28711e+00
      18,       0,    7209,    3604,1.62701e-02,5.83506e-03,3.35101e-02,       R,1.47191e+00
      18,       1,   10405,    5202,1.34099e-02,1.09060e-02,3.02270e-02,       Z,1.70595e+00
      19,       0,   10405,    5202,1.34113e-02,6.95138e-03,2.83665e-02,       Z,2.06611e+00
      20,       0,   10405,    5202,1.34124e-02,4.59347e-03,2.75480e-02,       R,2.36538e+00
      20,       1,   15749,    7874,1.09239e-02,9.03663e-03,2.47425e-02,       Z,2.38567e+00
      21,       0,   15749,    7874,1.09231e-02,5.79918e-03,2.31691e-02,       Z,2.64613e+00
      22,       0,   15749,    7874,1.09225e-02,3.85009e-03,2.24666e-02,       R,2.85506e+00
      22,       1,   23243,   11621,8.96472e-03,7.33204e-03,2.02599e-02,       Z,3.26478e+00
      23,       0,   23243,   11621,8.96461e-03,4.71221e-03,1.89964e-02,       R,3.78608e+00
      23,       1,   34815,   17407,7.36648e-03,6.95012e-03,1.72182e-02,       Z,4.37836e+00
      24,       0,   34815,   17407,7.36743e-03,4.45108e-03,1.58784e-02,       Z,5.35046e+00
      25,       0,   34815,   17407,7.36813e-03,2.95363e-03,1.52754e-02,       R,6.16870e+00
      25,       1,   49897,   24948,6.13029e-03,5.04309e-03,1.38672e-02,       Z,7.11315e+00
      26,       0,   49897,   24948,6.13083e-03,3.23682e-03,1.29947e-02,       Z,8.65007e+00
      27,       0,   49897,   24948,6.13123e-03,2.14839e-03,1.26058e-02,       R,9.94367e+00
      27,       1,   71927,   35963,5.07506e-03,4.05603e-03,1.14271e-02,       Z,1.13337e+01
      28,       0,   71927,   35963,5.07490e-03,2.62395e-03,1.07382e-02,       Z,1.36648e+01
      29,       0,   71927,   35963,5.07481e-03,1.75212e-03,1.04271e-02,       R,1.56234e+01
      29,       1,  104369,   52184,4.24766e-03,3.28344e-03,9.47892e-03,       Z,1.79842e+01
      30,       0,  104369,   52184,4.24777e-03,2.10161e-03,8.94193e-03,       Z,2.17464e+01
      31,       0,  104369,   52184,4.24785e-03,1.39413e-03,8.70459e-03,       R,2.49139e+01
      31,       1,  148793,   74396,3.55369e-03,2.71278e-03,7.89503e-03,       Z,2.51079e+01
      32,       0,  148793,   74396,3.55397e-03,1.71787e-03,7.46003e-03,       R,2.78847e+01
      32,       1,  208883,  104441,2.99862e-03,2.56711e-03,6.83621e-03,       Z,3.18377e+01
      33,       0,  208883,  104441,2.99873e-03,1.64025e-03,6.37899e-03,       Z,3.76191e+01
      34,       0,  208883,  104441,2.99882e-03,1.08604e-03,6.17594e-03,       R,4.24448e+01
      34,       1,  292283,  146141,2.52507e-03,1.94844e-03,5.64490e-03,       Z,4.81223e+01
      35,       0,  292283,  146141,2.52500e-03,1.26107e-03,5.32361e-03,       Z,5.71812e+01
      36,       0,  292283,  146141,2.52495e-03,8.42382e-04,5.17882e-03,       R,6.47636e+01
      36,       1,  410411,  205205,2.14108e-03,1.58138e-03,4.73709e-03,       Z,7.43316e+01
      37,       0,  410411,  205205,2.14110e-03,1.01269e-03,4.48825e-03,       R,8.92206e+01
      37,       1,  577419,  288709,1.80705e-03,1.53115e-03,4.10657e-03,       Z,8.99622e+01
      38,       0,  577419,  288709,1.80718e-03,9.74680e-04,3.83738e-03,       Z,1.03147e+02
      39,       0,  577419,  288709,1.80727e-03,6.44270e-04,3.71868e-03,       R,1.14021e+02
      39,       1,  793523,  396761,1.54019e-03,1.14417e-03,3.41008e-03,       Z,1.29910e+02
      40,       0,  793523,  396761,1.54024e-03,7.31283e-04,3.22924e-03,       R,1.52856e+02
      40,       1, 1106953,  553476,1.30207e-03,1.10079e-03,2.96734e-03,       Z,1.74297e+02
        }\tableDataOhFive

        \pgfplotstableread[col sep=comma]{
       k,     ell,    ndof,   nElem,        eta,         mu,        res,    case,   maxGradU
       0,       0,     193,      96,2.01791e-02,1.16608e-01,1.56407e+00,       Z,8.62011e-05
       1,       0,     193,      96,2.01787e-02,1.04957e-01,1.41328e+00,       R,3.22309e-04
       1,       1,     271,     135,1.63518e-02,1.05620e-01,1.41144e+00,       Z,5.05323e-04
       2,       0,     271,     135,1.63509e-02,9.50898e-02,1.27439e+00,       R,1.21357e-03
       2,       1,     427,     213,1.39869e-02,9.54661e-02,1.27340e+00,       Z,1.58483e-03
       3,       0,     427,     213,1.39853e-02,8.59721e-02,1.14944e+00,       Z,3.14044e-03
       4,       0,     427,     213,1.39833e-02,7.74454e-02,1.03816e+00,       R,5.07462e-03
       4,       1,     641,     320,1.11505e-02,7.79037e-02,1.03710e+00,       Z,5.90044e-03
       5,       0,     641,     320,1.11490e-02,7.02076e-02,9.36215e-01,       Z,9.27065e-03
       6,       0,     641,     320,1.11474e-02,6.32975e-02,8.45657e-01,       R,1.31199e-02
       6,       1,     977,     488,9.08072e-03,6.36269e-02,8.44965e-01,       Z,1.46976e-02
       7,       0,     977,     488,9.07952e-03,5.73971e-02,7.63037e-01,       Z,2.09227e-02
       8,       0,     977,     488,9.07829e-03,5.18047e-02,6.89498e-01,       R,2.77673e-02
       8,       1,    1389,     694,7.44356e-03,5.20647e-02,6.89001e-01,       Z,3.06736e-02
       9,       0,    1389,     694,7.44234e-03,4.70255e-02,6.22533e-01,       Z,4.15314e-02
      10,       0,    1389,     694,7.44112e-03,4.25023e-02,5.62875e-01,       R,5.32422e-02
      10,       1,    2183,    1091,6.00558e-03,4.27288e-02,5.62468e-01,       Z,5.35594e-02
      11,       0,    2183,    1091,6.00465e-03,3.86516e-02,5.08543e-01,       Z,6.62749e-02
      12,       0,    2183,    1091,6.00375e-03,3.49919e-02,4.60142e-01,       R,7.94878e-02
      12,       1,    3205,    1602,4.82654e-03,3.51736e-02,4.59843e-01,       Z,8.53811e-02
      13,       0,    3205,    1602,4.82592e-03,3.18761e-02,4.16119e-01,       Z,1.05598e-01
      14,       0,    3205,    1602,4.82532e-03,2.89158e-02,3.76870e-01,       R,1.26630e-01
      14,       1,    4861,    2430,3.91680e-03,2.90528e-02,3.76658e-01,       Z,1.35979e-01
      15,       0,    4861,    2430,3.91619e-03,2.63856e-02,3.41225e-01,       Z,1.68165e-01
      16,       0,    4861,    2430,3.91559e-03,2.39903e-02,3.09414e-01,       R,2.01693e-01
      16,       1,    7747,    3873,3.14967e-03,2.41028e-02,3.09243e-01,       Z,2.16545e-01
      17,       0,    7747,    3873,3.14927e-03,2.19428e-02,2.80520e-01,       Z,2.67869e-01
      18,       0,    7747,    3873,3.14888e-03,2.00019e-02,2.54725e-01,       R,3.21413e-01
      18,       1,   11017,    5508,2.60220e-03,2.00803e-02,2.54611e-01,       Z,3.42851e-01
      19,       0,   11017,    5508,2.60196e-03,1.83297e-02,2.31336e-01,       Z,4.22307e-01
      20,       0,   11017,    5508,2.60173e-03,1.67550e-02,2.10427e-01,       R,5.05159e-01
      20,       1,   16537,    8268,2.12684e-03,1.68219e-02,2.10334e-01,       Z,5.41014e-01
      21,       0,   16537,    8268,2.12659e-03,1.53998e-02,1.91457e-01,       Z,6.67015e-01
      22,       0,   16537,    8268,2.12636e-03,1.41190e-02,1.74490e-01,       R,7.98617e-01
      22,       1,   24715,   12357,1.75004e-03,1.41705e-02,1.74418e-01,       Z,8.54897e-01
      23,       0,   24715,   12357,1.74985e-03,1.30119e-02,1.59093e-01,       Z,1.05460e+00
      24,       0,   24715,   12357,1.74966e-03,1.19664e-02,1.45308e-01,       R,1.26348e+00
      24,       1,   36311,   18155,1.43946e-03,1.20076e-02,1.45250e-01,       Z,1.26873e+00
      25,       0,   36311,   18155,1.43935e-03,1.10598e-02,1.32788e-01,       Z,1.49029e+00
      26,       0,   36311,   18155,1.43925e-03,1.02025e-02,1.21566e-01,       R,1.71719e+00
      26,       1,   52341,   26170,1.20073e-03,1.02333e-02,1.21524e-01,       Z,1.81983e+00
      27,       0,   52341,   26170,1.20067e-03,9.45439e-03,1.11372e-01,       R,2.16342e+00
      27,       1,   75275,   37637,9.95314e-04,9.47820e-03,1.11340e-01,       Z,2.29649e+00
      28,       0,   75275,   37637,9.95216e-04,8.76989e-03,1.02153e-01,       Z,2.80098e+00
      29,       0,   75275,   37637,9.95123e-04,8.12643e-03,9.38601e-02,       R,3.32715e+00
      29,       1,  108139,   54069,8.35173e-04,8.14443e-03,9.38352e-02,       Z,3.55083e+00
      30,       0,  108139,   54069,8.35108e-04,7.55744e-03,8.63197e-02,       Z,4.34666e+00
      31,       0,  108139,   54069,8.35046e-04,7.02241e-03,7.95229e-02,       R,5.17859e+00
      31,       1,  153133,   76566,7.00348e-04,7.03712e-03,7.95023e-02,       Z,5.52283e+00
      32,       0,  153133,   76566,7.00319e-04,6.54731e-03,7.33295e-02,       R,6.77317e+00
      32,       1,  216925,  108462,5.90695e-04,6.55811e-03,7.33145e-02,       Z,6.79834e+00
      33,       0,  216925,  108462,5.90670e-04,6.10914e-03,6.77049e-02,       Z,8.13585e+00
      34,       0,  216925,  108462,5.90647e-04,5.69770e-03,6.26134e-02,       R,9.52109e+00
      34,       1,  301495,  150747,4.98491e-04,5.70650e-03,6.26013e-02,       Z,1.00692e+01
      35,       0,  301495,  150747,4.98455e-04,5.32800e-03,5.79627e-02,       Z,1.21024e+01
      36,       0,  301495,  150747,4.98422e-04,4.97984e-03,5.37410e-02,       R,1.42133e+01
      36,       1,  421393,  210696,4.23331e-04,4.98679e-03,5.37311e-02,       Z,1.51281e+01
      37,       0,  421393,  210696,4.23307e-04,4.66532e-03,4.98740e-02,       R,1.83255e+01
      37,       1,  587603,  293801,3.57978e-04,4.67079e-03,4.98660e-02,       Z,1.83963e+01
      38,       0,  587603,  293801,3.57966e-04,4.37350e-03,4.63374e-02,       Z,2.18084e+01
      39,       0,  587603,  293801,3.57956e-04,4.09859e-03,4.31112e-02,       R,2.53308e+01
      39,       1,  812197,  406098,3.04838e-04,4.10288e-03,4.31050e-02,       Z,2.68529e+01
      40,       0,  812197,  406098,3.04831e-04,3.84781e-03,4.01446e-02,       R,3.21507e+01
      40,       1, 1130281,  565140,2.58243e-04,3.85122e-03,4.01399e-02,       Z,3.41234e+01
        }\tableDataOhOne

        \pgfplotstableread[col sep=comma]{
       k,     ell,    ndof,   nElem,        eta,         mu,        res,    case,   maxGradU
       0,       0,     193,      96,1.00896e-02,5.83038e-02,1.64801e+00,       Z,2.15503e-05
       1,       0,     193,      96,1.00895e-02,5.53898e-02,1.56827e+00,       R,8.33586e-05
       1,       1,     271,     135,8.17602e-03,5.57044e-02,1.56738e+00,       Z,1.29800e-04
       2,       0,     271,     135,8.17589e-03,5.29229e-02,1.49084e+00,       R,3.21658e-04
       2,       1,     427,     213,6.99385e-03,5.30921e-02,1.49038e+00,       Z,4.17045e-04
       3,       0,     427,     213,6.99363e-03,5.04438e-02,1.41729e+00,       Z,8.50637e-04
       4,       0,     427,     213,6.99332e-03,4.79302e-02,1.34793e+00,       R,1.41585e-03
       4,       1,     641,     320,5.57624e-03,4.81157e-02,1.34745e+00,       Z,1.63383e-03
       5,       0,     641,     320,5.57601e-03,4.57215e-02,1.28114e+00,       Z,2.63337e-03
       6,       0,     641,     320,5.57573e-03,4.34493e-02,1.21821e+00,       R,3.82719e-03
       6,       1,     977,     488,4.54176e-03,4.35695e-02,1.21792e+00,       Z,4.25412e-03
       7,       0,     977,     488,4.54154e-03,4.14078e-02,1.15790e+00,       Z,6.19205e-03
       8,       0,     977,     488,4.54130e-03,3.93564e-02,1.10093e+00,       R,8.41455e-03
       8,       1,    1381,     690,3.73504e-03,3.94411e-02,1.10073e+00,       Z,9.20589e-03
       9,       0,    1381,     690,3.73480e-03,3.74907e-02,1.04646e+00,       Z,1.26931e-02
      10,       0,    1381,     690,3.73454e-03,3.56399e-02,9.94945e-01,       R,1.66050e-02
      10,       1,    2189,    1094,2.99783e-03,3.57094e-02,9.94791e-01,       Z,1.66950e-02
      11,       0,    2189,    1094,2.99761e-03,3.39503e-02,9.45740e-01,       Z,2.10957e-02
      12,       0,    2189,    1094,2.99739e-03,3.22810e-02,8.99177e-01,       R,2.58525e-02
      12,       1,    3195,    1597,2.41440e-03,3.23298e-02,8.99075e-01,       Z,2.75896e-02
      13,       0,    3195,    1597,2.41421e-03,3.07439e-02,8.54774e-01,       Z,3.47162e-02
      14,       0,    3195,    1597,2.41401e-03,2.92391e-02,8.12719e-01,       R,4.24089e-02
      14,       1,    4879,    2439,1.96332e-03,2.92729e-02,8.12652e-01,       Z,4.51411e-02
      15,       0,    4879,    2439,1.96316e-03,2.78438e-02,7.72664e-01,       Z,5.65589e-02
      16,       0,    4879,    2439,1.96301e-03,2.64878e-02,7.34702e-01,       R,6.88667e-02
      16,       1,    7637,    3818,1.59194e-03,2.65127e-02,7.34654e-01,       Z,7.33294e-02
      17,       0,    7637,    3818,1.59179e-03,2.52252e-02,6.98571e-01,       Z,9.16837e-02
      18,       0,    7637,    3818,1.59164e-03,2.40036e-02,6.64315e-01,       R,1.11461e-01
      18,       1,   10853,    5426,1.31303e-03,2.40205e-02,6.64284e-01,       Z,1.17881e-01
      19,       0,   10853,    5426,1.31293e-03,2.28608e-02,6.31734e-01,       Z,1.46466e-01
      20,       0,   10853,    5426,1.31283e-03,2.17605e-02,6.00834e-01,       R,1.77202e-01
      20,       1,   16455,    8227,1.06704e-03,2.17739e-02,6.00810e-01,       Z,1.88238e-01
      21,       0,   16455,    8227,1.06694e-03,2.07295e-02,5.71452e-01,       Z,2.33835e-01
      22,       0,   16455,    8227,1.06684e-03,1.97386e-02,5.43582e-01,       R,2.82884e-01
      22,       1,   24873,   12436,8.74459e-04,1.97480e-02,5.43566e-01,       Z,3.00185e-01
      23,       0,   24873,   12436,8.74386e-04,1.88075e-02,5.17091e-01,       Z,3.72616e-01
      24,       0,   24873,   12436,8.74315e-04,1.79152e-02,4.91958e-01,       R,4.50549e-01
      24,       1,   36239,   18119,7.22245e-04,1.79219e-02,4.91947e-01,       Z,4.52129e-01
      25,       0,   36239,   18119,7.22181e-04,1.70750e-02,4.68075e-01,       Z,5.36839e-01
      26,       0,   36239,   18119,7.22119e-04,1.62715e-02,4.45412e-01,       R,6.26230e-01
      26,       1,   52157,   26078,6.01692e-04,1.62764e-02,4.45405e-01,       Z,6.58259e-01
      27,       0,   52157,   26078,6.01651e-04,1.55138e-02,4.23882e-01,       R,7.88286e-01
      27,       1,   75027,   37513,4.99556e-04,1.55174e-02,4.23876e-01,       Z,8.34139e-01
      28,       0,   75027,   37513,4.99512e-04,1.47936e-02,4.03438e-01,       Z,1.02530e+00
      29,       0,   75027,   37513,4.99469e-04,1.41069e-02,3.84033e-01,       R,1.23032e+00
      29,       1,  107991,   53995,4.19094e-04,1.41095e-02,3.84029e-01,       Z,1.29510e+00
      30,       0,  107991,   53995,4.19060e-04,1.34578e-02,3.65603e-01,       Z,1.58807e+00
      31,       0,  107991,   53995,4.19026e-04,1.28393e-02,3.48109e-01,       R,1.90211e+00
      31,       1,  152657,   76328,3.51207e-04,1.28413e-02,3.48106e-01,       Z,2.01307e+00
      32,       0,  152657,   76328,3.51183e-04,1.22543e-02,3.31494e-01,       R,2.47548e+00
      32,       1,  216225,  108112,2.96244e-04,1.22558e-02,3.31492e-01,       Z,2.48331e+00
      33,       0,  216225,  108112,2.96225e-04,1.16987e-02,3.15717e-01,       Z,2.98890e+00
      34,       0,  216225,  108112,2.96206e-04,1.11699e-02,3.00740e-01,       R,3.52660e+00
      34,       1,  300325,  150162,2.50144e-04,1.11710e-02,3.00739e-01,       Z,3.71557e+00
      35,       0,  300325,  150162,2.50124e-04,1.06691e-02,2.86517e-01,       Z,4.49983e+00
      36,       0,  300325,  150162,2.50104e-04,1.01926e-02,2.73013e-01,       R,5.33652e+00
      36,       1,  420147,  210073,2.12338e-04,1.01935e-02,2.73012e-01,       Z,5.60154e+00
      37,       0,  420147,  210073,2.12323e-04,9.74113e-03,2.60189e-01,       R,6.79165e+00
      37,       1,  584281,  292140,1.79794e-04,9.74178e-03,2.60188e-01,       Z,6.81353e+00
      38,       0,  584281,  292140,1.79785e-04,9.31231e-03,2.48010e-01,       Z,8.10947e+00
      39,       0,  584281,  292140,1.79772e-04,8.90452e-03,2.36447e-01,       R,9.48173e+00
      39,       1,  809405,  404702,1.52820e-04,8.90502e-03,2.36446e-01,       Z,9.97210e+00
      40,       0,  809405,  404702,1.52812e-04,8.51777e-03,2.25465e-01,       R,1.19740e+01
      40,       1, 1127939,  563969,1.29402e-04,8.51815e-03,2.25464e-01,       Z,1.25926e+01
        }\tableDataOhOhFive

        \pgfplotstableread[col sep=comma]{
       k,     ell,    ndof,   nElem,        eta,         mu,        res,    case,   maxGradU
       0,       0,     193,      96,2.01791e-03,1.16608e-02,1.71524e+00,       Z,8.62011e-07
       1,       0,     193,      96,2.01791e-03,1.15442e-02,1.69859e+00,       R,3.42511e-06
       1,       1,     271,     135,1.63521e-03,1.16046e-02,1.69842e+00,       Z,5.30486e-06
       2,       0,     271,     135,1.63521e-03,1.14885e-02,1.68177e+00,       R,1.34817e-05
       2,       1,     427,     213,1.39880e-03,1.15197e-02,1.68168e+00,       Z,1.73811e-05
       3,       0,     427,     213,1.39880e-03,1.14046e-02,1.66512e+00,       Z,3.63002e-05
       4,       0,     427,     213,1.39879e-03,1.12906e-02,1.64872e+00,       R,6.19250e-05
       4,       1,     641,     320,1.11533e-03,1.13221e-02,1.64863e+00,       Z,7.10320e-05
       5,       0,     641,     320,1.11532e-03,1.12090e-02,1.63231e+00,       Z,1.16994e-04
       6,       0,     641,     320,1.11532e-03,1.10970e-02,1.61615e+00,       R,1.73961e-04
       6,       1,     977,     488,9.08473e-04,1.11158e-02,1.61610e+00,       Z,1.92167e-04
       7,       0,     977,     488,9.08471e-04,1.10048e-02,1.60005e+00,       Z,2.85239e-04
       8,       0,     977,     488,9.08468e-04,1.08949e-02,1.58417e+00,       R,3.95852e-04
       8,       1,    1381,     690,7.47200e-04,1.09072e-02,1.58414e+00,       Z,4.30244e-04
       9,       0,    1381,     690,7.47197e-04,1.07983e-02,1.56838e+00,       Z,6.03718e-04
      10,       0,    1381,     690,7.47193e-04,1.06905e-02,1.55278e+00,       R,8.05092e-04
      10,       1,    2191,    1095,6.00029e-04,1.06998e-02,1.55276e+00,       Z,8.09082e-04
      11,       0,    2191,    1095,6.00026e-04,1.05930e-02,1.53729e+00,       Z,1.04282e-03
      12,       0,    2191,    1095,6.00022e-04,1.04873e-02,1.52199e+00,       R,1.30436e-03
      12,       1,    3197,    1598,4.83420e-04,1.04933e-02,1.52197e+00,       Z,1.38299e-03
      13,       0,    3197,    1598,4.83417e-04,1.03886e-02,1.50680e+00,       Z,1.76842e-03
      14,       0,    3197,    1598,4.83414e-04,1.02850e-02,1.49178e+00,       R,2.19819e-03
      14,       1,    4849,    2424,3.93179e-04,1.02889e-02,1.49177e+00,       Z,2.32416e-03
      15,       0,    4849,    2424,3.93177e-04,1.01863e-02,1.47690e+00,       Z,2.95138e-03
      16,       0,    4849,    2424,3.93174e-04,1.00847e-02,1.46217e+00,       R,3.64865e-03
      16,       1,    7719,    3859,3.16335e-04,1.00874e-02,1.46217e+00,       Z,3.85694e-03
      17,       0,    7719,    3859,3.16333e-04,9.98691e-03,1.44758e+00,       Z,4.87648e-03
      18,       0,    7719,    3859,3.16330e-04,9.88740e-03,1.43315e+00,       R,6.00773e-03
      18,       1,   10943,    5471,2.61878e-04,9.88899e-03,1.43314e+00,       Z,6.31270e-03
      19,       0,   10943,    5471,2.61876e-04,9.79048e-03,1.41885e+00,       Z,7.92309e-03
      20,       0,   10943,    5471,2.61873e-04,9.69297e-03,1.40470e+00,       R,9.70424e-03
      20,       1,   16579,    8289,2.13176e-04,9.69417e-03,1.40469e+00,       Z,1.02348e-02
      21,       0,   16579,    8289,2.13174e-04,9.59764e-03,1.39068e+00,       Z,1.28246e-02
      22,       0,   16579,    8289,2.13172e-04,9.50210e-03,1.37681e+00,       R,1.56869e-02
      22,       1,   24895,   12447,1.75140e-04,9.50288e-03,1.37681e+00,       Z,1.65269e-02
      23,       0,   24895,   12447,1.75138e-04,9.40830e-03,1.36308e+00,       Z,2.06703e-02
      24,       0,   24895,   12447,1.75136e-04,9.31469e-03,1.34949e+00,       R,2.52462e-02
      24,       1,   36421,   18210,1.44430e-04,9.31522e-03,1.34948e+00,       Z,2.65930e-02
      25,       0,   36421,   18210,1.44428e-04,9.22256e-03,1.33603e+00,       Z,3.32167e-02
      26,       0,   36421,   18210,1.44426e-04,9.13084e-03,1.32271e+00,       R,4.05281e-02
      26,       1,   52453,   26226,1.20305e-04,9.13119e-03,1.32270e+00,       Z,4.06349e-02
      27,       0,   52453,   26226,1.20303e-04,9.04041e-03,1.30952e+00,       R,4.87472e-02
      27,       1,   75677,   37838,9.96930e-05,9.04066e-03,1.30952e+00,       Z,5.11147e-02
      28,       0,   75677,   37838,9.96917e-05,8.95080e-03,1.29646e+00,       Z,6.27195e-02
      29,       0,   75677,   37838,9.96903e-05,8.86185e-03,1.28354e+00,       R,7.54283e-02
      29,       1,  108613,   54306,8.37518e-05,8.86202e-03,1.28354e+00,       Z,7.91380e-02
      30,       0,  108613,   54306,8.37506e-05,8.77398e-03,1.27074e+00,       Z,9.73450e-02
      31,       0,  108613,   54306,8.37494e-05,8.68684e-03,1.25808e+00,       R,1.17308e-01
      31,       1,  153791,   76895,7.01301e-05,8.68696e-03,1.25808e+00,       Z,1.22652e-01
      32,       0,  153791,   76895,7.01291e-05,8.60070e-03,1.24554e+00,       R,1.50700e-01
      32,       1,  217791,  108895,5.91580e-05,8.60079e-03,1.24554e+00,       Z,1.58974e-01
      33,       0,  217791,  108895,5.91572e-05,8.51541e-03,1.23313e+00,       Z,1.99738e-01
      34,       0,  217791,  108895,5.91561e-05,8.43091e-03,1.22085e+00,       R,2.44865e-01
      34,       1,  302883,  151441,4.99424e-05,8.43097e-03,1.22085e+00,       Z,2.45593e-01
      35,       0,  302883,  151441,4.99416e-05,8.34733e-03,1.20869e+00,       Z,2.95869e-01
      36,       0,  302883,  151441,4.99404e-05,8.26454e-03,1.19665e+00,       R,3.50486e-01
      36,       1,  424553,  212276,4.23708e-05,8.26458e-03,1.19665e+00,       Z,3.66279e-01
      37,       0,  424553,  212276,4.23701e-05,8.18264e-03,1.18473e+00,       R,4.43718e-01
      37,       1,  590551,  295275,3.58552e-05,8.18267e-03,1.18473e+00,       Z,4.43679e-01
      38,       0,  590551,  295275,3.58541e-05,8.10157e-03,1.17293e+00,       Z,5.27977e-01
      39,       0,  590551,  295275,3.58534e-05,8.02129e-03,1.16126e+00,       R,6.19054e-01
      39,       1,  816693,  408346,3.04937e-05,8.02131e-03,1.16126e+00,       Z,6.45447e-01
      40,       0,  816693,  408346,3.04923e-05,7.94185e-03,1.14970e+00,       R,7.73927e-01
      40,       1, 1136315,  568157,2.58417e-05,7.94187e-03,1.14970e+00,       Z,8.08112e-01
        }\tableDataOhOhOne

        %
        %

        \addplot+ [marker1, adaptive, forget plot]
        table [col sep=comma, x=cumulativeNdof, y=res] {\tableDataOne};

        \addplot+ [marker2, adaptive, forget plot]
        table [col sep=comma, x=cumulativeNdof, y=res] {\tableDataOhFive};

        \addplot+ [marker3, adaptive, forget plot]
        table [col sep=comma, x=cumulativeNdof, y=res] {\tableDataOhOne};

        \addplot+ [marker4, adaptive, forget plot]
        table [col sep=comma, x=cumulativeNdof, y=res] {\tableDataOhOhFive};

        \addplot+ [marker5, adaptive, forget plot]
        table [col sep=comma, x=cumulativeNdof, y=res] {\tableDataOhOhOne};

        \drawslopetriangle[ST1]{0.5}{4e5}{3e-3} 
        \drawslopetriangle[ST2]{0.4}{6e5}{3.2e-2} 
    \end{loglogaxis}
\end{tikzpicture}

%% file: figures/Fig05a_Nonlinear_convergence_weighting.tex
\begin{tikzpicture}[>=stealth]
    \begin{loglogaxis}[%
            width            = 5.5cm,%
            xlabel           = {cumulative ndof},%
            ylabel           = {Least-squares functional \(N^k_\ell\)},%
            ymajorgrids      = true,%
            ymax             = 3e0,%
            font             = \footnotesize,%
            grid style       = {%
                densely dotted,%
                semithick%
            },%
            legend style     = {%
                legend pos = south west,%
                font = \footnotesize%
            },%
        ]

        \addlegendimage{marker1}
        \addlegendentry{emph.\ grad.}
        \addlegendimage{marker2}
        \addlegendentry{balanced w.}
        \addlegendimage{marker3}
        \addlegendentry{downsc. flux}
        \addlegendimage{marker4}
        \addlegendentry{split w.}

        \pgfplotstableread[col sep=comma]{
       k,     ell,    ndof,   nElem,        eta,         mu,        res,    case,   maxGradU
       0,       0,     193,      96,2.01791e-01,1.16608e+00,2.55303e-01,       Z,8.62011e-03
       1,       0,     193,      96,2.01405e-01,1.51358e-01,2.13723e-01,       R,1.60043e-02
       1,       1,     271,     135,1.63298e-01,1.91851e-01,1.79901e-01,       Z,2.95074e-02
       2,       0,     271,     135,1.63037e-01,7.16045e-02,1.67435e-01,       Z,4.14307e-02
       3,       0,     271,     135,1.62901e-01,2.90517e-02,1.65264e-01,       R,4.69924e-02
       3,       1,     427,     213,1.39329e-01,8.92647e-02,1.42920e-01,       Z,7.06245e-02
       4,       0,     427,     213,1.39178e-01,2.65099e-02,1.40792e-01,       R,8.57288e-02
       4,       1,     633,     316,1.12356e-01,8.63095e-02,1.16536e-01,       Z,1.20423e-01
       5,       0,     633,     316,1.12328e-01,2.82196e-02,1.13647e-01,       Z,1.47185e-01
       6,       0,     633,     316,1.12315e-01,1.15980e-02,1.13149e-01,       R,1.60334e-01
       6,       1,     937,     468,9.12346e-02,6.65240e-02,9.38424e-02,       Z,1.62853e-01
       7,       0,     937,     468,9.12420e-02,2.03783e-02,9.19383e-02,       Z,1.70365e-01
       8,       0,     937,     468,9.12472e-02,7.78354e-03,9.16595e-02,       R,1.73925e-01
       8,       1,    1309,     654,7.54564e-02,5.18938e-02,7.74173e-02,       Z,2.33313e-01
       9,       0,    1309,     654,7.54574e-02,1.63962e-02,7.59310e-02,       Z,2.67113e-01
      10,       0,    1309,     654,7.54574e-02,6.39634e-03,7.57034e-02,       R,2.84215e-01
      10,       1,    2137,    1068,6.08165e-02,4.51230e-02,6.24039e-02,       Z,3.76092e-01
      11,       0,    2137,    1068,6.08435e-02,1.33715e-02,6.12043e-02,       Z,4.36950e-01
      12,       0,    2137,    1068,6.08539e-02,5.40803e-03,6.10093e-02,       R,4.68905e-01
      12,       1,    3089,    1544,4.92090e-02,3.62065e-02,5.05179e-02,       Z,6.00485e-01
      13,       0,    3089,    1544,4.92085e-02,1.11692e-02,4.94582e-02,       Z,6.94138e-01
      14,       0,    3089,    1544,4.92082e-02,4.33988e-03,4.92970e-02,       R,7.44808e-01
      14,       1,    4589,    2294,4.01562e-02,2.87710e-02,4.11967e-02,       Z,9.64786e-01
      15,       0,    4589,    2294,4.01500e-02,9.09263e-03,4.03372e-02,       Z,1.12408e+00
      16,       0,    4589,    2294,4.01479e-02,3.57172e-03,4.02034e-02,       R,1.21258e+00
      16,       1,    7187,    3593,3.25030e-02,2.38362e-02,3.33021e-02,       Z,1.56129e+00
      17,       0,    7187,    3593,3.25053e-02,7.16600e-03,3.26503e-02,       Z,1.82386e+00
      18,       0,    7187,    3593,3.25067e-02,2.87589e-03,3.25452e-02,       R,1.97272e+00
      18,       1,   10419,    5209,2.68242e-02,1.85855e-02,2.74027e-02,       Z,2.52575e+00
      19,       0,   10419,    5209,2.68268e-02,5.54342e-03,2.69252e-02,       Z,2.95377e+00
      20,       0,   10419,    5209,2.68280e-02,2.17301e-03,2.68520e-02,       R,3.20017e+00
      20,       1,   15673,    7836,2.18321e-02,1.57425e-02,2.23725e-02,       Z,3.24298e+00
      21,       0,   15673,    7836,2.18306e-02,4.87122e-03,2.19192e-02,       Z,3.40512e+00
      22,       0,   15673,    7836,2.18303e-02,1.91498e-03,2.18486e-02,       R,3.49461e+00
      22,       1,   23249,   11624,1.78925e-02,1.26526e-02,1.83225e-02,       Z,4.42420e+00
      23,       0,   23249,   11624,1.78942e-02,3.92574e-03,1.79647e-02,       R,5.04005e+00
      23,       1,   34811,   17405,1.47452e-02,1.08718e-02,1.51273e-02,       Z,6.43375e+00
      24,       0,   34811,   17405,1.47492e-02,3.35415e-03,1.48135e-02,       Z,7.70067e+00
      25,       0,   34811,   17405,1.47509e-02,1.34287e-03,1.47634e-02,       R,8.44703e+00
      25,       1,   49807,   24903,1.22872e-02,8.27143e-03,1.25374e-02,       Z,1.06775e+01
      26,       0,   49807,   24903,1.22882e-02,2.48284e-03,1.23277e-02,       Z,1.25320e+01
      27,       0,   49807,   24903,1.22886e-02,9.69309e-04,1.22958e-02,       R,1.36208e+01
      27,       1,   71459,   35729,1.01768e-02,6.95582e-03,1.04136e-02,       Z,1.69306e+01
      28,       0,   71459,   35729,1.01764e-02,2.20749e-03,1.02143e-02,       Z,1.96559e+01
      29,       0,   71459,   35729,1.01763e-02,8.74291e-04,1.01829e-02,       R,2.12546e+01
      29,       1,  103867,   51933,8.51327e-03,5.64327e-03,8.68291e-03,       Z,2.69515e+01
      30,       0,  103867,   51933,8.51369e-03,1.70441e-03,8.54125e-03,       Z,3.14958e+01
      31,       0,  103867,   51933,8.51393e-03,6.78939e-04,8.51880e-03,       R,3.41693e+01
      31,       1,  148269,   74134,7.12601e-03,4.70829e-03,7.25922e-03,       Z,3.45545e+01
      32,       0,  148269,   74134,7.12693e-03,1.37864e-03,7.14859e-03,       R,3.63128e+01
      32,       1,  207617,  103808,6.02299e-03,4.05184e-03,6.15480e-03,       Z,4.56736e+01
      33,       0,  207617,  103808,6.02330e-03,1.26496e-03,6.04426e-03,       Z,5.25239e+01
      34,       0,  207617,  103808,6.02347e-03,4.99684e-04,6.02701e-03,       R,5.65371e+01
      34,       1,  290993,  145496,5.05900e-03,3.30732e-03,5.16659e-03,       Z,6.98422e+01
      35,       0,  290993,  145496,5.05879e-03,1.04959e-03,5.07584e-03,       R,8.05348e+01
      35,       1,  407169,  203584,4.29572e-03,2.87049e-03,4.38946e-03,       Z,1.02460e+02
      36,       0,  407169,  203584,4.29585e-03,9.01526e-04,4.31101e-03,       Z,1.23282e+02
      37,       0,  407169,  203584,4.29592e-03,3.59901e-04,4.29850e-03,       R,1.35684e+02
      37,       1,  572869,  286434,3.62967e-03,2.32595e-03,3.69334e-03,       Z,1.37730e+02
      38,       0,  572869,  286434,3.63005e-03,6.80692e-04,3.64041e-03,       Z,1.46105e+02
      39,       0,  572869,  286434,3.63022e-03,2.71917e-04,3.63200e-03,       R,1.50865e+02
      39,       1,  787563,  393781,3.09451e-03,1.91740e-03,3.14701e-03,       Z,1.87968e+02
      40,       0,  787563,  393781,3.09452e-03,5.72286e-04,3.10273e-03,       R,2.13507e+02
      40,       1, 1097597,  548798,2.61659e-03,1.74844e-03,2.67739e-03,       Z,2.65291e+02
        }\tableDataOne

        \pgfplotstableread[col sep=comma]{
       k,     ell,    ndof,   nElem,        eta,         mu,        res,    case,   maxGradU
       0,       0,     193,      96,9.60646e-02,8.06816e-01,2.55284e-01,       Z,9.01679e-03
       1,       0,     193,      96,2.01459e-01,1.54929e-01,2.09919e-01,       R,2.75277e-02
       1,       1,     271,     135,1.62861e-01,1.95103e-01,1.91307e-01,       Z,6.78911e-02
       2,       0,     271,     135,1.61705e-01,9.73864e-02,1.90339e-01,       Z,4.81053e-02
       3,       0,     271,     135,1.61894e-01,9.52332e-02,1.95929e-01,       R,5.20205e-02
       3,       1,     427,     213,1.38310e-01,1.27080e-01,1.90134e-01,       Z,1.08792e-01
       4,       0,     427,     213,1.38392e-01,1.28396e-01,2.00474e-01,       R,9.69161e-02
       4,       1,     641,     320,1.09456e-01,1.53809e-01,1.92033e-01,       Z,1.78676e-01
       5,       0,     641,     320,1.09937e-01,1.56749e-01,2.06895e-01,       Z,1.70231e-01
       6,       0,     641,     320,1.09855e-01,1.74688e-01,2.25024e-01,       R,1.71260e-01
       6,       1,    1019,     509,8.55644e-02,1.87784e-01,2.19230e-01,       Z,1.76531e-01
       7,       0,    1019,     509,8.58007e-02,2.01544e-01,2.41107e-01,       Z,1.76739e-01
       8,       0,    1019,     509,8.57667e-02,2.25158e-01,2.66721e-01,       Z,1.76486e-01
       9,       0,    1019,     509,8.57760e-02,2.52393e-01,2.95854e-01,       R,1.76595e-01
       9,       1,    1613,     806,6.97576e-02,2.57281e-01,2.93780e-01,       Z,3.13816e-01
      10,       0,    1613,     806,6.98129e-02,2.85307e-01,3.27250e-01,       Z,2.99466e-01
      11,       0,    1613,     806,6.98094e-02,3.19666e-01,3.65210e-01,       R,3.01260e-01
      11,       1,    2357,    1178,5.72907e-02,3.22145e-01,3.64392e-01,       Z,5.14428e-01
      12,       0,    2357,    1178,5.73334e-02,3.59831e-01,4.07464e-01,       R,5.00202e-01
      12,       1,    3413,    1706,4.68728e-02,3.61343e-01,4.06800e-01,       Z,8.30419e-01
      13,       0,    3413,    1706,4.68849e-02,4.04081e-01,4.55436e-01,       Z,8.23118e-01
      14,       0,    3413,    1706,4.68849e-02,4.53008e-01,5.10124e-01,       R,8.23705e-01
      14,       1,    5143,    2571,3.79220e-02,4.53847e-01,5.09738e-01,       Z,1.35254e+00
      15,       0,    5143,    2571,3.79254e-02,5.08322e-01,5.71201e-01,       Z,1.35260e+00
      16,       0,    5143,    2571,3.79264e-02,5.69938e-01,6.40203e-01,       R,1.35330e+00
      16,       1,    8095,    4047,3.06069e-02,5.70378e-01,6.40043e-01,       Z,2.19506e+00
      17,       0,    8095,    4047,3.06314e-02,6.39309e-01,7.17504e-01,       Z,2.21323e+00
      18,       0,    8095,    4047,3.06267e-02,7.16850e-01,8.04400e-01,       R,2.21466e+00
      18,       1,   11791,    5895,2.52295e-02,7.17060e-01,8.04317e-01,       Z,3.42841e+00
      19,       0,   11791,    5895,2.52326e-02,8.03921e-01,9.01795e-01,       Z,3.47959e+00
      20,       0,   11791,    5895,2.52329e-02,9.01442e-01,1.01112e+00,       R,3.48199e+00
      20,       1,   17589,    8794,2.05474e-02,9.01560e-01,1.01106e+00,       Z,5.55291e+00
      21,       0,   17589,    8794,2.05391e-02,1.01085e+00,1.13366e+00,       Z,5.65807e+00
      22,       0,   17589,    8794,2.05420e-02,1.13348e+00,1.27116e+00,       R,5.66404e+00
      22,       1,   26267,   13133,1.69073e-02,1.13354e+00,1.27113e+00,       Z,8.95155e+00
      23,       0,   26267,   13133,1.69138e-02,1.27102e+00,1.42531e+00,       Z,9.14293e+00
      24,       0,   26267,   13133,1.69125e-02,1.42521e+00,1.59821e+00,       R,9.15426e+00
      24,       1,   38249,   19124,1.39836e-02,1.42524e+00,1.59820e+00,       Z,9.41418e+00
      25,       0,   38249,   19124,1.39914e-02,1.59814e+00,1.79208e+00,       Z,9.42253e+00
      26,       0,   38249,   19124,1.39897e-02,1.79202e+00,2.00948e+00,       R,9.42309e+00
      26,       1,   55879,   27939,1.16120e-02,1.79204e+00,2.00947e+00,       Z,1.47847e+01
      27,       0,   55879,   27939,1.16113e-02,2.00944e+00,2.25325e+00,       Z,1.51095e+01
      28,       0,   55879,   27939,1.16117e-02,2.25322e+00,2.52661e+00,       R,1.51282e+01
      28,       1,   79695,   39847,9.65584e-03,2.25323e+00,2.52660e+00,       Z,2.35544e+01
      29,       0,   79695,   39847,9.65551e-03,2.52659e+00,2.83313e+00,       R,2.40957e+01
      29,       1,  115623,   57811,8.07827e-03,2.52659e+00,2.83313e+00,       Z,3.63136e+01
      30,       0,  115623,   57811,8.08148e-03,2.83311e+00,3.17684e+00,       Z,3.71305e+01
      31,       0,  115623,   57811,8.08075e-03,3.17683e+00,3.56225e+00,       R,3.71789e+01
      31,       1,  161421,   80710,6.81258e-03,3.17683e+00,3.56225e+00,       Z,5.85275e+01
      32,       0,  161421,   80710,6.81395e-03,3.56224e+00,3.99442e+00,       Z,5.99063e+01
      33,       0,  161421,   80710,6.81367e-03,3.99442e+00,4.47903e+00,       R,5.99869e+01
      33,       1,  231105,  115552,5.72127e-03,3.99442e+00,4.47903e+00,       Z,6.15426e+01
      34,       0,  231105,  115552,5.72058e-03,4.47902e+00,5.02242e+00,       R,6.16270e+01
      34,       1,  318659,  159329,4.83695e-03,4.47902e+00,5.02242e+00,       Z,9.58121e+01
      35,       0,  318659,  159329,4.83707e-03,5.02242e+00,5.63174e+00,       Z,9.80175e+01
      36,       0,  318659,  159329,4.83707e-03,5.63174e+00,6.31499e+00,       R,9.81457e+01
      36,       1,  448975,  224487,4.09767e-03,5.63174e+00,6.31499e+00,       Z,1.47648e+02
      37,       0,  448975,  224487,4.09891e-03,6.31499e+00,7.08112e+00,       R,1.50833e+02
      37,       1,  620141,  310070,3.47793e-03,6.31499e+00,7.08112e+00,       Z,1.54955e+02
      38,       0,  620141,  310070,3.47810e-03,7.08112e+00,7.94021e+00,       Z,1.55370e+02
      39,       0,  620141,  310070,3.47809e-03,7.94021e+00,8.90351e+00,       R,1.55394e+02
      39,       1,  868029,  434014,2.95422e-03,7.94021e+00,8.90351e+00,       Z,2.42490e+02
      40,       0,  868029,  434014,2.95403e-03,8.90351e+00,9.98369e+00,       Z,2.48087e+02
      41,       0,  868029,  434014,2.95419e-03,9.98369e+00,1.11949e+01,       R,2.48411e+02
      41,       1, 1187487,  593743,2.51468e-03,9.98369e+00,1.11949e+01,       Z,3.73280e+02
        }\tableDataTwo

        \pgfplotstableread[col sep=comma]{
       k,     ell,    ndof,   nElem,        eta,         mu,        res,    case,   maxGradU
       0,       0,     193,      96,4.54996e-02,5.56005e-01,2.56527e-01,       Z,9.21402e-03
       1,       0,     193,      96,2.02972e-01,1.56388e-01,3.29270e-01,       R,6.23118e-02
       1,       1,     271,     135,1.63709e-01,1.97114e-01,4.29268e-01,       Z,2.03996e-01
       2,       0,     271,     135,1.60311e-01,3.96790e-01,1.07123e+00,       Z,1.90154e-02
       3,       0,     271,     135,1.62342e-01,1.05835e+00,3.36292e+00,       R,6.80706e-01
       3,       1,     427,     213,1.38599e-01,1.06172e+00,3.37008e+00,       Z,1.09683e+00
       4,       0,     427,     213,1.36545e-01,3.36723e+00,1.15418e+01,       R,4.34832e-01
       4,       1,     615,     307,1.08537e-01,3.36825e+00,1.15432e+01,       Z,4.33294e-01
       5,       0,     615,     307,1.16765e-01,1.15426e+01,4.02283e+01,       Z,3.56972e+00
       6,       0,     615,     307,1.09789e-01,4.02282e+01,1.40727e+02,       R,2.53758e+00
       6,       1,     989,     494,8.54091e-02,4.02282e+01,1.40727e+02,       Z,3.14846e+00
       7,       0,     989,     494,9.99854e-02,1.40727e+02,4.92535e+02,       R,9.85007e+00
       7,       1,    1249,     624,8.44980e-02,1.40727e+02,4.92535e+02,       Z,1.25577e+01
       8,       0,    1249,     624,8.14751e-02,4.92535e+02,1.72391e+03,       Z,8.42874e+00
       9,       0,    1249,     624,8.85706e-02,1.72391e+03,6.03373e+03,       R,2.56954e+01
       9,       1,    1851,     925,7.39555e-02,1.72391e+03,6.03373e+03,       Z,3.13280e+01
      10,       0,    1851,     925,7.49541e-02,6.03373e+03,2.11181e+04,       R,2.18526e+01
      10,       1,    2493,    1246,6.06636e-02,6.03373e+03,2.11181e+04,       Z,2.18156e+01
      11,       0,    2493,    1246,6.79063e-02,2.11181e+04,7.39135e+04,       R,5.35179e+01
      11,       1,    3519,    1759,5.58090e-02,2.11181e+04,7.39135e+04,       Z,6.36517e+01
        }\tableDataThree

        \pgfplotstableread[col sep=comma]{
       k,     ell,    ndof,   nElem,        eta,         mu,        res,    case,   maxGradU
       0,       0,     193,      96,3.96416e-01,1.62597e+00,2.63724e-01,       Z,7.94203e-03
       1,       0,     193,      96,2.11806e-01,1.54617e-01,2.31688e-01,       Z,1.18378e-02
       2,       0,     193,      96,2.05213e-01,1.06488e-01,2.18439e-01,       Z,1.46085e-02
       3,       0,     193,      96,2.04843e-01,7.44515e-02,2.11775e-01,       Z,1.67011e-02
       4,       0,     193,      96,2.04752e-01,5.21035e-02,2.08434e-01,       Z,1.82525e-02
       5,       0,     193,      96,2.04693e-01,3.65399e-02,2.06763e-01,       Z,1.93848e-02
       6,       0,     193,      96,2.04650e-01,2.56655e-02,2.05925e-01,       Z,2.02029e-02
       7,       0,     193,      96,2.04620e-01,1.80479e-02,2.05504e-01,       R,2.07900e-02
       7,       1,     271,     135,1.66582e-01,1.20189e-01,1.81366e-01,       Z,2.82188e-02
       8,       0,     271,     135,1.64779e-01,7.52204e-02,1.72208e-01,       Z,3.42216e-02
       9,       0,     271,     135,1.64691e-01,4.95071e-02,1.68254e-01,       R,3.88302e-02
       9,       1,     437,     218,1.38503e-01,1.01937e-01,1.51718e-01,       Z,4.95674e-02
      10,       0,     437,     218,1.37734e-01,6.32914e-02,1.43864e-01,       Z,6.18173e-02
      11,       0,     437,     218,1.37669e-01,4.12609e-02,1.40555e-01,       R,7.14215e-02
      11,       1,     653,     326,1.10661e-01,9.17032e-02,1.24456e-01,       Z,8.67691e-02
      12,       0,     653,     326,1.10160e-01,5.77624e-02,1.16632e-01,       Z,1.07081e-01
      13,       0,     653,     326,1.10156e-01,3.80931e-02,1.13229e-01,       R,1.23124e-01
      13,       1,    1069,     534,8.54867e-02,7.92299e-02,9.94177e-02,       Z,1.47598e-01
      14,       0,    1069,     534,8.51726e-02,5.12155e-02,9.18009e-02,       Z,1.80778e-01
      15,       0,    1069,     534,8.51616e-02,3.41827e-02,8.83159e-02,       R,2.07239e-01
      15,       1,    1593,     796,7.01506e-02,5.91596e-02,7.97951e-02,       Z,2.45396e-01
      16,       0,    1593,     796,7.00281e-02,3.82170e-02,7.45291e-02,       Z,2.98846e-01
      17,       0,    1593,     796,7.00243e-02,2.54642e-02,7.21478e-02,       R,3.41873e-01
      17,       1,    2361,    1180,5.67807e-02,4.82471e-02,6.44427e-02,       Z,3.45998e-01
      18,       0,    2361,    1180,5.67247e-02,3.05602e-02,6.02033e-02,       Z,3.82974e-01
      19,       0,    2361,    1180,5.67269e-02,2.01307e-02,5.83493e-02,       R,4.11685e-01
      19,       1,    3495,    1747,4.63946e-02,3.83501e-02,5.24575e-02,       Z,4.77151e-01
      20,       0,    3495,    1747,4.63589e-02,2.45374e-02,4.91334e-02,       Z,5.54261e-01
      21,       0,    3495,    1747,4.63586e-02,1.62610e-02,4.76578e-02,       R,6.16068e-01
      21,       1,    5269,    2634,3.75910e-02,3.16299e-02,4.27465e-02,       Z,7.15746e-01
      22,       0,    5269,    2634,3.75698e-02,2.03848e-02,3.99405e-02,       Z,8.51180e-01
      23,       0,    5269,    2634,3.75685e-02,1.35506e-02,3.86781e-02,       R,9.61638e-01
      23,       1,    8203,    4101,3.03753e-02,2.59296e-02,3.45444e-02,       Z,1.11591e+00
      24,       0,    8203,    4101,3.03715e-02,1.64555e-02,3.22574e-02,       Z,1.34105e+00
      25,       0,    8203,    4101,3.03756e-02,1.08517e-02,3.12549e-02,       R,1.52661e+00
      25,       1,   12039,    6019,2.49139e-02,2.04874e-02,2.81248e-02,       Z,1.76797e+00
      26,       0,   12039,    6019,2.49101e-02,1.30554e-02,2.63646e-02,       Z,2.13431e+00
      27,       0,   12039,    6019,2.49112e-02,8.62997e-03,2.55894e-02,       R,2.43859e+00
      27,       1,   17971,    8985,2.02646e-02,1.68639e-02,2.30006e-02,       Z,2.78667e+00
      28,       0,   17971,    8985,2.02599e-02,1.08877e-02,2.15178e-02,       Z,3.34757e+00
      29,       0,   17971,    8985,2.02589e-02,7.25082e-03,2.08488e-02,       R,3.81462e+00
      29,       1,   27131,   13565,1.66707e-02,1.36046e-02,1.87948e-02,       Z,4.41354e+00
      30,       0,   27131,   13565,1.66704e-02,8.67957e-03,1.76316e-02,       Z,5.33663e+00
      31,       0,   27131,   13565,1.66713e-02,5.73853e-03,1.71188e-02,       R,6.10955e+00
      31,       1,   39381,   19690,1.37739e-02,1.10065e-02,1.54350e-02,       Z,6.16032e+00
      32,       0,   39381,   19690,1.37751e-02,6.96270e-03,1.45186e-02,       R,6.83664e+00
      32,       1,   57199,   28599,1.14952e-02,1.03002e-02,1.32536e-02,       Z,7.82864e+00
      33,       0,   57199,   28599,1.14952e-02,6.59662e-03,1.23004e-02,       Z,9.25790e+00
      34,       0,   57199,   28599,1.14956e-02,4.37591e-03,1.18734e-02,       R,1.04476e+01
      34,       1,   81001,   40500,9.57788e-03,7.71759e-03,1.07993e-02,       Z,1.19931e+01
      35,       0,   81001,   40500,9.57721e-03,4.99006e-03,1.01389e-02,       Z,1.43582e+01
      36,       0,   81001,   40500,9.57696e-03,3.32854e-03,9.84101e-03,       R,1.63398e+01
      36,       1,  118667,   59333,7.95874e-03,6.28138e-03,8.90544e-03,       Z,1.85568e+01
      37,       0,  118667,   59333,7.95924e-03,3.99461e-03,8.38552e-03,       Z,2.21894e+01
      38,       0,  118667,   59333,7.95973e-03,2.63801e-03,8.15808e-03,       R,2.52348e+01
      38,       1,  165681,   82840,6.72541e-03,5.00852e-03,7.43737e-03,       Z,2.90187e+01
      39,       0,  165681,   82840,6.72593e-03,3.17428e-03,7.04399e-03,       R,3.49570e+01
      39,       1,  237353,  118676,5.65643e-03,4.82897e-03,6.44810e-03,       Z,3.52233e+01
      40,       0,  237353,  118676,5.65648e-03,3.09550e-03,6.01851e-03,       Z,4.04594e+01
      41,       0,  237353,  118676,5.65657e-03,2.05561e-03,5.82656e-03,       R,4.47786e+01
      41,       1,  326373,  163186,4.78075e-03,3.65606e-03,5.33295e-03,       Z,5.11256e+01
      42,       0,  326373,  163186,4.78059e-03,2.36353e-03,5.03363e-03,       R,6.02789e+01
      42,       1,  462189,  231094,4.03487e-03,3.48714e-03,4.61301e-03,       Z,6.87518e+01
      43,       0,  462189,  231094,4.03501e-03,2.23575e-03,4.29943e-03,       Z,8.42619e+01
      44,       0,  462189,  231094,4.03513e-03,1.48417e-03,4.15918e-03,       R,9.74115e+01
      44,       1,  635997,  317998,3.43152e-03,2.59032e-03,3.80753e-03,       Z,1.12191e+02
      45,       0,  635997,  317998,3.43171e-03,1.64943e-03,3.60091e-03,       Z,1.36670e+02
      46,       0,  635997,  317998,3.43186e-03,1.09034e-03,3.51060e-03,       R,1.57374e+02
      46,       1,  895773,  447886,2.91252e-03,2.11749e-03,3.21240e-03,       Z,1.58622e+02
      47,       0,  895773,  447886,2.91251e-03,1.35527e-03,3.04812e-03,       R,1.76822e+02
      47,       1, 1220213,  610106,2.48107e-03,2.04054e-03,2.80964e-03,       Z,1.99958e+02
        }\tableDataFour

        %
        %

        \addplot+ [marker1, adaptive, forget plot]
        table [col sep=comma, x=cumulativeNdof, y=res] {\tableDataOne};

        \addplot+ [marker2, adaptive, forget plot]
        table [col sep=comma, x=cumulativeNdof, y=res] {\tableDataTwo};

        \addplot+ [marker3, adaptive, forget plot]
        table [col sep=comma, x=cumulativeNdof, y=res] {\tableDataThree};

        \addplot+ [marker4, adaptive, forget plot]
        table [col sep=comma, x=cumulativeNdof, y=res] {\tableDataFour};

        \drawslopetriangle[ST1]{0.5}{3e5}{3e-3} 
    \end{loglogaxis}
\end{tikzpicture}

%% file: figures/Fig05b_Nonlinear_condition_weighting.tex
\begin{tikzpicture}[>=stealth]
    \begin{loglogaxis}[%
            width            = 5.5cm,%
            xlabel           = {ndof},%
            ylabel           = {Least-squares functional \(N^k_\ell\)},%
            ymajorgrids      = true,%
            font             = \footnotesize,%
            grid style       = {%
                densely dotted,%
                semithick%
            },%
            legend style     = {%
                legend pos = north west,%
                font = \footnotesize%
            },%
        ]

        \addlegendimage{marker1}
        \addlegendentry{emph.\ grad.}
        \addlegendimage{marker2}
        \addlegendentry{balanced w.}
        \addlegendimage{marker3}
        \addlegendentry{downsc. flux}
        \addlegendimage{marker4}
        \addlegendentry{split w.}

        \pgfplotstableread[col sep=comma]{
       k,     ell,    ndof,    cost,   nElem,        eta,         mu,        res,    case,   maxGradU,       cond
       0,       0,     193,     193,      96,2.01791e-01,1.16608e+00,2.55303e-01,       Z,8.62011e-03,6.28126e+03
       1,       0,     193,     386,      96,2.01405e-01,1.51358e-01,2.13723e-01,       R,1.60043e-02,6.28126e+03
       1,       1,     271,     657,     135,1.63298e-01,1.91851e-01,1.79901e-01,       Z,2.95074e-02,2.34648e+04
       2,       0,     271,     928,     135,1.63037e-01,7.16045e-02,1.67435e-01,       Z,4.14307e-02,2.34648e+04
       3,       0,     271,    1199,     135,1.62901e-01,2.90517e-02,1.65264e-01,       R,4.69924e-02,2.34648e+04
       3,       1,     427,    1626,     213,1.39329e-01,8.92647e-02,1.42920e-01,       Z,7.06245e-02,7.27323e+04
       4,       0,     427,    2053,     213,1.39178e-01,2.65099e-02,1.40792e-01,       R,8.57288e-02,7.27323e+04
       4,       1,     633,    2686,     316,1.12356e-01,8.63095e-02,1.16536e-01,       Z,1.20423e-01,2.47481e+05
       5,       0,     633,    3319,     316,1.12328e-01,2.82196e-02,1.13647e-01,       Z,1.47185e-01,2.47481e+05
       6,       0,     633,    3952,     316,1.12315e-01,1.15980e-02,1.13149e-01,       R,1.60334e-01,2.47481e+05
       6,       1,     937,    4889,     468,9.12346e-02,6.65240e-02,9.38424e-02,       Z,1.62853e-01,6.49892e+05
       7,       0,     937,    5826,     468,9.12420e-02,2.03783e-02,9.19383e-02,       Z,1.70365e-01,6.49892e+05
       8,       0,     937,    6763,     468,9.12472e-02,7.78354e-03,9.16595e-02,       R,1.73925e-01,6.49892e+05
       8,       1,    1309,    8072,     654,7.54564e-02,5.18938e-02,7.74173e-02,       Z,2.33313e-01,2.58303e+06
       9,       0,    1309,    9381,     654,7.54574e-02,1.63962e-02,7.59310e-02,       Z,2.67113e-01,2.58303e+06
      10,       0,    1309,   10690,     654,7.54574e-02,6.39634e-03,7.57034e-02,       R,2.84215e-01,2.58303e+06
      10,       1,    2137,   12827,    1068,6.08165e-02,4.51230e-02,6.24039e-02,       Z,3.76092e-01,1.00592e+07
      11,       0,    2137,   14964,    1068,6.08435e-02,1.33715e-02,6.12043e-02,       Z,4.36950e-01,1.00592e+07
      12,       0,    2137,   17101,    1068,6.08539e-02,5.40803e-03,6.10093e-02,       R,4.68905e-01,1.00592e+07
      12,       1,    3089,   20190,    1544,4.92090e-02,3.62065e-02,5.05179e-02,       Z,6.00485e-01,1.72442e+07
      13,       0,    3089,   23279,    1544,4.92085e-02,1.11692e-02,4.94582e-02,       Z,6.94138e-01,1.72442e+07
      14,       0,    3089,   26368,    1544,4.92082e-02,4.33988e-03,4.92970e-02,       R,7.44808e-01,1.72442e+07
      14,       1,    4589,   30957,    2294,4.01562e-02,2.87710e-02,4.11967e-02,       Z,9.64786e-01,7.01696e+07
      15,       0,    4589,   35546,    2294,4.01500e-02,9.09263e-03,4.03372e-02,       Z,1.12408e+00,7.01696e+07
      16,       0,    4589,   40135,    2294,4.01479e-02,3.57172e-03,4.02034e-02,       R,1.21258e+00,7.01696e+07
      16,       1,    7187,   47322,    3593,3.25030e-02,2.38362e-02,3.33021e-02,       Z,1.56129e+00,2.82663e+08
      17,       0,    7187,   54509,    3593,3.25053e-02,7.16600e-03,3.26503e-02,       Z,1.82386e+00,2.82663e+08
      18,       0,    7187,   61696,    3593,3.25067e-02,2.87589e-03,3.25452e-02,       R,1.97272e+00,2.82663e+08
      18,       1,   10419,   72115,    5209,2.68242e-02,1.85855e-02,2.74027e-02,       Z,2.52575e+00,1.14756e+09
      19,       0,   10419,   82534,    5209,2.68268e-02,5.54342e-03,2.69252e-02,       Z,2.95377e+00,1.14756e+09
      20,       0,   10419,   92953,    5209,2.68280e-02,2.17301e-03,2.68520e-02,       R,3.20017e+00,1.14756e+09
      20,       1,   15673,  108626,    7836,2.18321e-02,1.57425e-02,2.23725e-02,       Z,3.24298e+00,2.78210e+09
      21,       0,   15673,  124299,    7836,2.18306e-02,4.87122e-03,2.19192e-02,       Z,3.40512e+00,2.78210e+09
      22,       0,   15673,  139972,    7836,2.18303e-02,1.91498e-03,2.18486e-02,       R,3.49461e+00,2.78210e+09
      22,       1,   23249,  163221,   11624,1.78925e-02,1.26526e-02,1.83225e-02,       Z,4.42420e+00,1.10895e+10
      23,       0,   23249,  186470,   11624,1.78942e-02,3.92574e-03,1.79647e-02,       R,5.04005e+00,1.10895e+10
      23,       1,   34811,  221281,   17405,1.47452e-02,1.08718e-02,1.51273e-02,       Z,6.43375e+00,4.45557e+10
      24,       0,   34811,  256092,   17405,1.47492e-02,3.35415e-03,1.48135e-02,       Z,7.70067e+00,4.45557e+10
      25,       0,   34811,  290903,   17405,1.47509e-02,1.34287e-03,1.47634e-02,       R,8.44703e+00,4.45557e+10
      25,       1,   49807,  340710,   24903,1.22872e-02,8.27143e-03,1.25374e-02,       Z,1.06775e+01,1.78414e+11
      26,       0,   49807,  390517,   24903,1.22882e-02,2.48284e-03,1.23277e-02,       Z,1.25320e+01,1.78414e+11
      27,       0,   49807,  440324,   24903,1.22886e-02,9.69309e-04,1.22958e-02,       R,1.36208e+01,1.78414e+11
      27,       1,   71459,  511783,   35729,1.01768e-02,6.95582e-03,1.04136e-02,       Z,1.69306e+01,2.92323e+11
      28,       0,   71459,  583242,   35729,1.01764e-02,2.20749e-03,1.02143e-02,       Z,1.96559e+01,2.92323e+11
      29,       0,   71459,  654701,   35729,1.01763e-02,8.74291e-04,1.01829e-02,       R,2.12546e+01,2.92323e+11
      29,       1,  103867,  758568,   51933,8.51327e-03,5.64327e-03,8.68291e-03,       Z,2.69515e+01,1.18448e+12
      30,       0,  103867,  862435,   51933,8.51369e-03,1.70441e-03,8.54125e-03,       Z,3.14958e+01,1.18448e+12
      31,       0,  103867,  966302,   51933,8.51393e-03,6.78939e-04,8.51880e-03,       R,3.41693e+01,1.18448e+12
      31,       1,  148269, 1114571,   74134,7.12601e-03,4.70829e-03,7.25922e-03,       Z,3.45545e+01,2.87800e+12
      32,       0,  148269, 1262840,   74134,7.12693e-03,1.37864e-03,7.14859e-03,       R,3.63128e+01,2.87800e+12
      32,       1,  207617, 1470457,  103808,6.02299e-03,4.05184e-03,6.15480e-03,       Z,4.56736e+01,1.15228e+13
      33,       0,  207617, 1678074,  103808,6.02330e-03,1.26496e-03,6.04426e-03,       Z,5.25239e+01,1.15228e+13
      34,       0,  207617, 1885691,  103808,6.02347e-03,4.99684e-04,6.02701e-03,       R,5.65371e+01,1.15228e+13
      34,       1,  290993, 2176684,  145496,5.05900e-03,3.30732e-03,5.16659e-03,       Z,6.98422e+01,1.88621e+13
      35,       0,  290993, 2467677,  145496,5.05879e-03,1.04959e-03,5.07584e-03,       R,8.05348e+01,1.88621e+13
      35,       1,  407169, 2874846,  203584,4.29572e-03,2.87049e-03,4.38946e-03,       Z,1.02460e+02,7.61739e+13
      36,       0,  407169, 3282015,  203584,4.29585e-03,9.01526e-04,4.31101e-03,       Z,1.23282e+02,7.61739e+13
      37,       0,  407169, 3689184,  203584,4.29592e-03,3.59901e-04,4.29850e-03,       R,1.35684e+02,7.61739e+13
      37,       1,  572869, 4262053,  286434,3.62967e-03,2.32595e-03,3.69334e-03,       Z,1.37730e+02,1.84783e+14
      38,       0,  572869, 4834922,  286434,3.63005e-03,6.80692e-04,3.64041e-03,       Z,1.46105e+02,1.84783e+14
      39,       0,  572869, 5407791,  286434,3.63022e-03,2.71917e-04,3.63200e-03,       R,1.50865e+02,1.84783e+14
      39,       1,  787563, 6195354,  393781,3.09451e-03,1.91740e-03,3.14701e-03,       Z,1.87968e+02,7.39395e+14
      40,       0,  787563, 6982917,  393781,3.09452e-03,5.72286e-04,3.10273e-03,       R,2.13507e+02,7.39395e+14
      40,       1, 1097597, 8080514,  548798,2.61659e-03,1.74844e-03,2.67739e-03,       Z,2.65291e+02,1.21045e+15
        }\tableDataOne

        \pgfplotstableread[col sep=comma]{
       k,     ell,    ndof,    cost,   nElem,        eta,         mu,        res,    case,   maxGradU,       cond
       0,       0,     193,     193,      96,9.60646e-02,8.06816e-01,2.55284e-01,       Z,9.01679e-03,6.27801e+03
       1,       0,     193,     386,      96,2.01459e-01,1.54929e-01,2.09919e-01,       R,2.75277e-02,6.27801e+03
       1,       1,     271,     657,     135,1.62861e-01,1.95103e-01,1.91307e-01,       Z,6.78911e-02,2.34620e+04
       2,       0,     271,     928,     135,1.61705e-01,9.73864e-02,1.90339e-01,       Z,4.81053e-02,2.34620e+04
       3,       0,     271,    1199,     135,1.61894e-01,9.52332e-02,1.95929e-01,       R,5.20205e-02,2.34620e+04
       3,       1,     427,    1626,     213,1.38310e-01,1.27080e-01,1.90134e-01,       Z,1.08792e-01,7.27217e+04
       4,       0,     427,    2053,     213,1.38392e-01,1.28396e-01,2.00474e-01,       R,9.69161e-02,7.27217e+04
       4,       1,     641,    2694,     320,1.09456e-01,1.53809e-01,1.92033e-01,       Z,1.78676e-01,2.46226e+05
       5,       0,     641,    3335,     320,1.09937e-01,1.56749e-01,2.06895e-01,       Z,1.70231e-01,2.46226e+05
       6,       0,     641,    3976,     320,1.09855e-01,1.74688e-01,2.25024e-01,       R,1.71260e-01,2.46226e+05
       6,       1,    1019,    4995,     509,8.55644e-02,1.87784e-01,2.19230e-01,       Z,1.76531e-01,6.44570e+05
       7,       0,    1019,    6014,     509,8.58007e-02,2.01544e-01,2.41107e-01,       Z,1.76739e-01,6.44570e+05
       8,       0,    1019,    7033,     509,8.57667e-02,2.25158e-01,2.66721e-01,       Z,1.76486e-01,6.44570e+05
       9,       0,    1019,    8052,     509,8.57760e-02,2.52393e-01,2.95854e-01,       R,1.76595e-01,6.44570e+05
       9,       1,    1613,    9665,     806,6.97576e-02,2.57281e-01,2.93780e-01,       Z,3.13816e-01,2.53967e+06
      10,       0,    1613,   11278,     806,6.98129e-02,2.85307e-01,3.27250e-01,       Z,2.99466e-01,2.53967e+06
      11,       0,    1613,   12891,     806,6.98094e-02,3.19666e-01,3.65210e-01,       R,3.01260e-01,2.53967e+06
      11,       1,    2357,   15248,    1178,5.72907e-02,3.22145e-01,3.64392e-01,       Z,5.14428e-01,1.03558e+07
      12,       0,    2357,   17605,    1178,5.73334e-02,3.59831e-01,4.07464e-01,       R,5.00202e-01,1.03558e+07
      12,       1,    3413,   21018,    1706,4.68728e-02,3.61343e-01,4.06800e-01,       Z,8.30419e-01,4.21938e+07
      13,       0,    3413,   24431,    1706,4.68849e-02,4.04081e-01,4.55436e-01,       Z,8.23118e-01,4.21938e+07
      14,       0,    3413,   27844,    1706,4.68849e-02,4.53008e-01,5.10124e-01,       R,8.23705e-01,4.21938e+07
      14,       1,    5143,   32987,    2571,3.79220e-02,4.53847e-01,5.09738e-01,       Z,1.35254e+00,1.68991e+08
      15,       0,    5143,   38130,    2571,3.79254e-02,5.08322e-01,5.71201e-01,       Z,1.35260e+00,1.68991e+08
      16,       0,    5143,   43273,    2571,3.79264e-02,5.69938e-01,6.40203e-01,       R,1.35330e+00,1.68991e+08
      16,       1,    8095,   51368,    4047,3.06069e-02,5.70378e-01,6.40043e-01,       Z,2.19506e+00,6.82820e+08
      17,       0,    8095,   59463,    4047,3.06314e-02,6.39309e-01,7.17504e-01,       Z,2.21323e+00,6.82820e+08
      18,       0,    8095,   67558,    4047,3.06267e-02,7.16850e-01,8.04400e-01,       R,2.21466e+00,6.82820e+08
      18,       1,   11791,   79349,    5895,2.52295e-02,7.17060e-01,8.04317e-01,       Z,3.42841e+00,1.13550e+09
      19,       0,   11791,   91140,    5895,2.52326e-02,8.03921e-01,9.01795e-01,       Z,3.47959e+00,1.13550e+09
      20,       0,   11791,  102931,    5895,2.52329e-02,9.01442e-01,1.01112e+00,       R,3.48199e+00,1.13550e+09
      20,       1,   17589,  120520,    8794,2.05474e-02,9.01560e-01,1.01106e+00,       Z,5.55291e+00,4.57347e+09
      21,       0,   17589,  138109,    8794,2.05391e-02,1.01085e+00,1.13366e+00,       Z,5.65807e+00,4.57347e+09
      22,       0,   17589,  155698,    8794,2.05420e-02,1.13348e+00,1.27116e+00,       R,5.66404e+00,4.57347e+09
      22,       1,   26267,  181965,   13133,1.69073e-02,1.13354e+00,1.27113e+00,       Z,8.95155e+00,1.83059e+10
      23,       0,   26267,  208232,   13133,1.69138e-02,1.27102e+00,1.42531e+00,       Z,9.14293e+00,1.83059e+10
      24,       0,   26267,  234499,   13133,1.69125e-02,1.42521e+00,1.59821e+00,       R,9.15426e+00,1.83059e+10
      24,       1,   38249,  272748,   19124,1.39836e-02,1.42524e+00,1.59820e+00,       Z,9.41418e+00,4.48723e+10
      25,       0,   38249,  310997,   19124,1.39914e-02,1.59814e+00,1.79208e+00,       Z,9.42253e+00,4.48723e+10
      26,       0,   38249,  349246,   19124,1.39897e-02,1.79202e+00,2.00948e+00,       R,9.42309e+00,4.48723e+10
      26,       1,   55879,  405125,   27939,1.16120e-02,1.79204e+00,2.00947e+00,       Z,1.47847e+01,1.79581e+11
      27,       0,   55879,  461004,   27939,1.16113e-02,2.00944e+00,2.25325e+00,       Z,1.51095e+01,1.79581e+11
      28,       0,   55879,  516883,   27939,1.16117e-02,2.25322e+00,2.52661e+00,       R,1.51282e+01,1.79581e+11
      28,       1,   79695,  596578,   39847,9.65584e-03,2.25323e+00,2.52660e+00,       Z,2.35544e+01,7.18229e+11
      29,       0,   79695,  676273,   39847,9.65551e-03,2.52659e+00,2.83313e+00,       R,2.40957e+01,7.18229e+11
      29,       1,  115623,  791896,   57811,8.07827e-03,2.52659e+00,2.83313e+00,       Z,3.63136e+01,1.17224e+12
      30,       0,  115623,  907519,   57811,8.08148e-03,2.83311e+00,3.17684e+00,       Z,3.71305e+01,1.17224e+12
      31,       0,  115623, 1023142,   57811,8.08075e-03,3.17683e+00,3.56225e+00,       R,3.71789e+01,1.17224e+12
      31,       1,  161421, 1184563,   80710,6.81258e-03,3.17683e+00,3.56225e+00,       Z,5.85275e+01,4.76338e+12
      32,       0,  161421, 1345984,   80710,6.81395e-03,3.56224e+00,3.99442e+00,       Z,5.99063e+01,4.76338e+12
      33,       0,  161421, 1507405,   80710,6.81367e-03,3.99442e+00,4.47903e+00,       R,5.99869e+01,4.76338e+12
      33,       1,  231105, 1738510,  115552,5.72127e-03,3.99442e+00,4.47903e+00,       Z,6.15426e+01,1.15455e+13
      34,       0,  231105, 1969615,  115552,5.72058e-03,4.47902e+00,5.02242e+00,       R,6.16270e+01,1.15455e+13
      34,       1,  318659, 2288274,  159329,4.83695e-03,4.47902e+00,5.02242e+00,       Z,9.58121e+01,4.61966e+13
      35,       0,  318659, 2606933,  159329,4.83707e-03,5.02242e+00,5.63174e+00,       Z,9.80175e+01,4.61966e+13
      36,       0,  318659, 2925592,  159329,4.83707e-03,5.63174e+00,6.31499e+00,       R,9.81457e+01,4.61966e+13
      36,       1,  448975, 3374567,  224487,4.09767e-03,5.63174e+00,6.31499e+00,       Z,1.47648e+02,7.55826e+13
      37,       0,  448975, 3823542,  224487,4.09891e-03,6.31499e+00,7.08112e+00,       R,1.50833e+02,7.55826e+13
      37,       1,  620141, 4443683,  310070,3.47793e-03,6.31499e+00,7.08112e+00,       Z,1.54955e+02,1.84930e+14
      38,       0,  620141, 5063824,  310070,3.47810e-03,7.08112e+00,7.94021e+00,       Z,1.55370e+02,1.84930e+14
      39,       0,  620141, 5683965,  310070,3.47809e-03,7.94021e+00,8.90351e+00,       R,1.55394e+02,1.84930e+14
      39,       1,  868029, 6551994,  434014,2.95422e-03,7.94021e+00,8.90351e+00,       Z,2.42490e+02,7.40164e+14
      40,       0,  868029, 7420023,  434014,2.95403e-03,8.90351e+00,9.98369e+00,       Z,2.48087e+02,7.40164e+14
      41,       0,  868029, 8288052,  434014,2.95419e-03,9.98369e+00,1.11949e+01,       R,2.48411e+02,7.40164e+14
      41,       1, 1187487, 9475539,  593743,2.51468e-03,9.98369e+00,1.11949e+01,       Z,3.73280e+02,1.20212e+15
        }\tableDataTwo

        \pgfplotstableread[col sep=comma]{
       k,     ell,    ndof,    cost,   nElem,        eta,         mu,        res,    case,   maxGradU,       cond
       0,       0,     193,     193,      96,4.54996e-02,5.56005e-01,2.56527e-01,       Z,9.21402e-03,6.27649e+03
       1,       0,     193,     386,      96,2.02972e-01,1.56388e-01,3.29270e-01,       R,6.23118e-02,6.27649e+03
       1,       1,     271,     657,     135,1.63709e-01,1.97114e-01,4.29268e-01,       Z,2.03996e-01,2.34608e+04
       2,       0,     271,     928,     135,1.60311e-01,3.96790e-01,1.07123e+00,       Z,1.90154e-02,2.34608e+04
       3,       0,     271,    1199,     135,1.62342e-01,1.05835e+00,3.36292e+00,       R,6.80706e-01,2.34608e+04
       3,       1,     427,    1626,     213,1.38599e-01,1.06172e+00,3.37008e+00,       Z,1.09683e+00,7.27168e+04
       4,       0,     427,    2053,     213,1.36545e-01,3.36723e+00,1.15418e+01,       R,4.34832e-01,7.27168e+04
       4,       1,     615,    2668,     307,1.08537e-01,3.36825e+00,1.15432e+01,       Z,4.33294e-01,7.36325e+04
       5,       0,     615,    3283,     307,1.16765e-01,1.15426e+01,4.02283e+01,       Z,3.56972e+00,7.36325e+04
       6,       0,     615,    3898,     307,1.09789e-01,4.02282e+01,1.40727e+02,       R,2.53758e+00,7.36325e+04
       6,       1,     989,    4887,     494,8.54091e-02,4.02282e+01,1.40727e+02,       Z,3.14846e+00,2.56918e+05
       7,       0,     989,    5876,     494,9.99854e-02,1.40727e+02,4.92535e+02,       R,9.85007e+00,2.56918e+05
       7,       1,    1249,    7125,     624,8.44980e-02,1.40727e+02,4.92535e+02,       Z,1.25577e+01,1.05125e+06
       8,       0,    1249,    8374,     624,8.14751e-02,4.92535e+02,1.72391e+03,       Z,8.42874e+00,1.05125e+06
       9,       0,    1249,    9623,     624,8.85706e-02,1.72391e+03,6.03373e+03,       R,2.56954e+01,1.05125e+06
       9,       1,    1851,   11474,     925,7.39555e-02,1.72391e+03,6.03373e+03,       Z,3.13280e+01,4.20126e+06
      10,       0,    1851,   13325,     925,7.49541e-02,6.03373e+03,2.11181e+04,       R,2.18526e+01,4.20126e+06
      10,       1,    2493,   15818,    1246,6.06636e-02,6.03373e+03,2.11181e+04,       Z,2.18156e+01,4.27135e+06
      11,       0,    2493,   18311,    1246,6.79063e-02,2.11181e+04,7.39135e+04,       R,5.35179e+01,4.27135e+06
      11,       1,    3519,   21830,    1759,5.58090e-02,2.11181e+04,7.39135e+04,       Z,6.36517e+01,1.70756e+07
        }\tableDataThree

        \pgfplotstableread[col sep=comma]{
       k,     ell,    ndof,    cost,   nElem,        eta,         mu,        res,    case,   maxGradU,       cond
       0,       0,     193,     193,      96,3.96416e-01,1.62597e+00,2.63724e-01,       Z,7.94203e-03,6.28740e+03
       1,       0,     193,     386,      96,2.11806e-01,1.54617e-01,2.31688e-01,       Z,1.18378e-02,6.28740e+03
       2,       0,     193,     579,      96,2.05213e-01,1.06488e-01,2.18439e-01,       Z,1.46085e-02,6.28740e+03
       3,       0,     193,     772,      96,2.04843e-01,7.44515e-02,2.11775e-01,       Z,1.67011e-02,6.28740e+03
       4,       0,     193,     965,      96,2.04752e-01,5.21035e-02,2.08434e-01,       Z,1.82525e-02,6.28740e+03
       5,       0,     193,    1158,      96,2.04693e-01,3.65399e-02,2.06763e-01,       Z,1.93848e-02,6.28740e+03
       6,       0,     193,    1351,      96,2.04650e-01,2.56655e-02,2.05925e-01,       Z,2.02029e-02,6.28740e+03
       7,       0,     193,    1544,      96,2.04620e-01,1.80479e-02,2.05504e-01,       R,2.07900e-02,6.28740e+03
       7,       1,     271,    1815,     135,1.66582e-01,1.20189e-01,1.81366e-01,       Z,2.82188e-02,2.34700e+04
       8,       0,     271,    2086,     135,1.64779e-01,7.52204e-02,1.72208e-01,       Z,3.42216e-02,2.34700e+04
       9,       0,     271,    2357,     135,1.64691e-01,4.95071e-02,1.68254e-01,       R,3.88302e-02,2.34700e+04
       9,       1,     437,    2794,     218,1.38503e-01,1.01937e-01,1.51718e-01,       Z,4.95674e-02,7.27531e+04
      10,       0,     437,    3231,     218,1.37734e-01,6.32914e-02,1.43864e-01,       Z,6.18173e-02,7.27531e+04
      11,       0,     437,    3668,     218,1.37669e-01,4.12609e-02,1.40555e-01,       R,7.14215e-02,7.27531e+04
      11,       1,     653,    4321,     326,1.10661e-01,9.17032e-02,1.24456e-01,       Z,8.67691e-02,2.49264e+05
      12,       0,     653,    4974,     326,1.10160e-01,5.77624e-02,1.16632e-01,       Z,1.07081e-01,2.49264e+05
      13,       0,     653,    5627,     326,1.10156e-01,3.80931e-02,1.13229e-01,       R,1.23124e-01,2.49264e+05
      13,       1,    1069,    6696,     534,8.54867e-02,7.92299e-02,9.94177e-02,       Z,1.47598e-01,1.04061e+06
      14,       0,    1069,    7765,     534,8.51726e-02,5.12155e-02,9.18009e-02,       Z,1.80778e-01,1.04061e+06
      15,       0,    1069,    8834,     534,8.51616e-02,3.41827e-02,8.83159e-02,       R,2.07239e-01,1.04061e+06
      15,       1,    1593,   10427,     796,7.01506e-02,5.91596e-02,7.97951e-02,       Z,2.45396e-01,4.24439e+06
      16,       0,    1593,   12020,     796,7.00281e-02,3.82170e-02,7.45291e-02,       Z,2.98846e-01,4.24439e+06
      17,       0,    1593,   13613,     796,7.00243e-02,2.54642e-02,7.21478e-02,       R,3.41873e-01,4.24439e+06
      17,       1,    2361,   15974,    1180,5.67807e-02,4.82471e-02,6.44427e-02,       Z,3.45998e-01,1.03678e+07
      18,       0,    2361,   18335,    1180,5.67247e-02,3.05602e-02,6.02033e-02,       Z,3.82974e-01,1.03678e+07
      19,       0,    2361,   20696,    1180,5.67269e-02,2.01307e-02,5.83493e-02,       R,4.11685e-01,1.03678e+07
      19,       1,    3495,   24191,    1747,4.63946e-02,3.83501e-02,5.24575e-02,       Z,4.77151e-01,4.18350e+07
      20,       0,    3495,   27686,    1747,4.63589e-02,2.45374e-02,4.91334e-02,       Z,5.54261e-01,4.18350e+07
      21,       0,    3495,   31181,    1747,4.63586e-02,1.62610e-02,4.76578e-02,       R,6.16068e-01,4.18350e+07
      21,       1,    5269,   36450,    2634,3.75910e-02,3.16299e-02,4.27465e-02,       Z,7.15746e-01,1.72944e+08
      22,       0,    5269,   41719,    2634,3.75698e-02,2.03848e-02,3.99405e-02,       Z,8.51180e-01,1.72944e+08
      23,       0,    5269,   46988,    2634,3.75685e-02,1.35506e-02,3.86781e-02,       R,9.61638e-01,1.72944e+08
      23,       1,    8203,   55191,    4101,3.03753e-02,2.59296e-02,3.45444e-02,       Z,1.11591e+00,6.91206e+08
      24,       0,    8203,   63394,    4101,3.03715e-02,1.64555e-02,3.22574e-02,       Z,1.34105e+00,6.91206e+08
      25,       0,    8203,   71597,    4101,3.03756e-02,1.08517e-02,3.12549e-02,       R,1.52661e+00,6.91206e+08
      25,       1,   12039,   83636,    6019,2.49139e-02,2.04874e-02,2.81248e-02,       Z,1.76797e+00,2.76558e+09
      26,       0,   12039,   95675,    6019,2.49101e-02,1.30554e-02,2.63646e-02,       Z,2.13431e+00,2.76558e+09
      27,       0,   12039,  107714,    6019,2.49112e-02,8.62997e-03,2.55894e-02,       R,2.43859e+00,2.76558e+09
      27,       1,   17971,  125685,    8985,2.02646e-02,1.68639e-02,2.30006e-02,       Z,2.78667e+00,4.56507e+09
      28,       0,   17971,  143656,    8985,2.02599e-02,1.08877e-02,2.15178e-02,       Z,3.34757e+00,4.56507e+09
      29,       0,   17971,  161627,    8985,2.02589e-02,7.25082e-03,2.08488e-02,       R,3.81462e+00,4.56507e+09
      29,       1,   27131,  188758,   13565,1.66707e-02,1.36046e-02,1.87948e-02,       Z,4.41354e+00,1.84491e+10
      30,       0,   27131,  215889,   13565,1.66704e-02,8.67957e-03,1.76316e-02,       Z,5.33663e+00,1.84491e+10
      31,       0,   27131,  243020,   13565,1.66713e-02,5.73853e-03,1.71188e-02,       R,6.10955e+00,1.84491e+10
      31,       1,   39381,  282401,   19690,1.37739e-02,1.10065e-02,1.54350e-02,       Z,6.16032e+00,4.47006e+10
      32,       0,   39381,  321782,   19690,1.37751e-02,6.96270e-03,1.45186e-02,       R,6.83664e+00,4.47006e+10
      32,       1,   57199,  378981,   28599,1.14952e-02,1.03002e-02,1.32536e-02,       Z,7.82864e+00,1.78595e+11
      33,       0,   57199,  436180,   28599,1.14952e-02,6.59662e-03,1.23004e-02,       Z,9.25790e+00,1.78595e+11
      34,       0,   57199,  493379,   28599,1.14956e-02,4.37591e-03,1.18734e-02,       R,1.04476e+01,1.78595e+11
      34,       1,   81001,  574380,   40500,9.57788e-03,7.71759e-03,1.07993e-02,       Z,1.19931e+01,7.20266e+11
      35,       0,   81001,  655381,   40500,9.57721e-03,4.99006e-03,1.01389e-02,       Z,1.43582e+01,7.20266e+11
      36,       0,   81001,  736382,   40500,9.57696e-03,3.32854e-03,9.84101e-03,       R,1.63398e+01,7.20266e+11
      36,       1,  118667,  855049,   59333,7.95874e-03,6.28138e-03,8.90544e-03,       Z,1.85568e+01,1.17837e+12
      37,       0,  118667,  973716,   59333,7.95924e-03,3.99461e-03,8.38552e-03,       Z,2.21894e+01,1.17837e+12
      38,       0,  118667, 1092383,   59333,7.95973e-03,2.63801e-03,8.15808e-03,       R,2.52348e+01,1.17837e+12
      38,       1,  165681, 1258064,   82840,6.72541e-03,5.00852e-03,7.43737e-03,       Z,2.90187e+01,4.75215e+12
      39,       0,  165681, 1423745,   82840,6.72593e-03,3.17428e-03,7.04399e-03,       R,3.49570e+01,4.75215e+12
      39,       1,  237353, 1661098,  118676,5.65643e-03,4.82897e-03,6.44810e-03,       Z,3.52233e+01,1.15286e+13
      40,       0,  237353, 1898451,  118676,5.65648e-03,3.09550e-03,6.01851e-03,       Z,4.04594e+01,1.15286e+13
      41,       0,  237353, 2135804,  118676,5.65657e-03,2.05561e-03,5.82656e-03,       R,4.47786e+01,1.15286e+13
      41,       1,  326373, 2462177,  163186,4.78075e-03,3.65606e-03,5.33295e-03,       Z,5.11256e+01,4.62162e+13
      42,       0,  326373, 2788550,  163186,4.78059e-03,2.36353e-03,5.03363e-03,       R,6.02789e+01,4.62162e+13
      42,       1,  462189, 3250739,  231094,4.03487e-03,3.48714e-03,4.61301e-03,       Z,6.87518e+01,7.56118e+13
      43,       0,  462189, 3712928,  231094,4.03501e-03,2.23575e-03,4.29943e-03,       Z,8.42619e+01,7.56118e+13
      44,       0,  462189, 4175117,  231094,4.03513e-03,1.48417e-03,4.15918e-03,       R,9.74115e+01,7.56118e+13
      44,       1,  635997, 4811114,  317998,3.43152e-03,2.59032e-03,3.80753e-03,       Z,1.12191e+02,3.04977e+14
      45,       0,  635997, 5447111,  317998,3.43171e-03,1.64943e-03,3.60091e-03,       Z,1.36670e+02,3.04977e+14
      46,       0,  635997, 6083108,  317998,3.43186e-03,1.09034e-03,3.51060e-03,       R,1.57374e+02,3.04977e+14
      46,       1,  895773, 6978881,  447886,2.91252e-03,2.11749e-03,3.21240e-03,       Z,1.58622e+02,7.38998e+14
      47,       0,  895773, 7874654,  447886,2.91251e-03,1.35527e-03,3.04812e-03,       R,1.76822e+02,7.38998e+14
      47,       1, 1220213, 9094867,  610106,2.48107e-03,2.04054e-03,2.80964e-03,       Z,1.99958e+02,1.20812e+15
        }\tableDataFour

        %
        %

        \addplot+ [marker1, adaptive, forget plot]
        table [col sep=comma, x=ndof, y=cond] {\tableDataOne};

        \addplot+ [marker2, adaptive, forget plot]
        table [col sep=comma, x=ndof, y=cond] {\tableDataTwo};

        \addplot+ [marker3, adaptive, forget plot]
        table [col sep=comma, x=ndof, y=cond] {\tableDataThree};

        \addplot+ [marker4, adaptive, forget plot]
        table [col sep=comma, x=ndof, y=cond] {\tableDataFour};

        \drawslopetriangleup[ST1]{3}{7e3}{1e8} 
    \end{loglogaxis}
\end{tikzpicture}

%% file: figures/Fig06_PorousMedium_rhs.tex
\begin{tikzpicture}
    \begin{axis}[%
        axis equal image,%
        width = 6.1cm,%
        ymin  = -1.1,%
        xmin  = -1.1,%
        ymax  = 1.1,%
        xmax  = 1.1,%
        font  = \footnotesize%
    ]
        \addplot[%
            mark      = none,%
            thick,%
            line cap  = round,%
            line join = round,%
        ] coordinates {
            (0, 0)
            (1, 0)
            (1, -1)
            (-1, -1)
            (-1, 1)
            (0, 1)
            (0, 0)
        };
        \fill[TUyellow]
            (axis cs: -0.6, 0.6) -- (axis cs: -0.6, 0.4) --
            (axis cs: -0.4, 0.4) -- (axis cs: -0.4, 0.6) -- cycle;
        \draw [decorate, decoration={brace, amplitude=2pt}, semithick, gray]
            (axis cs: -0.6, 0.4) -- (axis cs: -0.6, 0.6);
        \draw [decorate, decoration={brace, amplitude=2pt}, semithick, gray]
            (axis cs: -0.4, 0.4) -- (axis cs: -0.6, 0.4);
        \path (axis cs: -0.6, 0.5) node[left] {\scriptsize\color{gray}\(0.2\)};
        \path (axis cs: -0.5, 0.4) node[below] {\scriptsize\color{gray}\(0.2\)};
        \draw
            (axis cs: 0.3, 0.7) -- (axis cs: 0.3, 0.85) --
            (axis cs: 0.45, 0.85) -- (axis cs: 0.45, 0.7)
            node[midway,right] {\scriptsize\(f \equiv 0\)} -- cycle;
        \draw[fill=TUyellow]
            (axis cs: 0.3, 0.35) -- (axis cs: 0.3, 0.5) --
            (axis cs: 0.45, 0.5) -- (axis cs: 0.45, 0.35)
            node[midway,right] {\scriptsize\(f \equiv 1\)} -- cycle;
    \end{axis}
\end{tikzpicture}

%% file: figures/Fig06_PorousMedium_mesh.tex
\begin{tikzpicture}
    \begin{axis}[%
        axis equal image,%
        width=6.1cm,%
        xmin=-1.1, xmax=1.1,%
        ymin=-1.1, ymax=1.1,%
        font=\footnotesize%
    ]
        \addplot graphics [xmin=-1, xmax=1, ymin=-1, ymax=1]
            {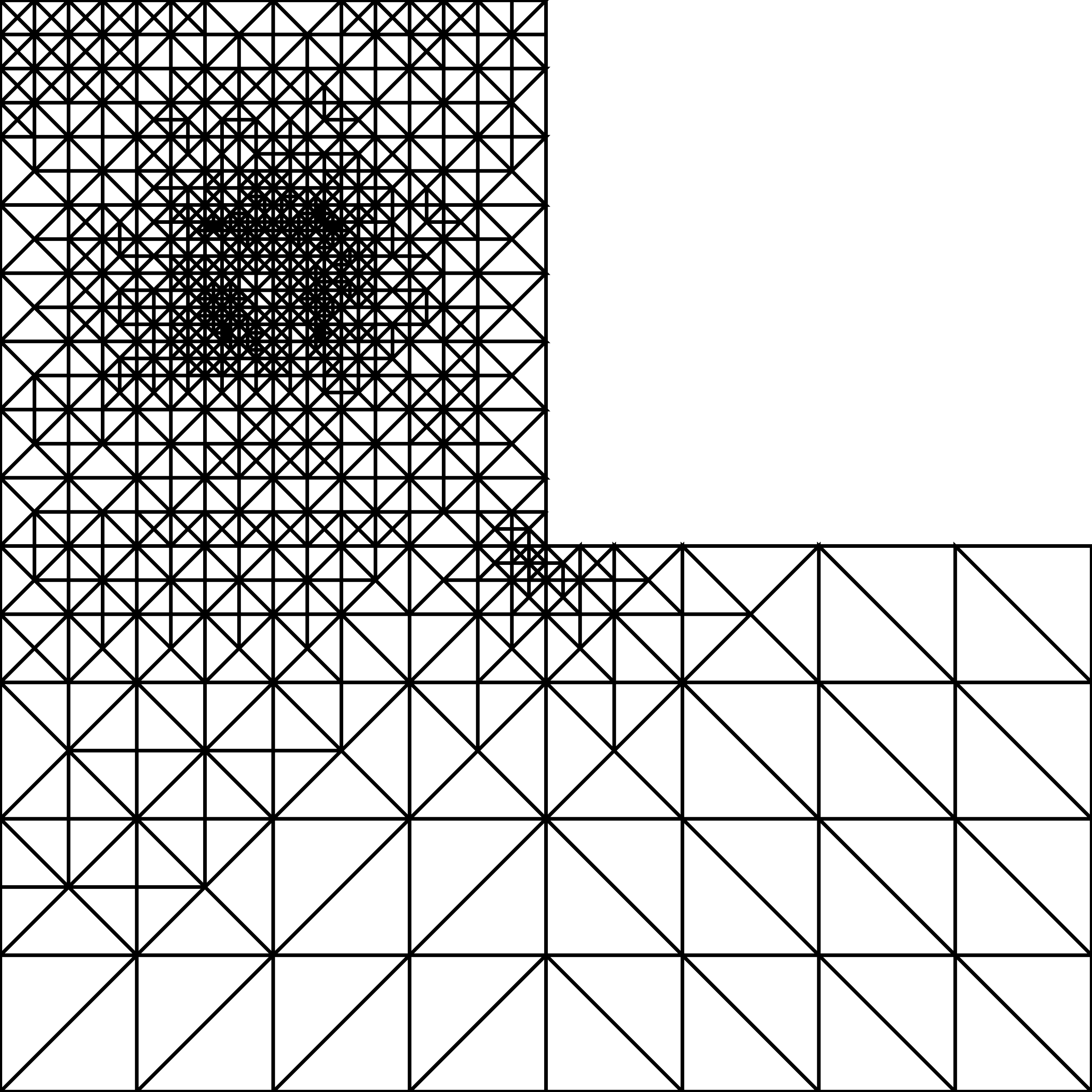};
    \end{axis}
\end{tikzpicture}

%% file: figures/Fig06_PorousMedium_potential.tex
\begin{tikzpicture}
    \begin{axis}[%
        axis equal image,%
        width=6.1cm,%
        xmin=-1.1, xmax=1.1,%
        ymin=-1.1, ymax=1.1,%
        font=\footnotesize,%
        point meta min=0.0,%
        point meta max=4.65942540e-03,%
        colorbar,%
        colorbar style={%
            font=\footnotesize,%
            width=2.5mm,%
            title style={yshift=-2mm},%
        },%
    ]
        \addplot graphics [xmin=-1, xmax=1, ymin=-1, ymax=1] {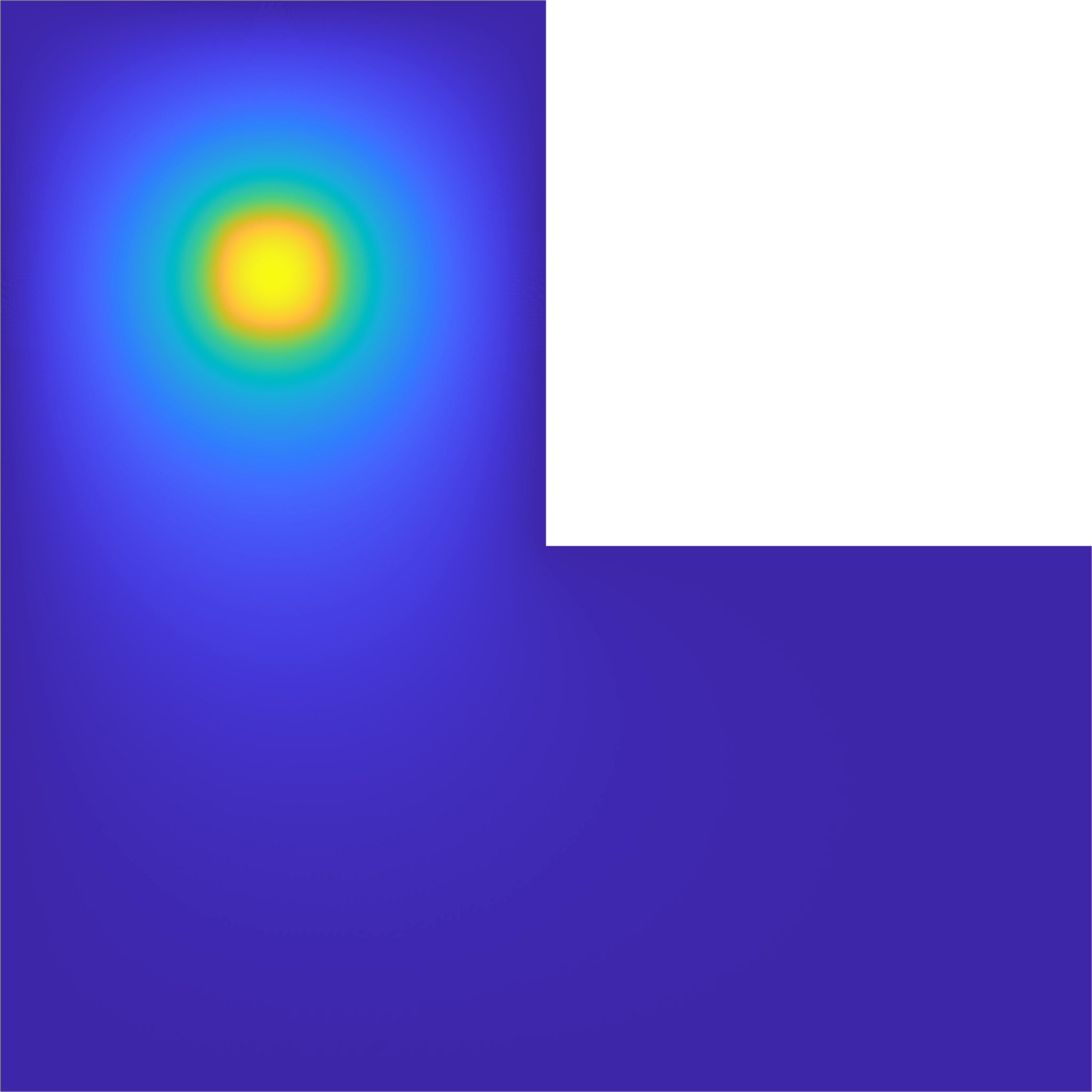};
    \end{axis}
\end{tikzpicture}

%% file: figures/Fig06_PorousMedium_flux.tex
\begin{tikzpicture}
    \begin{axis}[%
        axis equal image,%
        width=6.1cm,%
        xmin=-1.1, xmax=1.1,%
        ymin=-1.1, ymax=1.1,%
        font=\footnotesize,%
        point meta min=0.0,%
        point meta max=3.78967248e-02,%
        colorbar,%
        colorbar style={%
            font=\footnotesize,%
            width=2.5mm,%
            title style={yshift=-2mm},%
        },%
    ]
        \addplot graphics [xmin=-1, xmax=1, ymin=-1, ymax=1]
{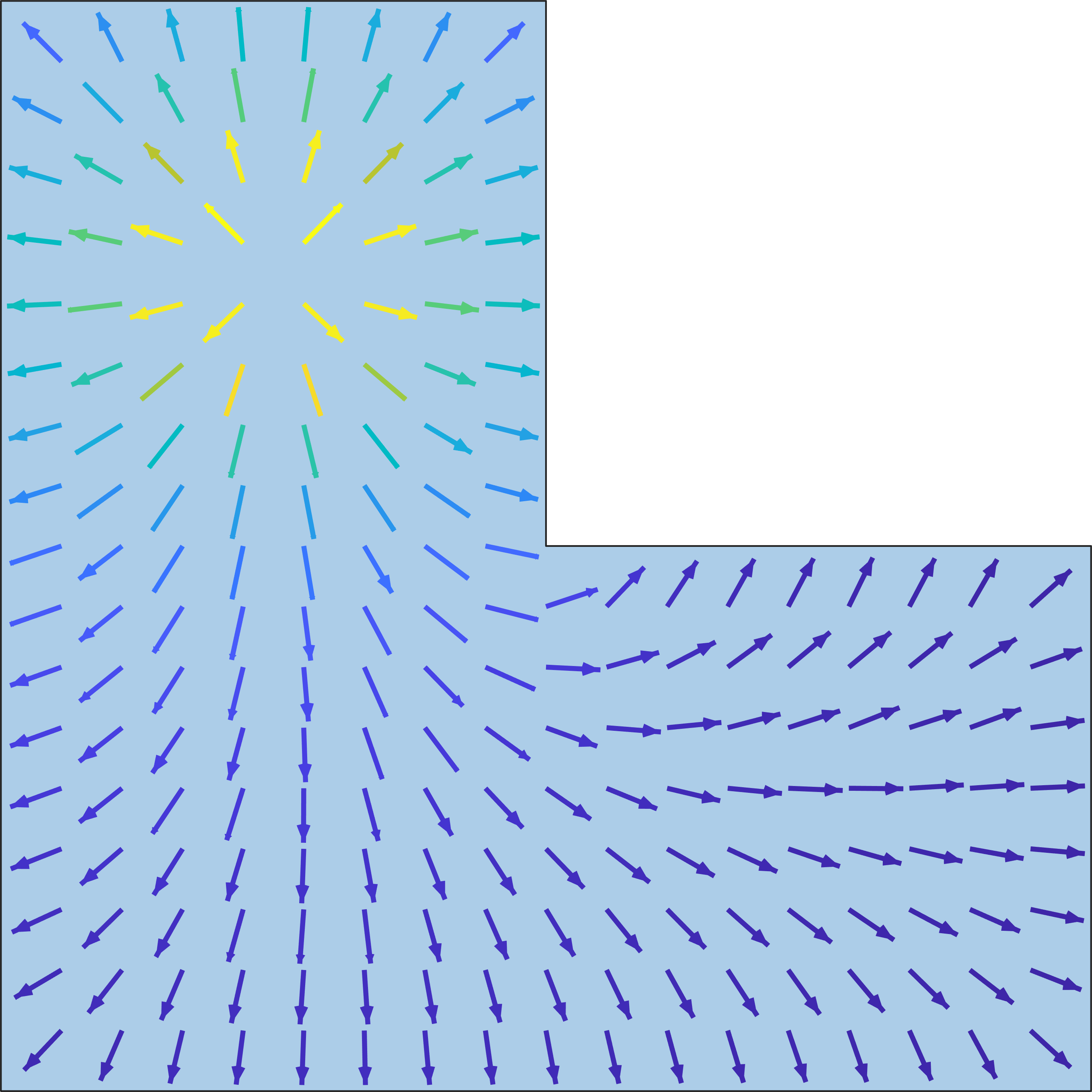};
    \end{axis}
\end{tikzpicture}

%% file: figures/Fig07a_PorousMedium_convergence_damping.tex
\begin{tikzpicture}[>=stealth]
    \begin{loglogaxis}[%
            width            = 5.5cm,%
            xlabel           = {cumulative ndof},%
            ylabel           = {Least-squares functional \(N^k_\ell\)},%
            ymax             = 0.35,%
            ymin             = 4e-5,%
            ymajorgrids      = true,%
            font             = \footnotesize,%
            grid style       = {%
                densely dotted,%
                semithick%
            },%
            legend style     = {%
                legend pos = south west,%
                font = \footnotesize%
            },%
        ]

        \addlegendimage{marker1}
        \addlegendentry{\(\delta = 1.0\phantom{0}\)}
        \addlegendimage{marker2}
        \addlegendentry{\(\delta = 0.5\phantom{0}\)}
        \addlegendimage{marker3}
        \addlegendentry{\(\delta = 0.1\phantom{0}\)}
        \addlegendimage{marker4}
        \addlegendentry{\(\delta = 0.05\)}
        \addlegendimage{marker5}
        \addlegendentry{\(\delta = 0.01\)}

        \pgfplotstableread[col sep=comma]{
       k,     ell,    ndof,    cost,   nElem,        eta,         mu,        res,    case,   maxGradU,       cond
       0,       0,     193,     193,      96,1.76237e-02,4.80764e-01,2.67715e-02,       Z,3.19183e-06,3.16582e+05
       1,       0,     193,     386,      96,1.76012e-02,2.01838e-02,2.42249e-02,       R,1.02740e-05,3.16582e+05
       1,       1,     225,     611,     112,1.48393e-02,4.81384e-01,2.57419e-02,       Z,1.29884e-05,1.00129e+06
       2,       0,     225,     836,     112,1.48714e-02,2.10175e-02,2.33022e-02,       R,2.93689e-05,1.00129e+06
       2,       1,     345,    1181,     172,1.10176e-02,2.32703e-02,2.17484e-02,       Z,3.04053e-05,1.76782e+06
       3,       0,     345,    1526,     172,1.10154e-02,1.87547e-02,1.96280e-02,       Z,5.17603e-05,1.76782e+06
       4,       0,     345,    1871,     172,1.10123e-02,1.62503e-02,1.79983e-02,       Z,7.53309e-05,1.76782e+06
       5,       0,     345,    2216,     172,1.10099e-02,1.42411e-02,1.67160e-02,       R,1.00032e-04,1.76782e+06
       5,       1,     441,    2657,     220,5.45785e-03,3.40434e-01,6.73712e-03,       Z,7.78977e-05,4.62117e+06
       6,       0,     441,    3098,     220,5.46326e-03,3.94374e-03,6.46451e-03,       Z,8.21723e-05,4.62117e+06
       7,       0,     441,    3539,     220,5.46831e-03,3.44950e-03,6.25401e-03,       Z,9.28305e-05,4.62117e+06
       8,       0,     441,    3980,     220,5.47298e-03,3.02830e-03,6.09144e-03,       Z,1.02934e-04,4.62117e+06
       9,       0,     441,    4421,     220,5.47731e-03,2.66763e-03,5.96583e-03,       Z,1.12453e-04,4.62117e+06
      10,       0,     441,    4862,     220,5.48132e-03,2.35748e-03,5.86866e-03,       Z,1.21378e-04,4.62117e+06
      11,       0,     441,    5303,     220,5.48505e-03,2.08971e-03,5.79336e-03,       Z,1.29716e-04,4.62117e+06
      12,       0,     441,    5744,     220,5.48851e-03,1.85763e-03,5.73490e-03,       R,1.37479e-04,4.62117e+06
      12,       1,     843,    6587,     421,4.45560e-03,1.90211e-01,6.24638e-03,       Z,1.49848e-04,1.69471e+07
      13,       0,     843,    7430,     421,4.45840e-03,4.37577e-03,5.91846e-03,       Z,1.69463e-04,1.69471e+07
      14,       0,     843,    8273,     421,4.46092e-03,3.89047e-03,5.65479e-03,       R,1.88423e-04,1.69471e+07
      14,       1,    1259,    9532,     629,3.58005e-03,3.04388e-02,5.32612e-03,       Z,1.88484e-04,2.35056e+07
      15,       0,    1259,   10791,     629,3.58390e-03,3.94024e-03,5.00618e-03,       Z,2.06975e-04,2.35056e+07
      16,       0,    1259,   12050,     629,3.58733e-03,3.49216e-03,4.74874e-03,       R,2.24759e-04,2.35056e+07
      16,       1,    1993,   14043,     996,2.94840e-03,1.12578e-01,4.83135e-03,       Z,2.27896e-04,5.58964e+07
      17,       0,    1993,   16036,     996,2.94947e-03,3.82671e-03,4.50515e-03,       R,2.48236e-04,5.58964e+07
      17,       1,    3137,   19173,    1568,2.40286e-03,8.64673e-02,4.73193e-03,       Z,2.50388e-04,9.86844e+07
      18,       0,    3137,   22310,    1568,2.40409e-03,4.07578e-03,4.35844e-03,       Z,2.72774e-04,9.86844e+07
      19,       0,    3137,   25447,    1568,2.40520e-03,3.63477e-03,4.04661e-03,       R,2.94294e-04,9.86844e+07
      19,       1,    4701,   30148,    2350,1.95592e-03,5.22247e-02,3.97627e-03,       Z,2.93042e-04,1.16060e+08
      20,       0,    4701,   34849,    2350,1.95643e-03,3.46170e-03,3.66593e-03,       Z,3.12382e-04,1.16060e+08
      21,       0,    4701,   39550,    2350,1.95693e-03,3.09995e-03,3.40428e-03,       R,3.30792e-04,1.16060e+08
      21,       1,    7251,   46801,    3625,1.55546e-03,5.63487e-02,3.07587e-03,       Z,3.29786e-04,2.35873e+08
      22,       0,    7251,   54052,    3625,1.55578e-03,2.65343e-03,2.84882e-03,       Z,3.46318e-04,2.35873e+08
      23,       0,    7251,   61303,    3625,1.55610e-03,2.38630e-03,2.65618e-03,       Z,3.61996e-04,2.35873e+08
      24,       0,    7251,   68554,    3625,1.55642e-03,2.15243e-03,2.49257e-03,       R,3.76845e-04,2.35873e+08
      24,       1,   10933,   79487,    5466,1.27120e-03,7.44506e-02,2.37738e-03,       Z,3.76063e-04,9.44373e+08
      25,       0,   10933,   90420,    5466,1.27137e-03,2.00888e-03,2.21570e-03,       R,3.89914e-04,9.44373e+08
      25,       1,   16415,  106835,    8207,1.03975e-03,4.99078e-02,2.18228e-03,       Z,3.90206e-04,2.21020e+09
      26,       0,   16415,  123250,    8207,1.03975e-03,1.91867e-03,2.02124e-03,       Z,4.03805e-04,2.21020e+09
      27,       0,   16415,  139665,    8207,1.03978e-03,1.73330e-03,1.88320e-03,       R,4.16659e-04,2.21020e+09
      27,       1,   25087,  164752,   12543,8.57826e-04,6.33657e-02,2.04078e-03,       Z,4.17623e-04,9.45823e+09
      28,       0,   25087,  189839,   12543,8.57987e-04,1.85166e-03,1.88476e-03,       Z,4.31062e-04,9.45823e+09
      29,       0,   25087,  214926,   12543,8.58152e-04,1.67807e-03,1.74936e-03,       R,4.45702e-04,9.45823e+09
      29,       1,   36521,  251447,   18260,7.05036e-04,3.08140e-02,1.65144e-03,       Z,4.45472e-04,3.60800e+10
      30,       0,   36521,  287968,   18260,7.05155e-04,1.49332e-03,1.52996e-03,       Z,4.59519e-04,3.60800e+10
      31,       0,   36521,  324489,   18260,7.05277e-04,1.35771e-03,1.42406e-03,       R,4.72768e-04,3.60800e+10
      31,       1,   53831,  378320,   26915,5.73596e-04,5.16900e-02,1.22805e-03,       Z,4.69312e-04,7.15173e+10
      32,       0,   53831,  432151,   26915,5.73684e-04,1.08581e-03,1.14363e-03,       Z,4.78888e-04,7.15173e+10
      33,       0,   53831,  485982,   26915,5.73776e-04,9.89278e-04,1.07014e-03,       R,4.90407e-04,7.15173e+10
      33,       1,   79435,  565417,   39717,4.73099e-04,2.34949e-02,1.06144e-03,       Z,4.90488e-04,1.42713e+11
      34,       0,   79435,  644852,   39717,4.73162e-04,9.50146e-04,9.87682e-04,       Z,5.01761e-04,1.42713e+11
      35,       0,   79435,  724287,   39717,4.73229e-04,8.66931e-04,9.23137e-04,       R,5.12375e-04,1.42713e+11
      35,       1,  115071,  839358,   57535,3.94404e-04,1.99087e-02,9.24222e-04,       Z,5.12524e-04,5.84801e+11
      36,       0,  115071,  954429,   57535,3.94424e-04,8.35833e-04,8.59752e-04,       Z,5.22655e-04,5.84801e+11
      37,       0,  115071, 1069500,   57535,3.94448e-04,7.63927e-04,8.03010e-04,       R,5.32187e-04,5.84801e+11
      37,       1,  167251, 1236751,   83625,3.28930e-04,3.77339e-02,8.32997e-04,       Z,5.32933e-04,2.28183e+12
      38,       0,  167251, 1404002,   83625,3.28966e-04,7.65288e-04,7.73685e-04,       R,5.42593e-04,2.28183e+12
      38,       1,  239061, 1643063,  119530,2.74886e-04,3.98282e-02,7.59972e-04,       Z,5.42707e-04,4.63888e+12
      39,       0,  239061, 1882124,  119530,2.74926e-04,7.08501e-04,7.04673e-04,       Z,5.51952e-04,4.63888e+12
      40,       0,  239061, 2121185,  119530,2.74968e-04,6.48812e-04,6.55580e-04,       R,5.61064e-04,4.63888e+12
      40,       1,  342357, 2463542,  171178,2.29749e-04,2.66512e-02,6.11468e-04,       Z,5.60682e-04,9.20755e+12
      41,       0,  342357, 2805899,  171178,2.29770e-04,5.66656e-04,5.68493e-04,       Z,5.68885e-04,9.20755e+12
      42,       0,  342357, 3148256,  171178,2.29792e-04,5.19981e-04,5.30251e-04,       R,5.76591e-04,9.20755e+12
      42,       1,  480575, 3628831,  240287,1.93994e-04,1.96721e-02,5.12206e-04,       Z,5.76857e-04,3.74795e+13
      43,       0,  480575, 4109406,  240287,1.94006e-04,4.74043e-04,4.76777e-04,       R,5.84355e-04,3.74795e+13
      43,       1,  676569, 4785975,  338284,1.63823e-04,2.15554e-02,4.94339e-04,       Z,5.84528e-04,7.39035e+13
      44,       0,  676569, 5462544,  338284,1.63837e-04,4.66400e-04,4.58777e-04,       Z,6.32554e-04,7.39035e+13
      45,       0,  676569, 6139113,  338284,1.63852e-04,4.28519e-04,4.26941e-04,       R,7.47702e-04,7.39035e+13
      45,       1,  942389, 7081502,  471194,1.38727e-04,1.34932e-02,4.26901e-04,       Z,7.51910e-04,1.47154e+14
      46,       0,  942389, 8023891,  471194,1.38741e-04,4.03728e-04,3.96564e-04,       R,8.69136e-04,1.47154e+14
      46,       1, 1325243, 9349134,  662621,1.17131e-04,1.53014e-02,3.84460e-04,       Z,9.42473e-04,5.99540e+14
        }\tableDataOne

        \pgfplotstableread[col sep=comma]{
       k,     ell,    ndof,    cost,   nElem,        eta,         mu,        res,    case,   maxGradU,       cond
       0,       0,     193,     193,      96,8.81184e-03,2.40382e-01,1.25869e-01,       Z,7.97958e-07,3.16582e+05
       1,       0,     193,     386,      96,8.80823e-03,1.20294e-01,6.56486e-02,       R,2.85170e-06,3.16582e+05
       1,       1,     225,     611,     112,7.42089e-03,2.08320e-01,1.10105e-01,       Z,3.50828e-06,1.00129e+06
       2,       0,     225,     836,     112,7.42656e-03,1.04448e-01,5.87041e-02,       R,8.51891e-06,1.00129e+06
       2,       1,     345,    1181,     172,5.50142e-03,1.04567e-01,5.83261e-02,       Z,8.78263e-06,1.76782e+06
       3,       0,     345,    1526,     172,5.50156e-03,5.29828e-02,3.49399e-02,       Z,1.58493e-05,1.76782e+06
       4,       0,     345,    1871,     172,5.50112e-03,2.78107e-02,2.56073e-02,       Z,2.42866e-05,1.76782e+06
       5,       0,     345,    2216,     172,5.50059e-03,1.60994e-02,2.21130e-02,       R,3.37815e-05,1.76782e+06
       5,       1,     441,    2657,     220,2.72555e-03,1.62692e-01,8.55830e-02,       Z,2.72969e-05,4.62117e+06
       6,       0,     441,    3098,     220,2.72697e-03,8.14419e-02,4.35584e-02,       Z,3.02761e-05,4.62117e+06
       7,       0,     441,    3539,     220,2.72837e-03,4.08489e-02,2.29862e-02,       Z,3.57057e-05,4.62117e+06
       8,       0,     441,    3980,     220,2.72974e-03,2.06110e-02,1.33819e-02,       Z,4.13742e-05,4.62117e+06
       9,       0,     441,    4421,     220,2.73105e-03,1.05935e-02,9.34337e-03,       Z,4.70802e-05,4.62117e+06
      10,       0,     441,    4862,     220,2.73233e-03,5.74648e-03,7.82759e-03,       Z,5.27861e-05,4.62117e+06
      11,       0,     441,    5303,     220,2.73355e-03,3.53269e-03,7.24263e-03,       Z,5.84612e-05,4.62117e+06
      12,       0,     441,    5744,     220,2.73474e-03,2.60368e-03,6.95828e-03,       R,6.40808e-05,4.62117e+06
      12,       1,     831,    6575,     415,2.25388e-03,9.45611e-02,4.96835e-02,       Z,6.81611e-05,1.75560e+07
      13,       0,     831,    7406,     415,2.25499e-03,4.73572e-02,2.57115e-02,       R,7.80482e-05,1.75560e+07
      13,       1,    1191,    8597,     595,1.83206e-03,4.73754e-02,2.56577e-02,       Z,7.80087e-05,1.75834e+07
      14,       0,    1191,    9788,     595,1.83353e-03,2.38535e-02,1.42835e-02,       Z,8.80091e-05,1.75834e+07
      15,       0,    1191,   10979,     595,1.83492e-03,1.22324e-02,9.37813e-03,       R,9.82556e-05,1.75834e+07
      15,       1,    1947,   12926,     973,1.49257e-03,6.13362e-02,3.26408e-02,       Z,9.94751e-05,4.87919e+07
      16,       0,    1947,   14873,     973,1.49310e-03,3.08012e-02,1.74196e-02,       Z,1.11093e-04,4.87919e+07
      17,       0,    1947,   16820,     973,1.49358e-03,1.56462e-02,1.05012e-02,       R,1.22766e-04,4.87919e+07
      17,       1,    3093,   19913,    1546,1.21026e-03,4.07212e-02,2.22175e-02,       Z,1.23203e-04,9.82524e+07
      18,       0,    3093,   23006,    1546,1.21059e-03,2.05560e-02,1.25422e-02,       Z,1.35530e-04,9.82524e+07
      19,       0,    3093,   26099,    1546,1.21091e-03,1.06313e-02,8.39866e-03,       R,1.48108e-04,9.82524e+07
      19,       1,    4671,   30770,    2335,9.73969e-04,3.04265e-02,1.70067e-02,       Z,1.47648e-04,1.39189e+08
      20,       0,    4671,   35441,    2335,9.74304e-04,1.54394e-02,1.00081e-02,       Z,1.59700e-04,1.39189e+08
      21,       0,    4671,   40112,    2335,9.74625e-04,8.11466e-03,7.09641e-03,       R,1.71657e-04,1.39189e+08
      21,       1,    7213,   47325,    3606,7.78799e-04,2.26476e-02,1.30106e-02,       Z,1.71390e-04,2.34423e+08
      22,       0,    7213,   54538,    3606,7.79047e-04,1.15617e-02,7.96363e-03,       Z,1.82990e-04,2.34423e+08
      23,       0,    7213,   61751,    3606,7.79284e-04,6.18361e-03,5.90341e-03,       Z,1.94479e-04,2.34423e+08
      24,       0,    7213,   68964,    3606,7.79510e-04,3.72183e-03,5.08890e-03,       R,2.05806e-04,2.34423e+08
      24,       1,   10807,   79771,    5403,6.35876e-04,3.84823e-02,2.05285e-02,       Z,2.05614e-04,9.45301e+08
      25,       0,   10807,   90578,    5403,6.36011e-04,1.93408e-02,1.09692e-02,       R,2.16564e-04,9.45301e+08
      25,       1,   16253,  106831,    8126,5.21902e-04,3.16848e-02,1.70455e-02,       Z,2.16686e-04,3.24634e+09
      26,       0,   16253,  123084,    8126,5.22015e-04,1.59550e-02,9.29456e-03,       Z,2.27652e-04,3.24634e+09
      27,       0,   16253,  139337,    8126,5.22126e-04,8.17769e-03,5.83815e-03,       R,2.38547e-04,3.24634e+09
      27,       1,   24949,  164286,   12474,4.30776e-04,3.18555e-02,1.70683e-02,       Z,2.39042e-04,8.83833e+09
      28,       0,   24949,  189235,   12474,4.30887e-04,1.60312e-02,9.25551e-03,       Z,2.50421e-04,8.83833e+09
      29,       0,   24949,  214184,   12474,4.30995e-04,8.20487e-03,5.76349e-03,       R,2.61644e-04,8.83833e+09
      29,       1,   36467,  250651,   18233,3.52861e-04,1.98306e-02,1.09826e-02,       Z,2.61598e-04,3.83249e+10
      30,       0,   36467,  287118,   18233,3.52935e-04,1.00509e-02,6.34651e-03,       Z,2.72565e-04,3.83249e+10
      31,       0,   36467,  323585,   18233,3.53007e-04,5.26237e-03,4.38698e-03,       R,2.83313e-04,3.83249e+10
      31,       1,   54179,  377764,   27089,2.85290e-04,2.69899e-02,1.43993e-02,       Z,2.83241e-04,5.26695e+10
      32,       0,   54179,  431943,   27089,2.85343e-04,1.35670e-02,7.67649e-03,       Z,2.93697e-04,5.26695e+10
      33,       0,   54179,  486122,   27089,2.85395e-04,6.90854e-03,4.59191e-03,       R,3.03932e-04,5.26695e+10
      33,       1,   79311,  565433,   39655,2.35321e-04,1.26547e-02,7.17704e-03,       Z,3.03785e-04,2.46169e+11
      34,       0,   79311,  644744,   39655,2.35371e-04,6.44986e-03,4.32480e-03,       Z,3.13651e-04,2.46169e+11
      35,       0,   79311,  724055,   39655,2.35421e-04,3.43732e-03,3.15811e-03,       R,3.23294e-04,2.46169e+11
      35,       1,  115209,  839264,   57604,1.96537e-04,9.41569e-03,5.54661e-03,       Z,3.23342e-04,6.28120e+11
      36,       0,  115209,  954473,   57604,1.96564e-04,4.84465e-03,3.54546e-03,       Z,3.32808e-04,6.28120e+11
      37,       0,  115209, 1069682,   57604,1.96591e-04,2.65678e-03,2.75647e-03,       R,3.42632e-04,6.28120e+11
      37,       1,  167369, 1237051,   83684,1.63939e-04,1.94159e-02,1.03741e-02,       Z,3.42774e-04,2.52729e+12
      38,       0,  167369, 1404420,   83684,1.63964e-04,9.76583e-03,5.57947e-03,       R,3.52671e-04,2.52729e+12
      38,       1,  239657, 1644077,  119828,1.36832e-04,2.23300e-02,1.18337e-02,       Z,3.52726e-04,2.52730e+12
      39,       0,  239657, 1883734,  119828,1.36855e-04,1.12117e-02,6.23694e-03,       Z,3.62441e-04,2.52730e+12
      40,       0,  239657, 2123391,  119828,1.36877e-04,5.69070e-03,3.64538e-03,       R,3.71920e-04,2.52730e+12
      40,       1,  342773, 2466164,  171386,1.14567e-04,1.34706e-02,7.31058e-03,       Z,3.71849e-04,1.03510e+13
      41,       0,  342773, 2808937,  171386,1.14582e-04,6.79857e-03,4.06365e-03,       Z,3.81025e-04,1.03510e+13
      42,       0,  342773, 3151710,  171386,1.14597e-04,3.51188e-03,2.64654e-03,       R,3.89971e-04,1.03510e+13
      42,       1,  480805, 3632515,  240402,9.67986e-05,1.01774e-02,5.62207e-03,       Z,3.89912e-04,4.14380e+13
      43,       0,  480805, 4113320,  240402,9.68086e-05,5.15730e-03,3.23892e-03,       R,3.98574e-04,4.14380e+13
      43,       1,  676187, 4789507,  338093,8.17629e-05,1.09048e-02,5.95661e-03,       Z,3.98691e-04,4.14528e+13
      44,       0,  676187, 5465694,  338093,8.17727e-05,5.51237e-03,3.36251e-03,       Z,4.07245e-04,4.14528e+13
      45,       0,  676187, 6141881,  338093,8.17824e-05,2.86367e-03,2.25005e-03,       R,4.15577e-04,4.14528e+13
      45,       1,  943257, 7085138,  471628,6.91625e-05,7.71207e-03,4.34930e-03,       Z,4.15778e-04,1.65809e+14
      46,       0,  943257, 8028395,  471628,6.91709e-05,3.92764e-03,2.60568e-03,       Z,4.87106e-04,1.65809e+14
      47,       0,  943257, 8971652,  471628,6.91791e-05,2.09000e-03,1.89439e-03,       R,5.65855e-04,1.65809e+14
      47,       1, 1325565,10297217,  662782,5.84207e-05,8.04217e-03,4.45580e-03,       Z,6.01168e-04,6.63058e+14
        }\tableDataOhFive

        \pgfplotstableread[col sep=comma]{
       k,     ell,    ndof,    cost,   nElem,        eta,         mu,        res,    case,   maxGradU,       cond
       0,       0,     193,     193,      96,1.76237e-03,4.80764e-02,2.25049e-01,       Z,3.19183e-08,3.16582e+05
       1,       0,     193,     386,      96,1.76233e-03,4.32692e-02,2.02625e-01,       R,1.24720e-07,3.16582e+05
       1,       1,     225,     611,     112,1.48430e-03,4.58733e-02,2.14735e-01,       Z,1.49959e-07,1.00129e+06
       2,       0,     225,     836,     112,1.48437e-03,4.12876e-02,1.93365e-01,       R,3.90762e-07,1.00129e+06
       2,       1,     345,    1181,     172,1.09950e-03,4.12997e-02,1.93347e-01,       Z,4.01499e-07,1.76782e+06
       3,       0,     345,    1526,     172,1.09951e-03,3.71729e-02,1.74133e-01,       Z,7.71457e-07,1.76782e+06
       4,       0,     345,    1871,     172,1.09951e-03,3.34602e-02,1.56878e-01,       Z,1.25255e-06,1.76782e+06
       5,       0,     345,    2216,     172,1.09950e-03,3.01204e-02,1.41388e-01,       R,1.83820e-06,1.76782e+06
       5,       1,     425,    2641,     212,5.55145e-04,2.49430e-02,1.17062e-01,       Z,1.53333e-06,4.62118e+06
       6,       0,     425,    3066,     212,5.55179e-04,2.24532e-02,1.05473e-01,       Z,1.83648e-06,4.62118e+06
       7,       0,     425,    3491,     212,5.55216e-04,2.02130e-02,9.50563e-02,       Z,2.23375e-06,4.62118e+06
       8,       0,     425,    3916,     212,5.55256e-04,1.81973e-02,8.56958e-02,       Z,2.72520e-06,4.62118e+06
       9,       0,     425,    4341,     212,5.55299e-04,1.63837e-02,7.72872e-02,       Z,3.25535e-06,4.62118e+06
      10,       0,     425,    4766,     212,5.55344e-04,1.47522e-02,6.97363e-02,       Z,3.82204e-06,4.62118e+06
      11,       0,     425,    5191,     212,5.55391e-04,1.32844e-02,6.29586e-02,       R,4.42322e-06,4.62118e+06
      11,       1,     793,    5984,     396,4.59946e-04,1.93529e-02,9.10933e-02,       Z,4.64205e-06,1.69456e+07
      12,       0,     793,    6777,     396,4.59979e-04,1.74249e-02,8.21608e-02,       Z,5.52732e-06,1.69456e+07
      13,       0,     793,    7570,     396,4.60013e-04,1.56906e-02,7.41417e-02,       R,6.47447e-06,1.69456e+07
      13,       1,    1173,    8743,     586,3.67860e-04,1.58786e-02,7.49974e-02,       Z,6.47646e-06,2.34342e+07
      14,       0,    1173,    9916,     586,3.67908e-04,1.43000e-02,6.77101e-02,       Z,7.48465e-06,2.34342e+07
      15,       0,    1173,   11089,     586,3.67959e-04,1.28803e-02,6.11740e-02,       R,8.54891e-06,2.34342e+07
      15,       1,    1937,   13026,     968,2.95333e-04,1.57168e-02,7.42697e-02,       Z,8.59140e-06,4.56264e+07
      16,       0,    1937,   14963,     968,2.95367e-04,1.41549e-02,6.70603e-02,       Z,9.75735e-06,4.56264e+07
      17,       0,    1937,   16900,     968,2.95402e-04,1.27503e-02,6.05941e-02,       R,1.09796e-05,4.56264e+07
      17,       1,    3039,   19939,    1519,2.42578e-04,1.41015e-02,6.68286e-02,       Z,1.10771e-05,8.88720e+07
      18,       0,    3039,   22978,    1519,2.42596e-04,1.27030e-02,6.03976e-02,       Z,1.24607e-05,8.88720e+07
      19,       0,    3039,   26017,    1519,2.42613e-04,1.14455e-02,5.46351e-02,       R,1.39051e-05,8.88720e+07
      19,       1,    4521,   30538,    2260,1.97722e-04,1.21547e-02,5.78840e-02,       Z,1.38676e-05,1.43302e+08
      20,       0,    4521,   35059,    2260,1.97744e-04,1.09527e-02,5.23818e-02,       Z,1.53283e-05,1.43302e+08
      21,       0,    4521,   39580,    2260,1.97766e-04,9.87227e-03,4.74567e-02,       R,1.68402e-05,1.43302e+08
      21,       1,    6945,   46525,    3472,1.56722e-04,1.00890e-02,4.84291e-02,       Z,1.68098e-05,2.32018e+08
      22,       0,    6945,   53470,    3472,1.56739e-04,9.09560e-03,4.39077e-02,       Z,1.83379e-05,2.32018e+08
      23,       0,    6945,   60415,    3472,1.56756e-04,8.20294e-03,3.98658e-02,       R,1.99105e-05,2.32018e+08
      23,       1,   10431,   70846,    5215,1.27658e-04,9.41393e-03,4.53389e-02,       Z,1.98943e-05,9.45245e+08
      24,       0,   10431,   81277,    5215,1.27669e-04,8.48862e-03,4.11357e-02,       Z,2.14917e-05,9.45245e+08
      25,       0,   10431,   91708,    5215,1.27680e-04,7.65728e-03,3.73799e-02,       R,2.31283e-05,9.45245e+08
      25,       1,   15791,  107499,    7895,1.04806e-04,8.64992e-03,4.18517e-02,       Z,2.31446e-05,2.21199e+09
      26,       0,   15791,  123290,    7895,1.04818e-04,7.80211e-03,3.80168e-02,       Z,2.48358e-05,2.21199e+09
      27,       0,   15791,  139081,    7895,1.04830e-04,7.04061e-03,3.45932e-02,       R,2.66045e-05,2.21199e+09
      27,       1,   24063,  163144,   12031,8.67408e-05,8.59808e-03,4.15953e-02,       Z,2.66312e-05,9.48764e+09
      28,       0,   24063,  187207,   12031,8.67524e-05,7.75548e-03,3.77863e-02,       Z,2.84721e-05,9.48764e+09
      29,       0,   24063,  211270,   12031,8.67640e-05,6.99870e-03,3.43864e-02,       R,3.03501e-05,9.48764e+09
      29,       1,   35599,  246869,   17799,7.05851e-05,7.11174e-03,3.48857e-02,       Z,3.03547e-05,3.82328e+10
      30,       0,   35599,  282468,   17799,7.05926e-05,6.42051e-03,3.17955e-02,       Z,3.22729e-05,3.82328e+10
      31,       0,   35599,  318067,   17799,7.06000e-05,5.80012e-03,2.90439e-02,       R,3.42250e-05,3.82328e+10
      31,       1,   53283,  371350,   26641,5.70104e-05,6.31538e-03,3.12946e-02,       Z,3.42235e-05,1.44093e+11
      32,       0,   53283,  424633,   26641,5.70175e-05,5.70484e-03,2.85780e-02,       Z,3.62066e-05,1.44093e+11
      33,       0,   53283,  477916,   26641,5.70246e-05,5.15706e-03,2.61615e-02,       R,3.82208e-05,1.44093e+11
      33,       1,   78885,  556801,   39442,4.69217e-05,5.43030e-03,2.73456e-02,       Z,3.82110e-05,2.85920e+11
      34,       0,   78885,  635686,   39442,4.69279e-05,4.91070e-03,2.50634e-02,       Z,4.02446e-05,2.85920e+11
      35,       0,   78885,  714571,   39442,4.69341e-05,4.44493e-03,2.30385e-02,       R,4.23060e-05,2.85920e+11
      35,       1,  114555,  829126,   57277,3.91088e-05,4.80219e-03,2.45581e-02,       Z,4.23113e-05,1.02124e+12
      36,       0,  114555,  943681,   57277,3.91129e-05,4.34735e-03,2.25841e-02,       Z,4.44046e-05,1.02124e+12
      37,       0,  114555, 1058236,   57277,3.91171e-05,3.94000e-03,2.08365e-02,       R,4.65234e-05,1.02124e+12
      37,       1,  166345, 1224581,   83172,3.26509e-05,4.84441e-03,2.46601e-02,       Z,4.65464e-05,2.30460e+12
      38,       0,  166345, 1390926,   83172,3.26545e-05,4.38421e-03,2.26565e-02,       R,4.87134e-05,2.30460e+12
      38,       1,  237199, 1628125,  118599,2.72645e-05,5.12111e-03,2.58167e-02,       Z,4.87167e-05,9.35431e+12
      39,       0,  237199, 1865324,  118599,2.72678e-05,4.63141e-03,2.36667e-02,       Z,5.09110e-05,9.35431e+12
      40,       0,  237199, 2102523,  118599,2.72710e-05,4.19247e-03,2.17593e-02,       R,5.31282e-05,9.35431e+12
      40,       1,  339675, 2442198,  169837,2.28137e-05,4.31563e-03,2.22775e-02,       Z,5.31190e-05,1.84448e+13
      41,       0,  339675, 2781873,  169837,2.28164e-05,3.90931e-03,2.05234e-02,       Z,5.53485e-05,1.84448e+13
      42,       0,  339675, 3121548,  169837,2.28190e-05,3.54557e-03,1.89721e-02,       R,5.75984e-05,1.84448e+13
      42,       1,  479563, 3601111,  239781,1.92130e-05,3.75873e-03,1.98479e-02,       Z,5.75996e-05,7.49146e+13
      43,       0,  479563, 4080674,  239781,1.92150e-05,3.41038e-03,1.83670e-02,       R,5.98701e-05,7.49146e+13
      43,       1,  676221, 4756895,  338110,1.62144e-05,3.46990e-03,1.86108e-02,       Z,5.98766e-05,1.47772e+14
      44,       0,  676221, 5433116,  338110,1.62161e-05,3.15209e-03,1.72753e-02,       Z,6.21722e-05,1.47772e+14
      45,       0,  676221, 6109337,  338110,1.62178e-05,2.86822e-03,1.60992e-02,       R,6.44856e-05,1.47772e+14
      45,       1,  940347, 7049684,  470173,1.37302e-05,2.96878e-03,1.64935e-02,       Z,6.44879e-05,2.94285e+14
      46,       0,  940347, 7990031,  470173,1.37317e-05,2.70438e-03,1.54077e-02,       Z,6.68204e-05,2.94285e+14
      47,       0,  940347, 8930378,  470173,1.37332e-05,2.46871e-03,1.44545e-02,       R,6.91688e-05,2.94285e+14
      47,       1, 1318035,10248413,  659017,1.15952e-05,2.64572e-03,1.51227e-02,       Z,6.91736e-05,1.19817e+15
        }\tableDataOhOne

        \pgfplotstableread[col sep=comma]{
       k,     ell,    ndof,    cost,   nElem,        eta,         mu,        res,    case,   maxGradU,       cond
       0,       0,     193,     193,      96,8.81184e-04,2.40382e-02,2.37520e-01,       Z,7.97958e-09,3.16582e+05
       1,       0,     193,     386,      96,8.81179e-04,2.28363e-02,2.25672e-01,       R,3.15449e-08,3.16582e+05
       1,       1,     225,     611,     112,7.42155e-04,2.34655e-02,2.31827e-01,       Z,3.78157e-08,1.00129e+06
       2,       0,     225,     836,     112,7.42164e-04,2.22924e-02,2.20266e-01,       R,9.95566e-08,1.00129e+06
       2,       1,     345,    1181,     172,5.49730e-04,2.22980e-02,2.20259e-01,       Z,1.02249e-07,1.76782e+06
       3,       0,     345,    1526,     172,5.49732e-04,2.11834e-02,2.09277e-01,       Z,1.98424e-07,1.76782e+06
       4,       0,     345,    1871,     172,5.49732e-04,2.01248e-02,1.98853e-01,       Z,3.25302e-07,1.76782e+06
       5,       0,     345,    2216,     172,5.49731e-04,1.91193e-02,1.88957e-01,       R,4.81898e-07,1.76782e+06
       5,       1,     425,    2641,     212,2.77555e-04,1.39817e-02,1.38184e-01,       Z,4.03780e-07,4.62118e+06
       6,       0,     425,    3066,     212,2.77560e-04,1.32833e-02,1.31307e-01,       Z,4.89816e-07,4.62118e+06
       7,       0,     425,    3491,     212,2.77565e-04,1.26198e-02,1.24776e-01,       Z,5.99655e-07,4.62118e+06
       8,       0,     425,    3916,     212,2.77572e-04,1.19895e-02,1.18576e-01,       Z,7.38535e-07,4.62118e+06
       9,       0,     425,    4341,     212,2.77579e-04,1.13908e-02,1.12689e-01,       Z,8.90258e-07,4.62118e+06
      10,       0,     425,    4766,     212,2.77587e-04,1.08222e-02,1.07100e-01,       Z,1.05444e-06,4.62118e+06
      11,       0,     425,    5191,     212,2.77595e-04,1.02820e-02,1.01794e-01,       R,1.23069e-06,4.62118e+06
      11,       1,     793,    5984,     396,2.29908e-04,1.23760e-02,1.22419e-01,       Z,1.28838e-06,1.69456e+07
      12,       0,     793,    6777,     396,2.29913e-04,1.17583e-02,1.16349e-01,       Z,1.54353e-06,1.69456e+07
      13,       0,     793,    7570,     396,2.29918e-04,1.11716e-02,1.10586e-01,       R,1.81929e-06,1.69456e+07
      13,       1,    1173,    8743,     586,1.83848e-04,1.12430e-02,1.11280e-01,       Z,1.81981e-06,2.34342e+07
      14,       0,    1173,    9916,     586,1.83856e-04,1.06822e-02,1.05775e-01,       Z,2.11620e-06,2.34342e+07
      15,       0,    1173,   11089,     586,1.83865e-04,1.01495e-02,1.00548e-01,       R,2.43213e-06,2.34342e+07
      15,       1,    1937,   13026,     968,1.47551e-04,1.09070e-02,1.08000e-01,       Z,2.44335e-06,4.56264e+07
      16,       0,    1937,   14963,     968,1.47557e-04,1.03631e-02,1.02663e-01,       Z,2.79121e-06,4.56264e+07
      17,       0,    1937,   16900,     968,1.47564e-04,9.84651e-03,9.75986e-02,       R,3.15921e-06,4.56264e+07
      17,       1,    3031,   19931,    1515,1.21361e-04,1.03164e-02,1.02221e-01,       Z,3.18526e-06,8.88720e+07
      18,       0,    3031,   22962,    1515,1.21364e-04,9.80234e-03,9.71826e-02,       Z,3.60199e-06,8.88720e+07
      19,       0,    3031,   25993,    1515,1.21367e-04,9.31408e-03,9.24019e-02,       R,4.04078e-06,8.88720e+07
      19,       1,    4543,   30536,    2271,9.84201e-05,9.45772e-03,9.38107e-02,       Z,4.03081e-06,1.43362e+08
      20,       0,    4543,   35079,    2271,9.84242e-05,8.98681e-03,8.92024e-02,       Z,4.47986e-06,1.43362e+08
      21,       0,    4543,   39622,    2271,9.84284e-05,8.53958e-03,8.48301e-02,       R,4.94862e-06,1.43362e+08
      21,       1,    6999,   46621,    3499,7.79492e-05,8.41880e-03,8.36400e-02,       Z,4.93979e-06,2.32018e+08
      22,       0,    6999,   53620,    3499,7.79528e-05,8.00008e-03,7.95495e-02,       Z,5.41816e-06,2.32018e+08
      23,       0,    6999,   60619,    3499,7.79565e-05,7.60244e-03,7.56692e-02,       Z,5.91451e-06,2.32018e+08
      24,       0,    6999,   67618,    3499,7.79602e-05,7.22483e-03,7.19886e-02,       R,6.42839e-06,2.32018e+08
      24,       1,   10543,   78161,    5271,6.34808e-05,7.50125e-03,7.46923e-02,       Z,6.42660e-06,9.45245e+08
      25,       0,   10543,   88704,    5271,6.34830e-05,7.12875e-03,7.10623e-02,       R,6.95559e-06,9.45245e+08
      25,       1,   15879,  104583,    7939,5.22758e-05,7.39405e-03,7.36581e-02,       Z,6.96113e-06,2.21159e+09
      26,       0,   15879,  120462,    7939,5.22784e-05,7.02701e-03,7.00834e-02,       Z,7.51252e-06,2.21159e+09
      27,       0,   15879,  136341,    7939,5.22810e-05,6.67848e-03,6.66934e-02,       R,8.08581e-06,2.21159e+09
      27,       1,   24321,  160662,   12160,4.29677e-05,6.99191e-03,6.97560e-02,       Z,8.09330e-06,9.48207e+09
      28,       0,   24321,  184983,   12160,4.29702e-05,6.64522e-03,6.63856e-02,       Z,8.69334e-06,9.48207e+09
      29,       0,   24321,  209304,   12160,4.29727e-05,6.31602e-03,6.31900e-02,       R,9.31002e-06,9.48207e+09
      29,       1,   36135,  245439,   18067,3.49559e-05,6.25909e-03,6.26337e-02,       Z,9.31348e-06,3.82670e+10
      30,       0,   36135,  281574,   18067,3.49574e-05,5.94938e-03,5.96323e-02,       Z,9.95006e-06,3.82670e+10
      31,       0,   36135,  317709,   18067,3.49590e-05,5.65533e-03,5.67875e-02,       R,1.06027e-05,3.82670e+10
      31,       1,   53905,  371614,   26952,2.82733e-05,5.62116e-03,5.64528e-02,       Z,1.05940e-05,1.44262e+11
      32,       0,   53905,  425519,   26952,2.82749e-05,5.34366e-03,5.37718e-02,       Z,1.12532e-05,1.44262e+11
      33,       0,   53905,  479424,   26952,2.82765e-05,5.08022e-03,5.12316e-02,       R,1.19272e-05,1.44262e+11
      33,       1,   79499,  558923,   39749,2.33475e-05,5.16104e-03,5.20136e-02,       Z,1.19319e-05,2.88853e+11
      34,       0,   79499,  638422,   39749,2.33490e-05,4.90687e-03,4.95666e-02,       Z,1.26253e-05,2.88853e+11
      35,       0,   79499,  717921,   39749,2.33504e-05,4.66560e-03,4.72491e-02,       R,1.33331e-05,2.88853e+11
      35,       1,  115401,  833322,   57700,1.94574e-05,4.75791e-03,4.81384e-02,       Z,1.33348e-05,1.15561e+12
      36,       0,  115401,  948723,   57700,1.94585e-05,4.52422e-03,4.58971e-02,       Z,1.40584e-05,1.15561e+12
      37,       0,  115401, 1064124,   57700,1.94595e-05,4.30241e-03,4.37751e-02,       R,1.47957e-05,1.15561e+12
      37,       1,  167185, 1231309,   83592,1.62432e-05,4.44826e-03,4.51747e-02,       Z,1.48006e-05,2.30191e+12
      38,       0,  167185, 1398494,   83592,1.62440e-05,4.23033e-03,4.30918e-02,       R,1.55564e-05,2.30191e+12
      38,       1,  237693, 1636187,  118846,1.35944e-05,4.36624e-03,4.43935e-02,       Z,1.55550e-05,9.35655e+12
      39,       0,  237693, 1873880,  118846,1.35952e-05,4.15248e-03,4.23518e-02,       Z,1.63235e-05,9.35655e+12
      40,       0,  237693, 2111573,  118846,1.35960e-05,3.94962e-03,4.04198e-02,       R,1.71079e-05,9.35655e+12
      40,       1,  339555, 2451128,  169777,1.13829e-05,3.94462e-03,4.03714e-02,       Z,1.71069e-05,1.84481e+13
      41,       0,  339555, 2790683,  169777,1.13835e-05,3.75236e-03,3.85456e-02,       Z,1.79031e-05,1.84481e+13
      42,       0,  339555, 3130238,  169777,1.13842e-05,3.56993e-03,3.68190e-02,       R,1.87119e-05,1.84481e+13
      42,       1,  476345, 3606583,  238172,9.61315e-06,3.58738e-03,3.69835e-02,       Z,1.87106e-05,7.48836e+13
      43,       0,  476345, 4082928,  238172,9.61368e-06,3.41340e-03,3.53416e-02,       R,1.95304e-05,7.48836e+13
      43,       1,  668275, 4751203,  334137,8.12870e-06,3.42239e-03,3.54267e-02,       Z,1.95350e-05,1.47800e+14
      44,       0,  668275, 5419478,  334137,8.12910e-06,3.25689e-03,3.38707e-02,       Z,2.03715e-05,1.47800e+14
      45,       0,  668275, 6087753,  334137,8.12959e-06,3.09992e-03,3.24008e-02,       R,2.12198e-05,1.47800e+14
      45,       1,  929573, 7017326,  464786,6.88098e-06,3.11356e-03,3.25282e-02,       Z,2.12204e-05,5.99521e+14
      46,       0,  929573, 7946899,  464786,6.88124e-06,2.96398e-03,3.11329e-02,       Z,2.20810e-05,5.99521e+14
      47,       0,  929573, 8876472,  464786,6.88147e-06,2.82214e-03,2.98157e-02,       R,2.29529e-05,5.99521e+14
      47,       1, 1298345,10174817,  649172,5.81940e-06,2.83283e-03,2.99141e-02,       Z,2.29519e-05,1.18131e+15
        }\tableDataOhOhFive

        \pgfplotstableread[col sep=comma]{
       k,     ell,    ndof,    cost,   nElem,        eta,         mu,        res,    case,   maxGradU,       cond
       0,       0,     193,     193,      96,1.76237e-04,4.80764e-03,2.47503e-01,       Z,3.19183e-10,3.16582e+05
       1,       0,     193,     386,      96,1.76237e-04,4.75956e-03,2.45032e-01,       R,1.27372e-09,3.16582e+05
       1,       1,     225,     611,     112,1.48431e-04,4.78473e-03,2.46277e-01,       Z,1.52324e-09,1.00129e+06
       2,       0,     225,     836,     112,1.48431e-04,4.73689e-03,2.43817e-01,       R,4.04461e-09,1.00129e+06
       2,       1,     345,    1181,     172,1.09945e-04,4.73794e-03,2.43816e-01,       Z,4.15258e-09,1.76782e+06
       3,       0,     345,    1526,     172,1.09945e-04,4.69056e-03,2.41380e-01,       Z,8.12738e-09,1.76782e+06
       4,       0,     345,    1871,     172,1.09945e-04,4.64366e-03,2.38969e-01,       Z,1.34383e-08,1.76782e+06
       5,       0,     345,    2216,     172,1.09945e-04,4.59723e-03,2.36582e-01,       R,2.00760e-08,1.76782e+06
       5,       1,     425,    2641,     212,5.55098e-05,3.25479e-03,1.67475e-01,       Z,1.68848e-08,4.62118e+06
       6,       0,     425,    3066,     212,5.55098e-05,3.22225e-03,1.65802e-01,       Z,2.07167e-08,4.62118e+06
       7,       0,     425,    3491,     212,5.55099e-05,3.19003e-03,1.64146e-01,       Z,2.55195e-08,4.62118e+06
       8,       0,     425,    3916,     212,5.55099e-05,3.15814e-03,1.62506e-01,       Z,3.17081e-08,4.62118e+06
       9,       0,     425,    4341,     212,5.55100e-05,3.12656e-03,1.60883e-01,       Z,3.85519e-08,4.62118e+06
      10,       0,     425,    4766,     212,5.55101e-05,3.09530e-03,1.59277e-01,       Z,4.60463e-08,4.62118e+06
      11,       0,     425,    5191,     212,5.55102e-05,3.06436e-03,1.57686e-01,       R,5.41869e-08,4.62118e+06
      11,       1,     793,    5984,     396,4.59763e-05,3.38604e-03,1.74224e-01,       Z,5.66067e-08,1.69456e+07
      12,       0,     793,    6777,     396,4.59764e-05,3.35218e-03,1.72484e-01,       Z,6.82445e-08,1.69456e+07
      13,       0,     793,    7570,     396,4.59765e-05,3.31867e-03,1.70762e-01,       R,8.09526e-08,1.69456e+07
      13,       1,    1173,    8743,     586,3.67625e-05,3.33115e-03,1.71398e-01,       Z,8.09747e-08,2.34342e+07
      14,       0,    1173,    9916,     586,3.67626e-05,3.29785e-03,1.69686e-01,       Z,9.47721e-08,2.34342e+07
      15,       0,    1173,   11089,     586,3.67627e-05,3.26488e-03,1.67992e-01,       R,1.09630e-07,2.34342e+07
      15,       1,    1937,   13026,     968,2.94997e-05,3.33877e-03,1.71789e-01,       Z,1.10103e-07,4.56264e+07
      16,       0,    1937,   14963,     968,2.94998e-05,3.30539e-03,1.70073e-01,       Z,1.26565e-07,4.56264e+07
      17,       0,    1937,   16900,     968,2.94998e-05,3.27235e-03,1.68376e-01,       R,1.44150e-07,4.56264e+07
      17,       1,    3031,   19931,    1515,2.42772e-05,3.35761e-03,1.72759e-01,       Z,1.45257e-07,9.00319e+07
      18,       0,    3031,   22962,    1515,2.42772e-05,3.32405e-03,1.71034e-01,       Z,1.65211e-07,9.00319e+07
      19,       0,    3031,   25993,    1515,2.42772e-05,3.29083e-03,1.69327e-01,       R,1.86419e-07,9.00319e+07
      19,       1,    4533,   30526,    2266,1.96421e-05,3.29426e-03,1.69502e-01,       Z,1.85971e-07,1.44389e+08
      20,       0,    4533,   35059,    2266,1.96421e-05,3.26133e-03,1.67810e-01,       Z,2.07914e-07,1.44389e+08
      21,       0,    4533,   39592,    2266,1.96422e-05,3.22874e-03,1.66136e-01,       R,2.31034e-07,1.44389e+08
      21,       1,    6999,   46591,    3499,1.56020e-05,3.18818e-03,1.64048e-01,       Z,2.30695e-07,2.32063e+08
      22,       0,    6999,   53590,    3499,1.56020e-05,3.15631e-03,1.62411e-01,       Z,2.54610e-07,2.32063e+08
      23,       0,    6999,   60589,    3499,1.56021e-05,3.12477e-03,1.60791e-01,       Z,2.79650e-07,2.32063e+08
      24,       0,    6999,   67588,    3499,1.56021e-05,3.09354e-03,1.59186e-01,       R,3.05810e-07,2.32063e+08
      24,       1,   10603,   78191,    5301,1.26620e-05,3.10008e-03,1.59522e-01,       Z,3.05741e-07,9.45245e+08
      25,       0,   10603,   88794,    5301,1.26620e-05,3.06910e-03,1.57931e-01,       R,3.32936e-07,9.45245e+08
      25,       1,   15905,  104699,    7952,1.04405e-05,3.10097e-03,1.59570e-01,       Z,3.33236e-07,2.22476e+09
      26,       0,   15905,  120604,    7952,1.04406e-05,3.06998e-03,1.57978e-01,       Z,3.61852e-07,2.22476e+09
      27,       0,   15905,  136509,    7952,1.04406e-05,3.03930e-03,1.56402e-01,       R,3.91584e-07,2.22476e+09
      27,       1,   24451,  160960,   12225,8.55582e-06,3.07098e-03,1.58031e-01,       Z,3.91904e-07,9.51411e+09
      28,       0,   24451,  185411,   12225,8.55585e-06,3.04029e-03,1.56455e-01,       Z,4.23289e-07,9.51411e+09
      29,       0,   24451,  209862,   12225,8.55588e-06,3.00991e-03,1.54895e-01,       R,4.55933e-07,9.51411e+09
      29,       1,   36343,  246205,   18171,6.98109e-06,2.99567e-03,1.54163e-01,       Z,4.56087e-07,3.82810e+10
      30,       0,   36343,  282548,   18171,6.98110e-06,2.96573e-03,1.52626e-01,       Z,4.90037e-07,3.82810e+10
      31,       0,   36343,  318891,   18171,6.98112e-06,2.93610e-03,1.51104e-01,       R,5.25134e-07,3.82810e+10
      31,       1,   53843,  372734,   26921,5.65181e-06,2.90787e-03,1.49652e-01,       Z,5.24715e-07,1.41778e+11
      32,       0,   53843,  426577,   26921,5.65183e-06,2.87881e-03,1.48160e-01,       Z,5.60505e-07,1.41778e+11
      33,       0,   53843,  480420,   26921,5.65185e-06,2.85005e-03,1.46684e-01,       R,5.97400e-07,1.41778e+11
      33,       1,   79619,  560039,   39809,4.65073e-06,2.85835e-03,1.47111e-01,       Z,5.97627e-07,2.85412e+11
      34,       0,   79619,  639658,   39809,4.65075e-06,2.82979e-03,1.45645e-01,       Z,6.35860e-07,2.85412e+11
      35,       0,   79619,  719277,   39809,4.65076e-06,2.80152e-03,1.44194e-01,       R,6.75199e-07,2.85412e+11
      35,       1,  116469,  835746,   58234,3.87207e-06,2.80769e-03,1.44511e-01,       Z,6.75273e-07,1.16716e+12
      36,       0,  116469,  952215,   58234,3.87209e-06,2.77964e-03,1.43072e-01,       Z,7.15788e-07,1.16716e+12
      37,       0,  116469, 1068684,   58234,3.87208e-06,2.75187e-03,1.41647e-01,       R,7.57399e-07,1.16716e+12
      37,       1,  168977, 1237661,   84488,3.23085e-06,2.76311e-03,1.42225e-01,       Z,7.57632e-07,2.30591e+12
      38,       0,  168977, 1406638,   84488,3.23082e-06,2.73551e-03,1.40808e-01,       R,8.00626e-07,2.30591e+12
      38,       1,  240001, 1646639,  120000,2.70242e-06,2.73548e-03,1.40807e-01,       Z,8.00746e-07,9.35548e+12
      39,       0,  240001, 1886640,  120000,2.70237e-06,2.70816e-03,1.39405e-01,       Z,8.45015e-07,9.35548e+12
      40,       0,  240001, 2126641,  120000,2.70230e-06,2.68111e-03,1.38017e-01,       R,8.90387e-07,9.35548e+12
      40,       1,  343197, 2469838,  171598,2.25865e-06,2.67605e-03,1.37757e-01,       Z,8.90257e-07,1.84459e+13
      41,       0,  343197, 2813035,  171598,2.25875e-06,2.64932e-03,1.36386e-01,       Z,9.36587e-07,1.84459e+13
      42,       0,  343197, 3156232,  171598,2.25862e-06,2.62286e-03,1.35029e-01,       R,9.84000e-07,1.84459e+13
      42,       1,  481883, 3638115,  240941,1.90560e-06,2.62069e-03,1.34917e-01,       Z,9.83921e-07,7.49158e+13
      43,       0,  481883, 4119998,  240941,1.90552e-06,2.59452e-03,1.33575e-01,       R,1.03233e-06,7.49158e+13
      43,       1,  676799, 4796797,  338399,1.60884e-06,2.59805e-03,1.33756e-01,       Z,1.03242e-06,1.47794e+14
      44,       0,  676799, 5473596,  338399,1.60860e-06,2.57211e-03,1.32426e-01,       Z,1.08200e-06,1.47794e+14
      45,       0,  676799, 6150395,  338399,1.60853e-06,2.54642e-03,1.31109e-01,       R,1.13263e-06,1.47794e+14
      45,       1,  940019, 7090414,  470009,1.35673e-06,2.54669e-03,1.31122e-01,       Z,1.13269e-06,5.96463e+14
      46,       0,  940019, 8030433,  470009,1.35649e-06,2.52126e-03,1.29818e-01,       Z,1.18445e-06,5.96463e+14
      47,       0,  940019, 8970452,  470009,1.35645e-06,2.49609e-03,1.28528e-01,       R,1.23726e-06,5.96463e+14
      47,       1, 1314507,10284959,  657253,1.13687e-06,2.49423e-03,1.28433e-01,       Z,1.23721e-06,1.18113e+15
        }\tableDataOhOhOne

        %
        %

        \addplot+ [marker1, adaptive, forget plot]
        table [col sep=comma, x=cumulativeNdof, y=res] {\tableDataOne};

        \addplot+ [marker2, adaptive, forget plot]
        table [col sep=comma, x=cumulativeNdof, y=res] {\tableDataOhFive};

        \addplot+ [marker3, adaptive, forget plot]
        table [col sep=comma, x=cumulativeNdof, y=res] {\tableDataOhOne};

        \addplot+ [marker4, adaptive, forget plot]
        table [col sep=comma, x=cumulativeNdof, y=res] {\tableDataOhOhFive};

        \addplot+ [marker5, adaptive, forget plot]
        table [col sep=comma, x=cumulativeNdof, y=res] {\tableDataOhOhOne};

        \drawslopetriangle[ST1]{0.5}{2e5}{4e-4} 
        \drawswappedslopetriangle[ST2]{0.2}{6e6}{6.5e-2} 
    \end{loglogaxis}
\end{tikzpicture}

%% file: figures/Fig07b_PorousMedium_gradient_norm_damping.tex
\begin{tikzpicture}[>=stealth]
    \begin{loglogaxis}[%
            width            = 5.5cm,%
            xlabel           = {cumulative ndof},%
            ylabel           = {\(\Vert \nabla u^k_\ell \Vert_{L^\infty(\Omega)}\)},%
            ymajorgrids      = true,%
            font             = \footnotesize,%
            grid style       = {%
                densely dotted,%
                semithick%
            },%
            legend style     = {%
                legend pos = south east,%
                font = \footnotesize%
            },%
        ]

        \addlegendimage{marker1}
        \addlegendentry{\(\delta = 1.0\phantom{0}\)}
        \addlegendimage{marker2}
        \addlegendentry{\(\delta = 0.5\phantom{0}\)}
        \addlegendimage{marker3}
        \addlegendentry{\(\delta = 0.1\phantom{0}\)}
        \addlegendimage{marker4}
        \addlegendentry{\(\delta = 0.05\)}
        \addlegendimage{marker5}
        \addlegendentry{\(\delta = 0.01\)}

        \pgfplotstableread[col sep=comma]{
       k,     ell,    ndof,    cost,   nElem,        eta,         mu,        res,    case,   maxGradU,       cond
       0,       0,     193,     193,      96,1.76237e-02,4.80764e-01,2.67715e-02,       Z,3.19183e-06,3.16582e+05
       1,       0,     193,     386,      96,1.76012e-02,2.01838e-02,2.42249e-02,       R,1.02740e-05,3.16582e+05
       1,       1,     225,     611,     112,1.48393e-02,4.81384e-01,2.57419e-02,       Z,1.29884e-05,1.00129e+06
       2,       0,     225,     836,     112,1.48714e-02,2.10175e-02,2.33022e-02,       R,2.93689e-05,1.00129e+06
       2,       1,     345,    1181,     172,1.10176e-02,2.32703e-02,2.17484e-02,       Z,3.04053e-05,1.76782e+06
       3,       0,     345,    1526,     172,1.10154e-02,1.87547e-02,1.96280e-02,       Z,5.17603e-05,1.76782e+06
       4,       0,     345,    1871,     172,1.10123e-02,1.62503e-02,1.79983e-02,       Z,7.53309e-05,1.76782e+06
       5,       0,     345,    2216,     172,1.10099e-02,1.42411e-02,1.67160e-02,       R,1.00032e-04,1.76782e+06
       5,       1,     441,    2657,     220,5.45785e-03,3.40434e-01,6.73712e-03,       Z,7.78977e-05,4.62117e+06
       6,       0,     441,    3098,     220,5.46326e-03,3.94374e-03,6.46451e-03,       Z,8.21723e-05,4.62117e+06
       7,       0,     441,    3539,     220,5.46831e-03,3.44950e-03,6.25401e-03,       Z,9.28305e-05,4.62117e+06
       8,       0,     441,    3980,     220,5.47298e-03,3.02830e-03,6.09144e-03,       Z,1.02934e-04,4.62117e+06
       9,       0,     441,    4421,     220,5.47731e-03,2.66763e-03,5.96583e-03,       Z,1.12453e-04,4.62117e+06
      10,       0,     441,    4862,     220,5.48132e-03,2.35748e-03,5.86866e-03,       Z,1.21378e-04,4.62117e+06
      11,       0,     441,    5303,     220,5.48505e-03,2.08971e-03,5.79336e-03,       Z,1.29716e-04,4.62117e+06
      12,       0,     441,    5744,     220,5.48851e-03,1.85763e-03,5.73490e-03,       R,1.37479e-04,4.62117e+06
      12,       1,     843,    6587,     421,4.45560e-03,1.90211e-01,6.24638e-03,       Z,1.49848e-04,1.69471e+07
      13,       0,     843,    7430,     421,4.45840e-03,4.37577e-03,5.91846e-03,       Z,1.69463e-04,1.69471e+07
      14,       0,     843,    8273,     421,4.46092e-03,3.89047e-03,5.65479e-03,       R,1.88423e-04,1.69471e+07
      14,       1,    1259,    9532,     629,3.58005e-03,3.04388e-02,5.32612e-03,       Z,1.88484e-04,2.35056e+07
      15,       0,    1259,   10791,     629,3.58390e-03,3.94024e-03,5.00618e-03,       Z,2.06975e-04,2.35056e+07
      16,       0,    1259,   12050,     629,3.58733e-03,3.49216e-03,4.74874e-03,       R,2.24759e-04,2.35056e+07
      16,       1,    1993,   14043,     996,2.94840e-03,1.12578e-01,4.83135e-03,       Z,2.27896e-04,5.58964e+07
      17,       0,    1993,   16036,     996,2.94947e-03,3.82671e-03,4.50515e-03,       R,2.48236e-04,5.58964e+07
      17,       1,    3137,   19173,    1568,2.40286e-03,8.64673e-02,4.73193e-03,       Z,2.50388e-04,9.86844e+07
      18,       0,    3137,   22310,    1568,2.40409e-03,4.07578e-03,4.35844e-03,       Z,2.72774e-04,9.86844e+07
      19,       0,    3137,   25447,    1568,2.40520e-03,3.63477e-03,4.04661e-03,       R,2.94294e-04,9.86844e+07
      19,       1,    4701,   30148,    2350,1.95592e-03,5.22247e-02,3.97627e-03,       Z,2.93042e-04,1.16060e+08
      20,       0,    4701,   34849,    2350,1.95643e-03,3.46170e-03,3.66593e-03,       Z,3.12382e-04,1.16060e+08
      21,       0,    4701,   39550,    2350,1.95693e-03,3.09995e-03,3.40428e-03,       R,3.30792e-04,1.16060e+08
      21,       1,    7251,   46801,    3625,1.55546e-03,5.63487e-02,3.07587e-03,       Z,3.29786e-04,2.35873e+08
      22,       0,    7251,   54052,    3625,1.55578e-03,2.65343e-03,2.84882e-03,       Z,3.46318e-04,2.35873e+08
      23,       0,    7251,   61303,    3625,1.55610e-03,2.38630e-03,2.65618e-03,       Z,3.61996e-04,2.35873e+08
      24,       0,    7251,   68554,    3625,1.55642e-03,2.15243e-03,2.49257e-03,       R,3.76845e-04,2.35873e+08
      24,       1,   10933,   79487,    5466,1.27120e-03,7.44506e-02,2.37738e-03,       Z,3.76063e-04,9.44373e+08
      25,       0,   10933,   90420,    5466,1.27137e-03,2.00888e-03,2.21570e-03,       R,3.89914e-04,9.44373e+08
      25,       1,   16415,  106835,    8207,1.03975e-03,4.99078e-02,2.18228e-03,       Z,3.90206e-04,2.21020e+09
      26,       0,   16415,  123250,    8207,1.03975e-03,1.91867e-03,2.02124e-03,       Z,4.03805e-04,2.21020e+09
      27,       0,   16415,  139665,    8207,1.03978e-03,1.73330e-03,1.88320e-03,       R,4.16659e-04,2.21020e+09
      27,       1,   25087,  164752,   12543,8.57826e-04,6.33657e-02,2.04078e-03,       Z,4.17623e-04,9.45823e+09
      28,       0,   25087,  189839,   12543,8.57987e-04,1.85166e-03,1.88476e-03,       Z,4.31062e-04,9.45823e+09
      29,       0,   25087,  214926,   12543,8.58152e-04,1.67807e-03,1.74936e-03,       R,4.45702e-04,9.45823e+09
      29,       1,   36521,  251447,   18260,7.05036e-04,3.08140e-02,1.65144e-03,       Z,4.45472e-04,3.60800e+10
      30,       0,   36521,  287968,   18260,7.05155e-04,1.49332e-03,1.52996e-03,       Z,4.59519e-04,3.60800e+10
      31,       0,   36521,  324489,   18260,7.05277e-04,1.35771e-03,1.42406e-03,       R,4.72768e-04,3.60800e+10
      31,       1,   53831,  378320,   26915,5.73596e-04,5.16900e-02,1.22805e-03,       Z,4.69312e-04,7.15173e+10
      32,       0,   53831,  432151,   26915,5.73684e-04,1.08581e-03,1.14363e-03,       Z,4.78888e-04,7.15173e+10
      33,       0,   53831,  485982,   26915,5.73776e-04,9.89278e-04,1.07014e-03,       R,4.90407e-04,7.15173e+10
      33,       1,   79435,  565417,   39717,4.73099e-04,2.34949e-02,1.06144e-03,       Z,4.90488e-04,1.42713e+11
      34,       0,   79435,  644852,   39717,4.73162e-04,9.50146e-04,9.87682e-04,       Z,5.01761e-04,1.42713e+11
      35,       0,   79435,  724287,   39717,4.73229e-04,8.66931e-04,9.23137e-04,       R,5.12375e-04,1.42713e+11
      35,       1,  115071,  839358,   57535,3.94404e-04,1.99087e-02,9.24222e-04,       Z,5.12524e-04,5.84801e+11
      36,       0,  115071,  954429,   57535,3.94424e-04,8.35833e-04,8.59752e-04,       Z,5.22655e-04,5.84801e+11
      37,       0,  115071, 1069500,   57535,3.94448e-04,7.63927e-04,8.03010e-04,       R,5.32187e-04,5.84801e+11
      37,       1,  167251, 1236751,   83625,3.28930e-04,3.77339e-02,8.32997e-04,       Z,5.32933e-04,2.28183e+12
      38,       0,  167251, 1404002,   83625,3.28966e-04,7.65288e-04,7.73685e-04,       R,5.42593e-04,2.28183e+12
      38,       1,  239061, 1643063,  119530,2.74886e-04,3.98282e-02,7.59972e-04,       Z,5.42707e-04,4.63888e+12
      39,       0,  239061, 1882124,  119530,2.74926e-04,7.08501e-04,7.04673e-04,       Z,5.51952e-04,4.63888e+12
      40,       0,  239061, 2121185,  119530,2.74968e-04,6.48812e-04,6.55580e-04,       R,5.61064e-04,4.63888e+12
      40,       1,  342357, 2463542,  171178,2.29749e-04,2.66512e-02,6.11468e-04,       Z,5.60682e-04,9.20755e+12
      41,       0,  342357, 2805899,  171178,2.29770e-04,5.66656e-04,5.68493e-04,       Z,5.68885e-04,9.20755e+12
      42,       0,  342357, 3148256,  171178,2.29792e-04,5.19981e-04,5.30251e-04,       R,5.76591e-04,9.20755e+12
      42,       1,  480575, 3628831,  240287,1.93994e-04,1.96721e-02,5.12206e-04,       Z,5.76857e-04,3.74795e+13
      43,       0,  480575, 4109406,  240287,1.94006e-04,4.74043e-04,4.76777e-04,       R,5.84355e-04,3.74795e+13
      43,       1,  676569, 4785975,  338284,1.63823e-04,2.15554e-02,4.94339e-04,       Z,5.84528e-04,7.39035e+13
      44,       0,  676569, 5462544,  338284,1.63837e-04,4.66400e-04,4.58777e-04,       Z,6.32554e-04,7.39035e+13
      45,       0,  676569, 6139113,  338284,1.63852e-04,4.28519e-04,4.26941e-04,       R,7.47702e-04,7.39035e+13
      45,       1,  942389, 7081502,  471194,1.38727e-04,1.34932e-02,4.26901e-04,       Z,7.51910e-04,1.47154e+14
      46,       0,  942389, 8023891,  471194,1.38741e-04,4.03728e-04,3.96564e-04,       R,8.69136e-04,1.47154e+14
      46,       1, 1325243, 9349134,  662621,1.17131e-04,1.53014e-02,3.84460e-04,       Z,9.42473e-04,5.99540e+14
        }\tableDataOne

        \pgfplotstableread[col sep=comma]{
       k,     ell,    ndof,    cost,   nElem,        eta,         mu,        res,    case,   maxGradU,       cond
       0,       0,     193,     193,      96,8.81184e-03,2.40382e-01,1.25869e-01,       Z,7.97958e-07,3.16582e+05
       1,       0,     193,     386,      96,8.80823e-03,1.20294e-01,6.56486e-02,       R,2.85170e-06,3.16582e+05
       1,       1,     225,     611,     112,7.42089e-03,2.08320e-01,1.10105e-01,       Z,3.50828e-06,1.00129e+06
       2,       0,     225,     836,     112,7.42656e-03,1.04448e-01,5.87041e-02,       R,8.51891e-06,1.00129e+06
       2,       1,     345,    1181,     172,5.50142e-03,1.04567e-01,5.83261e-02,       Z,8.78263e-06,1.76782e+06
       3,       0,     345,    1526,     172,5.50156e-03,5.29828e-02,3.49399e-02,       Z,1.58493e-05,1.76782e+06
       4,       0,     345,    1871,     172,5.50112e-03,2.78107e-02,2.56073e-02,       Z,2.42866e-05,1.76782e+06
       5,       0,     345,    2216,     172,5.50059e-03,1.60994e-02,2.21130e-02,       R,3.37815e-05,1.76782e+06
       5,       1,     441,    2657,     220,2.72555e-03,1.62692e-01,8.55830e-02,       Z,2.72969e-05,4.62117e+06
       6,       0,     441,    3098,     220,2.72697e-03,8.14419e-02,4.35584e-02,       Z,3.02761e-05,4.62117e+06
       7,       0,     441,    3539,     220,2.72837e-03,4.08489e-02,2.29862e-02,       Z,3.57057e-05,4.62117e+06
       8,       0,     441,    3980,     220,2.72974e-03,2.06110e-02,1.33819e-02,       Z,4.13742e-05,4.62117e+06
       9,       0,     441,    4421,     220,2.73105e-03,1.05935e-02,9.34337e-03,       Z,4.70802e-05,4.62117e+06
      10,       0,     441,    4862,     220,2.73233e-03,5.74648e-03,7.82759e-03,       Z,5.27861e-05,4.62117e+06
      11,       0,     441,    5303,     220,2.73355e-03,3.53269e-03,7.24263e-03,       Z,5.84612e-05,4.62117e+06
      12,       0,     441,    5744,     220,2.73474e-03,2.60368e-03,6.95828e-03,       R,6.40808e-05,4.62117e+06
      12,       1,     831,    6575,     415,2.25388e-03,9.45611e-02,4.96835e-02,       Z,6.81611e-05,1.75560e+07
      13,       0,     831,    7406,     415,2.25499e-03,4.73572e-02,2.57115e-02,       R,7.80482e-05,1.75560e+07
      13,       1,    1191,    8597,     595,1.83206e-03,4.73754e-02,2.56577e-02,       Z,7.80087e-05,1.75834e+07
      14,       0,    1191,    9788,     595,1.83353e-03,2.38535e-02,1.42835e-02,       Z,8.80091e-05,1.75834e+07
      15,       0,    1191,   10979,     595,1.83492e-03,1.22324e-02,9.37813e-03,       R,9.82556e-05,1.75834e+07
      15,       1,    1947,   12926,     973,1.49257e-03,6.13362e-02,3.26408e-02,       Z,9.94751e-05,4.87919e+07
      16,       0,    1947,   14873,     973,1.49310e-03,3.08012e-02,1.74196e-02,       Z,1.11093e-04,4.87919e+07
      17,       0,    1947,   16820,     973,1.49358e-03,1.56462e-02,1.05012e-02,       R,1.22766e-04,4.87919e+07
      17,       1,    3093,   19913,    1546,1.21026e-03,4.07212e-02,2.22175e-02,       Z,1.23203e-04,9.82524e+07
      18,       0,    3093,   23006,    1546,1.21059e-03,2.05560e-02,1.25422e-02,       Z,1.35530e-04,9.82524e+07
      19,       0,    3093,   26099,    1546,1.21091e-03,1.06313e-02,8.39866e-03,       R,1.48108e-04,9.82524e+07
      19,       1,    4671,   30770,    2335,9.73969e-04,3.04265e-02,1.70067e-02,       Z,1.47648e-04,1.39189e+08
      20,       0,    4671,   35441,    2335,9.74304e-04,1.54394e-02,1.00081e-02,       Z,1.59700e-04,1.39189e+08
      21,       0,    4671,   40112,    2335,9.74625e-04,8.11466e-03,7.09641e-03,       R,1.71657e-04,1.39189e+08
      21,       1,    7213,   47325,    3606,7.78799e-04,2.26476e-02,1.30106e-02,       Z,1.71390e-04,2.34423e+08
      22,       0,    7213,   54538,    3606,7.79047e-04,1.15617e-02,7.96363e-03,       Z,1.82990e-04,2.34423e+08
      23,       0,    7213,   61751,    3606,7.79284e-04,6.18361e-03,5.90341e-03,       Z,1.94479e-04,2.34423e+08
      24,       0,    7213,   68964,    3606,7.79510e-04,3.72183e-03,5.08890e-03,       R,2.05806e-04,2.34423e+08
      24,       1,   10807,   79771,    5403,6.35876e-04,3.84823e-02,2.05285e-02,       Z,2.05614e-04,9.45301e+08
      25,       0,   10807,   90578,    5403,6.36011e-04,1.93408e-02,1.09692e-02,       R,2.16564e-04,9.45301e+08
      25,       1,   16253,  106831,    8126,5.21902e-04,3.16848e-02,1.70455e-02,       Z,2.16686e-04,3.24634e+09
      26,       0,   16253,  123084,    8126,5.22015e-04,1.59550e-02,9.29456e-03,       Z,2.27652e-04,3.24634e+09
      27,       0,   16253,  139337,    8126,5.22126e-04,8.17769e-03,5.83815e-03,       R,2.38547e-04,3.24634e+09
      27,       1,   24949,  164286,   12474,4.30776e-04,3.18555e-02,1.70683e-02,       Z,2.39042e-04,8.83833e+09
      28,       0,   24949,  189235,   12474,4.30887e-04,1.60312e-02,9.25551e-03,       Z,2.50421e-04,8.83833e+09
      29,       0,   24949,  214184,   12474,4.30995e-04,8.20487e-03,5.76349e-03,       R,2.61644e-04,8.83833e+09
      29,       1,   36467,  250651,   18233,3.52861e-04,1.98306e-02,1.09826e-02,       Z,2.61598e-04,3.83249e+10
      30,       0,   36467,  287118,   18233,3.52935e-04,1.00509e-02,6.34651e-03,       Z,2.72565e-04,3.83249e+10
      31,       0,   36467,  323585,   18233,3.53007e-04,5.26237e-03,4.38698e-03,       R,2.83313e-04,3.83249e+10
      31,       1,   54179,  377764,   27089,2.85290e-04,2.69899e-02,1.43993e-02,       Z,2.83241e-04,5.26695e+10
      32,       0,   54179,  431943,   27089,2.85343e-04,1.35670e-02,7.67649e-03,       Z,2.93697e-04,5.26695e+10
      33,       0,   54179,  486122,   27089,2.85395e-04,6.90854e-03,4.59191e-03,       R,3.03932e-04,5.26695e+10
      33,       1,   79311,  565433,   39655,2.35321e-04,1.26547e-02,7.17704e-03,       Z,3.03785e-04,2.46169e+11
      34,       0,   79311,  644744,   39655,2.35371e-04,6.44986e-03,4.32480e-03,       Z,3.13651e-04,2.46169e+11
      35,       0,   79311,  724055,   39655,2.35421e-04,3.43732e-03,3.15811e-03,       R,3.23294e-04,2.46169e+11
      35,       1,  115209,  839264,   57604,1.96537e-04,9.41569e-03,5.54661e-03,       Z,3.23342e-04,6.28120e+11
      36,       0,  115209,  954473,   57604,1.96564e-04,4.84465e-03,3.54546e-03,       Z,3.32808e-04,6.28120e+11
      37,       0,  115209, 1069682,   57604,1.96591e-04,2.65678e-03,2.75647e-03,       R,3.42632e-04,6.28120e+11
      37,       1,  167369, 1237051,   83684,1.63939e-04,1.94159e-02,1.03741e-02,       Z,3.42774e-04,2.52729e+12
      38,       0,  167369, 1404420,   83684,1.63964e-04,9.76583e-03,5.57947e-03,       R,3.52671e-04,2.52729e+12
      38,       1,  239657, 1644077,  119828,1.36832e-04,2.23300e-02,1.18337e-02,       Z,3.52726e-04,2.52730e+12
      39,       0,  239657, 1883734,  119828,1.36855e-04,1.12117e-02,6.23694e-03,       Z,3.62441e-04,2.52730e+12
      40,       0,  239657, 2123391,  119828,1.36877e-04,5.69070e-03,3.64538e-03,       R,3.71920e-04,2.52730e+12
      40,       1,  342773, 2466164,  171386,1.14567e-04,1.34706e-02,7.31058e-03,       Z,3.71849e-04,1.03510e+13
      41,       0,  342773, 2808937,  171386,1.14582e-04,6.79857e-03,4.06365e-03,       Z,3.81025e-04,1.03510e+13
      42,       0,  342773, 3151710,  171386,1.14597e-04,3.51188e-03,2.64654e-03,       R,3.89971e-04,1.03510e+13
      42,       1,  480805, 3632515,  240402,9.67986e-05,1.01774e-02,5.62207e-03,       Z,3.89912e-04,4.14380e+13
      43,       0,  480805, 4113320,  240402,9.68086e-05,5.15730e-03,3.23892e-03,       R,3.98574e-04,4.14380e+13
      43,       1,  676187, 4789507,  338093,8.17629e-05,1.09048e-02,5.95661e-03,       Z,3.98691e-04,4.14528e+13
      44,       0,  676187, 5465694,  338093,8.17727e-05,5.51237e-03,3.36251e-03,       Z,4.07245e-04,4.14528e+13
      45,       0,  676187, 6141881,  338093,8.17824e-05,2.86367e-03,2.25005e-03,       R,4.15577e-04,4.14528e+13
      45,       1,  943257, 7085138,  471628,6.91625e-05,7.71207e-03,4.34930e-03,       Z,4.15778e-04,1.65809e+14
      46,       0,  943257, 8028395,  471628,6.91709e-05,3.92764e-03,2.60568e-03,       Z,4.87106e-04,1.65809e+14
      47,       0,  943257, 8971652,  471628,6.91791e-05,2.09000e-03,1.89439e-03,       R,5.65855e-04,1.65809e+14
      47,       1, 1325565,10297217,  662782,5.84207e-05,8.04217e-03,4.45580e-03,       Z,6.01168e-04,6.63058e+14
        }\tableDataOhFive

        \pgfplotstableread[col sep=comma]{
       k,     ell,    ndof,    cost,   nElem,        eta,         mu,        res,    case,   maxGradU,       cond
       0,       0,     193,     193,      96,1.76237e-03,4.80764e-02,2.25049e-01,       Z,3.19183e-08,3.16582e+05
       1,       0,     193,     386,      96,1.76233e-03,4.32692e-02,2.02625e-01,       R,1.24720e-07,3.16582e+05
       1,       1,     225,     611,     112,1.48430e-03,4.58733e-02,2.14735e-01,       Z,1.49959e-07,1.00129e+06
       2,       0,     225,     836,     112,1.48437e-03,4.12876e-02,1.93365e-01,       R,3.90762e-07,1.00129e+06
       2,       1,     345,    1181,     172,1.09950e-03,4.12997e-02,1.93347e-01,       Z,4.01499e-07,1.76782e+06
       3,       0,     345,    1526,     172,1.09951e-03,3.71729e-02,1.74133e-01,       Z,7.71457e-07,1.76782e+06
       4,       0,     345,    1871,     172,1.09951e-03,3.34602e-02,1.56878e-01,       Z,1.25255e-06,1.76782e+06
       5,       0,     345,    2216,     172,1.09950e-03,3.01204e-02,1.41388e-01,       R,1.83820e-06,1.76782e+06
       5,       1,     425,    2641,     212,5.55145e-04,2.49430e-02,1.17062e-01,       Z,1.53333e-06,4.62118e+06
       6,       0,     425,    3066,     212,5.55179e-04,2.24532e-02,1.05473e-01,       Z,1.83648e-06,4.62118e+06
       7,       0,     425,    3491,     212,5.55216e-04,2.02130e-02,9.50563e-02,       Z,2.23375e-06,4.62118e+06
       8,       0,     425,    3916,     212,5.55256e-04,1.81973e-02,8.56958e-02,       Z,2.72520e-06,4.62118e+06
       9,       0,     425,    4341,     212,5.55299e-04,1.63837e-02,7.72872e-02,       Z,3.25535e-06,4.62118e+06
      10,       0,     425,    4766,     212,5.55344e-04,1.47522e-02,6.97363e-02,       Z,3.82204e-06,4.62118e+06
      11,       0,     425,    5191,     212,5.55391e-04,1.32844e-02,6.29586e-02,       R,4.42322e-06,4.62118e+06
      11,       1,     793,    5984,     396,4.59946e-04,1.93529e-02,9.10933e-02,       Z,4.64205e-06,1.69456e+07
      12,       0,     793,    6777,     396,4.59979e-04,1.74249e-02,8.21608e-02,       Z,5.52732e-06,1.69456e+07
      13,       0,     793,    7570,     396,4.60013e-04,1.56906e-02,7.41417e-02,       R,6.47447e-06,1.69456e+07
      13,       1,    1173,    8743,     586,3.67860e-04,1.58786e-02,7.49974e-02,       Z,6.47646e-06,2.34342e+07
      14,       0,    1173,    9916,     586,3.67908e-04,1.43000e-02,6.77101e-02,       Z,7.48465e-06,2.34342e+07
      15,       0,    1173,   11089,     586,3.67959e-04,1.28803e-02,6.11740e-02,       R,8.54891e-06,2.34342e+07
      15,       1,    1937,   13026,     968,2.95333e-04,1.57168e-02,7.42697e-02,       Z,8.59140e-06,4.56264e+07
      16,       0,    1937,   14963,     968,2.95367e-04,1.41549e-02,6.70603e-02,       Z,9.75735e-06,4.56264e+07
      17,       0,    1937,   16900,     968,2.95402e-04,1.27503e-02,6.05941e-02,       R,1.09796e-05,4.56264e+07
      17,       1,    3039,   19939,    1519,2.42578e-04,1.41015e-02,6.68286e-02,       Z,1.10771e-05,8.88720e+07
      18,       0,    3039,   22978,    1519,2.42596e-04,1.27030e-02,6.03976e-02,       Z,1.24607e-05,8.88720e+07
      19,       0,    3039,   26017,    1519,2.42613e-04,1.14455e-02,5.46351e-02,       R,1.39051e-05,8.88720e+07
      19,       1,    4521,   30538,    2260,1.97722e-04,1.21547e-02,5.78840e-02,       Z,1.38676e-05,1.43302e+08
      20,       0,    4521,   35059,    2260,1.97744e-04,1.09527e-02,5.23818e-02,       Z,1.53283e-05,1.43302e+08
      21,       0,    4521,   39580,    2260,1.97766e-04,9.87227e-03,4.74567e-02,       R,1.68402e-05,1.43302e+08
      21,       1,    6945,   46525,    3472,1.56722e-04,1.00890e-02,4.84291e-02,       Z,1.68098e-05,2.32018e+08
      22,       0,    6945,   53470,    3472,1.56739e-04,9.09560e-03,4.39077e-02,       Z,1.83379e-05,2.32018e+08
      23,       0,    6945,   60415,    3472,1.56756e-04,8.20294e-03,3.98658e-02,       R,1.99105e-05,2.32018e+08
      23,       1,   10431,   70846,    5215,1.27658e-04,9.41393e-03,4.53389e-02,       Z,1.98943e-05,9.45245e+08
      24,       0,   10431,   81277,    5215,1.27669e-04,8.48862e-03,4.11357e-02,       Z,2.14917e-05,9.45245e+08
      25,       0,   10431,   91708,    5215,1.27680e-04,7.65728e-03,3.73799e-02,       R,2.31283e-05,9.45245e+08
      25,       1,   15791,  107499,    7895,1.04806e-04,8.64992e-03,4.18517e-02,       Z,2.31446e-05,2.21199e+09
      26,       0,   15791,  123290,    7895,1.04818e-04,7.80211e-03,3.80168e-02,       Z,2.48358e-05,2.21199e+09
      27,       0,   15791,  139081,    7895,1.04830e-04,7.04061e-03,3.45932e-02,       R,2.66045e-05,2.21199e+09
      27,       1,   24063,  163144,   12031,8.67408e-05,8.59808e-03,4.15953e-02,       Z,2.66312e-05,9.48764e+09
      28,       0,   24063,  187207,   12031,8.67524e-05,7.75548e-03,3.77863e-02,       Z,2.84721e-05,9.48764e+09
      29,       0,   24063,  211270,   12031,8.67640e-05,6.99870e-03,3.43864e-02,       R,3.03501e-05,9.48764e+09
      29,       1,   35599,  246869,   17799,7.05851e-05,7.11174e-03,3.48857e-02,       Z,3.03547e-05,3.82328e+10
      30,       0,   35599,  282468,   17799,7.05926e-05,6.42051e-03,3.17955e-02,       Z,3.22729e-05,3.82328e+10
      31,       0,   35599,  318067,   17799,7.06000e-05,5.80012e-03,2.90439e-02,       R,3.42250e-05,3.82328e+10
      31,       1,   53283,  371350,   26641,5.70104e-05,6.31538e-03,3.12946e-02,       Z,3.42235e-05,1.44093e+11
      32,       0,   53283,  424633,   26641,5.70175e-05,5.70484e-03,2.85780e-02,       Z,3.62066e-05,1.44093e+11
      33,       0,   53283,  477916,   26641,5.70246e-05,5.15706e-03,2.61615e-02,       R,3.82208e-05,1.44093e+11
      33,       1,   78885,  556801,   39442,4.69217e-05,5.43030e-03,2.73456e-02,       Z,3.82110e-05,2.85920e+11
      34,       0,   78885,  635686,   39442,4.69279e-05,4.91070e-03,2.50634e-02,       Z,4.02446e-05,2.85920e+11
      35,       0,   78885,  714571,   39442,4.69341e-05,4.44493e-03,2.30385e-02,       R,4.23060e-05,2.85920e+11
      35,       1,  114555,  829126,   57277,3.91088e-05,4.80219e-03,2.45581e-02,       Z,4.23113e-05,1.02124e+12
      36,       0,  114555,  943681,   57277,3.91129e-05,4.34735e-03,2.25841e-02,       Z,4.44046e-05,1.02124e+12
      37,       0,  114555, 1058236,   57277,3.91171e-05,3.94000e-03,2.08365e-02,       R,4.65234e-05,1.02124e+12
      37,       1,  166345, 1224581,   83172,3.26509e-05,4.84441e-03,2.46601e-02,       Z,4.65464e-05,2.30460e+12
      38,       0,  166345, 1390926,   83172,3.26545e-05,4.38421e-03,2.26565e-02,       R,4.87134e-05,2.30460e+12
      38,       1,  237199, 1628125,  118599,2.72645e-05,5.12111e-03,2.58167e-02,       Z,4.87167e-05,9.35431e+12
      39,       0,  237199, 1865324,  118599,2.72678e-05,4.63141e-03,2.36667e-02,       Z,5.09110e-05,9.35431e+12
      40,       0,  237199, 2102523,  118599,2.72710e-05,4.19247e-03,2.17593e-02,       R,5.31282e-05,9.35431e+12
      40,       1,  339675, 2442198,  169837,2.28137e-05,4.31563e-03,2.22775e-02,       Z,5.31190e-05,1.84448e+13
      41,       0,  339675, 2781873,  169837,2.28164e-05,3.90931e-03,2.05234e-02,       Z,5.53485e-05,1.84448e+13
      42,       0,  339675, 3121548,  169837,2.28190e-05,3.54557e-03,1.89721e-02,       R,5.75984e-05,1.84448e+13
      42,       1,  479563, 3601111,  239781,1.92130e-05,3.75873e-03,1.98479e-02,       Z,5.75996e-05,7.49146e+13
      43,       0,  479563, 4080674,  239781,1.92150e-05,3.41038e-03,1.83670e-02,       R,5.98701e-05,7.49146e+13
      43,       1,  676221, 4756895,  338110,1.62144e-05,3.46990e-03,1.86108e-02,       Z,5.98766e-05,1.47772e+14
      44,       0,  676221, 5433116,  338110,1.62161e-05,3.15209e-03,1.72753e-02,       Z,6.21722e-05,1.47772e+14
      45,       0,  676221, 6109337,  338110,1.62178e-05,2.86822e-03,1.60992e-02,       R,6.44856e-05,1.47772e+14
      45,       1,  940347, 7049684,  470173,1.37302e-05,2.96878e-03,1.64935e-02,       Z,6.44879e-05,2.94285e+14
      46,       0,  940347, 7990031,  470173,1.37317e-05,2.70438e-03,1.54077e-02,       Z,6.68204e-05,2.94285e+14
      47,       0,  940347, 8930378,  470173,1.37332e-05,2.46871e-03,1.44545e-02,       R,6.91688e-05,2.94285e+14
      47,       1, 1318035,10248413,  659017,1.15952e-05,2.64572e-03,1.51227e-02,       Z,6.91736e-05,1.19817e+15
        }\tableDataOhOne

        \pgfplotstableread[col sep=comma]{
       k,     ell,    ndof,    cost,   nElem,        eta,         mu,        res,    case,   maxGradU,       cond
       0,       0,     193,     193,      96,8.81184e-04,2.40382e-02,2.37520e-01,       Z,7.97958e-09,3.16582e+05
       1,       0,     193,     386,      96,8.81179e-04,2.28363e-02,2.25672e-01,       R,3.15449e-08,3.16582e+05
       1,       1,     225,     611,     112,7.42155e-04,2.34655e-02,2.31827e-01,       Z,3.78157e-08,1.00129e+06
       2,       0,     225,     836,     112,7.42164e-04,2.22924e-02,2.20266e-01,       R,9.95566e-08,1.00129e+06
       2,       1,     345,    1181,     172,5.49730e-04,2.22980e-02,2.20259e-01,       Z,1.02249e-07,1.76782e+06
       3,       0,     345,    1526,     172,5.49732e-04,2.11834e-02,2.09277e-01,       Z,1.98424e-07,1.76782e+06
       4,       0,     345,    1871,     172,5.49732e-04,2.01248e-02,1.98853e-01,       Z,3.25302e-07,1.76782e+06
       5,       0,     345,    2216,     172,5.49731e-04,1.91193e-02,1.88957e-01,       R,4.81898e-07,1.76782e+06
       5,       1,     425,    2641,     212,2.77555e-04,1.39817e-02,1.38184e-01,       Z,4.03780e-07,4.62118e+06
       6,       0,     425,    3066,     212,2.77560e-04,1.32833e-02,1.31307e-01,       Z,4.89816e-07,4.62118e+06
       7,       0,     425,    3491,     212,2.77565e-04,1.26198e-02,1.24776e-01,       Z,5.99655e-07,4.62118e+06
       8,       0,     425,    3916,     212,2.77572e-04,1.19895e-02,1.18576e-01,       Z,7.38535e-07,4.62118e+06
       9,       0,     425,    4341,     212,2.77579e-04,1.13908e-02,1.12689e-01,       Z,8.90258e-07,4.62118e+06
      10,       0,     425,    4766,     212,2.77587e-04,1.08222e-02,1.07100e-01,       Z,1.05444e-06,4.62118e+06
      11,       0,     425,    5191,     212,2.77595e-04,1.02820e-02,1.01794e-01,       R,1.23069e-06,4.62118e+06
      11,       1,     793,    5984,     396,2.29908e-04,1.23760e-02,1.22419e-01,       Z,1.28838e-06,1.69456e+07
      12,       0,     793,    6777,     396,2.29913e-04,1.17583e-02,1.16349e-01,       Z,1.54353e-06,1.69456e+07
      13,       0,     793,    7570,     396,2.29918e-04,1.11716e-02,1.10586e-01,       R,1.81929e-06,1.69456e+07
      13,       1,    1173,    8743,     586,1.83848e-04,1.12430e-02,1.11280e-01,       Z,1.81981e-06,2.34342e+07
      14,       0,    1173,    9916,     586,1.83856e-04,1.06822e-02,1.05775e-01,       Z,2.11620e-06,2.34342e+07
      15,       0,    1173,   11089,     586,1.83865e-04,1.01495e-02,1.00548e-01,       R,2.43213e-06,2.34342e+07
      15,       1,    1937,   13026,     968,1.47551e-04,1.09070e-02,1.08000e-01,       Z,2.44335e-06,4.56264e+07
      16,       0,    1937,   14963,     968,1.47557e-04,1.03631e-02,1.02663e-01,       Z,2.79121e-06,4.56264e+07
      17,       0,    1937,   16900,     968,1.47564e-04,9.84651e-03,9.75986e-02,       R,3.15921e-06,4.56264e+07
      17,       1,    3031,   19931,    1515,1.21361e-04,1.03164e-02,1.02221e-01,       Z,3.18526e-06,8.88720e+07
      18,       0,    3031,   22962,    1515,1.21364e-04,9.80234e-03,9.71826e-02,       Z,3.60199e-06,8.88720e+07
      19,       0,    3031,   25993,    1515,1.21367e-04,9.31408e-03,9.24019e-02,       R,4.04078e-06,8.88720e+07
      19,       1,    4543,   30536,    2271,9.84201e-05,9.45772e-03,9.38107e-02,       Z,4.03081e-06,1.43362e+08
      20,       0,    4543,   35079,    2271,9.84242e-05,8.98681e-03,8.92024e-02,       Z,4.47986e-06,1.43362e+08
      21,       0,    4543,   39622,    2271,9.84284e-05,8.53958e-03,8.48301e-02,       R,4.94862e-06,1.43362e+08
      21,       1,    6999,   46621,    3499,7.79492e-05,8.41880e-03,8.36400e-02,       Z,4.93979e-06,2.32018e+08
      22,       0,    6999,   53620,    3499,7.79528e-05,8.00008e-03,7.95495e-02,       Z,5.41816e-06,2.32018e+08
      23,       0,    6999,   60619,    3499,7.79565e-05,7.60244e-03,7.56692e-02,       Z,5.91451e-06,2.32018e+08
      24,       0,    6999,   67618,    3499,7.79602e-05,7.22483e-03,7.19886e-02,       R,6.42839e-06,2.32018e+08
      24,       1,   10543,   78161,    5271,6.34808e-05,7.50125e-03,7.46923e-02,       Z,6.42660e-06,9.45245e+08
      25,       0,   10543,   88704,    5271,6.34830e-05,7.12875e-03,7.10623e-02,       R,6.95559e-06,9.45245e+08
      25,       1,   15879,  104583,    7939,5.22758e-05,7.39405e-03,7.36581e-02,       Z,6.96113e-06,2.21159e+09
      26,       0,   15879,  120462,    7939,5.22784e-05,7.02701e-03,7.00834e-02,       Z,7.51252e-06,2.21159e+09
      27,       0,   15879,  136341,    7939,5.22810e-05,6.67848e-03,6.66934e-02,       R,8.08581e-06,2.21159e+09
      27,       1,   24321,  160662,   12160,4.29677e-05,6.99191e-03,6.97560e-02,       Z,8.09330e-06,9.48207e+09
      28,       0,   24321,  184983,   12160,4.29702e-05,6.64522e-03,6.63856e-02,       Z,8.69334e-06,9.48207e+09
      29,       0,   24321,  209304,   12160,4.29727e-05,6.31602e-03,6.31900e-02,       R,9.31002e-06,9.48207e+09
      29,       1,   36135,  245439,   18067,3.49559e-05,6.25909e-03,6.26337e-02,       Z,9.31348e-06,3.82670e+10
      30,       0,   36135,  281574,   18067,3.49574e-05,5.94938e-03,5.96323e-02,       Z,9.95006e-06,3.82670e+10
      31,       0,   36135,  317709,   18067,3.49590e-05,5.65533e-03,5.67875e-02,       R,1.06027e-05,3.82670e+10
      31,       1,   53905,  371614,   26952,2.82733e-05,5.62116e-03,5.64528e-02,       Z,1.05940e-05,1.44262e+11
      32,       0,   53905,  425519,   26952,2.82749e-05,5.34366e-03,5.37718e-02,       Z,1.12532e-05,1.44262e+11
      33,       0,   53905,  479424,   26952,2.82765e-05,5.08022e-03,5.12316e-02,       R,1.19272e-05,1.44262e+11
      33,       1,   79499,  558923,   39749,2.33475e-05,5.16104e-03,5.20136e-02,       Z,1.19319e-05,2.88853e+11
      34,       0,   79499,  638422,   39749,2.33490e-05,4.90687e-03,4.95666e-02,       Z,1.26253e-05,2.88853e+11
      35,       0,   79499,  717921,   39749,2.33504e-05,4.66560e-03,4.72491e-02,       R,1.33331e-05,2.88853e+11
      35,       1,  115401,  833322,   57700,1.94574e-05,4.75791e-03,4.81384e-02,       Z,1.33348e-05,1.15561e+12
      36,       0,  115401,  948723,   57700,1.94585e-05,4.52422e-03,4.58971e-02,       Z,1.40584e-05,1.15561e+12
      37,       0,  115401, 1064124,   57700,1.94595e-05,4.30241e-03,4.37751e-02,       R,1.47957e-05,1.15561e+12
      37,       1,  167185, 1231309,   83592,1.62432e-05,4.44826e-03,4.51747e-02,       Z,1.48006e-05,2.30191e+12
      38,       0,  167185, 1398494,   83592,1.62440e-05,4.23033e-03,4.30918e-02,       R,1.55564e-05,2.30191e+12
      38,       1,  237693, 1636187,  118846,1.35944e-05,4.36624e-03,4.43935e-02,       Z,1.55550e-05,9.35655e+12
      39,       0,  237693, 1873880,  118846,1.35952e-05,4.15248e-03,4.23518e-02,       Z,1.63235e-05,9.35655e+12
      40,       0,  237693, 2111573,  118846,1.35960e-05,3.94962e-03,4.04198e-02,       R,1.71079e-05,9.35655e+12
      40,       1,  339555, 2451128,  169777,1.13829e-05,3.94462e-03,4.03714e-02,       Z,1.71069e-05,1.84481e+13
      41,       0,  339555, 2790683,  169777,1.13835e-05,3.75236e-03,3.85456e-02,       Z,1.79031e-05,1.84481e+13
      42,       0,  339555, 3130238,  169777,1.13842e-05,3.56993e-03,3.68190e-02,       R,1.87119e-05,1.84481e+13
      42,       1,  476345, 3606583,  238172,9.61315e-06,3.58738e-03,3.69835e-02,       Z,1.87106e-05,7.48836e+13
      43,       0,  476345, 4082928,  238172,9.61368e-06,3.41340e-03,3.53416e-02,       R,1.95304e-05,7.48836e+13
      43,       1,  668275, 4751203,  334137,8.12870e-06,3.42239e-03,3.54267e-02,       Z,1.95350e-05,1.47800e+14
      44,       0,  668275, 5419478,  334137,8.12910e-06,3.25689e-03,3.38707e-02,       Z,2.03715e-05,1.47800e+14
      45,       0,  668275, 6087753,  334137,8.12959e-06,3.09992e-03,3.24008e-02,       R,2.12198e-05,1.47800e+14
      45,       1,  929573, 7017326,  464786,6.88098e-06,3.11356e-03,3.25282e-02,       Z,2.12204e-05,5.99521e+14
      46,       0,  929573, 7946899,  464786,6.88124e-06,2.96398e-03,3.11329e-02,       Z,2.20810e-05,5.99521e+14
      47,       0,  929573, 8876472,  464786,6.88147e-06,2.82214e-03,2.98157e-02,       R,2.29529e-05,5.99521e+14
      47,       1, 1298345,10174817,  649172,5.81940e-06,2.83283e-03,2.99141e-02,       Z,2.29519e-05,1.18131e+15
        }\tableDataOhOhFive

        \pgfplotstableread[col sep=comma]{
       k,     ell,    ndof,    cost,   nElem,        eta,         mu,        res,    case,   maxGradU,       cond
       0,       0,     193,     193,      96,1.76237e-04,4.80764e-03,2.47503e-01,       Z,3.19183e-10,3.16582e+05
       1,       0,     193,     386,      96,1.76237e-04,4.75956e-03,2.45032e-01,       R,1.27372e-09,3.16582e+05
       1,       1,     225,     611,     112,1.48431e-04,4.78473e-03,2.46277e-01,       Z,1.52324e-09,1.00129e+06
       2,       0,     225,     836,     112,1.48431e-04,4.73689e-03,2.43817e-01,       R,4.04461e-09,1.00129e+06
       2,       1,     345,    1181,     172,1.09945e-04,4.73794e-03,2.43816e-01,       Z,4.15258e-09,1.76782e+06
       3,       0,     345,    1526,     172,1.09945e-04,4.69056e-03,2.41380e-01,       Z,8.12738e-09,1.76782e+06
       4,       0,     345,    1871,     172,1.09945e-04,4.64366e-03,2.38969e-01,       Z,1.34383e-08,1.76782e+06
       5,       0,     345,    2216,     172,1.09945e-04,4.59723e-03,2.36582e-01,       R,2.00760e-08,1.76782e+06
       5,       1,     425,    2641,     212,5.55098e-05,3.25479e-03,1.67475e-01,       Z,1.68848e-08,4.62118e+06
       6,       0,     425,    3066,     212,5.55098e-05,3.22225e-03,1.65802e-01,       Z,2.07167e-08,4.62118e+06
       7,       0,     425,    3491,     212,5.55099e-05,3.19003e-03,1.64146e-01,       Z,2.55195e-08,4.62118e+06
       8,       0,     425,    3916,     212,5.55099e-05,3.15814e-03,1.62506e-01,       Z,3.17081e-08,4.62118e+06
       9,       0,     425,    4341,     212,5.55100e-05,3.12656e-03,1.60883e-01,       Z,3.85519e-08,4.62118e+06
      10,       0,     425,    4766,     212,5.55101e-05,3.09530e-03,1.59277e-01,       Z,4.60463e-08,4.62118e+06
      11,       0,     425,    5191,     212,5.55102e-05,3.06436e-03,1.57686e-01,       R,5.41869e-08,4.62118e+06
      11,       1,     793,    5984,     396,4.59763e-05,3.38604e-03,1.74224e-01,       Z,5.66067e-08,1.69456e+07
      12,       0,     793,    6777,     396,4.59764e-05,3.35218e-03,1.72484e-01,       Z,6.82445e-08,1.69456e+07
      13,       0,     793,    7570,     396,4.59765e-05,3.31867e-03,1.70762e-01,       R,8.09526e-08,1.69456e+07
      13,       1,    1173,    8743,     586,3.67625e-05,3.33115e-03,1.71398e-01,       Z,8.09747e-08,2.34342e+07
      14,       0,    1173,    9916,     586,3.67626e-05,3.29785e-03,1.69686e-01,       Z,9.47721e-08,2.34342e+07
      15,       0,    1173,   11089,     586,3.67627e-05,3.26488e-03,1.67992e-01,       R,1.09630e-07,2.34342e+07
      15,       1,    1937,   13026,     968,2.94997e-05,3.33877e-03,1.71789e-01,       Z,1.10103e-07,4.56264e+07
      16,       0,    1937,   14963,     968,2.94998e-05,3.30539e-03,1.70073e-01,       Z,1.26565e-07,4.56264e+07
      17,       0,    1937,   16900,     968,2.94998e-05,3.27235e-03,1.68376e-01,       R,1.44150e-07,4.56264e+07
      17,       1,    3031,   19931,    1515,2.42772e-05,3.35761e-03,1.72759e-01,       Z,1.45257e-07,9.00319e+07
      18,       0,    3031,   22962,    1515,2.42772e-05,3.32405e-03,1.71034e-01,       Z,1.65211e-07,9.00319e+07
      19,       0,    3031,   25993,    1515,2.42772e-05,3.29083e-03,1.69327e-01,       R,1.86419e-07,9.00319e+07
      19,       1,    4533,   30526,    2266,1.96421e-05,3.29426e-03,1.69502e-01,       Z,1.85971e-07,1.44389e+08
      20,       0,    4533,   35059,    2266,1.96421e-05,3.26133e-03,1.67810e-01,       Z,2.07914e-07,1.44389e+08
      21,       0,    4533,   39592,    2266,1.96422e-05,3.22874e-03,1.66136e-01,       R,2.31034e-07,1.44389e+08
      21,       1,    6999,   46591,    3499,1.56020e-05,3.18818e-03,1.64048e-01,       Z,2.30695e-07,2.32063e+08
      22,       0,    6999,   53590,    3499,1.56020e-05,3.15631e-03,1.62411e-01,       Z,2.54610e-07,2.32063e+08
      23,       0,    6999,   60589,    3499,1.56021e-05,3.12477e-03,1.60791e-01,       Z,2.79650e-07,2.32063e+08
      24,       0,    6999,   67588,    3499,1.56021e-05,3.09354e-03,1.59186e-01,       R,3.05810e-07,2.32063e+08
      24,       1,   10603,   78191,    5301,1.26620e-05,3.10008e-03,1.59522e-01,       Z,3.05741e-07,9.45245e+08
      25,       0,   10603,   88794,    5301,1.26620e-05,3.06910e-03,1.57931e-01,       R,3.32936e-07,9.45245e+08
      25,       1,   15905,  104699,    7952,1.04405e-05,3.10097e-03,1.59570e-01,       Z,3.33236e-07,2.22476e+09
      26,       0,   15905,  120604,    7952,1.04406e-05,3.06998e-03,1.57978e-01,       Z,3.61852e-07,2.22476e+09
      27,       0,   15905,  136509,    7952,1.04406e-05,3.03930e-03,1.56402e-01,       R,3.91584e-07,2.22476e+09
      27,       1,   24451,  160960,   12225,8.55582e-06,3.07098e-03,1.58031e-01,       Z,3.91904e-07,9.51411e+09
      28,       0,   24451,  185411,   12225,8.55585e-06,3.04029e-03,1.56455e-01,       Z,4.23289e-07,9.51411e+09
      29,       0,   24451,  209862,   12225,8.55588e-06,3.00991e-03,1.54895e-01,       R,4.55933e-07,9.51411e+09
      29,       1,   36343,  246205,   18171,6.98109e-06,2.99567e-03,1.54163e-01,       Z,4.56087e-07,3.82810e+10
      30,       0,   36343,  282548,   18171,6.98110e-06,2.96573e-03,1.52626e-01,       Z,4.90037e-07,3.82810e+10
      31,       0,   36343,  318891,   18171,6.98112e-06,2.93610e-03,1.51104e-01,       R,5.25134e-07,3.82810e+10
      31,       1,   53843,  372734,   26921,5.65181e-06,2.90787e-03,1.49652e-01,       Z,5.24715e-07,1.41778e+11
      32,       0,   53843,  426577,   26921,5.65183e-06,2.87881e-03,1.48160e-01,       Z,5.60505e-07,1.41778e+11
      33,       0,   53843,  480420,   26921,5.65185e-06,2.85005e-03,1.46684e-01,       R,5.97400e-07,1.41778e+11
      33,       1,   79619,  560039,   39809,4.65073e-06,2.85835e-03,1.47111e-01,       Z,5.97627e-07,2.85412e+11
      34,       0,   79619,  639658,   39809,4.65075e-06,2.82979e-03,1.45645e-01,       Z,6.35860e-07,2.85412e+11
      35,       0,   79619,  719277,   39809,4.65076e-06,2.80152e-03,1.44194e-01,       R,6.75199e-07,2.85412e+11
      35,       1,  116469,  835746,   58234,3.87207e-06,2.80769e-03,1.44511e-01,       Z,6.75273e-07,1.16716e+12
      36,       0,  116469,  952215,   58234,3.87209e-06,2.77964e-03,1.43072e-01,       Z,7.15788e-07,1.16716e+12
      37,       0,  116469, 1068684,   58234,3.87208e-06,2.75187e-03,1.41647e-01,       R,7.57399e-07,1.16716e+12
      37,       1,  168977, 1237661,   84488,3.23085e-06,2.76311e-03,1.42225e-01,       Z,7.57632e-07,2.30591e+12
      38,       0,  168977, 1406638,   84488,3.23082e-06,2.73551e-03,1.40808e-01,       R,8.00626e-07,2.30591e+12
      38,       1,  240001, 1646639,  120000,2.70242e-06,2.73548e-03,1.40807e-01,       Z,8.00746e-07,9.35548e+12
      39,       0,  240001, 1886640,  120000,2.70237e-06,2.70816e-03,1.39405e-01,       Z,8.45015e-07,9.35548e+12
      40,       0,  240001, 2126641,  120000,2.70230e-06,2.68111e-03,1.38017e-01,       R,8.90387e-07,9.35548e+12
      40,       1,  343197, 2469838,  171598,2.25865e-06,2.67605e-03,1.37757e-01,       Z,8.90257e-07,1.84459e+13
      41,       0,  343197, 2813035,  171598,2.25875e-06,2.64932e-03,1.36386e-01,       Z,9.36587e-07,1.84459e+13
      42,       0,  343197, 3156232,  171598,2.25862e-06,2.62286e-03,1.35029e-01,       R,9.84000e-07,1.84459e+13
      42,       1,  481883, 3638115,  240941,1.90560e-06,2.62069e-03,1.34917e-01,       Z,9.83921e-07,7.49158e+13
      43,       0,  481883, 4119998,  240941,1.90552e-06,2.59452e-03,1.33575e-01,       R,1.03233e-06,7.49158e+13
      43,       1,  676799, 4796797,  338399,1.60884e-06,2.59805e-03,1.33756e-01,       Z,1.03242e-06,1.47794e+14
      44,       0,  676799, 5473596,  338399,1.60860e-06,2.57211e-03,1.32426e-01,       Z,1.08200e-06,1.47794e+14
      45,       0,  676799, 6150395,  338399,1.60853e-06,2.54642e-03,1.31109e-01,       R,1.13263e-06,1.47794e+14
      45,       1,  940019, 7090414,  470009,1.35673e-06,2.54669e-03,1.31122e-01,       Z,1.13269e-06,5.96463e+14
      46,       0,  940019, 8030433,  470009,1.35649e-06,2.52126e-03,1.29818e-01,       Z,1.18445e-06,5.96463e+14
      47,       0,  940019, 8970452,  470009,1.35645e-06,2.49609e-03,1.28528e-01,       R,1.23726e-06,5.96463e+14
      47,       1, 1314507,10284959,  657253,1.13687e-06,2.49423e-03,1.28433e-01,       Z,1.23721e-06,1.18113e+15
        }\tableDataOhOhOne

        %
        %

        \addplot+ [marker1, adaptive, forget plot]
        table [col sep=comma, x=cumulativeNdof, y=maxGradU] {\tableDataOne};

        \addplot+ [marker2, adaptive, forget plot]
        table [col sep=comma, x=cumulativeNdof, y=maxGradU] {\tableDataOhFive};

        \addplot+ [marker3, adaptive, forget plot]
        table [col sep=comma, x=cumulativeNdof, y=maxGradU] {\tableDataOhOne};

        \addplot+ [marker4, adaptive, forget plot]
        table [col sep=comma, x=cumulativeNdof, y=maxGradU] {\tableDataOhOhFive};

        \addplot+ [marker5, adaptive, forget plot]
        table [col sep=comma, x=cumulativeNdof, y=maxGradU] {\tableDataOhOhOne};

    \end{loglogaxis}
\end{tikzpicture}

%% file: figures/Fig08a_PorousMedium_convergence_weighting.tex
\begin{tikzpicture}[>=stealth]
    \begin{loglogaxis}[%
            width            = 5.5cm,%
            xlabel           = {cumulative ndof},%
            ylabel           = {Least-squares functional \(N^k_\ell\)},%
            ymajorgrids      = true,%
            ymax             = 2e-1,%
            font             = \footnotesize,%
            grid style       = {%
                densely dotted,%
                semithick%
            },%
            legend style     = {%
                legend pos = north east,%
                font = \footnotesize%
            },%
        ]

        \addlegendimage{marker1}
        \addlegendentry{emph.\ grad.}
        \addlegendimage{marker2}
        \addlegendentry{balanced w.}
        \addlegendimage{marker3}
        \addlegendentry{downsc. flux}
        \addlegendimage{marker4}
        \addlegendentry{split w.}

        \pgfplotstableread[col sep=comma]{
       k,     ell,    ndof,   nElem,        eta,         mu,        res,    case,   maxGradU
       0,       0,     193,      96,1.76237e-02,4.80764e-01,2.67715e-02,       Z,3.19183e-06
       1,       0,     193,      96,1.76012e-02,2.01838e-02,2.42249e-02,       R,1.02740e-05
       1,       1,     225,     112,1.48393e-02,4.81384e-01,2.57419e-02,       Z,1.29884e-05
       2,       0,     225,     112,1.48714e-02,2.10175e-02,2.33022e-02,       R,2.93689e-05
       2,       1,     345,     172,1.10176e-02,2.32703e-02,2.17484e-02,       Z,3.04053e-05
       3,       0,     345,     172,1.10154e-02,1.87547e-02,1.96280e-02,       Z,5.17603e-05
       4,       0,     345,     172,1.10123e-02,1.62503e-02,1.79983e-02,       Z,7.53309e-05
       5,       0,     345,     172,1.10099e-02,1.42411e-02,1.67160e-02,       R,1.00032e-04
       5,       1,     441,     220,5.45785e-03,3.40434e-01,6.73712e-03,       Z,7.78977e-05
       6,       0,     441,     220,5.46326e-03,3.94374e-03,6.46451e-03,       Z,8.21723e-05
       7,       0,     441,     220,5.46831e-03,3.44950e-03,6.25401e-03,       Z,9.28305e-05
       8,       0,     441,     220,5.47298e-03,3.02830e-03,6.09144e-03,       Z,1.02934e-04
       9,       0,     441,     220,5.47731e-03,2.66763e-03,5.96583e-03,       Z,1.12453e-04
      10,       0,     441,     220,5.48132e-03,2.35748e-03,5.86866e-03,       Z,1.21378e-04
      11,       0,     441,     220,5.48505e-03,2.08971e-03,5.79336e-03,       Z,1.29716e-04
      12,       0,     441,     220,5.48851e-03,1.85763e-03,5.73490e-03,       R,1.37479e-04
      12,       1,     843,     421,4.45560e-03,1.90211e-01,6.24638e-03,       Z,1.49848e-04
      13,       0,     843,     421,4.45840e-03,4.37577e-03,5.91846e-03,       Z,1.69463e-04
      14,       0,     843,     421,4.46092e-03,3.89047e-03,5.65479e-03,       R,1.88423e-04
      14,       1,    1259,     629,3.58005e-03,3.04388e-02,5.32612e-03,       Z,1.88484e-04
      15,       0,    1259,     629,3.58390e-03,3.94024e-03,5.00618e-03,       Z,2.06975e-04
      16,       0,    1259,     629,3.58733e-03,3.49216e-03,4.74874e-03,       R,2.24759e-04
      16,       1,    1993,     996,2.94840e-03,1.12578e-01,4.83135e-03,       Z,2.27896e-04
      17,       0,    1993,     996,2.94947e-03,3.82671e-03,4.50515e-03,       R,2.48236e-04
      17,       1,    3137,    1568,2.40286e-03,8.64673e-02,4.73193e-03,       Z,2.50388e-04
      18,       0,    3137,    1568,2.40409e-03,4.07578e-03,4.35844e-03,       Z,2.72774e-04
      19,       0,    3137,    1568,2.40520e-03,3.63477e-03,4.04661e-03,       R,2.94294e-04
      19,       1,    4701,    2350,1.95592e-03,5.22247e-02,3.97627e-03,       Z,2.93042e-04
      20,       0,    4701,    2350,1.95643e-03,3.46170e-03,3.66593e-03,       Z,3.12382e-04
      21,       0,    4701,    2350,1.95693e-03,3.09995e-03,3.40428e-03,       R,3.30792e-04
      21,       1,    7251,    3625,1.55546e-03,5.63487e-02,3.07587e-03,       Z,3.29786e-04
      22,       0,    7251,    3625,1.55578e-03,2.65343e-03,2.84882e-03,       Z,3.46318e-04
      23,       0,    7251,    3625,1.55610e-03,2.38630e-03,2.65618e-03,       Z,3.61996e-04
      24,       0,    7251,    3625,1.55642e-03,2.15243e-03,2.49257e-03,       R,3.76845e-04
      24,       1,   10933,    5466,1.27120e-03,7.44506e-02,2.37738e-03,       Z,3.76063e-04
      25,       0,   10933,    5466,1.27137e-03,2.00888e-03,2.21570e-03,       R,3.89914e-04
      25,       1,   16415,    8207,1.03975e-03,4.99078e-02,2.18228e-03,       Z,3.90206e-04
      26,       0,   16415,    8207,1.03975e-03,1.91867e-03,2.02124e-03,       Z,4.03805e-04
      27,       0,   16415,    8207,1.03978e-03,1.73330e-03,1.88320e-03,       R,4.16659e-04
      27,       1,   25087,   12543,8.57826e-04,6.33657e-02,2.04078e-03,       Z,4.17623e-04
      28,       0,   25087,   12543,8.57987e-04,1.85166e-03,1.88476e-03,       Z,4.31062e-04
      29,       0,   25087,   12543,8.58152e-04,1.67807e-03,1.74936e-03,       R,4.45702e-04
      29,       1,   36521,   18260,7.05036e-04,3.08140e-02,1.65144e-03,       Z,4.45472e-04
      30,       0,   36521,   18260,7.05155e-04,1.49332e-03,1.52996e-03,       Z,4.59519e-04
      31,       0,   36521,   18260,7.05277e-04,1.35771e-03,1.42406e-03,       R,4.72768e-04
      31,       1,   53831,   26915,5.73596e-04,5.16900e-02,1.22805e-03,       Z,4.69312e-04
      32,       0,   53831,   26915,5.73684e-04,1.08581e-03,1.14363e-03,       Z,4.78888e-04
      33,       0,   53831,   26915,5.73776e-04,9.89278e-04,1.07014e-03,       R,4.90407e-04
      33,       1,   79435,   39717,4.73099e-04,2.34949e-02,1.06144e-03,       Z,4.90488e-04
      34,       0,   79435,   39717,4.73162e-04,9.50146e-04,9.87682e-04,       Z,5.01761e-04
      35,       0,   79435,   39717,4.73229e-04,8.66931e-04,9.23137e-04,       R,5.12375e-04
      35,       1,  115071,   57535,3.94404e-04,1.99087e-02,9.24222e-04,       Z,5.12524e-04
      36,       0,  115071,   57535,3.94424e-04,8.35833e-04,8.59752e-04,       Z,5.22655e-04
      37,       0,  115071,   57535,3.94448e-04,7.63927e-04,8.03010e-04,       R,5.32187e-04
      37,       1,  167251,   83625,3.28930e-04,3.77339e-02,8.32997e-04,       Z,5.32933e-04
      38,       0,  167251,   83625,3.28966e-04,7.65288e-04,7.73685e-04,       R,5.42593e-04
      38,       1,  239061,  119530,2.74886e-04,3.98282e-02,7.59972e-04,       Z,5.42707e-04
      39,       0,  239061,  119530,2.74926e-04,7.08501e-04,7.04673e-04,       Z,5.51952e-04
      40,       0,  239061,  119530,2.74968e-04,6.48812e-04,6.55580e-04,       R,5.61064e-04
      40,       1,  342357,  171178,2.29749e-04,2.66512e-02,6.11468e-04,       Z,5.60682e-04
      41,       0,  342357,  171178,2.29770e-04,5.66656e-04,5.68493e-04,       Z,5.68885e-04
      42,       0,  342357,  171178,2.29792e-04,5.19981e-04,5.30251e-04,       R,5.76591e-04
      42,       1,  480575,  240287,1.93994e-04,1.96721e-02,5.12206e-04,       Z,5.76857e-04
      43,       0,  480575,  240287,1.94006e-04,4.74043e-04,4.76777e-04,       R,5.84355e-04
      43,       1,  676569,  338284,1.63823e-04,2.15554e-02,4.94339e-04,       Z,5.84528e-04
      44,       0,  676569,  338284,1.63837e-04,4.66400e-04,4.58777e-04,       Z,6.32554e-04
      45,       0,  676569,  338284,1.63852e-04,4.28519e-04,4.26941e-04,       R,7.47702e-04
      45,       1,  942389,  471194,1.38727e-04,1.34932e-02,4.26901e-04,       Z,7.51910e-04
      46,       0,  942389,  471194,1.38741e-04,4.03728e-04,3.96564e-04,       R,8.69136e-04
      46,       1, 1325243,  662621,1.17131e-04,1.53014e-02,3.84460e-04,       Z,9.42473e-04
        }\tableDataOne

        \pgfplotstableread[col sep=comma]{
       k,     ell,    ndof,   nElem,        eta,         mu,        res,    case,   maxGradU
       0,       0,     193,      96,3.83765e-03,2.24371e-01,2.68192e-02,       Z,3.20455e-06
       1,       0,     193,      96,1.76092e-02,2.02283e-02,1.86884e-02,       R,6.92304e-05
       1,       1,     225,     112,1.40392e-02,2.25542e-01,1.92694e-02,       Z,1.06035e-04
       2,       0,     225,     112,1.49157e-02,1.21990e-02,2.11419e-02,       R,2.22957e-04
       2,       1,     345,     172,1.10390e-02,1.57935e-02,3.30442e-02,       Z,2.37078e-04
       3,       0,     345,     172,1.10203e-02,3.11523e-02,1.10049e-01,       Z,3.41849e-04
       4,       0,     345,     172,1.10152e-02,1.09497e-01,3.93624e-01,       Z,4.18875e-04
       5,       0,     345,     172,1.10152e-02,3.93470e-01,1.41498e+00,       R,4.73358e-04
       5,       1,     441,     220,8.21264e-03,4.24293e-01,1.41497e+00,       Z,4.49899e-04
       6,       0,     441,     220,5.46677e-03,1.41495e+00,5.08818e+00,       Z,3.39112e-04
       7,       0,     441,     220,5.48198e-03,5.08817e+00,1.82973e+01,       Z,2.98044e-04
       8,       0,     441,     220,5.50064e-03,1.82973e+01,6.57981e+01,       Z,2.71977e-04
       9,       0,     441,     220,5.51452e-03,6.57981e+01,2.36614e+02,       Z,2.55331e-04
      10,       0,     441,     220,5.52433e-03,2.36614e+02,8.50875e+02,       Z,2.44766e-04
      11,       0,     441,     220,5.54627e-03,8.50875e+02,3.05979e+03,       R,2.38128e-04
      11,       1,     843,     421,4.13713e-03,8.50875e+02,3.05979e+03,       Z,2.82448e-04
      12,       0,     843,     421,5.49266e-03,3.05979e+03,1.10032e+04,       R,3.48199e-04
      12,       1,    1273,     636,6.32554e-03,3.05979e+03,1.10032e+04,       R,3.47697e-04
      12,       2,    1765,     882,8.77631e-03,3.05979e+03,1.10032e+04,       R,3.49350e-04
      12,       3,    2741,    1370,1.16910e-02,3.05979e+03,1.10032e+04,       R,3.61375e-04
      12,       4,    4185,    2092,1.57102e-02,3.05979e+03,1.10032e+04,       R,3.63060e-04
      12,       5,    6725,    3362,2.15666e-02,3.05979e+03,1.10032e+04,       R,3.62964e-04
      12,       6,   11397,    5698,2.93051e-02,3.05979e+03,1.10032e+04,       R,3.62230e-04
      12,       7,   19701,    9850,3.93216e-02,3.05979e+03,1.10032e+04,       R,3.62304e-04
      12,       8,   35017,   17508,5.42982e-02,3.05979e+03,1.10032e+04,       R,3.62380e-04
      12,       9,   62957,   31478,7.39439e-02,3.05979e+03,1.10032e+04,       R,3.62428e-04
      12,      10,  114297,   57148,1.00615e-01,3.05979e+03,1.10032e+04,       R,3.62912e-04
      12,      11,  208925,  104462,1.36780e-01,3.05979e+03,1.10032e+04,       R,3.62952e-04
      12,      12,  384913,  192456,1.87376e-01,3.05979e+03,1.10032e+04,       R,3.62975e-04
      12,      13,  710037,  355018,2.55759e-01,3.05979e+03,1.10032e+04,       R,3.62985e-04
      12,      14, 1313645,  656822,3.49155e-01,3.05979e+03,1.10032e+04,       Z,3.63023e-04
        }\tableDataTwo

        \pgfplotstableread[col sep=comma]{
       k,     ell,    ndof,   nElem,        eta,         mu,        res,    case,   maxGradU
       0,       0,     193,      96,8.35136e-04,1.04670e-01,2.68301e-02,       Z,3.20733e-06
       1,       0,     193,      96,1.76169e-02,2.02361e-02,2.89589e-02,       R,1.01316e-03
       1,       1,     225,     112,1.43463e-02,1.07099e-01,3.49659e-02,       Z,1.87614e-03
       2,       0,     225,     112,1.52211e-02,3.14780e-02,4.20552e-01,       R,7.26745e-04
       2,       1,     347,     173,1.12831e-02,3.30943e-02,4.48512e-01,       Z,8.00017e-04
       3,       0,     347,     173,1.08143e-02,4.48381e-01,8.99728e+00,       Z,1.71979e-03
       4,       0,     347,     173,1.13506e-02,8.99727e+00,1.81058e+02,       Z,1.68842e-03
       5,       0,     347,     173,1.12201e-02,1.81058e+02,3.64356e+03,       R,2.33023e-03
       5,       1,     447,     223,9.11745e-03,1.81058e+02,3.64356e+03,       Z,2.36025e-03
       6,       0,     447,     223,9.64935e-03,3.64356e+03,7.33221e+04,       R,2.11819e-03
       6,       1,     607,     303,1.21677e-02,3.64356e+03,7.33221e+04,       R,2.71755e-03
       6,       2,     783,     391,1.64161e-02,3.64356e+03,7.33221e+04,       R,3.52731e-03
       6,       3,    1111,     555,2.31663e-02,3.64356e+03,7.33221e+04,       R,4.02228e-03
       6,       4,    1619,     809,3.21667e-02,3.64356e+03,7.33221e+04,       R,4.05056e-03
       6,       5,    2591,    1295,4.58869e-02,3.64356e+03,7.33221e+04,       R,4.20678e-03
       6,       6,    4263,    2131,6.35945e-02,3.64356e+03,7.33221e+04,       R,4.20669e-03
       6,       7,    7475,    3737,8.53735e-02,3.64356e+03,7.33221e+04,       R,4.42122e-03
       6,       8,   13055,    6527,1.15370e-01,3.64356e+03,7.33221e+04,       R,4.65698e-03
       6,       9,   23411,   11705,1.58243e-01,3.64356e+03,7.33221e+04,       R,4.69263e-03
       6,      10,   42555,   21277,2.16069e-01,3.64356e+03,7.33221e+04,       R,4.75783e-03
       6,      11,   77819,   38909,2.95827e-01,3.64356e+03,7.33221e+04,       R,4.90764e-03
       6,      12,  143083,   71541,4.02623e-01,3.64356e+03,7.33221e+04,       R,5.04898e-03
       6,      13,  265059,  132529,5.47128e-01,3.64356e+03,7.33221e+04,       R,5.23535e-03
       6,      14,  488563,  244281,7.46453e-01,3.64356e+03,7.33221e+04,       R,5.22665e-03
       6,      15,  903239,  451619,1.01682e+00,3.64356e+03,7.33221e+04,       R,5.22131e-03
       6,      16, 1669943,  834971,1.38340e+00,3.64356e+03,7.33221e+04,       Z,5.18623e-03
        }\tableDataThree

        \pgfplotstableread[col sep=comma]{
       k,     ell,    ndof,   nElem,        eta,         mu,        res,    case,   maxGradU
       0,       0,     193,      96,2.08537e-02,5.22953e-01,2.67610e-02,       Z,3.18886e-06
       1,       0,     193,      96,1.76058e-02,2.01745e-02,2.45855e-02,       Z,8.90874e-06
       2,       0,     193,      96,1.75796e-02,1.72045e-02,2.30120e-02,       R,1.60775e-05
       2,       1,     225,     112,1.55083e-02,5.23480e-01,2.45746e-02,       Z,1.79351e-05
       3,       0,     225,     112,1.48710e-02,1.95725e-02,2.27618e-02,       Z,3.21293e-05
       4,       0,     225,     112,1.49004e-02,1.72159e-02,2.13664e-02,       R,4.84195e-05
       4,       1,     345,     172,1.10411e-02,1.99124e-02,1.97476e-02,       Z,4.94148e-05
       5,       0,     345,     172,1.10365e-02,1.63790e-02,1.83006e-02,       Z,6.85179e-05
       6,       0,     345,     172,1.10325e-02,1.46050e-02,1.71512e-02,       Z,8.87240e-05
       7,       0,     345,     172,1.10289e-02,1.31390e-02,1.62050e-02,       R,1.09526e-04
       7,       1,     441,     220,5.08737e-03,3.70246e-01,7.03863e-03,       Z,8.64959e-05
       8,       0,     441,     220,5.47820e-03,4.42067e-03,6.73487e-03,       Z,8.94852e-05
       9,       0,     441,     220,5.48099e-03,3.91557e-03,6.49878e-03,       Z,9.78448e-05
      10,       0,     441,     220,5.48359e-03,3.48981e-03,6.30791e-03,       Z,1.05865e-04
      11,       0,     441,     220,5.48608e-03,3.11563e-03,6.15388e-03,       Z,1.13527e-04
      12,       0,     441,     220,5.48846e-03,2.78593e-03,6.02988e-03,       Z,1.20968e-04
      13,       0,     441,     220,5.49075e-03,2.49512e-03,5.93028e-03,       R,1.28032e-04
      13,       1,     831,     415,4.63990e-03,2.06925e-01,6.33443e-03,       Z,1.39361e-04
      14,       0,     831,     415,4.50775e-03,4.45144e-03,6.03940e-03,       Z,1.55852e-04
      15,       0,     831,     415,4.50948e-03,4.01851e-03,5.80277e-03,       R,1.71945e-04
      15,       1,    1243,     621,3.61696e-03,3.30704e-02,5.49530e-03,       Z,1.72057e-04
      16,       0,    1243,     621,3.61650e-03,4.13790e-03,5.19681e-03,       Z,1.87768e-04
      17,       0,    1243,     621,3.61927e-03,3.72970e-03,4.95317e-03,       R,2.03322e-04
      17,       1,    2005,    1002,2.96173e-03,1.22468e-01,5.01396e-03,       Z,2.06893e-04
      18,       0,    2005,    1002,2.94173e-03,4.06048e-03,4.70641e-03,       Z,2.24945e-04
      19,       0,    2005,    1002,2.94201e-03,3.67374e-03,4.44968e-03,       R,2.42439e-04
      19,       1,    3229,    1614,2.37120e-03,9.54555e-02,4.63486e-03,       Z,2.43921e-04
      20,       0,    3229,    1614,2.35818e-03,3.99018e-03,4.32523e-03,       Z,2.61921e-04
      21,       0,    3229,    1614,2.35866e-03,3.62560e-03,4.06196e-03,       R,2.79294e-04
      21,       1,    4761,    2380,1.95196e-03,5.43818e-02,4.02740e-03,       Z,2.79695e-04
      22,       0,    4761,    2380,1.94977e-03,3.52400e-03,3.75705e-03,       Z,2.96773e-04
      23,       0,    4761,    2380,1.95014e-03,3.21133e-03,3.52456e-03,       R,3.13198e-04
      23,       1,    7149,    3574,1.57167e-03,5.18314e-02,3.28262e-03,       Z,3.12458e-04
      24,       0,    7149,    3574,1.57402e-03,2.88067e-03,3.06809e-03,       Z,3.27628e-04
      25,       0,    7149,    3574,1.57437e-03,2.63339e-03,2.88231e-03,       R,3.42159e-04
      25,       1,   10763,    5381,1.28389e-03,9.03677e-02,2.74644e-03,       Z,3.41793e-04
      26,       0,   10763,    5381,1.28288e-03,2.42842e-03,2.56681e-03,       Z,3.55346e-04
      27,       0,   10763,    5381,1.28298e-03,2.22318e-03,2.41095e-03,       R,3.68279e-04
      27,       1,   16161,    8080,1.05000e-03,5.67047e-02,2.38823e-03,       Z,3.68495e-04
      28,       0,   16161,    8080,1.04895e-03,2.14555e-03,2.23108e-03,       Z,3.81006e-04
      29,       0,   16161,    8080,1.04899e-03,1.96910e-03,2.09349e-03,       R,3.93259e-04
      29,       1,   24761,   12380,8.68790e-04,6.84501e-02,2.29187e-03,       Z,3.93943e-04
      30,       0,   24761,   12380,8.65763e-04,2.12206e-03,2.13588e-03,       Z,4.06614e-04
      31,       0,   24761,   12380,8.65949e-04,1.95246e-03,1.99773e-03,       R,4.18687e-04
      31,       1,   36009,   18004,7.10269e-04,4.03204e-02,1.87156e-03,       Z,4.19573e-04
      32,       0,   36009,   18004,7.11154e-04,1.73119e-03,1.74824e-03,       Z,4.31849e-04
      33,       0,   36009,   18004,7.11262e-04,1.59702e-03,1.63863e-03,       R,4.43555e-04
      33,       1,   53437,   26718,5.74991e-04,6.09442e-02,1.41452e-03,       Z,4.42019e-04
      34,       0,   53437,   26718,5.76466e-04,1.29172e-03,1.32539e-03,       Z,4.51929e-04
      35,       0,   53437,   26718,5.76539e-04,1.19343e-03,1.24629e-03,       R,4.62101e-04
      35,       1,   78153,   39076,4.76076e-04,2.38738e-02,1.24449e-03,       Z,4.62110e-04
      36,       0,   78153,   39076,4.75927e-04,1.14989e-03,1.16558e-03,       R,4.71917e-04
      36,       1,  113845,   56922,3.96056e-04,1.87745e-02,1.16821e-03,       Z,4.72138e-04
      37,       0,  113845,   56922,3.95938e-04,1.09907e-03,1.09188e-03,       Z,4.82516e-04
      38,       0,  113845,   56922,3.95952e-04,1.01756e-03,1.02335e-03,       R,4.92494e-04
      38,       1,  165227,   82613,3.30922e-04,3.92324e-02,1.04786e-03,       Z,4.93139e-04
      39,       0,  165227,   82613,3.30703e-04,9.94305e-04,9.79634e-04,       Z,5.03152e-04
      40,       0,  165227,   82613,3.30734e-04,9.22116e-04,9.18109e-04,       R,5.12669e-04
      40,       1,  235657,  117828,2.77085e-04,4.68081e-02,9.23527e-04,       Z,5.12842e-04
      41,       0,  235657,  117828,2.76936e-04,8.81028e-04,8.64085e-04,       Z,5.22025e-04
      42,       0,  235657,  117828,2.76970e-04,8.18493e-04,8.10194e-04,       R,5.30995e-04
      42,       1,  337059,  168529,2.31223e-04,2.71349e-02,7.77994e-04,       Z,5.30721e-04
      43,       0,  337059,  168529,2.31296e-04,7.42817e-04,7.28794e-04,       R,5.39197e-04
      43,       1,  474421,  237210,1.95102e-04,2.42015e-02,7.08368e-04,       Z,5.39357e-04
      44,       0,  474421,  237210,1.95129e-04,6.80962e-04,6.63156e-04,       Z,5.47606e-04
      45,       0,  474421,  237210,1.95140e-04,6.33795e-04,6.21965e-04,       R,5.55448e-04
      45,       1,  667197,  333598,1.65025e-04,2.03774e-02,6.37915e-04,       Z,5.55580e-04
      46,       0,  667197,  333598,1.64985e-04,6.16211e-04,5.97378e-04,       R,5.85521e-04
      46,       1,  928549,  464274,1.39737e-04,1.57490e-02,5.99117e-04,       Z,5.88662e-04
      47,       0,  928549,  464274,1.39739e-04,5.82593e-04,5.60794e-04,       Z,6.90533e-04
      48,       0,  928549,  464274,1.39752e-04,5.43102e-04,5.25678e-04,       R,7.91353e-04
      48,       1, 1304055,  652027,1.17877e-04,1.63165e-02,5.15052e-04,       Z,8.51239e-04
        }\tableDataFour

        %
        %

        \addplot+ [marker1, adaptive, forget plot]
        table [col sep=comma, x=cumulativeNdof, y=res] {\tableDataOne};

        \addplot+ [marker2, adaptive, forget plot]
        table [col sep=comma, x=cumulativeNdof, y=res] {\tableDataTwo};

        \addplot+ [marker3, adaptive, forget plot]
        table [col sep=comma, x=cumulativeNdof, y=res] {\tableDataThree};

        \addplot+ [marker4, adaptive, forget plot]
        table [col sep=comma, x=cumulativeNdof, y=res] {\tableDataFour};

        \drawslopetriangle[ST1]{0.5}{1.5e5}{5e-4} 
    \end{loglogaxis}
\end{tikzpicture}

%% file: figures/Fig08b_PorousMedium_condition_weighting.tex
\begin{tikzpicture}[>=stealth]
    \begin{loglogaxis}[%
            width            = 5.5cm,%
            xlabel           = {ndof},%
            ylabel           = {condition number},%
            ymajorgrids      = true,%
            font             = \footnotesize,%
            grid style       = {%
                densely dotted,%
                semithick%
            },%
            legend style     = {%
                legend pos = north west,%
                font = \footnotesize%
            },%
        ]

        \addlegendimage{marker1}
        \addlegendentry{emph.\ grad.}
        \addlegendimage{marker2}
        \addlegendentry{balanced w.}
        \addlegendimage{marker3}
        \addlegendentry{downsc. flux}
        \addlegendimage{marker4}
        \addlegendentry{split w.}

        \pgfplotstableread[col sep=comma]{
       k,     ell,    ndof,    cost,   nElem,        eta,         mu,        res,    case,   maxGradU,       cond
       0,       0,     193,     193,      96,1.76237e-02,4.80764e-01,2.67715e-02,       Z,3.19183e-06,1.38285e+05
       1,       0,     193,     386,      96,1.76012e-02,2.01838e-02,2.42249e-02,       R,1.02740e-05,1.38285e+05
       1,       1,     225,     611,     112,1.48393e-02,4.81384e-01,2.57419e-02,       Z,1.29884e-05,3.54673e+05
       2,       0,     225,     836,     112,1.48714e-02,2.10175e-02,2.33022e-02,       R,2.93689e-05,3.54673e+05
       2,       1,     345,    1181,     172,1.10176e-02,2.32703e-02,2.17484e-02,       Z,3.04053e-05,6.48398e+05
       3,       0,     345,    1526,     172,1.10154e-02,1.87547e-02,1.96280e-02,       Z,5.17603e-05,6.48398e+05
       4,       0,     345,    1871,     172,1.10123e-02,1.62503e-02,1.79983e-02,       Z,7.53309e-05,6.48398e+05
       5,       0,     345,    2216,     172,1.10099e-02,1.42411e-02,1.67160e-02,       R,1.00032e-04,6.48398e+05
       5,       1,     441,    2657,     220,5.45785e-03,3.40434e-01,6.73712e-03,       Z,7.78977e-05,1.57838e+06
       6,       0,     441,    3098,     220,5.46326e-03,3.94374e-03,6.46451e-03,       Z,8.21723e-05,1.57838e+06
       7,       0,     441,    3539,     220,5.46831e-03,3.44950e-03,6.25401e-03,       Z,9.28305e-05,1.57838e+06
       8,       0,     441,    3980,     220,5.47298e-03,3.02830e-03,6.09144e-03,       Z,1.02934e-04,1.57838e+06
       9,       0,     441,    4421,     220,5.47731e-03,2.66763e-03,5.96583e-03,       Z,1.12453e-04,1.57838e+06
      10,       0,     441,    4862,     220,5.48132e-03,2.35748e-03,5.86866e-03,       Z,1.21378e-04,1.57838e+06
      11,       0,     441,    5303,     220,5.48505e-03,2.08971e-03,5.79336e-03,       Z,1.29716e-04,1.57838e+06
      12,       0,     441,    5744,     220,5.48851e-03,1.85763e-03,5.73490e-03,       R,1.37479e-04,1.57838e+06
      12,       1,     843,    6587,     421,4.45560e-03,1.90211e-01,6.24638e-03,       Z,1.49848e-04,5.85406e+06
      13,       0,     843,    7430,     421,4.45840e-03,4.37577e-03,5.91846e-03,       Z,1.69463e-04,5.85406e+06
      14,       0,     843,    8273,     421,4.46092e-03,3.89047e-03,5.65479e-03,       R,1.88423e-04,5.85406e+06
      14,       1,    1259,    9532,     629,3.58005e-03,3.04388e-02,5.32612e-03,       Z,1.88484e-04,8.42979e+06
      15,       0,    1259,   10791,     629,3.58390e-03,3.94024e-03,5.00618e-03,       Z,2.06975e-04,8.42979e+06
      16,       0,    1259,   12050,     629,3.58733e-03,3.49216e-03,4.74874e-03,       R,2.24759e-04,8.42979e+06
      16,       1,    1993,   14043,     996,2.94840e-03,1.12578e-01,4.83135e-03,       Z,2.27896e-04,1.99747e+07
      17,       0,    1993,   16036,     996,2.94947e-03,3.82671e-03,4.50515e-03,       R,2.48236e-04,1.99747e+07
      17,       1,    3137,   19173,    1568,2.40286e-03,8.64673e-02,4.73193e-03,       Z,2.50388e-04,3.53323e+07
      18,       0,    3137,   22310,    1568,2.40409e-03,4.07578e-03,4.35844e-03,       Z,2.72774e-04,3.53323e+07
      19,       0,    3137,   25447,    1568,2.40520e-03,3.63477e-03,4.04661e-03,       R,2.94294e-04,3.53323e+07
      19,       1,    4701,   30148,    2350,1.95592e-03,5.22247e-02,3.97627e-03,       Z,2.93042e-04,4.11772e+07
      20,       0,    4701,   34849,    2350,1.95643e-03,3.46170e-03,3.66593e-03,       Z,3.12382e-04,4.11772e+07
      21,       0,    4701,   39550,    2350,1.95693e-03,3.09995e-03,3.40428e-03,       R,3.30792e-04,4.11772e+07
      21,       1,    7251,   46801,    3625,1.55546e-03,5.63487e-02,3.07587e-03,       Z,3.29786e-04,8.39699e+07
      22,       0,    7251,   54052,    3625,1.55578e-03,2.65343e-03,2.84882e-03,       Z,3.46318e-04,8.39699e+07
      23,       0,    7251,   61303,    3625,1.55610e-03,2.38630e-03,2.65618e-03,       Z,3.61996e-04,8.39699e+07
      24,       0,    7251,   68554,    3625,1.55642e-03,2.15243e-03,2.49257e-03,       R,3.76845e-04,8.39699e+07
      24,       1,   10933,   79487,    5466,1.27120e-03,7.44506e-02,2.37738e-03,       Z,3.76063e-04,3.36194e+08
      25,       0,   10933,   90420,    5466,1.27137e-03,2.00888e-03,2.21570e-03,       R,3.89914e-04,3.36194e+08
      25,       1,   16415,  106835,    8207,1.03975e-03,4.99078e-02,2.18228e-03,       Z,3.90206e-04,9.48598e+08
      26,       0,   16415,  123250,    8207,1.03975e-03,1.91867e-03,2.02124e-03,       Z,4.03805e-04,9.48598e+08
      27,       0,   16415,  139665,    8207,1.03978e-03,1.73330e-03,1.88320e-03,       R,4.16659e-04,9.48598e+08
      27,       1,   25087,  164752,   12543,8.57826e-04,6.33657e-02,2.04078e-03,       Z,4.17623e-04,3.83432e+09
      28,       0,   25087,  189839,   12543,8.57987e-04,1.85166e-03,1.88476e-03,       Z,4.31062e-04,3.83432e+09
      29,       0,   25087,  214926,   12543,8.58152e-04,1.67807e-03,1.74936e-03,       R,4.45702e-04,3.83432e+09
      29,       1,   36521,  251447,   18260,7.05036e-04,3.08140e-02,1.65144e-03,       Z,4.45472e-04,1.52701e+10
      30,       0,   36521,  287968,   18260,7.05155e-04,1.49332e-03,1.52996e-03,       Z,4.59519e-04,1.52701e+10
      31,       0,   36521,  324489,   18260,7.05277e-04,1.35771e-03,1.42406e-03,       R,4.72768e-04,1.52701e+10
      31,       1,   53831,  378320,   26915,5.73596e-04,5.16900e-02,1.22805e-03,       Z,4.69312e-04,2.52395e+10
      32,       0,   53831,  432151,   26915,5.73684e-04,1.08581e-03,1.14363e-03,       Z,4.78888e-04,2.52395e+10
      33,       0,   53831,  485982,   26915,5.73776e-04,9.89278e-04,1.07014e-03,       R,4.90407e-04,2.52395e+10
      33,       1,   79435,  565417,   39717,4.73099e-04,2.34949e-02,1.06144e-03,       Z,4.90488e-04,6.17517e+10
      34,       0,   79435,  644852,   39717,4.73162e-04,9.50146e-04,9.87682e-04,       Z,5.01761e-04,6.17517e+10
      35,       0,   79435,  724287,   39717,4.73229e-04,8.66931e-04,9.23137e-04,       R,5.12375e-04,6.17517e+10
      35,       1,  115071,  839358,   57535,3.94404e-04,1.99087e-02,9.24222e-04,       Z,5.12524e-04,2.48468e+11
      36,       0,  115071,  954429,   57535,3.94424e-04,8.35833e-04,8.59752e-04,       Z,5.22655e-04,2.48468e+11
      37,       0,  115071, 1069500,   57535,3.94448e-04,7.63927e-04,8.03010e-04,       R,5.32187e-04,2.48468e+11
      37,       1,  167251, 1236751,   83625,3.28930e-04,3.77339e-02,8.32997e-04,       Z,5.32933e-04,9.81977e+11
      38,       0,  167251, 1404002,   83625,3.28966e-04,7.65288e-04,7.73685e-04,       R,5.42593e-04,9.81977e+11
      38,       1,  239061, 1643063,  119530,2.74886e-04,3.98282e-02,7.59972e-04,       Z,5.42707e-04,1.63150e+12
      39,       0,  239061, 1882124,  119530,2.74926e-04,7.08501e-04,7.04673e-04,       Z,5.51952e-04,1.63150e+12
      40,       0,  239061, 2121185,  119530,2.74968e-04,6.48812e-04,6.55580e-04,       R,5.61064e-04,1.63150e+12
      40,       1,  342357, 2463542,  171178,2.29749e-04,2.66512e-02,6.11468e-04,       Z,5.60682e-04,3.97408e+12
      41,       0,  342357, 2805899,  171178,2.29770e-04,5.66656e-04,5.68493e-04,       Z,5.68885e-04,3.97408e+12
      42,       0,  342357, 3148256,  171178,2.29792e-04,5.19981e-04,5.30251e-04,       R,5.76591e-04,3.97408e+12
      42,       1,  480575, 3628831,  240287,1.93994e-04,1.96721e-02,5.12206e-04,       Z,5.76857e-04,1.59242e+13
      43,       0,  480575, 4109406,  240287,1.94006e-04,4.74043e-04,4.76777e-04,       R,5.84355e-04,1.59242e+13
      43,       1,  676569, 4785975,  338284,1.63823e-04,2.15554e-02,4.94339e-04,       Z,5.84528e-04,2.60539e+13
      44,       0,  676569, 5462544,  338284,1.63837e-04,4.66400e-04,4.58777e-04,       Z,6.32554e-04,2.60539e+13
      45,       0,  676569, 6139113,  338284,1.63852e-04,4.28519e-04,4.26941e-04,       R,7.47702e-04,2.60539e+13
      45,       1,  942389, 7081502,  471194,1.38727e-04,1.34932e-02,4.26901e-04,       Z,7.51910e-04,6.36697e+13
      46,       0,  942389, 8023891,  471194,1.38741e-04,4.03728e-04,3.96564e-04,       R,8.69136e-04,6.36697e+13
      46,       1, 1325243, 9349134,  662621,1.17131e-04,1.53014e-02,3.84460e-04,       Z,9.42473e-04,2.54750e+14
        }\tableDataOne

        \pgfplotstableread[col sep=comma]{
       k,     ell,    ndof,    cost,   nElem,        eta,         mu,        res,    case,   maxGradU,       cond
       0,       0,     193,     193,      96,3.83765e-03,2.24371e-01,2.68192e-02,       Z,3.20455e-06,3.16814e+05
       1,       0,     193,     386,      96,1.76092e-02,2.02283e-02,1.86884e-02,       R,6.92304e-05,3.16814e+05
       1,       1,     225,     611,     112,1.40392e-02,2.25542e-01,1.92694e-02,       Z,1.06035e-04,1.00305e+06
       2,       0,     225,     836,     112,1.49157e-02,1.21990e-02,2.11419e-02,       R,2.22957e-04,1.00305e+06
       2,       1,     345,    1181,     172,1.10390e-02,1.57935e-02,3.30442e-02,       Z,2.37078e-04,1.77130e+06
       3,       0,     345,    1526,     172,1.10203e-02,3.11523e-02,1.10049e-01,       Z,3.41849e-04,1.77130e+06
       4,       0,     345,    1871,     172,1.10152e-02,1.09497e-01,3.93624e-01,       Z,4.18875e-04,1.77130e+06
       5,       0,     345,    2216,     172,1.10152e-02,3.93470e-01,1.41498e+00,       R,4.73358e-04,1.77130e+06
       5,       1,     441,    2657,     220,8.21264e-03,4.24293e-01,1.41497e+00,       Z,4.49899e-04,4.62325e+06
       6,       0,     441,    3098,     220,5.46677e-03,1.41495e+00,5.08818e+00,       Z,3.39112e-04,4.62325e+06
       7,       0,     441,    3539,     220,5.48198e-03,5.08817e+00,1.82973e+01,       Z,2.98044e-04,4.62325e+06
       8,       0,     441,    3980,     220,5.50064e-03,1.82973e+01,6.57981e+01,       Z,2.71977e-04,4.62325e+06
       9,       0,     441,    4421,     220,5.51452e-03,6.57981e+01,2.36614e+02,       Z,2.55331e-04,4.62325e+06
      10,       0,     441,    4862,     220,5.52433e-03,2.36614e+02,8.50875e+02,       Z,2.44766e-04,4.62325e+06
      11,       0,     441,    5303,     220,5.54627e-03,8.50875e+02,3.05979e+03,       R,2.38128e-04,4.62325e+06
      11,       1,     843,    6146,     421,4.13713e-03,8.50875e+02,3.05979e+03,       Z,2.82448e-04,1.69556e+07
      12,       0,     843,    6989,     421,5.49266e-03,3.05979e+03,1.10032e+04,       R,3.48199e-04,1.69556e+07
      12,       1,    1273,    8262,     636,6.32554e-03,3.05979e+03,1.10032e+04,       R,3.47697e-04,6.44018e+07
      12,       2,    1765,   10027,     882,8.77631e-03,3.05979e+03,1.10032e+04,       R,3.49350e-04,2.55086e+08
      12,       3,    2741,   12768,    1370,1.16910e-02,3.05979e+03,1.10032e+04,       R,3.61375e-04,1.00951e+09
      12,       4,    4185,   16953,    2092,1.57102e-02,3.05979e+03,1.10032e+04,       R,3.63060e-04,4.09812e+09
      12,       5,    6725,   23678,    3362,2.15666e-02,3.05979e+03,1.10032e+04,       R,3.62964e-04,1.62313e+10
      12,       6,   11397,   35075,    5698,2.93051e-02,3.05979e+03,1.10032e+04,       R,3.62230e-04,7.17736e+10
      12,       7,   19701,   54776,    9850,3.93216e-02,3.05979e+03,1.10032e+04,       R,3.62304e-04,2.63050e+11
      12,       8,   35017,   89793,   17508,5.42982e-02,3.05979e+03,1.10032e+04,       R,3.62380e-04,1.04391e+12
      12,       9,   62957,  152750,   31478,7.39439e-02,3.05979e+03,1.10032e+04,       R,3.62428e-04,2.81238e+12
      12,      10,  114297,  267047,   57148,1.00615e-01,3.05979e+03,1.10032e+04,       R,3.62912e-04,1.13612e+13
      12,      11,  208925,  475972,  104462,1.36780e-01,3.05979e+03,1.10032e+04,       R,3.62952e-04,2.29773e+13
      12,      12,  384913,  860885,  192456,1.87376e-01,3.05979e+03,1.10032e+04,       R,3.62975e-04,1.39234e+14
      12,      13,  710037, 1570922,  355018,2.55759e-01,3.05979e+03,1.10032e+04,       R,3.62985e-04,4.09578e+14
      12,      14, 1313645, 2884567,  656822,3.49155e-01,3.05979e+03,1.10032e+04,       Z,3.63023e-04,1.47819e+15
        }\tableDataTwo

        \pgfplotstableread[col sep=comma]{
       k,     ell,    ndof,    cost,   nElem,        eta,         mu,        res,    case,   maxGradU,       cond
       0,       0,     193,     193,      96,8.35136e-04,1.04670e-01,2.68301e-02,       Z,3.20733e-06,1.38269e+05
       1,       0,     193,     386,      96,1.76169e-02,2.02361e-02,2.89589e-02,       R,1.01316e-03,1.38269e+05
       1,       1,     225,     611,     112,1.43463e-02,1.07099e-01,3.49659e-02,       Z,1.87614e-03,3.54664e+05
       2,       0,     225,     836,     112,1.52211e-02,3.14780e-02,4.20552e-01,       R,7.26745e-04,3.54664e+05
       2,       1,     347,    1183,     173,1.12831e-02,3.30943e-02,4.48512e-01,       Z,8.00017e-04,6.48394e+05
       3,       0,     347,    1530,     173,1.08143e-02,4.48381e-01,8.99728e+00,       Z,1.71979e-03,6.48394e+05
       4,       0,     347,    1877,     173,1.13506e-02,8.99727e+00,1.81058e+02,       Z,1.68842e-03,6.48394e+05
       5,       0,     347,    2224,     173,1.12201e-02,1.81058e+02,3.64356e+03,       R,2.33023e-03,6.48394e+05
       5,       1,     447,    2671,     223,9.11745e-03,1.81058e+02,3.64356e+03,       Z,2.36025e-03,1.54279e+06
       6,       0,     447,    3118,     223,9.64935e-03,3.64356e+03,7.33221e+04,       R,2.11819e-03,1.54279e+06
       6,       1,     607,    3725,     303,1.21677e-02,3.64356e+03,7.33221e+04,       R,2.71755e-03,5.59893e+06
       6,       2,     783,    4508,     391,1.64161e-02,3.64356e+03,7.33221e+04,       R,3.52731e-03,2.27275e+07
       6,       3,    1111,    5619,     555,2.31663e-02,3.64356e+03,7.33221e+04,       R,4.02228e-03,8.91504e+07
       6,       4,    1619,    7238,     809,3.21667e-02,3.64356e+03,7.33221e+04,       R,4.05056e-03,3.58695e+08
       6,       5,    2591,    9829,    1295,4.58869e-02,3.64356e+03,7.33221e+04,       R,4.20678e-03,1.43296e+09
       6,       6,    4263,   14092,    2131,6.35945e-02,3.64356e+03,7.33221e+04,       R,4.20669e-03,5.93037e+09
       6,       7,    7475,   21567,    3737,8.53735e-02,3.64356e+03,7.33221e+04,       R,4.42122e-03,2.30293e+10
       6,       8,   13055,   34622,    6527,1.15370e-01,3.64356e+03,7.33221e+04,       R,4.65698e-03,9.21967e+10
       6,       9,   23411,   58033,   11705,1.58243e-01,3.64356e+03,7.33221e+04,       R,4.69263e-03,2.58123e+11
       6,      10,   42555,  100588,   21277,2.16069e-01,3.64356e+03,7.33221e+04,       R,4.75783e-03,1.06538e+12
       6,      11,   77819,  178407,   38909,2.95827e-01,3.64356e+03,7.33221e+04,       R,4.90764e-03,4.68866e+12
       6,      12,  143083,  321490,   71541,4.02623e-01,3.64356e+03,7.33221e+04,       R,5.04898e-03,1.67004e+13
       6,      13,  265059,  586549,  132529,5.47128e-01,3.64356e+03,7.33221e+04,       R,5.23535e-03,7.45453e+13
       6,      14,  488563, 1075112,  244281,7.46453e-01,3.64356e+03,7.33221e+04,       R,5.22665e-03,3.02627e+14
       6,      15,  903239, 1978351,  451619,1.01682e+00,3.64356e+03,7.33221e+04,       R,5.22131e-03,1.54388e+15
       6,      16, 1669943, 3648294,  834971,1.38340e+00,3.64356e+03,7.33221e+04,       Z,5.18623e-03,4.23920e+15
        }\tableDataThree

        \pgfplotstableread[col sep=comma]{
       k,     ell,    ndof,    cost,   nElem,        eta,         mu,        res,    case,   maxGradU,       cond
       0,       0,     193,     193,      96,2.08537e-02,5.22953e-01,2.67610e-02,       Z,3.18886e-06,1.38288e+05
       1,       0,     193,     386,      96,1.76058e-02,2.01745e-02,2.45855e-02,       Z,8.90874e-06,1.38288e+05
       2,       0,     193,     579,      96,1.75796e-02,1.72045e-02,2.30120e-02,       R,1.60775e-05,1.38288e+05
       2,       1,     225,     804,     112,1.55083e-02,5.23480e-01,2.45746e-02,       Z,1.79351e-05,3.54675e+05
       3,       0,     225,    1029,     112,1.48710e-02,1.95725e-02,2.27618e-02,       Z,3.21293e-05,3.54675e+05
       4,       0,     225,    1254,     112,1.49004e-02,1.72159e-02,2.13664e-02,       R,4.84195e-05,3.54675e+05
       4,       1,     345,    1599,     172,1.10411e-02,1.99124e-02,1.97476e-02,       Z,4.94148e-05,6.48399e+05
       5,       0,     345,    1944,     172,1.10365e-02,1.63790e-02,1.83006e-02,       Z,6.85179e-05,6.48399e+05
       6,       0,     345,    2289,     172,1.10325e-02,1.46050e-02,1.71512e-02,       Z,8.87240e-05,6.48399e+05
       7,       0,     345,    2634,     172,1.10289e-02,1.31390e-02,1.62050e-02,       R,1.09526e-04,6.48399e+05
       7,       1,     441,    3075,     220,5.08737e-03,3.70246e-01,7.03863e-03,       Z,8.64959e-05,1.57838e+06
       8,       0,     441,    3516,     220,5.47820e-03,4.42067e-03,6.73487e-03,       Z,8.94852e-05,1.57838e+06
       9,       0,     441,    3957,     220,5.48099e-03,3.91557e-03,6.49878e-03,       Z,9.78448e-05,1.57838e+06
      10,       0,     441,    4398,     220,5.48359e-03,3.48981e-03,6.30791e-03,       Z,1.05865e-04,1.57838e+06
      11,       0,     441,    4839,     220,5.48608e-03,3.11563e-03,6.15388e-03,       Z,1.13527e-04,1.57838e+06
      12,       0,     441,    5280,     220,5.48846e-03,2.78593e-03,6.02988e-03,       Z,1.20968e-04,1.57838e+06
      13,       0,     441,    5721,     220,5.49075e-03,2.49512e-03,5.93028e-03,       R,1.28032e-04,1.57838e+06
      13,       1,     831,    6552,     415,4.63990e-03,2.06925e-01,6.33443e-03,       Z,1.39361e-04,6.03233e+06
      14,       0,     831,    7383,     415,4.50775e-03,4.45144e-03,6.03940e-03,       Z,1.55852e-04,6.03233e+06
      15,       0,     831,    8214,     415,4.50948e-03,4.01851e-03,5.80277e-03,       R,1.71945e-04,6.03233e+06
      15,       1,    1243,    9457,     621,3.61696e-03,3.30704e-02,5.49530e-03,       Z,1.72057e-04,8.56033e+06
      16,       0,    1243,   10700,     621,3.61650e-03,4.13790e-03,5.19681e-03,       Z,1.87768e-04,8.56033e+06
      17,       0,    1243,   11943,     621,3.61927e-03,3.72970e-03,4.95317e-03,       R,2.03322e-04,8.56033e+06
      17,       1,    2005,   13948,    1002,2.96173e-03,1.22468e-01,5.01396e-03,       Z,2.06893e-04,1.79242e+07
      18,       0,    2005,   15953,    1002,2.94173e-03,4.06048e-03,4.70641e-03,       Z,2.24945e-04,1.79242e+07
      19,       0,    2005,   17958,    1002,2.94201e-03,3.67374e-03,4.44968e-03,       R,2.42439e-04,1.79242e+07
      19,       1,    3229,   21187,    1614,2.37120e-03,9.54555e-02,4.63486e-03,       Z,2.43921e-04,3.56623e+07
      20,       0,    3229,   24416,    1614,2.35818e-03,3.99018e-03,4.32523e-03,       Z,2.61921e-04,3.56623e+07
      21,       0,    3229,   27645,    1614,2.35866e-03,3.62560e-03,4.06196e-03,       R,2.79294e-04,3.56623e+07
      21,       1,    4761,   32406,    2380,1.95196e-03,5.43818e-02,4.02740e-03,       Z,2.79695e-04,4.11675e+07
      22,       0,    4761,   37167,    2380,1.94977e-03,3.52400e-03,3.75705e-03,       Z,2.96773e-04,4.11675e+07
      23,       0,    4761,   41928,    2380,1.95014e-03,3.21133e-03,3.52456e-03,       R,3.13198e-04,4.11675e+07
      23,       1,    7149,   49077,    3574,1.57167e-03,5.18314e-02,3.28262e-03,       Z,3.12458e-04,8.37268e+07
      24,       0,    7149,   56226,    3574,1.57402e-03,2.88067e-03,3.06809e-03,       Z,3.27628e-04,8.37268e+07
      25,       0,    7149,   63375,    3574,1.57437e-03,2.63339e-03,2.88231e-03,       R,3.42159e-04,8.37268e+07
      25,       1,   10763,   74138,    5381,1.28389e-03,9.03677e-02,2.74644e-03,       Z,3.41793e-04,3.36394e+08
      26,       0,   10763,   84901,    5381,1.28288e-03,2.42842e-03,2.56681e-03,       Z,3.55346e-04,3.36394e+08
      27,       0,   10763,   95664,    5381,1.28298e-03,2.22318e-03,2.41095e-03,       R,3.68279e-04,3.36394e+08
      27,       1,   16161,  111825,    8080,1.05000e-03,5.67047e-02,2.38823e-03,       Z,3.68495e-04,9.53725e+08
      28,       0,   16161,  127986,    8080,1.04895e-03,2.14555e-03,2.23108e-03,       Z,3.81006e-04,9.53725e+08
      29,       0,   16161,  144147,    8080,1.04899e-03,1.96910e-03,2.09349e-03,       R,3.93259e-04,9.53725e+08
      29,       1,   24761,  168908,   12380,8.68790e-04,6.84501e-02,2.29187e-03,       Z,3.93943e-04,3.84423e+09
      30,       0,   24761,  193669,   12380,8.65763e-04,2.12206e-03,2.13588e-03,       Z,4.06614e-04,3.84423e+09
      31,       0,   24761,  218430,   12380,8.65949e-04,1.95246e-03,1.99773e-03,       R,4.18687e-04,3.84423e+09
      31,       1,   36009,  254439,   18004,7.10269e-04,4.03204e-02,1.87156e-03,       Z,4.19573e-04,1.52648e+10
      32,       0,   36009,  290448,   18004,7.11154e-04,1.73119e-03,1.74824e-03,       Z,4.31849e-04,1.52648e+10
      33,       0,   36009,  326457,   18004,7.11262e-04,1.59702e-03,1.63863e-03,       R,4.43555e-04,1.52648e+10
      33,       1,   53437,  379894,   26718,5.74991e-04,6.09442e-02,1.41452e-03,       Z,4.42019e-04,2.52366e+10
      34,       0,   53437,  433331,   26718,5.76466e-04,1.29172e-03,1.32539e-03,       Z,4.51929e-04,2.52366e+10
      35,       0,   53437,  486768,   26718,5.76539e-04,1.19343e-03,1.24629e-03,       R,4.62101e-04,2.52366e+10
      35,       1,   78153,  564921,   39076,4.76076e-04,2.38738e-02,1.24449e-03,       Z,4.62110e-04,6.17417e+10
      36,       0,   78153,  643074,   39076,4.75927e-04,1.14989e-03,1.16558e-03,       R,4.71917e-04,6.17417e+10
      36,       1,  113845,  756919,   56922,3.96056e-04,1.87745e-02,1.16821e-03,       Z,4.72138e-04,2.48403e+11
      37,       0,  113845,  870764,   56922,3.95938e-04,1.09907e-03,1.09188e-03,       Z,4.82516e-04,2.48403e+11
      38,       0,  113845,  984609,   56922,3.95952e-04,1.01756e-03,1.02335e-03,       R,4.92494e-04,2.48403e+11
      38,       1,  165227, 1149836,   82613,3.30922e-04,3.92324e-02,1.04786e-03,       Z,4.93139e-04,9.93384e+11
      39,       0,  165227, 1315063,   82613,3.30703e-04,9.94305e-04,9.79634e-04,       Z,5.03152e-04,9.93384e+11
      40,       0,  165227, 1480290,   82613,3.30734e-04,9.22116e-04,9.18109e-04,       R,5.12669e-04,9.93384e+11
      40,       1,  235657, 1715947,  117828,2.77085e-04,4.68081e-02,9.23527e-04,       Z,5.12842e-04,1.62218e+12
      41,       0,  235657, 1951604,  117828,2.76936e-04,8.81028e-04,8.64085e-04,       Z,5.22025e-04,1.62218e+12
      42,       0,  235657, 2187261,  117828,2.76970e-04,8.18493e-04,8.10194e-04,       R,5.30995e-04,1.62218e+12
      42,       1,  337059, 2524320,  168529,2.31223e-04,2.71349e-02,7.77994e-04,       Z,5.30721e-04,3.97342e+12
      43,       0,  337059, 2861379,  168529,2.31296e-04,7.42817e-04,7.28794e-04,       R,5.39197e-04,3.97342e+12
      43,       1,  474421, 3335800,  237210,1.95102e-04,2.42015e-02,7.08368e-04,       Z,5.39357e-04,1.59233e+13
      44,       0,  474421, 3810221,  237210,1.95129e-04,6.80962e-04,6.63156e-04,       Z,5.47606e-04,1.59233e+13
      45,       0,  474421, 4284642,  237210,1.95140e-04,6.33795e-04,6.21965e-04,       R,5.55448e-04,1.59233e+13
      45,       1,  667197, 4951839,  333598,1.65025e-04,2.03774e-02,6.37915e-04,       Z,5.55580e-04,2.60531e+13
      46,       0,  667197, 5619036,  333598,1.64985e-04,6.16211e-04,5.97378e-04,       R,5.85521e-04,2.60531e+13
      46,       1,  928549, 6547585,  464274,1.39737e-04,1.57490e-02,5.99117e-04,       Z,5.88662e-04,6.36731e+13
      47,       0,  928549, 7476134,  464274,1.39739e-04,5.82593e-04,5.60794e-04,       Z,6.90533e-04,6.36731e+13
      48,       0,  928549, 8404683,  464274,1.39752e-04,5.43102e-04,5.25678e-04,       R,7.91353e-04,6.36731e+13
      48,       1, 1304055, 9708738,  652027,1.17877e-04,1.63165e-02,5.15052e-04,       Z,8.51239e-04,2.54780e+14
        }\tableDataFour

        %
        %

        \addplot+ [marker1, adaptive, forget plot]
        table [col sep=comma, x=ndof, y=cond] {\tableDataOne};

        \addplot+ [marker2, adaptive, forget plot]
        table [col sep=comma, x=ndof, y=cond] {\tableDataTwo};

        \addplot+ [marker3, adaptive, forget plot]
        table [col sep=comma, x=ndof, y=cond] {\tableDataThree};

        \addplot+ [marker4, adaptive, forget plot]
        table [col sep=comma, x=ndof, y=cond] {\tableDataFour};

        \drawslopetriangleup[ST1]{3}{1.5e4}{1e8} 
    \end{loglogaxis}
\end{tikzpicture}

%% file: literature.bib
@article {Bringmann2024,
    AUTHOR = {Bringmann, Philipp},
     TITLE = {Review and computational comparison of adaptive least-squares
              finite element schemes},
   JOURNAL = {Comput. Math. Appl.},
  FJOURNAL = {Computers \& Mathematics with Applications. An International
              Journal},
    VOLUME = {172},
      YEAR = {2024},
     PAGES = {1--15},
      ISSN = {0898-1221,1873-7668},
   MRCLASS = {65N30 (65N50)},
  MRNUMBER = {4782072},
MRREVIEWER = {Huipo\ Liu},
       DOI = {10.1016/j.camwa.2024.07.022},
       URL = {https://doi.org/10.1016/j.camwa.2024.07.022},
}

@article {Bringmann2023,
    AUTHOR = {Bringmann, Philipp},
     TITLE = {How to prove optimal convergence rates for adaptive
              least-squares finite element methods},
   JOURNAL = {J. Numer. Math.},
  FJOURNAL = {Journal of Numerical Mathematics},
    VOLUME = {31},
      YEAR = {2023},
    NUMBER = {1},
     PAGES = {43--58},
      ISSN = {1570-2820,1569-3953},
   MRCLASS = {65N12 (65N15 65N30 65N50 65Y20)},
  MRNUMBER = {4557622},
       DOI = {10.1515/jnma-2021-0116},
       URL = {https://doi.org/10.1515/jnma-2021-0116},
}

@incollection {BringmannCarstensenTran2022,
    AUTHOR = {Bringmann, Philipp and Carstensen, Carsten and Tran, Ngoc
              Tien},
     TITLE = {Adaptive least-squares, discontinuous {P}etrov-{G}alerkin, and
              hybrid high-order methods},
 BOOKTITLE = {Non-standard discretisation methods in solid mechanics},
    SERIES = {Lect. Notes Appl. Comput. Mech.},
    VOLUME = {98},
     PAGES = {107--147},
 PUBLISHER = {Springer, Cham},
      YEAR = {2022},
      ISBN = {978-3-030-92671-7; 978-3-030-92672-4},
   MRCLASS = {65N30 (49J10 65K10 65N12 65N15)},
  MRNUMBER = {4433562},
       DOI = {10.1007/978-3-030-92672-4\_5},
       URL = {https://doi.org/10.1007/978-3-030-92672-4_5},
}

@article {CarstensenBringmannHellwigWriggers2018,
    AUTHOR = {Carstensen, C. and Bringmann, P. and Hellwig, F. and Wriggers,
              P.},
     TITLE = {Nonlinear discontinuous {P}etrov-{G}alerkin methods},
   JOURNAL = {Numer. Math.},
  FJOURNAL = {Numerische Mathematik},
    VOLUME = {139},
      YEAR = {2018},
    NUMBER = {3},
     PAGES = {529--561},
      ISSN = {0029-599X,0945-3245},
   MRCLASS = {65N30 (47H05 49M15 65N12 65N15)},
  MRNUMBER = {3814605},
MRREVIEWER = {T.\ A.\ Angelov},
       DOI = {10.1007/s00211-018-0947-5},
       URL = {https://doi.org/10.1007/s00211-018-0947-5},
}

@misc {BringmannPraetorius2025_code,
    AUTHOR = {Bringmann, Philipp and Praetorius, Dirk},
     TITLE = {{octAFEM} -- Numerical investigation of an adaptive least-squares finite element method for the solution of Zarantonello-linearized first-order systems of quasilinear {PDEs}},
      NOTE = {MATLAB software package, available under DOI: \href{https://doi.org/10.24433/CO.0796231.v1}{10.24433/CO.0796231.v1}.},
       DOI = {10.24433/CO.0796231.v1},
      YEAR = {2026},
HOWPUBLISHED = {\url{https://www.codeocean.com/}},
}

@incollection {BringmannMiraciPraetorius2024,
    AUTHOR    = {Bringmann, Philipp and Mira\c{c}i, Ani and Praetorius, Dirk},
    TITLE     = {Iterative solvers in adaptive {FEM}},
    EDITOR    = {F. Chouly and S.P.A. Bordas and R. Becker and P. Omnes},
    BOOKTITLE = {Error Control, Adaptive Discretizations, and Applications. Part 2},
    SUBTITLE  = {Adaptivity yields quasi-optimal computational runtime},
    SERIES    = {Advances in Applied Mechanics (AAMS)},
    VOLUME    = {59},
    PAGES     = {147--212},
    PUBLISHER = {Elsevier},
    YEAR      = {2024},
    ISSN      = {0065-2156},
    ISBN      = {978-0-443-29448-8},
    DOI       = {10.1016/bs.aams.2024.08.002},
    % NOTE      = {Preprint under \href{https://arxiv.org/abs/2404.07126}{\ttfamily\mdseries arxiv.org:2404.07126}},
}

@article {BrunnerInnerbergerMiraciPraetoriusStreitbergerHeid2024,
    AUTHOR = {Brunner, Maximilian and Innerberger, Michael and Mira\c ci,
              Ani and Praetorius, Dirk and Streitberger, Julian and Heid,
              Pascal},
     TITLE = {Adaptive {FEM} with quasi-optimal overall cost for
              nonsymmetric linear elliptic {PDE}s},
   JOURNAL = {IMA J. Numer. Anal.},
  FJOURNAL = {IMA Journal of Numerical Analysis},
    VOLUME = {44},
      YEAR = {2024},
    NUMBER = {3},
     PAGES = {1560--1596},
      ISSN = {0272-4979,1464-3642},
   MRCLASS = {65N30 (35J15)},
  MRNUMBER = {4755062},
MRREVIEWER = {Huipo\ Liu},
       DOI = {10.1093/imanum/drad039},
       URL = {https://doi.org/10.1093/imanum/drad039},
}

@article {HaberlPraetoriusSchimankoVohralik2021,
    AUTHOR = {Haberl, Alexander and Praetorius, Dirk and Schimanko, Stefan
              and Vohral\'ik, Martin},
     TITLE = {Convergence and quasi-optimal cost of adaptive algorithms for
              nonlinear operators including iterative linearization and
              algebraic solver},
   JOURNAL = {Numer. Math.},
  FJOURNAL = {Numerische Mathematik},
    VOLUME = {147},
      YEAR = {2021},
    NUMBER = {3},
     PAGES = {679--725},
      ISSN = {0029-599X,0945-3245},
   MRCLASS = {65N30 (35J15 65N12 65N15 65N50 68Q25)},
  MRNUMBER = {4224933},
MRREVIEWER = {Mohammad\ Asadzadeh},
       DOI = {10.1007/s00211-021-01176-w},
       URL = {https://doi.org/10.1007/s00211-021-01176-w},
}

@article {FuhrerPraetorius2020,
    AUTHOR = {F\"uhrer, Thomas and Praetorius, Dirk},
     TITLE = {A short note on plain convergence of adaptive least-squares
              finite element methods},
   JOURNAL = {Comput. Math. Appl.},
  FJOURNAL = {Computers \& Mathematics with Applications. An International
              Journal},
    VOLUME = {80},
      YEAR = {2020},
    NUMBER = {6},
     PAGES = {1619--1632},
      ISSN = {0898-1221,1873-7668},
   MRCLASS = {65N30 (65N12 65N22 65N50)},
  MRNUMBER = {4138307},
MRREVIEWER = {Yunying\ Zheng},
       DOI = {10.1016/j.camwa.2020.07.022},
       URL = {https://doi.org/10.1016/j.camwa.2020.07.022},
}

@article {FuhrerPraetorius2018,
    AUTHOR = {F\"uhrer, Thomas and Praetorius, Dirk},
     TITLE = {A linear {U}zawa-type {FEM}-{BEM} solver for nonlinear
              transmission problems},
   JOURNAL = {Comput. Math. Appl.},
  FJOURNAL = {Computers \& Mathematics with Applications. An International
              Journal},
    VOLUME = {75},
      YEAR = {2018},
    NUMBER = {8},
     PAGES = {2678--2697},
      ISSN = {0898-1221,1873-7668},
   MRCLASS = {65N22 (35J25 35J62 65N30 65N38)},
  MRNUMBER = {3787479},
MRREVIEWER = {Jennifer\ Pestana},
       DOI = {10.1016/j.camwa.2017.12.035},
       URL = {https://doi.org/10.1016/j.camwa.2017.12.035},
}

@unpublished {BertrandBrodbeckRickenSchneider2025,
    AUTHOR = {Fleurianne Bertrand and Maximilian Brodbeck and Tim Ricken and Henrik Schneider},
     TITLE = {Least-squares finite element methods for nonlinear problems: {A} unified framework}, 
      YEAR = {2025},
    EPRINT = {2503.18739},
ARCHIVEPREFIX = {arXiv},
       URL ={https://arxiv.org/abs/2503.18739}, 
      NOTE = {Preprint}
}

@unpublished {KotheLoscherSteinbach2023,
    AUTHOR = {K{\"o}the, Christian and L{\"o}scher, Richard and Steinbach, Olaf},
     TITLE = {Adaptive least-squares space-time finite element methods},
      YEAR = {2023},
ARCHIVEPREFIX = {arXiv},
       DOI = {10.48550/arXiv.2309.14300},
      NOTE = {Preprint}
}

@software {quiver2,
    AUTHOR = {Vargas Aguilera, Carlos Adrian},
     TITLE = {{quiver2.m}},
      YEAR = {2009},
   VERSION = {1.2},
       URL = {https://de.mathworks.com/matlabcentral/fileexchange/24600-quiver2-m-v1-2-nov-2009},
      NOTE = {MATLAB Central File Exchange}
}

@thesis {Riveros2023,
    AUTHOR  = {Riveros Neira, Andr{\'e}s Sebasti{\'a}},
    TITLE   = {Elementos finitos m{\'i}nimos cuadrados para ecuaciones fuertemente mon{\'o}tonas},
    TYPE    = {Master Thesis (Supervisor: Prof.~Michael Karkulik)},
    SCHOOL  = {Universidad T{\'e}cnica Federico Santa Mar{\'i}a},
    YEAR    = {2023},
    ADDRESS = {Chile},
}

@article {HoonhoutLoscherSteinbachUrzuaTorres2026,
    AUTHOR = {Hoonhout, Daniel and L\"oscher, Richard and Steinbach, Olaf
              and Urz\'ua--Torres, Carolina},
     TITLE = {Stable least-squares space-time boundary element methods for
              the wave equation},
   JOURNAL = {Adv. Comput. Math.},
  FJOURNAL = {Advances in Computational Mathematics},
    VOLUME = {52},
      YEAR = {2026},
    NUMBER = {1},
     PAGES = {Paper No. 7},
      ISSN = {1019-7168,1572-9044},
   MRCLASS = {65M38 (65M12)},
  MRNUMBER = {5022352},
       DOI = {10.1007/s10444-026-10282-y},
       URL = {https://doi.org/10.1007/s10444-026-10282-y},
}

@article{KotheSteinbach2026,
    AUTHOR = {K{\"o}the, Christian and Steinbach, Olaf},
     TITLE = {Adaptive {Least}-{Squares} ({Space}-{Time}) {Finite} {Element} {Methods} for {Convection}-{Diffusion} {Problems}},
   JOURNAL = {Comput. Methods Appl. Math.},
  FJOURNAL = {Computational Methods in Applied Mathematics},
      YEAR = {2026},
     PAGES = {1--29},
       DOI = {10.1515/cmam-2025-0156},
      NOTE = {Published online}
}

@article {FuhrerGonzalezKarkulik2025,
    AUTHOR = {F{\"u}hrer, Thomas and Gonz{\'a}lez, Roberto and Karkulik, Michael},
     TITLE = {Well-posedness of first-order acoustic wave equations and space-time finite element approximation},
   JOURNAL = {IMA J. Numer. Anal.},
  FJOURNAL = {IMA Journal of Numerical Analysis},
      YEAR = {2025},
     PAGES = {1--30},
       DOI = {10.1093/imanum/drae104},
      NOTE = {Published online}
}

@article{GallistlTran2025,
    AUTHOR = {Gallistl, Dietmar and Tran, Ngoc Tien},
     TITLE = {Minimal residual discretization of a class of fully nonlinear elliptic PDE},
   JOURNAL = {IMA J. Numer. Anal.},
  FJOURNAL = {IMA Journal of Numerical Analysis},
      YEAR = {2025},
     PAGES = {1--21},
       DOI = {10.1093/imanum/draf075},
      NOTE = {Published online}
}

@article {DieningGehringStorn2025,
    AUTHOR = {Diening, Lars and Gehring, Lukas and Storn, Johannes},
     TITLE = {Adaptive {Mesh} {Refinement} for {Arbitrary} {Initial} {Triangulations}},
   JOURNAL = {Found. Comput. Math.},
  FJOURNAL = {Foundations of Computational Mathematics},
      YEAR = {2025},
     PAGES = {1--26},
      ISSN = {1615-3375},
   MRCLASS = {65N50 (65Y20)},
       DOI = {10.1007/s10208-025-09698-7},
       URL = {https://doi.org/10.1007/s10208-025-09698-7},
}

@article {MonsuurSmeetsStevenson2025,
    AUTHOR = {Monsuur, Harald and Smeets, Robin and Stevenson, Rob},
     TITLE = {Quasi-optimal least squares: {Inhomogeneous} boundary conditions, and application with machine learning},
      YEAR = {2025},
   JOURNAL = {IMA J. Numer. Anal.},
  FJOURNAL = {IMA Journal of Numerical Analysis},
     PAGES = {1--38},
       DOI = {10.1093/imanum/draf051},
      NOTE = {Published online}
}

@article {MeissnerHuynhKuberryBochev2025,
    AUTHOR = {Meissner, Teddy and Huynh, Edward and Kuberry, Paul and
              Bochev, Pavel},
     TITLE = {A deep least-squares method for the {S}tokes equations},
   JOURNAL = {Comput. Math. Appl.},
  FJOURNAL = {Computers \& Mathematics with Applications. An International
              Journal},
    VOLUME = {196},
      YEAR = {2025},
     PAGES = {1--12},
      ISSN = {0898-1221,1873-7668},
   MRCLASS = {65M70 (35Q84 65N 68T07)},
  MRNUMBER = {4931131},
       DOI = {10.1016/j.camwa.2025.07.004},
       URL = {https://doi.org/10.1016/j.camwa.2025.07.004},
}

@article {LiZhang2025,
    AUTHOR = {Li, Ziyan and Zhang, Shun},
     TITLE = {Non-intrusive least-squares functional a posteriori error
              estimator: linear and nonlinear problems with plain
              convergence},
   JOURNAL = {Comput. Math. Appl.},
  FJOURNAL = {Computers \& Mathematics with Applications. An International
              Journal},
    VOLUME = {191},
      YEAR = {2025},
     PAGES = {275--295},
      ISSN = {0898-1221,1873-7668},
   MRCLASS = {65N50 (65N12)},
  MRNUMBER = {4910387},
       DOI = {10.1016/j.camwa.2025.05.011},
       URL = {https://doi.org/10.1016/j.camwa.2025.05.011},
}

@article {EggerEngertsbergerDomenigRoppertKaltenbacher2025,
    AUTHOR = {Egger, H. and Engertsberger, F. and Domenig, L. and Roppert,
              K. and Kaltenbacher, M.},
     TITLE = {On nonlinear magnetic field solvers using local quasi-{N}ewton
              updates},
   JOURNAL = {Comput. Math. Appl.},
  FJOURNAL = {Computers \& Mathematics with Applications. An International
              Journal},
    VOLUME = {183},
      YEAR = {2025},
     PAGES = {20--31},
      ISSN = {0898-1221,1873-7668},
   MRCLASS = {78M10 (65H10 78A25 78A30 78A60)},
  MRNUMBER = {4857640},
       DOI = {10.1016/j.camwa.2025.01.033},
       URL = {https://doi.org/10.1016/j.camwa.2025.01.033},
}

@article {BrennerSungTanZhang2024,
    AUTHOR = {Brenner, Susanne C. and Sung, Li-yeng and Tan, Zhiyu and
              Zhang, Hongchao},
     TITLE = {A nonlinear least-squares convexity enforcing {$C^0$} interior
              penalty method for the {M}onge-{A}mp\`ere equation on strictly
              convex smooth planar domains},
   JOURNAL = {Commun. Am. Math. Soc.},
  FJOURNAL = {Communications of the American Mathematical Society},
    VOLUME = {4},
      YEAR = {2024},
     PAGES = {607--640},
      ISSN = {2692-3688},
   MRCLASS = {65N12 (35J96 65K10 65N15 65N30 90C30 90C90)},
  MRNUMBER = {4806850},
       DOI = {10.1090/cams/39},
       URL = {https://doi.org/10.1090/cams/39},
}

@article {Storn2024,
    AUTHOR = {Storn, Johannes},
     TITLE = {Solving minimal residual methods in {$W^{-1,p'}$} with large
              exponents {$p$}},
   JOURNAL = {J. Sci. Comput.},
  FJOURNAL = {Journal of Scientific Computing},
    VOLUME = {99},
      YEAR = {2024},
    NUMBER = {2},
     PAGES = {Paper No. 35, 18},
      ISSN = {0885-7474,1573-7691},
   MRCLASS = {65N12 (49M20 49M29 65N15 65N20)},
  MRNUMBER = {4721138},
       DOI = {10.1007/s10915-024-02494-5},
       URL = {https://doi.org/10.1007/s10915-024-02494-5},
}

@article {BertrandSchneider2024,
    AUTHOR = {Bertrand, Fleurianne and Schneider, Henrik},
     TITLE = {Least-squares finite element method for the simulation of
              sea-ice motion},
   JOURNAL = {Comput. Math. Appl.},
  FJOURNAL = {Computers \& Mathematics with Applications. An International
              Journal},
    VOLUME = {172},
      YEAR = {2024},
     PAGES = {38--46},
      ISSN = {0898-1221,1873-7668},
   MRCLASS = {65M60 (74M99 76M10 86A05)},
  MRNUMBER = {4782920},
       DOI = {10.1016/j.camwa.2024.07.023},
       URL = {https://doi.org/10.1016/j.camwa.2024.07.023},
}

@article {GantnerStevenson2024_rates,
    AUTHOR = {Gantner, Gregor and Stevenson, Rob},
     TITLE = {Improved rates for a space-time {FOSLS} of parabolic {PDE}s},
   JOURNAL = {Numer. Math.},
  FJOURNAL = {Numerische Mathematik},
    VOLUME = {156},
      YEAR = {2024},
    NUMBER = {1},
     PAGES = {133--157},
      ISSN = {0029-599X,0945-3245},
   MRCLASS = {65M12 (35F35 65M15 65M50 65M60)},
  MRNUMBER = {4700409},
       DOI = {10.1007/s00211-023-01387-3},
       URL = {https://doi.org/10.1007/s00211-023-01387-3},
}

@article {BalciDieningStorn2023,
    AUTHOR = {Balci, Anna Kh. and Diening, Lars and Storn, Johannes},
     TITLE = {Relaxed {K}a\v canov scheme for the {$p$}-{L}aplacian with
              large exponent},
   JOURNAL = {SIAM J. Numer. Anal.},
  FJOURNAL = {SIAM Journal on Numerical Analysis},
    VOLUME = {61},
      YEAR = {2023},
    NUMBER = {6},
     PAGES = {2775--2794},
      ISSN = {0036-1429,1095-7170},
   MRCLASS = {65N30 (35J70 35J92 65N22)},
  MRNUMBER = {4668389},
MRREVIEWER = {D.\ J.\ Liu},
       DOI = {10.1137/22M1528550},
       URL = {https://doi.org/10.1137/22M1528550},
}

@article {DieningStornTscherpel2023,
    AUTHOR = {Diening, Lars and Storn, Johannes and Tscherpel, Tabea},
     TITLE = {Interpolation operator on negative {S}obolev spaces},
   JOURNAL = {Math. Comp.},
  FJOURNAL = {Mathematics of Computation},
    VOLUME = {92},
      YEAR = {2023},
    NUMBER = {342},
     PAGES = {1511--1541},
      ISSN = {0025-5718,1088-6842},
   MRCLASS = {65D05 (65M15 65N15 65N30)},
  MRNUMBER = {4570332},
       DOI = {10.1090/mcom/3824},
       URL = {https://doi.org/10.1090/mcom/3824},
}

@article {FuhrerHeuerKarkulik2022,
    AUTHOR = {F\"uhrer, Thomas and Heuer, Norbert and Karkulik, Michael},
     TITLE = {M{INRES} for second-order {PDE}s with singular data},
   JOURNAL = {SIAM J. Numer. Anal.},
  FJOURNAL = {SIAM Journal on Numerical Analysis},
    VOLUME = {60},
      YEAR = {2022},
    NUMBER = {3},
     PAGES = {1111--1135},
      ISSN = {0036-1429,1095-7170},
   MRCLASS = {65N30 (65N12)},
  MRNUMBER = {4425909},
MRREVIEWER = {Da\ Xu},
       DOI = {10.1137/21M1457023},
       URL = {https://doi.org/10.1137/21M1457023},
}

@article {LiDemkowicz2022,
    AUTHOR = {Li, Jiaqi and Demkowicz, Leszek},
     TITLE = {An {$L^p$}-{DPG} method with application to 2{D}
              convection-diffusion problems},
   JOURNAL = {Comput. Methods Appl. Math.},
  FJOURNAL = {Computational Methods in Applied Mathematics},
    VOLUME = {22},
      YEAR = {2022},
    NUMBER = {3},
     PAGES = {649--662},
      ISSN = {1609-4840,1609-9389},
   MRCLASS = {65N30 (65N50)},
  MRNUMBER = {4445530},
       DOI = {10.1515/cmam-2021-0158},
       URL = {https://doi.org/10.1515/cmam-2021-0158},
}

@article {HoustonRoggendorfVanDerZee2022,
    AUTHOR = {Houston, Paul and Roggendorf, Sarah and van der Zee,
              Kristoffer G.},
     TITLE = {Gibbs phenomena for {${\rm L}^q$}-best approximation in finite
              element spaces},
   JOURNAL = {ESAIM Math. Model. Numer. Anal.},
  FJOURNAL = {ESAIM. Mathematical Modelling and Numerical Analysis},
    VOLUME = {56},
      YEAR = {2022},
    NUMBER = {1},
     PAGES = {177--211},
      ISSN = {2822-7840,2804-7214},
   MRCLASS = {65N30 (41A50)},
  MRNUMBER = {4376273},
MRREVIEWER = {Gerrit\ Welper},
       DOI = {10.1051/m2an/2021086},
       URL = {https://doi.org/10.1051/m2an/2021086},
}

@article {ErnGudiSmearsVohralik2022,
    AUTHOR = {Ern, Alexandre and Gudi, Thirupathi and Smears, Iain and
              Vohral\'ik, Martin},
     TITLE = {Equivalence of local- and global-best approximations, a simple
              stable local commuting projector, and optimal {$hp$}
              approximation estimates in {$\textbf{H({\rm div})}$}},
   JOURNAL = {IMA J. Numer. Anal.},
  FJOURNAL = {IMA Journal of Numerical Analysis},
    VOLUME = {42},
      YEAR = {2022},
    NUMBER = {2},
     PAGES = {1023--1049},
      ISSN = {0272-4979,1464-3642},
   MRCLASS = {65N30 (65N15)},
  MRNUMBER = {4410735},
       DOI = {10.1093/imanum/draa103},
       URL = {https://doi.org/10.1093/imanum/draa103},
}

@article {GantnerStevenson2021,
    AUTHOR = {Gantner, Gregor and Stevenson, Rob},
     TITLE = {Further results on a space-time {FOSLS} formulation of
              parabolic {PDE}s},
   JOURNAL = {ESAIM Math. Model. Numer. Anal.},
  FJOURNAL = {ESAIM. Mathematical Modelling and Numerical Analysis},
    VOLUME = {55},
      YEAR = {2021},
    NUMBER = {1},
     PAGES = {283--299},
      ISSN = {2822-7840,2804-7214},
   MRCLASS = {65M60 (35K20 65M12 65M15)},
  MRNUMBER = {4216839},
MRREVIEWER = {Lu\ Zhang},
       DOI = {10.1051/m2an/2020084},
       URL = {https://doi.org/10.1051/m2an/2020084},
}

@article {FuhrerKarkulik2021,
    AUTHOR = {F\"uhrer, Thomas and Karkulik, Michael},
     TITLE = {Space-time least-squares finite elements for parabolic
              equations},
   JOURNAL = {Comput. Math. Appl.},
  FJOURNAL = {Computers \& Mathematics with Applications. An International
              Journal},
    VOLUME = {92},
      YEAR = {2021},
     PAGES = {27--36},
      ISSN = {0898-1221,1873-7668},
   MRCLASS = {65M60 (65M12 65M15)},
  MRNUMBER = {4242919},
       DOI = {10.1016/j.camwa.2021.03.004},
       URL = {https://doi.org/10.1016/j.camwa.2021.03.004},
}

@article {HeidPraetoriusWihler2021,
	AUTHOR   = {Heid, Pascal and Praetorius, Dirk and Wihler, Thomas P.},
	TITLE    = {Energy contraction and optimal convergence of adaptive iterative linearized finite element methods},
	JOURNAL  = {Comput. Methods Appl. Math.},
	FJOURNAL = {Computational Methods in Applied Mathematics},
	VOLUME   = {21},
	YEAR     = {2021},
	NUMBER   = {2},
	PAGES    = {407--422},
	ISSN     = {1609-4840},
	MRCLASS  = {65N30 (35J62 47H05 47J25 65N12 65Y20)},
	MRNUMBER = {4235817},
	DOI      = {10.1515/cmam-2021-0025},
	URL      = {https://doi.org/10.1515/cmam-2021-0025},
}

@article {CarstensenMa2021,
    AUTHOR = {Carstensen, Carsten and Ma, Rui},
     TITLE = {Collective marking for arbitrary order adaptive least-squares
              finite element methods with optimal rates},
   JOURNAL = {Comput. Math. Appl.},
  FJOURNAL = {Computers \& Mathematics with Applications. An International
              Journal},
    VOLUME = {95},
      YEAR = {2021},
     PAGES = {271--281},
      ISSN = {0898-1221,1873-7668},
   MRCLASS = {65N30 (65N12 65N50)},
  MRNUMBER = {4271577},
       DOI = {10.1016/j.camwa.2020.12.005},
       URL = {https://doi.org/10.1016/j.camwa.2020.12.005},
}

@article {CarstensenErnPuttkammer2021,
    AUTHOR = {Carstensen, Carsten and Ern, Alexandre and Puttkammer, Sophie},
     TITLE = {Guaranteed lower bounds on eigenvalues of elliptic operators
              with a hybrid high-order method},
   JOURNAL = {Numer. Math.},
  FJOURNAL = {Numerische Mathematik},
    VOLUME = {149},
      YEAR = {2021},
    NUMBER = {2},
     PAGES = {273--304},
      ISSN = {0029-599X,0945-3245},
   MRCLASS = {65N25 (65N30)},
  MRNUMBER = {4332791},
MRREVIEWER = {Xia\ Ji},
       DOI = {10.1007/s00211-021-01228-1},
       URL = {https://doi.org/10.1007/s00211-021-01228-1},
}

@article {Fuhrer2021,
    AUTHOR = {F\"uhrer, Thomas},
     TITLE = {Ultraweak formulation of linear {PDE}s in nondivergence form
              and {DPG} approximation},
   JOURNAL = {Comput. Math. Appl.},
  FJOURNAL = {Computers \& Mathematics with Applications. An International
              Journal},
    VOLUME = {95},
      YEAR = {2021},
     PAGES = {67--84},
      ISSN = {0898-1221,1873-7668},
   MRCLASS = {65N30},
  MRNUMBER = {4271567},
MRREVIEWER = {Yongxiang\ Liu},
       DOI = {10.1016/j.camwa.2020.07.007},
       URL = {https://doi.org/10.1016/j.camwa.2020.07.007},
}

@article{HeidWihler2020_convergence,
    AUTHOR     = {Heid, Pascal and Wihler, Thomas P.},
    TITLE      = {On the convergence of adaptive iterative linearized {G}alerkin methods},
    JOURNAL    = {Calcolo},
    FJOURNAL   = {Calcolo. A Quarterly on Numerical Analysis and Theory of Computation},
    VOLUME     = {57},
    YEAR       = {2020},
    NUMBER     = {3},
    %	PAGES  = {\#24},
    ISSN       = {0008-0624},
    MRCLASS    = {65N30 (35J62 47H05 47H10 47J25 49M15 65N12 65N50)},
    MRNUMBER   = {4131951},
    MRREVIEWER = {Riccardo Sacco},
    DOI        = {10.1007/s10092-020-00368-4},
    URL        = {https://doi.org/10.1007/s10092-020-00368-4},
}

@article {HeidWihler2020_linearization,
	AUTHOR     = {Heid, Pascal and Wihler, Thomas P.},
	TITLE      = {Adaptive iterative linearization {G}alerkin methods for nonlinear problems},
	JOURNAL    = {Math. Comp.},
	FJOURNAL   = {Mathematics of Computation},
	VOLUME     = {89},
	YEAR       = {2020},
	NUMBER     = {326},
	PAGES      = {2707--2734},
	ISSN       = {0025-5718},
	MRCLASS    = {65N30 (47H05 47H10 47J25 49M15 65J15 65N12 65N50)},
	MRNUMBER   = {4136544},
	MRREVIEWER = {Huai Zhang},
	DOI        = {10.1090/mcom/3545},
	URL        = {https://doi.org/10.1090/mcom/3545},
}

@article {MugaVanDerZee2020,
    AUTHOR = {Muga, Ignacio and van der Zee, Kristoffer G.},
     TITLE = {Discretization of linear problems in {B}anach spaces: residual
              minimization, nonlinear {P}etrov-{G}alerkin, and monotone
              mixed methods},
   JOURNAL = {SIAM J. Numer. Anal.},
  FJOURNAL = {SIAM Journal on Numerical Analysis},
    VOLUME = {58},
      YEAR = {2020},
    NUMBER = {6},
     PAGES = {3406--3426},
      ISSN = {0036-1429,1095-7170},
   MRCLASS = {65N30 (41A65 46B20 65J05 65N12 65N15)},
  MRNUMBER = {4178383},
MRREVIEWER = {Andreas\ Petersson},
       DOI = {10.1137/20M1324338},
       URL = {https://doi.org/10.1137/20M1324338},
}

@article {DieningFornasierTomasiWank2020,
    AUTHOR = {Diening, L. and Fornasier, M. and Tomasi, R. and Wank, M.},
     TITLE = {A relaxed {K}a\v canov iteration for the {$p$}-{P}oisson
              problem},
   JOURNAL = {Numer. Math.},
  FJOURNAL = {Numerische Mathematik},
    VOLUME = {145},
      YEAR = {2020},
    NUMBER = {1},
     PAGES = {1--34},
      ISSN = {0029-599X,0945-3245},
   MRCLASS = {65N30 (35J70 35J92 65N12 65N22)},
  MRNUMBER = {4091593},
MRREVIEWER = {Michael\ Neilan},
       DOI = {10.1007/s00211-020-01107-1},
       URL = {https://doi.org/10.1007/s00211-020-01107-1},
}

@article {CaiChenLiuLiu2020,
    AUTHOR = {Cai, Zhiqiang and Chen, Jingshuang and Liu, Min and Liu,
              Xinyu},
     TITLE = {Deep least-squares methods: an unsupervised learning-based
              numerical method for solving elliptic {PDE}s},
   JOURNAL = {J. Comput. Phys.},
  FJOURNAL = {Journal of Computational Physics},
    VOLUME = {420},
      YEAR = {2020},
     PAGES = {109707, 13},
      ISSN = {0021-9991,1090-2716},
   MRCLASS = {65N99 (68T05)},
  MRNUMBER = {4128623},
       DOI = {10.1016/j.jcp.2020.109707},
       URL = {https://doi.org/10.1016/j.jcp.2020.109707},
}

@article {Carstensen2020,
    AUTHOR = {Carstensen, Carsten},
     TITLE = {Collective marking for adaptive least-squares finite element
              methods with optimal rates},
   JOURNAL = {Math. Comp.},
  FJOURNAL = {Mathematics of Computation},
    VOLUME = {89},
      YEAR = {2020},
    NUMBER = {321},
     PAGES = {89--103},
      ISSN = {0025-5718,1088-6842},
   MRCLASS = {65N30 (65N12 65N15)},
  MRNUMBER = {4011536},
MRREVIEWER = {Mohammad\ Asadzadeh},
       DOI = {10.1090/mcom/3474},
       URL = {https://doi.org/10.1090/mcom/3474},
}

@article {QiuZhang2020,
    AUTHOR = {Qiu, Weifeng and Zhang, Shun},
     TITLE = {Adaptive first-order system least-squares finite element
              methods for second-order elliptic equations in nondivergence
              form},
   JOURNAL = {SIAM J. Numer. Anal.},
  FJOURNAL = {SIAM Journal on Numerical Analysis},
    VOLUME = {58},
      YEAR = {2020},
    NUMBER = {6},
     PAGES = {3286--3308},
      ISSN = {0036-1429,1095-7170},
   MRCLASS = {65N30 (65N12 65N15 65N50)},
  MRNUMBER = {4173220},
MRREVIEWER = {Liang\ Ge},
       DOI = {10.1137/19M1271099},
       URL = {https://doi.org/10.1137/19M1271099},
}

@article {RaissiPerdikarisKarniadakis2019,
    AUTHOR = {Raissi, M. and Perdikaris, P. and Karniadakis, G. E.},
     TITLE = {Physics-informed neural networks: a deep learning framework
              for solving forward and inverse problems involving nonlinear
              partial differential equations},
   JOURNAL = {J. Comput. Phys.},
  FJOURNAL = {Journal of Computational Physics},
    VOLUME = {378},
      YEAR = {2019},
     PAGES = {686--707},
      ISSN = {0021-9991,1090-2716},
   MRCLASS = {65M70 (68T05)},
  MRNUMBER = {3881695},
       DOI = {10.1016/j.jcp.2018.10.045},
       URL = {https://doi.org/10.1016/j.jcp.2018.10.045},
}

@article {Westphal2019,
    AUTHOR = {Westphal, Chad R.},
     TITLE = {A {N}ewton div-curl least-squares finite element method for
              the elliptic {M}onge-{A}mp\`ere equation},
   JOURNAL = {Comput. Methods Appl. Math.},
  FJOURNAL = {Computational Methods in Applied Mathematics},
    VOLUME = {19},
      YEAR = {2019},
    NUMBER = {3},
     PAGES = {631--643},
      ISSN = {1609-4840,1609-9389},
   MRCLASS = {65N30 (35J96 65N12)},
  MRNUMBER = {3977490},
MRREVIEWER = {Xiaobing\ Henry\ Feng},
       DOI = {10.1515/cmam-2018-0196},
       URL = {https://doi.org/10.1515/cmam-2018-0196},
}

@article {CarstensenStorn2018,
    AUTHOR = {Carstensen, Carsten and Storn, Johannes},
     TITLE = {Asymptotic exactness of the least-squares finite element
              residual},
   JOURNAL = {SIAM J. Numer. Anal.},
  FJOURNAL = {SIAM Journal on Numerical Analysis},
    VOLUME = {56},
      YEAR = {2018},
    NUMBER = {4},
     PAGES = {2008--2028},
      ISSN = {0036-1429,1095-7170},
   MRCLASS = {65N30 (65N12 65N15)},
  MRNUMBER = {3820383},
MRREVIEWER = {Ignacio\ Romero},
       DOI = {10.1137/17M1125972},
       URL = {https://doi.org/10.1137/17M1125972},
}

@article {CantinHeuer2018,
    AUTHOR = {Cantin, Pierre and Heuer, Norbert},
     TITLE = {A {DPG} framework for strongly monotone operators},
   JOURNAL = {SIAM J. Numer. Anal.},
  FJOURNAL = {SIAM Journal on Numerical Analysis},
    VOLUME = {56},
      YEAR = {2018},
    NUMBER = {5},
     PAGES = {2731--2750},
      ISSN = {0036-1429,1095-7170},
   MRCLASS = {65N30 (47H05 65J15 65N12)},
  MRNUMBER = {3851046},
MRREVIEWER = {Rajen\ Kumar\ Sinha},
       DOI = {10.1137/18M1166663},
       URL = {https://doi.org/10.1137/18M1166663},
}

@article {CarstensenHellwig2018,
    AUTHOR = {Carstensen, Carsten and Hellwig, Friederike},
     TITLE = {Optimal convergence rates for adaptive lowest-order
              discontinuous {P}etrov-{G}alerkin schemes},
   JOURNAL = {SIAM J. Numer. Anal.},
  FJOURNAL = {SIAM Journal on Numerical Analysis},
    VOLUME = {56},
      YEAR = {2018},
    NUMBER = {2},
     PAGES = {1091--1111},
      ISSN = {0036-1429,1095-7170},
   MRCLASS = {65N30 (65N12 65N15)},
  MRNUMBER = {3790080},
MRREVIEWER = {Igor\ Bock},
       DOI = {10.1137/17M1146671},
       URL = {https://doi.org/10.1137/17M1146671},
}

@article {CarstensenParkBringmann2017,
    AUTHOR = {Carstensen, Carsten and Park, Eun-Jae and Bringmann, Philipp},
     TITLE = {Convergence of natural adaptive least squares finite element
              methods},
   JOURNAL = {Numer. Math.},
  FJOURNAL = {Numerische Mathematik},
    VOLUME = {136},
      YEAR = {2017},
    NUMBER = {4},
     PAGES = {1097--1115},
      ISSN = {0029-599X,0945-3245},
   MRCLASS = {65N30 (65N12 65N15 65N50 65Y20)},
  MRNUMBER = {3671598},
MRREVIEWER = {Seydi\ Battal Gazi Karakoc},
       DOI = {10.1007/s00211-017-0866-x},
       URL = {https://doi.org/10.1007/s00211-017-0866-x},
}

@article {CongreveWihler2017,
    AUTHOR = {Congreve, Scott and Wihler, Thomas P.},
     TITLE = {Iterative {G}alerkin discretizations for strongly monotone
              problems},
   JOURNAL = {J. Comput. Appl. Math.},
  FJOURNAL = {Journal of Computational and Applied Mathematics},
    VOLUME = {311},
      YEAR = {2017},
     PAGES = {457--472},
      ISSN = {0377-0427,1879-1778},
   MRCLASS = {65N30 (65N50)},
  MRNUMBER = {3552717},
MRREVIEWER = {C.\ Ilioi},
       DOI = {10.1016/j.cam.2016.08.014},
       URL = {https://doi.org/10.1016/j.cam.2016.08.014},
}

@article {Gallistl2017,
    AUTHOR = {Gallistl, Dietmar},
     TITLE = {Variational formulation and numerical analysis of linear
              elliptic equations in nondivergence form with {C}ordes
              coefficients},
   JOURNAL = {SIAM J. Numer. Anal.},
  FJOURNAL = {SIAM Journal on Numerical Analysis},
    VOLUME = {55},
      YEAR = {2017},
    NUMBER = {2},
     PAGES = {737--757},
      ISSN = {0036-1429,1095-7170},
   MRCLASS = {65N30 (31B30 35J30 65N12 65N15 65N50)},
  MRNUMBER = {3628316},
MRREVIEWER = {Srinivasan\ Kesavan},
       DOI = {10.1137/16M1080495},
       URL = {https://doi.org/10.1137/16M1080495},
}

@article {MonneslandLeeGunzburgerYoon2016,
    AUTHOR = {Monnesland, Irene Sonja and Lee, Eunjung and Gunzburger, Max
              and Yoon, Ryeongkyung},
     TITLE = {A least-squares finite element method for a nonlinear {S}tokes
              problem in glaciology},
   JOURNAL = {Comput. Math. Appl.},
  FJOURNAL = {Computers \& Mathematics with Applications. An International
              Journal},
    VOLUME = {71},
      YEAR = {2016},
    NUMBER = {11},
     PAGES = {2421--2431},
      ISSN = {0898-1221,1873-7668},
   MRCLASS = {65N30 (76D07 86A40)},
  MRNUMBER = {3501329},
MRREVIEWER = {H.\ P.\ Dikshit},
       DOI = {10.1016/j.camwa.2015.11.001},
       URL = {https://doi.org/10.1016/j.camwa.2015.11.001},
}

@article {CarstensenPark2015,
    AUTHOR = {Carstensen, Carsten and Park, Eun-Jae},
     TITLE = {Convergence and optimality of adaptive least squares finite
              element methods},
   JOURNAL = {SIAM J. Numer. Anal.},
  FJOURNAL = {SIAM Journal on Numerical Analysis},
    VOLUME = {53},
      YEAR = {2015},
    NUMBER = {1},
     PAGES = {43--62},
      ISSN = {0036-1429,1095-7170},
   MRCLASS = {65N30 (65N12 65N15 65N50 65Y20)},
  MRNUMBER = {3296614},
MRREVIEWER = {J\'an\ Lov\'i\v sek},
       DOI = {10.1137/130949634},
       URL = {https://doi.org/10.1137/130949634},
}

@article {MullerStarkeSchwarzSchroder2014,
    AUTHOR = {M\"uller, Benjamin and Starke, Gerhard and Schwarz, Alexander
              and Schr\"oder, J\"org},
     TITLE = {A first-order system least squares method for hyperelasticity},
   JOURNAL = {SIAM J. Sci. Comput.},
  FJOURNAL = {SIAM Journal on Scientific Computing},
    VOLUME = {36},
      YEAR = {2014},
    NUMBER = {5},
     PAGES = {B795--B816},
      ISSN = {1064-8275,1095-7197},
   MRCLASS = {65N30 (74G15 74S05)},
  MRNUMBER = {3264569},
MRREVIEWER = {Yves\ Renard},
       DOI = {10.1137/130937573},
       URL = {https://doi.org/10.1137/130937573},
}

@article {CarstensenGedicke2014,
    AUTHOR = {Carstensen, Carsten and Gedicke, Joscha},
     TITLE = {Guaranteed lower bounds for eigenvalues},
   JOURNAL = {Math. Comp.},
  FJOURNAL = {Mathematics of Computation},
    VOLUME = {83},
      YEAR = {2014},
    NUMBER = {290},
     PAGES = {2605--2629},
      ISSN = {0025-5718,1088-6842},
   MRCLASS = {65N25 (15A42 65N15 65N30 65N50)},
  MRNUMBER = {3246802},
MRREVIEWER = {J\'{a}n\ Lov\'{\i}\v{s}ek},
       DOI = {10.1090/S0025-5718-2014-02833-0},
       URL = {https://doi.org/10.1090/S0025-5718-2014-02833-0},
}

@article {BelenkiDieningKreuzer2012,
    AUTHOR = {Belenki, Liudmila and Diening, Lars and Kreuzer, Christian},
     TITLE = {Optimality of an adaptive finite element method for the
              {$p$}-{L}aplacian equation},
   JOURNAL = {IMA J. Numer. Anal.},
  FJOURNAL = {IMA Journal of Numerical Analysis},
    VOLUME = {32},
      YEAR = {2012},
    NUMBER = {2},
     PAGES = {484--510},
      ISSN = {0272-4979,1464-3642},
   MRCLASS = {65N30 (65N15)},
  MRNUMBER = {2911397},
MRREVIEWER = {Snorre\ H.\ Christiansen},
       DOI = {10.1093/imanum/drr016},
       URL = {https://doi.org/10.1093/imanum/drr016},
}

@article {PayetteReddy2011,
    AUTHOR = {Payette, G. S. and Reddy, J. N.},
     TITLE = {On the roles of minimization and linearization in
              least-squares finite element models of nonlinear
              boundary-value problems},
   JOURNAL = {J. Comput. Phys.},
  FJOURNAL = {Journal of Computational Physics},
    VOLUME = {230},
      YEAR = {2011},
    NUMBER = {9},
     PAGES = {3589--3613},
      ISSN = {0021-9991,1090-2716},
   MRCLASS = {76M10},
  MRNUMBER = {2780480},
       DOI = {10.1016/j.jcp.2011.02.002},
       URL = {https://doi.org/10.1016/j.jcp.2011.02.002},
}

@article {Siebert2011,
    AUTHOR = {Siebert, Kunibert G.},
     TITLE = {A convergence proof for adaptive finite elements without lower
              bound},
   JOURNAL = {IMA J. Numer. Anal.},
  FJOURNAL = {IMA Journal of Numerical Analysis},
    VOLUME = {31},
      YEAR = {2011},
    NUMBER = {3},
     PAGES = {947--970},
      ISSN = {0272-4979,1464-3642},
   MRCLASS = {65N30 (65N12 65N50)},
  MRNUMBER = {2832786},
MRREVIEWER = {Christos\ A.\ Xenophontos},
       DOI = {10.1093/imanum/drq001},
       URL = {https://doi.org/10.1093/imanum/drq001},
}

@article {CianchiMazya2011,
    AUTHOR = {Cianchi, Andrea and Maz'ya, Vladimir G.},
     TITLE = {Global {L}ipschitz regularity for a class of quasilinear
              elliptic equations},
   JOURNAL = {Comm. Partial Differential Equations},
  FJOURNAL = {Communications in Partial Differential Equations},
    VOLUME = {36},
      YEAR = {2011},
    NUMBER = {1},
     PAGES = {100--133},
      ISSN = {0360-5302,1532-4133},
   MRCLASS = {35J62 (35B65 35D30 35J25)},
  MRNUMBER = {2763349},
MRREVIEWER = {Jens\ Habermann},
       DOI = {10.1080/03605301003657843},
       URL = {https://doi.org/10.1080/03605301003657843},
}

@article {ManteuffelMcCormickSchmidtWestphal2006,
    AUTHOR = {Manteuffel, T. A. and McCormick, S. F. and Schmidt, J. G. and
              Westphal, C. R.},
     TITLE = {First-order system least squares for geometrically nonlinear
              elasticity},
   JOURNAL = {SIAM J. Numer. Anal.},
  FJOURNAL = {SIAM Journal on Numerical Analysis},
    VOLUME = {44},
      YEAR = {2006},
    NUMBER = {5},
     PAGES = {2057--2081},
      ISSN = {0036-1429,1095-7170},
   MRCLASS = {74G15 (65N30 74B20 74S05)},
  MRNUMBER = {2263040},
MRREVIEWER = {Penny\ J.\ Davies},
       DOI = {10.1137/050628027},
       URL = {https://doi.org/10.1137/050628027},
}

@article {Guermond2004,
    AUTHOR = {Guermond, J. L.},
     TITLE = {A finite element technique for solving first-order {PDE}s in
              {$L^P$}},
   JOURNAL = {SIAM J. Numer. Anal.},
  FJOURNAL = {SIAM Journal on Numerical Analysis},
    VOLUME = {42},
      YEAR = {2004},
    NUMBER = {2},
     PAGES = {714--737},
      ISSN = {0036-1429,1095-7170},
   MRCLASS = {65N30},
  MRNUMBER = {2084233},
MRREVIEWER = {Bruce\ A.\ Finlayson},
       DOI = {10.1137/S0036142902417054},
       URL = {https://doi.org/10.1137/S0036142902417054},
}

@article {BanschMorinNochetto2002,
    AUTHOR = {B\"ansch, Eberhard and Morin, Pedro and Nochetto, Ricardo H.},
     TITLE = {An adaptive {U}zawa {FEM} for the {S}tokes problem:
              convergence without the inf-sup condition},
   JOURNAL = {SIAM J. Numer. Anal.},
  FJOURNAL = {SIAM Journal on Numerical Analysis},
    VOLUME = {40},
      YEAR = {2002},
    NUMBER = {4},
     PAGES = {1207--1229},
      ISSN = {0036-1429,1095-7170},
   MRCLASS = {65N30 (65N50 76D07 76M10)},
  MRNUMBER = {1951892},
MRREVIEWER = {Murli\ M.\ Gupta},
       DOI = {10.1137/S0036142901392134},
       URL = {https://doi.org/10.1137/S0036142901392134},
}

@article {BochevCaiManteuffelMcCormick1998,
    AUTHOR = {Bochev, P. and Cai, Z. and Manteuffel, T. A. and McCormick, S.
              F.},
     TITLE = {Analysis of velocity-flux first-order system least-squares
              principles for the {N}avier-{S}tokes equations. {I}},
   JOURNAL = {SIAM J. Numer. Anal.},
  FJOURNAL = {SIAM Journal on Numerical Analysis},
    VOLUME = {35},
      YEAR = {1998},
    NUMBER = {3},
     PAGES = {990--1009},
      ISSN = {0036-1429,1095-7170},
   MRCLASS = {76M25 (65N12 65N30 76D05)},
  MRNUMBER = {1619922},
       DOI = {10.1137/S0036142996313592},
       URL = {https://doi.org/10.1137/S0036142996313592},
}

@article {Dorfler1996,
    AUTHOR = {D\"orfler, Willy},
     TITLE = {A convergent adaptive algorithm for {P}oisson's equation},
   JOURNAL = {SIAM J. Numer. Anal.},
  FJOURNAL = {SIAM Journal on Numerical Analysis},
    VOLUME = {33},
      YEAR = {1996},
    NUMBER = {3},
     PAGES = {1106--1124},
      ISSN = {0036-1429},
   MRCLASS = {65N50 (65N55)},
  MRNUMBER = {1393904},
MRREVIEWER = {S.\ F.\ McCormick},
       DOI = {10.1137/0733054},
       URL = {https://doi.org/10.1137/0733054},
}

@article {CarstensenStephan1995,
    AUTHOR = {Carstensen, Carsten and Stephan, Ernst P.},
     TITLE = {Adaptive coupling of boundary elements and finite elements},
   JOURNAL = {RAIRO Mod\'el. Math. Anal. Num\'er.},
  FJOURNAL = {RAIRO Mod\'elisation Math\'ematique et Analyse Num\'erique},
    VOLUME = {29},
      YEAR = {1995},
    NUMBER = {7},
     PAGES = {779--817},
      ISSN = {0764-583X},
   MRCLASS = {65N35 (65D07 65N30 65R20)},
  MRNUMBER = {1364401},
MRREVIEWER = {Stephen\ W.\ Brady},
       DOI = {10.1051/m2an/1995290707791},
       URL = {https://doi.org/10.1051/m2an/1995290707791},
}

@article {Park1995,
    AUTHOR = {Park, Eun-Jae},
     TITLE = {Mixed finite element methods for nonlinear second-order
              elliptic problems},
   JOURNAL = {SIAM J. Numer. Anal.},
  FJOURNAL = {SIAM Journal on Numerical Analysis},
    VOLUME = {32},
      YEAR = {1995},
    NUMBER = {3},
     PAGES = {865--885},
      ISSN = {0036-1429},
   MRCLASS = {65N30 (35J60 65N15)},
  MRNUMBER = {1335659},
MRREVIEWER = {Jeff\ S.\ McGough},
       DOI = {10.1137/0732040},
       URL = {https://doi.org/10.1137/0732040},
}

@article {BochevGunzburger1993,
    AUTHOR = {Bochev, Pavel B. and Gunzburger, Max D.},
     TITLE = {A least-squares finite element method for the
              {N}avier-{S}tokes equations},
   JOURNAL = {Appl. Math. Lett.},
  FJOURNAL = {Applied Mathematics Letters. An International Journal of Rapid
              Publication},
    VOLUME = {6},
      YEAR = {1993},
    NUMBER = {2},
     PAGES = {27--30},
      ISSN = {0893-9659,1873-5452},
   MRCLASS = {76M10 (65N30 76D05)},
  MRNUMBER = {1347770},
       DOI = {10.1016/0893-9659(93)90007-A},
       URL = {https://doi.org/10.1016/0893-9659(93)90007-A},
}

@incollection {DouglasPaesLemeGiorgi1993,
    AUTHOR = {Douglas, Jr., Jim and Paes-Leme, Paulo Jorge and Giorgi,
              Tiziana},
     TITLE = {Generalized {F}orchheimer flow in porous media},
 BOOKTITLE = {Boundary value problems for partial differential equations and
              applications},
    SERIES = {RMA Res. Notes Appl. Math.},
    VOLUME = {29},
     PAGES = {99--111},
 PUBLISHER = {Masson, Paris},
      YEAR = {1993},
      ISBN = {2-225-84334-1},
   MRCLASS = {76S05 (35Q35 76M10)},
  MRNUMBER = {1260441},
       DOI = {10.1007/978-3-0348-8564-5},
       URL = {https://doi.org/10.1007/978-3-0348-8564-5},
}

@article {Zarantonello1960,
    AUTHOR = {Zarantonello, E.H.},
     TITLE = {Solving functional equations by contractive averaging},
   JOURNAL = {Technical Report},
    VOLUME = {160},
      YEAR = {1960},
    NUMBER = {},
     PAGES = {},
      NOTE = {Mathematics Research Center, Univ. of Wisconsin, Madison}
}

@book {BoffiBrezziFortin2013,
    AUTHOR = {Boffi, Daniele and Brezzi, Franco and Fortin, Michel},
     TITLE = {Mixed finite element methods and applications},
    SERIES = {Springer Series in Computational Mathematics},
    VOLUME = {44},
 PUBLISHER = {Springer, Heidelberg},
      YEAR = {2013},
     PAGES = {xiv+685},
      ISBN = {978-3-642-36518-8; 978-3-642-36519-5},
   MRCLASS = {65-02 (65M60 65N30)},
  MRNUMBER = {3097958},
MRREVIEWER = {Beny\ Neta},
       DOI = {10.1007/978-3-642-36519-5},
       URL = {https://doi.org/10.1007/978-3-642-36519-5},
}

@book {BochevGunzburger2009,
    AUTHOR = {Bochev, Pavel B. and Gunzburger, Max D.},
     TITLE = {Least-squares finite element methods},
    SERIES = {Applied Mathematical Sciences},
    VOLUME = {166},
 PUBLISHER = {Springer, New York},
      YEAR = {2009},
     PAGES = {xxii+660},
      ISBN = {978-0-387-30888-3},
   MRCLASS = {65-02 (35A35 65M60 65N30 74S05 76M10)},
  MRNUMBER = {2490235},
MRREVIEWER = {Tsu-Fen\ Chen},
       DOI = {10.1007/b13382},
       URL = {https://doi.org/10.1007/b13382},
}

@article {Stevenson2008,
    AUTHOR = {Stevenson, Rob},
     TITLE = {The completion of locally refined simplicial partitions
              created by bisection},
   JOURNAL = {Math. Comp.},
  FJOURNAL = {Mathematics of Computation},
    VOLUME = {77},
      YEAR = {2008},
    NUMBER = {261},
     PAGES = {227--241},
      ISSN = {0025-5718,1088-6842},
   MRCLASS = {65N50},
  MRNUMBER = {2353951},
       DOI = {10.1090/S0025-5718-07-01959-X},
       URL = {https://doi.org/10.1090/S0025-5718-07-01959-X},
}

@article {Traxler1997,
    AUTHOR = {Traxler, C. T.},
     TITLE = {An algorithm for adaptive mesh refinement in {$n$} dimensions},
   JOURNAL = {Computing},
  FJOURNAL = {Computing. Archives for Scientific Computing},
    VOLUME = {59},
      YEAR = {1997},
    NUMBER = {2},
     PAGES = {115--137},
      ISSN = {0010-485X,1436-5057},
   MRCLASS = {65N50},
  MRNUMBER = {1475530},
       DOI = {10.1007/BF02684475},
       URL = {https://doi.org/10.1007/BF02684475},
}

@article {Maubach1995,
    AUTHOR = {Maubach, Joseph M.},
     TITLE = {Local bisection refinement for {$n$}-simplicial grids
              generated by reflection},
   JOURNAL = {SIAM J. Sci. Comput.},
  FJOURNAL = {SIAM Journal on Scientific Computing},
    VOLUME = {16},
      YEAR = {1995},
    NUMBER = {1},
     PAGES = {210--227},
      ISSN = {1064-8275},
   MRCLASS = {65M50},
  MRNUMBER = {1311687},
MRREVIEWER = {Patrick\ M.\ Knupp},
       DOI = {10.1137/0916014},
       URL = {https://doi.org/10.1137/0916014},
}

@article {PehlivanovCareyLazarov1994,
    AUTHOR = {Pehlivanov, A. I. and Carey, G. F. and Lazarov, R. D.},
     TITLE = {Least-squares mixed finite elements for second-order elliptic
              problems},
   JOURNAL = {SIAM J. Numer. Anal.},
  FJOURNAL = {SIAM Journal on Numerical Analysis},
    VOLUME = {31},
      YEAR = {1994},
    NUMBER = {5},
     PAGES = {1368--1377},
      ISSN = {0036-1429},
   MRCLASS = {65N30},
  MRNUMBER = {1293520},
MRREVIEWER = {Lutz\ Angermann},
       DOI = {10.1137/0731071},
       URL = {https://doi.org/10.1137/0731071},
}

@book {Zeidler1990,
    AUTHOR = {Zeidler, Eberhard},
     TITLE = {Nonlinear functional analysis and its applications. {II}/{B}},
 PUBLISHER = {Springer-Verlag, New York},
      YEAR = {1990},
     PAGES = {i--xvi and 469--1202},
      ISBN = {0-387-97167-X},
   MRCLASS = {47-02 (35-01 35J60 47Hxx 58-01 65Jxx)},
  MRNUMBER = {1033498},
MRREVIEWER = {Jean\ Mawhin},
       DOI = {10.1007/978-1-4612-0985-0},
       URL = {https://doi.org/10.1007/978-1-4612-0985-0},
}

@book {GiraulRaviart1986,
    AUTHOR = {Girault, Vivette and Raviart, Pierre-Arnaud},
     TITLE = {Finite element methods for {N}avier-{S}tokes equations},
    SERIES = {Springer Series in Computational Mathematics},
    VOLUME = {5},
      NOTE = {Theory and algorithms},
 PUBLISHER = {Springer-Verlag, Berlin},
      YEAR = {1986},
     PAGES = {x+374},
      ISBN = {3-540-15796-4},
   MRCLASS = {65N30 (65-02 76-08)},
  MRNUMBER = {851383 (88b:65129)},
MRREVIEWER = {Max D. Gunzburger},
}

@article {Jespersen1977,
    AUTHOR = {Jespersen, Dennis C.},
     TITLE = {A least squares decomposition method for solving elliptic
              equations},
   JOURNAL = {Math. Comp.},
  FJOURNAL = {Mathematics of Computation},
    VOLUME = {31},
      YEAR = {1977},
    NUMBER = {140},
     PAGES = {873--880},
      ISSN = {0025-5718,1088-6842},
   MRCLASS = {65N30},
  MRNUMBER = {461948},
MRREVIEWER = {Joachim\ A.\ Nitsche},
       DOI = {10.2307/2006118},
       URL = {https://doi.org/10.2307/2006118},
}

@article {ghps2021,
    AUTHOR = {Gantner, Gregor and Haberl, Alexander and Praetorius, Dirk and
              Schimanko, Stefan},
     TITLE = {Rate optimality of adaptive finite element methods with
              respect to overall computational costs},
   JOURNAL = {Math. Comp.},
  FJOURNAL = {Mathematics of Computation},
    VOLUME = {90},
      YEAR = {2021},
    NUMBER = {331},
     PAGES = {2011--2040},
      ISSN = {0025-5718},
   MRCLASS = {65N30 (65N22 65N50 65Y20)},
  MRNUMBER = {4280291},
       DOI = {10.1090/mcom/3654},
       URL = {https://doi.org/10.1090/mcom/3654},
}

@article{kpp2013,
	Author = {Karkulik, Michael and Pavlicek, David and Praetorius, Dirk},
	Doi = {10.1007/s00365-013-9192-4},
	Fjournal = {Constructive Approximation. An International Journal for Approximations and Expansions},
	Issn = {0176-4276},
	Journal = {Constr. Approx.},
	Mrclass = {65N50 (65N30 65Y20)},
	Mrnumber = {3097045},
	Number = {2},
	Pages = {213--234},
	Title = {On 2{D} newest vertex bisection: optimality of mesh-closure and {$H^1$}-stability of {$L_2$}-projection},
	Url = {https://doi.org/10.1007/s00365-013-9192-4},
	Volume = {38},
	Year = {2013}}

@article {affkp2013,
	AUTHOR = {Aurada, Markus and Feischl, Michael and F\"{u}hrer, Thomas and
	Karkulik, Michael and Praetorius, Dirk},
	TITLE = {Efficiency and optimality of some weighted-residual error
	estimator for adaptive 2{D} boundary element methods},
	JOURNAL = {Comput. Methods Appl. Math.},
	FJOURNAL = {Computational Methods in Applied Mathematics},
	VOLUME = {13},
	YEAR = {2013},
	NUMBER = {3},
	PAGES = {305--332},
	ISSN = {1609-4840},
	MRCLASS = {65N38 (65N15 65N50)},
	MRNUMBER = {3094620},
	MRREVIEWER = {Paul Andrew Martin},
	DOI = {10.1515/cmam-2013-0010},
	URL = {https://doi.org/10.1515/cmam-2013-0010},
}

@article{gmz2012,
	Author = {Garau, Eduardo M. and Morin, Pedro and Zuppa, Carlos},
	Doi = {10.4208/nmtma.2012.m1023},
	Fjournal = {Numerical Mathematics. Theory, Methods and Applications},
	Issn = {1004-8979},
	Journal = {Numer. Math. Theory Methods Appl.},
	Mrclass = {65N30 (35J25 35J62 65N12 65N50)},
	Mrnumber = {2911871},
	Mrreviewer = {Meng Zhao Qin},
	Number = {2},
	Pages = {131--156},
	Title = {Quasi-optimal convergence rate of an {AFEM} for quasi-linear problems of monotone type},
	Url = {https://doi.org/10.4208/nmtma.2012.m1023},
	Volume = {5},
	Year = {2012},
%	Bdsk-Url-1 = {https://doi.org/10.4208/nmtma.2012.m1023}
}
